\documentclass[spanningrule]{cambridge7A}

\usepackage{amssymb}
\usepackage{amsmath}
\usepackage{amsthm}

\usepackage[numbers]{natbib}
\usepackage[figuresright]{rotating}
\usepackage[refpage]{nomencl}  

\usepackage{floatpag}
\rotfloatpagestyle{empty}
\usepackage{graphicx}
\usepackage{txfonts}
\usepackage{caption}  
\usepackage{subcaption}  
\usepackage{placeins}  

\usepackage[within=section]{newfloat}  

\def\exmode{0} 

\usepackage{xr}
\usepackage{amsxtra}  
\usepackage{bm}   

\makeatletter
\let\old@font@info\@font@info
\def\@font@info#1{%
\expandafter\ifx\csname\detokenize{#1}\endcsname\relax
  \old@font@info{#1}%
\fi
\expandafter\xdef\csname\detokenize{#1}\endcsname{}%
}
\makeatother

\usepackage[shortlabels]{enumitem}
\setlist[enumerate]{leftmargin=*,left=0pt .. 3\parindent,align=left}

\newlist{holgerenum}{enumerate}{3}
\setlist[holgerenum]{label=(\alph*),ref=(\alph*),left=0pt .. 2\parindent,align=left,itemsep=0.5em,leftmargin=0pt,itemindent=1.7\parindent,topsep=0.5em}

\newlist{holgerenum2}{enumerate}{3}  
\setlist[holgerenum2]{label=(\alph*),ref=(\alph*),left=0pt .. 2\parindent,align=left,itemsep=0.5em,leftmargin=0pt,itemindent=1.7\parindent,topsep=0.5em}

\makeatletter
\newcommand{\advanceenumeratelevel}{\advance\@enumdepth\@ne}
\makeatother

\makeatletter                           
\providecommand*{\toclevel@schapter}{0}
\makeatother                            

\swapnumbers
  \theoremstyle{plain}
  \newtheorem{theorem}{Theorem}[section]
  \newtheorem{thm}[theorem]{Theorem}
  
  \newtheorem*{lem*}{Lemma}
    \newtheorem*{claim*}{Claim}
    \newtheorem*{fact*}{fact}
    \newtheorem{lem}[theorem]{Lemma}
  \newtheorem{prop}[theorem]{Proposition}
  
   \newtheorem{cor}[theorem]{Corollary}

 \newtheorem{COAL-norm}[theorem]{Norm Equivalence Theorem}
 \newtheorem{CUNOS-coord}[theorem]{Jacobson Co-ordinatization Theorem}
 \newtheorem{Groups-desc}[theorem]{Descent Theorem}

  \newtheorem*{theorem*}{Theorem}
  \newtheorem*{lemma*}{Lemma}
  \newtheorem*{proposition*}{Proposition}
  \newtheorem*{corollary*}{Corollary}
  \newtheorem*{conjecture*}{Conjecture}

  \theoremstyle{definition}
  \newtheorem{definition}[theorem]{Definition}
  \newtheorem{example}[theorem]{Example}

  \newtheorem*{definition*}{Definition}
  \newtheorem*{example*}{Example}
  \newtheorem*{notation*}{Notation}
  \newtheorem*{exer*}{Exercise}
  
  \theoremstyle{remark}

\newtheorem{rmk*}[theorem]{Remark}  
      \newtheorem{remark*}[theorem]{Remark}

\newtheoremstyle{exercise}
  {0.8\baselineskip}
  {3pt}
  {\small}
  {}
  {\bf\small}
  {.}
  {.5em}
   {}     
  {}

\theoremstyle{exercise}

\newenvironment{sol}[1]{\medskip\noindent\textbf{\ref{#1}} \setcounter{equation}{0}}{}

\counterwithout*{section}{chapter}  

\numberwithin{equation}{theorem} 
\renewcommand{\thechapter}{\Roman{chapter}}

\makeatletter
\newenvironment{ermk}%
{\smallskip \noindent{\emph{Remark}.} }%
{\smallskip \global\@ignoretrue}
\makeatother

\makeatletter
{\smallskip \noindent{\emph{Remarks}.} }%
{\smallskip \global\@ignoretrue}
\makeatother

\makeatletter
\newenvironment{borel}[1]%
{\medskip \refstepcounter{theorem}\noindent{\bf \thetheorem\ }{{\bf #1.}}}%
{\smallskip \global\@ignoretrue}
\makeatother

\makeatletter
{\medskip \refstepcounter{theorem}\noindent{\bf \thetheorem\ }{{\bf #1} {(#2)}}\textbf{.}\ }%
{\smallskip \global\@ignoretrue}
\makeatother

\makeatletter
\newenvironment{borel*}%
{\medskip \refstepcounter{theorem}\noindent{\bf \thetheorem.}}%
{\smallskip \global\@ignoretrue}
\makeatother

\newcommand{\chapspace}{\vspace{-75pt}}

\DeclareFloatingEnvironment[
    name={Table of Identities},
    placement=hbtp,
    within=section
    ]{tableids}
\ifnum\exmode=0
	\newcommand{\solnchap}[1]{\chapter*{#1}\markboth{#1}{#1}}
	\newcommand{\solnsec}[1]{\addcontentsline{toc}{section}{#1}\section*{Solutions for #1}\markright{#1}\setcounter{equation}{0}}
\else
	\newcommand{\solnchap}[1]{\vspace{-100pt}\section{{#1}}}
	\newcommand{\solnsec}[1]{\begin{center} \textbf{{#1}}\\ \end{center}}
\fi

\renewcommand{\thesection}{\arabic{section}}

\renewcommand{\indexname}{Subject index}

\usepackage{makeidx}         
\usepackage{graphicx}        
\usepackage{multicol}        

\renewcommand{\theequation}{\arabic{equation}}
\renewcommand{\thechapter}{\Roman{chapter}}
\numberwithin{figure}{section}
\renewcommand{\thefigure}{\thesection\alph{figure}}

\ifnum\exmode=1  
	\usepackage[notcite,notref]{showkeys}   
\fi

\usepackage{varioref}
\usepackage{dynkin-diagrams}
\usepackage{booktabs}
\usepackage{tikz}
\usetikzlibrary{chains}
\usepackage[all]{xy}

\renewcommand{\:}{\colon}

\newcommand{\IA}{\mathbb{A}}
\newcommand{\IC}{\mathbb{C}}
\newcommand{\ID}{\mathbb{D}} 
\newcommand{\IF}{\mathbb{F}}
\newcommand{\IH}{\mathbb{H}}
\newcommand{\IN}{\mathbb{N}}
\newcommand{\IO}{\mathbb{O}}

\newcommand{\IQ}{\mathbb{Q}}
\newcommand{\IR}{\mathbb{R}}
\newcommand{\IS}{\mathbb{S}} 
\newcommand{\IZ}{\mathbb{Z}}
\newcommand{\vart}{\vartheta}
\newcommand{\vep}{\varepsilon}  
\newcommand{\vph}{\varphi}
\newcommand{\vrh}{\varrho}
\newcommand{\bu}{\bullet}
\newcommand{\la}{\langle}
\newcommand{\ra}{\rangle}
\newcommand{\raq}{\rangle_{\mathrm{quad}}}

\newcommand{\ras}{\rangle_{\mathrm{sesq}}}
\newcommand{\pt}{\,.\,}
\newcommand{\Eins}{\mathrm{\mathbf{1}}}

\newcommand{\mcC}{\mathcal{C}}
\newcommand{\mcD}{\mathcal{D}}

\newcommand{\mcF}{\mathcal{F}}

\newcommand{\msA}{\mathscr{A}}
\newcommand{\msB}{\mathscr{B}}
\newcommand{\msO}{\mathscr{O}}

\newcommand{\pmsB}{\mathscr{B}^\prime}
\newcommand{\hac}{\hat{C}}
\newcommand{\haj}{\hat{J}}
\newcommand{\han}{\hat{N}}
\newcommand{\bfa}{\mathbf{a}}
\newcommand{\bfe}{\mathbf{e}}
\newcommand{\bfh}{\mathbf{h}}
\newcommand{\bfi}{\mathbf{i}}
\newcommand{\bfj}{\mathbf{j}}
\newcommand{\bfk}{\mathbf{k}}

\newcommand{\bfp}{\mathbf{p}}
\newcommand{\bfq}{\mathbf{q}}
\newcommand{\bfs}{\mathbf{s}}
\newcommand{\bft}{\mathbf{t}}
\newcommand{\bfu}{\mathbf{u}}
\newcommand{\bfx}{\mathbf{x}}

\newcommand{\bfG}{\mathbf{G}}
\newcommand{\bfH}{\mathbf{H}}

\newcommand{\bfM}{\mathbf{M}}

\newcommand{\bfO}{\mathbf{O}}
\newcommand{\bfS}{\mathbf{S}}
\newcommand{\bfT}{\mathbf{T}}

\newcommand{\bfX}{\mathbf{X}}

\newcommand{\bfGL}{\mathbf{GL}}
\newcommand{\bfGm}{\mathbf{G}_{\mathbf{m}}}
\newcommand{\bfOr}{\bfO}
\newcommand{\bfmu}{\bm{\mu}}
\newcommand{\set}{\mathbf{set}}
\newcommand{\bfAut}{\mathbf{Aut}}
\newcommand{\bfEld}{\mathbf{Elid}}
\newcommand{\bfElf}{\mathbf{Elfr}}
\newcommand{\bfHer}{\mathbf{Her}}
\newcommand{\bfHyp}{\mathbf{Hyp}}
\newcommand{\bfIsom}{\mathbf{Isom}}

\newcommand{\bfSplid}{\mathbf{Splid}}
\newcommand{\bfaist}{\mathop{\mathbf{altist}_k}\nolimits}
\newcommand{\bfpainv}{\mathop{\mathbf{paltinv}_k}\nolimits}
\newcommand{\hgamma}{\hat{\gamma}}
\newcommand{\heta}{\hat{\eta}}
\newcommand{\hiota}{\hat{\iota}}
\newcommand{\htau}{\hat{\tau}}
\newcommand{\hp}{\hat{p}}
\newcommand{\hC}{\hat{C}}
\newcommand{\hmsC}{\hat{\mathscr{C}}}

\newcommand{\hDelta}{\hat{\Delta}}
\newcommand{\hGamma}{\hat{\Gamma}}
\newcommand{\mfa}{\mathfrak{a}}
\newcommand{\mfb}{\mathfrak{b}}

\newcommand{\mfd}{\mathfrak{d}}

\newcommand{\mfg}{\mathfrak{g}}

\newcommand{\mfm}{\mathfrak{m}}
\newcommand{\mfn}{\mathfrak{n}}
\newcommand{\mfo}{\mathfrak{o}}
\newcommand{\mfp}{\mathfrak{p}}
\newcommand{\mfq}{\mathfrak{q}}

\newcommand{\mfD}{\mathfrak{D}}
\newcommand{\mfE}{\mathfrak{E}}

\newcommand{\mfP}{\mathfrak{P}}
\newcommand{\mfS}{\mathfrak{S}}

\newcommand{\pA}{A^\prime}
\newcommand{\pB}{B^\prime}
\newcommand{\pC}{C^\prime}
\newcommand{\pD}{D^\prime}
\newcommand{\pE}{E^\prime}
\newcommand{\pF}{F^\prime}
\newcommand{\pH}{H^\prime}
\newcommand{\pJ}{J^\prime}
\newcommand{\pK}{K^\prime}
\newcommand{\pL}{L^\prime}
\newcommand{\pM}{M^\prime}
\newcommand{\pN}{N^\prime}
\newcommand{\pQ}{Q^\prime}
\newcommand{\pR}{R^\prime}
\newcommand{\pS}{S^\prime}
\newcommand{\pT}{T^\prime}
\newcommand{\pX}{X^\prime}
\newcommand{\pY}{Y^\prime}
\usepackage[mathscr]{eucal}
\newcommand{\pmfp}{\mfp^\prime}
\newcommand{\ppJ}{J^{\prime\prime}}

\newcommand{\ppL}{L^{\prime\prime}}
\newcommand{\ppM}{M^{\prime\prime}}

\newcommand{\ppN}{N^{\prime\prime}}
\newcommand{\ppR}{R^{\prime\prime}}
\newcommand{\ppY}{Y^{\prime\prime}}
\newcommand{\pc}{c^\prime}
\newcommand{\pd}{d^\prime}
\newcommand{\pe}{e^\prime}
\newcommand{\pf}{f^\prime}
\newcommand{\pg}{g^\prime}
\newcommand{\pj}{j^\prime}
\newcommand{\pk}{k^\prime}
\newcommand{\ppk}{k^{\prime\prime}}

\newcommand{\pp}{p^\prime}
\newcommand{\pq}{q^\prime}
\newcommand{\pr}{r^\prime}
\newcommand{\pu}{u^\prime}
\newcommand{\pv}{v^\prime}
\newcommand{\pw}{w^\prime}
\newcommand{\px}{x^\prime}

\newcommand{\py}{y^\prime}
\newcommand{\palpha}{\alpha^\prime}
\newcommand{\pbeta}{\beta^\prime}
\newcommand{\peta}{\eta^\prime}
\newcommand{\pgamma}{\gamma^\prime}
\newcommand{\piota}{\iota^\prime}
\newcommand{\plambda}{\lambda^\prime}
\newcommand{\pmu}{\mu^\prime}
\newcommand{\ptau}{\tau^\prime}
\newcommand{\ppi}{\pi^\prime}
\newcommand{\prho}{\rho^\prime}

\newcommand{\pvep}{\vep^\prime}
\newcommand{\pvph}{\vph^\prime}
\newcommand{\pphi}{\phi^\prime}
\newcommand{\ppsi}{\psi^\prime}

\newcommand{\pvrh}{\vrh^\prime}
\newcommand{\pxi}{\xi^\prime}
\newcommand{\pGamma}{\Gamma^\prime}

\newcommand{\pOmega}{\Omega^\prime}

\newcommand{\psharp}{\sharp^\prime}
\newcommand{\pmfS}{\mathfrak{S}^\prime}

\newcommand{\qmsB}{^q{\hspace{-1pt}}\msB}

\newcommand{\lptau}{^p{\hspace{-1pt}}\tau}
\newcommand{\pwtau}{^{pw}{\hspace{-1pt}}\tau}

\newcommand{\rhB}{^\vrh{\hspace{-1pt}}B}

\newcommand{\pprhB}{^{\pvrh}{\hspace{-1pt}}\pB}
\newcommand{\rhmsB}{^\vrh{\hspace{-1pt}}\msB}

\newcommand{\pprhmsB}{^{\pvrh}{\hspace{-2pt}}\pmsB}
\newcommand{\rhK}{^\vrh{\hspace{-1pt}}K}

\newcommand{\prhpK}{^{\pvrh}{\hspace{-1pt}}\pK}
\newcommand{\rhk}{^\vrh{\hspace{-1pt}}k}
\newcommand{\rhp}{^\vrh{\hspace{-1pt}}p}
\newcommand{\rhpk}{^{\pvrh}{\hspace{-1pt}}k}

\newcommand{\rhtau}{^\vrh{\hspace{-1pt}}\tau}
\newcommand{\rhiota}{^\vrh{\hspace{-1pt}}\iota}
\newcommand{\rhsharp}{^\vrh{\hspace{-1pt}}\sharp}
\newcommand{\prhsharp}{^{\pvrh}{\hspace{-1pt}}\sharp}
\newcommand{\bPos}{\overline{\mathrm{Pos}}}
\newcommand{\ort}{\mathrm{ort}}
\newcommand{\spl}{\mathrm{spl}}

\newcommand{\Alt}{\mathop\mathrm{Alt}\nolimits}
\newcommand{\Ann}{\mathop\mathrm{Ann}\nolimits}

\newcommand{\Aut}{\mathop\mathrm{Aut}\nolimits}

\newcommand{\Cay}{\mathop\mathrm{Cay}\nolimits}
\newcommand{\Cent}{\mathop\mathrm{Cent}\nolimits}
\newcommand{\Cenout}{\mathop\mathrm{Cent}_\mathrm{out}\nolimits}

\newcommand{\Con}{\mathop\mathrm{Con}\nolimits}

\newcommand{\Cosp}{\mathop\mathrm{Cosp}\nolimits}
\newcommand{\Cot}{\mathop\mathrm{Cop}\nolimits}
\newcommand{\Cox}{\mathop\mathrm{DiCo}\nolimits}

\DeclareMathOperator{\Der}{Der}

\newcommand{\Diag}{\mathop\mathrm{Diag}\nolimits}

\newcommand{\End}{\mathop\mathrm{End}\nolimits}
\newcommand{\Fix}{\mathop\mathrm{Fix}\nolimits}
\DeclareMathOperator{\Ga}{Ga}
\newcommand{\GL}{\mathop \mathrm{GL}\nolimits}
\DeclareMathOperator{\Her}{Her}
\newcommand{\Hom}{\mathop\mathrm{Hom}\nolimits}
\DeclareMathOperator{\Hur}{Hur}
\newcommand{\Hyp}{\mathop\mathrm{Hyp}\nolimits}
\newcommand{\Img}{\mathop\mathrm{Im}\nolimits}

\newcommand{\Int}{\mathop\mathrm{Int}\nolimits}
\newcommand{\Isom}{\mathop\mathrm{Isom}\nolimits}

\newcommand{\Jcub}{J_{\mathrm{cub}}}

\newcommand{\Jquad}{J^{\mathrm{quad}}}

\newcommand{\Kir}{\mathop\mathrm{Kir}\nolimits}

\newcommand{\Map}{\mathop\mathrm{Map}\nolimits}
\newcommand{\Mat}{\mathop\mathrm{Mat}\nolimits}
\newcommand{\Min}{\mathop\mathrm{Min}\nolimits}
\newcommand{\Mon}{\mathop\mathrm{Mon}\nolimits}

\newcommand{\Mult}{\mathop\mathrm{Mult}\nolimits}
\newcommand{\Nil}{\mathop\mathrm{Nil}\nolimits}
\newcommand{\Nrd}{\mathop\mathrm{Nrd}\nolimits}
\newcommand{\Nuc}{\mathop\mathrm{Nuc}\nolimits}

\newcommand{\Pic}{\mathop\mathrm{Pic}\nolimits}
\newcommand{\Pocu}{\mathop\mathrm{Pocu}\nolimits}
\newcommand{\Pol}{\mathop\mathrm{Pol}\nolimits}
\newcommand{\Pos}{\mathop\mathrm{Pos}\nolimits}
\newcommand{\Quad}{\mathop\mathrm{Quad}\nolimits}
\newcommand{\Rad}{\mathop\mathrm{Rad}\nolimits}
\newcommand{\Rat}{\mathop\mathrm{Rat}\nolimits}
\newcommand{\Rex}{\mathop\mathrm{Rex}\nolimits}

\newcommand{\Spec}{\mathop\mathrm{Spec}\nolimits}

\newcommand{\Sq}{\mathop\mathrm{Sq}\nolimits}

\newcommand{\Str}{\mathop\mathrm{Str}\nolimits}
\newcommand{\Sym}{\mathop\mathrm{Sym}\nolimits}

\newcommand{\Ter}{\mathop\mathrm{Ter}\nolimits}
\newcommand{\Tri}{\mathop\mathrm{Tri}\nolimits}
\newcommand{\Zar}{\mathop\mathrm{Zar}\nolimits}
\newcommand{\Zer}{\mathop\mathrm{Zer}\nolimits}
\newcommand{\Zor}{\mathop\mathrm{Zor}\nolimits}
\newcommand{\tie}{\tilde{e}}

\newcommand{\tiw}{\tilde{w}}
\newcommand{\til}{\tilde{L}}
\newcommand{\tin}{\tilde{N}}

\newcommand{\tir}{\tilde{R}}
\newcommand{\tis}{\tilde{S}}
\newcommand{\tit}{\tilde{T}}
\newcommand{\tiX}{\tilde{X}}
\newcommand{\tiY}{\tilde{Y}}
\newcommand{\biw}{\mathop{\bigwedge}\nolimits}
\newcommand{\can}{\mathop\mathrm{can}\nolimits}
\newcommand{\ch}{\mathop\mathrm{char}\nolimits}
\newcommand{\cl}{\mathop\mathrm{cl}\nolimits}

\newcommand{\cub}{\mathrm{cub}}
\newcommand{\diag}{\mathop\mathrm{diag}\nolimits}
\DeclareMathOperator{\disc}{disc}

\newcommand{\grad}{\mathop\mathrm{grad}\nolimits}

\newcommand{\op}{\mathop\mathrm{op}\nolimits}
\newcommand{\rk}{\mathop\mathrm{rk}\nolimits}

\newcommand{\Ker}{\mathop\mathrm{Ker}\nolimits}
\renewcommand{\ker}{\Ker}  

\newcommand{\tr}{\mathop\mathrm{tr}\nolimits}
\newcommand{\Tr}{\tr}

\newcommand{\alg}{\mathchar45\mathbf{alg}}
\newcommand{\Falg}{F\alg}
\newcommand{\kalg}{k\alg}

\newcommand{\pRalg}{R^\prime\mathchar45\mathbf{alg}}
\newcommand{\Kalg}{K\alg}

\newcommand{\Ralg}{R\mathchar45\mathbf{alg}}
\newcommand{\Salg}{S\mathchar45\mathbf{alg}}

\newcommand{\kcocujo}{k\mathchar45\mathbf{cocujo}}
\newcommand{\kcopa}{k\mathchar45\mathbf{copa}}

\newcommand{\kosp}{k\mathchar45\mathbf{cosp}}

\newcommand{\kcumap}{k\mathchar45\mathbf{cumap}}
\newcommand{\kcuno}{k\mathchar45\mathbf{cuno}}

\newcommand{\kholaw}{k3\mathchar45\mathbf{holaw}}

\newcommand{\kpoc}{k\mathchar45\mathbf{pocu}}
\newcommand{\kpolaw}{k\mathchar45\mathbf{polaw}}
\newcommand{\kracuno}{k\mathchar45\mathbf{racuno}}

\newcommand{\modstub}[1]{#1\mathchar45\mathbf{mod}}
\newcommand{\kmod}{\modstub{k}}

\newcommand{\prooftype}[1]{\ifvmode\else\par\fi\medskip\noindent \emph{#1.}\ \ }

\DeclareMathOperator{\car}{char}

\DeclareMathOperator{\Lie}{Lie}

\newcommand{\PGL}{\mathbf{PGL}}
\newcommand{\PGSp}{\mathbf{PGSp}}

\newcommand{\Inv}{\mathrm{Inv}}

\DeclareMathOperator{\rmO}{O}
\newcommand{\Or}{\rmO}    

\renewcommand{\bfO}{\mathbf{O}}

\renewcommand{\sl}{\mathfrak{sl}}

\newcommand{\opgrp}[1]{#1^{\mathrm{op}}}

\newcommand{\plalg}[1]{#1^{(+)}}  

\newcommand{\Id}{\mathbf{1}}

\newcommand{\qbar}{\overline{q}}
\newcommand{\Qbar}{\overline{Q}}

\newcommand{\aclosure}[1]{\bar{#1}}  
\newcommand{\aclosk}{\aclosure{k}} 
\newcommand{\aclosF}{\aclosure{F}}

\newcommand{\kx}{k^\times}  
\newcommand{\Fx}{F^\times}

\newcommand{\trans}{\mathsf{T}}   

\newcommand{\emptyslot}{{\--}}  

\newcommand{\skiprootsystem}{\mathsf}
\newcommand{\rsF}{\skiprootsystem{F}_4}

\newcommand{\rsA}{\skiprootsystem{A}}

\newcommand{\rsE}{\skiprootsystem{E}}

\newcommand{\bform}[1]{\la #1 \ra}
\newcommand{\qform}[1]{\la #1 \raq}

\newcommand{\Proofend}{\hfill $\square$}

\newcommand{\step}[1]{\medskip $#1^\circ$.} 
\newcommand{\case}[1]{\smallskip \emph{Case #1:}} 

\newcommand{\stacks}[1]{\cite[\href{https://stacks.math.columbia.edu/tag/#1}{Tag #1}]{stacks-project}}

\newcommand{\noexsec}{This section contains no exercises.}

\ifnum\exmode=0  
	\newcommand{\todo}[1]{}  
\else
	\newcommand{\todo}[1]{\footnote{\emph{\textcolor{red}{To do:}} #1}}
\fi

\newcommand{\stbtmat}[4]{\left( \begin{smallmatrix} #1 & #2 \\ #3 & #4 \end{smallmatrix} \right)}  

\newcommand{\mfsl}{\mathfrak{sl}}
\newcommand{\mfsp}{\mathfrak{sp}}

\DeclareMathOperator{\LMDer}{LMDer}
\DeclareMathOperator{\InDer}{InDer}

\DeclareMathOperator{\ComDer}{ComDer}

\newcommand{\dref}[2]{\ref{#2}}

\usepackage[colorlinks=true, urlcolor=blue, backref=page, citecolor=green!60!black, linkcolor=red!80!black,plainpages=false,breaklinks=true,pdfstartview=FitH]{hyperref}

\allowdisplaybreaks[2]  

\makeindex             

\begin{document}

\author{Skip Garibaldi, Holger P. Petersson and Michel L. Racine}
\title{Solutions to the exercises\\ from the book\\ \textbf{Albert algebras over commutative rings}}
\maketitle

\frontmatter

\tableofcontents

\mainmatter

\backmatter

%
%
\renewcommand{\theequation}{s\arabic{equation}}

\ifnum\exmode=0
   \chapter*{Introduction}
\else
   \chapter{Solutions} 
\fi

\chapspace

What follows are errata, addenda, and solutions to the exercises in the book \emph{Albert algebras over commutative rings} \cite{GPR:book}.

The current version of this document can be found on the arXiv at \href{https://arxiv.org/abs/2406.02933}{2406.02933}.  A shorter text containing just the errata and addenda is also available at Skip's website \href{https://www.garibaldibros.com}{garibaldibros.com}; it may be updated more frequently then the text on the arxiv.  There you can also find a current ``draft'' version of the book with the errata corrected.  Please note that the draft version of the book preserves the numbering of theorems, definitions, and so on, but it has different page numbering.

\subsection*{Equation numbering} In order to reduce ambiguity in cross-references, equations are numbered differently in the printed book from how they are numbered here.  In the book, the equations are numbered (1), (2), etc., whereas here in the addenda we number them (a1), (a2), etc., and in the  solutions we number them  (s1), (s2), etc.

%
%
\newcommand{\typo}[4]{\medskip\noindent \emph{Page #1, item #2, line #3 {\small (fixed #4)}}: }
\newcommand{\typop}[3]{\medskip\noindent \emph{Page #1, #2 {\small (fixed #3)}}: }
\newcommand{\rcite}[2]{#1\footnote{For the convenience of the reader: [#1] in the printed book = \cite{#2}.}}

\chapter*{Errata for the published edition of \emph{Albert algebras over commutative rings}}

\chapspace

Page numbers refer to the published/paper edition of the book, published by Cambridge in November 2024.  These typos have been corrected in the online/``draft'' version you can find on the web.  The online/draft versions have a date on their front cover, and the ``fixed'' date says that the typo has been corrected in all online/draft versions starting with that date. 

\ifdefined\errata
 This errata sheet was compiled on \today.
\fi

Thank you to Alberto Elduque, Darij Grinberg, and Eric Rains for identifying some of the following.

If you find other typos, please do tell us!  You can send us an email or use the contact form on Skip's website at \href{https://www.garibaldibros.com/wp/contact/}{this link}.

\section*{Notation and conventions}

\typop{xxii}{line -17}{20 Feb 2025} The $R$ in $M\times M\rightarrow R$ should be replaced by $k$.

\section*{Chapter I: Prologue: the ancient protagonists}

\typo{38}{6.4}{4}{16 Dec 2024} The $U$-operator maps $U\: J \to \End(J)$.  That is, the codomain is $\End(J)$ not $J$.

\section*{Chapter II: Foundations}

\typo{94}{12.16}{-2}{9 Nov 2025} The $(Df)$ on this line should be $(Df)_R$.

\typo{95}{12.17}{1 after (8)}{9 Nov 2025} Insert the following sentence before ``Note'':
To verify each of (1)--(8), one can use (12.15.4) or (12.15.5) to expand both sides and equate coefficients of $\varepsilon$.

\typo{103}{12.36}{-1}{26 Nov 2025} The lower bound on $p$ should read ``$1 \le p$'' instead of ``$0 \le p$''.
 
\typo{104}{12.39}{2}{2 Jan 2025} Insert after the first sentence: (See for example [\rcite{27}{MR0360549}, \S{VI.3.6}] for background on discrete valuations on fields.) 

\section*{Chapter IV: Composition algebras}

\typo{157}{19.24}{5}{16 Dec 2024} Replace $\Her_r(\ID)$ with $\Her_r(\IO)$.

\section*{Chapter V: Jordan algebras}

\typo{275}{29.11}{3}{20 Feb 2025} Replace ``that that'' with a single ``that''.

\section*{Chapter VI: Cubic Jordan algebras}

\typo{344}{34.25}{-2}{19 Jan 2026} Replace ``in (c)'' with ``in (b)''.

\section*{Chapter VII: The two Tits constructions}

\typo{458}{42.18}{1}{16 Jan 2025} The first sentence should read: 
Let $J = J(D, \mu)$ be an Albert algebra arising from a first Tits construction over a field $F$ of characteristic $\ne 3$.

\typo{495}{45.1}{-2}{31 Dec 2025} Delete the phrase ``of the second kind''.

\typo{509}{46.6 proof}{3}{16 Dec 2024} Should refer to Lemma 46.5, not Proposition.

\section*{Chapter VIII: Lie algebras}

\typop{527}{line 1}{20 Feb 2025} The first sentence should read: Let $V$ be the subspace of $E = \IR^{8}$ whose points have coordinates $(\xi_i)$ such that $\xi_6 = \xi_7 = - \xi_8$.

(That is, replace every $\epsilon$ in the printed version with a $\xi$, to align with the notation in \cite{Bou:g4}.)

\typo{529}{47.17}{4 of the remark}{20 Feb 2025} Replace ``very'' with ``every".

\typop{565}{51.29}{20 Feb 2025} The statement of part (a) should also include the hypothesis $2 \in \kx$, in order for the proof provided to be sufficient.  The proof needs that the map $\mfg_0 \to \mfo(n_C)$ is bijective from 51.26(b), which relies on $2 \in \kx$.

The claim in part (a), that $\Der(J)$ is finitely generated projective of rank 52, does hold in greater generality.  Here is a proof using the material from the next chapter, Chapter IX.  

First suppose that $k$ is a field, in which case we only need to prove that $\dim(\Der(J)) = 52$.
By 52.1(b), $\Lie(\bfAut(J)) = \Der(J)$.  By Theorem 53.4, $\bfAut(J)$ is semisimple of type $\rsF$, so $\dim \Lie(\bfAut(J)) = \dim \bfAut(J)$   because $\bfAut(J)$ is smooth and  $\dim \bfAut(J) = 52$  because type $\rsF$. 

Next suppose that $J$ is the split Albert algebra $\Her_3(\Zor(k))$ over a principal ideal domain $k$.  Since $\Zor(k)$ is a finitely generated free module, so is $\End(\Zor(k))$; since additionally $k$ is a principal ideal domain we conclude that the submodule $\Der(J)$ is also finitely generated and free \stacks{0AUW}.  The field of fractions $F$ of $k$ is a flat extension of $k$, so $\Der(J)_F \cong \Der(J_F)$ by Prop.~50.4, and this has dimension 52 by the previous paragraph, hence $\Der(J)$ has rank 52 at the prime 0.  Since $k$ is connected, $\Der(J)$ has constant rank 52.

If $J$ is any Albert algebra over a principal ideal domain $k$, then there is some faithfully flat $K \in \kalg$ such that $J_K$ is split.  Then $\Der(J)_K \cong \Der(J_K) \cong \Der(\Her_3(\Zor(K))) \cong \Der(\Her_3(\Zor(k))) \otimes K$ is finitely generated projective of rank 52 by the previous paragraph, and it follows that $\Der(J)$ is finitely generated projective of rank 52.

\typo{557}{(51.10.2)}{1}{16 Dec 2024} The right side of the equation should read $\{(a_1,a_1^\trans)\mid a_1 \in \Mat_3(k)\}$.

\typop{566}{alternative proof of 51.29(b), line 2}{20 Feb 2025} Replace $\Zor(k)$ with $\Her_3(\Zor(k))$.

\section*{Chapter IX: Group schemes}

\typo{567}{52.1, paragraph 2}{-2}{22 Feb 2025} Replace $G$ with $\bfG$.

\typo{588}{proof of 54.6}{3}{20 Feb 2025} Replace ``(1\,3\,2)'' by ``(1\,2)''.

\typo{593}{54.11}{2}{20 Feb 2025} Replace ``structure'' by ``system''.

\typo{593}{54.12}{2}{20 Feb 2025} Replace ``structure'' by ``system''.

\typop{603}{55.12}{29 Mar 2025} Replace the statement of the exercise with the following, clearer, version:

Let $F$ be a field.  Show that the maps $(A,\tau) \mapsto H(A,\tau)$ and $J \mapsto \bfAut(J)$ define bijections between the isomorphism classes of 
\begin{enumerate}[(i)]
\item Azumaya algebras $(A, \tau)$ of degree 3 with unitary involution over $F$, as in 44.23,
\item rank 9 Freudenthal $F$-algebras $J$, and
\item adjoint semi-simple $F$-group schemes of type $\rsA_2$.
\end{enumerate}

\typop{619}{bottom of page}{20 Feb 2025} Insert new paragraph: Results for groups of type $\rsE_7$ --- analogous to the results in this section for type $\rsE_6$ --- can be found in sections 16 and 17 of [\rcite{95}{GPR}].

\section*{Subject Index}

\typop{648}{middle of left column}{1 Dec 2025} Insert line: balanced pair, 470


\renewcommand{\todo}[1]{\footnote{\emph{\textcolor{red}{To do:}} #1}}

\renewcommand{\theequation}{a\arabic{equation}}
\renewcommand{\thesection}{A\arabic{section}}

\chapter*{Addenda to the published edition of \emph{Albert algebras over commutative rings}}
\chaptermark{Addenda to \emph{Albert algebras over commutative rings}}

This chapter contains material that could have been in the book and does not require substantial extra background. 

\setcounter{section}{0}
\setcounter{equation}{0}

\section{Albert algebras are exceptional}

In his first paper on the subject of Jordan algebras, Albert \cite{MR1503142} showed that (1) the euclidean Albert algebra $\Her_3(\IO)$ over the reals is a (linear) Jordan algebra and (2) that it is exceptional.  We have shown (1) in Theorem 5.10.  More generally, for every octonion algebra $C$ over every ring $k$, we have shown that $\Her_3(C)$ is a Jordan algebra in Theorem 36.5.  (Note that Albert algebras are defined in the book to be cubic Jordan algebras, hence Jordan algebras.)
 In this addendum we would like to extend (2), the exceptionality of Albert algebras, to an arbitrary ring $k$.  We  prove the following slightly stronger result.

\begin{thm} \label{t.ALBEX}
Every Albert algebra over every non-zero ring is i-exceptional. 
\end{thm}

We will define the term i-exceptional in a moment; it is stronger than the notion of exceptional defined in 29.9.  When the ring is a field of characteristic different from 2, Theorem \ref{t.ALBEX} was first proved in \cite{AlbertPaige}.  A proof of Theorem \ref{t.ALBEX} for every field can be found in Jacobson's Arkansas notes \cite[\S2.5]{MR634508}.  That proof considers separately the cases of fields of characteristic 2 and different from 2.  We provide a proof that does not contain special considerations involving 2.

In the excluded case of the zero ring, all algebras are zero and so the notion of being exceptional or i-exceptional do not make sense.

Here is a beautiful application of the theorem.

\begin{thm}
A central simple Jordan algebra over a field is exceptional if and only if it is an Albert algebra. 
\end{thm}

\begin{proof}
The central simple exceptional Jordan algebras over a field have been classified into types in \S15 of 
\cite{MR946263}.  One finds that each type is either special or an Albert algebra.  Albert algebras are exceptional by Theorem \ref{t.ALBEX}.
\end{proof}

We now provide the promised definition.

\begin{definition}
A Jordan algebra $J$ is \emph{i-special} if $J$ is a homomorphic image of a special algebra and \emph{i-exceptional} if it is not.
\end{definition}

Clearly any special Jordan algebra is i-special. However P.M.~Cohn has given an example of an i-special Jordan algebra that is exceptional \cite[\S{I.3}, Thm.~2]{MR0251099}. Therefore being i-exceptional is stronger than being exceptional.   We also have the following:

\begin{lem} \label{l.EXBASE}
Let $J$ be a Jordan $k$-algebra.  If $J_R$ is i-exceptional for some flat $R \in \kalg$, then $J$ is i-exceptional.
\end{lem}

\begin{proof}
We prove the contrapositive.
Suppose that $J$ is not i-exceptional, i.e.,  there is a special Jordan algebra $J'$ and a (surjective) homomorphism $f\:J' \rightarrow J$. Hence there exists a unital associative algebra $A$ and an injective Jordan homomorphism $j\: J' \rightarrow A^{(+)}$. It follows that $f_R\: J'_R \rightarrow J_R$ is still surjective and, by flatness, $j_R\: J'_R \rightarrow (A_R)^{(+)}$ is still injective, i.e., $J_R$ is i-special.
\end{proof}

As a first step in the proof of Theorem \ref{t.ALBEX}, we note the following, which could have appeared in 37.7.

\begin{lem} \label{l.TP} 
{For $J = \Her_3(C)$, $C$ a multiplicative conic alternative algebra, and $a[ij]$, $b[ji]$, $c[il] \in J$ with $\{i, j, k\} = \{1, 2, 3\}$,} 
\begin{align}
\label{TPijl} \{a[ij]\, b[ji]\, c[il]\} = a(bc)[il]. 
\end{align}
\end{lem}

We remark that if $C$ is associative, then (\ref{TPijl}) follows immediately from matrix multiplication. 

\begin{proof} 
Using (33.6.2), (36.4.6), (36.4.7),
\begin{align*}
\{a[ij]\, b[ji]\, c[il]\} &= T_J(a[ij], \bar{b}[ij])c[il] - (c[il] \times a[ij]) \times b[ji] \\
&= n_C(a,\, \bar{b})c[il] - \bar{a}c[jl] \times b[ji] \\
&= (n_C(a,\, \bar{b})c - \bar{b}(\bar{a}c))[il].   
\end{align*}
Replace $\bar{a}$ and $\bar{b}$ using (16.5.4) and expand. By (16.5.5), Proposition 16.10 and alternativity, we obtain (\ref{TPijl}).
\end{proof}

%
%
%

The following polynomial map was introduced by Glennie \cite{MR086708}:
\begin{multline*}
g_9(x,y,z) := U_xz\circ\{x\,U_zy^2\,y\} - U_yz\circ\{x\,U_zx^2\,y\} \\
 - U_xU_z\{x\,U_yz\,y\} + U_yU_z\{x\,U_xz\,y\}. 
\end{multline*}
It is homogeneous of degree 9 and skew-symmetric in $x$ and $y$.

\begin{lem} \label{l.GLENEX}
Let $C$ be a multiplicative conic alternative algebra.  If $C$ is not associative, then $g_9$ is not identically zero on $\Her_3(C)$.
\end{lem}

\begin{proof}
Consider the following substitution,
\[
x = 1[12], \quad y = 1[23],\quad z = a[12] + b[23] + c[31], \quad x, y, z \in \Her_3(C), 
\]
for $a$, $b$, $c \in C$. Using (16.12.2), (17.4.2), the Peirce decomposition rules (32.15), the identities of 37.7 and Lemma \ref{l.TP}, one computes 
\[
x^2 = e_{11} + e_{22}, \quad y^2 = e_{22} + e_{33}, \quad U_xz = \bar{a}[12], \quad U_yz = \bar{b}[23].
\]
Futhermore
\begin{align*}
U_zx^2 &= n_C(a)e_{11} + n_C(a)e_{22} + (n_C(b) + n_C(c))e_{33} + \bar{a}\,\bar{c}[23] + \bar{b}\,\bar{a}[31]. \\
U_zy^2 &= n_C(b)e_{22} + n_C(b)e_{33} + (n_C(a) + n_C(c))e_{11} + \bar{c}\,\bar{b}[12] + \bar{b}\,\bar{a}[31]. \\
\{x\, U_zy^2\, y\} &= t_C(ab)e_{22} + \bar{c}\,\bar{b}[23] + n_C(b)1[31].  \\
\{x\, U_zx^2\, y\} &= t_C(ab)e_{22} + \bar{a}\,\bar{c}[12] + n_C(a)1[31]. \\
U_xz\circ\{x\,U_zy^2\,y\} &= t_C(ab)\bar{a}[12] + n_C(b)a[23] + (bc)a[31]. \\
U_yz\circ\{x\,U_zx^2\,y\} &= n_C(a)b[12] + t_C(ab)\bar{b}[23] + b(ca)[31]. \\
\{x\,U_yz\, y\} &= \bar{b}[12],\\
U_z\{x\,U_yz\, y\} &= n_C(b)t_C(c)e_{33} + aba[12] + n_C(b)\bar{a}[23] + (c\bar{b})\bar{a}[31]. \\
U_xU_z\{x\,U_yz\, y\} &= \overline{aba}[12]. \\
\{x\, U_xz\, y\} &= \bar{a}[23]. \\
U_z\{x\, U_xz\, y\} &= n_C(a)t_C(c)e_{11} + n_C(a)\bar{b}[12] + bab[23] + \bar{b}(\bar{a}c)[31]. \\
U_yU_z\{x\, U_xz\, y\} &= \overline{bab}[23]. 
\end{align*}
Therefore
\begin{multline*}
g_3(x,y,z) = (t_C(ab)\bar{a} - n_C(a)b - \overline{aba})[12] + (n_C(b)a - t_C(ab)\bar{b} + \overline{bab})[23] \\ 
  + [b,c,a][31]. 
\end{multline*}
The [12] and [23] entries are 0 by (17.4.2). If $C$ is not associative, one can choose $a$, $b$, $c \in C$ such that the associator $[a,b,c] \neq 0$, in which case the [31] entry is not 0.
\end{proof}

\begin{lem}  \label{l.GLEN}
The Glennie polynomial $g_9(x, y, z)$ vanishes on every i-special Jordan algebra.
\end{lem}

\begin{proof}
A special Jordan algebra is a subalgebra of $A^{(+)}$ for some associative algebra $A$. Let $u$, $v$, $w \in A$. Substituting $x = u$, $y = v$, $z = w$  in $g_9(x,y,z)$, we obtain
\begin{multline*}
uwuuwv^2wv + uwuvwv^2wu + uwv^2wvuwu + vwv^2wuuwu \\
-vwvuwu^2wv - vwvvwu^2wu -uwu^2wvvwv - vwu^2wuvwv \\
-uwuvwvvwu - uwvvwvuwu + vwuuwuvwv + vwvuwuuwv 
= 0. 
\end{multline*}

If $J$ is an i-special Jordan algebra, then by definition there is a surjection $f\: J' \to J$ such that $J'$ is special.  Then $g_9$ on $J$ can be computed as the composition $g_9 f$ on $J'$, which is identically zero by the preceding paragraph.
\end{proof} 

In fancier language, we have shown that $g_9$ is an \emph{s-identity}, i.e., it is not identically zero for every Jordan algebra (Lemma \ref{l.GLENEX}) and it vanishes on every special Jordan algebra (Lemma \ref{l.GLEN}).

%
%

\begin{proof}[Proof of Theorem \ref{t.ALBEX}]
If $J$ is an Albert algebra, by Corollary 39.32, there exists a faithfully flat extension $R \in \kalg$ such that $J_R \cong \Her_3(\Zor(R))$.  Since $\Zor(R)$ is not associative, Lemmas \ref{l.GLENEX} and \ref{l.GLEN} show that $J_R$ is not i-special, i.e., is i-exceptional, and it follows (Lemma \ref{l.EXBASE}) that $J$ is i-exceptional.
\end{proof}

\section{Descent of algebras to Dedekind domains}

Chapter I of the book, especially section 6, addressed the question of what properties of a $\IZ$-submodule of a real Jordan or Albert algebra correspond to it being a Jordan or Albert algebra over $\IZ$.  That part of the text, however, takes a naive view of the objects involved and therefore does not obviously give statements about the more sophisticated notions of cubic Jordan and Freudenthal algebras defined later in the text.  We rectify that here.

The general setting we describe is some integral domain $k$ contained in a field $F$ and we wish to give properties of a $k$-submodule $M$ of some algebraic object over $F$ to guarantee that $M$ is also of the same type of object.  We need the following basic tool.

\begin{lem} \label{permanence}
Let $M$, $N$ be $k$-modules for some integral domain $k$, and let $f \colon M \to N$ be a polynomial law.  If there exists a field $F$ containing $k$ such that $f \otimes F = 0$, then $f = 0$ as a polynomial law.
\end{lem}

\begin{proof}
If $k$ is itself a field, then $F$ is a faithfully flat $k$-algebra.  Since $0 \otimes F = f \otimes F$, $f = 0$ by Exercise 25.35(b).  If $k$ is not a field, then it is infinite, so $F$ is infinite.  The Principle of Permanence \cite[\S{IV.2.3}, Scholium]{MR643362} then says that $f_R = 0$ for every $k$-algebra $R$ and therefore $f = 0$.
\end{proof}

\subsection*{Jordan algebras}

Recall from Definition 29.1 that a Jordan algebra is a para-quadratic algebra such that two identities concerning the $U$-operator hold strictly, meaning as polynomial laws.

\begin{prop}[Jordan algebras, first version] \label{jordan.1}
Let $A$ be a para-quadratic $k$-algebra where $k$ is an integral domain.  If there is a field $F$ containing $k$ such that $A \otimes_k F$ is a Jordan $F$-algebra, then $A$ is a Jordan $k$-algebra.
\end{prop}

\begin{proof}
We must verify that the polynomial laws $A \times A \times A \to A$ given by sending $(x, y, z)$ to 
\[
U_{U_x y}z - U_x U_y U_x z\quad \text{and} \quad U_x V_{y,x}z - V_{x,y} U_x z
\]
are zero as polynomial laws.  This follows from Lemma \ref{permanence}.
\end{proof}

We can rephrase this in terms of lattices.  Let $M$ be a vector space over a field $F$.  For an integral domain $k$ contained in $F$, a \emph{$k$-lattice} $\Lambda$ in $M$ is a finitely generated $k$-submodule of $M$ such that the natural $F$-linear map $\Lambda \otimes_k F \to M$ is an isomorphism of $F$-modules.

\begin{prop}[Jordan algebras, second version] \label{jordan.2}
Suppose $J$ is a Jordan $F$-algebra for a field $F$, $k$ is an integral domain contained in $F$, and $\Lambda$ is a $k$-lattice $J$ that contains $1_J$.  If 
\begin{enumerate}[(1)]
\item \label{jordan.2.U} $\Lambda$ is closed under the $U$-operator for $J$, or
\item \label{jordan.2.sq} $2 \in \kx$ and $x^2 \in \Lambda$ for every $x \in \Lambda$,
\end{enumerate}
then the restriction of the $U$-operator of $J$ turns $\Lambda$ into a Jordan $k$-algebra.
\end{prop}

\begin{proof}
Assume \ref{jordan.2.U}, i.e., the restriction of $U$ to $\Lambda$ defines a function $\Lambda \to \End_k(\Lambda)$.  Since $U$ on $J$ is a quadratic map in the sense of 11.1, so is the restriction of $U$ to $\Lambda$.  Therefore, $\Lambda$ together with $1_J$ and the restriction of $U$ is a para-quadratic $k$-algebra.  Now apply Proposition \ref{jordan.1}.

Now assume \ref{jordan.2.sq}. Since $2 \in \kx$, we view $J$ as a linear Jordan algebra and write $xy$ for the product of $x, y \in J$.  Again because $2 \in \kx$, for $x, y \in \Lambda$, $xy := \frac12 ((x+y)^2 - x^2 - y^2)$ belongs to $\Lambda$.  It follows that $\Lambda$ is closed under the $U$-operator (27.10.1), and we are done by \ref{jordan.2.U}.
\end{proof}

\subsection*{Jordan algebras of degree 3}

Next we wish to treat Jordan algebras of degree 3, for which we will restrict to the case where $k$ is a Dedekind domain.  In that case, because a $k$-lattice $\Lambda$ is necessarily torsion-free (and by definition finitely generated), it follows that $\Lambda$ is projective \cite[Tags 00NX, 0AUW]{stacks-project}.

Recall that a cubic Jordan algebra can be defined as a Jordan algebra $J$ together with a cubic form $N$ on $J$ as in 34.1, or as a cubic norm structure meaning a module $J$ together with a base point $1_J$, quadratic map $\sharp$, and cubic form $N$ as in 33.1 and 33.4; the two notions are interchangeable by Corollary 34.6.  Certain identities regarding these maps are required to hold as polynomial laws.

A Jordan algebra of degree 3 is defined in 38.1(b) as a cubic Jordan algebra $J$ such that the set map $J_K \to \wedge^3 J_K$ given by $x \mapsto 1_{J_K} \wedge x \wedge x^2$ is not the zero map for every algebraically closed field $K \in \kalg$.

\begin{prop}[Compare Theorem 6.11] \label{cubic.jordan}
Let $F$ be a field of characteristic $\ne 2$, and suppose $J$ is a Jordan algebra of degree 3 over $F$. 
Let $k$ be an infinite Dedekind domain contained in $F$.  For a $k$-lattice $\Lambda$ in $J$ containing $1_J$, the following are equivalent:
\begin{enumerate}[(i)]
\item \label{611.1} The restriction of the $U$-operator and $N$ turn $\Lambda$ into a cubic Jordan algebra.
\item \label{611.2} $x^2 \in \Lambda$ for all $x \in \Lambda$.
\item \label{611.3} $x^\sharp \in \Lambda$ for all $x \in \Lambda$.
\end{enumerate}
\end{prop}

\begin{proof}
It is trivial that \ref{611.1} implies \ref{611.2} and \ref{611.3}, so the work is in proving the opposite implications.
The idea is to imitate the proof of Theorem 6.11 in the book, substituting items according to the following table, where $K$ denotes the fraction field of $k$:
\[
\begin{array}{c|cccc}
\text{book} & \IZ & \IQ & \IR & \Her_3(\ID) \\
\text{here} & k & K & F & J
\end{array}
\]
We also replace various identities regarding $\Her_3(\ID)$ from Chapter I with the corresponding properties for cubic Jordan algebras from Chapter VI, for example: (5a.8) is replaced by (34.1.4); Exercise 6.16 is replaced by (33a.6); and Exercise 6.12(b) is replaced by (33.6.1) and (33.24).  We now step through the various arguments leading up to Theorem 6.11, and comment on what further changes are required.

The proof of Proposition 3.5 goes through, using that a Dedekind domain is integrally closed.  Lemma 3.7, Proposition 3.8, and Exercises 3.19 and 6.15 go through.

Lemma 6.7 shows that \ref{611.2} implies that the trace $T$, quadratic trace $S$, and norm $N$ of $J$ restrict to maps $\Lambda \otimes K \to K$.  In the proof, we use that $K$ has characteristic $\ne 2$ so that we can work with linear Jordan algebras, and so that \ref{611.2} gives that $\Lambda \otimes K$ is a $K$-algebra.  For the Zariski-density argument following equation (1) in the proof of Lemma 6.7, 
we require $J$ to be a Jordan algebra of degree 3 and not merely a cubic Jordan algebra in order that the set of elements in $k$ with minimal polynomial of degree 3 is not empty.\footnote{Note that this property is lacking for the example in Exercise 34.27 that shows that cubic Jordan algebras are not stable under faithfully flat descent.}
We require $k$ to be infinite so that the Zariski-open set constructed there has a $k$-point.

Remark 6.8 observes that every $K$-subspace $B$ of $J$ that is closed under squaring is a $K$-subalgebra, because $\car K \ne 2$.

Proposition 6.9 shows that, if $\Lambda$ is closed under taking powers $x \mapsto x^n$ for $n \in \IN$, then every $x \in \Lambda$ is integral over $k$ and has $T(x)$, $S(x)$, $N(x) \in k$.  In the proof, when Cor.~5.14 is invoked to produce elements $c_1, c_2 \in \IR[x]$, we instead produce elements in $\aclosF[x]$, where $\aclosF$ denotes the algebraic closure of $F$.

Assuming \ref{611.2}, the argument in the proof of Theorem 6.11 gives \ref{611.3}, where 
Lemma 6.10 is replaced by Lemma \ref{lem6.10} below.  

Assuming \ref{611.3}, the argument in the proof of Theorem 6.11 shows that the adjoint $\sharp$ and norm $N$ on $J$ restrict to $\Lambda$ to give $\Lambda$ the structure of a cubic array as defined in 33.1.  To be a cubic array requires that  the element $1_J$ be unimodular, but that is automatic by Exercise 12.44 since $\Lambda$ is a projective $k$-module.  

For this cubic array to be a cubic norm structure, certain identities must hold strictly, i.e., as polynomial laws.  For each such polynomial law $f$ that we wish to prove is zero, we use Lemma \ref{permanence} to reduce to checking it over $F$, where it holds by hypothesis.  The equivalence between cubic norm structures and cubic Jordan algebras (Cor.~34.6) completes the proof of \ref{611.1}.
\end{proof}

Here is the promised generalization of Lemma 6.10 about $\IZ$ to the case of a Dedekind domain.

\begin{lem}[Compare Lemma 6.10]  \label{lem6.10}
Let $k$ be a Dedekind domain with quotient field $F$ and suppose $X \subseteq F$ is closed under squaring (i.e., $\xi \in X$ $\Rightarrow$ $\xi^2 \in X$). If the $k$-submodule of $F$ generated by $X$ is finitely generated, then $X \subseteq k$.
\end{lem}

The proof is only cosmetically different from the one in the book.

\begin{proof}
For any $\xi \in F$, we denote by
\[
\mfd_\xi := \{b \in k\mid b\xi \in k\} \subseteq k
\]
the \emph{denominator ideal} of $\xi$. These ideals satisfy the following obvious relations for all $\xi,\eta \in F$ and $a \in k$:
\begin{equation}
\label{OBRE} 
\mfd_\xi = k \Longleftrightarrow \xi \in k, \quad \mfd_\xi\mfd_\eta \subseteq \mfd_\xi \cap \mfd_\eta \subseteq \mfd_{\xi + \eta}, \quad \mfd_\xi \subseteq \mfd_{a\xi}.
\end{equation}
Slightly less obvious is
\begin{align}
\label{PRIDE} \mfd_\xi = \prod_{\mfp,v_\mfp(\xi)<0} \mfp^{-v_\mfp(\xi)}.
\end{align}
In order to see this, let $b \in k$. Then $b \in \mfd_\xi$ if and only if $v_\mfp(b) \ge -v_\mfp(\xi)$ for all $\mfp$, which holds trivially for $v_\mfp(\xi) \ge 0$ and thus is  equivalent to its validity for all $\mfp$ having $v_\mfp(\xi) < 0$, that is, to $b$ belonging to the right-hand side of \eqref{PRIDE}.  This verifies \eqref{PRIDE}.
As an immediate consequence, we conclude
\begin{align}
\label{DESQ} \mfd_{\xi^2} = \mfd_\xi^2.
\end{align} 

Now let $\xi_1,\dots,\xi_n$ be a finite sequence in $X$ generating the $k$-submodule of $F$ generated by $X$. Then any $\xi \in X$ may be written as $\xi = \sum_{i=1}^n a_i\xi_i$, for some $a_1,\dots, a_n \in k$, which by \eqref{OBRE} implies
\[
\mfd_\xi \supseteq \prod_{i=1}^n\mfd_{a_i\xi_i} \supseteq \mfd := \prod_{i=1}^n \mfd_{\xi_i},
\]
the latter being a fractional ideal of $F$ depending only on $X$. Thus $\mfd_\xi$ divides $\mfd$ for all $\xi \in X$. But $\mfd$ has only finitely many divisors, so $\{\mfd_\xi\mid \xi \in X\}$ is a finite set. On the other hand, the fractional ideals of $F$ form a free, hence torsion-free, abelian group, and for any $\xi \in X$ having $\mfd_\xi \ne (1) = k$, the hypothesis of the lemma would yield a sequence $(\xi^{2^n})_{n\ge 0}$ of elements in $X$ such that the $\mfd_{\xi^{2^n}}$, $n = 0,1,2, \ldots$, by \eqref{DESQ} would be mutually distinct, a contradiction. Thus $\mfd_\xi = k$ for all $\xi \in X$, and \eqref{OBRE} implies $X \subseteq k$.
\end{proof}

%
%
%
%
%
%

\subsection*{Albert algebras}

Next we wish to give a result for Albert algebras that is similar to Propositions \ref{jordan.1}, which we will do in Theorem \ref{albert.descent}.  Recall that an Albert $k$-algebra $J$ is a cubic Jordan algebra whose underlying module is finitely generated projective\footnote{The definition in 39.8 has that the module is projective and its rank is locally constant, which is a priori a weaker condition, yet Corollary 39.11 proves that the module is finitely generated.  Conversely, a finitely generated projective module has locally constant rank \stacks{00NX}.} and such that $J_K$ is simple of rank 27 for every field $K \in \kalg$.  

The heavy lifting for this is done by the following lemma, which encapsulates arguments by Racine and van der Blij--Springer.

\begin{lem} \label{racine.lem}
    Let $F$ be a field that is complete with respect to a discrete valuation and let $J$ be a cubic Jordan algebra over the valuation ring $k$ of $F$.  Then $J$ is the split Albert $k$-algebra if and only if 
    \begin{enumerate}
        \item[(i)]  $J \otimes_k F$ is the split Albert $F$-algebra and
        \item[(ii)] the trace bilinear form on $J$ is regular.
    \end{enumerate}
\end{lem}

\begin{proof}
    Since the trace bilinear form on an Albert algebra is regular (39.19(b)) and split Albert algebras remain so after base change, the two conditions are clearly necessary.

    Conversely, given the two conditions and since $k$ is a Dedekind domain, \cite[p.~115, \S{IV.4}, Prop.~5]{MR0399187} says that $J$ is isomorphic to $\Her_3(C)$ for $C$ a maximal order in the split octonions over $F$.  (We note that while the isomorphism statements for $J$ in \cite{MR0399187} are as quadratic Jordan algebras, such isomorphisms are isomorphisms also as cubic Jordan algebras by Cor.~38.18.)  Since $k$ is a complete discrete valuation ring, all maximal orders in the split octonions are isomorphic 
    \cite[(3.4)]{MR0152555}, so we may assume that $C$ is the split octonions over $k$, ergo $J$ is itself the split Albert algebra.  
\end{proof}

\begin{thm} \label{albert.descent}
    Let $k$ be a Dedekind domain and suppose $J$ is a cubic Jordan algebra over $k$ that is regular.  If there is a field $F \supseteq k$ such that $J \otimes_k F$ is an Albert algebra, then $J$ is an Albert algebra.
\end{thm}

Recall from 33.3 that $J$ is defined to be regular if its underlying $k$-module is finitely generated projective and the bilinear trace form is regular as a symmetric bilinear form.

\begin{proof}
Because the rank of $J$ is locally constant and $k$ is an integral domain, the rank of $J$ at every prime equals $\dim_F(J_F) = 27$.  Consequently, it suffices to prove that $J_K$ is an Albert algebra for every field $K \in \kalg$.  So suppose that $K = k(\mfp)$ for some prime ideal $\mfp$ of $k$. 

If $\mfp = 0$, then $K$ is the fraction field of $k$. 
Since $F$ is faithfully flat over $K$, the base change $J_K$ is an Albert algebra by Cor.~39.32.  

Assume $\mfp \ne 0$, and write 
$E$ for the completion of $F$ with respect to the discrete valuation defined by $\mfp$.  Now $J_E$ is an Albert algebra so there is a finite\footnote{In fact, separable of dimension dividing 6, by Proposition \ref{albert.split}.} extension $E'$ of $E$ that splits $J_E$ by Cor.~54.19 and Prop.~55.5.  The discrete valuation on $E$ extends to one on $E'$ and we take $R$ to be the valuation ring of $E'$ with respect to that extension.  Since $E'$ is finite over $E$, it is complete with respect to the valuation.  Now $J_R$ is regular because $J$ is, and $J_{E'}$ is the split Albert algebra, so $J_R$ is the split Albert algebra by Lemma \ref{racine.lem}.  For $\ell$ the residue field of $R$, then, $J_\ell$ is the split Albert algebra.  But $\ell$ is a finite extension of $K$, so we deduce that $J_K$ is also an Albert algebra.

For any $K \in \kalg$, there is some prime ideal $\mfp$ (possibly zero) such that $K \supseteq k(\mfp)$.  Since $J_{k(\mfp)}$ is an Albert algebra, so is $J_K$.
\end{proof}

\borel{Application} Pick $k$ to be an infinite Dedekind domain of characteristic $\ne 2$.  Let $F$ be any field containing $k$, and let $A$ be an Albert algebra over $F$.  For example, you could take $k$ to be the ring of integers in a number field with a real embedding, $F = \IR$, and $A$ to be the euclidean Albert algebra.

If you pick any $k$-lattice $J$ in $A$ that contains $1_J$ and is closed under the $\sharp$ operator in $A$, then by Theorem \ref{cubic.jordan} $J$ is a cubic Jordan algebra.  If furthermore the bilinear trace form on $J$ (the restriction of the bilinear form on $A$) is regular, then Theorem \ref{albert.descent} gives that $J$ is an Albert algebra over $k$ as defined in 39.18.

\section{Splitting and reducing fields for Freudenthal algebras}

In this section, we consider a Freudenthal algebra $J$ over a field $F$ and prove the existence of separable extensions $K \supseteq L \supseteq F$ of small dimension such that $J_K$ is split (as defined in 39.20) and $J_L$ is reduced (as defined in 41.1).

\subsection*{Dimensions 1 and 3}
If $\dim_F J = 1$, then $J \cong \plalg{k}$, so $J$ is split and we may take $K = F$.

If $\dim _F J = 3$, then $J \cong \plalg{E}$ for a cubic \'etale $F$-algebra $E$ that is uniquely determined up to isomorphism, see Cor.~55.2 or Exercise 39.42(b).  For such a $J$, being reduced is equivalent to being split --- i.e., $J \cong \plalg{(F \times F \times F)}$ --- which is equivalent to $E \cong F \times F \times F$.  If $E \cong F \times K$ for a field $K$, then $E \otimes_F K$ is split.  Otherwise, $E$ is a field, and the normal closure $K$ of $E$ has dimension 3 or 6 and $E \otimes_F K$ is split.  In summary, one can take $L = K$ and $K$ of dimension dividing 6.

\subsection*{Dimensions 9,  15, and 27}

\begin{lem} \label{freud.red}
Suppose $J$ is a Freudenthal $F$-algebra of rank $\ge 9$.  If $J$ is reduced, then there is a separable extension $K$ of $J$ of degree dividing 2 such that $J_K$ is split.
\end{lem}

\begin{proof}
The hypothesis that $J$ is reduced says that $J \cong \Her_3(C, \Gamma)$ for some composition algebra $C$ (Prop.~39.17) of rank $\ge 2$.  If $C$ is split, then $J$ is split (Prop.~40.6), and we take $K = F$.  Otherwise, there is a separable quadratic extension $K$ of $F$ contained in $C$ and $C \otimes_F K$ is split (Prop.~22.20, Cor.~22.16).  
\end{proof}

Because every Freudenthal algebra of dimension 15 is reduced (Th.~46.8), we obtain:
\begin{cor} 
If $J$ is a Freduenthal $F$-algebra of dimension 15, then there is a separable extension $K$ of $F$ of dimension dividing 2 such that $J_K$ is split. \Proofend
\end{cor}

\begin{prop} \label{albert.split}
Let $J$ be a Freudenthal $F$-algebra of  dimension 9 or 27.  There exist separable extensions $K \supseteq L \supseteq F$ such that $[K:L]$ divides 2, $[L:F]$ divides 3, $J_K$ is split, and $J_L$ is reduced.
\end{prop}

\begin{proof}[Proof \#1]
If $J$ is reduced, take $L = F$.  Otherwise, $J$ is a division algebra (Prop.~39.17), so it contains a a separable cubic subfield $L$ (Cor.~46.7). 
Since $L \otimes_F L$ is not a field, $\plalg{(L \otimes_F L)}$ contains a nonzero element of norm zero, ergo so does $J_L$, and we conclude that $J_L$ is reduced (Prop.~39.17).  

Applying 
Lemma \ref{freud.red} to $J_L$ provides $K$.
\end{proof}

\begin{proof}[Proof \#2 (sketch)]
To produce $K$, use Exercise 55.12 or Prop.~55.5 to translate the problem into producing a field $K$ that splits an absolutely simple algebraic group of type $\rsA_2$ or $\rsF$.  Then apply the main result of \cite{Ti:deg}.
\end{proof}

\begin{example}
We sketch an example of a field $F$ and an Albert $F$-algebra $J$ that shows that the dimensions in Prop.~\ref{albert.split} are best possible. 

Pick a field $k$.  There exists a field $F \supseteq k$ and an Albert $F$-algebra $J$ such that the Rost invariant $r(J)$ of $J$ as described in 55.10 has order 6 as an element of the abelian group  $H^3(F, \IZ/6\IZ(2))$.  Here are two ways to produce $F$ and $J$:
\begin{enumerate}[(1)]
    \item Do it explicitly, using the second Tits construction and the formulas for $r(J)$ provided in \cite{MR99j:17045}, \cite{MR1337184}, \cite{PR:char3}, and \cite[\S40]{MR2000a:16031}.

    \item Take $\bfG$ to be the automorphism group of the split Albert algebra, which is a split semi-simple group scheme of type $\rsF$.  Find a versal $\bfG$-torsor $\bfX$, which exists by \cite[p.~12]{GMS} and is defined over some extension field $F$ of $k$.  Take $J$ to be the Albert $F$-algebra corresponding to $\bfX$ by Prop.~55.5.  The fact that $r(J)$ has order 6 follows from \cite[pp.~31, 135]{GMS}.
\end{enumerate}
Suppose now that $K$ is a finite separable extension of $F$ such that $J_K$ is split.  Then the restriction $\mathrm{res}_{K/F}(r(J)) \in H^3(K, \IZ/6\IZ(2))$ is zero, consequently the corestriction satisfies
\begin{equation} \label{rescores}
[K:F]\, r(J) = \mathrm{cor}_{K/F} \, \mathrm{res}_{K/F}(r(J)) = \mathrm{cor}_{K/F}(0) = 0
\end{equation}
in $H^3(F, \IZ/6\IZ(2))$, ergo 6 divides $[K:F]$.  A similar argument shows that every separable extension $L$ of $F$ such that $J_L$ is reduced has dimension divisible by 3.
\end{example}

The argument in the example works with essentially no changes also for dimension 9. 

The restriction/corestriction argument around \eqref{rescores} gives a converse to Proposition \ref{albert.split}: \emph{If $J$ is an Albert $F$-algebra that is not reduced and $L$ is a separable extension of $F$ such that $J_L$ is reduced, then 3 divides $[L:F]$.}  (Compare Exercise 12.41.)

Using the same restriction/corestriction argument or combining Th.~41.26 and Springer's theorem on quadratic forms under odd-degree extensions \cite[Cor.~18.5, 18.6]{MR2427530} completes the converse: \emph{If $J$ is an Albert $F$-algebra in which every nilpotent element is zero and $K$ is a separable extension of $F$ such that $J_K$ does contain nonzero nilpotents, then 2 divides $[K:F]$.}

\subsection*{Dimension 6}
It remains to treat the case of a Freudenthal algebra $J$ of rank 6.  We will assume that $J$ is regular, equivalently (by Cor.~39.15) that $\car F \ne 2$.
By Th.~46.8, any such $J$ is reduced, i.e., of the form $\Her_3(F, \Gamma)$ for some $\Gamma = \diag( \gamma_1, \gamma_2, \gamma_3 ) \in \bfGL_3(F)$.  The algebra $J$ that is defined to be split in 39.20 has $\Gamma = \diag(1,1,1)$, whereas the $J$ that has $\bfAut(J)$ split has $\Gamma = \diag(1, -1, -1)$ (Prop.~55.6).  Example 39.35 provides an example of a $J$ that is not split by any separable extension of $F$.  Consequently, rather than splitting $J$, we instead split $\bfAut(J)$.

\begin{lem}
Let $J$ be a Freudenthal $F$-algebra of dimension 6 and suppose that $\car F \ne 2$.  Then there exists a separable extension $K$ of $F$ of dimension dividing 2 such that $J_K \cong \Her_3(F, \diag(1, -1, -1))$.
\end{lem}

\begin{proof}
There is a 3-dimensional regular quadratic form $Q_J$ defined in (41.5.1).  Adjoining a square root (which would be a separable extension) is sufficient to make $Q_J$ isotropic, which implies the necessary isomorphism by Th.~41.26.
\end{proof}


\section{Describing Albert algebras by the Tits constructions}

The classical treatment of the two Tits constructions as set forth in Jacobson \cite[IX.12]{MR0251099} and McCrimmon \cite{MR0238916, MR0271181}, culminates in the following result (\cite[Thm.~IX.22]{MR0251099}, \cite[Thm.~8]{MR0271181}): If $J$ is an Albert algebra over a field $F$ and $A$ is a central simple associative $F$-algebra of degree $3$ making $A^{(+)}$ a subalgebra of $J$, then there exists a scalar $\lambda \in F^\times$ such that the inclusion $A^{(+)} \hookrightarrow J$ can be extended to an isomorphism from the first Tits construction $J(A,\lambda)$ onto $J$. As an immediate corollary, basically the same conclusion can be drawn after replacing $A^{(+)}$ by the Jordan algebra of symmetric elements of a central simple associative $F$-algebra of degree $3$ with unitary involution.

These results are quick to derive, given the approach to the Tits constructions as described in Chap.~VII, so we give a proof here. 
The approach adopted here is different from that taken by Jacobson and McCrimmon.  They relied on properties of the special universal envelope of $A^{(+)}$, where $A$ is a central simple associative algebra of degree 3, whereas we combine results from the book with \cite{MR0038335}.

Throughout we let $F$ be an arbitrary field.


\begin{prop} \label{p.UNITIN}
For a central simple associative $F$-algebra $(B,\tau)$ of degree $3$ with unitary involution, the following conditions are equivalent. 
\begin{enumerate}[(i)]
\item [\emph{(i)}] The center of $B$ is not a field.

\item [\emph{(ii)}] $B$ is not simple.

\item [\emph{(iii)}] There exists a central simple associative $F$-algebra $A$ of degree $3$ such that $H(B,\tau) \cong A^{(+)}$.
\end{enumerate}
\end{prop}

\begin{proof}
Put $K := \Cent(B)$, a quadratic \'etale $F$-algebra.

(i)$\Rightarrow$(ii). The center of a unital simple algebra is a field (Cor.~9.20).

(ii)$\Rightarrow$(iii). Since $(B,\tau)$ is simple but $B$ is not, Prop.~10.5 yields a simple associative $F$-algebra $A$ having $(B,\tau) \cong (A \times A^{\op},\vep_A)$ where $\vep_A$ is the exchange involution. We conclude $H(B,\tau) \cong A^{(+)}$ from (29.8.1) and with $L := \Cent(A)$ have $L \times L = K$, hence $L = F$. Thus $A$ is central, while a dimension count shows that $A$ has degree $3$.

(iii)$\Rightarrow$(i). Suppose (iii) and not (i), i.e., that $K$ is a field. From Exc.~44.29~(a) we deduce $A_K^{(+)} \cong H(B,\tau)_K = B^{(+)}$, where $A_K$ and then $B$ are both central simple over $K$ (Cor.~9.23~(a)). By \cite[Cor. of Thm.~21]{MR0038335}, therefore, either $A_K \cong B$ or $A_K^{\op} \cong B$. Replacing $A$ by $A^{\op}$ if necessary, we may assume $A_K \cong B$. The action of $\tau$ on $A_K$ induced by any such isomorphism fixes $A$ since $A^{(+)} \cong H(B,\tau)$ and is the conjugation on $K$, hence an automorphism on all of $A_K$, a contradiction to $\tau$ being an anti-automorphism on $B$.
\end{proof}

\begin{cor} \label{c.CENTIN}
For $i = 1,2$, let $(B_i,\tau_i)$ be central simple associative $F$-algebras of degree $3$ with unitary involutions and $K_i := \Cent(B_i)$. If $H(B_1,\tau_1) \cong H(B_2,\tau_2)$, then $K_1 \cong K_2$ over $F$.
\end{cor}

\begin{proof}
If $K_1$ is split, then Prop.~\ref{p.UNITIN} yields a central simple associative $F$-algebra $A_1$ of degree $3$ such that $H(B_2,\tau_2) \cong A_1^{(+)}$. Hence $K_2$ is split as well. By symmetry, we may therefore assume that $K_1,K_2$ are both fields. Then 
\[
H((B_1)_{K_2}, (\tau_1)_{K_2}) \cong B_2^{(+)},
\]
so $(K_1)_{K_2}$ by Prop.~\ref{p.UNITIN} is a split quadratic \'etale $K_2$-algebra. By Exc.~23.40, this implies $K_1 \cong K_2$ over $F$.  
\end{proof}

\begin{lem} \label{l.INVSYS}
Let $\msB = (K,B,\tau,p)$ be an involutorial system of the second kind over $F$ such that $H(\msB)$ is a Freudenthal $F$-algebra of dimension $9$. Then $(B,\tau)$ is a central simple associative $F$-algebra of degree $3$ with unitary involution and $K = \Cent(B)$.
\end{lem}

\begin{proof}
$B^{(+)} = H(\msB)_K$ (by Exc.~44.29~(a)) is a Freudenthal $K$-algebra of constant rank $9$ as a $K$-module. If $K$ is a field, we therefore conclude that $B^{(+)}$ is simple. Hence so is $B$, and Exc.~42.21 implies that $B$ is, in fact,  a central simple associative $K$-algebra of degree $3$ on which $\tau$ acts as a $K/F$-involution. On the other hand, if $K = F \times F$ is split, Exc.~44.31~(a) yields a pointed cubic alternative $F$-algebra $(A,q)$ such that $B = A \times (A^q)^{\op}$, $\tau = \vep_A$ is the switch and $p = (q,q)$. Letting $K$ act on $F$ through the first factor, we obtain $F \in \Kalg$ and, for the same reason as before, $A = B_F$ is a central simple associative $F$-algebra of degree $3$. In both cases, therefore, $(B,\tau)$ is a central simple associative $F$-algebra of degree $3$ with unitary involution having $K = \Cent(B)$.  \end{proof}

We can now prove the results announced at the beginning of this section, albeit in the reverse order.

\begin{thm} \label{t.SETIAL}
Let $J$ be an Albert $F$-algebra and $(B,\tau)$ a central simple associative algebra of degree $3$ with unitary involution over $F$ making $J_0 := H(B,\tau)$  a Freudenthal subalgebra of $J$. Then there exists an admissible scalar $(p,\mu)$ for $(B,\tau)$ such that the inclusion $J_0 \hookrightarrow J$ can be extended to an isomorphism from the second Tits construction $J(B,\tau,p,\mu)$ onto $J$.   
\end{thm}

\begin{proof}
We first assume that $F$ is finite. Then $J$ is split (Exc.~40.17~(b)). But $J^\prime := J(B,\tau,1,1)$ is a split Albert $F$-algebra as well, so there is an isomorphism $\Psi\:J^\prime \overset{\sim}\to J$, and this isomorphism maps the subalgebra $J_0 \subseteq J^\prime$ onto a subalgebra $J_0^\prime \subseteq J$. Since $F$ is finite, both $J_0$ and $J^\prime_0$ are reduced simple Freudenthal subalgebras of $J$, and $\Psi$ induces an isomorphism $\vph\:J_0 \overset{\sim}\to J_0^\prime$. By the Skolem-Noether theorem (Exc.~41.34), $\vph$ can be extended to an automorphism $\Phi$ of $J$. Hence $\Phi^{-1}\circ \Psi\:J^\prime \overset{\sim}\to J$ is an isomorphism of the desired kind.

Next assume that $F$ is infinite. Applying Thm.~45.10, we find an \'etale element of $J$ relative to $J_0$\footnote{The full force of the difficult Thm.~45.10 will not be needed here. Instead, by a Zariski density argument, we are reduced to the case that $F$ is algebraically closed. Then $J = \Her_3(C)$, $C := \Zor(F)$, is split, and, again by the Skolem-Noether theorem, we may assume $J_0 = \Her_3(D)$, $D$ being the diagonal of $\Zor(F)$. Since $D^\perp \subseteq \Zor(F)$ is a free right $D$-module of rank $3$, \'etale elements of $J$ relative to $J_0$ exist, thanks to Exc.~45.17~(b).}. Hence Cor.~44.17 yields an involutorial system  $\msA = (K,A,\sigma,p)$ of the second kind over $F$ as well as an admissible scalar $\mu \in K^\times$ such that
\begin{align}
\label{HASI} H(A,\sigma) = H(\msA) = H(B,\tau)   
\end{align}
and the inclusion $J_0 \hookrightarrow J$ can be extended to an isomorphism from $J(\msA,\mu)$ onto $J$. By Lemma~\ref{l.INVSYS}, $(A,\sigma)$ is a central simple associative $F$-algebra of degree $3$ with unitary involution and $K = \Cent(A)$. From Cor.~\ref{c.CENTIN} we deduce $K = \Cent(B)$, while \eqref{HASI} implies
\[
A^{(+)} = H(A,\sigma) \otimes_F K = H(B,\tau) \otimes_F K = B^{(+)},
\]
hence $A = B$ by \cite[Cor. of Thm.~21]{MR0038335}, after replacing $A$ by $A^{\op}$ if necessary. But $\sigma$ acts on $H(A,\sigma) = H(B,\tau)$ as the identity and on $K$ by conjugation. Hence $\sigma = \tau$, and the proof is complete.
\end{proof}

\begin{cor} \label{c.FITIAL}
Let $J$ be an Albert $F$-algebra and $A$ a central sinple associative algebra of degree $3$ over $F$ making $A^{(+)}$ a subalgebra of $J$. Then there exists a $\lambda \in F^\times$ such that the inclusion $A^{(+)} \hookrightarrow J$ can be extended to an isomorphism from the first Tits construction $J(A,\lambda)$ onto $J$.
\end{cor}

\begin{proof}
Put $(B,\tau) := (A \times A^{\op},\vep_A)$, $\vep_A$ being the exchange involution, and identify $A^{(+)} = H(B,\tau)$ canonically. By Thm.~\ref{t.SETIAL}, there exists an admissible scalar $(p,\mu)$ for $(B,\tau)$ such that the inclusion $A^{(+)} \hookrightarrow J$ can be extended to an isomorphism $\Psi\:J(B,\tau,p,\mu)\overset{\sim}\to J$. Write $p = (q,q)$ with $q \in A^\times$, and $\mu = (\lambda,\lambda^\prime)$ with $\lambda,\lambda^\prime \in F^\times$. Applying Thm.~44.19, we find an isomorphism $\Phi\:J(B,\tau,p,\mu)\overset{\sim}\to J(A,\lambda)$ extending the identity of $A^{(+)}$. Thus $\Psi\circ \Phi^{-1}\:J(A,\lambda) \overset{\sim}\to J$ is an isomorphism of the desired kind.    
\end{proof}
    \renewcommand{\theequation}{s\arabic{equation}}

\solnchap{Solutions for Chapter~\ref{c.PROLOGUE}}

\solnsec{Section~\ref{s.GRACA}}

\begin{sol}{pr.RiceBrown} 
Suppose $A$ is a finite-dimensional real algebra of odd dimension $> 1$.  Then there exists $x, y \in A$ that are linearly independent.  Put $L_x$ for the linear transformations that is left multiplication by $x$ and similarly for $y$.  For $\bft$ an indeterminate, $\det(L_x + \bft L_y)$ is a polynomial in $\IR[\bft]$ of odd degree $\dim A$, and therefore it has a nontrivial root in $\IR$.  That is, there is some $\alpha \in \IR$ such that $L_{x + \alpha y} = L_x + \alpha L_y$ is not invertible.  Since $x + \alpha y \ne 0$, $A$ is not division.

(This observation, under the additional hypothesis that $A$ is unital, can be found in the entertaining paper \cite{RiceBrown:cant}.)
\end{sol}

\begin{sol}{pr.ASSGC} \label{sol.ASSGC} Putting $[x_1,x_2,x_3] =: (a,u)$ for some $a \in \IC$, $u \in \IC^3$ and using (\ref{ss.DEFGC}.\ref{PROGC}), we obtain
\begin{align*}
(a,u) =\,\,& [(a_1,u_1),(a_2,u_2),(a_3,u_3)] \\
=\,\,&\big((a_1,u_1)(a_2,u_2)\big)(a_3,u_3) - (a_1,u_1)\big((a_2,u_2)8a_3,u_3)\big) \\
=\,\,&\big((a_1a_2 - \bar u_1^\trans u_2,u_2\bar a_1 + u_1a_2 + \bar u_1 \times \bar u_2)\big)(a_3,u_3) \\
\,\,&- (a_1,u_1)\big((a_2a_3 - \bar u_2^\trans u_3,u_3\bar a_2 + u_2a_3 + \bar u_2 \times \bar u_3)\big),
\end{align*}
and an application of (\ref{ss.CROPRO}.\ref{DECRO}) combined with the Grassmann identity (\ref{ss.CROPRO}.\ref{GRAS}) yields
\begin{align*}
a =\,\,&a_1a_2a_3 -\bar u_1^\trans u_2a_3 - a_1\bar u_2^\trans u_3 - \bar a_2\bar u_1^\trans u_3 - (u_1 \times u_2)^\trans u_3 - a_1a_2a_3 \\
\,\,&+ a_1\bar u_2^\trans u_3 + \bar u_1^\trans u_3\bar a_2 +\bar u_1^\trans u_2a_3 + \overline{u_1^\trans (u_2 \times u_3)} \\
=\,\,&\overline{\det(u_1,u_2,u_3)} - \det(u_1,u_2,u_3), \\
x =\,\,&u_3\overline{a_1a_2} - u_3\bar u_2^\trans u_1 + u_2\bar a_1a_3 + u_1a_2a_3 + (\bar u_1 \times \bar u_2)a_3 + (\bar u_2 \times \bar u_3)a_1 \\
\,\,&+ (\bar u_1 \times \bar u_3)\bar a_2 + (u_1 \times u_2) \times \bar u_3  - u_3\bar a_2\bar a_1 -u_2a_3\bar a_1 - (\bar u_2 \times \bar u_3)\bar a_1 \\
\,\,&- u_1a_2a_3 + u_1\bar u_2^\trans u_3 - (\bar u_1 \times \bar u_3)a_2 - (\bar u_1 \times \bar u_2)\bar a_3 - \bar u_1 \times (u_2 \times u_3) \\
=\,\,&-u_3\bar u_2^\trans u_1 + (\bar u_1 \times \bar u_2)(a_3 - \bar a_3) + (\bar u_2 \times \bar u_3)(a_1 - \bar a_1) + (\bar u_3 \times \bar u_1)(a_2 - \bar a_2) \\
\,\,&+ (u_1 \times u_2) \times \bar u_3 + u_1\bar u_2^\trans u_3 + (u_2 \times u_3) \times \bar u_1 \\
=\,\,&\sum (\bar u_i \times \bar u_j)(a_k - \bar a_k) + \sum (u_i \times u_j) \times \bar u_k. 
\end{align*}
Thus the associator formula holds, While the $\IC$-component of the right-hand side is trivially alternating in $x_1,x_2,x_3$, so is its $\IC^3$-component since it is the sum over all cyclic permutations of $(123)$ of a trilinear expression, namely, $(u_1 \times u_2) \times \bar u_3 + (\bar u_1 \times \bar u_2)(a_3 - \bar a_3)$, that vanishes for $x_1 = x_2$.
\end{sol}

\begin{sol}{pr.DIVO} \label{sol.DIVO} Let $x = (a,u), y = (b,v)$ with $a,b \in \IC$, $u,v \in \IC^3$ and assume $x \neq 0 = xy$. We must show $y = 0$. First of all, (\ref{ss.DEFGC}.\ref{PROGC}) implies
\begin{align}
\label{ZEDI}  ab = \bar u^\trans v, \quad v\bar a + ub + \bar u \times \bar v = 0.
\end{align}
Multiplying the second equation with $\bar u^\trans$ from the left and applying the first combined with (\ref{ss.CROPRO}.\ref{DECRO}), (\ref{ss.NOTR}.\ref{EUNGC}), we obtain
\[
0 =\bar u^\trans (\bar u \times \bar v) + \bar u^\trans ub + \bar u^\trans v\bar a = (\|u\|^2 + \|a\|^2)b = \|x\|^2b,
\]
hence $b = 0$. Thus \eqref{ZEDI} reduces to
\begin{align}
\label{ZEDIB} u^\trans \bar v = 0, \quad \bar u \times \bar v = -v\bar a.
\end{align}
Here the second equation gives $0 = \bar v^\trans (\bar u \times \bar v) = -\|v\|^2\bar a$, so if $a \neq 0$ then $v = 0$. Hence we are left with the case $a = 0$, which implies $u \neq 0$. On the other hand, \eqref{ZEDIB} yields $(\bar u \times \bar v) \times u = \bar vu^\trans \bar u - \bar uu^\trans \bar v = \|u\|^2\bar v$, hence again $v = 0$.
\end{sol}

\begin{sol}{pr.INVAS}\label{sol.INVAS} a) For $a,b \in \IC$, $u,v \in \IC^3$, $x = (a,u), y = (b,v) \in \IO$, we apply (\ref{ss.DEFGC}.\ref{PROGC}), (\ref{ss.NOTR}.\ref{CONGC}) and conclude
\begin{align*}
\overline{xy} =\,\,&\overline{(ab - \bar u^\trans v,v\bar a + ub + \bar u \times \bar v)} = (\overline{ab - \bar u^\trans v},-v\bar a -ub - \bar u \times \bar v) \\
=\,\,&\big(\bar b\bar a - \overline{(-v)}^\trans (-u),(-u)\bar{\bar b} + (-v)\bar a + (-\bar v) \times (-\bar u)\big) \\
=\,\,&\big(\bar b,-v)\big)\big(\bar a,-u\big) = \bar y\bar x,
\end{align*}
which proves (a). 

(b)  Let $x = (a,u)$, $y = (b,v)$, with $a,b \in \IC$, $u,v \in \IC^3$. From (\ref{ss.DEFGC}.\ref{PROGC}) and (\ref{ss.DEFGC}.\ref{TRAGC}) we conclude that $t_\IO(xy) = ab + \bar a\bar b - \bar u^\trans v - u^\trans \bar v = ab + \bar a\bar b - \bar u^\trans v - \bar v^\trans u$ is symmetric in $x,y$, giving the first equation. Moreover, the  $\IC$-component of the associator $[x_1,x_2,x_3]$ in Exc.~\ref{pr.ASSGC} is purely imaginary. Hence the second equation follows from (\ref{ss.NOTR}.\ref{TRAGC}). Combined with (\ref{ss.NOTR}.\ref{TRANO}) and (a) it implies
\begin{align*}
n_\IO(xy,z) =\,\, &t_\IO\big((xy)\bar z\big) = t_\IO\big(x(y\bar z)\big) = n_\IO(x,\overline{y\bar z}) = n_\IO(x,z\bar y),
\end{align*}
hence the first of the remaining equations. The second one follows analogously.

(c) The first set of equations follows immediately from (\ref{ss.NOTR}.\ref{QUAGC}) and the definition of the conjugation. As to the remaining ones, we apply alternativity, (\ref{ss.NOTR}.\ref{BLAGC}), (\ref{ss.NOTR}.\ref{QUAGC}), (\ref{ss.NOTR}.\ref{COGC}) to deduce
\begin{align*}
xyx =\,\,&x(x \circ y) - x^2y = t_\IO(x)xy + t_\IO(y)x^2 - n_\IO(x,y)x - t_\IO(x)xy + n_\IO(x)y \\
=\,\,&t_\IO(x)t_\IO(y)x - n_\IO(x)t_\IO(y)1_\IO - n_\IO(x,y)x + n_\IO(x)y \\
=\,\,&n_\IO\big(x,t_\IO(y)1_\IO - y\big)x - n_\IO(x)\big(t_\IO(y)1_\IO - y\big) = n_\IO(x,\bar y)x - n_\IO(x)\bar y,
\end{align*}
which completes the proof of (c).
\end{sol}

\begin{sol}{pr.MIMOU}\label{sol.MIMOU}
We begin by proving \eqref{CROMA}. If $\tr$ stands for the trace form of $\Mat_3(\IC)$, then for any $A \in \Mat_3(\IC)$, we find a scalar $\lambda(A) \in \IC$ such that
\begin{align}
\label{CROMAT} \tr\Big(\big((u \times v)w^\trans  + (v \times w)u^\trans  + (w \times u)v^\trans \big)A\Big) = \lambda(A)\det (u,v,w) &&(u,v,w \in \IC^3)
\end{align} 
since the left-hand side is an alternating trilinear function of its arguments $u,v,w \in \IC^3$. Specializing $u := e_1$, $v := e_2$, $w := e_3$ and observing $e_ie_i^\trans  = e_{ii}$ in terms of the ordinary matrix units for $1 \leq i \leq 3$, we deduce $\lambda(A) = \tr(A)$, hence
\[
\tr\Big(\big((u \times v)w^\trans  + (v \times w)u^\trans  + (w \times u)v^\trans \big)A\Big) = \tr\Big(\big(\det(u,v,w)\Eins_3\big)A\Big).
\]
Since the symmetric bilinear form $(X,Y) \mapsto \tr(XY)$ on $\Mat_3(\IC)$ is easily seen to be non-degenerate, \eqref{CROMA} holds.

It remains to prove the Moufang identities. Since the conjugation of $\IO$ is an algebra involution (Exc.~\ref{pr.INVAS}~(a)), the third Moufang identity immediately follows from the first.

We begin by proving the first Moufang identity, which by Exc.~\ref{pr.INVAS}~(c) comes down to showing
\begin{align}
\label{FIEX.MIMOU} x\big(y(xz)\big) = n_\IO(x,\bar y)xz - n_\IO(x)\bar yz. 
\end{align}
Since this relation is bilinear in $(y,z)$, it suffices to put $x = a \oplus u$, $a \in \IC$, $u \in \IC^3$, and to consider the following cases.

\case{1}
$y = (b,0)$, $z = (c,0)$, $b,c \in \IC$. Then  (\ref{ss.DEFGC}.\ref{PROGC}) yields
\begin{align*}
x\big(y(xz)\big) =\,\,&x\Big(y\big((a,u)(c,0)\big)\Big) = x\big((b,0)(ac,uc)\big) = (a,u)(abc,u\bar bc) \\
=\,\,&\big(a^2bc - \bar u^\trans u\bar bc,u(ab + \bar a\bar b)c\big).
\end{align*}
On the other hand, applying (\ref{ss.NOTR}.\ref{NOGC}), (\ref{ss.NOTR}.\ref{BLINOGC}) and (\ref{ss.NOTR}.\ref{CONGC}) we obtain
\begin{align*}
n_\IO(x,\bar y)xz - n_\IO(x)\bar yz =\,\,&n_\IO\big((a,u),(\bar b,v)\big)(a,u)(c,0) - n_\IO\big((a,u)\big)(\bar b,0)(c,0) \\
=\,\,&(\bar a\bar b + ab)(ac,uc) - (\bar aa + \bar u^\trans u)(\bar bc,0) \\
=\,\,&\big(\bar aa\bar bc + a^2bc - \bar aa\bar bc - \bar u^\trans u\bar bc,u(abc + \bar a\bar bc)\big) \\
=\,\,&\big(a^2bc - \bar u^\trans u\bar bc,u(ab + \bar a\bar b)c\big).
\end{align*}
Thus \eqref{FIEX.MIMOU} holds. 

\case{2} 
$y = (b,0)$, $z = (0,w)$, $b \in \IC$, $w \in \IC^3$. Basically arguing as before, we obtain
\begin{align*}
x\big(y(xz)\big) =\,\,&x\Big(y\big((a,u)(0,w)\big)\Big) = x\big((b,0)(-\bar u^\trans w,w\bar a + \bar u \times \bar w)\big) \\
=\,\,&(a,u)\big(-b\bar u^\trans w,w\bar a\bar b + (\bar u \times \bar w)\bar b\big) \\
=\,\,&\big(-ab\bar u^\trans w - \bar u^\trans w\bar a\bar b -\bar u^\trans (\bar u \times \bar w)\bar b, \\
\,\,&\big(w\bar a^2\bar b + (\bar u \times \bar w)\bar a\bar b - ub\bar u^\trans w + (\bar u \times \bar w)ab + \bar u \times (u \times w)b\big). 
\end{align*}
Here (\ref{ss.CROPRO}.\ref{DECRO}) and the Grassmann identity (\ref{ss.CROPRO}.\ref{GRAS}) imply $\bar u^\trans (\bar u \times \bar w) = 0$ as well as $\bar u \times (u \times w) = (w \times u) \times \bar u = u\bar u^\trans w - w\bar u^\trans u$, hence
\begin{align*}
x\big(y(xz)\big) =\,\,&-\bar u^\trans w(ab + \bar a\bar b) \oplus \big(w\bar a^2\bar b + (\bar u \times \bar w)(ab + \bar a\bar b) - w(\bar u^\trans u)b\big) 
\end{align*} 
On the other hand,
\begin{align*}
n_\IO(x,\bar y)xz - n_\IO(x)\bar yz =\,\,&n_\IO\big((a,u),(\bar b,0)\big)(a,u)(0,w) - n_\IO\big((a,u)\big)(\bar b,0)(0,w) \\
=\,\,&(\bar a\bar b + ab)(-\bar u^\trans w,w\bar a + \bar u \times \bar w) - (\bar aa + \bar u^\trans u)(0,wb) \\
=\,\,&-\big((\bar u^\trans w)(ab + \bar a\bar b),w(\bar aab + \bar a^2\bar b) \\
\,\,&+ (\bar u \times \bar w)(ab + \bar a\bar b) - w\bar aab - w\bar u^\trans ub\big) \\
=\,\,&\big(-(\bar u^\trans w)(ab + \bar a\bar b),w\bar a^2\bar b + (\bar u \times \bar w)(ab + \bar a\bar b) - w\bar u^\trans ub\big)
\end{align*}
Comparing with the preceding expression for $x(y(xz))$ settles Case~$2$. 

\case{3} 
$y = (0,v)$, $z = (c,0)$, $v \in \IC^3$, $c \in \IC$. Then
\begin{align*}
x\big(y(xz)\big) =\,\,& x\Big(y\big((a,u)(c,0)\big)\Big) = x\big((0,v)(ac,uc)\big) = (a,u)\big(-\bar v^\trans uc,vac + (\bar v \times \bar u)\bar c\big) \\
=\,\,&\big(-a\bar v^\trans uc - \bar u^\trans vac - \bar u^\trans (\bar v\times \bar u)\bar c, \\
\,\,&v\bar aac + (\bar v \times \bar u)\bar a\bar c - u\bar v^\trans uc + (\bar u \times \bar v)\bar a\bar c + \bar u \times (v \times u)c\big).
\end{align*}
Here $\bar u^\trans (\bar v \times \bar u) = 0$ and $\bar u \times (v \times u) = (u \times v) \times \bar u = v\bar u^\trans u - u\bar u^\trans v$, which implies
\begin{align*}
x\big(y(xz)\big) =\,\,&\big(-ac(\bar u^\trans v + \bar v^\trans u),v(\bar aa + \bar u^\trans u)c - u(\bar u^\trans v + \bar v^\trans u)c\big).
\end{align*}
On the other hand,
\begin{align*}
n_\IO(x,\bar y)xz - n_\IO(x)\bar yz =\,\,&-n_\IO\big((a,u),(0,v)\big)(a,u)(c,0) + n_\IO\big((a,u)\big)(0,v)(c,0) \\
=\,\,&-(\bar u^\trans v + \bar v^\trans u)(ac,uc) + (\bar aa + \bar u^\trans u)(0,vc) \\
=\,\,&\big(-ac(\bar u^\trans v + \bar v^\trans u),-u(\bar u^\trans v + \bar v^\trans u)c + v(\bar aa + \bar u^\trans u)c\big),
\end{align*}
and the discussion of Case~$3$ is complete. 

\case{4}
$y = (0,v)$, $z = (0,w)$, $v,w \in \IC^3$. This is the most delicate case of them all. 
\begin{align*}
x\big(y(xz)\big) =\,\,&x\Big(y\big((a,u)(0,w)\big)\Big) = x\big((0,v)(-\bar u^\trans w,w\bar a + \bar u \times \bar w)\big) \\
=\,\,&(a,u)\big(-\bar v^\trans w\bar a - \bar v^\trans (\bar u \times \bar w),-v\bar u^\trans w + (\bar v \times \bar w)a + \bar v \times (u \times w)\big).
\end{align*}
Since $\bar v \times (u \times w) = (w \times u) \times \bar v = u\bar v^\trans w - w\bar v^\trans u$, we conclude
\begin{align*}
x\big(y(xz)\big) =\,\,&(a,u)\big(-\bar v^\trans w\bar a - \bar v^\trans (\bar u \times \bar w),u\bar v^\trans w - w\bar v^\trans u - v\bar u^\trans w + (\bar v \times \bar w)a\big) \\
=\,\,&\big(-\bar aa\bar v^\trans w - a\bar v^\trans (\bar u \times \bar w) - \bar u^\trans u\bar v^\trans w + \bar u^\trans w\bar v^\trans u + \bar u^\trans v\bar u^\trans w - \bar u^\trans (\bar v \times \bar w)a\big), \\
\,\,&u\bar v^\trans w\bar a - w\bar v^\trans u\bar a - v\bar u^\trans w\bar a + (\bar v \times \bar w)\bar aa - u\bar v^\trans w\bar a - u\bar v^\trans (\bar u \times \bar w) \\
\,\,&+ (\bar u \times \bar u)v^\trans \bar w - (\bar u \times \bar w)v^\trans \bar u - (\bar u \times \bar v)u^\trans \bar w + \bar u \times (v \times w)\bar a\big).
\end{align*}
Here $\bar u \times (v \times w) = (w \times v) \times \bar u = v\bar u^\trans w - w\bar u^\trans v$, so that we obtain
\begin{align}
\label{LHLM} x\big(y(xz)\big) =\,\,&\big(\bar u^\trans v\bar u^\trans w + \bar u^\trans w\bar v^\trans u - \bar u^\trans u\bar v^\trans w - \bar aa\bar v^\trans w,-w\bar v^\trans u\bar a - w\bar u^\trans v\bar a \\
\,\,&+ (\bar v \times \bar w)\bar aa - u\bar v^\trans (\bar u \times \bar w) - (\bar u \times \bar v)u^\trans \bar w - (\bar u \times \bar w)v^\trans \bar u\big). \notag
\end{align}
On the other hand,
\begin{align*}
n_\IO(x,\bar y)xz - n_\IO(x)\bar yz =\,\,&-n_\IO\big((a,u),(0,v)\big)(a,u)(0,w) + n_\IO\big((a,u)\big)(0,v)(0,w) \\
=\,\,&-(\bar u^\trans v + \bar v^\trans u)(-\bar u^\trans w,w\bar a + \bar u \times \bar w) + (\bar aa + \bar u^\trans u)(-\bar v^\trans w,\bar v \times \bar w) \\
=\,\,&\big(\bar u^\trans v\bar u^\trans w + \bar v^\trans u\bar u^\trans w - \bar aa\bar v^\trans w - \bar u^\trans u\bar v^\trans w,-w\bar u^\trans v\bar a - w\bar v^\trans u\bar a \\
\,\,&- (\bar u \times \bar w)\bar u^\trans v - (\bar u \times \bar w)\bar v^\trans u + (\bar v \times \bar w)\bar aa + \bar u^\trans u(\bar v \times \bar w)\big).
\end{align*}
Comparing this with \eqref{LHLM} we see that the $\IC$-components are the same. Since $u^\trans \bar w = \bar ww^\trans u$ and $v^\trans \bar u = \bar u^\trans v$, equality of the $\IC^3$-components is equivalent to
\begin{align*}
-u\bar v^\trans (\bar u \times \bar w) - (\bar u \times \bar v)\bar w^\trans u - (\bar u \times \bar w)\bar u^\trans v = -(\bar u \times \bar w)\bar u^\trans v - (\bar u \times \bar w)\bar v^\trans u + (\bar v \times \bar w)\bar u^\trans u,
\end{align*}
which in turn is equivalent to 
\[
\big(\det(\bar u,\bar v,\bar w)\Eins_3\big)u = -u\det(\bar v,\bar u,\bar w) = -u\bar v^\trans (\bar u \times \bar w) = \big((\bar u \times \bar v)\bar w^\trans  + (\bar v \times \bar w)\bar u^\trans  + (\bar w \times \bar u)\bar v^\trans \big)u.
\]
But this relation holds by \eqref{CROMA}, and the discussion of Case~$4$ is complete.  

We now turn to the second Moufang identity. Arguing as before, it will be enough to consider the following cases. 

\case{1}
$y = (b,0)$, $z = (c,0)$, $b,c \in \IC$. Then (\ref{ss.DEFGC}.\ref{PROGC}) yields
\begin{align*}
(xy)(zx) =\,\,&\big((a,u)(b,0)\big)\big((c,0)(a,u)\big) = (ab,ub)(ca,u\bar c) \\
=\,\,&(a^2bc - \bar b\bar c\bar u^\trans u,u\bar a\bar b\bar c+ uabc), \\
x(yz)x =\,\,&\big((a,u)(bc,0)\big)x = (abc,ubc)(a,u) \\
=\,\,&(a^2bc - \overline{bc}\bar u^\trans u,u\overline{abc} + uabc),
\end{align*}
hence the assertion. 

\case{2}
$y = (b,0)$, $z = (0,w)$, $b \in \IC$, $w \in \IC^3$. Then (\ref{ss.DEFGC}.\ref{PROGC}) and (\ref{ss.CROPRO}.\ref{DECRO}) yield
\begin{align*}
(xy)(zx) =\,\,&\big((a,u)(b,0)\big)\big((0,w)(a,u)\big) = (ab,ub)(-\bar w^\trans u,wa + \bar w \times \bar u) \\
=\,\,&\big(-ab\bar w^ tu - a\bar b\bar u^ tw - \bar b\bar u^\trans (\bar w \times \bar u), \\
\,\,&wa\bar a\bar b + (\bar w \times \bar u)\bar a\bar b - ub\bar w^ tu + (\bar u \times \bar w)\bar a\bar b + \bar u \times (w \times u)\bar b\big) \\
=\,\,&\big(-ab\bar w^\trans u - a\bar b\bar u^\trans w,wa\bar a\bar b - ub\bar w^\trans u + \bar u \times (w \times u)\bar b\big), \\
x(yz)x =\,\,&\Big((a,u)\big((b,0)(0,w)\big)\Big)x = \big((a,u)(0,w\bar b)\big)x \\
=\,\,&\big(-\bar b\bar u^\trans w,w\bar a\bar b + (\bar u \times \bar w)b\big)(a,u) \\
=\,\,&\big(-a\bar b\bar u^\trans w - ab\bar w^\trans u - \bar b(u \times w)^\trans u, \\
\,\,&-ub\bar w^\trans u + wa \bar a\bar b + (\bar u \times \bar w)ab + (\bar w \times \bar u)ab + (u \times w) \times \bar u\bar b\big) \\
=\,\,&\big(-a\bar b\bar u^\trans w - ab\bar w^\trans u,-ub\bar w^\trans u + wa\bar a\bar b + \bar u \times (w \times u)\bar b\big),
\end{align*} 
and a comparison gives the assertion.

Since $\iota_\IO$ is an algebra involution by Exc.~\ref{pr.INVAS}~(a), the case $y = (0,v)$, $z = (c,0)$ with $v \in \IC^ 3$, $c \in \IC$ immediately reduces to Case~$2$. Thus we are left with 

\case{3}
$y = (0,v)$, $z = (0,w)$, $v,w \in \IC^3$. Then we put
\begin{align}
\label{DEMO} (xy)(zx) = (a_0,u_0), \quad x(yz)x = (a_1,u_1) &&(a_i \in \IC,\;u_i \in \IC^3,\;i = 1,2)
\end{align}
and must show $a_0 = a_1$, $u_0 = u_1$. Since
\begin{align}
\label{LEFA} xy =\,\,&(a,u)(0,v) = (-\bar u^\trans v,v\bar a + \bar u \times \bar v), \\
\label{REFA} zx =\,\,&(0,w)(a,u) = (-\bar w^\trans u,wa + \bar w \times \bar u),
\end{align}
we conclude
\begin{align*}
a_0 =\,\,&\bar u^\trans v\bar w^\trans u - a^2\bar v^\trans w - a\bar v^\trans (\bar w \times \bar u) - a(u \times v)^\trans w - (u \times v)^\trans (\bar w \times \bar u), 
\end{align*}
where the last summand may be written as
\begin{align*}
(u \times v)^\trans (\bar w \times \bar u) =\,\,&\bar u^\trans \big((u \times v) \times \bar w\big) =\bar u^\trans (v\bar w^\trans u - u\bar w^\trans v) = \bar u^\trans v\bar w^\trans u - \bar u^\trans u\bar w^\trans v.
\end{align*}
Thus
\begin{align}
\label{ANU} a_0 =\,\,&\bar u^\trans u\bar w^\trans v - a^2\bar v^\trans w - a\big(\det(u,v,w) + \overline{\det(u,v,w)}\,\big).
\end{align}
Using \eqref{LEFA}, \eqref{REFA} to compute the $\IC^3$-component of $(xy)(zx)$, we obtain
\begin{align*}
u_0 =\,\,&-wa\bar v^\trans u - (\bar w \times \bar u)\bar v^\trans u - v\bar a\bar w^\trans u - (\bar u \times \bar v)\bar w^\trans u + (\bar v \times \bar w)a\bar a \\
\,\,&+ \bar v \times(w \times u)a + (u \times v) \times \bar w\bar a + (u \times v) \times (w \times u).
\end{align*}
Here the last three terms may be written as
\begin{align*}
\bar v \times (w \times u)a =\,\,& (u \times w) \times \bar va = w\bar v^\trans ua  - u\bar v^\trans wa, \\
(u \times v) \times \bar w\bar a =\,\,&v\bar w^\trans u\bar a - u\bar w^\trans v\bar a, \\
(u \times v) \times (w \times u) =\,\,&v(w \times u)^\trans u - u(w \times u)^\trans v = -u\det(u,v,w).
\end{align*}
Thus
\begin{align}
\label{UNU} u_0 =\,\,&-u\bar v^\trans wa - u\bar w^\trans v\bar a + (\bar v \times \bar w)a\bar a - \big((\bar u \times \bar v)\bar w^\trans  + (\bar w \times \bar u)\bar v^\trans \big)u - u\det(u,v,w).
\end{align}
We now apply Exc.~\ref{pr.INVAS}~(c) to tackle the expression 
\begin{align}
\label{REMO} x(yz)x = n_\IO(x,\overline{yz})x - n_\IO(x)\overline{yz}.
\end{align} 
Since $yz = (-\bar v^\trans w,\bar v \times \bar w)$, we deduce from (\ref{ss.NOTR}.\ref{CONGC}) that
\begin{align}
\label{COYZ} \overline{yz} = (-\bar w^\trans v,-\bar v \times \bar w). 
\end{align}
When combined with (\ref{ss.NOTR}.\ref{BLINOGC}), this implies
\begin{align}
\label{BNYZ} n_\IO(x,\overline{yz}) =\,\,&n_\IO(\overline{yz},x) = -\bar v^\trans wa - \bar w^\trans v\bar a - (v \times w)^\trans u - (\bar v \times \bar w)^\trans \bar u \notag \\
=\,\,&-\bar v^\trans wa - \bar w^\trans v\bar a - \det(u,v,w) - \overline{\det(u,v,w)},
\end{align}
hence, in view of \eqref{REMO}, \eqref{COYZ}, \eqref{ANU}
\begin{align*}
a_1 =\,\,&n_\IO(x,\overline{yz})a + n_\IO(x)\bar w^\trans v \\
=\,\,&-\bar v^\trans wa^2 - \bar w^\trans v\bar aa - \det(u,v,w)a - \overline{\det(u,v,w)}a + \bar w^\trans va\bar a + \bar w^\trans v\bar u^\trans u = a_0,
\end{align*}
giving our first claim. It remains to prove the second, i.e., $u_0 = u_1$. Again by \eqref{REMO}, \eqref{COYZ}, we obtain
\begin{align*}
u_1 =\,\,&n_\IO(x,\overline{yz})u + n_\IO(x)(\bar v \times \bar w) \\
=\,\,&-u\bar v^\trans wa - u\bar w^\trans v\bar a - u\det(u,v,w) - u\overline{\det(u,v,w)} +(\bar v \times \bar w)\bar aa + (\bar v \times \bar w)\bar u^\trans u,
\end{align*}
and comparing with \eqref{UNU} yields
\begin{align*}
u_1 - u_0 =\,\,&\big((\bar u \times \bar v)\bar w^\trans  + (\bar v \times \bar w)\bar u^\trans  + (\bar w \times \bar u)\bar v^\trans   - \det(\bar u,\bar v,\bar w)\Eins_3\big)u,
\end{align*}
which is zero by \eqref{CROMA}, and the assertion follows.
\end{sol}

\begin{sol}{rotation} \label{sol.rotation}
(a): Since
\[
n_\IH(v) = \cos^2(\theta/2) + \left\| \sin(\theta/2) u \right\|^2 = \cos^2(\theta/2) + \sin^2(\theta/2) = 1,
\]
$v$ is a versor.  

Conversely, if $v = \alpha \oplus t \in \IH$ is a versor, then since 
\[
1 = n_\IH(v) = \alpha^2 + \| t \|^2,
\]
$\alpha = \cos(\theta/2)$ for some angle $\theta$.  If $\alpha = 1$, then we can take $\theta = 0$ and $u$ any unit vector and we have verified the claim.  Similarly if $\alpha = -1$.  Thus we may assume that $\sin(\theta/2) \ne 0$ and set $u = (1/\sin(\theta/2)) t$ to obtain the required expression
\[
\cos(\theta/2) \oplus \sin(\theta/2) u = \alpha \oplus t = v.
\]

(b): We use $v^{-1} = \bar{v} = \cos(\theta/2) \oplus -\sin(\theta/2)u$ and expand the product $vsv^{-1}$ using the multiplication formula (\ref{ss.HAQU}.\ref{PROH}) and trigonometric double angle identities to find:
\[
vsv^{-1} = \frac{1-\cos \theta}{2} (u\cdot s)u + \frac{1+\cos \theta}2 s + (\sin \theta) u \times s - \frac{1 - \cos \theta}{2} (u \times s) \times u.
\]
Applying now $(u \times s) \times u = (u \cdot u) s - (u \cdot s) u$ simplifies the expression to the one given by Rodrigues' Formula.
\end{sol}

\solnsec{Section~\ref{s.CASH}}

\begin{sol}{pr.MUTA}\label{sol.MUTA} Put $N := \{(s,t) \in \IZ \times \IZ\mid 1 \leq s,t \leq 7,\;s \neq t\}$ and note that $M$ (resp. $N$) has $21$ (resp. $42$) elements. Defining maps $\varphi_\pm\:M \to N$ by
\begin{align}
\label{CTPM} \varphi_+(r,i) := (r + i,r + 3i), \quad \varphi_-(r,i) := (r + 3i,r + i) &&((r,i) \in M,\;\text{integers} \bmod 7),
\end{align}
it suffices to show that (i) the assignments \eqref{CTPM} for $\varphi_\pm$ do indeed take values in $N$, (ii) $\varphi_\pm$ are both injective, and (iii) their images in $N$ are disjoint. We prove these assertions one at a time.

(i) Let $(r,i) \in M$ and assume $r + i \equiv r + 3i \bmod 7$. Then $2i \equiv 0 \bmod 7$, hence $i \equiv 0 \bmod 7$, a contradiction.

(ii) Let $(r,i),(s,j) \in M$. If $\varphi_+(r,i) = \varphi_+(s,j)$, then $r + i \equiv s + j \bmod 7$ and $r + 3i \equiv s + 3j \bmod 7$. Taking differences, we conclude $2i \equiv 2j \bmod 7$, hence $i \equiv j \bmod 7$ and then $i = j$. This implies $r \equiv s \bmod 7$, hence $r = s$, and we have show that $\varphi_+$ is injective. But $\varphi_-$ agrees with the switch $(s,t) \mapsto (t,s)$ of $N$, which is bijective, composed with $\varphi_+$. Thus $\varphi_-$ is injective as well, which completes the proof of (ii).

(iii) Let $(r,i),(s,j) \in M$ and suppose $\varphi_+(r,i) = \varphi_-(s,j)$. Then $r + i \equiv s + 3j \bmod 7$, $r + 3i \equiv s + j \bmod 7$, and taking differences yields $2i \equiv -2j \bmod 7$, hence $i + j \equiv 0 \bmod 7$, a contradiction. This proves (iii).
\end{sol}

\begin{sol}{pr.CHACS}\label{sol.CHACS}
(i) $\Rightarrow$ (iii). Combining (\ref{ss.NOTR}.\ref{QUAGC}) with (\ref{ss.DECS}.\ref{SQCS}), we deduce $t_\IO(u_r) = 0$, $n_\IO(u_r) = 1$ for $1 \leq r \leq 7$, so $u_r \in \IO^0$ has euclidean norm $1$, while \eqref{HAMI} follows immediately from (\ref{ss.DECS}.\ref{MICS}). Thus (iii) holds. 

(iii) $\Rightarrow$ (iv). For $1 \leq r \leq 7$, $i = 1,2,4$, indices $\bmod 7$, we combine \eqref{HAMI} with (\ref{ss.NOTR}.\ref{TRANO}), (\ref{ss.NOTR}.\ref{PUGC}) and obtain $n_\IO(u_{r+i},u_{r+3i})  = t_\IO(u_{r+i}\bar u_{r+3i}) = -t_\IO(u_{r+i}u_{r+3i}) = -t_\IO(u_r) = 0$. By Exc.~\ref{pr.MUTA}, therefore, $(u_r)_{1\leq r\leq 7}$ forms an orthonormal basis of $\IO^0$, giving the final assertion of the problem. The second set of equations in \eqref{EXMI} is just \eqref{HAMI} for $r \geq 4$, $i = 4$, and implies $n_\IO(u_1u_2,u_3) = n_\IO(u_4,u_3) = 0$. Thus (iv) holds.

(iv) $\Rightarrow$ (i). Systematically counting indices $\bmod\;7$, we establish the desired implication by proving the following intermediate assertions. 

\step{1}
\emph{$n_\IO(u_r) = 1$ for $1 \leq r \leq 7$}. By (iv) this is clear for $r \leq 3$, and for $r \geq 4$ follows from \eqref{EXMI} by induction since the norm of $\IO$ by Theorem~\ref{t.DIGC} permits composition. 

\step{2}
\emph{$n_\IO(u_r,u_s) = 0$ for $1 \leq r,s \leq 4$ distinct.} By (iv), we may assume $s = 4$, $r \leq 2$. Then \ref{r.LICO} implies $n_\IO(u_r,u_4) = n_\IO(u_r,u_1u_2) = n_\IO(u_r)t_\IO(u_{3-r}) = 0$, as claimed. 

\step{3}
\emph{Let $1 \leq r \leq 7$ and $i = 1,2,4$ such that $u_r = u_{r+i}u_{r+3i}$. If $u_{r+i}$ has trace zero, then $u_r$ has trace zero if and only if $u_{r+i},u_{r+3i}$ are orthogonal. Moreover, if this is so and $u_{r+i}$ as well as $u_{r+3i}$ both have trace zero, then $u_r = -u_{r+3i}u_{r+i}$.} Since $\bar u_{r+i} = -u_{r+i}$, we deduce from (\ref{ss.NOTR}.\ref{TRANO}) that $t_\IO(u_r) = -n_\IO(u_{r+i},u_{r+3i})$, giving the first assertion. Since, under the final conditions stated, $u_r,u_{r+i},u_{r+3i}$ are all skew relative to the conjugation of $\IO$, and since this is an algebra involution, we conclude $u_r = -\bar u_r = -\bar u_{r+3i}\bar u_{r+i} = -u_{r+3i}u_{r+i}$, as claimed. 

\step{4}
\emph{$u_r \in \IO^0$ for $1 \leq r \leq 7$.} By (iv), this is automatic for $r \leq 3$, while it follows immediately from $2^\circ$, $3^\circ$ (with $i = 4$) for $r = 4,5,6$. It remains to discuss the case $r = 7$. We have $u_7 = u_4u_5$, and $3^\circ$ combined with \ref{r.LICO} yields $n_\IO(u_4,u_5) = n_\IO(u_1u_2,u_2u_3) = -n_\IO(u_2u_1,u_2u_3) = -n_\IO(u_2)n_\IO(u_1,u_3) = 0$. Hence the first part of  $3^\circ$ yields $t_\IO(u_7) = 0$. 

\step{5}
\emph{$(u_r)_{1\leq r\leq 7}$ is an orthonormal basis of $\IO^0$ and $u_r^2 = -u_0$, $u_ru_s = -u_su_r$ for $1 \leq r,s \leq 7$ distinct.} The first part follows immediately from $1^\circ$, $3^\circ$, $4^\circ$ combined with Exc.~\ref{pr.MUTA}, while the second part is now a consequence of $1^\circ$ and (\ref{ss.NOTR}.\ref{QUAGC}), (\ref{ss.NOTR}.\ref{BLAGC}). 

\step{6}
\emph{Let $1 \leq r \leq 7$ and $i = 1,2,4$ such that $u_{r+i}u_{r+3i} = u_r$. Then}  
\begin{align}
\label{RTWI} u_{s+j}u_{s+3j} =\,\,&u_s &&(s = r + i,\;j \in \{1,2,4\},\;j \equiv 2i \bmod 7), \\
\label{RFOI} u_{s+j}u_{s+3j} =\,\,&u_s &&(s = r + 3i,\;j \in \{1,2,4\},\;j \equiv 4i \bmod 7).
\end{align}
Since $u_{r+3i}u_{r+i} = -u_r$ by $5^\circ$ and $\IO$ is an alternative algebra, $u_{r+3i}u_r = -u_{r+3i}^2u_{r+i} = u_{r+i}$, and \eqref{RTWI} follows. \eqref{RFOI} is proved similarly. 

The proof of the implication (iv) $\Rightarrow$ (i) will be complete once we have shown that the hypothesis of $6^\circ$ holds for all $r = 1,\dots, 7$ and $i = 4$. Hence it suffices to establish  

\step{7}
\emph{$u_{r+4}u_{r+5} = u_r$ for $1 \leq r \leq 7$.} By \eqref{EXMI}, this is clear for $4 \leq r \leq 7$, and $6^\circ$ yields
\begin{align}
\label{ESONE} u_{s+1}u_{s+3} =\,\,&u_s &&(1 \leq s \leq 4), \\
\label{ESTWO} u_{s+2}u_{s+6} =\,\,&u_s &&(2 \leq s \leq 5).
\end{align}
Now let $1 \leq r \leq 2$. Then \eqref{EXMI} combined with $5^\circ$, the middle Moufang identity (cf. Exc.~\ref{pr.MIMOU}) and \eqref{ESONE} yields
\begin{align*}
u_{r+4}u_{r+5} =\,\,&(u_{r+1}u_{r+2})(u_{r+2}u_{r+3}) = (u_{r+2}u_{r+1})(u_{r+3}u_{r+2}) = u_{r+2}(u_{r+1}u_{r+3})u_{r+2} \\
=\,\,&u_{r+2}u_ru_{r+2} = -u_{r+2}(u_{r+2}u_r) = -u_{r+2}^2u_r = u_r,
\end{align*}
so we are left with the case $r = 3$. Applying \eqref{ESONE} for $s = 1$ and \eqref{ESTWO} for $s = 3$, we obtain
\begin{align*}
u_7u_8 =\,\,&(u_4u_5)u_1 = (u_4u_5)(u_2u_4) = u_4(u_5u_2)u_4 = u_4(u_{3+2}u_{3+6})u_4 \\
=\,\,&u_4u_3u_4 = -u_4(u_4u_3) = -u_4^ 2u_3 = u_3,
\end{align*}
as claimed.

Summing up, we have shown not only that conditions (i), (iii), (iv) are equivalent but also the final statement of the problem. It remains to establish (i) $\Leftrightarrow$ (ii). The first equation of \eqref{ELGRO} is the same as (\ref{ss.DECS}.\ref{SQCS}). As to the second, the alternative laws combined with the fact that $\IO$ is a division algebra yield the following chain of equivalent conditions.     
\begin{align}
\label{FELGRO} (u_ru_{r+1})u_{r+3} = -1_\IO\,\,&\Longleftrightarrow (u_ru_{r+1})u_{r+3}^2 = -u_{r+3} \Longleftrightarrow u_ru_{r+1} = u_{r+3} \\
\,\,&\Longleftrightarrow u_{(r+3) + 4}u_{(r+3)+3\cdot 4} =  u_{r+3}. \notag
\end{align} 
Similarly, the third equation of \eqref{ELGRO} may be rephrased as
\begin{align}
\label{SELGRO}  u_r(u_{r+1}u_{r+3}) = -1_\IO \,\,&\Longleftrightarrow u_r^2(u_{r+1}u_{r+3}) = -u_r \Longleftrightarrow u_{r+1}u_{r+3} = u_r.
\end{align}
Hence (i) implies (ii). Conversely, suppose (ii) holds. Then the first equation of \eqref{ELGRO} combined with (\ref{ss.NOTR}.\ref{QUAGC}) shows that the elements $u_r$, $1 \leq r \leq 7$, belong to $\IO^ 0$ and have length $1$, while \eqref{FELGRO}, \eqref{SELGRO} yield $u_ru_{r+1} = u_{r+3}$, $u_{r+1}u_{r+3} = u_r$. Combining with (\ref{ss.NOTR}.\ref{PUGC}), (\ref{ss.NOTR}.\ref{TRANO}) we deduce $n_\IO(u_r,u_{r+1}) = t_\IO(u_r\bar u_{r+1}) = -t_\IO(u_ru_{r+1}) = t_\IO(u_{r+3}) = 0$ and, similarly, $n_\IO(u_{r+1},u_{r+3}) = 0$. In particular, $u_1,u_2,u_3$ form an orthonormal system in $\IO^0$ such that $n_\IO(u_1u_2,u_3) = n_\IO(u_3,u_4) = 0$, and we have $u_r = u_{r-3}u_{r-2}$ for $4 \leq r \leq 7$. Since, therefore, condition (iv) holds, so does (ii).  
 
\end{sol}

\begin{sol}{pr.OPLU}\label{sol.OPLU}
$\IO^+$ is clearly a commutative real algebra with identity element $1_{\IO^+} = 1_\IO$. Moreover, the squarings in $\IO$ and $\IO^+$ coincide. Hence any $\varphi \in \Aut(\IO^+)$ fixes $1_\IO$ and preserves squares in $\IO$. From (\ref{ss.NOTR}.\ref{QUAGC}) we therefore deduce $t(x)x = n(x)1_\IO$ for all $x \in \IO$, where $t := t_\IO \circ \varphi - t_\IO$, $n := n_\IO \circ \varphi - n_\IO$. Thus $\IO$ is the \emph{union} of the subspaces $\Ker(t)$ and $\IR1_\IO$. This implies first $t = 0$ and then $n = 0$, so $\varphi$ is an orthogonal transformation of $\IO$ fixing $1_\IO$, hence stabilizing $\IO^0 = (\IR1_\IO)^\perp$. The assignment $\varphi \mapsto \varphi\vert_{\IO^0}$ clearly gives an injective group homomorphism from $\Aut(\IO^+)$ to $\Or(\IO^0)$. Conversely, the unique linear extension fixing $1_\IO$ of $\psi \in \Or(\IO^0)$ belongs to the orthogonal group of $\IO$ and leaves the trace invariant, hence is an automorphism of $\IO^+$. Moreover, the maps in question are inverse isomorphisms not only of abstract groups but, in fact, of topological ones since both are restrictions of linear maps and hence automatically continuous. Finally, $\Aut(\IO)$ is trivially a subgroup of $\Aut(\IO^+)$ which is closed since it may be realized via
\[
\Aut(\IO) = \bigcap_{x,y \in \IO} \{\varphi \in \GL(\IO)\mid \varphi(xy) = \varphi(x)\varphi(y)\}
\] 
as the intersection of a family of closed subsets of $\GL(\IO)$.
\end{sol}

\solnsec{Section~\ref{s.SUBGCZ}}

\begin{sol}{pr.CHALAT} \label{sol.CHALAT}
Put $n := \dim_\IR(V)$. 

(i) $\Rightarrow$ (iii). Obvious.

(iii) $\Rightarrow$ (ii). Since $L_0$ spans $V$ as a real vector space, so does $L$. Since $L$ is an additive subgroup of $L_1$, it is finitely generated and torsion-free, hence free of finite rank with $\rk(L) \leq \rk(L_1) = n$.

(ii) $\Rightarrow$ (i). Let $\vep = (e_1,\dots,e_r)$, $r \leq n$, be a $\IZ$-basis of $L$. Since $L$ spans $V$, so does $\vep$, and we conclude that $\vep$ is an $\IR$-basis of $V$, forcing $L \subseteq V$ to be a lattice.
\end{sol}

\begin{sol}{pr.LASTRU} \label{sol.LASTRU}
$M + \sqrt{2}M \subseteq \ID$ is a free abelian subgroup of rank $2\dim_\IR(\ID)$ which is obviously closed under multiplication but not a lattice, hence not a $\IZ$-structure of $A$. 
\end{sol}

\begin{sol}{pr.REDONE}\label{sol.REDONE} Since $M$ is clearly a discrete additive subgroup of $\ID$, it follows that $M^\prime := M \cap \IR$ is a discrete additive subgroup of $\IR$ containing $1$. Thus $M^\prime = \IZ v$ for some $0 \neq v \in \IR$. Hence $M^\prime$ is a $\IZ$-structure of $\IR$. As such it contains $\IZ$. On the other hand, $v^2 \in M^\prime$ implies $v^2 = rv$ for some integer $r$, hence $v = r \in \IZ$, and we also conclude $M^\prime \subseteq \IZ$, as desired.
\end{sol}

\solnsec{Section~\ref{s.MAQO}}

\begin{sol}{pr.HURCON} \label{sol.HURCON}
By Thm.~\ref{t.HURO}, $(1_\IH,\bfi,\bfj,\bfh)$ is an $\IR$-basis of $\IH$ that  is associated with $\Hur(\IH)$. Hence $\Hur(\IH)$ consists of all linear combinations
\begin{align*}
m1_\IH + n\bfi + p\bfj + q\bfh = (m + \frac{q}{2})1_\IH + (n + \frac{q}{2})\bfi + (p + \frac{q}{2})\bfj + \frac{q}{2}\bfk 
\end{align*}
with $m,n,p,q \in \IZ$. Since the coefficients of the linear combination on the right either all belong to $\IZ$ or to $\frac{1}{2} + \IZ$ depending on whether $q$ is even or odd, the assertion follows. 
\end{sol}

\begin{sol}{pr.MAXHUR} \label{sol.MAXHUR}
Write $x \in L$ in the form $x = \alpha 1_\IH + \beta\bfi + \gamma\bfj + \delta\bfh$ with $\alpha,\beta,\gamma,\delta \in \IR$. We must show $\alpha,\beta,\gamma,\delta \in \IZ$. Replacing $x$ by $x - \lfloor x\rfloor$, $\lfloor x\rfloor := \lfloor\alpha\rfloor1_\IH + \lfloor\beta\rfloor\bfi + \lfloor\gamma\rfloor\bfj + \lfloor\delta\rfloor\bfh \in \Hur(\IH) \subseteq L$, if necessary, we may assume $0 \leq \alpha,\beta,\gamma,\delta < 1$ and then must show $\alpha = \beta = \gamma = \delta = 0$. By hypothesis and (\ref{ss.HUZH},\ref{ENTH}), the quantities
\begin{align}
\label{NHX} n_\IH(x) =\,\,&\alpha^2 + \beta^2 + \gamma^2 + (\alpha + \beta + \gamma + \delta)\delta, \\
\label{THX} t_\IH(x) =\,\,& 2\alpha + \delta, \\
\label{NHIX} n_\IH(\bfi,x) =\,\,&2\beta + \delta, \\
\label{NHJX} n_\IH(\bfj,x) =\,\,&2\gamma + \delta, \\
\label{NHHX} n_\IH(\bfh,x) =\,\,&\alpha + \beta + \gamma + 2\delta 
\end{align}
are all integers. Multiplying \eqref{NHHX} by $2$ and subtracting the sum of \eqref{THX}, \eqref{NHIX}, \eqref{NHJX} from the result, we conclude $\delta \in \IZ$, hence $\delta = 0$, and then $2\alpha,2\beta,2\gamma,\alpha^2 + \beta^2 + \gamma^2 \in \IZ$ by \eqref{NHX}--\eqref{NHJX}. Thus $\alpha,\beta,\gamma \in \{0,\frac{1}{2}\}$, forcing $\alpha = \beta = \gamma = 0$ as well. 
\end{sol}

\begin{sol}{pr.DETFO} \label{sol.DETFO}
We begin by deriving a general formula for the determinant of a block matrix. Let $A \in \GL_p(\IR)$, $B \in \Mat_{p,q}(\IR)$, $C \in \Mat_{q,p}(\IR)$ and $D \in \Mat_q(\IR)$. Then we claim
\begin{align}
\label{DEBLO} \det(\left(\begin{matrix}
A & B \\
C & D
\end{matrix}\right)) = \det(A)\det(D - CA^{-1}B).
\end{align}
To see this, we make use of the well known fact that \eqref{DEBLO} holds for $B = 0$ or $C = 0$. Hence
\begin{align*}
\det(\left(\begin{matrix}
A & B \\
C & D
\end{matrix}\right)) =\,\,&\det(A)\det(\left(\begin{matrix}
A & B \\
C & D
\end{matrix}\right)\left(\begin{matrix}
A^ {-1} & 0 \\
0 & \Eins_q
\end{matrix}\right)) = \det(A)\det(\left(\begin{matrix}
\Eins_p & B \\
CA^ {-1} & D
\end{matrix}\right)) \\
=\,\,&\det(A)\det(\left(\begin{matrix}
\Eins_p & 0 \\
-CA^ {-1} & \Eins_q
\end{matrix}\right)\left(\begin{matrix}
\Eins_p & B \\
CA^ {-1} & D
\end{matrix}\right)) = \det(A)\det(\left(\begin{matrix}
\Eins_p & B \\
0 & -CA^{-1}B + D
\end{matrix}\right)) \\
=\,\,&\det(A)\det(D - CA^ {-1}B),
\end{align*}
as claimed. Turning to the matrix $S$ of our problem, we now conclude
\begin{align*}
\det(S) =\,\,&\det(r\cdot\Eins_p)\det\big(s\cdot\Eins_q - T_2(r^{-1}\cdot\Eins_p)T_1\big) = r^p\det(s\cdot\Eins_q - r^{-1}T_2T_1) \\
=\,\,&r^{p-q}\det(rs - T_2T_1) = r^{p-q}\ch_{T_2T_1}(rs),
\end{align*} 
as claimed. Next suppose $T_1 = T$, $T_2 = T^\trans$, with $T \in \Mat_{p,q}(\IR)$ being comprised of column vectors $v_1,\dots,v_q \in \IR$ that are mutually orthogonal but may have different euclidean lengths. Then
\[
T_2T_1 = T^\trans T = \left(\begin{matrix}
v_1^\trans  \\
\vdots \\
v_q^\trans 
\end{matrix}\right)(v_1,\dots,v_q) = (v_i^\trans v_j) _{1\leq i,j\leq q} = \diag(\|v_1\|^2,\dots,\|v_q\|^2).
\]
Thus the preceding determinant formula reduces to
\[
\det(S) = r^{p-q}\prod_{i=1}^q(rs - \|v_i\|^ 2). 
\]
In particular, if $\|v_i\| = \sqrt{a}$ for some $a > 0$ and all $i = 1,\dots,q$, the second assertion stated in our problem drops out. Finally, assuming $\|v_i\|^2 < rs$ for all $i = 1,\dots,q$, we claim that the matrix $S$ is positive definite. Since it is symmetric, it will be enough to show that all its principal minors are positive. Let $1 \leq i \leq p + q$. For $i \leq p$, the the $i$-th principal minor of $S$ is $1$, while for $p < i \leq p + q$, it is the determinant of the matrix
\[
S^\prime := \left(\begin{matrix}
r\cdot\Eins_p & T^\prime \\
(T^\prime)^\trans  & s\cdot\Eins_{i-p}
\end{matrix}\right),
\]
where $T^\prime = (v_1,\dots,v_{i-p}) \in \Mat_{p,i-p}(\IR)$. But by what we have just seen, $S^\prime$ has determinant 
\[
\det(S^\prime) = r^{p-q}\prod_{j=1}^{p-i}(rs - \|v_j\|^2) > 0, 
\]
and the assertion follows. Since the final statement of the problem is a special case of this assertion, the solution is complete.
\end{sol}

\begin{sol}{pr.UNIT} \label{sol.UNIT}
(a) If (iii) holds, then (\ref{ss.NOTR}.\ref{CONOTR}) implies $u\bar u = n_\ID(u)1_\ID = 1_\ID = \bar uu$. Thus (iii) implies (i) and (ii). Conversely, suppose (i) holds. Then $1 = n_\ID(uv) = n_\ID(u)n_\ID(v)$, and since both factors on the right are integers, (iii) holds. This shows that (i) and (iii) are equivalent. The equivalence of (ii) and (iii) is proved similarly. Moreover, for any $v$ satisfying (i) we apply Exc.~\ref{pr.INVAS}~(c) and obtain $\bar u = 1_\ID\bar u = (vu)\bar u = n_\ID(u)v = v$ (by (iii)). Reading this in $\ID^{\op}$ shows that $v$ as in (ii) is unique and equal to $\bar u$. Since the norm of $\ID$ permits composition by Thm.~\ref{t.DIGC}, the set of units in $M$ is closed under multiplication and contains $1_\ID$. Finally, by (\ref{ss.NOTR}.\ref{NOCOGC}), if $n_\ID(u) = 1$, then $n_\ID(\bar u) = 1$, so $M^\times$ is also closed under taking inverses.

(b) $(1,i)$ is an orthonormal basis of $\IC$ associated with $\Ga(\IC)$. Hence $\Ga(\IC)^\times = \{\pm 1,\pm i\}$. Similarly, $(1_\IH,\bfi,\bfj,\bfk)$ is an orthonormal basis of $\IH$ associated with $\Ga(\IH)$. Hence $\Ga(\IH)^\times = \{\pm 1_\IH,\pm\bfi,\pm\bfj,\pm\bfk\}$. And, finally, any Cartan-Schouten basis $E = (u_r)_{0\leq r\leq 7}$ of $\IO$ is an orthonormal basis  associated with $\Ga_\vep(\IO)$. Hence $\Ga_E(\IO)^\times = \{\pm u_r\mid 0 \leq r \leq 7\}$. Now we claim
\begin{align}
\label{HUNI} \Hur(\IH)^\times = \{\pm 1_\IH,\pm\bfi,\pm\bfj,\pm\bfk,\frac{1}{2}(\pm 1_\IH \pm\bfi \pm\bfj \pm\bfk)\}.
\end{align} 
To prove this, we note that the first eight elements of the right-hand side belong to $\Ga(\IH)^\times \subseteq \Hur(\IH)^\times$. The remaining sixteen may all be obtained from adding a unit of $\Ga(\IH)$ to $\bfh$, hence belong to $\Hur(\IH)$. Since they also have norm $1$, the right-hand side of \eqref{HUNI} belongs to the left. Conversely, suppose $u = \alpha_0 1_\IH + \alpha_1\bfi + \alpha_2\bfj + \alpha_3\bfk$ with $\alpha_i \in \IR$, $0 \leq i \leq 3$, is invertible in $\Hur(\IH)$. By Exc.~\ref{pr.HURCON}, there are two cases. 

\case{1} 
$\alpha_i \in \IZ$, $0 \leq i \leq 3$. Then $u \in \Ga(\IH)^\times \subseteq \Hur(\IH)^\times$. 

\case{2}
$\alpha_i = \frac{1}{2} + \beta_i$, $\beta_i \in \IZ$, $0 \leq i \leq 3$. Then
\[
1 = n_\IH(u) = \sum_{i=0}^3 (\frac{1}{4} + \beta_i + \beta_i^2) = 1 + \sum_{i=0}^3 (\beta_i^2 + \beta_i),
\]  
and we conclude $\sum_{i=0}^3(\beta_i^2 + \beta_i) = 0$, where the summands are all non-negative. Thus $\beta_i^2 = -\beta_i$, i.e., $\beta_i = -1,0$ for $0 \leq i \leq 3$. But this means $\alpha_i = \pm\frac{1}{2}$ for $0 \leq i \leq 3$, and $u$ belongs to the right-hand side of \eqref{HUNI}. 
\end{sol}

\begin{sol}{pr.HURDE} \label{sol.HURDE}
That the $\vep_i$, $1 \leq i \leq 4$, form an orthonormal basis of $\IH$ relative to $2\la x,y \ra$ follows immediately from the fact that the vectors $1_\IH,\bfi,\bfj,\bfk$ of $\IH$ are orthogonal of norm $1$. Moreover, the latter vectors can be linearly expressed by the former according to the formulas
\[
\vep_1 - \vep_2 = \bfj, \quad \vep_2 - \vep_3 = -\bfh, \quad \vep_3 - \vep_4 = 1_\IH, \quad \vep_3 + \vep_4 = \bfi.
\]
In particular, the second equation of the problem holds. Given $\xi_i \in \IR$, $1 \leq i \leq 4$, we also obtain
\begin{align*}
\sum_{i=1}^4\xi_i\vep_i =\,\,&\frac{1}{2}\big(\xi_1\bfj - \xi_1\bfk - \xi_2\bfj - \xi_2\bfk + \xi_3 1_\IH + \xi_3\bfi - \xi_4 1_\IH + \xi_4\bfi\big) \\
=\,\,&\frac{1}{2}\big((\xi_3 - \xi_4)1_\IH + (\xi_3 + \xi_4)\bfi + (\xi_1 - \xi_2)\bfj - (\xi_1 + \xi_2)\bfk\big).
\end{align*}
Using (\ref{ss.HUZH}.\ref{HIJK}), we therefore deduce
\begin{align*}
\sum_{i=1}^4\xi_i\vep_i =\,\,&\frac{1}{2}\big((\xi_3 - \xi_4)1_\IH + (\xi_3 + \xi_4)\bfi + (\xi_1 - \xi_2)\bfj - (\xi_1 + \xi_2)(2\bfh - 1_\IH - \bfi - \bfj)\big) \\
=\,\,&\frac{1}{2}\big((\xi_1 + \xi_2 + \xi_3 - \xi_4)1_\IH + (\xi_1 + \xi_2 + \xi_3 + \xi_4)\bfi + 2\xi_1\bfj\big) - (\xi_1 + \xi_2)\bfh \\
=\,\,&\frac{1}{2}(\xi_1 + \xi_2 + \xi_3 - \xi_4)1_\IH + \frac{1}{2}(\xi_1 + \xi_2 + \xi_3 + \xi_4)\bfi + \xi_1\bfj - (\xi_1 + \xi_2)\bfh,
\end{align*}
which by Thm.~\ref{t.HURO} yields the following chain of equivalent conditions.
\begin{align*}
\sum_{i=1}^4\xi_i\vep_i \in \Hur(\IH) \quad \,\,&\Longleftrightarrow \quad \xi_1,\xi_1 + \xi_2 \in \IZ, \quad \xi_1 + \xi_2 + \xi_3 \pm \xi_4 \in 2\IZ \\
\,\,&\Longleftrightarrow \quad \xi_1,\xi_2 \in \IZ, \quad \sum_{i=1}^4\xi_i \in 2\IZ, \quad 2\xi_4 \in 2\IZ \\
\,\,&\Longleftrightarrow \quad \xi_1,\xi_2,\xi_3,\xi_4 \in \IZ, \quad \sum_{i=1}^4\xi_i \in 2\IZ.
\end{align*}
Thus the second equation of the problem holds. Finally, we note $n_\IH(\vep_i) = \frac{1}{2}$ for $1 \leq i \leq 4$, so for $\xi_i \in \IZ$, $\sum \xi_i \in 2\IZ$, we obtain
\begin{align*}
\sum_{i=1}^4\xi_i\vep_i \in \Hur(\IH)^\times \quad \,\,&\Longleftrightarrow \quad \frac{1}{2}\sum_{i=1}^4\xi_i^2 = 1 \quad \Longleftrightarrow \quad \sum_{i=1}^4\xi_i^2 = 2,
\end{align*} 
which amounts to $\xi_i = \pm 1$ for precisely two indices $i = 1,2,3,4$ and $\xi_i = 0$ for the remaining ones. But this means the units of $\Hur(\IH)$ have the form stated at the very end of the problem.
\end{sol}

\begin{sol}{pr.ROOC} \label{sol.ROOC}
(a) This follows immediately from the fact that Cartan-Schouten bases are orthonormal relative to the scalar product $\la x,y\ra$. (Exc.~\ref{pr.CHACS}). For example, $n_\IO(\vep_1) = \frac{1}{4}(n_\IO(u_0) + n_\IO(u_2) = \frac{1}{2}$, which implies $2\la \vep_1,\vep_1\ra = 2n_\IO(\vep_1) = 1$.

(b) We write $L$ for the additive subgroup of $\IO$ generated by the quantities assembled in \eqref{EEIGHT}. We claim
\begin{align}
\label{LEIGHT} \Cox(\IO) \subseteq L \subseteq \frac{1}{2}\sum_{i=1}^8 \IZ\vep_i,
\end{align} 
where the second inclusion is obvious. In order to prove the first, we note that
\begin{align}
\label{UZER} u_0 =\,\,&-\vep_1 + \vep_2, \quad u_1 = -\vep_3 + \vep_4, \quad u_2 = \vep_1 + \vep_2, \quad u_3 = -\vep_3 - \vep_4, \\ 
u_4 =\,\,&-\vep_5 + \vep_6, \quad u_5 =  \vep_5 + \vep_6,, \quad u_6 = \vep_7 + \vep_8, \quad u_7 = -\vep_7 + \vep_8 \notag
\end{align}
all belong to $L$. Moreover, \eqref{UZER} and (\ref{ss.COZO}.\ref{BFEL})--(\ref{ss.COZO}.\ref{FOUL})  imply
\begin{align}
\label{PEP} \bfp =\,\,&\frac{1}{2}(-\vep_1 + \vep_2 - \vep_3 + \vep_4 + \vep_1 + \vep_2  - \vep_3 - \vep_4) = \vep_2 - \vep_3 \in L, \\
\label{ONEPEP} u_1\bfp =\,\,&\frac{1}{2}(\vep_1 - \vep_2 - \vep_3 + \vep_4 - \vep_5 + \vep_6 - \vep_7 + \vep_8) \in L, \\
\label{TWOPEP} u_2\bfp =\,\,&\frac{1}{2}(\vep_1 - \vep_2+ \vep_1 + \vep_2 + \vep_5 - \vep_6 + \vep_5 + \vep_6) = \vep_1 + \vep_5 \in L, \\
\label{FOURPEP} u_4\bfp =\,\,&\frac{1}{2}(\vep_3 - \vep_4+ \vep_1 + \vep_2 - \vep_5 + \vep_6 - \vep_7 - \vep_8) \in L.
\end{align}
Thus the basis $E^\prime$ exhibited in (\ref{t.COGC}.\ref{BECO}), which  is associated with $\Cox(\IO)$ by Thm.~\ref{t.COGC}, is entirely contained in $L$, and the first inclusion of \eqref{LEIGHT} holds. From this and Exc.~\ref{pr.CHALAT} we deduce that $L$ is a lattice in $\IO$. Since $\Cox(\IO)$ is unimodular by Thm.~\ref{t.COGC}, the proof of (b)  will be complete once we have shown that the lattice $L$ is integral quadratic (Cor.~\ref{c.MAXIN}). Since $n_\IO(\vep_i) = \frac{1}{2}$ for $1 \leq i \leq 8$ by (a), equation \eqref{EEIGHT} yields
\begin{align}
\label{UNEI} n_\IO(\pm\vep_i\,\pm\,\vep_j) = n_\IO\big(\frac{1}{2}\sum_{i=1}^8 s_i\vep_i\big) = 1
\end{align}
for $1 \leq i < j \leq 8$ and all families $s = (s_i) \in \{\pm 1\}^8$. Writing $M := \{1,2,\dots,8\}$ and $M_s^\pm := \{i \in M\mid s_i = \pm 1\}$, it remains to show
\begin{align}
\label{MIXEI} n_\IO(\pm\vep_i\,\pm\,\vep_j,\pm\vep_k\,\pm\,\vep_l) \in \IZ, \quad n_\IO\big(\frac{1}{2}\sum_{i=1}^8s_i\vep_i,\frac{1}{2}\sum_{i=1}^8t_i\vep_i\big) \in \IZ
\end{align}
for all $1 \leq i < j \leq 8$, $1 \leq k < l \leq 8$ and all families $s = (s_i),t = (t_i) \in \{\pm 1\}^8$ such that $\vert M_s^-\vert$ and $\vert M_t^-\vert$ are both even. The first relation of \eqref{MIXEI} is obvious. In order to prove the second, we let $s = (s_i) \in \{\pm 1\}^8$ be arbitrary and note $\vert M_s^+\vert + \vert M_s^-\vert = 8$, so $\vert M_s^\pm\vert$ are either both even or both odd. Moreover, $-s := (-s_i) \in \{\pm 1\}^8$ and $M_{-s}^\pm = M_s^\mp$. Now, letting also $t = (t_i) \in \{\pm 1\}^8$ be arbitrary and setting $st := (s_it_i) \in \{\pm 1\}^8$, we obtain 
\[
n_\IO\big(\frac{1}{2}\sum_{i=1}^8s_i\vep_i,\frac{1}{2}\sum_{i=1}^8t_i\vep_i\big) =\frac{1}{4}\sum_{i=1}^8s_it_i = \frac{1}{4}\big(\vert M_{st}^+\vert - \vert M_{st}^-\vert\big) = 2 - \frac{1}{2}\vert M_{st}^-\vert,
\]  
and must show that if $\vert M_s^\pm\vert$ and $\vert M_t^\pm\vert$ are even, then so is $\vert M_{st}^\pm\vert$. Suppose not. Then
\[
\vert M_{st}^+\vert = \vert M_s^+ \cup M_t^+\vert + \vert M_s^- \cup M_t^-\vert
\]
is odd. Replacing $s$ by $-s$, $t$ by $-t$ if necessary, we may assume that $\vert M_s^+ \cup M_t^+\vert$ is odd and $\vert M_s^- \cup M_t^-\vert$ is even. But then
\[
\vert M_s^+ \vert = \vert M_s^+ \cup M_t^+\vert + \vert M_s^+ \cup M_t^-\vert
\]
shows that $\vert M_s^+ \cup M_t^-\vert$ is odd, whence 
\[
\vert M_t^-\vert = \vert M_s^+ \cup M_t^-\vert + \vert M_s^- \cup M_t^-\vert 
\]
is odd as well, a contradiction.

(c) Write $L^\prime$ for the right-hand side of \eqref{COXCON}. Inspecting \eqref{UZER}--\eqref{FOURPEP}, we conclude $\Cox(\IO) \subseteq L^ \prime \subseteq \frac{1}{2}\sum \IZ\vep_i$, so $L^ \prime$ is a lattice in $\IO$ by Exc.~\ref{pr.CHALAT}. As before, we will be through once we have shown that $L^ \prime$ is integral quadratic. Let $x = \sum\xi_i\vep_i \in L^\prime$ with $\xi_1,\dots,\xi_8 \in \IR$. Then \eqref{COXCON} implies $2\xi_1 \in \IZ$, $\xi_i = \xi_1 + n_i$, $n_i \in \IZ$ for $1 \leq i \leq 8$, and
\[
\sum_{i=1}^8\xi_i = \sum_{i=1}^8(\xi_1 + n_i) = 4(2\xi_1) + \sum_{i=1}^8n_i
\] 
belongs to $2\IZ$. Thus $\sum n_i \in 2\IZ$, and $N := \{i \in M \mid n_i \equiv 1 \bmod 2\}$ contains an even number of elements. Now $n_\IO(x) = \frac{1}{2}\sum\xi_i^2$ and 
\begin{align*}
\sum_{i=1}^8\xi_i^2 =\,\,& \sum_{i=1}^8(\xi_1 + n_i)^2 = \sum_{i=1}^8(\xi_1^2 + 2\xi_1n_i + n_i^2) = 2(2\xi_1)^2 + 2\xi_1\sum_{i=1}^8n_i + \sum_{i=1}^8n_i^2 \\
\equiv\,\,&\sum_{i=1}^8n_i^2 \equiv \sum_{i \in N} n_i^2 \equiv \vert N\vert \equiv 0 \bmod 2. 
\end{align*}
Hence $n_\IO(x) \in \IZ$, as claimed.

(d) By (b), (c) and \eqref{UNEI}, the quantities in \eqref{EEIGHT} are all units in $\Cox(\IO)$. Conversely, suppose $u = \sum_{i=1}^8\xi_i\vep_i$, $\xi_i \in \IR$, $1 \leq i \leq 8$, is invertible in $\Cox(\IO)$. Then (a) implies
\begin{align}
\label{UNO} \sum_{i=1}^8\xi_i^2 = 2n_\IO(u) = 2.
\end{align} 
If $\xi_1$ is an integer, then all $\xi_i$, $1 \leq i \leq 8$, are, and \eqref{UNO} shows that $u$ belongs to the first type of \eqref{EEIGHT}. On the other hand, if $\xi_1$ belongs to $\frac{1}{2} + \IZ$, then all $\xi_i$, $1 \leq i \leq 8$, do, and \eqref{UNO} shows $\xi_i^2 = \frac{1}{4}$, hence $\xi_i = \pm\frac{1}{2}$ for all  $i = 1,\dots,8$. Since $\sum\xi_i \in 2\IZ$ by \eqref{COXCON} the number if indices $i = 1,\dots,8$ having $\xi_i$ negative must be even, so $u$ is of the second type in \eqref{EEIGHT}. 

The units of $\Cox(\IO)$ belonging to the first type of \eqref{EEIGHT} are in bijective correspondence with the quadruples of subsets of $M$ consisting of two elements, hence are $4\cdot\binom{8}{2} = 4\cdot\frac{8\cdot 7}{2} = 4\cdot 28 = 112$ in number. On the other hand, the units of $\Cox(\IO)$ belonging to the second type of \eqref{EEIGHT} are in bijective correspondence with the subsets of $M$ consisting of an even number of elements, hence are $\frac{2^8}{2} = 2^7 = 128$ in number. All in all, therefore, $\Cox(\IO)$ has exactly $112 + 128 = 240$ units.  

\end{sol}

\begin{sol}{pr.KILA} \label{sol.KILA}
The inclusions $\Ga(\IO) \subseteq \Kir(\IO) \subseteq \frac{1}{2}\Ga(\IO)$ follow immediately from the definition and show by Exc.~\ref{pr.CHALAT} that $L := \Kir(\IO)$ is a lattice in $\IO$.

Since $(u_r)_{0\leq r\leq 7}$ is an orthonormal basis of $\IO$, we have $n_\IO(\Ga(\IO),v_i) \subseteq \IZ$ for $1 \leq i \leq 4$, and a straightforward verification shows
\begin{align}
\label{NOVI} n_\IO(v_i) = 1, \quad n_\IO(v_1,v_2) = 0, \quad n_\IO(v_j,v_k) = 1 
\end{align}
for $1 \leq i,j,l \leq 4$, $j \neq l$, $\{j,l\} \neq \{1,2\}$. Combining this with the definition of $L$, we conclude $n_\IO(L) \subseteq \IZ$, so $L$ is a unital integral quadratic lattice in $\IO$.

Next we consider the family
\[
E_1 := (1_\IO,u_1,u_2,u_3,v_1,v_2,v_3,v_4)
\] 
of elements in $L$. From $u_4 = 2v_1 - 1_\IO - u_1 - u_2$, $u_5 = 2v_4 - 1_\IO - u_2 - u_3$, $u_7 = -2v_3 + 1_\IO + u_1 + u_3$, $u_6 = -2v_2 + u_3 + u_5 - u_7$ we deduce that $E_1$ spans $\IO$ as a real vector space, hence is a basis, and that the additive map from $\IZ^8$ to $L$ determined by $E_1$ is surjective, hence bijective. Thus $E_1$ is an $\IR$-basis of $\IO$  associated with $L$. 

We now proceed to determine the discriminant of $L$ by computing the determinant of $Dn_\IO(E_1)$. A straightforward verification using \eqref{NOVI} shows
\[
Dn_\IO(E_1) = \left(\begin{matrix}
2\cdot\Eins_4 & B \\
B^\trans  & D
\end{matrix}\right),
\]
where
\[
B = \left(\begin{matrix}
1 & 0 & 1 & 1 \\
1 & 0 & 1 & 0 \\
1 & 0 & 0 & 1 \\
0 & 1 & 1 & 1
\end{matrix}\right), \quad D = \left(\begin{matrix}
2 & 0 & 1 & 1 \\
0 & 2 & 1 & 1 \\
1 & 1 & 2 & 1 \\
1 & 1 & 1 & 2
\end{matrix}\right).
\]
Now we use \eqref{DEBLO} in the solution to Exc.~\ref{pr.DETFO} and obtain
\begin{align*}
\det\big(Dn_\IO(E_1)\big) =\,\,&\det(2\cdot\Eins_4)\det\big(D - B^\trans (\frac{1}{2}\cdot\Eins_4)B\big) = \det(2\cdot D - B^ tB) \\
=\,\,&\det(\left(\begin{matrix}
1 & 0 & 0 & 0 \\
0 & 3 & 1 & 1 \\
0 & 1 & 1 & 0 \\
0 & 1 & 0 & 1
\end{matrix}\right)) = 3 - 1 - 1 = 1. 
\end{align*} Thus $L$ is a unimodular integral quadratic lattice in $\IO$.

It remains to show that $L$ is \emph{not} a $\IZ$-structure of $\IO$, equivalently, that it is not closed under multiplication. To this end, we consider the product $v_1v_3$. A short computation gives
\begin{align}
\label{PROTH} v_1v_3 = \frac{1}{2}(u_1 + u_2 + u_3 + u_5).
\end{align} 
Assuming $v_1v_3 \in \Kir(\IO)$, we would find integers $\alpha_r,\beta_i$, $0 \leq r \leq 3$, $1 \leq i \leq 4$, such that
\begin{align}
\label{LICOV} v_1v_3 = \sum_{r=0}^3 \alpha_ru_r + \sum_{i=1}^4\beta_iv_i.
\end{align}
Comparing coefficients of $u_4$, we obtain $\beta_1 = 0$; comparing coefficients of $u_1$, we obtain $\frac{1}{2} = \alpha_1 + \frac{1}{2}\beta_3$, hence $\alpha_1 = \frac{1 - \beta_3}{2}$, \emph{and $\beta_3$ is odd.} Comparing coefficients of $u_2$, we obtain $\frac{1}{2} = \alpha_2 + \frac{1}{2}\beta_4$, hence $\alpha_2 = \frac{1 - \beta_4}{2}$, \emph{and $\beta_4$ is odd.} Comparing coefficients of $u_3$, we obtain $\frac{1}{2} = \alpha_3 + \frac{1}{2}(\beta_2 + \beta_3 + \beta_4)$, hence $\alpha_3 = -\frac{1}{2}(\beta_3 + \beta_4) + \frac{1 - \beta_2}{2}$, \emph{and $\beta_2$ is odd} since both $\beta_3,\beta_4$ are odd. And, finally, comparing coefficients of $u_5$, we obtain $\frac{1}{2} = \frac{1}{2}(\beta_2 + \beta_4)$, hence $\beta_2 + \beta_4 = 1$, which is a contradiction since $\beta_2,\beta_4$ are both odd. Thus $v_1v_3$ does not belong to $\Kir(\IO)$, as desired.
\end{sol}

\begin{sol}{pr.ALCO} \label{sol.ALCO}
The idea of the solution is to imitate our construction of the Dickson-Coxeter octonions with a different model of the Hamiltonian quaternions inside $\IO$ and a different vector $\bfq \in \IO$ in place of $\bfp$ such that the resulting $\IZ$-structure agrees with $R$. In doing so, repeated use will be made of the multiplcation rules for Cartan-Shouten bases as depicted in Fig.~\vref{fig.fano}. 

On a more formal level, these multiplication rules as described in (\ref{ss.DECS}.\ref{SQCS}) and (\ref{ss.DECS}.\ref{MICS}) show that the linear bijection $\vph\:\IO \to \IO$ determined by 
\begin{align*}
\vph(1_\IO) = 1_\IO, \quad \vph(u_r) = u_{r+2} &&(1 \leq r \leq 7,\;\text{indices} \bmod 7)
\end{align*}
is an automorphism of $\IO$. With $\IH \subseteq \IO$ as given in (\ref{ss.COZO}.\ref{HACS}) we therefore conclude that
\[
B := \vph(\IH) = \IR 1_\IO + \IR u_3 + \IR u_4 + \IR u_6 \subseteq \IO
\]
is a model of the Hamiltonian quaternions matching $1_\IH= 1_\IO$, $\bfi = u_3$, $\bfj = u_4$, $\bfk = u_6$ and
\[
\Ga(B) := \vph\big(\Ga(\IH)\big) = \IZ 1_\IO + \IZ u_3 + \IZ u_4 + \IZ u_6.
\]
Observing (\ref{ss.COZO}.\ref{BFEL}), we also put
\begin{align}
\label{KUFIP} \bfq := \vph(\bfp) = \frac{1}{2}(1_\IO + u_3 + u_4 + u_5)
\end{align}
and deduce from Thm.~\ref{t.COGC} that
\[
R^\prime := \vph\big(\Cox(\IO)\big) = \Ga(B) + \Ga(B)\bfq
\]
is a $\IZ$-structure of $\IO$ isomorphic to $\Cox(\IO)$. It therefore remains to show $R^\prime = R$.

We first note that $\Ga(\IO)$ is contained in $R$ (by definition) but also in $R^\prime$ (by Thm.~\ref{t.COGC}). Moreover, invoking (\ref{ss.COZO}.\ref{ONEL})--(\ref{ss.COZO}.\ref{FOUL}), we deduce
\begin{align*}
u_3\bfq =\,\,&\vph(u_1\bfp) = \frac{1}{2}(-1_\IO + u_2 + u_3 + u_6), \\
u_4\bfq =\,\,&\vph(u_2\bfp) = \frac{1}{2}(-1_\IO + u_4 - u_6 + u_7), \\
u_6\bfq =\,\,&\vph(u_4\bfp) = \frac{1}{2}(-u_1 - u_3 + u_4 + u_6).
\end{align*}
Now an inspection shows
\begin{align*}
v_1 + v_4\,\,&\equiv \frac{1}{2}(1_\IO + u_3 + u_4 + u_5) \equiv \bfq \bmod \Ga(\IO), \\
v_2 + v_4\,\,&\equiv \frac{1}{2}(1_\IO + u_2 + u_3 + u_6) \equiv u_3\bfq \bmod \Ga(\IO), \\ 
v_1 + v_2 + v_3 + v_4\,\,&\equiv \frac{1}{2}(1_\IO + u_4 + u_6 + u_7) \equiv u_4\bfq \bmod \Ga(\IO), \\
v_1 + v_2 + v_4\,\,&\equiv \frac{1}{2}(u_1 + u_3 + u_4 + u_6) \equiv u_6\bfq \bmod \Ga(\IO).
\end{align*}
Since 
\[
\vph(E^\prime) = (1_\IO,u_3,u_4,u_6,\bfq,u_3\bfq,u_4\bfq,u_6\bfq)
\]
by (\ref{t.COGC}.\ref{BECO}) is a basis of $\IO$ associated with $R^\prime$, the preceding relations imply $R^\prime \subseteq R$. But they also yield
\begin{align*}
u_4\bfq - u_6\bfq\,\,&\equiv v_3 \bmod \Ga(\IO), \\
u_6\bfq - u_3\bfq\,\,&\equiv v_1 \bmod \Ga(\IO), \\
u_6\bfq - \bfq\,\,&\equiv v_2 \bmod \Ga(\IO), \\
\bfq + u_3\bfq - u_6\bfq\,\,&\equiv v_4 \bmod \Ga(\IO)
\end{align*}
and hence imply $R \subseteq R^\prime$. This shows $R = R^\prime$ and the proof is complete.
\end{sol}

\begin{sol}{pr.HUCO} \label{sol.HUCO}
Let $(u_r)_{0\leq r\leq 7}$ be a Cartan-Schouten basis of $\IO$. In the language of Exc.~\ref{pr.ROOC} we have
\[
\bfq := \frac{1}{2}(1_\IO + u_3 + u_4 + u_6) = \frac{1}{2}(-\vep_1 + \vep_2 - \vep_3 - \vep_4 - \vep_5 + \vep_6 + \vep_7 + \vep_8),
\]
whence part (d) of that exercise shows that $\bfq$ is invertible in $\Cox(\IO)$. Since a straightforward verification yields $u_3u_4 = u_6$, $u_4u_6 = u_3$, $u_6u_3 = u_4$, we conclude from Thm.~\ref{t.HURO} that the assignments $1_\IH \mapsto 1_\IO$, $\bfi \mapsto u_3$, $\bfj \mapsto u_4$, $\bfh \mapsto \bfq$ determine an additive embedding $\Hur(\IH) \hookrightarrow \Cox(\IO)$ that preserves multiplication, hence is an embedding of $\IZ$-algebras.
\end{sol}

\solnsec{Section~\ref{s.EUAL}}

\begin{sol}{pr.FORMU} \label{sol.FORMU}
Computing the square of the matrix (\ref{ss.CASEN}.\ref{ELTH}) in a straighforward manner yields (\ref{FUIDRE.fig}.\ref{SQUARE}). Linearizing and dividing by $2$  immediately yields (\ref{FUIDRE.fig}.\ref{SYMREA}). Taking the trace of (\ref{FUIDRE.fig}.\ref{SYMREA}) and applying (\ref{ss.NOTRAD}.\ref{TRADE}), we obtain
\begin{align*}
T(x\bu y) =\,\,&\sum\big(\alpha_i\beta_i + \frac{1}{2}n_\IO(u_j,v_j) + \frac{1}{2}n_\IO(u_l,v_l)\big) \\
=\,\,&\sum\alpha_i\beta_i + \frac{1}{2}\big(\sum n_\IO(u_i,v_i) + \sum n_\IO(u_i,v_i)\big), 
\end{align*}
and (\ref{FUIDRE.fig}.\ref{TRSYRE}) follows. Writing
\begin{align*}
z = \sum\big(\gamma_ie_{ii} + w_i[jl]\big) &&(\gamma_i \in \IR,\;\;w_i \in \IO,\;\;1 \leq i \leq 3),
\end{align*}
we first note by (\ref{ss.NOTR}.\ref{TRANO}) and Exc.~\ref{pr.INVAS}~(b) that, for all $u,v,w \in \IO$, the expression
\begin{align}
\label{TCYC} n_\IO(\overline{uv},w) = t_\IO(uvw)\;\;\text{is invariant under cyclic permutations of the variables.}
\end{align}
Combining (\ref{FUIDRE.fig}.\ref{SYMREA}) with (\ref{FUIDRE.fig}.\ref{TRSYRE}), we therefore obtain 
\begin{align*}
T\big((x\bu y)\bu z\big) =\,\,&\sum\Big(\big(\alpha_i\beta_i + \frac{1}{2}n_\IO(u_j,v_j) + \frac{1}{2}n_\IO(u_l,v_l)\big)\gamma_i \\
\,\,&+ \frac{1}{2}n_\IO\big((\alpha_j + \alpha_l)v_i + (\beta_j + \beta_l)u_i + \overline{u_jv_l + v_ju_l},w_i\big)\Big) \\
=\,\,&\sum\Big(\alpha_i\beta_i\gamma_i + \frac{1}{2}(\gamma_j + \gamma_l)n_\IO(u_i,v_i) + \frac{1}{2}(\alpha_j + \alpha_l)n_\IO(v_i,w_i) \\
\,\,&+ \frac{1}{2}(\beta_j + \beta_l)n_\IO(w_i,u_i) + \frac{1}{2}\big(t_\IO(u_iv_jw_l) + t_\IO(w_iv_ju_l)\big)\Big).
\end{align*}
This expression obviously being symmetric in $x,z$, we end up with (\ref{ss.FUIDRE}.\ref{TASRE}).

We now turn to (\ref{ss.FUIDRE}.\ref{BST}). Combining (\ref{ss.NOTRAD}.\ref{TRADE}) with (\ref{ss.FUIDRE}.\ref{TRSYRE}) we deduce
\begin{align*}
T(x)T(y) - T(x\bu y) =\,\,&(\sum\alpha_i)(\sum\beta_i) - \sum\big(\alpha_i\beta_i + n_\IO(u_i,v_i)\big) \\
=\,\,&\sum\big(\alpha_i\beta_j + \beta_i\alpha_j - n_\IO(u_i,v_i)\big),
\end{align*}
which by (\ref{ss.NOTRAD}.\ref{BQUADE}) agrees with $S(x,y)$. Similarly, by (\ref{ss.FUIDRE}.\ref{SQUARE}), (\ref{ss.NOTRAD}.\ref{TRADE}), \ref{ss.NOTRAD}.\ref{QUADE}) and (\ref{ss.NOTRAD}.\ref{ADDE}) we obtain
\begin{align*}
x^2 - T(x)x + S(x)\Eins_3 =\,\,&\sum\Big(\big(\alpha_i^2 + n_\IO(u_j) + n_\IO(u_l)\big)e_{ii} + \big((\alpha_j + \alpha_l)u_i + \overline{u_ju_l}\big)[jl]\Big) \\
\,\,&- \big(\sum\alpha_i\big)\sum\big(\alpha_ie_{ii} + u_i[jl]\big) + \sum\big(\alpha_j\alpha_l - n_\IO(u_i)\big)\big(\sum e_{ii}\big) \\
=\,\,&\sum\Big(\big(\alpha_i^2 + n_\IO(u_j) + n_\IO(u_l) - \alpha_i^2 - \alpha_i\alpha_j - \alpha_l\alpha_i \\
\,\,&+ \alpha_i\alpha_j + \alpha_j\alpha_l + \alpha_l\alpha_i - n_\IO(u_i) - n_\IO(u_j) - n_\IO(u_l)\big)e_{ii} \\
\,\,&+ \big((\alpha_j + \alpha_l)u_i + \overline{u_ju_l} - \alpha_iu_i - \alpha_ju_i - \alpha_lu_i\big)[jl]\Big) \\
=\,\,&\sum\Big(\big(\alpha_j\alpha_l - n_\IO(u_i)\big)e_{ii} + \big(-\alpha_iu_i + \overline{u_ju_l}\big)[jl]\Big) = x^\sharp,  
\end{align*}
which completes the proof of (\ref{ss.FUIDRE}.\ref{SHASQ}).

Differentiating (\ref{ss.NOTRAD}.\ref{NORDE}) at $x$ in the direction $y$ and observing \eqref{TCYC}, we obtain
\begin{align*}
DN(x)(y) =\,\,&\beta_1\alpha_2\alpha_3 + \alpha_1\beta_2\alpha_3 + \alpha_1\alpha_2\beta_3 - \sum\big(\beta_in_\IO(u_i) + \alpha_in_\IO(u_i,v_i)\big) \\
\,\,&+ t_\IO(v_1u_2u_3) + t_\IO(u_1v_2u_3)  + t_\IO(u_1u_2v_3) \\
=\,\,&\sum\big(\alpha_j\alpha_l\beta_i - n_\IO(u_i)\beta_i - \alpha_in_\IO(u_i,v_i) + t_\IO(u_ju_lv_i)\big),
\end{align*}
while combining (\ref{ss.FUIDRE}.\ref{TRSYRE}) with (\ref{ss.NOTRAD}.\ref{ADDE}) and \eqref{TCYC} yields
\begin{align*}
T(x^\sharp\bu y) =\,\,&\sum\Big(\big(\alpha_j\alpha_l - n_\IO(u_i)\big)\beta_i + n_\IO\big(-\alpha_iu_i + \overline{u_ju_l}, v_i\big)\Big) \\
=\,\,&\sum\big(\alpha_j\alpha_l\beta_i - n_\IO(u_i)\beta_i - \alpha_in_\IO(u_i,v_i + t_\IO(u_ju_lv_i)\big),
\end{align*}
which proves the first equation of (\ref{ss.FUIDRE}.\ref{DIREN}), while the second one follows from (\ref{ss.FUIDRE}.\ref{SHASQ}).

In order to prove (\ref{FUIDRE.fig}.\ref{CUBRE}), we apply (\ref{FUIDRE.fig}.\ref{SYMREA}) and (\ref{ss.NOTRAD}.\ref{ADDE}) to compute
\begin{align*}
x\bu x^\sharp =\,\,&\sum\Big(\big(\alpha_i[\alpha_j\alpha_l - n_\IO(u_i)] + \frac{1}{2}n_\IO(u_j,-\alpha_ju_j + \overline{u_lu_i}) + \frac{1}{2}n_\IO(u_l,-\alpha_lu_l + \overline{u_iu_j})\big)e_{ii} \\
\,\,&+ \frac{1}{2}\big((\alpha_j + \alpha_l)[-\alpha_iu_i + \overline{u_ju_l}] + [\alpha_l\alpha_i - n_\IO(u_j) + \alpha_i\alpha_j - n_\IO(u_l)]u_i \\
\,\,&+ \overline{u_j[-\alpha_lu_l + \overline{u_iu_j}] + [-\alpha_ju_j + \overline{u_lu_i}]u_l}\big)[jl]\Big).
\end{align*}
Since the conjugation of $\IO$ is an involution and $(u_iu_j)\bar u_j = n_\IO(u_j)u_i$, $\bar u_l(u_lu_i) = n_\IO(u_l)u_i$ by Exc.~\ref{pr.INVAS}~(c), we conclude
\begin{align*}
x\bu x^\sharp = \sum\big(\alpha_1\alpha_2\alpha_3 - \sum_{r=1}^3\alpha_rn_\IO(u_r) + t_\IO(u_1u_2u_3)\big)e_{ii} = N(x)\Eins_3,
\end{align*}
and plugging in (\ref{FUIDRE.fig}.\ref{SHASQ}) gives (\ref{FUIDRE.fig}.\ref{CUBRE}). 

Finally, differentiating $x^3 = (x\bu x)\bu x$ at $x$ in the direction $y$, we find $(y\bu x)\bu x + (x\bu y)\bu x + (x\bu x)\bu y = 2x\bu(x\bu y) + x^2\bu y$. Performing the same procedure on the right-hand side of (\ref{FUIDRE.fig}.\ref{CUBRE}), we therefore get (\ref{FUIDRE.fig}.\ref{DICURE}). 
\end{sol} 

\begin{sol}{pr.HEMAO} \label{sol.HEMAO}
We have $\Her_1(\IO) \cong \IR^+$, while $\Her_2(\IO)$ sits in the euclidean Albert algebra via
\[
\Her_2(\IO) \cong \Big\{\left(\begin{matrix}
\alpha_1 & u_3 & 0 \\
\bar u_3 & \alpha_2 & 0 \\
0 & 0 & 0
\end{matrix}\right)\mid \alpha_1,\alpha_2 \in \IR,\;u_3 \in \IO\Big\}
\]
as a non-unital subalgebra.  By Thm.~\ref{t.HEUJO}, therefore, $\Her_n(\IO)$ is a Jordan algebra for $1 \leq n \leq 3$. It remains to show that it is \emph{not} a Jordan algebra for $n \geq 4$. We begin by differentiating the Jordan identity $u(u^2v) = u^2(u v)$ at $u$ in the direction $w$ and obtain 
\[
w(u^2v) + 2u\big((uw)v\big) = 2(uw)(uv) + u^2(wv).
\]
Differentiating at $u$ again, this time in the direction $x$, changing notation and rearranging terms, we end up with the fully linearized Jordan identity, which therefore holds in arbitrary real Jordan algebras. 

Now consider the commutative real algebra $J := \Her_n(\IO)$ for some integer $n \geq 4$. For $1 \leq i,j,l,m \leq n$ satisfying $i \ne j$, $l\ne m$ and $u,v \in \IO$, we will make use of the obvious relations
\begin{align*}
e_{ii}\bu u[lm] =\,\,&\frac{1}{2}(\delta_{il} + \delta_{im})u[lm] = u[lm]\bu e_{ii}, \\
u[ij] =\,\,&\bar u[ji], \\ 
u[ij]\bu v[jl] =\,\,&\frac{1}{2}(uv)[il] &&(j \ne l \ne i), \\
u[ij]\bu v[lm] =\,\,&0 &&(\{i,j\} \cap \{l,m\} = \emptyset). 
\end{align*}
Assuming now $J$ is a Jordan algebra, we combine the preceding relations with the fully linearized Jordan identity and conclude
\begin{align*}
\frac{1}{4}\big((uv)w\big)[14] =\,\,&\big(u[12]\bu v[23]\big)\bu\big((e_{33} + e_{44})\bu w[34]\big) \\
=\,\,&-\big(v[23]\bu(e_{33} + e_{44})\big)\bu\big(u[12]\bu w[34]\big) - \big((e_{33} + e_{44})\bu u[12]\big)\bu\big(v[23]\bu w[34]\big) \\
\,\,&+ u[12]\bu\Big(\big(v[23]\bu(e_{33} + e_{44})\big)\bu w[34]\Big) + v[23]\bu\Big(\big((e_{33} + e_{44})\bu u[12]\big)\bu w[34]\Big) \\
\,\,&+ (e_{33} + e_{44})\bu\big((u[12]\bu v[23]) \bu w[34]\big) \\
=\,\,&\frac{1}{8}\big(u(vw)\big)[14] + \frac{1}{8}\big((uv)w\big)[14].
\end{align*}
Thus we obtain the contradiction that $\IO$ is associative. Hence $J$ cannot be a Jordan algebra.
\end{sol}

\begin{sol}{pr.INVRE} \label{sol.INVRE}
(a) By hypothesis, $\deg(\mu_x) = \dim_\IR(\IR[x]) = 3$. Hence $\mu := \det(\bft\Eins_{\IR[x]} - L_x^0$) is monic of degree $3$ as well. Moreover, $\mu(x) = \mu(L_x^0)\Eins_3 = 0$, which proves the first assertion. The remaining ones now follow from (\ref{ss.MIPOD}.\ref{SQSQ}).

(b) Using the polynomial map
\[
f\:J \longrightarrow \bigwedge^3(J), \quad x \longmapsto \Eins_3 \wedge x \wedge x^2,
\]
we consider the set $X$ of all septuples
\begin{align} \label{SEPTUP}
(\alpha_0,\alpha_1,\alpha_2,\beta_0,\beta_1,\beta_2,x) \in \IR^6 \times J
\end{align}
such that, setting $u := \sum_{r=0}^2\alpha_rx^r$, $v := \sum_{r=0}^2\beta_rx^r$, the quantities $f(x), f(u), f(v), f(u\bu v)$ are all different from zero. Note that this condition is equivalent to $\IR[x] = \IR[u] = \IR[v] = \IR[u\bu v]$ having dimension $3$. We claim that $X \subseteq \IR^6 \times J$, which is obviously Zariski-open, is also not empty, hence Zariski-dense. Indeed, let $\xi_i,\eta_i,\zeta_i \in \IR$ for $i = 1,2,3$ such that the sets
\[
\{\xi_1,\xi_2,\xi_3\}, \quad \{\eta_1,\eta_2,\eta_3\}, \quad \{\zeta_1,\zeta_2,\zeta_3\}, \quad \{\eta_1\zeta_1,\eta_2\zeta_2,\eta_3\zeta_3\}
\]
all have cardinality $3$. Setting $x :=\sum\xi_ie_{ii}$, $u := \sum\eta_ie_{ii}$, $v := \sum\zeta_ie_{ii}$, we conclude $\IR[x] = \IR[u] = \IR[v] = \IR[u\bu v] = \sum\IR e_{ii}$, allowing us to write $u = \sum_{r=0}^2\alpha_rx^r$, $v = \sum_{r=0}^2\beta_rx^r$ for some $\alpha_r,\beta_r \in \IR$, $0\le r \le 2$, and it follows that the quantity \eqref{SEPTUP} belongs to $X$.

Hence it suffices to prove the first claim under the additional hypothesis that $\IR[x] = \IR[u] = \IR[v] = \IR[u\bu v]$ is three-dimensional. But then, since $\IR[x]$ is commutative associative, (a) gives
\[
N(u\bu v) = \det(L_{u\bu v}^0) = \det(L_u^0L_v^0) = \det(L_u^0)\det(L_v^0) = N(u)N(v). 
\]
It remains to show that the equation $N(x\bu y) = N(x)N(y)$ does not hold for all $x,y \in J$. To this end, put
\begin{align*}
x := \sum\alpha_ie_{ii}, \quad y := \sum v_i[jl] &&(\alpha_i \in \IR,\;v_i \in \ID,\;1 \leq i \leq 3).
\end{align*}
Then (\ref{ss.NOTRAD}.\ref{NORDE}) and (\ref{FUIDRE.fig}.\ref{SYMREA}) imply
\begin{align*}
N(x) =\,\,&\alpha_1\alpha_2\alpha_3, \\ 
N(y) =\,\,&t_\ID(v_1v_2v_3), \\
x\bu y =\,\,&\frac{1}{2}\sum(\alpha_j + \alpha_l)v_i[jl]
\end{align*}
hence
\[
N(x\bu y) = \frac{1}{8}(\alpha_2 + \alpha_3)(\alpha_3 + \alpha_1)(\alpha_1 + \alpha_2)t_\ID(v_1v_2v_3),
\]
which in general will be distinct from $N(x)N(y) = \alpha_1\alpha_2\alpha_3t_\ID(v_1v_2v_3)$. 

(c) If $x$ is invertible in $J$, then, by definition, it is so in $\IR[x]$, and $x^{-1} \in \IR[x]$ satisfies $x\bu x^{-1} = \Eins_3$. Now (b) implies $N(x)N(x^{-1}) =1$, hence $N(x) \neq 0$ as well as $N(x^{-1}) = N(x)^{-1}$. Conversely, let $N(x) \ne 0$. Combining (\ref{FUIDRE.fig}.\ref{SHASQ}) with (\ref{FUIDRE.fig}.\ref{CUBRE}), we obtain $x^\sharp \in \IR[x]$ and $x\bu x^\sharp = N(x)\Eins_3$. But this implies $x\bu y = \Eins_3$ with $y := N(x)^{-1}x^\sharp \in \IR[x]$, so $x$ is invertible in $J$, and the inverse has the desired form.

(d) By Zariski density and (c), it suffices to prove  both identities for $x$ invertible. Taking norms on both sides of the displayed equation in (c), we obtain $N(x)^{-1} = N(x^{-1}) = N(x)^{-3}N(x^\sharp)$, hence $N(x^\sharp) = N(x)^2$. In particular, $x^\sharp$ is invertible with inverse $x^{\sharp {-1}} = N(x^\sharp)^{-1}x^{\sharp\sharp} = N(x)^{-2}x^{\sharp\sharp}$. Hence
\[
x = (x^{-1})^{-1} = \big(N(x)^{-1}x^\sharp\big)^{-1} = N(x)x^{\sharp {-1}} = N(x)^{-1}x^{\sharp\sharp}, 
\] 
and  also the adjoint identity follows. 
\end{sol}

\begin{sol}{pr.AURE} \label{sol.AURE}
First let $\vph$ be an automorphism of $J$. Then $\vph$ fixes $\Eins_3$ and in order to prove $N(\vph(x)) = N(x)$ for $x \in J$, we may assume, by Zariski densitiy, that $\Eins_3 \wedge x \wedge x^2 \neq 0$ in $\biw^3(J)$. Since $\vph$ preserves powers, (\ref{FUIDRE.fig}.\ref{CUBRE}) implies
\begin{align*}
T(x)\vph(x)^2 - S(x)\vph(x) + N(x)\Eins_3 =\,\,&\vph(x^3) = \vph(x)^3  \\
=\,\,&T\big(\vph(x)\big)\vph(x)^2 - S\big(\vph(x)\big)\vph(x) + N\big(\vph(x)\big)\Eins_3,
\end{align*} 
and since $\Eins_3,\vph(x),\vph(x)^2$ are linearly independent, we conclude $N(\vph(x)) = N(x)$, as claimed. 

Conversely, suppose $\vph$ preserves $\Eins_3$ and $N$. Then $\vph$ also preserves the directional derivative of $N$ at $x$ in the direction $y$, i.e., in view (\ref{FUIDRE.fig}.\ref{DIREN}) we have
\begin{align}
\label{TEPH} T\big(\vph(x)^\sharp\bu\vph(y\big) = DN\big(\vph(x)\big)\big(\vph(y)\big) = DN(x)(y) = T(x^\sharp\bu y).
\end{align}
Setting $x = \Eins_3$ implies $T \circ \vph = T$, while setting $y = \Eins_3$ implies $S \circ \vph  = S$. Taking traces in (\ref{FUIDRE.fig}.\ref{SHASQ}), on the other hand, we conclude $T(x^2) = T(x)^2 - 2S(x)$. Replacing $x$ by $\vph(x)$ and applying what we have just proved yields $T(\vph(x)^2) = T(x^2)$, hence $T(\vph(x)\bu\vph(y)) = T(x\bu y)$ after linearization. Substituting $x^\sharp$ for $x$ and combining with \eqref{TEPH}, we deduce $T(\vph(x)^\sharp\bu\vph(y)) = T(x^\sharp\bu y) = T(\vph(x^\sharp)\bu\vph(y))$. But (\ref{FUIDRE.fig}.\ref{TRSYRE}) shows that the bilinear form $T(x\bu y)$ is positive definite, hence non-degenerate, and we conclude that $\vph$ preserves adjoints. By (\ref{FUIDRE.fig}.\ref{SHASQ}), therefore, $\vph$ preserves squares, and since the algebra $J$ is commutative, it must be an isomorphism.

As to the second part, the linear map $\vph$ is clearly bijective and preserves the unit. Assume that $\vph$ is an automorphism of $J$. Then (\ref{FUIDRE.fig}.\ref{SQUARE}) yields
\[
e_{22} + e_{33} = \vph(1_\ID[23]^2) = \vph(1_\ID[23])^2 = u^{-1}[23]^2 = n_\ID(u)^{-1}(e_{22} + e_{33}), 
\]
and we conclude $n_\ID(u) = 1$. Conversely, let this be so. By the first part of the problem,  it suffices to show that $\vph$ preserves the norm. To this end, we compute
\begin{align*}
N\big(\vph(x)\big) =\,\,&N\big(\sum\alpha_ie_{ii} + (u^{-1}u_1)[23] + (u_2u^{-1})[31] + (uu_3u)[12]\big) \\
=\,\,&\alpha_1\alpha_2\alpha_3 - \alpha_1n_\ID(u^{-1}u_1) - \alpha_2n_\ID(u_2u^{-1}) - \alpha_3n_\ID(uu_3u)\ + \\ 
\qquad &t_\ID\big((u^{-1}u_1)(u_2u^{-1})(uu_3u)\big).
\end{align*} 
But $u$ has norm $1$, $n_\ID$ permits composition, and the Moufang identities (Exc.~\ref{pr.MIMOU}) combined with Exc.~\ref{pr.INVAS}~(b), (c) give
\begin{align*}
t_\ID\big((u^{-1}u_1)(u_2u^{-1})(uu_3u)\big) =\,\,&t_\ID\Big(\big(u^{-1}(u_1u_2)u^{-1}\big)\big(uu_3u\big)\Big) \\
=\,\,&t_\ID\Big(u^{-1}\big((u_1u_2)[u^{-1}(uu_3u)]\big)\Big) \\
=\,\,&t_\ID\Big(u^{-1}\big((u_1u_2)(u_3u)\big)\Big) = t_\ID\Big(\big(u^{-1}(u_1u_2)\big)(u_3u)\Big) \\
=\,\,&t_\ID\Big(\big([u^{-1}(u_1u_2)]u_3\big)u\Big)= t_\ID\Big(u\big([u^{-1}(u_1u_2)]u_3\big)\Big) \\
=\,\,&t_\ID\Big(\big(u[\bar u(u_1u_2)]\big)u_3\Big) = t_\ID(u_1u_2u_3).
\end{align*}
Hence $N(\vph(x)) = N(x)$, and the problem is solved.
\end{sol}

\solnsec{Section~\ref{s.ZUNREJ}}

\begin{sol}{pr.UORE} \label{sol.UORE}
(a) This follows by a straightforward computation. 

(b) Linearizing (\ref{FUIDRE.fig}.\ref{SHASQ}) and applying (\ref{FUIDRE.fig}.\ref{BST}), we obtain
\[
x \times y = 2x\bu y -T(x)y - T(y)x + T(x)T(y)\Eins_3 - T(x\bu y)\Eins_3,
\]
and combining this with (\ref{FUIDRE.fig}.\ref{SHASQ}), (\ref{FUIDRE.fig}.\ref{DIREN}), (\ref{FUIDRE.fig}.\ref{DICURE}) yields
\begin{align*}
x^\sharp \times y =\,\,&2x^\sharp\bu y - T(x^\sharp)y - T(y)x^\sharp + T(x^\sharp)T(y)\Eins_3 - T(x^\sharp\bu y)\Eins_3 \\
=\,\,&2x^2\bu y - 2T(x)x\bu y + 2S(x)y - S(x)y - T(y)x^\sharp + S(x)T(y)\Eins_3 - DN(x)(y)\Eins_3 \\
=\,\,&2x^2\bu y - 2T(x)x\bu y + S(x)y - T(y)x^2 + T(x)T(y)x -S(x)T(y)\Eins_3 \\
\,&+ S(x)T(y)\Eins_3 - x^2\bu y + 2T(x)x\bu y + T(y)x^2 - S(x)y -S(x,y)x -2x\bu(x\bu y) \\
=\,\,&-U_xy + T(x\bu y)x,
\end{align*}
as claimed.
\end{sol}

\begin{sol}{pr.OPLSP} \label{sol.OPLSP}
Let $A,B$ be commutative real algebras and $\vph\:A \to B$ a linear map. If $\vph$ is a homomorphism, then it clearly preserves squares. Conversely, let this be so. Then $\vph(x^2) = \vph(x)^2$ for all $x \in A$, and linearizing gives 
\[
\vph(xy) = \frac{1}{2}\Big(\vph\big((x + y)^2 - x^2 - y^2\big)\Big) = \frac{1}{2}\Big(\big(\vph(x) + \vph(y)\big)^2 - \vph(x)^2 - \vph(y)^2\Big) = \vph(x)\vph(y).
\]
Thus $\vph$ is a homomorphism.

Now define the \emph{left multiplication} of $\IO$ as the linear map $L\:\IO^+ \to \End_\IR(\IO)^+$ sending $x \in \IO^+$ to $L_x\:\IO \to \IO$, $y \mapsto xy$. Since $L_{x^2} = L_x^2$ for all $x \in \IO$ by \ref{ss.ALTGC}, $L$ preserves squares  and hence is a homomorphism of commutative algebras. It is also injective since $L_x = 0$ implies $x = L_x1_\IO = 0$. Thus $L$ maps $\IO^+$ isomorphically onto a subalgebra of $\End_\IR(\IO)^+$, and we conclude from \ref{ss.SPEXRE}~(a) that $\IO^+$ is a special Jordan algebra. As in Exc.~\ref{pr.UORE}~(a), the verification of the formula for the $U$-operator is straightforward: for $x,y \in \IO$ we obtain, using alternativity \ref{ss.ALTGC},
\begin{align*}
U_xy =\,\,&2x\bu (x\bu y) - x^2\bu y = \frac{1}{2}\big(x(xy + yx) + (xy + yx)x - x^2y - yx^2\big) \\
=\,\,&\frac{1}{2}\big(x(xy) - x^2y + x(yx) + (xy)x + (yx)x - yx^2\big) =xyx,
\end{align*} 
as claimed. Finally, the Jordan algebra $\IO^+$ is \emph{not} euclidean since, given a Cartan-Shouten basis $(u_i)_{0\leq i\leq 7}$ of $\IO$, we have $1_\IO^2 + u_i^2 = 0$ for $i = 1,\dots,7$.
\end{sol}

\begin{sol}{pr.GAULI} \label{sol.GAULI}
For $0 \leq i \leq n$ we have $e_0\bu e_i = \frac{1}{2}(1_Ae_i + e_i1_A) = e_i$. Now let $1 \leq i,j \leq n$, $i \neq j$. By orthnormality, $e_i$ has trace $0$ and norm $1$, which by (\ref{ss.NOTR}.\ref{QUAGC}) implies $e_i\bu e_i = e_i^2 = -1_A = -e_0$, while (\ref{ss.NOTR}.\ref{BLAGC}) implies
\[
e_i\bu e_j = \frac{1}{2}e_i \circ e_j = \frac{1}{2}\big(t_\ID(e_i)e_j + t_\ID(e_j)e_i - n_\ID(e_i,e_j)\big) = 0.
\]
Thus $\Lambda$ is closed under the multiplication of $\ID^+$ and hence is a linear $\IZ$-structure of that algebra. But note that, in general, $\Lambda$ will not be closed under the multiplication of $\ID$. Finally, by Ex.~\ref{e.GAI} and Exc.~\ref{pr.CHACS}, the Gaussian integers of $\ID$ are the $\IZ$-linear span of some orthonormal basis of $\ID$, yielding the final assertion of the problem.
\end{sol} 

\begin{sol}{pr.IDREJO} \label{sol.IDREJO}
Exc.~\ref{pr.INVRE}~(b) implies $N(c) = N(c^2) = N(c)^2$, hence $N(c) \in \{0,1\}$. If $N(c) = 1$, then $c$ is invertible (Exc.~\ref{pr.INVRE}~(c)), and the relations $\Eins_3 - c \in \IR[c]$, $c(\Eins_3 - c) = 0$ yield $c = \Eins_3$, a contradiction. Thus $N(c) = 0$. Using this and (\ref{FUIDRE.fig}.\ref{CUBRE}), we deduce $c = c^3 = (T(c) - S(c))c$, hence
\begin{align}
\label{IDCU} S(c) = T(c) - 1.
\end{align}
Now (\ref{FUIDRE.fig}.\ref{SHASQ}) yields $c^\sharp = (1 - T(c))c + S(c)\Eins_3 =S(c)(\Eins_3 - c)$, hence $S(c) = T(c^\sharp) = S(c)(3 - T(c)) = -S(c)^2 + 2S(c)$. Thus $S(c)^2 = S(c)$, forcing $S(c) \in \{0,1\}$, and \eqref{IDCU} implies the assertion.
\end{sol}

\begin{sol}{pr.EXPNO} \label{sol.EXPNO}
Setting
\[
x = \sum\big(\alpha_ie_{ii} + u_i[jl]\big), \quad y = \sum\big(\beta_ie_{ii} + v_i[jl]\big)
\]
with $\alpha_i, \beta_i \in \IR$, $u_i,v_i \in \ID$, $1 \leq i \leq 3$ and applying (\ref{ss.NOTRAD}.\ref{NORDE}), (\ref{ss.NOTRAD}.\ref{TRADE}), (\ref{FUIDRE.fig}.\ref{TRSYRE}) and the fact that the expression $t_\IO(uvw) = n_\IO(\overline{uv},w)$ remains invariant under cyclic permutations of the variables, we compute
\begin{align*}
N(x + y) =\,\,&(\alpha_1 + \beta_1)(\alpha_2 + \beta_2)(\alpha_3 + \beta_3) - \sum(\alpha_i + \beta_i)n_\IO(u_i + v_i) \\
\,\,&+ t_\IO\big((u_1 + v_1)(u_2 + v_2)(u_3 + v_3)\big) \\
=\,\,&\alpha_1\alpha_2\alpha_3 + \sum\alpha_j\alpha_l\beta_i + \sum\alpha_i\beta_j\beta_l + \beta_1\beta_2\beta_3 \\
\,\,&-\sum\big(\alpha_in_\IO(u_i) + \alpha_in_\IO(u_i,v_i) + \alpha_in_\IO(v_i) + \beta_in_\IO(u_i) + \beta_in_\IO(u_i,v_i) \\
\,\,&+ \beta_in_\IO(v_i)\big) + t_\IO(u_1u_2u_3) + \sum t_\IO(u_ju_lv_i) + \sum t_\IO(u_iv_jv_l) + t_\IO(v_1v_2v_3) \\
=\,\,&\big(\alpha_1\alpha_2\alpha_3 - \sum\alpha_in_\IO(u_i) + t_\IO(u_1u_2u_3)\big) \\
\,\,&+\sum\Big(\big(\alpha_j\alpha_l - n_\IO(u_i)\big)\beta_i + n_\IO(-\alpha_iu_i + \overline{u_ju_l},v_i)\Big) \\
\,\,&+\sum\Big(\alpha_i\big(\beta_j\beta_l - n_\IO(v_i)\big) + n_\IO(u_i,-\beta_iv_i + \overline{v_jv_l})\Big) \\
\,\,&+\big(\beta_1\beta_2\beta_3 - \sum\beta_in_\IO(v_i) + t_\IO(v_1v_2v_3)\big) \\
=\,\,&N(x) + T(x^\sharp\bu y) + T(x\bu y^\sharp) + N(y),
\end{align*}
as desired.
\end{sol}



\solnchap{Solutions for Chapter~\ref{c.FOUND}}

\solnsec{Section~\ref{s.LAN}}

\begin{sol}{pr.HOSI} \label{sol.HOSI}
For $x \in A$ and $\vep = \pm$, we put
\[
M_\vep(x) := \{y \in A \mid f(xy) = \vep f(x)f(y)\},
\]
which is a submodule of $A$. By hypothesis we have $A = M_+(x) \cup M_-(x)$, so some $\vep = \pm$ has $M_\vep(x) = A$. Now put
\[
N_\vep := \{x \in A \mid \forall y \in A\,: f(xy) = \vep f(x)f(y)\} = \{x \in A \mid M_\vep(x) = A\},
\]
which is again a submodule of $A$, and from what we have just proved, we obtain $A = N_+ \cup N_-$. This gives $N_\vep = A$ for some $\vep = \pm$, and if $\vep = +$ (resp. $\vep = -$), $f$ (resp. $-f$) is a homomorphism from $A$ to $B$.
\end{sol}

\begin{sol}{pr.NILRAD} \label{sol.NILRAD}
For any subset $X \subseteq A$ and any
homomorphism $\varphi\:A \to B$ of $k$-algebras, one makes the
following observations by straightforward induction:
\begin{enumerate}[(a)]
\item $\varphi(\Mon_m(X)) = \Mon_m(\varphi(X))$ for all integers
$m > 0$.

\item $\Mon_m(Y) \subseteq \Mon_{mn}(X)$ for all $Y
\subseteq \Mon_n(X)$ and all integers $m,n > 0$.
\end{enumerate}
If $x \in A$ is nilpotent, then (a) implies that $\varphi(x) \in B$
is nilpotent. Hence if $A$ is nil, so are $I$ and $A^\prime := A/I$.
Conversely, suppose $I$ and $A^\prime$ are nil. Writing $\pi\:A \to
A^\prime$ for the canonical epimorphism, any $x \in A$ makes $\pi(x)
\in A^\prime$ nilpotent, so by (a), $\Mon_n(\{x\})$ meets $I$ for
some integer $n > 0$. But since $I$ is nil, we conclude that some $y
\in \Mon_n(\{x\})$ is nilpotent, leading to a positive integer $m$
such that $0 \in \Mon_{mn}(\{x\})$ by (b). Thus $x$ is nilpotent,
forcing $A$ to be a nil algebra. Now write $\mfn := \Nil(A)$ for the
sum of all nil ideals in $A$. Given $x \in \mfn$, there are finitely
many nil ideals $\mfn_1,\dots,\mfn_r \subseteq A$ such that $x \in
\mfn_1 + \cdots + \mfn_r$. Thus we only have to show that the sum of
finitely many nil ideals in $A$ is nil. Arguing by induction, we are
actually reduced to the case $r = 2$. But then the isomorphism
$(\mfn_1 + \mfn_2)/\mfn_1 \cong \mfn_2/(\mfn_1 \cap \mfn_2)$
combines with what we have proved earlier to yield the assertion. To establish the final statement of the problem, we require the following standard result. 

\begin{lem*} 
Every finitely generated nil ideal $I \subseteq k$ is nilpotent, i.e., there exists a positive integer $n$ such that $I^n = \{0\}$.
\end{lem*}

\begin{proof}
Assume $\alpha_1,\dots,\alpha_m \in I$ generate $I$ as an ideal and let $p \in \IZ$, $p > 0$, satisfy $\alpha_i^p = 0$ for all $i = 1,\dots,m$. For any positive integer $n$, the ideal $I^n \subseteq k$ is generated by the expressions $\alpha_{i_1}\cdots\alpha_{i_n}$, where $1 \leq i_1 \leq \cdots \leq i_n \leq m$. If for all $i = 1,\dots,m$,
\[
q_i := \vert\{j \in \IZ\vert 1 \leq j \leq n,\;i_j = i\}\vert \leq p - 1,
\] 
then $n = \sum_{i=1}^m q_i \leq m(p - 1)$. Thus for $n > m(p - 1)$, some $i = 1,\dots,m$ has $q_i \geq p$. But this implies first $\alpha_{i_1}\cdots\alpha_{i_n} = 0$ and then $I^n = \{0\}$. 
\end{proof}

Since $\Nil(k)A \subseteq A$ is an ideal, it suffices to show that it consists entirely of nilpotent elements, so let $x \in \Nil(k)A$. Then $x \in IA$ for some finitely generated nil ideal $I \subseteq A$, and the lemma implies $x^n \in I^nA = \{0\}$ for some positive integer $n$. The assertion follows.
\end{sol}

\begin{sol}{pr.ID} \label{sol.ID}
(a) Since $x^n$ is a linear combination of
powers $x^i$, $d \leq i < n$, the $k$-module $k_d[x]$ is finitely
generated. Moreover, right multiplication by $x$ in the power
associative algebra $A$ gives a linear map $\varphi\:k_d[x] \to
k_d[x]$ whose image is $k_{d+1}[x]$. But $k_{d+1}[x] = k_d[x]$ since
$\alpha_d \in k^\times$. Therefore, $\varphi$ is surjective,
hence bijective (by Prop.~\ref{p.SURBI}), and so is $\varphi^d$. This proves existence and
uniqueness of $c \in k_d[x]$ satisfying $cx^d = x^d$, and the
relation $c^2x^d = c(cx^d) = cx^d$ shows that $c$ is an idempotent.

(b) Let $x \in A$ be a pre-image of $c^\prime$ under $\vph$. Then
$x^2 - x$ is nilpotent, so some integer $n > 0$ has $0 = (x^2 - x)^n
= x^{2n} - \cdots + (-1)^nx^n$. Now (a) yields an idempotent $c \in
k_n[x]$ satisfying $cx^n = x^n$. Writing $\bar u := \vph(u)$ for $u
\in A$, we conclude $\bar cc^\prime = \bar cc^{\prime n} = \bar
c\bar x^n = \bar x^n = c^\prime$. On the other hand, $c =
\sum_{i\geq n}\alpha_ix^i$ for some scalars $\alpha_i \in k$, $i
\geq n$, so with $\alpha := \sum_{i\geq n}\alpha_i$ we obtain $\bar
c = \bar\alpha c^\prime = \bar\alpha c^{\prime 2} = \bar cc^\prime =
c^\prime$.
\end{sol}

\solnsec{Section~\ref{s.UNAL}}

\begin{sol}{pr.IDDIRSUM} \label{sol.IDDIRSUM}
If $I_j$ are ideals in $A_j$ for $1 \leq j
\leq n$, then $I := I_1 \times \cdots \times I_n$ is an ideal in $A$
since $AI + IA = \prod_j (A_jI_j + I_jA_j) \subseteq \prod_j I_j = I$.
Conversely, let $I$ be an ideal in $A$ and $e_j$ the identity
element of $A_j$ for $1 \leq j \leq n$. If we identify $A_j \subseteq A$ through the $j$-th factor, then $I_j = e_jI = Ie_j$ is
clearly an ideal in $A_j$ such that $I = \prod_j I_j$.

The assertion \emph{fails} without the assumption that $A$ be
unital: suppose the $A_j$ have trivial multiplication. Then so has
$A$, and \emph{any} decomposition of $A$ into the direct product of $n$
submodules not related to the $A_j$ in any way is a decomposition
into a direct product of ideals.

To complete the solution, let us assume that the $A_j$, $1 \leq j \leq n$, are simple $k$-algebras. By the first part of the exercise, the ideals of $A$ have the form $\prod_{j=1}^n I_j$, where either $I_j = A_j$ or $I_j = \{0\}$, for all $ j = 1,\dots,n$. Hence the $A_j$, $1 \leq j \leq n$, are precisely the minimal ideals of $A$. This description being independent of the decomposition chosen, the assertion follows.
\end{sol}

\begin{sol}{pr.ALGEL} \label{sol.ALGEL}(a) The polynomials $\mu_1,\dots,\mu_r$ have greatest common divisor $1$. Hence $f_1,\dots,f_r \in F[\bft]$ with the desired properties exist. We therefore conclude $\sum_{i=1}^r c_i = 1_A$. Moreover, for $1 \leq i,j \leq r$, $i \neq j$, the polynomial $\mu_i\mu_j$ is a multiple of $\mu_x$, hence kills $x$. This proves $c_ic_j = 0$ and then $c_i = c_i1_A =\sum_{l=1}^r c_ic_l = c_i^2$. Thus $(c_1,\dots,c_r)$ is a complete orthogonal system of idempotents in $R := F[x]$. Now put
\begin{align}
\label{MINIL} v := x - \sum_{i=1}^r \alpha_ic_i = \sum_{i=1}^r (x - \alpha_i1_A)c_i \in R. 
\end{align}
Then the first relation of \eqref{MIDEC} (in the exercise) trivially holds and 
\[
v^n = \sum_{i=1}^r (x - \alpha_i1_A)^nc_i = \Big(\sum_{i=1}^r(\bft - \alpha_i)^n\mu_if_i\Big)(x).
\]
For $1 \leq i \leq r$,
\[
(\bft - \alpha_i)^n\mu_if_i = (\bft - \alpha_i)^{n_i}\mu_i(\bft - \alpha_i)^{n-n_i}f_i = \mu_x(\bft - \alpha_i)^{n-n_i}f_i
\]
is a multiple of $\mu_x$ and hence kills $x$. This shows $v^n = 0$, and the proof of \eqref{MIDEC} (in the exercise) is complete.

(b) For $I \subseteq \{1,\dots,r\}$, the quantity $c_I$ is clearly an idempotent in $R$. Conversely, let $c$ be an idempotent in $R$. Then $c = f(x)$ for some $f \in F[\bft]$. Since $R$ is commutative associative, we have
\begin{align}
\label{NILCAS} \Nil(R) = \{u \in R \mid \text{$u$ is nilpotent}\}.
\end{align}
Hence \eqref{MIDEC} (in the exercise) implies, for all $m \in \IN$,
\[
x^m = \sum_{i=1}^m \alpha_i^mc_i + v_m, \quad v_m \in \Nil(R)
\]
and then
\[
g(x) = \sum_{i=1}^r g(\alpha_i)c_i + v_g, \quad v_g \in \Nil(R),
\]
for all $g \in F[\bft]$. In particular, setting $\beta_i := f(\alpha_i)$ for $1 \leq i \leq r$,
\[
c = f(x) = \sum_{i=1}^r \beta_ic_i + v_f, \quad c^2 = f^2(x) = \sum_{i=1}^r \beta_i^2c_i + v_{f^2}. 
\]
But $c = c^2$, so $\sum(\beta_i^2 - \beta_i)c_i = v_{f^2} - v_f$ is nilpotent, which obviously implies $\beta_i^2 = \beta_i$, hence $\beta_i \in \{0,1\}$ for $1 \leq i \leq r$. Put
\[
I := \{i \in \IZ\mid 1 \leq i \leq r,\;\beta_i = 1\}.
\]
Then, setting $w := v_f \in \Nil(R)$,
\[
c_I + w = c = c^2 = c_I + 2c_Iw + w^2,
\]
and we conclude
\begin{align}
\label{CECI} (1_A - 2c_I)w = w^2, \quad (1_A - 2c_I)^2 = 1_A.
\end{align}
In particular, $1_A - 2c_I$ is invertible in $R$. Let $m \in \IZ$ satisfy $m \geq 1$ and $w^m = 0 \neq w^{m-1}$. Assuming $m \geq 2$, we conclude from \eqref{CECI} (in the solution) that $(1 - 2c_I)w^{m-1} = w^m = 0$, hence $w^{m-1} = 0$, a contradiction. Thus $m = 1$, forcing $w = 0$ and $c = c_I$.

(c) (i) $\Rightarrow$ (ii). Let $1 \leq i \leq r$ and put $v := ((\bft - \alpha_i)\mu_i)(x) \in R$. Then
\[
v^{n_i} = \big((\bft - \alpha_i)^{n_i}\mu_i^{n_i}\big)(x) = (\mu_x\mu_i^{n_i-1})(x) = 0,
\]
and (i) implies $v = 0$.Thus $\mu_x = (\bft - \alpha_i)^{n_i}\mu_i$ divides $(\bft - \alpha_i)\mu_i$, forcing $n_i = 1$ since $\mu_i$ is not divisible by $\bft - \alpha_i$.

(ii) $\Rightarrow$ (iii). This follows immediately from \eqref{MIDEC} (in the exercise).

(iii) $\Rightarrow$ (i). Let $v \in \Nil(R)$. Then
\[
v = f(x) = \sum_{i=1}^rf(\alpha_i)c_i
\]
for some $f \in F[\bft]$ is nilpotent, which can happen only if $f(\alpha_1) = \cdots = f(\alpha_r) = 0$ and hence amounts to $v = 0$.

(d) Write $\mu_x = \sum_{j=0}^d \gamma_j\bft^j$ with $d = \sum_{i=1}^r n_i$ and $\gamma_0,\dots,\gamma_d \in F$. Then the condition $\alpha_i \neq 0$ for all $i = 1,\dots,r$ is equivalent to $\gamma_0 = \mu_x(0) \neq 0$, hence implies $xy = 1_A$ with $y = -\gamma_0^{-1}\sum_{j=1}^d\gamma_jx^{j-1} \in R$, so $x$ is invertible in $R$. Conversely, let this be so and assume $\gamma_0 = 0$. Then $x\sum_{j=1}^d\gamma_jx^{j-1} = \mu_x(x) = 0$, and since $x$ is invertible in $R$, the polynomial $\sum_{j=1}^d\gamma_j\bft^{j-1} \in  F[\bft]$ of degree at most $d - 1$ kills $x$, in contradiction to $\mu_x$ being the minimum polynomial of $x$. Thus $\gamma_0 \neq 0$. 

(e) Let $\beta_i \in F^\times$  satisfy $\beta_i^2 = \alpha_i$ for $1 \leq i \leq r$. Setting $z := \sum_{i=1}^r \beta_ic_i$, we obtain $z^{-2}x = 1_A + w$, for some $w \in \Nil(R)$. Thus we may assume $x = 1_A + w$ and $w^m = 0 \neq w^{m-1}$ for some positive integer $m$. The case $m = 1$ being obvious, we may actually assume $m \geq 2$. Let
\begin{align*}
y = \sum_{j=0}^{m-1} \gamma_jw^j \in R &&(\gamma_0\dots,\gamma_{m-1} \in F).
\end{align*} 
The minimum polynomial of $w$ is $\bft^m$, and we conclude that $1_A,w,\dots,w^{m-1}$ are linearly independent over $F$. Hence $y^2 = x$ if and only if
\[
1_A + w = \sum_{j=0}^{m-1}\big(\sum_{l=0}^j \gamma_l\gamma_{j-l}\big)w^j,
\]                                                                       
which in turn is equivalent to
\begin{align*}
\gamma_0^2 =2\gamma_0\gamma_1 = 1, \quad \sum_{l=0}^j \gamma_l\gamma_{j-l} = 0 &&(2 \leq j < m).
\end{align*}
This system of quadratic equations can be solved recursively by
\begin{align*}
\gamma_0 = 1, \quad \gamma_1 = \frac{1}{2}, \quad \gamma_j = -\frac{1}{2}\sum_{l=1}^{j-1} \gamma_l\gamma_{j-l} &&(2 \leq j < m).
\end{align*}
This proves the assertion.

\begin{rmk*} 
The equation $y^2 = x = 1_A + w$, $w \in R$ nilpotent, can be solved more concisely  by 
\[
y = \sum_{l=0}^\infty \binom{\frac{1}{2}}{l}w^l,
\]
combined with the observation that $\binom{\frac{1}{2}}{l} \in \IQ$ has exact denominator a power of $2$ (depending on $l$).
\end{rmk*}
\end{sol}

\begin{sol}{pr.CENIDSPLIDE} \label{sol.CENIDSPLIDE}
We write $\mfD$ (resp. $\mfE$) for the set
of decompositions of $A$ into the direct sum of $n$ complementary
ideals (resp. of complete orthogonal systems of $n$ central
idempotents) in $A$. Let $(e_j)_{1\leq j\leq n} \in \mfE$. For $1
\leq j,l \leq n$, centrality implies $(Ae_j)(Ae_l) =A^2(e_je_l)$,
which is zero for $j \neq l$ and $Ae_j$ for $j = l$. Moreover, $x =
\sum xe_j \in \sum Ae_j$ for every $x \in A$, and given $x_j \in A$
for $1 \leq j \leq n$ satisfying $\sum x_je_j = 0$, we multiply with
$e_l$, $1 \leq l \leq n$, and deduce $x_le_l = 0$. All this boils
down to $A= (Ae_1) \oplus \cdots \oplus (Ae_n)$ as a direct sum of
ideals. Thus the assignment $(e_j) \mapsto (Ae_j)$ gives a map from
$\mfE$ to $\mfD$.

Conversely, let $(I_j)_{1\leq j\leq n} \in \mfD$. Then each $I_j$,
$1 \leq j \leq n$, is a unital $k$-algebra in its own right, and
setting $e_j := 1_{A_j}$, we deduce $\sum e_j = 1_A$, $e_je_l =
\delta_{jl}e_j$ for $1 \leq j,l \leq n$, so $(e_j)_{1\leq j\leq n}$
is a complete orthogonal system of idempotents in $A$, which are
clearly central since $\Cent(A) = \Cent(I_1) \oplus \cdots \oplus
\Cent(I_n)$ as a direct sum of ideals. Thus the assignment
$(I_j)_{1\leq j\leq n} \mapsto (e_j)_{1\leq j\leq n}$ determines a
map from $\mfD$ to $\mfE$, and it is readily checked that the two
maps constructed are inverse to each other.
\end{sol}

\begin{sol}{pr.PRIMID} \label{sol.PRIMID}
(a) (i) $\Rightarrow$ (ii). Let $d \in Rc$ be a non-zero idempotent. Then $cd = d$ and $c = d + \pd$, where $\pd := c - d\in Rc$ is an idempotent orthogonal to $d$. Since $c$ is primitive, this implies $\pd = 0$, i.e., $d = c$.

(ii) $\Rightarrow$ (iii). Suppose $x \in Rc$ is not nilpotent. Since $F$ is a field, we conclude from Exc.~\ref{pr.ID}~(a) that $F_1[x] \subseteq Rc$ contains a non-zero idempotent. By (ii), therefore, $c \in F_1[x]$, and there are $n \in \IZ$, $n > 0$, as well as $\alpha_1,\dots,\alpha_n \in F$ satisfying $c = \sum_{i=1}^n\alpha_ix^i = xy$, $y := \alpha_1c + \sum_{i=1}^{n-1}\alpha_{i+1}x^i \in Rc$. Hence $x$ is invertible in $Rc$.

(iii) $\Rightarrow$ (i). Assume $c_1,c_2 \in R$ are non-zero orthogonal idempotents having $c = c_1 + c_2$. Then $cc_i = c_i$, hence $c_i \in Rc$, for $i = 1,2$. But $c_i$, being a non-zero idempotent, cannot be nilpotent. By (iii), therefore, $c_i$ is invertible in $Rc$, which contradicts the relation $c_1c_2 = 0$. 

\smallskip

(b) Non-zero orthogonal idempotents are linearly independent. Hence, as $R$ is finite-dimensional over $F$, there exists a unique maximal number $n$ such that $c$ splits into the orthogonal sum of $n$ non-zero idempotents: $c = \sum_{i=1}^n c_i$. Assume $c_n$, say, is not primitive. Then $c_n = d_n + d_{n+1}$, for some non-zero orthogonal idempotents $d_n,d_{n+1}$. For $1 \leq i < n$ and $j = n,n + 1$, we conclude $c_id_j = c_ic_nd_j = 0$, and $c = \sum_{i=1}^{n-1} c_i + d_n + d_{n+1}$ decomposes into the orthogonal sum of $n + 1$ non-zero idempotents, a contradiction. Thus $c_n$, hence (by symmetry) each $c_i$, $1 \leq i \leq n$, is primitive. 

\smallskip

(c) Let $c,d \in R$ be distinct primitive idempotents. Then $cd \in Rc \cap Rd$ is an idempotent. Assuming $cd \neq 0$, condition (ii) in (a) would imply $c = cd = d$, a contradiction. Hence $cd = 0$. Summing up, the primitive idempotents of $R$ are mutually orthogonal and thus finite in number since $R$ is finite-dimensional over $F$. We may therefore form the orthogonal system $(c_1,\dots,c_n)$ of \emph{all} primitive idempotents in $R$. Put $c := \sum_{i=1}^n c_i$ and $R_0 := \{x \in R \mid cx = 0\}$. Then Exc.~\ref{pr.CENIDSPLIDE} implies
\begin{align}
\label{ARCER} R = Rc \oplus R_0 = \bigoplus_{i=1}^n Rc_i \oplus R_0
\end{align}
as direct sums of ideals. Assuming $R_0 \neq \{0\}$, Exc.~\ref{pr.ID} would yield a non-zero idempotent in $R_0$, hence also a primitive one by (b). Using \eqref{ARCER}, one checks that a primitive idempotent in $R_0$ stays primitive in $R$, and we have arrived at a contradiction. Thus $R_0 = \{0\}$, i.e., $R$ is unital with identity element $c$. From condition (iii) in (a) we conclude that $K_i := Rc_i$ ($1 \leq i \leq n$) are finite algebraic field extensions of $F$, and \eqref{ARCER} yields $R \cong \prod_{i=1}^n K_i$.
\end{sol}

\begin{sol}{pr.INFORS} \label{sol.INFORS} Let $(c_i)_{i\in I}$ be an orthogonal system of non-zero idempotents in $R$, indexed by an infinite set $I$. Pick a countably infinite ascending chain
\[
I_0 \subset I_1 \subset \cdots \subset I_n \subset I_{n+1} \subset \cdots 
\]
of finite subsets of $I$ and put $\mfa_n := \sum_{i\in I_n} Rc_i$ for all $n \in \IN$. Then
\[
\mfa_0 \subset \mfa_1 \subset \cdots \subset \mfa_n \subset \mfa_{n+1} \subset \cdots
\]
is a countably infinite ascending chain of ideals in $R$. If $R$ were finitely generated as a $k$-algebra, it would be noetherian (since $k$ is) and the aforementioned chain would become stationary eventually, a contradiction.
\end{sol}

\begin{sol}{pr.NUCMAT} \label{sol.NUCMAT}
For $1 \leq i,j,l,m,p,q \leq n$, a
straightforward verification shows
\begin{align}
\label{ASMA} [ae_{ij},be_{lm},ce_{pq}]
=\,\,&\delta_{jl}\delta_{mp}[a,b,c]e_{iq}, \\
\label{COMMA} [ae_{ij},be_{lm}] =\,\,&\delta_{jl}(ab)e_{im} -
\delta_{im}(ba)e_{lj}.
\end{align}
The inclusion $\Mat_n(\Nuc(A)) \subseteq \Nuc(\Mat_n(A))$ follows
immediately from \eqref{ASMA}. Conversely, let $x = (a_{ij}) = \sum
a_{ij}e_{ij} \in \Nuc(\Mat_n(A))$. Then \eqref{ASMA} for $m = p$
implies $0 = [x,be_{lm},c_{mp}] = \sum_i [a_{il},b,c]e_{iq}$, hence
$[a_{il},A,A] = \{0\}$. Similarly, $0 = [be_{lm},ce_{mq},x] = \sum_j
[b,c,a_{qj}]e_{lj}$, forcing $[A,A,a_{qj}] = \{0\}$. And, finally,
$0 = [be_{lm},x,ce_{pq}] = [b,a_{mp},c]e_{lq}$ implies $[b,a_{mp},c]
= 0$. Summing up, we conclude $a_{ij} \in \Nuc(A)$ for $1 \leq i,j
\leq n$, solving the first part of the problem.

Concerning the second one, let $x = (a_{ij}) = \sum a_{ij}e_{ij} \in
\Cent(\Mat_n(A))$. By the first part, $a_{ij} \in \Nuc(A)$ for $1
\leq i,j \leq n$. Moreover, \eqref{COMMA} for $l = m$ yields $0 =
[x,e_{ll}] = \sum_i a_{il}e_{il} - \sum_j a_{lj}e_{lj}$, hence
$a_{ij} = 0$ for $1 \leq i,j \leq n$, $i \neq j$. Now let $l \neq m$
in \eqref{COMMA}. Then $[x,be_{lm}] = (a_{ll}b)e_{lm} -
(ba_{mm})e_{lm}$, forcing $a_{ll} = a_{mm} \in \Cent(A)$, and we
conclude $x \in \Cent(A)\Eins_n$. Conversely, it is clear that every
element of $\Cent(A)\Eins_n$ belongs to $\Cent(\Mat_n(A))$.
\end{sol}

\begin{sol}{pr.FIDIA} \label{sol.FIDIA}
Let $D$ be a finite-dimensional non-associative division algebra over $F$ and suppose $D$ has dimension $n > 1$. Picking linearly independent vectors $x,y \in D$, the polynomial $\det(L_{x+\bft y}) \in F[\bft]$ has degree $n$, hence a root $\alpha \in F$. But then left multiplication by the non-zero element $x + \alpha  y \in D$ is not bijective, contradicting the property of $D$ being a division algebra. Thus $n = 1$, and we conclude $D = Fe$ for some non-zero element $e \in D$. Therefore $e^2 = \alpha e$ for some $\alpha \in F^\times$, and it is readily checked that $F \to D$, $\xi \mapsto \alpha^{-1}\xi e$ is an isomorphism of $F$-algebras.  
\end{sol}

\solnsec{Section~\ref{s.SCEX}}

\begin{sol}{pr.CEBA} \label{sol.CEBA} Let $m,n \in \IZ$ such that $u := mx + ny
\in \Nuc(A)$. Then $0 = [u,x,x] = m(x^2x - xx^2) + n((yx)x - yx^2) =
2m(x - x) - 2ny = -2ny$, hence $n = 0$, and $0 = [u,x,y] = m(x^2y -
x(xy)) = 2my$, hence $m = 0$. Thus $\Cent(A) = \Nuc(A) = \IZ1_A$; in
particular, $A$ is central. On the other hand, passing to the base
change $A_R$, $R = \IZ/2\IZ$, which by \ref{ss.REDID} agrees with $A/2A$, and setting $\bar u := u_R$ for $u \in
A$, we conclude that $1_{\bar A},\bar x,\bar y$ is an $R$-basis of
$A_R$ with multiplication table $\bar x^2 = \bar x\bar y = \bar
y\bar x = \bar y^2 = 0$. It follows from this at once that $A_R$ is
commutative associative, so while $\Cent(A_R) = \Nuc(A_R) = A_R$ is
a free $R$-module of rank $3$, $(\Cent(A))_R = (\Nuc(A))_R =
R1_{\bar A}$ is a free $R$-module of rank $1$.
\end{sol}

\begin{sol}{pr.DIVKT}
\ref{DIVKT.1}: Suppose $0 \ne x \in A$.  If $y \in A$ is such that $0 = xy$, then $y = 0$ by hypothesis on $A$.  That is, the linear transformation $L_x\: A \to A$ is injective.  Since $A$ is finite-dimensional, $L_x$ is bijective.  Similarly, $R_x$ is bijective.  This verifies the claim.

\ref{DIVKT.2}: Since $F(\bft) \subset F((\bft))$, $A_{F(\bft)} \subset A_{F((\bft))}$, so it suffices to treat the case $K = F((\bft))$.

Since $A$ is a division algebra, it has no zero divisors.  Applying Exc.~\ref{pr.NZKT} with $k = F$ and $R = F[[\bft]]$, we find that $A_R$ has no zero divisors.

Every element of $A_K$ can be written as $a\bft^{-m}$ for some $a \in A_R$ and $m \ge 0$.  Suppose $b\bft^{-n} \in A_K$ for $b \in A_R$ and $n \ge 0$ is such that $0 = (a\bft^{-m})(b\bft^{-n})$.  Then $0 = (ab)\bft^{-m-n}$, i.e., $ab = 0$ in $A_R$, so $a = 0$ or $b = 0$, ergo $a\bft^{-m} = 0$ or $b\bft^{-n} = 0$.
\end{sol}

\begin{sol}{pr.COVERS} \label{sol.COVERS} Put $\mfp$ for the prime ideal that is the kernel of the unit morphism $\vart\:k \to K$.
By hypothesis, there is a prime ideal $\mfp^\prime \in \Spec(k^\prime)$ such that $\mfp = \vart^{-1}(\mfp^\prime)$.  This provides a commutative diagram as in (\ref{ss.LOPROM}.\ref{KAYPE}) and in particular the field $k(\mfp)$ is contained in both $K$ and $\pk(\mfp')$.  It suffices to take $K'$ to be any field containing both $K$ and $\pk(\mfp')$ (called a ``compositum'' of the two fields), which exists by an argument as in \cite[p.~552]{MR1009787}.
\end{sol}

\begin{sol}{pr.FIB} \label{sol.FIB} We proceed in several steps. 

\step{1}
Let $\can = \can_{k,k(\mfp)}\:k \to k(\mfp)$, $\alpha \mapsto \alpha(\mfp)$, be the natural map. Then the diagram
\[
\xymatrix{
k \ar[r]_{\can} \ar[d]_{\vph} & k(\mfp) \ar[d]^{\vph(\mfp)} \\
k^\prime \ar[r]_{\can(\mfp)} & k^\prime(\mfp)}
\]
commutes since
\begin{align*}
\vph(\mfp) \circ \can(\alpha) =\,\,&(\vph \otimes \Eins_{k(\mfp)})\big(1 \otimes \alpha(\mfp)\big) = 1_{k^\prime} \otimes \alpha(\mfp) = 1_{k^\prime} \otimes (\alpha\cdot 1_{k(\mfp)}) \\
=\,\,&(\alpha\cdot 1_{k^\prime}) \otimes 1_{k(\mfp)} = \vph(\alpha) \otimes 1_{k(\mfp)} = \vph(\alpha)(\mfp) = \can(\mfp) \circ \vph(\alpha)
\end{align*}
for all $\alpha \in k$. Applying the contra-variant functor $\Spec$ to this diagram, we conclude that also
\[
\xymatrix{
\Spec\big(k^\prime(\mfp)\big) \ar[rr]_(0.55){\Spec(\can(\mfp))} \ar[d]_{\Spec(\vph(\mfp))} && \Spec(k^\prime) \ar[d]^{\Spec(\vph)} \\
\Spec\big(k(\mfp)\big) \ar[rr]_(0.55){\Spec(\can)} && \Spec(k)}
\]
commutes. But $\Spec(k(\mfp))$ consists of a single point sent by $\Spec(\can)$ to $\can^{-1}(\{0\}) = \Ker(\can) = \mfp$. Hence the image of $\Spec(\can(\mfp))$ belongs to the fiber $\Spec(\vph)^{-1}(\mfp)$. Thus \emph{$\Spec(\can(\mfp))$ induces canonically a continuous map}
\[
\Psi\:\Spec\big(k^\prime(\mfp)\big) \longrightarrow \Spec(\vph)^{-1}(\mfp).
\]

\step{2}
Let $\mfp^\prime \in \Spec(\vph)^{-1}(\mfp)$. Then $\vph^{-1}(\mfp^\prime) = \mfp$, and consulting (\ref{ss.LOPROM}.\ref{KAYPE}), we find a unique homomorphism 
\[
\sigma_{\mfp^\prime}\:k^\prime(\mfp) \longrightarrow k^\prime(\mfp^\prime)
\]
of $k(\mfp)$-algebras making a commutative diagram
\begin{align}
\vcenter{\label{SIGMAPE} \xymatrix{
k \ar[r]^{\vph} \ar[d]_{\can_\mfp} & k^\prime \ar[r]^{\can(\mfp)} \ar[d]_{\can_{\mfp^\prime}} & k^\prime(\mfp) \ar[ldd]^{\sigma_{\mfp^\prime}} \\
k_\mfp \ar[r]^{\pvph} \ar[d]_{\vrh(\mfp)} & k^\prime_{\mfp^\prime} \ar[d]_{\vrh(\mfp^\prime)} \\
k(\mfp) \ar[r]_{\bar\vph} & k^\prime(\mfp^\prime).}}
\end{align}
Thus $\Ker(\sigma_{\mfp^\prime}) \in \Spec(k^\prime(\mfp))$, and we obtain a map
\[
\Phi\:\Spec(\vph)^{-1}(\mfp) \longrightarrow \Spec\big(k^\prime(\mfp)\big)
\]
given by
\begin{align*}
\Phi(\mfp^\prime) := \Ker(\sigma_{\mfp^\prime}) &&(\mfp^\prime \in \Spec(\vph)^{-1}(\mfp)).
\end{align*} 
The definitions and \eqref{SIGMAPE} imply
\[
\Psi \circ \Phi(\mfp^\prime) = \can(\mfp)^{-1}\big(\sigma_{\mfp^\prime}^{-1}(\{0\})\big) = \big(\sigma_{\mfp^\prime} \circ \can(\mfp)^{-1}\big)(\{0\}) = \{\alpha^\prime \in k^\prime \mid \alpha^\prime(\mfp^\prime) = 0\} = \mfp^\prime
\]
for all $\pmfp \in \Spec(\vph)^{-1}(\mfp)$, hence $\Psi \circ \Phi = \Eins_{\Spec(\vph)^{-1}(\mfp)}$. 

\step{3}
Let $\mfq \in \Spec(k^\prime(\mfp))$. Then
\[
\Phi \circ \Psi(\mfq) = \Phi(\mfp^\prime) = \Ker(\sigma_{\mfp^\prime}), \quad \mfp^\prime := \can(\mfp)^{-1}(\mfq) \in \Spec(\vph)^{-1}(\mfp).
\]
By (\ref{ss.LOPROM}.\ref{KAYPE}) again, we obtain a diagram
\[
\xymatrix{
k \ar[r]^{\vph} \ar[d]_{\can_\mfp} & k^\prime \ar[r]^{\can(\mfp)} \ar[d]_{\can_{\mfp^\prime}} & k^\prime(\mfp) \ar@/^7pc/[ldd]^{\sigma_{\mfp^\prime}} \ar[d]_{\can_\mfq} \\
k_\mfp \ar[r]^{\pvph} \ar[d]_{\vrh(\mfp)} & k^\prime_{\mfp^\prime} \ar[r]^{\can(\mfp)^\prime} \ar[d]_{\vrh(\mfp^\prime)} & \big(k^\prime(\mfp)\big)_\mfq \ar[d]_{\vrh(\mfq)} \\
k(\mfp) \ar[r]^{\bar\vph} & k^\prime(\mfp^\prime) \ar[r]^{\overline{\can(\mfp)}} & k(\mfq)}
\]
where the four squares commute.  But so does the diagram
\[
\xymatrix{
k^\prime \ar[r]^{\can(\mfp)} \ar[d]_{\vrh(\pmfp) \circ \can_{\mfp^\prime}} & k^\prime(\mfp) \ar[d]^{\vrh(\mfq) \circ \can_\mfq} \ar[ld]^{\sigma_{\mfp^\prime}} \\
k^\prime(\mfp^\prime) \ar[r]_{\overline{\can(\mfp)}} & k^\prime(\mfq),}
\]
which by \eqref{SIGMAPE} is clear for the upper triangle but holds true also for the lower one because $\overline{\can(\mfp)} \circ \sigma_{\mfp^\prime}$ and $\vrh(\mfq) \circ \can_\mfq$ are both $k(\mfp)$-linear satisfying $\overline{\can(\mfp)} \circ \sigma_{\mfp^\prime} \circ       \can(\mfp) = \vrh(\mfq) \circ \can_\mfq \circ \can(\mfp)$, hence $\overline{\can(\mfp)} \circ \sigma_{\mfp^\prime} = \vrh(\mfq) \circ \can_\mfq$. Since $\overline{\can(\mfp)}$ as a field homomorphism is injective, we therefore conclude
\[
\Phi \circ \Psi (\mfq) = \Ker(\sigma_{\mfp^\prime}) = \Ker(\overline{\can(\mfp)} \circ \sigma_{\mfp^\prime}) = \Ker\big(\vrh(\mfq) \circ \can_\mfq\big) = \mfq.
\]
Summing up we have shown that \emph{$\Psi$ is bijective with inverse $\Phi$}. 

\step{4}
By $1^\circ$ and $3^\circ$, it remains to show that $\Phi$ is continuous, equivalently, that $\Psi$ is open. In order to see this, we first deduce from (\ref{ss.ITSCA}.\ref{ITBA}) and \ref{ss.REDID} the natural identifications
\[
k^\prime(\mfp) = (k^\prime \otimes k_\mfp) \otimes_{k_\mfp} (k_\mfp/\mfp_\mfp) = k^\prime_\mfp \otimes_{k_\mfp}(k_\mfp/\mfp_\mfp) =k^\prime_\mfp/\mfp_\mfp k^\prime_\mfp
\] 
such that, with the natural maps $i\:k^\prime \to k^\prime_\mfp$, $\pi\:k^\prime_\mfp \to k^\prime_\mfp/\mfp_\mfp k^\prime_\mfp$, the triangle 
\[
\xymatrix{
k^\prime \ar[d]_{i} \ar[rd]^{\psi} \\
k^\prime_\mfp \ar[r]_(0.3){\pi} & k^\prime_\mfp/\mfp_\mfp k^\prime_\mfp = k^\prime \otimes k(\mfp)}
\]
commutes. The proof will be complete once we have shown
\begin{align}
\label{PSIDEPI} \Psi\Big(D\big(\pi (f/s)\big)\Big) = D(f) \cap \Spec(\vph)^{-1}(\mfp)
\end{align}
for all $f \in k^ \prime$, $s \in k \setminus \mfp$. Let $\mfq \in \Spec(k^\prime(\mfp)) = \Spec(k^ \prime_\mfp/\mfp_\mfp k^ \prime_\mfp)$. Then the chain of equivalent conditions
\begin{align}
\label{QUDEPI} \mfq \in D\big(\pi(f/s)\big)\;\;&\Longleftrightarrow \pi(f/s) \notin \mfq \;\;\Longleftrightarrow f/s \notin \pi^{-1}(\mfq) \\
\;\;&\Longleftrightarrow f \notin i^{-1}\big(\pi^{-1}(\mfq)\big) = \psi^{-1}(\mfq)\;\;\Longleftrightarrow \can(\mfp)^{-1}(\mfq) \in D(f) \notag
\end{align}
shows that the left-hand side of \eqref{PSIDEPI} is contained in the right. Conversely, let $\mfp^\prime \in D(f) \cap \Spec(\vph)^{-1}(\mfp)$. Then $\mfq := \Psi^{-1}(\mfp^\prime) \in \Spec(k^\prime(\mfp))$ and $\can(\mfp)^{-1}(\mfq) = \mfp^\prime \in D(f)$. Consulting \eqref{QUDEPI}, we conclude $\mfq \in D(\pi(f/s))$, hence $\mfp^\prime = \Psi(\mfq) \in \Psi(D(\pi(f/s)))$, and the proof of \eqref{PSIDEPI} is complete.
\end{sol}

\begin{sol}{pr.LOCHOM} \label{sol.LOCHOM} (a) We may assume $\chi = 0$.  Pick $\mfp \in U$.  Let $u_1,\dots,u_m$ be a finite set of generators for the $k$-module $M$. For each $1 \leq i \leq m$, we have $\psi(u_i)/1 = \psi(u_i)_\mfp = \psi_\mfp(u_{i\mfp}) = 0$ in $N_\mfp$. Hence there are $f_i \in k \setminus \mfp$ such that $f_i\psi(u_i) = 0$ in $N$. But this means $\psi(u_i)/1 = 0$ in $N_{f_i}$, which in turn, by \ref{ss.PROPSP}, amounts to $\psi_\mfq = 0$ for all $\mfq \in \cap_{i=1}^m D(f_i) = D(f)$ for $f = f_1 f_2 \cdots f_m$. Hence we have shown $\mfp \in D(f) \subseteq U$, and $U$ is Zariski-open in $X$.

(b) Again let $u_1,\dots,u_m$ be a finite set of generators for $A$ as a $k$-module. Given $i = 1,\dots,m$, we define $k$-linear maps $\psi_i,\chi_i\:A \to B$ by  $\psi_i(v) := \vph(u_iv)$, $\chi_i(v) := \vph(u_i)\vph(v)$, $v \in A$. By (a),
\[
U_i := \{\mfp \in X \mid \psi_{i\mfp} = \chi_{i\mfp}\}
\]
is Zariski-open in $X$. Hence so is $V = \bigcap_{i=1}^mU_i$. 
\end{sol}

\begin{sol}{pr.IDEPAR} \label{sol.IDEPAR}
In the solution to this problem, free use
will be made of the formalism for the Zariski topology as explained,
for example, in \cite[II, \S4.3]{MR0360549}. 

(a) For $\mfp \in X$, the chain of equivalent conditions
\begin{align*}
\mfp \in D(\varepsilon) &\Longleftrightarrow \varepsilon \notin \mfp
\Longleftrightarrow \varepsilon_\mfp \notin \mfp_\mfp
\Longleftrightarrow \varepsilon_\mfp \in k_\mfp^\times \\
&\Longleftrightarrow \varepsilon_\mfp = 1_\mfp \Longleftrightarrow (1
- \varepsilon)_\mfp = 0
\end{align*}
gives the first displayed equality, while the second one follows
from the first and the fact that local rings are \emph{connected},
i.e., contain only the trivial idempotents $0$ and $1$. The third
equality is now obvious. To establish the final assertion of (a), we
note for $\mfp \in X$ that $(e_{i\mfp})_{i\in I}$ is an
orthogonal system of idempotents in $k_\mfp$, whence $(\sum
e_i)_\mfp \neq 0$ if and only if $e_{i\mfp} \neq 0$ for some $i =
1,\dots,r$. This proves $\bigcup D(\varepsilon_i) = D(\sum
\varepsilon_i)$, and the union on the left is disjoint since
$D(\varepsilon_i) \cap D(\varepsilon_j) =
D(\varepsilon_i\varepsilon_j) = D(0) = \emptyset$ for $i,j\in I$, $i \neq j$. 

\smallskip

(b) Given a complete system $(\varepsilon_i)_{i\in I}$ of 
orthogonal idempotents in $k$, the subsets $U_i := D(\varepsilon_i)
\subseteq X$, $i\in I$, are open in $X$ and by (a) satisfy
$\bigcup U_i = D(\sum e_i) = D(1) = X$, where the union on the very
left is disjoint and the $U_i$ are empty for almost all $i \in I$. Thus the assignment $(\varepsilon_i) \mapsto (U_i)$
gives a map of the right kind, which is injective since the $U_i$ by
(a) determine the $\varepsilon_i$ locally, hence globally. It
remains to show that the map in question is surjective, so let
$(U_i)_{i\in I}$ be a decomposition of $X$ into disjoint
open subsets almost all of which are empty.

First, view $X$ as a geometric space with structure sheaf $\msO$.
Since the $U_i$ form a disjoint open cover of $X$, we may apply the
sheaf axioms to find, for each $i \in I$, a unique global
section $\varepsilon_i \in H^0(X,\msO) = k$ which restricts to the
identity element of $H^0(U_i,\msO)$ and to zero on $H^0(U_j,\msO)$ for all $j \in I$, $j \neq i$. It is then readily checked that
the $\varepsilon_i$ form a complete orthogonal system of idempotents
in $k$, forcing the $D(\varepsilon_i)$ to be an open disjoint cover
of $X$. On the other hand, given $\mfp \in U_i$, we obtain
$\varepsilon_{i\mfp} = 1_\mfp$, which implies $U_i \subseteq
D(\varepsilon_i)$ by (a). Summing up, we have equality and
surjectivity is proved.

Adopting a more direct approach based on the proof of \cite[II.4, Prop. 15]{MR0360549}, let $i \in I$. Then $U_i$ and
$Y_i := \bigcup_{j\neq i}U_j$ are closed subsets of $X$, hence may
be written in the form $U_i = V(\mfa_i)$, $Y_i = V(\mfb_i)$ for some
ideals $\mfa_i,\mfb_i \subseteq k$. Since $V(\mfa_i + \mfb_i) =
V(\mfa_i) \cap V(\mfb_i) = U_i \cap Y_i = \emptyset$, we conclude
$\mfa_i+ \mfb_i = k$, hence $\alpha_i + \beta_i = 1$ for some
$\alpha_i \in \mfa_i$, $\beta_i \in \mfb_i$. This implies $\alpha_i
= \alpha_i(\alpha_i + \beta_i) = \alpha_i^2 + \alpha_i\beta_i \equiv
\alpha_i^2 \bmod \mfa_i\mfb_i$, so $\alpha_i \in \mfa_i$ becomes an
idempotent modulo $\mfa_i\mfb_i$. On the other hand, $X = U_i \cup
Y_i = V(\mfa_i) \cup V(\mfb_i) = V(\mfa_i\mfb_i)$, so $\mfa_i\mfb_i
\subseteq k$ is a nil ideal. Now Exc.~\ref{pr.ID}~(b) yields an
idempotent $\varepsilon_i^\prime \in \mfa_i$ such that
$\varepsilon_i^\prime \equiv \alpha_i \bmod \mfa_i\mfb_i$. Thus
$\varepsilon_i := 1 - \varepsilon_i^\prime \in k$ is an
idempotent as well satisfying $\varepsilon_i \equiv \beta_i \bmod
\mfa_i\mfb_i$, hence $\pvep_i \in \mfb_i$. We conclude $U_i = V(\mfa_i) \subseteq
V(\varepsilon_i^\prime) = D(\varepsilon_i)$, $Y_i = V(\mfb_i)
\subseteq V(\varepsilon_i) = D(\varepsilon_i^\prime)$ and
$D(\varepsilon_i) \cap D(\varepsilon_i^\prime) = \emptyset$ by (a).
Therefore $U_i = D(\varepsilon_i)$, and it remains to show that
$(\varepsilon_i)_{i\in I}$ is a complete orthogonal system of
idempotents. For $i,j \in I$, $i \neq j$, we deduce
$D(\varepsilon_i\varepsilon_j) = D(\varepsilon_i) \cap
D(\varepsilon_j) = U_i \cap U_j = \emptyset$, so the idempotent
$\varepsilon_i\varepsilon_j$ is also nilpotent, hence zero. Thus the
$\varepsilon_i$ are orthogonal. For the same reason, $\vep_i = 0$ for almost all $i \in I$ Now (a) gives
\[
X = \bigcup U_i = \bigcup D(\varepsilon_i) = \bigcup V(1 -
\varepsilon_i) = V\big(\prod (1 - \varepsilon_i)\big) = V(1 -
\sum\varepsilon_i),
\]
forcing the idempotent $1 - \sum\varepsilon_i$ to be nilpotent,
hence zero.
\end{sol}

\begin{sol}{pr.SURPRO} \label{sol.SURPRO}
 By \cite[\S{II.3.3}, Thm.~1]{MR0360549}, it suffices to prove that the map $M_\mfm \to N_\mfm$ induced by $f$ is an isomorphism for every $\mfm$, so we may replace $k$ by $k_\mfm$ and assume $k$ is local.  Then $M$ and $N$ are free modules, and $M$ has finite rank.  Since they have the same rank, they are isomorphic.  Let $g\:N \to M$ be an isomorphism. Then $g \circ f\:M \to M$ is a $k$-linear surjection, and Prop.~\ref{p.SURBI} shows that it is, in fact, a bijection. Hence so is $f$.  
\end{sol}

\begin{sol}{pr.RANDEC} \label{sol.RANDEC} 
Setting $X := \Spec(k)$, the subsets
\begin{align}
\label{RANKI} U_i := \{\mfp \in \Spec(k)\mid \rk_\mfp(M) = i\} \subseteq X &&(i \in \IN)
\end{align}
form an open disjoint cover of $X$ and since $X$ is quasi-compact \cite[II.4, Prop.~12]{MR0360549}, $U_i$ is empty for almost all $i \in \IN$. This proves (a) while, in order to prove (b), we apply Exc.~\ref{pr.IDEPAR} to obtain a complete orthogonal system $(\vep_i)_{i\in\IN}$ of idempotents in $k$ such that $U_i = D(\vep_i)$ for all $i \in \IN$. Given $j \in \IN$, let $\pi_j\:k \to k_j$ be the canonical projection induced by \eqref{KKEI} and let $\mfp_j \in \Spec(k_j)$. Then
\[
\mfp := \Spec(\pi_j)(\mfp_j) = \pi_j^{-1}(\mfp_j) = \mfp_j \times \prod_{i\neq j} k_i \in X,
\]  
and (\ref{ss.LOPROM}.\ref{MPRLO}) shows
\[
M_{j\mfp_j} = M_\mfp \otimes_{k_\mfp} k_{j\mfp_j}.
\]
Since $M_\mfp$ is free of rank $j$ over $k_\mfp$, it follows that $M_{j\mfp_j}$ is free of rank $j$ over $k_{j\mfp_j}$. This proves existence. In order to establish uniqueness, let $(\eta_i)_{i\in\IN}$ be any complete orthogonal system of idempotents in $k$, giving rise, in analogy to \eqref{KKEI}, \eqref{MMEI}, to the decompositions
\begin{align}
\label{KLEI} k =\,\,&\prod_{i\in\IN} l_i, \quad l_i = k\eta_i &&(i \in \IN), \\
\label{MNEI} M =\,\,&\prod_{i\in\IN} N_i, \quad N_i = N \otimes l_i &&(i \in \IN),
\end{align}
and making each $N_i$ a projective $l_i$-module of rank $i$. Again by Exc.~\ref{pr.IDEPAR}, the $V_i := D(\eta_i)$ are empty for almost all $i \in I$ and form an open disjoint cover of $X$, so it suffices to show $V_i \subseteq U_i$ for all $i \in \IN$. Fixing $j \in \IN$ and $\mfq \in V_j$, we apply Exc.~\ref{pr.IDDIRSUM}, and find a decomposition $\mfq = \prod_{i\in\IN} \mfq_i$, where $\mfq_i = \mfq \cap l_i$ is an ideal in $l_i$. But $\eta_j \notin \mfq$ while $\eta_i\eta_j = 0 \in \mfq$ for all $i \neq j$, forcing $\eta_i \in \mfq \cap l_i = \mfq_i$ and then $\mfq_i = l_i$. Summing up we have shown
\[
\mfq = \mfq_j \times \prod_{i\neq j} l_i,
\]
so $\mfq_j$ is a prime ideal in $l_j$ and $\Spec(\rho_j)(\mfq_j) = \mfq$, where $\rho_j\:k \to l_j$ is the canonical projection induced by \eqref{KLEI}. Since
\[
N_{j\mfq_j} = M_\mfq \otimes_{k_\mfq} l_{j\mfq_j}
\]
has rank $j$ over $l_{j\mfq_j}$, the free $k_\mfq$-module $M_\mfq$ must have rank $j$ over $k_\mfq$. This proves $\mfq \in U_j$ by \eqref{RANKI}, and the solution is complete.
\end{sol}

\begin{sol}{pr.LOCSIALG} \label{sol.LOCSIALG}
If $\varphi\:A \to A^\prime$ is a unital
algebra homomorphism, then so is $\varphi(\mfp)\:A(\mfp) \to
A^\prime(\mfp)$, for each $\mfp \in \Spec(k)$. But since $A(\mfp)$
is simple by hypothesis, $\varphi(\mfp)$ is injective and a
dimension argument shows that it is, in fact, an isomorphism. Hence,
by Nakayama's lemma \cite[X, Prop.~4.5]{MR1878556} or \stacks{07RC}, $\varphi_\mfp\:A_\mfp
\to A^\prime_\mfp$ is an isomorphism, whence $\varphi$ is one
\cite[II.3, Thm.~1]{MR0360549}.
\end{sol}

\begin{sol}{pr.GRANTO} \label{sol.GRANTO}
(a) The right-hand side is clearly
contained in the left, so it remains to establish the reverse
inclusion. To this end, put $(A^\prime,G^\prime) = (A,G) \otimes
F^\prime$ and write $a^\prime \in \Cent(A^\prime,G^\prime)$ in the
form $a^\prime = \sum a_i \otimes \alpha_i^\prime$ with $a_i \in A$
and $\alpha_i^\prime \in F^\prime$ linearly independent over $F$.
For $x,y \in A$, we have $\sum [a_i,x,y] \otimes \alpha_i^\prime =
[a^\prime,x \otimes 1_{F^\prime},y \otimes 1_{F^\prime}] = 0$ and,
similarly, $\sum [x,a_i,y] \otimes \alpha_i^\prime = \sum [x,y,a_i]
\otimes \alpha_i^\prime = \sum [x,a_i] \otimes \alpha_i^\prime = 0$.
Moreover, for $g \in G$, $\sum g(a_i) \otimes \alpha_i^\prime = (g
\otimes \Eins_{F^\prime})(a^\prime) = a^\prime = \sum a_i \otimes
\alpha_i^\prime$. Since the $\alpha_i^\prime$ are linearly
independent over $F$, we therefore obtain $[a_i,x,y] = [x,a_i,y] =
[x,y,a_i] = [a_i,x] = 0$, $g(a_i) = a_i$, hence $a_i \in \Cent(A,G)$
for all $i$, which implies $a^\prime \in \Cent(A^\prime,G^\prime)$,
as claimed. 

(b) By Prop.~\ref{ss.p-CENTMAL}, the right multiplication of $A$
induces an isomorphism $R:\Cent(A) \to \End_{\Mult(A)}(A)$. Hence it
suffices to show that the elements of $\Cent(A,G)$ correspond to the
$\Mult(A,G)$-linear maps $A \to A$ under this isomorphism, which
follows from the relations $R_{g(a)}g = gR_a$ for all $a \in
\Cent(A)$, $g \in G$. 

(c) If $(A,G)$ is a simple $F$-algebra with a group of
antomorphisms, then $M(A,G)$ acts irreducibly on $A$, so Schur's
lemma and (b) yield the assertion.

(d) Since $A$, by simplicity of $(A,G)$, is an irreducible faithful
$M(A,G)$-module, $M(A,G)$ is a primitive left-artinian $F$-algebra, hence simple \cite[Thm.~4.2]{MR1009787}. The first assertion now follows from the double centralizer theorem \cite[Thm.~4.10]{MR1009787}. It immediately  implies necessity in the second
part. Conversely, suppose $A \neq\{0\}$ and $M(A,G) = \End_F(A)$.
Then $A$ is an irreducible $M(A,G)$-module, forcing $(A,G)$ to be
simple and $\Cent(A,G)$ to be a field. Also, by (b), $\End_F(A) =
M(A,G) = \End_{\Cent(A,G)}(A)$ and a dimension count yields
sufficiency as well.
\end{sol}

\begin{sol}{pr.GRANTOCS} \label{sol.GRANTOCS}
(i) $\Rightarrow$ (ii). Let $F^\prime$ be
any extension field of $F$ and put $(A^\prime,G^\prime) = (A,G)
\otimes F^\prime$ as an $F^\prime$-algebra with a group of
antomorphisms. A moment's reflection shows $M(A^\prime,G^\prime) =
M(A,G) \otimes F^\prime$. Now the assertion follows from the second
part of Exc.~\ref{pr.GRANTO}~(d).

(ii) $\Rightarrow$ (iii). Obvious.

(iv) $\Rightarrow$ (i). This follows immediately from
Exc.~\ref{pr.GRANTO}~(a), (c).
\end{sol}

\begin{sol}{pr.TPGRANTO} \label{sol.TPGRANTO}
We follow standard arguments, reproduced in
\cite[I Satz 5.5]{MR34:4310}, for example. Let $J \subseteq A \otimes A^\prime$ be a non-zero
ideal of $(A,G) \otimes (A^\prime,G^\prime)$ and consider the least
positive integer $n$ such that there exist elements $x_1,\dots,x_n
\in A, x_1^\prime,\dots,x_n^\prime \in A^\prime$ satisfying
\begin{align}
\label{MIRE} 0 \neq x := \sum_{i=1}^{n}x_i \otimes x_i^\prime \in J.
\end{align}
Then the elements $x_1^\prime,\dots,x_n^\prime$ are linearly
independent over $F$. Furthermore, a straightforward verification
shows that the totality of elements $u_1 \in A$ such that
\[
\sum_{i=1}^{n} u_i \otimes x_i^\prime \in J
\]
for some $u_2,\dots,u_n \in A$ forms an ideal $I$ of $(A,G)$ which
is non-zero since $x_1 \in I$ and $x_1 \neq 0$ by  minimality of
$n$. Simplicity of $(A,G)$ therefore implies $I = A$, so in
\eqref{MIRE} we may assume $x_1 = 1$. But then, for all $u,v \in A$,
\[
\sum_{i=2}^{n} [x_i,u,v] \otimes x_i^\prime = [x,u \otimes 1,v
\otimes 1] \in J.
\]
Again minimality of $n$ yields $[x,u\otimes 1,v \otimes 1] = 0$
and then $[x_i,u,v] = 0$ for $2 \leq i \leq n$. Similarly,
$[u,x_i,v] = [u,v,x_i] = [x_i,u] = 0$, forcing $x_i \in \Cent(A)$
for $2 \leq i \leq n$. Applying $g \otimes \Eins$ for any $g \in G$
to \eqref{MIRE} and again observing $x_1 = 1$, we conclude
\[
\sum_{i=2}^{n} \big(x_i - g(x_i)\big) \otimes x_i^\prime = x -(g
\otimes \Eins)(x) \in J\,,
\]
and minimality of $n$ implies $\sum_{i=2}^{n} (x_i - g(x_i))
\otimes x_i^\prime = 0$, hence $x_i \in \Cent(A,G)$ for $2 \leq i
\leq n$. But $(A,G)$ was assumed to be central simple, so
$\Cent(A,G) = k1_A$, and we conclude $x = 1_A \otimes x^\prime \in
J$ for some $x^\prime \in A^\prime$. Now put
\[
I^\prime = \{ u^\prime \in A^\prime \mid 1_A \otimes u^\prime \in
J\}.
\]
Again it is straightforward to check that $I^\prime$ is an ideal in
$(A^\prime, G^\prime)$ which is non-zero since it contains $x^\prime
\neq 0$. Simplicity of $(A^\prime, G^\prime)$ therefore implies
$I^\prime = A^\prime$, hence $1_{A \otimes A^\prime} = 1_A \otimes
1_{A^\prime} \in J$ and $J = A$, as claimed.
\end{sol}


\solnsec{Section~\ref{s.INV}}

\begin{sol}{pr.PERMAT} \label{sol.PERMAT} (a) We must show $P_\pi = \Eins_m \Rightarrow \pi = 1_{S_m}$, $P_{\pi\vrh} = P_\pi P_\vrh$ and $P_{\pi}^{-1} = P_\pi^\trans$ for all $\pi,\vrh \in S_m$. The first statement is obvious, while the second one follows from
\[
P_\pi P_\vrh = \sum_{i,j=1}^m e_{\pi(i),i}e_{\vrh(j),j} = \sum_{i=1}^me_{\pi(i),\vrh^{-1}(i)} = \sum_{i=1}^me_{\pi\vrh(i),i} = P_{\pi\vrh}.
\]
And, finally, $P_\pi^\trans = \sum e_{i,\pi(i)} = \sum e_{\pi^{-1}(i),i} = P_{\pi^{-1}} = P_\pi^{-1}$. Thus (a) is proved. In addition, we have
\begin{align} \label{PERMAU}
P_\pi e_{rs}P_\pi^\trans = e_{\pi(r),\pi(s)} &&(1\le r,s\le m)
\end{align}
since the left-hand side equals $\sum_{i,j} e_{\pi(i),i}e_{rs}e_{j,\pi(j)} = e_{\pi(r),\pi(r)}$, as desired.

(b) In terms of matrix units, the quantities $I,J$ of \ref{ss.SPSPL} may be written as
\[
I = \sum_{i=1}^n(e_{i,i+n} - e_{i+n,i}), \quad J = \sum_{i=1}^n(e_{2i-1,2i} - e_{2i,2i-1}).
\]
Hence \eqref{PERMAU} implies  $P_\pi IP_\pi^\trans = J$, and setting $\tau := \tau_\ort^J$, we conclude
\[
\tau(P_\pi SP_\pi^\trans) = J^{-1}(P_\pi S P_\pi^\trans)^\trans J = J^{-1}P_\pi S^\trans P_\pi^\trans J = P_\pi I^{-1}S^\trans IP_\pi^\trans = P_\pi\tau_\spl(S)P_\pi^\trans
\]
for all $S \in \Mat_{2n}(k)$, and the assertion follows.
\end{sol}


\solnsec{Section~\ref{s.QUAMA}}

\begin{sol}{pr.ASLIFO}  Let $t\:A \to k$ be an associative
linear form and suppose the corresponding bilinear form $\sigma :=
\sigma_t$ (cf. \ref{ss.ASBILIFO}) is regular (cf.
\ref{ss.REBI}). By hypothesis, there exists a unique element $e \in
A$ such that $t(x) = \sigma(e,x) = t(ex)$ for all $x \in A$. Now the
relations $t((ex)y) = t(e(xy)) = t(xy)$ and $t((xe)y) = t(y(xe)) =
t((yx)e) = t(e(yx)) = t(yx) = t(xy)$ for all $x,y \in A$ combined
with regularity of $\sigma$ imply $ex = x = xe$ for all $x \in
A$. Hence $A$ is unital with $1_A = e$.
\end{sol}

\begin{sol}{skip.pr.bilind} Let $\alpha_1,\dots,\alpha_n \in k$ satisfy $\sum_{i=1}^n\alpha_im_i = 0$ and put $x := (\alpha_1,\dots,\alpha_n) \in k^n$. Then 
\[
\sum_{i=1}^n\alpha_i\sigma(m_i,m_j) = \sigma(\sum_{i=1}^n\alpha_im_i,m_j) = 0
\]
for $1\le j\le n$, hence $xS = 0$. Applying the adjoint of $S$ yields $0 = xSS^\sharp = \det(S)x$, and since $\det(S)$ is not a zero divisor in $k$, we conclude $x = 0$, i.e., $\alpha_i = 0$ for $1\le i\le n$, as desired.
\end{sol}

\begin{sol}{pr.SPLHYPLA} \label{sol.SPLHYPLA}
The implications (iii) $\Rightarrow$ (i)
$\Rightarrow$ (ii) being obvious, it suffices to prove (ii)
$\Rightarrow$ (iii), which we do by following a suggestion of O.
Loos. Letting $(e_+,e_-)$ be a hyperbolic pair of $\bfh_L$ and writing $e_\pm = (v_\pm^\ast,v_\pm)$, $v_\pm \in L$, $v_\pm^\ast \in L^\ast$, we apply (\ref{ss.HYSP}.\ref{HYPA}) and
obtain $\la v_+^\ast + v_-^\ast,v_+ + v_- \ra = \bfh_L(e_+ +
\,e_-) = 1$. Therefore $v_+ + \,v_- \in L$ is a unimodular vector,
hence a basis of $L$ over $k$, and (iii) follows.
\end{sol}

\begin{sol}{pr.HYPLA} \label{sol.HYPLA}
(i) $\Rightarrow$ (ii). Clear.

(ii) $\Rightarrow$ (iii). We have $q = \la S\ra_{\mathrm{quad}}$ with $S := \left(\begin{smallmatrix}
1 & 0 \\
0 & -1
\end{smallmatrix}\right)$, hence $Dq = \la S + S^\trans\ra$ with $S + S^\trans = \left(\begin{smallmatrix}
2 & 0 \\
0 & -2
\end{smallmatrix}\right)$ having determinant $-4$. Since $q$ is regular, $2 \in k$ is a unit.

(iii) $\Rightarrow$ (i). With the canonical basis $e_1,e_2$ of $k^2$, we have $q(e_1) = 1$, $q(e_2) = -1$, $q(e_1,e_2) = 0$. One checks that $(u_+,u_-)$ with $u_\pm := \frac{1}{2}(e_1 \pm e_2)$ is a hyperbolic pair relative to $q$, and since $e_1 = u_+ + u_-$, $e_2 = u_+ - u_-$, it follows that $q$ is a hyperbolic plane.
\end{sol}

\begin{sol}{pr.HYPDISC} \label{sol.HYPDISC}
Because the base ring is a field of characteristic $\ne 2$, $q$ is isomorphic to a form $\la \alpha_1, \alpha_2 \raq$ for some $\alpha_1, \alpha_2 \in F^\times$.  By hypothesis on the discriminant, $\alpha_1 \cong -\alpha_2$ and we find $q \cong \la \alpha_1 \ra \otimes \la 1, -1 \raq$.  It suffices now to note that $\la 1, -1 \raq$ is a split hyperbolic plane by Exercise \ref{pr.HYPLA}, so tensoring with $\la \alpha_1 \ra$ does not change its isomorphism class (\ref{e.HYSP}).
\end{sol}

\begin{sol}{pr.HYPPA} \label{sol.HYPPA}
 The implications (ii) $\Rightarrow$ (i) being obvious, let us assume that $(u_1,u_2)$ is a hyperbolic pair in $\bfh$. Writing $u_j = \alpha_{1j}e_1 + \alpha_{2j}e_2$ ($\alpha_{ij} \in k$,\;$i,j = 1,2$), we conclude
\begin{align}
\label{COHY} \alpha_{1j}\alpha_{2j} = 0, \quad \alpha_{11}\alpha_{22} + \alpha_{21}\alpha_{12} = 1 &&(j = 1,2).
\end{align}
Thus $(\vep_+,\vep_-)$, $\vep_+ := \alpha_{11}\alpha_{22}$, $\vep_- := \alpha_{21}\alpha_{12}$, is a complete orthogonal system of idempotents in $k$, giving rise to a decomposition of $k$ as in (ii) such that \eqref{HYPM} holds. Moreover, setting $\gamma_+ := \alpha_{11}\vep_+ \in k_+^\times$ with inverse $\gamma_+^{-1} = \alpha_{22}\vep_+$, $\gamma_-:= \alpha_{21}\vep_- \in k_-^\times$ with inverse $\gamma_-^{-1} = \alpha_{12}\vep_-$, we apply \eqref{COHY} and obtain
\begin{align*}
u_{1+} =\,\,&\vep_+u_1  = \alpha_{11}\alpha_{22}(\alpha_{11}e_1 + \alpha_{21}e_2) = \gamma_+e_{1+}, \\
u_{2+} =\,\,&\vep_+u_2 = \alpha_{11}\alpha_{22}(\alpha_{12}e_1 + \alpha_{22}e_2) = \gamma_+^{-1}e_{2+}, \\
u_{1-} =\,\,&\vep_-u_1 = \alpha_{21}\alpha_{12}(\alpha_{11}e_1 + \alpha_{21}e_2) = \gamma_-e_{2-}, \\
u_{2-} =\,\,&\vep_-u_2 = \alpha_{21}\alpha_{12}(\alpha_{12}e_1 + \alpha_{22}e_2) = \gamma_-^{-1}e_{1-},
\end{align*}
hence \eqref{GAHY}, \eqref{GOPHY}.
\end{sol}

\begin{sol}{pr.QUALI} \label{sol.QUALI}
 Uniqueness is clear, so it suffices to prove existence. Setting
\[
\Quad_k(M,P) := \{Q \mid Q\:M \to P\; \text{is $k$-quadratic}\}
\] 
as a $k$-module, the map
\[
h\:N \longrightarrow \Quad_k(M,P), \quad y \longmapsto h(y) := g(\emptyslot,y)
\]
is $k$-linear and hence induces a $k$-linear map
\[
h_1\:N \longrightarrow \Quad_R(M_R,P_R), \quad y \longmapsto h_1(y) := h(y)_R,
\]
$h(y)_R$ being the $R$-quadratic extension of $h(y)$ is the sense of Corollary~\ref{c.BAQUA}. Applying (\ref{ss.SCEM}.\ref{REX}), we obtain an $R$-linear map
\[
h_1^\prime\:N_R \longrightarrow \Quad_R(M_R,P_R)
\]
such that
\begin{align}
\label{HAPR} h_1^\prime(y_R) = h(y)_R
\end{align}
for all $y \in N$. Now define $g_R\:M_R \times N_R \to P_R$ by
\begin{align}
\label{GEAR} g_R(u,v) := h_1^\prime(v)(u) 
\end{align}
for all $u \in M_R$, $v \in N_R$. Combining \eqref{GEAR} with \eqref{HAPR} and (\ref{c.BAQUA}.\ref{BAQUA}), one checks that \eqref{QUALIR} is commutative.
\end{sol}

\begin{sol}{pr.MULQUA} \label{sol.MULQUA}
We argue by induction on $n$. The case $n = 1$ has been settled in Cor.~\ref{c.BAQUA}. Now suppose $n > 1$ and assume that the assertion holds for $n - 1$. For any $u_n \in M_n$, the map 
\begin{align}
\label{EFUN} F_{u_n}\:M_1 \times \cdots \times M_{n-1} \longrightarrow M, \quad (u_1,\dots,u_{n-1}) \longmapsto F(u_1,\dots,u_{n-1},u_n)
\end{align}
is obviously $k$-$(n-1)$-quadratic. Hence by the induction hypothesis, it has  a unique $R$-$(n-1)$-quadratic extension
\[
(F_{u_n})_R\:M_{1R} \times \cdots \times M_{n-1,R} \longrightarrow M_R.
\]
Now let $x_i \in M_{iR}$ be arbitrary for $1 \leq i < n$ and put $x := (x_1,\dots,x_{n-1})$. We claim that \emph{the map
\begin{align}
\label{EFEX} F_x\:M_n \longrightarrow M_R, \quad u_n \longmapsto F_x(u_n) := (F_{u_n})_R(x).
\end{align}
is $k$-quadratic.} In order to see this, we put 
\[
V := \Quad(M_1,\dots,M_{n-1};M)
\]
and conclude from \eqref{EFUN} that the map $Q\:M_n \to V$, $u_n \mapsto F_{u_n}$ is $k$-quadratic. On the other hand, \eqref{EFEX} implies $F_x(u_n) = Q(u_n)_R(x)$, so $F_x$ is $Q$ composed with two $k$-linear maps: the natural map
\[
\Quad(M_1,\dots,M_{n-1};M) \longrightarrow \Quad(M_{1R},\dots,M_{n-1,R};M_R), \quad G \longmapsto G_R,
\] 
given and $k$-linear by the induction hypothesis, and the evaluation at $x$. This proves our claim.

By Cor.~\ref{c.BAQUA} we therefore obtain an $R$-quadratic map 
\[
(F_x)_R\:M_{nR} \longrightarrow (M_R)_R := M \otimes R \otimes R,
\]
which composes with the $R$-linear map
\begin{align}
\label{MUR} \mu\:M \otimes R \otimes R \longrightarrow M \otimes R = M_R, \quad u \otimes r \otimes s \longmapsto u \otimes rs
\end{align}
to yield a map $F_R\:M_{1R} \times \cdots \times M_{n-1,R} \times M_{nR} \to M_R$ via
\begin{align}
\label{EFR} F_R(x_1,\dots,x_{n-1},x_n) := \mu\big((F_{x_1,\dots,x_{n-1}})_R(x_n)\big)
\end{align}
for $x_i \in M_{iR}$, $1 \leq i \leq n$. By construction, $F_R$ is $R$-quadratic in $x_n$. Hence, writing $x_n = \sum r_lu_{lnR}$, $r_l \in R$, $u_{ln} \in M_n$, we obtain, after a straightforward computation,
\begin{align*}
(F_R)_{x_n} =\,\,&\sum r_l^2(F_{u_{ln}})_R + \sum_{l<m}r_lr_m\big((F_{u_{ln} + u_{mn}})_R - (F_{u_{ln}})_R - (F_{u_{mn}})_R\big), 
\end{align*}
which shows that the left-hand side of \eqref{EFR} is $R$-quadratic in $x_i$, $1 \leq i < n$. Thus $F_R$ is an $R$-$n$-quadratic map. Next we show that \eqref{MULQUA} commutes. Indeed, given $u_i \in M_i$, $1 \leq i \leq n$, we apply \eqref{EFR}, Cor.~\ref{c.BAQUA}, \eqref{MUR}, \eqref{EFEX}, the induction hypothesis and \eqref{EFUN} to obtain
\begin{align*}
F_R(u_{1R},\dots,u_{n-1,R},u_{nR}) =\,\,&\mu\big((F_{u_{1R},\dots,u_{n-1,R}})_R(u_{nR})\big) = \mu\Big(\big(F_{u_{1R},\dots,u_{n-1,R}}(u_n)\big)_R\Big) \\
=\,\,&F_{u_{1R},\dots,u_{n-1,R}}(u_n) = (F_{u_{nR}})_R(u_{1R},\dots,u_{n-1,R}) \\
=\,\,&F_{u_n}(u_1,\dots,u_{n-1})_R = F(u_1,\dots,u_{n-1},u_n)_R,
\end{align*}
and commutativity of \eqref{MULQUA} is proved.

It remains to establish uniqueness. Let $G,H\:M_{1R} \times \cdots \times M_{nR} \to M_R$ be $R$-$n$-quadratic maps both rendering the diagram
\begin{align}
\vcenter{\label{QUAUN} \xymatrix{M_1 \times \cdots \times M_n \ar[r]^(0.65){F} \ar[d]_{\can} & M \ar[d]^{\can} \\
M_{1R} \times \cdots \times M_{nR} \ar[r]_(0.65){G,H} & M_R.}}
\end{align}
commutative. For $u_n \in M_n$, using obvious notation, both 
\[
G_{u_{nR}}, H_{u_{nR}}\:M_{1R} \times \cdots \times M_{n-1,R} \longrightarrow M_R
\]
are $R$-$(n-1)$-quadratic maps rendering the diagram
\begin{align*}
\xymatrix{M_1 \times \cdots \times M_n \ar[rr]^(0.65){F} \ar[d]_{\can} && M \ar[d]^{\can} \\
M_{1R} \times \cdots \times M_{nR} \ar[rr]_(0.65){G_{u_{nR}},H_{u_{nR}}} && M_R.}
\end{align*}
commutative. From the uniqueness part of the induction hypothesis, we therefore conclude $G_{u_{nR}} = H_{u_{nR}}$. Now let $x_i \in M_{iR}$ for $1 \leq i \leq n$ and write $x_n = \sum r_lu_{nlR}$, $r_l \in R$, $u_{ln} \in M_n$. Then what we have just shown implies
\begin{align*}
G(x_1,\dots,x_{n-1},x_n) =\,\,&\sum r_l^2G(x_1,\dots,x_{n-1},u_{nlR}) + \sum_{l<m}\big(G(x_1,\dots,x_{n-1},u_{lnR} + u_{mnR}) \\
\,\,&- G(x_1,\dots,x_{n-1},u_{lnR}) - G(x_1,\dots x_{n-1},u_{mnR})\big) \\
=\,\,&\sum r_l^2G_{u_{nlR}}(x_1,\dots,x_{n-1}) + \sum_{l<m}\big(G_{(u_{ln} + u_{mn})_R}(x_1,\dots,x_{n-1}) \\
\,\,&-G_{u_{lnR}}(x_1,\dots,x_{n-1}) - G_{u_{mnR}}(x_1,\dots,x_{n-1}) \\
=\,\,&H(x_1,\dots,x_{n-1},x_n)
\end{align*}
This proves uniqueness and completes the induction. Finally, a straightforward verification using uniqueness of $F_R$ shows that the assignment $F \mapsto F_R$ is $k$-linear.
\end{sol}

\begin{sol}{pr.QUABIL} \label{sol.QUABIL} 
We first assume that $M$ is free and let $(e_i)_{i \in I}$ be a basis of $M$.  By the Axiom of Choice, we may choose a well-ordering on $I$.  Let $B\:M \times M \to N$ be the unique bilinear map from $M$ to $N$ given by the values
\[
B(e_i,e_j) =
\begin{cases}
Q(e_i) & \text{for}\;\; i = j \in I, \\
Q(e_i,e_j) & \text{for}\;\; i,j \in I, i < j, \\
0 & \text{for}\;\; i,j \in I, i > j
\end{cases}
\]
on the basis vectors. Then, for any $x = \sum\xi_ie_i \in M$, $\xi_i \in k$, we have
\begin{multline*}
B(x,x) = B(\sum\xi_ie_i,\sum\xi_ie_i) = \sum_{i,j \in I}\xi_i\xi_j B(e_i,e_j) \\ = \sum\xi_i^2 Q(e_i) + \sum_{i<j}\xi_i\xi_jQ(e_i,e_j) = Q(x),
\end{multline*}
as claimed. Now let $M$ be an arbitrary projective module. Then there exists a $k$-module $M^\prime$ making $M_0 := M \oplus M^\prime$ free. Let $Q^\prime\:M^\prime \to N$ be the zero quadratic map on $M^\prime$ with values in $N$ and consider the quadratic map $Q_0\:M_0 \to N$ canonically induced by $Q$ and $\pQ$. By the special case just treated we find a bilinear map $B_0\:M_0 \times M_0 \to N$ such that $Q_0(x_0) = B_0(x_0,x_0)$ for all $x_0 \in M_0$. But this implies $Q(x) = B(x,x)$ for all $x \in M$, where $B\:M \times M \to N$ is the restriction of $B_0$ to $M \times M$.
\end{sol}

\begin{sol}{pr.SEQUAFI} \label{sol.SEQUAFI}
(i) $\Rightarrow$ (iv) $\Rightarrow$ (iii). Obvious.

(iii) $\Rightarrow$ (ii). If $x \in V$ satisfies $q(x) = q(x,y) = 0$
for all $y \in V$, then by linearity $x_K \in V_K$ satisfies
$(q_K)(x_K) = (q_K)(x_K,y^\prime) = 0$ for all $y^\prime \in V_K$,
and (iii) gives $x_K = 0$, hence $x = 0$. Thus $q$ is
non-degenerate, and it remains to show $\dim_F(\Rad(Dq)) \leq 1$.
Since the radical of a symmetric or skew-symmetric bilinear form
\emph{over a field} is compatible with base change, it suffices to
show $\dim_K(\Rad(Dq_K)) \leq 1$. But $q_K$ is non-degenerate by
hypothesis, so its restriction to $\Rad(Dq_K)$ is anisotropic. On
the other hand, $K$ being algebraically closed, all quadratic forms
of (possibly infinite) dimension $> 1$ over $K$ are isotropic. The assertion follows.

(ii) $\Rightarrow$ (i). By hypothesis, there exists an element $u
\in V$ satisfying $\Rad(Dq) = Fu$, and if $u \neq 0$, then $q(u)
\neq 0$ since $q$ is non-degenerate. Now let $L$ be any field
containing $F$. Then $\Rad(Dq_L) = Lu_L$. If $x^\prime \in V_L$
satisfies $(q_L)(x^\prime) = (q_L)(x^\prime,y^\prime) = 0$ for all
$y^\prime \in V_L$ then, in particular, $x^\prime \in \Rad(Dq_L) =
Lu_L$, so some $a \in L$ has $x^\prime = a(u_L)$. This implies $0 =
q_L(x^\prime) = a^2q(u)$, hence $a = 0$ or $q(u) = 0$, which in turn
yields $a = 0$ or $u = 0$, hence $x^\prime = 0$. Therefore $q_L$ is
non-degenerate, and $q$ is non-singular.
\end{sol}

\begin{sol}{pr.LIQUA} \label{sol.LIQUA}
(a) Let $x \in V$. The assertion is obvious if $x = 0$ or $x$ is anisotropic. Hence we may assume that $x$ is isotropic. Then there is a hyperbolic plane $H \subseteq V$ containing $x$, and the assumption $n \geq 3$ implies that $H^\perp \neq \{0\}$ is a quadratic space in its own right. Let $y \in H^\perp$ be anisotropic. Then so are $-y$ and $x + y$, and $x = (x + y) - y$ is a decomposition of the desired kind.

(b) We may assume $q(x_1,x_2) = 0$. Let $H \subseteq V$ be a hyperbolic subspace of dimension $4$. By hypothesis, $q$ is regular and isotropic, hence universal, on $H^\perp \neq \{0\}$.  There are two cases, depending on $x_1$ and $x_2$ are linearly dependent.

Suppose first that $x_1$ and $x_2$ are linearly dependent. Then $x_2 = \alpha x_1$ for some $\alpha \in F^\times$, and we are in characteristic $2$. The Witt Extension Theorem \cite[Thm.~8.3]{MR2427530} allows us to assume $x_1 \in H$. Pick $y_1 \in H$ such that $q(x_1,y_1) \neq 0$, and $y_3 \in H^\perp$ such that $y := y_1 + y_3 \in V$ is anisotropic. Then $q(x_1,y \neq 0 \neq q(y,x_2)$. 

Suppose now that $x_1$ and $x_2$ are linearly independent. Write $H = H_1 \perp H_2$ with hyperbolic planes $H_1,H_2$. Since $q$ is universal on each of the subspaces $H_1, H_2$, the Witt Extension Theorem again justifies the assumption $x_i \in H_i$ for $i = 1,2$. Pick $y_i \in H_i$ such that $q(x_i,y_i) \neq 0$ for $i = 1,2$ and $y_3 \in H^\perp$ such that $y := y_1 + y_2 + y_3$ is anisotropic. Then $q(x_1,y) = q(x_1,y_1) \neq 0 \neq q(y_2,x_2) = q(y,x_2)$, as desired. 

\medskip

We now show that \emph{in} (a) \emph{the assumption $n \geq 3$ cannot be avoided.} Indeed, let $(V,q)$ be the hyperbolic plane over $\IF_2$, the field with two elements. $V$ contains precisely two isotropic vectors relative to $q$, denoted by $e_1,e_2$, which form a hyperbolic pair, and precisely one anisotropic one, denoted by $a = e_1 + e_2$.  Hence any finite sum of anisotropic vectors in $(V,q)$ is either zero or equal to $a$, so neither $e_1$ nor $e_2$ can be written in this form.

Finally, we show that \emph{in} (b) \emph{the assumption of a hyperbolic subspace of dimension at least 6 cannot be avoided.} Indeed, with $(V,q)$ as in the previous paragraph, put $(W,Q) = (V,q) \perp (V,q)$. Then $W$ contains precisely four anisotropic vectors relative to $Q$, namely $(a,e_i)$, $(e_j,a)$ for $i,j = 1,2$. Moreover, $Q((a,e_i),(a,e_j)) = 1 - \delta_{ij}$ and $Q((a,e_i),(e_j,a)) = 0$. Hence $x_1 := (a,e_1)$ and $x_2 := (e_1,a)$ are anisotropic as well as orthogonal, and there is a unique anisotropic vector $y \in W$ satisfying $Q(x_1,y) \neq 0$, namely, $y = (a,e_2)$. But $Q(y,x_2) = 0$, so the conclusion of (b) does not hold.
\end{sol}

\begin{sol}{skip.pr.reg} \label{skip.sol.reg}
Suppose first that $k$ is local, i.e., has a unique maximal ideal $\mfm$.  Then $M \cong k^n$ for some $n$ and $Dq$ is identified with a symmetric $S \in \Mat_n(k)$.  By hypothesis, $q_{k(\mfm)}$ is non-singular, so by the remark after Exc.~ \ref{pr.SEQUAFI} it is regular.  Therefore $\det S$ is not zero in $k/\mfm$, i.e., $\det S$ is invertible in $k$ and the map $M \to M^*$ induced by $Dq$ is an isomorphism.

For general $k$, the map $M \to M^*$ induced by $Dq$  has the property that it is an isomorphism after base change to $k_\mfm$ for every maximal ideal $\mfm$ of $k$, so it is an isomorphism by \cite[II.3, Thm.~1]{MR0360549}.
\end{sol}

\begin{sol}{skip.pr.antiregular} \label{skip.sol.antiregular}
(a): For each $v \in \Rad(Dq)$, we have $q(v) = \frac12 Dq(v,v) = 0$.

(b): Replace $(M, q)$ with $(M', q\vert_{M'})$ and so assume that $M$ is free of finite rank, $Dq$ is identically zero, and $q$ is anisotropic.  Pick a basis $m_1, \ldots, m_r$ of $M$ and define $\alpha_i := q(m_i)$.  Then for $c_1, \ldots, c_r \in k$ we have $q(\sum c_i m_i) = \sum q(c_i m_i) = c_i^2 \alpha_i$, as required.  Finally, if $\alpha_i = \alpha_j \beta^2$ for some $i \ne j$ and $\beta \in \kx$, then $q(m_i + \beta m_j) = \alpha_i + \alpha_j \beta^2 = 0$, whence $m_i + \beta m_j = 0$, a contradiction.
\end{sol}

\begin{sol}{skip.pr.cgp} \label{skip.sol.cgp}
(i): We may identify each $x \in \Hom_{\kmod}(M/\Rad(q), \Rad(q))$ with an element of $\Hom_{\kmod}(M, M)$ by composing it with the natural maps $M \to M/\Rad(q)$ and $\Rad(q) \to M$.  For such an $x$, define
 $g_x := x + \Id_M \in \Hom_{\kmod}(M, M)$.  Then, for $m \in M$,
\[
q(g_x(m)) = q(m + x(m)) = q(m) + q(x(m)) + q(m, x(m)) = q(m),
\]
so $q \circ g_x = q$.

Also, $x\in \Hom_{\kmod}(M, M)$ vanishes on $\Rad(q)$ so $x^2 = 0$ and we have $g_x (\Id_M - x) = \Id_M - x^2 = \Id_M$ and similarly $(\Id_M - x) g_x = \Id_M$.  That is, $g_x$ is invertible, so it belongs to $\rmO(Q)$.

Evidently the map $x \mapsto g_x$ is injective.  For $x, y \in \Hom_{\kmod}(M, \Rad(q))$, we have $g_x g_y = g_{xy + x + y}$, so the image $R$ of $\Hom_{\kmod}(M, \Rad(q))$ in $\rmO(Q)$ is a subgroup.

For $\eta \in \rmO(Q)$ and $m \in M$, we have 
\[
\eta g_x \eta^{-1}(m) = \eta x \eta^{-1}(m) + m = g_{\eta x \eta^{-1}}(m),
\]
so $R$ is a normal subgroup.

(ii): Every $\eta \in \rmO(Q)$ maps $\Rad(q)$ to itself, whence restriction gives a well-defined homomorphism $\phi$.  Now suppose that $\eta \in \ker \phi$ and set $x := \eta - \Id_M$.  Since $\eta$ maps to the identity in $\rmO(\Qbar)$, it follows that $\eta m - m \in \Rad(q)$ for all $m \in M$, i.e., $x \colon M \to \Rad(q)$.  Since $\eta$ maps to the identity in $\GL(\Rad(q))$, it follows that for $m \in \Rad(q)$ we have $xm = \eta m - m = 0$.  That is, $x \colon M/\Rad(q) \to \Rad(q)$ and $\eta$ is in $R$.

(iii): The quotient map $M \to M/\Rad(q)$ identifies the quadratic module $(M', q\vert_{M'})$ with $(M/\Rad(q), \qbar)$, providing a recipe that takes an element of $\GL(\Rad(q)) \times \rmO(\Qbar)$ and identifies it with an element of $\GL(\Rad(q)) \times \rmO(M', q\vert_{M'}) \subseteq \rmO(Q)$.  This shows that the map $\phi$ is surjective.
\end{sol}

\begin{sol}{skip.pr.nonLG} \label{skip.sol.nonLG}
We first claim: \emph{If a ring $k$ is LG, then for every $c \in k$ the polynomial $x^2 + c$ represents a unit.}  By the LG property, it suffices to verify this in case $k$ is a field.  In that case: If $c \ne -1$, take $x = 1$.  If $c = -1$, take $x = 0$. 

(a): For $k = \IZ$, take $f = x^2 + 2$.  Then for every $z \in \IZ$, $f(z) = z^2 + 2 \ge 2$, so $f(z)$ is not a unit.  The previous paragraph shows that $\IZ$ is not LG.

(b): For $k = R[\bft]$, we have $k^\times = R^\times$ since $R$ is an integral domain.  Then $x^2 + \bft$ represents polynomials of degree $\ge 1$, hence non-units, and we again conclude that $k$ is not LG.

\medskip

\noindent{\textbf{Extension.}}  Part (b) of this exercise can be generalized in the following way.  \emph{Suppose $k$ is a ring with a discrete valuation $v$.  If there exists an element $c \in k$ such that $v(c)$ is odd and negative, then $k$ is not an LG ring.}

A discrete valuation on a ring $k$ is defined in \cite[\S{VI.3.1}, Def.~1]{MR643362}.  It is a $v \: k \to \IZ \cup \{ \infty \}$ such that following axioms hold for all $x, y \in k$:
\begin{enumerate}[(I)]
\item $v(xy) = v(x) + v(y)$.
\item \label{DVR.2} $v(x+y) \ge \min \{ v(x), v(y) \}$.
\item $v(1) = 0$ and $v(0) = \infty$.
\end{enumerate}
For such a $v$, we have equality in \ref{DVR.2} if $v(x) \ne v(y)$ \cite[\S{VI.3.1}, Prop.~1]{MR643362}.

The italicized claim generalizes part (b) by taking the valuation to be $v(f) = -\deg f$ (the ``$1/\bft$-adic valuation'') and $c = \bft$, which has $v(c) = -1$.
\begin{proof}
Assume $c \in k$ has $v(c)$ odd and negative. Assume further for some $x \in k$ that $x^2 + c$ is a unit. Since $v(x^2) = 2v(x)$ is even and hence distinct from $v(c)$, we conclude $0 = v(x^2 + c) = \min\{2v(x),v(c)\}$, which cannot be $v(c)$ since this is negative, hence must be $2v(x)$, forcing the minimum in question to be $0 = v(c)$, a contradiction.
\end{proof}
\end{sol}

\vfill

\solnsec{Section~\ref{s.POLA}}

\begin{sol}{pr.ZERLIN} \label{sol.ZERLIN}  Setting $\bfT= (\bft_1,\dots,\bft_n)$ as usual, $x = (x_1,\dots,x_n)$, $y = (x_{n+1},\dots,x_{n+p})$, $x_i \in M_R$ for $1 \leq i \leq n + p$ and applying Thm.~\ref{t.BILEXP}, we obtain
\begin{align*}
\sum_{\nu\in\IN^n}\bfT^\nu(\Pi^{(\nu,0)} f)_R(x,y) =\,\,&\sum_{(\nu,\nu^\prime)\in\IN^{n+p}}\bft_1^{\nu_1}\cdots\bft_n^{\nu_n}0^{\nu_{n+1}}\cdots 0^{\nu_{n+p}} (\Pi^{(\nu,\nu^\prime)} f)_R(x,y) \\
=\,\,&f_{R[\bfT]}(\sum_{i=1}^n\bft_ix_i + \sum_{i=1}^p 0x_{n+i}) = f_R(\sum_{i=1}^n \bft_ix_i) \\
=\,\,&\sum_{\nu\in\IN^n}\bfT^\nu(\Pi^\nu f)_R(x), 
\end{align*}
and comparing coefficients, the assertion follows.
\end{sol}

\begin{sol}{pr.HOCO} \label{sol.HOCO} For $d \in \IN$ define $f_d := \Pi^{(d)} f$, which by Thm.~\ref{t.BILEXP} and \ref{ss.LIPOLA}~(a) is a homogeneous polynomial law $M \to N$ of degree $d$ over $k$. Moreover, the family $(f_d)_{d\geq 0}$ is locally finite, and (\ref{t.BILEXP}.\ref{EFERSU}) for $n = 1$, $r_1 = 1_R$ implies $f = \sum_{d\geq 0} f_d$. This proves existence. Conversely, let $(f_d)_{d\geq 0}$ be any family of polynomial laws over $k$ with the desired properties. For $R \in \kalg$, $r \in R$, $x \in M_R$ we deduce
\[
f_R(rx) = \sum_{d\geq 0} f_{dR}(rx) = \sum_{d\geq 0} r^df_{dR}(x),
\]
and the uniqueness property of Thm.~\ref{t.BILEXP} yields $f_d = \Pi^{(d)}f$ for all $d \in \IN$. This shows uniqueness and solves the first part of the problem.

It remains to exhibit an example of a polynomial law $f$ having $f_d \neq 0$ for all $d \in \IN$. To this end, consider a free $k$-module $M$ of countably infinite rank and let $(e_i)_{i\in\IN}$ be a basis of $M$ over $k$. For $d \in \IN$ and $R \in \kalg$, define a set map $f_{dR}\:M_R \to R$ by
\[
f_{dR}(\sum_{i\in\IN}r_ie_{iR}) := r_d^d,
\]
where $(r_i)_{i\in\IN}$ is an arbitrary sequence of finite support in $R$. The family $(f_{dR})_{R\in\kalg}$ is clearly a homogeneous polynomial law $f_d\:M \to k$ of degree $d$ over $k$, and the family $(f_d)_{d\geq 0}$ is locally finite. Hence, by \ref{ss.POFIN}, $f := \sum_{d\geq 0}f_d\:M \to k$ exists as a polynomial law over $k$ and has the desired property.
\end{sol}

\begin{sol}{pr.COPOLA} \label{sol.COPOLA} The family of set maps $\hat{w}_R\:M_R \to N_R$ for $w \in N$ clearly varies functorially with $R \in \kalg$ and hence is a polynomial law $\hat{w}\:M \to N$ such that $\hat{w}_R(rx) = w_R = r^0\hat{w}_R(x)$ for all $r \in R$, $x \in M_R$. Thus $\hat{w}$ is homogeneous of degree $0$. Conversely, let $f\:M \to N$ be any homogeneous polynomial law of degree $0$ over $k$. Then $w := f_k(0) \in N$, and for all $R \in \kalg$, $x \in M_R$, we apply (\ref{ss.COPOLA}.\ref{FUBA}) to obtain $w_R = f_k(0)_R = f_R(0) = f_R(0x)= 0^0f_R(x) = f_R(x)$, hence $\hat{w}_R = f_R$ and then $\hat{w} = f$.
\end{sol}

\begin{sol}{pr.MULI} \label{sol.MULI} (a) The first part is obvious. To establish the second, we begin by treating the case $n = 1$, so let $f\:M \to N$ be a homogeneous polynomial law of degree $1$. Then \ref{ss.LIPOLA}~(c) implies $\Pi^\nu f = 0$ for all $\nu \in \IN^2 \setminus \{(1,0),(0,1)\}$, and (\ref{ss.BILI}.\ref{BILI}) yields
\[
f_k(\alpha x + \beta y) = \alpha(\Pi^{(1,0)} f)_k(x,y) + \beta(\Pi^{(0,1)}f)_k(x,y)
\]
for all $\alpha,\beta \in k$ and all $x,y \in M$. Specializing $\alpha = 1$, $\beta = 0$ and $\alpha = 0$, $\beta = 1$, we conclude $(\Pi^{(1,0)}f)_k(x,y) = f_k(x)$, $(\Pi^{(0,1)}f)_k(x,y) = f_k(y)$. Hence $f_k\:M \to N$ is a $k$-linear map. Since $f \otimes R\:M_R \to N_R$, for $R \in \kalg$, is a homogeneous polynomial law of degree $1$ over $R$, it follows, therefore, that $f_R\:M_R \to N_R$ is an $R$-linear map. Moreover, for $r \in R$, $x \in M$, we deduce $f_R(x \otimes r) = f_R(rx_R) = rf_R(x_R) = rf_k(x)_R = f_k(x) \otimes r = (f_k \otimes \Eins_R)(x \otimes r)$, hence $f_R = (f_k)_R$. Thus $f$ is the polynomial law derived from the scalar extensions of the linear map $f_k$. This settles the case $n = 1$.

Now let $n$ be arbitrary and $\mu\:M_1 \times \cdots \times M_n \to N$ be a multi-homogeneous polynomial law of multi-degree $\hat{1} = (1,\dots,1)$. For $1 \leq i \leq n$ and fixed $v_j \in M_j$ ($1 \leq j \leq n$, $j \neq i$), it is straightforward to verify, using (\ref{ss.STANID}.\ref{VEONVE}), that $f_i\:M_i \to N$ defined by
\[
f_{iR}(x) := \mu_R(v_{1R},\dots,v_{i-1,R},x,v_{i+1,R},\dots,v_{nR})
\]
for $R \in \kalg$, $x \in M_R$, is a homogeneous polynomial law of degree $1$. Therefore, by the special case treated above, $f_{ik}\:M_i \to N$ is a $k$-linear map, forcing $\mu_k\:M_1 \times \cdots \times M_n \to N$ to be $k$-multi-linear. Applying this to the extended polynomial law $\mu \otimes R\:M_{1R} \times \cdots \times M_{nR} \to N_R$ over $R$, we conclude that $\mu_R\:M_{1R} \times \cdots \times M_{nR} \to N_R$ is an $R$-multi-linear map. Hence, given $r_{ij} \in R$, $v_{ij} \in M_j$ for $1 \leq j \leq n$ and finitely many indices $i$, we obtain
\begin{align*}
\mu_R(\sum_ir_{i1}v_{i1R},\dots,\sum_ir_{in}v_{inR}) =\,\,&\sum_{i_1,\dots,i_n}r_{i_11}\cdots r_{i_nn}\mu_k(v_{i_11},\dots,v_{i_nn})_R \\
=\,\,&(\mu_k \otimes R)(\sum_ir_{i1}v_{i1R},\dots,\sum_ir_{in}v_{inR}).
\end{align*}
This proves $\mu_R = \mu_k \otimes R$ and completes the solution of (a).

(b) Let $Q\:M \to N$ be a quadratic map. We show that $\tilde{Q}\:M \to N$ is a polynomial law over $k$ by letting $\vph\:R \to S$ be a morphism in $\kalg$ and $x,x^\prime \in M$, $r,r^\prime \in R$. The Cor.~\ref{c.BAQUA} implies
\begin{align*}
(\Eins_N \otimes \vph) \circ Q_R(x \otimes r) =\,\,&(\Eins_N \otimes \vph) \circ Q_R(rx_R) = (\Eins_N \otimes \vph)\big(r^2Q(x)_R\big) 
= (\Eins_N \otimes \vph)\big(Q(x) \otimes r^2\big) \\
=\,\,& Q(x) \otimes \vph(r)^2 = \vph(r)^2Q(x)_S = \vph(r)^2Q_S(x_S) \\
=\,\,&Q_S\big(\vph(r)x_S\big) = Q_S\big(x \otimes \vph(r)\big) = Q_S \circ (\Eins_M \otimes \vph)(x \otimes r)
\end{align*}
and, similarly, $(\Eins_N \otimes \vph) \circ Q_R(x \otimes r,x^\prime \otimes r^\prime) = Q_S \circ (\Eins_M \otimes \vph)(x \otimes r,x^\prime \otimes r^\prime)$. This proves that $\tilde{Q}$ is indeed a polynomial law over $k$ and obviously homogeneous of degree $2$. Conversely, let $f\:M \to N$ be a homogeneous polynomial law of degree $2$ over $k$ and consider the set map $Q := f_k\:M \to N$. By \ref{ss.LIPOLA}~(c) and (\ref{ss.BILI}.\ref{BILI}), we have 
\[
Q(\alpha x + \beta y) = \alpha^2(\Pi^{(2,0)}f)_k(x,y) + \alpha\beta(\Pi^{(1,1)}f)_k(x,y) + \beta^2(\Pi^{(0,2)}f)_k(x,y)
\] 
for all $\alpha,\beta \in k$ and all $x,y \in M$. After specializing $\alpha = 1$, $\beta = 0$ and $\alpha = 0$, $\beta = 1$ and observing (\ref{ss.TODEHO}.\ref{DEEFPI}), this may be rewritten as 
\[
Q(\alpha x + \beta y) = \alpha^2Q(x) + \alpha\beta(Df)_k(x,y) + \beta^2Q(y),
\]
where $(Df)_k = (\Pi^{(1,1)}f)_k$ is $k$-bilinear by \ref{ss.LIPOLA}~(a) and (a). This shows that $Q$ is a quadratic map with $DQ = (Df)_k$. For $R \in \kalg$, the $R$-quadratic map $Q_R\:M_R \to N_R$ is uniquely determined by (\ref{c.BAQUA}.\ref{BAQUA}), which proves $Q_R = f_R$, hence $\tilde{Q} = f$.
\end{sol}

\begin{sol}{pr.HIDIDE} \label{sol.HIDIDE} (a) Applying Thm.~\ref{t.BILEXP}, we deduce
\begin{align*}
\sum_{\nu_1,\dots,\nu_p\geq 0}\bft_1^{\nu_1}\cdots\bft_p^{\nu_p}\Big(\,\,&\sum_{\nu_{p+1},\dots,\nu_n\geq 0}\bft_{p+1}^{\nu_{p+1}}\cdots\bft_n^{\nu_n}(\Pi^{(\nu_1,\dots,\nu_n)}f)_R(x_1,\dots,x_p,x,\dots,x)\Big) \\
=\,\,&\sum_{\nu\in\IN^n}\bfT^\nu(\Pi^\nu f)(x_1,\dots,x_p,x,\dots,x) = f_{R[\bfT]}\big(\sum_{j=1}^p\bft_jx_j + (\sum_{j=p+1}^n \bft_j)x\big) \\
=\,\,&\sum_{\nu_1,\dots,\nu_p,\nu_{p+1}\geq 0}\bft_1^{\nu_1}\cdots\bft_p^{\nu_p}(\sum_{j=p+1}^n\bft_j)^{\nu_{p+1}} (\Pi^{(\nu_1,\dots,\nu_p,\nu_{p+1})}f)_R(x_1,\dots,x_p,x) \\
=\,\,&\sum_{\nu_1,\dots,\nu_p\geq 0}\bft_1^{\nu_1}\cdots\bft_p^{\nu_p}\Big(\sum_{i\geq 0}(\sum_{j=p+1}^n\bft_j)^i(\Pi^{(\nu_1,\dots,\nu_p,i)}f)_R(x_1,\dots,x_p,x)\Big).
\end{align*}
Comparing coefficients and using the multi-nomial formula (Abramowitz-Segun \cite[p.~823]{MR0167642}), this implies, for all $\nu_1,\dots,\nu_p \in \IN$,
\begin{align*}
\sum_{\nu_{p+1},\dots,\nu_n\geq 0}&\bft_{p+1}^{\nu_{p+1}}\cdots\bft_n^{\nu_n}(\Pi^{(\nu_1,\dots,\nu_n)}f)_R(x_1,\dots,x_p,x,\dots,x) \\
&= \sum_{i\geq 0}(\sum_{j=p+1}^n\bft_j)^i(\Pi^{(\nu_1,\dots,\nu_p,i)}f)_R(x_1,\dots,x_p,x) \\
&= \sum_{i\geq 0}\Big(\sum_{\nu_{p+1},\dots,\nu_n\geq 0,\nu_{p+1}+\cdots +\nu_n = i}\frac{i!}{\nu_{p+1}!\cdots\nu_n!}\bft_{p+1}^{\nu_{p+1}}\cdots\bft_n^{\nu_n}(\Pi^{(\nu_1,\dots,\nu_p,i)}f)_R(x_1,\dots,x_p,x)\Big) \\
&= \sum_{\nu_{p+1},\dots,\nu_n\geq 0}\bft_{p+1}^{\nu_{p+1}}\cdots\bft_n^{\nu_n}\frac{(\nu_{p+1} + \cdots + \nu_n)!}{\nu_{p+1}! \cdots \nu_n!}(\Pi ^{(\nu_1,\dots,\nu_p,\nu_{p+1} + \cdots + \nu_n)}f)_R(x_1,\dots,x_p,x),
\end{align*}
and comparing coefficients again completes the proof of (a).

(b) Let $\vph\:R \to S$ be a morphism in $\kalg$ and $x \in M_R$. Since $D^pf$ is a polynomial law over $k$, an application of (\ref{ss.STANID}.\ref{ONEMON}) yields
\begin{align*}
(\partial_y^{[p]}f)_S \circ (\Eins_M \otimes \vph)(x) =\,\,&(D^pf)_S\big((\Eins_M \otimes \vph)(x),y_S\big) = (D^pf)_S\big((\Eins_M \otimes \vph)(x),(\Eins_M \otimes \vph)(y_R)\big) \\
=\,\,&(D^pf)_S \circ (\Eins_{M\times M} \otimes \vph)(x,y_R) = (\Eins_N \otimes \vph) \circ (D^pf)_R(x,y_R) \\
=\,\,&(\Eins_N \otimes \vph) \circ (\partial_y^{[p]})_R(x).
\end{align*}
Thus $\partial_y^{[p]}f$ is a polynomial law over $k$ and the map $\partial_y^{[p]}\:\Pol_k(M,N) \to \Pol_k(M,N)$ is clearly $k$-linear.

Next combine Prop.~\ref{p.ITDIDE} for $y_1 = \cdots = y_p = y$ with (a) to conclude
\begin{align*}
\big((\partial_y)^pf\big)_R(x) =\,\,&(\partial_y\cdots\partial_yf)_R(x) = \sum_{i\geq 0}(\Pi^{(i,1,\dots,1)}f)_R(x,y_R,\dots,y_R) = p!\sum_{i\geq 0}(\Pi^{(i,p)}f)_R(x,y_R) \\
=\,\,&p!(D^pf)_R(x,y_R) = p!(\partial_y^{[p]}f)_R(x),
\end{align*}
which implies $(\partial_y)^p = p!\partial_y^{[p]}$ as linear maps $\Pol_k(M,N) \to \Pol_k(M,N)$. This completes the solution of (b).

(c) Let $R \in \kalg$, $\bft,\bft_1,\bft_2$ be independent variables and $\bfT= (\bft_1,\bft_2)$. For $x,y_1,y_2 \in M_R$ we combine the Taylor expansion formula (\ref{ss.TODE}.\ref{TAYL}) with Thm.~\ref{t.BILEXP} and obtain
\begin{align*}
\sum_{n\geq 0}\bft^n(D^nf)_{R[\bfT]}(x,\bft_1y_1 + \bft_2y_2) =\,\,&f_{R[\bft,\bfT]}\big(x + \bft(\bft_1y_1 + \bft_2y_2)\big) = f_{R[\bft,\bfT]}\big(x + (\bft\bft_1)y_1 + (\bft\bft_2)y_2\big) \\
=\,\,&\sum_{\nu_0,\nu_1,\nu_2\geq 0}1^{\nu_0}(\bft\bft_1)^{\nu_1}(\bft\bft_2)^{\nu_2}(\Pi^{(\nu_0,\nu_1,\nu_2)}f)_R(x,y_1,y_2) \\
=\,\,&\sum_{\nu_1,\nu_2\geq 0}\bft^{\nu_1+\nu_2}\bft_1^{\nu_1}\bft_2^{\nu_2}\sum_{i\geq 0}(\Pi^{(i,\nu_1,\nu_2)}f)_R(x,y_1,y_2) \\
=\,\,&\sum_{n\geq 0}\bft^n\sum_{j=0}^n\bft_1^j\bft_2^{n-j}\sum_{i\geq 0}(\Pi^{(i,j,n-j)}f)_R(x,y_1,y_2)
\end{align*}
Comparing coefficients of $\bft^n$, we therefore end up with the formula
\begin{align}
\label{DIFENX} (D^nf)_{R[\bfT]}(x,\bft_1y_1 + \bft_2y_2) = \sum_{j=0}^n\bft_1^j\bft_2^{n-j}\sum_{i\geq 0}(\Pi^{(i,j,n-j)}f)_R(x,y_1,y_2)
\end{align}
This formula will now be applied in the special case $n = 2$. From Exc.~\ref{pr.ZERLIN} combined with \ref{ss.LIPOLA}~(b) we deduce, for all $i \in \IN$,
\begin{align*}
(\Pi^{(i,2,0)}f)_R(x,y_1,y_2) =\,\,&(\Pi^{(i,2)}f)_R(x,y_1), \\ 
(\Pi^{(i,0,2)}f)_R(x,y_1,y_2) =\,\,&(\Pi^{(i,2,0)}f)_R(x,y_2,y_1) = (\Pi^{(i,2)}f)_R(x,y_2). 
\end{align*}
Hence \eqref{DIFENX} for $n = 2$ implies
\begin{align*}
(D^2f)_{R[\bfT]}(x,\bft_1y_1 + \bft_2y_2) =\,\,&\bft_1^2\sum_{i\geq 0}(\Pi^{(i,2)}f)_R(x,y_1) + \bft_1\bft_2\sum_{i\geq 0}(\Pi^{(i,1,1)}f)_R(x,y_1,y_2) \\
\,\,&+ \bft_2^2\sum_{i\geq 0}(\Pi^{(i,2)}f)_R(x,y_2) \\
=\,\,&\bft_1^2(D^2f)_R(x,y_1) + \bft_1\bft_2\sum_{i\geq 0}(\Pi^{(i,1,1)}f)_R(x,y_1,y_2) + \bft_2^2(D^2f)_R(x,y_2).
\end{align*}
In particular, applying Prop.~\ref{p.ITDIDE},
\begin{align*}
(\partial_{y+z}^{[2]}f)_R(x) =\,\,&(D^2f)_R(x,y_R + z_R) \\
=\,\,&(D^2f)_R(x,y_R) + \sum_{i\geq 0}(\Pi^{(i,1,1)}f)_R(x,y_R,z_R) + (D^2f)_R(x,z_R) \\
=\,\,&\Big((\partial_y^{[2]}f)_R + (\partial_y\partial_zf)_R + (\partial_z^{[2]}f)_R\Big)(x),
\end{align*}
which completes the proof of (c).
\end{sol}

\begin{sol}{pr.POLINF} \label{sol.POLINF} (a) Let $R \in \Kalg$ and $x \in V_R$. We must show $f_R(x) = 0$. To this end, pick $K$-bases $(v_i)_{i\in I}$ of $V$, $(w_j)_{j\in J}$ of $W$. and write
\begin{align}
\label{EXSUMI} x = \sum_{i\in I}r_iv_{iR}, \quad f_R(x) = \sum_{j\in J}s_jw_{jR}
\end{align}
with $r_i,s_j \in R$, $i \in I$, $j \in J$. Then
\[
I_0 := \{i \in I \mid r_i \neq 0\} \subseteq I, \quad J_0 := \{j \in J \mid s_j \neq 0\} \subseteq J
\]
are finite subsets, forcing
\begin{align}
\label{VENAU} V_0 := \sum_{i\in I_0}Kv_i \subseteq V, \quad W_0 := \sum_{j\in J_0}Kw_j \subseteq W
\end{align}
to be finite-dimensional subspaces of $V,W$, respectively. Writing $\vep_0\:V_0 \to V$ for the inclusion and $\pi_0\:W \to W_0$ for the projection along $\sum_{j\notin J_0}Kw_j$,
\[
f_0 := \pi_0 \circ f \circ \vep_0\:V_0 \longrightarrow W_0
\]
is a polynomial law over $K$ such that $f_{0K} = \pi_{0K} \circ f_K \circ \vep_{0K} = 0$ as a set map from $V_0$ to $W_0$. Since $V_0$ and $W_0$ are both finite-dimensional, we conclude from Corollary~\ref{c.POMALA} that $f_0 = 0$ as a polynomial law over $K$. On the other hand, $x \in V_{0R}$, $f_R(x) \in W_{0R}$ by \eqref{EXSUMI}, \eqref{VENAU}, which implies $f_R(x) = f_{0R}(x) = 0$, as claimed.

(b) By (a) we may assume $f_K \neq 0$ and must show $\vph_K = 0$. Let $(v_i)_{i\in I}$, $(w_j)_{j\in J}$ be $K$-bases of $V,W$, respectively, and choose $x \in V$ such that
\begin{align}
\label{EFKAX} f_K(x) \neq 0.
\end{align}
Write
\begin{align}
\label{EXSUAL} x = \sum_{i\in I}\alpha_iv_i, \quad f_K(x) = \sum_{j\in J}\beta_jw_j &&(\alpha_i,\beta_j \in K, i \in I, j \in J).
\end{align}
For any $y \in V$, we have
\begin{align}
\label{YPSUGA} y = \sum_{i\in I}\gamma_iv_i, \quad f_K(y) = \sum_{j\in J}\delta_jw_j &&(\gamma_i,\delta_j \in K, i \in I, j \in J)
\end{align}
and must show $\vph_K(y) = 0$. Consider the finite subsets
\[
I_0 := \{i \in I \mid \alpha_i \neq 0\; \text{or}\; \gamma_i \neq 0\} \subseteq I, \quad J_0 := \{j \in J \mid \beta_j \neq 0\} \subseteq J,
\]
making
\begin{align}
\label{VEKAVE} V_0 := \sum_{i\in I_0}Kv_i, \quad W_0 = \sum_{j\in J_0}Kw_j
\end{align}
finite-dimensional subspaces of $V,W$, respectively, such that
\begin{align}
\label{EXYPVE} x,y \in V_0, \quad f_K(x) \in W_0.
\end{align}
Writing as before $\vep_0\:V_0 \to V$ for the inclusion and $\pi_0\:W \to W_0$ for the projection along $\sum_{j\notin J_0}Kw_j$, it follows that $f_0 := \pi_0 \circ f \circ \vep_0\:V_0 \to W_0$ is a polynomial law over $K$ having $f_{0K}(x) \neq 0$; in particular, $f_{0K} \neq 0$. On the other hand, $\vph_0 := \vph \circ \vep_0\:V_0 \to K$ is a scalar polynomial law over $K$, and one checks that $(\vph_0f_0)_K = 0$, which implies $\vph_{0K}(z_0) = 0$ for all $z_0 \in V_0$ having $f_{0K}(z_0) \neq 0$. But Corollary~\ref{c.POMALA} guarantees that $f_{0K}\:V_0 \to W_0$ is a non-zero polynomial map. By Zariski density, therefore, we obtain $\vph_K(z_0) = \vph_{0K}(z_0) = 0$ for all $z_0 \in V_0$; in particular, \eqref{EXYPVE} yields $\vph_K(y) = 0$, as desired.
\end{sol}

\begin{sol}{pr.POLALOC} \label{sol.POLALOC} 
(a) $f = 0$ clearly implies $f_\mfp := f \otimes k_\mfp = 0$ for all $\mfp \in \Spec(k)$. Before proving the converse, we show
\begin{enumerate}[label=($\ast$)]
\item \label{sol.POLALOC.ast} If $f_{k_\mfp} = 0$ as a set map $M_\mfp \to N_\mfp$ for all $\mfp \in \Spec(k)$, then $f_k = 0$ as a set map $M \to N$.
\end{enumerate}
Indeed, by the very definition of a polynomial law, the diagram
\[
\xymatrix{
M \ar[r]^{f_k} \ar[d]_{\mathrm{can_M}} & N \ar[d]^{\mathrm{can_N}} \\
M_\mfp \ar[r]_{f_{k_\mfp} = 0} & N_\mfp}
\]
commutes for all $\mfp \in \Spec(k)$. Hence for all $x \in M$ we have $f_k(x)/1 = 0$ in $N_\mfp$, which means $s_\mfp f_k(x) = 0$ for some $s_\mfp \in k \setminus \mfp$. It follows that the ideal generated by the $s_\mfp$, $\mfp \in \Spec(k)$, is all of $k$, since it cannot be contained in any prime ideal of $k$. Thus there exists a family $(\alpha_\mfp)_{\mfp \in \Spec(k)}$ with finite support such that $\sum_\mfp\alpha_\mfp s_\mfp = 1$. But this implies $f_k(x) = \sum\alpha_\mfp s_\mfp f_k(x) = 0$, and \ref{sol.POLALOC.ast} is proved.

Now assume $f_\mfp = 0$ as a polynomial law $M_\mfp \to N_\mfp$ over $k_\mfp$ for all $\mfp \in \Spec(k)$. Let $R \in \kalg$ and $\mfP \in \Spec(R)$. Then $\mfp := i^{-1}(\mfP) \in \Spec(k)$, where $i$ stands for the unit homomorphism $k \to R$, and we derive from (\ref{ss.LOPROM}.\ref{KAYPE}) the commutative diagram
\[
\xymatrix{
k \ar[r]^{i} \ar[d]_{\mathrm{can}} & R \ar[d]^{\mathrm{can}_R} \\
k_\mfp \ar[r]_{j} & R_\mfP}
\]
with a local homomorphism $j$ of local $k$-algebras, making $R_\mfP$ a $k_\mfp$-algebra compatible with its $k$-algebra structure. Since $f_\mfp = 0$ as a polynomial law $M_\mfp \to N_\mfp$ over $k_\mfp$ by hypothesis, we therefore conclude $(f_R)_{R_\mfP} = f_{R_\mfP} = (f_\mfp)_{R_\mfP} = 0$ as a set map $(M_R)_\mfP \to (N_R)_\mfP$. Now \ref{sol.POLALOC.ast} for $R$ in place of $k$ and $f_R$ in place of $f$ shows $f_R = 0$ as a set map $M_R \to N_R$. 

(b) We must show that the conjunction of (i) $f_k = 0$ and (ii) $D^pf = 0$ for $1 \le p \le \lfloor\frac{d}{2}\rfloor$ implies $f = 0$. For $d = 0,1$, condition (ii) is void, so we must show that $f_k = 0$ implies $f = 0$. But this is clear since $f$ for $d = 0$ is a constant polynomial law (Exc.~\ref{pr.COPOLA}), while for $d = 1$ it is basically a linear map (Exc.~\ref{pr.MULI}~(a)). We may therefore assume $d \ge 2$. In order to prove $f = 0$, it suffices to show, by Corollary~\ref{c.VANCRI}, that $(\Pi^\nu f)_k = 0$ for all $\nu \in \IN^n$, $n \in \IZ$, $n > 0$. By \ref{ss.LIPOLA}~(c),  we may assume $\vert\nu\vert = d$. For $n = 1$, we therefore have $\nu = d$, and (\ref{t.BILEXP}.\ref{EFERBT}) for $R = k$ and $n = 1$ yields $(\Pi^df)_k = f_k = 0$. For $n > 1$, we put
\[
m := \max\,\{i \in \IZ \mid 1 \le i \le n,\, \nu_i \ne 0\}.
\]
Hence $\nu = (\mu,0,\dots,0)$, $\mu = (\nu_1,\dots,\nu_m)$, and if $m < n$, then Exc.~\ref{pr.ZERLIN} implies $(\Pi^\nu f)_R(x,y) = (\Pi^\mu f)_R(x)$ for all $R \in \kalg$, $x \in M_R^m$, $y \in M_R^{n-m}$. We are thus reduced to the case $\nu_n \ne 0$. Simplifying notation, we conclude $\nu = (\varrho,p)$ for some $\varrho \in \IN^{n-1}$, $p \in \IN$ having $p \ge 1$ and $\vert\varrho\vert + p = n$. Applying \ref{ss.LIPOLA}~(b) with an appropriate permutation, we may assume $1 \le p \le \lfloor\frac{d}{2}\rfloor$. Now Lemma~\ref{l.TODEFI} combined with $D^pf = 0$ yields
\[
0 = (D^pf)_{k[\bfT]}(\sum_{i=1}^{n-1}\bft_ix_i,x_n) = \sum_{\lambda\in\IN^{n-1}}\bfT^\lambda (\Pi^{(\lambda,p)}f)_k(x_1,\dots,x_{n-1},x_n)
\]
for all $x_1,\dots,x_{n-1},x_n \in M$. Applying (\ref{ss.UFRBACH}.\ref{EMART}), we conclude $(\Pi^{(\lambda,p)}f)_k = 0$ for all $\lambda \in \IN^{n-1}$. Specializing $\lambda$ to $\varrho$, we obtain $(\Pi^\nu f)_k = 0$, as claimed.
\end{sol}

\begin{sol}{pr.UNI} \label{sol.UNI} Let $F \in \kalg$ be a field. By hypothesis, $f_F(m_F) = f_k(m)_F \in N$ is unimodular, hence not zero. Since $f_F(0) = f_k(0)_F = 0$, we conclude $m_F \ne 0$, whence $m \in M$ is unimodular (Lemma~\ref{l.FAITH}).
\end{sol}

\begin{sol}{pr.ANISO} \label{sol.ANISO} (a) Let $m \in M_R$ be nonzero, so $m = \sum_{j \ge j_0} m_j \bft^j$ for some $j_0 \ge 0$ with $m_{j_0} \ne 0$.  

Suppose that $j_0 = 0$.  Then the homomorphism $R \to k$ defined by sending $\bft \mapsto 0$ sends $m \mapsto m_0$ and $f_R(m) \mapsto f(m_0)$.  Since $f$ does not represent zero on $M$, $f(m_0) \ne 0$, whence $f_R(m) \ne 0$ and indeed $f_R(m)$ has lowest degree term $f(m_0)t^0$.

Suppose $j_0 > 0$.  Then
\[
f_R(m) =f_R(t^{j_0} (mt^{-j_0})) = t^{dj_0} (f(m_{j_0}) t^0 + \cdots),
\]
proving the claim.

(b)  Suppose there is a nonzero $m \in M$ such that $f_k(m) = 0$.  Because $M$ is torsion-free, $m \otimes 1$ is not zero in $M_R$ \cite[II.2, Prop.~4]{MR0360549} and $f_R(m \otimes 1) = f_k(m) \otimes 1_R = 0$.

Conversely, suppose $f_R$ represents zero.  Then by combining denominators, there is a nonzero $m \in M$ and $s \in S$ such that $f_R(m \otimes 1/s) = 0$.  Whence 
\[
0 = (1 \otimes s^{-d}) f_R(m \otimes 1) = (1 \otimes s^{-d}) (f_k(m) \otimes 1),
\]
and $f_k(m) \otimes 1 = 0$ in $R$.  Since $S$ contains no zero divisors, $f_k(m) = 0$.

(c): By hypothesis, there is a nonzero $m \in M$ such that $f_k(m) = 0$.  Identify $M = k^r$ for some $r \in \IN$.  The ideal of $k$ generated by the coordinates $m_i$ of $m$ for all $i $ is generated by a $g \ne 0$, and $g$ divides $m_i$ in $k$ for all $i$.  Defining $m' \in M$ via $gm'_i = m_i$ for each $i$, we find that $gm' = m$ and the coordinates of $m'$ generate the unit ideal in $k$ so $m'$ is unimodular.  Finally, 
$0 = f_k(m) = g^d f_k(m')$, whence $f_k(m') = 0$. 
\end{sol}

\begin{sol}{pr.DVR.FORMS} \label{sol.DVR.FORMS} The argument is similar to that of Exercise \ref{pr.ANISO}~(a).  Put $v$ for the value function on $R$, so in particular $v(\pi) = 1$.  If $M$ is a finitely generated free $R$-module, then the unimodular elements are exactly those $m \in M$ whose image $\bar{m} \in M \otimes k$ is not zero.  With respect to a basis of $M$, the unimodular elements are the ones with at least one coordinate in $R^\times$.  For any nonzero $m \in M$, we have $m = \pi^e m'$ for uniquely determined $e \ge 0$ and $m' \in M$ unimodular.

For $m \in M_i$, if $\bar{m} \in M_i \otimes k$ is not zero, then $\overline{f_i(m)} = \overline{f_i}(\bar{m}) \ne 0$, whence $v(f_i(m)) = 0$.  Otherwise, if $m \ne 0$, because $M_i$ is free and finitely generated it can be written as $m = \pi^e m'$ for some $e > 0$ and $m' \in M$ whose image in $M \otimes k$ is not zero.  Then $f_i(m) = \pi^{de} f_i(m')$, so $v(f_i(m)) = de$.  In summary, for all $m \in M_i$, $v(f_i(m))$ is divisible by $d$.  (If $m = 0$, then $f_i(m) = 0$ and $v(f_i(m)) = \infty$, which is divisible by $d$.)  Note that this argument shows that $f_i$ does not represent zero.

To prove the claim, pick $m_i \in M_i$ for all $i$ such that $m := \sum m_i$ is not zero.  By the preceding paragraph, $v(\pi^i f_i(m_i)) = i + v(f_i(m_i)) \equiv i \mod d$.  Since these represent distinct residue classes mod $d$, the one of smallest valuation, call it $i_0$, is the only one with that valuation.  It follows that $v(m) = i_0 + v(f_{i_0}(m_{i_0}))$, the value of the smallest term.  In particular, $\sum_i \pi^i f_i(m_i) \ne 0$.
\end{sol}

\begin{sol}{pr.EXOCUB} \label{sol.EXOCUB} By (\ref{ss.TODEHO}.\ref{DEEFPI}) we have $D^2f = \Pi^{(1,2)}f$, and (\ref{ss.LIPOLA}.\ref{PERLIN}) implies $(D^2f)(x,y) = \Pi^{(1,2)}f)(x,y) = (\Pi^{(2,1)}f)(y,x) = (Df)(y,x)$, hence
\begin{align}
\label{EFYEX} f(y,x) = (D^2f)(x,y).
\end{align}
Combining (\ref{ss.TODEHO}.\ref{HOSYMDIF}) for $d = 3$ with (\ref{ss.TODE}.\ref{TAYL}) and \eqref{EFEXY}, \eqref{EFYEX}, we obtain \eqref{EFEXTE}. In order to derive \eqref{EFEXPLY}, we first let $x,y \in M$. Then \eqref{EFYEX}, Exc.~\ref{pr.HIDIDE}~(c) and Prop.~\ref{p.ITDIDE} imply
\begin{align*}
f(x + y,z) =\,\,&(D^2f)(z,x + y) = (\partial_{x+y}^{[2]}f)(z) \\
=\,\,&(\partial_x^{[2]}f)(z) + (\partial_x\partial_yf)(z) + (\partial_y^{[2]}f)(z) \\
=\,\,&(D^2f)(z,x) + (\Pi^{(1,1,1)}f)(x,y,z) + (D^2f)(z,y) \\
=\,\,&f(x,z) + f(x,y,z) + f(y,z), 
\end{align*}
as claimed. The general case $x,y \in M_R$ now follows by looking at $f \otimes R$ rather than $f$. Combining \eqref{EFEXTE} for $\bft = 1$ and \eqref{EFEXPLY}, we now obtain
\begin{align*}
f(x + y + z) \,\,&- f(x + y) - f(y + z) - f(z + x) + f(x) + f(y) + f(z) \\
=\,\,&f(x + y) + f(x + y,z) + f(z,x + y) + f(z) - f(x + y) - f(y) - f(y,z) \\
\,\,&- f(z,y) - f(z) -f(z) - f(z,x) - f(x,z) - f(x)  + f(x) + f(y) + f(z) \\
=\,\,&f(x,z) + f(x,y,z) + f(y,z) + f(z,x) + f(z,y) + f(z) - f(y) - f(y,z) -f(z,y) \\
\,\,&- f(z) - f(z) - f(z,x) - f(x,z) - f(x) + f(x) + f(y) + f(z) \\
=\,\,&f(x,y,z),
\end{align*}
and this is \eqref{EFEXYZ}. Before dealing with \eqref{EFLICO}, we prove
\begin{align}
\label{EFRIXI} f(\sum_{i=1}^nr_ix_i,x) = \sum_{i=1}^nr_i^2f(x_i,x) + \sum_{1\leq i<j\leq n}r_ir_jf(x_i,x_j,x)
\end{align}
for $x \in M_R$ by induction on $n$. For $n = 1$ there is nothing to prove. For $n > 1$, the induction hypothesis and \eqref{EFEXPLY} imply
\begin{align*}
f(\sum_{i=1}^nr_ix_i,x) =\,\,&f(\sum_{i=1}^{n-1}r_ix_i,x) + f(\sum_{i=1}^{n-1}r_ix_i,r_nx_n,x) + f(r_nx_n,x) \\
=\,\,&\sum_{i=1}^{n-1}r_i^2f(x_i,x) + \sum_{1\leq i<j<n}r_ir_jf(x_i,x_j,x) + \sum_{i=1}^{n-1}r_ir_nf(x_i,x_n,x) + r_n^2f(x_n,x) \\
=\,\,&\sum_{i=1}^nr_i^2f(x_i,x) + \sum_{1\leq i<j\leq n}r_ir_jf(x_i,x_j,x),
\end{align*}
as claimed. Now we can prove \eqref{EFLICO}, again by induction on $n$. The case $n = 1$ again being obvious, let us assume $n > 1$ and that \eqref{EFLICO} holds for $n - 1$. Then this and \eqref{EFRIXI} yield
\begin{align*}
f(\sum_{i=1}^nr_ix_i) =\,\,&f(\sum_{i=1}^{n-1}r_ix_i) + f(\sum_{i=1}^{n-1}r_ix_i,r_nx_n) + f(r_nx_n,\sum_{i=1}^{n-1}r_ix_i) + f(r_nx_n) \\
=\,\,&\sum_{i=1}^{n-1}r_i^3f(x_i) + \sum_{1\leq i,j<n,i\neq j}r_i^2r_jf(x_i,x_j) + \sum_{1\leq i<j<l<n}r_ir_jr_lf(x_i,x_j,x_l) \\
\,\,&+ \sum_{i=1}^{n-1}r_i^2r_nf(x_i,x_n) + \sum_{1\leq i<j<n}r_ir_jr_nf(x_i,x_j,x_n) + \sum_{j=1}^{n-1}r_n^2r_jf(x_n,x_j) + r_n^3f(x_n) \\
=\,\,&\sum_{i=1}^nr_i^3f(x_i) + \sum_{1\leq i,j\leq n, i\neq j}r_i^2r_jf(x_i,x_j) + \sum_{1\leq i<j<l\leq n}r_ir_jr_lf(x_i,x_j,x_l).
\end{align*}
This completes the induction and the proof of (a).

(b) (i) $\Rightarrow$ (ii). From (i) we deduce $g_F(e_i) = 0$ for $1 \leq i \leq n$, while \eqref{EFEXYZ} implies $g_F(x,y,z) = 0$ for all $x,y,z \in F^n$. Moreover, from Euler's differential equation (\ref{ss.TODEHO}.\ref{EULDIFEQ} we conclude $g_F(x,x) = (Dg)_F(x,x) = 3g_F(x) = 0$. It remains to show that $F$ consists of two elements and $g_F(x_0,y_0) \neq 0$ for some $x_0,y_0 \in F^n$. Replacing $f$ by $g$ and specializing $\bft \mapsto \alpha \in F^\times$ in \eqref{EFEXTE}, we conclude $g_F(x,y) + \alpha g_F(y,x) = 0$. Assuming $\vert F\vert > 2$, this implies $g_F(x,y) = 0$ for all $x,y \in F^n$. But If $g_F(x,y) = 0$ for all $x,y \in F^n$, then \eqref{BASGEF} and \eqref{EFLICO} for $g$ in place of $f$ and $e_{iR}$ in place of $x_i$ for $1 \leq i \leq n$ would lead to the contradiction $g = 0$ as a polynomial law over $F$. Thus (ii) holds.

(ii) $\Rightarrow$ (iii). By \eqref{BASGEF} and \eqref{EFRIXI} for $g$ in place of $f$, the map $F^n \times F^n \to F$, $(x,y) \mapsto g_F(x,y)$ is an alternating $F$-bilinear form, forcing
\[
S := \big(g_F(e_i,e_j)\big)_{1\leq i,j\leq n} \in \Mat_n(F)
\] 
to be an alternating matrix and $S \neq 0$ by (ii). Now let
\begin{align*}
x = \left(\begin{matrix}
r_1 \\
\vdots \\
r_n
\end{matrix}\right) = \sum_{i=1}^ nr_ie_{iR} \in R^ n = (F^ n)_R &&(r_1,\dots,r_n \in R).
\end{align*}
Then \eqref{EFLICO} and \eqref{BASGEF} imply
\begin{align*}
g_R(x) =\,\,&g_R(\sum_{i=1}^nr_ie_{iR}) = \sum_{i=1}^nr_i^ 3g_F(e_i) + \sum_{1\leq i,j\leq n, i\neq j}r_i^ 2r_jg_F(e_i,e_j) + \sum_{1\leq i<j<l\leq n}r_ir_jr_lg_F(e_i,e_j,e_l) \\
=\,\,&\sum_{1\leq i,j\leq n}r_ig_F(e_i,e_j)r_j^2 = x^\trans Sx^2,   
\end{align*}
as claimed.

(iii) $\Rightarrow$ (i). Since $F \cong \IF_2$ consists of two elements, the elements of $F^ n$ as an $F$-algebra are all idempotents. Hence $g_F(x) = x^\trans Sx = 0$ for all $x \in F^n$, whence the set map $g_F\:F^n \to F$ is identically zero. On the other hand, since $S \neq 0$, there are elements $x_0,y_0 \in F^ n$ such that $x_0^\trans Sy_0 \neq 0$, and passing from $F \cong \IF^2$ to the separable quadratic extension $K = F(\theta) \cong \IF_4$ with $\theta \in K$ satisfying $\theta^2 = \theta + 1$, we deduce
\begin{align*}
g_K(x_0 + \theta y_0) =\,\,&(x_0 + \theta y_0)^\trans S(x_0^2 + \theta^2y_0^2) \\
=\,\,&(x_0 + \theta y_0)^\trans S(x_0 + \theta^2y_0) = (x_0 + \theta y_0)^\trans S(x_0 + \theta y_0) + (x_0 + \theta y_0)^\trans Sy_0 \\
=\,\,&x_0^\trans Sy_0 \neq 0.
\end{align*}
Thus $g_K \neq 0$, so $g\:F^n \to F$ is a non-zero cubic form.
\end{sol}

\begin{sol}{pr.SPRINGER.CUBIC} \label{sol.SPRINGER.CUBIC} Identifying $V \subseteq V_K$ canonically and picking a basis $(1,\xi)$ of $K$ as a vector space over $F$, we have $V_K = V \oplus \xi V$ as a direct sum of $F$-subspaces. Arguing indirectly, we assume that $f \otimes K$ does represent zero. Then $f_K(v + \xi w) = 0$ for some $v,w \in V$ not both zero. If $w = 0$, then $f_F(v) = 0$, hence $v = 0$ (since $f$ does not represenr zero), a contradiction. Thus $w \ne 0$ and then $f_F(w) \ne 0$. Consulting eqn. \eqref{EFEXTE} of Exc.~\ref{pr.EXOCUB}, we conclude that $\mu := f_{F[\bft]}(v + \bft w) \in F[\bft]$ is a cubic polynomial having the root $\xi$ in $K$ but no root  in $F$ (since $f_F$ is does not represent zero). Thus $\mu$ is irreducible and hence the minimum polynomial of $\xi$. On the other hand, the minimum polynomial of any element in the quadratic field extension $K$ of $F$ has degree at most $2$. This contradiction proves our claim.
\end{sol}

\begin{sol}{pr.THICHA} \label{sol.THICHA} Passing to $R[\vep]$, $\vep^4 = 0$, the Taylor expansion (\ref{ss.TODE}.\ref{TAYLEPS}) yields
\begin{align*}
\big(D^0(g \circ f)\big)(x,y) +\,\,& \vep\big(D(g \circ f)\big)(x,y) + \vep^2\big(D^2(g \circ f)\big)(x,y) + \vep^3\big(D^3(g \circ f)\big)(x,y) \\
=\,\,&(g \circ f)(x + \vep y) \\
=\,\,&g\big(f(x) + \vep(Df)(x,y) + \vep^2(D^2f)(x,y) + \vep^3(D^3f)(x,y)\big) \\
=\,\,&g\big(f(x)\big) + \vep(Dg)\big(f(x),(Df)(x,y) + \vep(D^2f)(x,y) + \vep^2(D^3f)(x,y)\big) \\	
\,\,&+ \vep^2(D^2g)\big(f(x),(Df)(x,y) + \vep(D^2f)(x,y) + \vep^2(D^3f)(x,y)\big) \\
\,\,&+ \vep^3(D^3g)\big(f(x),(Df)(x,y) + \vep(D^2f)(x,y) + \vep^2(D^3f)(x,y)\big).
\end{align*}
Comparing coefficients of $\vep^3$ by combining \eqref{BILSEC} with \eqref{EFLICO} of Exercise~\ref{pr.EXOCUB}, the assertion follows.
\end{sol}

\begin{sol}{pr.CUMA} \label{sol.CUMA} \ref{CUMA.a} If $f\:M \to N$ and $f^\prime\:M^\prime \to N^\prime$ are polynomial laws over $k$, a \emph{morphism} from $f$ to $f^\prime$ in $\kpolaw$ is defined as a pair of linear maps $\vph\:M \to M^\prime$, $\psi\:N \to N^\prime$ making a commutative diagram
\begin{align}
\vcenter{\label{MORPOL}\xymatrix{
M \ar[r]^{f} \ar[d]_{\vph} & N \ar[d]^{\psi} \\
M^\prime \ar[r]_{f^\prime} & N^\prime}}
\end{align}
of polynomial laws over $k$. 

\ref{CUMA.c} Since $g$ is quadratic in the first variable, the expression $g(x,y,z)$ is surely trilinear. Moreover, for $x,y,z \in M$, we repeatedly use condition (iii) to compute
\begin{align*}
f(x + y + z) -\,\,& f(x + y) - f(y + z) - f(z + x) + f(x) + f(y) + f(z) \\
=\,\,&g(x + y,z) + g(z,x + y) + f(z) -f(y) - g(y,z) - g(z.y) -f(z) \\
\,\,&- f(z) -g(z,x) - g(x,z) - f(x) + f(x) +f(y) + f(z) \\
=\,\,&g(x,y,z),
\end{align*}
and since the very left-hand side is totally symmetric in $x,y,z$, so is the right. This proves the first part of (a). As to the second, \eqref{TRIRAD} follows immediately by linearizing $g(y,x) = 0$ for $x \in \Rad(f,g)$ with respect to $y \in M$ and using the total symmetry of \eqref{TRICU}. Hence (iii) implies that $\Rad(f,g)$ is a submodule of $M$. Finally, for a linear map $\pi$ as indicated, the maps $f_1\:M _1 \to N$ and $g_1\:M_1 \times M_1 \to N$ given by $f_1(\pi(x)) := f(x)$ and $g_1(\pi(x),\pi(y)) := g(x,y)$, respectively, for $x,y \in M$ in view of \eqref{CURAD} are well-defined, and it is straightforward to check that $(f_1,g_1)\:M_1 \to N_1$ is a cubic map with the desired properties. This completes the proof of (c). 

\ref{CUMA.d} Combining (iii) with the fact that $g$ is quadratic-linear, \eqref{CUEXP} follows by a straightforward induction. \eqref{CUEXP} and the commutativity of \eqref{ARCUB} immediately imply uniqueness of the $R$-cubic extension. To prove its existence, let $R \in \kalg$ and assume first that $M$ is a free $k$-module, with basis $(e_i)_{i\in I}$. Write $g_R\:M_R \times M_R \to N_R$ for the $R$-quadratic-linear extension of $g$ in the sense of Exercise~\ref{pr.QUALI}. In the spirit of \eqref{CUEXP}, we define a set map $f_R\:M_R \to N_R$ by
\begin{align}
\label{CUMA.EFAREX} f_R(\sum_{i\in I} r_ie_{iR}) :=\,\,& \sum_{i} r_i^3f(e_i)_R + \sum_{i\neq j} r_i^2r_j g(e_i,e_j)_R + \sum_{|i<j<l} r_ir_jr_lg(e_i,e_j,e_l)_R
\end{align}
for all families $(r_i)_{i\in I}$ with finite support in $R$. Putting $x := \sum r_ie_{iR}$, $y := \sum s_ie_{iR}$ with another such family $(s_i)_{i\in I}$, we obtain
\begin{align}
\label{GEAREX} g_R(x,y) = \sum_{i,l} r_i^2s_lg(e_i,e_l)_R + \sum_{i<j,l} r_ir_jr_lg(e_i,e_j,e_l)_R,
\end{align}and using \eqref{CUMA.EFAREX}, \eqref{GEAREX}, (iv), one checks that $(f,g)_R := (f_R,g_R)$ is indeed an $R$-cubic extension of $(f,g)$.

Now assume that $M$ is arbitrary. Since the functor $\emptyslot \otimes R$ is right exact, we obtaina commutative diagram
\begin{align}
\vcenter{\label{ELPEM} \xymatrix{0 \ar[r] & L \ar[r]_{i} \ar[d]_{\can} & \pM \ar[r]_{\pi} \ar[d]^{\can} & M \ar[r] \ar[d]^{\can} & 0 \\
& L_R \ar[r]_{i_R} & \pM_R \ar[r]_{\pi_R} & M_R \ar[r] & 0}}
\end{align}
of $k$- (resp. $R$-)modules, with exact rows and $M$ a free $k$-module. By the special case just treated, the cubic map
\begin{align}
\label{PEFPEG} (\pf,\pg) := \big(f \circ \pi,g \circ(\pi \times \pi)\big)\:\pM \longrightarrow N
\end{align}
over $k$ extends to a cubic map
\begin{align}
\label{PEFPEGR} (\pf,\pg)_R := (\pf_R,\pg_R)\:\pM_R \longrightarrow N_R
\end{align}
over $R$. We claim
\begin{align}
\label{KERPIR} \Ker(\pi_R) \subseteq \Rad\big((\pf,\pg)_R\big).
\end{align}
Exactness of the second row in \eqref{ELPEM} implies that any element of $\Ker(\pi_R)$ can be written as $i_R(z)$, for some $z = \sum r_ju_{jR} \in L_R$, $r_j \in R$, $u_j \in L$. Applying \eqref{CUEXP} to $(\pf,\pg)_R$, we conclude
\begin{align*}
\pf_R\big(i_R(z)\big) =\,\,&\pf_R\big(\sum r_ji(u_j)_R\big) \\
=\,\,& \sum_jr_j^3\pf\big(i(u_j)\big)_R + \sum_{j\neq l}r_j^2r_l\pg\big(i(u_j),i(u_l)\big)_R + \sum_{j<l<m}r_jr_lr_m\pg\big(i(u_j),i(u_l),i(u_m)\big)_R,
\end{align*}
where \eqref{PEFPEGR} implies $\pf(i(u_j)) = (f \circ \pi \circ i)(u_j) = 0$, $\pg(i(u_j),i(u_l)) = g((\pi \circ i)(u_j)$, $(\pi \circ i)(u_l)) = 0$ and, similarly, $\pg(i(u_j),i(u_l),i(u_m)) = 0$. Since $\pg$ is $R$-quadratic-linear, we also have $\pg_R(i_R(z),y) = \pg_R(y,i_R(z)) = 0$ for all $y \in M$. Summing up, we conclude $i_R(z) \in \Rad((\pf,\pg)_R)$, and \eqref{KERPIR} follows. By (c), therefore, we find a cubic map
\begin{align}
\label{EFGER} (f,g)_R := (f_R,g_R)\:M_R \longrightarrow N
\end{align}
over $R$ such that
\begin{align}
\label{FARPIAR} f_R \circ \pi_R = \pf_R, \quad g_R \circ (\pi_R \times \pi_R) = \pg_R.
\end{align}
Since $\pi$ and $\pi_R$ are both surjective, one checks that \eqref{ARCUB} commutes, so $(f,g)_R$ is indeed an $R$-cubic extension of $(f,g)$.

\ref{CUMA.e} If $(f,g)\:M \to N$ and $(f^\prime,g^\prime)\:M^\prime \to N^\prime$ are cubic maps over $k$, a \emph{morphism} from $(f,g)$ to $(f^\prime,g^\prime)$ in $\kcumap$ is defined as a pair of linear maps $\vph\:M \to M^\prime$, $\psi\:N \to N^\prime$ making commutative diagrams 
\begin{align}
\vcenter{\label{MORCUMAP}\xymatrix{
M \ar[r]^{f} \ar[d]_{\vph} & N \ar[d]^{\psi} && M \times M \ar[r]^{g} \ar[d]_{\vph \times \vph} & N \ar[d]^{\psi} \\
M^\prime \ar[r]_{f^\prime} & N^\prime, && M^\prime \times M^\prime \ar[r]_{g^\prime} & N^\prime}}
\end{align}
of set maps.

Now let $f\:M \to N$ be a homogeneous polynomial law of degree $3$ over $k$. We combine equations \eqref{EFEXTE}, \eqref{EFEXPLY} of Exc.~\ref{pr.EXOCUB} with Euler's differential equation (\ref{ss.TODEHO}.\ref{EULDIFEQ}) to conclude that $(f_k,(Df)_k)$ is a cubic map from $M$ to $N$. If $(\vph,\psi)\:f \to f^\prime$ is a morphism in $\kholaw$, then the chain rules (\ref{ss.DICA}.\ref{DICH}), (\ref{ss.DICA}.\ref{DICL}) applied to the commutative diagram \eqref{MORPOL} show that $(\vph,\psi)$ is also a morphism from $(f_k,(Df)_k)$ to $(f^\prime_k,(Df^\prime)_k)$ in $\kcumap$. Thus we have obtained a functor from $\kholaw$ to $\kcumap$. 

Conversely, let $(f,g)\:M \to N$ be a cubic map over $k$. We claim that \emph{the family of set maps
\begin{align}
\label{EFSTAR} (f \ast g)_R := f_R\:M_R \longrightarrow N_R
\end{align}
defined for $R \in \kalg$ as the first component of the $R$-cubic extension of $(f,g)$ is a polynomial law over $k$.} Indeed, given a morphism$\vph\:R \to S$ in $\kalg$, we replace $f$ by $f_R$ and $R$ by $S$ in \eqref{ARCUB}, observe $(f_R)_S = f_S$ and conclude from (\ref{ss.ITSCA}.\ref{TENTH}) that (\ref{ss.COPOLA}.\ref{FUPOLA}) commutes, as desired. Note that the polynomial law $f \ast g\:M \to N$ by definition is homogeneous of degree $3$, hence an object of the category $\kholaw$. If $(\vph,\psi)\:(f,g) \to (\pf,\pg)$ is a morphism of cubic maps $(f,g)\:M \to N$, $(\pf,\pg)\:\pM \to \pN$ over $k$, one checks that $\pf_R(\vph(x_R)) = \psi_R(f_R(x_R))$, $\pg_R(\vph_R(x_R),\vph(y_R)) = \psi_R(g_R(x_R,y_R))$ for all $R \in \kalg$, $x,y \in M$, and we conclude from \eqref{EFSTAR} that $(\vph_R,\psi_R)\:(f,g)_R \to (\pf,\pg)_R$ is a morphism of cubic maps over $R$. Hence $(\vph,\psi)\:f \ast g \to \pf \ast \pg$ is a morphism in $\kholaw$., thus defining a functor from $\kcumap$ to $\kholaw$. It remains to show that the two functors are inverse to one another. Let $f\:M \to N$ be a homogeneous poynomial law of degree $3$ over $k$. For $R \in \kalg$, the diagrams
\[
\xymatrix{M \ar[r]_{f_k} \ar[d]_{\can} & N \ar[d]^{\can} & M \times M \ar[r]_(.6){(Df)_k} \ar[d]_{\can} & N \ar[d]^{\can} \\
M_R \ar[r]_{f_R} & N_R, & M_R \times M_R \ar[r]_(.6){(Df)_R} & N_R}
\]
commute, and we conclude $f_R = (f_k)_R$, $(Df)_R = ((Df)_k)_R$, so $(f_R,(Df)_R)$ is the $R$-cubic extension of $(f_k,(Df)_k)$. Hence $f_k \ast (Df)_k = f$ as polynomial laws over $k$. Conversely, let $(f,g)\:M \to N$ be a cubic map over $k$ and write $R := k[\vep]$, $\vep^2 = 0$, for the $k$-algebra of dual numbers. Since $g$ is quadratic-linear, condition (iii) in the definition of a cubic map implies
\[
f_R(x_1 + \vep x_2) = f(x_1 + \vep g(x_1,x_2))
\]
for all $x_1,x_2 \in M$, and comparing with the Taylor expansion (\ref{ss.TODE}.\ref{TAYLDUAL}) we conclude $((f \ast g)_k,(D(f \ast g))_k) = (f,g)$. This completes the proof.
\end{sol}

\begin{sol}{pr.CUTRIL} \label{sol.CUTRIL} Arguing  as in the solution to Exercise~\ref{pr.QUABIL}, we first assume that $M$ is free, let $(e_i)_{i \in I}$ be a basis of $M$, and choose a total order on $I$.
By Exercise~\ref{pr.CUMA}~\ref{CUMA.e}, we may write $F = f \ast g$, where $(f,g)\:M \to N$ is a cubic map in the sense of Exercise~\ref{pr.CUMA}. Let $T\:M \times M \times M \to N$ be the unique trilinear map from $M \times M \times M$ to $N$ given by the values
\[
T(e_i,e_j,e_l) =
\begin{cases}
f(e_i) & \text{for}\;\; i = j = l \in I, \\
g(e_i,e_l) & \text{for}\;\; i,j,l \in I, i = j \neq l, \\
g(e_i,e_j,e_l) & \text{for}\;\; i,j,l \in I, i < j < l, \\
0 & \text{otherwise}
\end{cases}
\]
on the basis vectors. Then, for any $x = \sum\xi_ie_i \in M$, $\xi_i \in k$, we have
\begin{align*}
T(x,x,x) =\,\,&T(\sum \xi_ie_i,\sum \xi_ie_i,\sum \xi_ie_i) = \sum_{i,j,l\in I} \xi_i\xi_j\xi_lT(e_i,e_j,e_l) \\
=\,\,&\sum_{i\in I}\xi_i^3T(e_i,e_i,e_i) + \sum_{i,l\in I,i \neq l}\xi_i^2\xi_l\big(T(e_i,e_i,e_l) + T(e_i,e_l,e_i) + T(e_l,e_i,e_i)\big) \\
\,\,&+ \sum_{i,j,l\in I,i<j<l}\xi_i\xi_j\xi_l\big(T(e_i,e_j,e_l) + T(e_j,e_l,e_i) + T(e_l,e_i,e_j) \\
\,\,&+ T(e_j,e_i,e_l) + T(e_i,e_l,e_j) + T(e_l,e_j,e_i)\big) \\
=\,\,&\sum \xi_i^3f(e_i) + \sum_{i\neq j} g(e_i,e_j) + \sum_{i<j<l}\xi_i\xi_j\xi_lg(e_i,e_j,e_l) = f(\sum\xi_ie_i)
\end{align*}
by Exercise~(\ref{pr.CUMA}.\ref{CUEXP}) Since the identification $F = f \ast g$ is compatible with base change, the assertion follows if $M$ is free. Now let $M$ be an arbitrary projective module. Then there exists a $k$-module $M^\prime$ making $M_0 := M \oplus M^\prime$ free. Let $F^\prime\:M^\prime \to N$ be the zero polynomial law on $M^\prime$ with values in $N$ and consider the homogeneous polynomial law $F_0\:M_0 \to N$ of degree $3$ canonically induced by $F$ and $\pF$. By the special case just treated we find a trilinear map $T_0\:M_0 \times M_0 \times M_0 \to N$ such that $F_0(x_0) = T_0(x_0,x_0,x_0)$ for all $x_0 \in M_{0R}$, $R \in \kalg$. But this implies $F(x) = T(x,x,x)$ for all $x \in M_R$, where $T\:M \times M \times M \to N$ is the restriction of $T_0$ to $M \times M \times M$. 
\end{sol}

\begin{sol}{pr.RESCAIN} \label{sol.RESCAIN} Let $f\:M \to N$ be a polynomial law over $K$. Letting $S \in \Kalg$ and
setting $R := \,_kS$, the map $\Eins_S\:R \to S$ is a morphism of
$k$-algebras, giving rise to a \emph{surjective} morphism
$\varphi\:R_K \to S$ of $K$-algebras and hence, in view of
(\ref{ss.RESC}.\ref{EMAR}), to a commutative diagram with exact columns
\[
\xymatrix{ (_kM)_R \ar[r]_{(_kf)_R}\ar[d]_{\Eins_M \otimes \varphi} & (_kN)_R \ar[d]^{\Eins_N \otimes \varphi} \\ 
M_S \ar[r]_{f_S}\ar[d] & N_S\ar[d] \\ 
0 & 0, }
\]
which shows that $f_S$ is uniquely determined by $(_kf)_R$, so $f$ is uniquely determined by $_kf$, as claimed.
\end{sol}

\begin{sol}{pr.RESCACAL} \label{sol.RESCACAL} Let $R \in \kalg$.
	
(a) Using (\ref{ss.RESC}.\ref{KAFAR}), we obtain $(g \circ f)_R
= (g \circ f)_{R_K} = (g_{R_K}) \circ (f_{R_K}) = (_kg)_R \circ
(_kf)_R = (_kg \circ\,_kf)_R$.
	
(b) For $r,r_i \in R$, $x_i \in M$, we recall $(_kN)_R = N_{R_K}$ as
$R_K$-modules and compute $(_kf)_R(r\sum x_i \otimes r_i) =
(_kf)_R(\sum x_i \otimes rr_i) = f_{R_K}(\sum x_i \otimes_K (rr_i
\otimes 1_K)) = f_{R_K}(\sum x_i \otimes_K ((r \otimes 1_K)(r_i	\otimes 1_K))) = f_{R_K}((r \otimes 1_K)\sum x_i \otimes_K (r_i
\otimes 1_K)) = (r \otimes 1_K)^df_{R_K}(\sum x_i \otimes_K (r_i
\otimes 1_K)) = r^d(_kf)_R(\sum x_i \otimes r_i)$, where the last
equality is justified by (\ref{ss.RESC}.\ref{KAFAR}) for $a = 1_K$.
	
(c) For a finite chain $\bfT= (\bft_1,\dots,\bft_n)$ of independent indeterminates, we make use of the identification
\begin{align}
\label{ERTKA} R[\bfT]_K = R_K[\bfT]
\end{align}
via
\begin{align}
\label{ELERTKA} (r\bfT^\nu) \otimes a = (r \otimes a)\bfT^\nu, \quad \bfT^\nu \otimes 1_K = \bfT^\nu &&(r \in R,\;a \in K,\;\nu \in \IN^n)
\end{align}
Combining this identification with \ref{ss.RESC}, we find
\begin{align}
\label{KAEMTE} ((_kM)_R)_{R[\bfT]} = (_kM)_{R[\bfT]} = M_{R[\bfT]_K}
= M_{R_K[\bfT]} = (M_{R_K})_{R_K[\bfT]}
\end{align}
and claim
\begin{align}
\label{ZETENU}z \otimes_R \bfT^\nu = z \otimes_{R_K} \bfT^\nu &&(z
\in (_kM)_R = M_{R_K}, \nu \in \IN_0^n).
\end{align}
In order to prove \eqref{ZETENU}, additivity in the first variable
allows us to assume $z = x \otimes r = x \otimes_K (r \otimes 1_K)$
with $x \in M$, $r \in R$, and following the chain of
identifications presented in \eqref{KAEMTE}, we conclude $z
\otimes_R \bfT^\nu = (x \otimes r) \otimes_R \bfT^\nu = x
\otimes(r\bfT^\nu) = x \otimes_K ((r\bfT^\nu) \otimes 1_K) = x
\otimes_K((r \otimes 1_K)\bfT^\nu) =(x \otimes_K (r \otimes 1_K))
\otimes_{R_K} \bfT^\nu = z \otimes_{R_K} \bfT^\nu$, and
\eqref{ZETENU} follows. Now let $x_1,\ldots,x_n \in (_kM)_R =
M_{R_K}$. Then we combine \eqref{KAEMTE} with the identifications of (\ref{ss.COMUIN}.\ref{EMKATE}) and Theorem~\ref{t.BILEXP} to obtain
\begin{align*}
\sum_{\nu \in \IN_0^n}
\Big(\big(\Pi^\nu\,_kf\big)_R(x_1,\ldots,x_n)\Big) \otimes_R
\bfT^\nu =\,\,&(_kf)_{R[\bfT]}(\sum_{j=1}^n x_j \otimes_R \bft_j) \\
=\,\,&f_{R_K[\bfT]}(\sum_{j=1}^n x_j \otimes_{R_K} \bft_j) \\
=\,\,&\sum_{\nu \in \IN_0^n} (\Pi^\nu\,f)_{R_K}(x_1,\ldots,x_n)
\otimes_{R_K} \bfT^\nu \\
=\,\,&\sum_{\nu \in \IN_0^n}
\Big(\big(_k(\Pi^\nu\,f)\big)_R(x_1,\ldots,x_n)\Big) \otimes_R
\bfT^\nu,
\end{align*}
and a comparison yields $(\Pi^\nu\,(_kf))_R(x_1,\ldots,x_n) =
(_k(\Pi^\nu f))_R(x_1,\ldots,x_n)$, hence (c).

(d) Let $x \in (_kM)_R = M_{R_K}$. Then $(_kf_i)_R(x) =
(f_i)_{R_K}(x) = 0$ for almost all $i \in I$, showing that the
family $(_kf_i)_{i \in I}$ is locally finite. Moreover, $(\sum_{i
\in I} (_kf_i)_R)(x) = \sum_{i \in I} (_kf_i)_R(x) = \sum_{i \in I}
(f_i)_{R_K}(x) = (\sum_{i \in I} f_i)_{R_K}(x) = (_k(\sum_{i \in I}
f_i))_R(x)$, and (d) follows.

(e) Combining (c) with (d) and (\ref{ss.TODE}.\ref{ENDER}) yields
\[
 _k(D^n f) =\,_k(\sum_{p\geq 0} \Pi^{(p,n)} f) = \sum_{p\geq	0}\,_k(\Pi^{(p,n)} f) = \sum_{p\geq 0}\Pi^{(p,n)}(_kf) = D^n(_kf)
\]
for all $n \in \IN^0$.

(f)  For $R \in \kalg$, $x \in (_kM)_R = M_{R_K}$,
we apply (e) and obtain
\begin{align*}
(\partial_y (_kf))_R(x) =\,\,&(D(_kf))_R(x,y_R) = (_k(Df))_R(x,y_R) \\
=\,\,&(Df)_{R_K}(x,y_{R_K}) = (\partial_y f)_{R_K}(x) =
(_k(\partial_y f))_R(x),
\end{align*}
as claimed. 
\end{sol}

\begin{sol}{pr.LADJ} \label{sol.LADJ} \emph{Existence}. Let $R \in \Kalg \subseteq \kalg$. There is a natural surjection
\begin{align}
\label{OMER} \omega_R\:(_kN)_R \longrightarrow N_R, \quad \omega_R(y
\otimes r) = y \otimes_K r &&(y \in N, r \in R).
\end{align}
Combining this with (\ref{ss.ITSCA}.\ref{ITBA}), we obtain a set map
\begin{align}
\label{GER} \xymatrix{g_R := \omega_R \circ f_R\:(M_K)_R = M_R
\ar[r]_(.7){f_R} & (_kN)_R \ar[r]_{\omega_R} & N_R.}
\end{align}
We claim that the totality of set maps $g_R$, $R \in \Kalg$, is a polynomial law over $K$. To see this, suppose $\varphi\:R \to S$ is a morphism in $\Kalg$. In the diagram
\[
\xymatrix{(M_K)_R \ar[r]_{\Eins} \ar[d]_{\Eins_{M_K} \otimes_K \varphi} & M_R \ar[d]^{\Eins_M \otimes \varphi} \ar[r]_{f_R} &
(_kN)_R \ar[d]^{\Eins_N \otimes \varphi} \ar[r]_{\omega_R} & N_R
\ar[d]^{\Eins_N \otimes_K \varphi} \\ 
(M_K)_S \ar[r]_{\Eins} & M_S \ar[r]_{f_S} & (_kN)_S
\ar[r]_{\omega_S} & N_S,}
\]
the left-hand square commutes by (\ref{ss.ITSCA}.\ref{ITBA}), the
middle square commutes because $f$ is a polynomial law over $k$,
while the right-hand square commutes by virtue of \eqref{OMER}.	Combining we obtain a commutative diagram
\[
\xymatrix{(M_K)_R \ar[r]_{g_R} \ar[d]_{\Eins_{M_K} \otimes_K \varphi} & N_R \ar[d]^{\Eins_N \otimes_K \varphi} \\ 
(M_K)_S \ar[r]_{g_S} & N_S,}	
\]
which shows that $g$ is indeed a polynomial law over $K$. To
complete the existence part of the proof, it remains to show that
the diagram \eqref{FIEGE} commutes, so let $R \in \kalg$. Then
$(_kN)_R = N_{R_K}$ by (\ref{ss.RESC}.\ref{EMAR}), leading to a
composition of additive maps as
\[
\xymatrix{ (_kN)_R \ar[rr]_{\Eins_{(_kN)} \otimes \can_{R,K}} &&
(_kN)_{R_K} \ar[rr]_{\omega_{R_K}} && N_{R_K} = (_kN)_R,}
\]
where $\can_{R,K}\:R \to R_K$ is the natural map of (\ref{ss.SCEM}.\ref{KTOR}). For $y \in N$, $r \in R$, we obtain $\omega_{R_K} \circ (\Eins_{(_kN)} \otimes \can_{R,K})(x \otimes r) = \omega_{R_K}(x \otimes (r \otimes 1_K)) = x \otimes_K (r \otimes 1_K) = x \otimes r$, so
\begin{align}
\label{OMEINS} \omega_{R_K} \circ (\Eins_{(_kN)} \otimes \can_{R,K}) =
\Eins_{(_kN)_R}.
\end{align}
We must show $(_kg)_R \circ (_k(\can_{M,K}))_R = f_R$. Under the natural identification $(_k(M_K))_R = (M_K)_{R_K} = M_{R_K}$, we have $(_k(\can_{M,K}))_R = \Eins_M \otimes \can_{R,K}$, so we must show $g_{R_K} \circ (\Eins_M \otimes \can_{R,K}) = f_R$. But
$f$ is a polynomial law over $k$, so the diagram
\[
\xymatrix{M_R \ar[r]_{f_R} \ar[d]_{\Eins_M \otimes \can_{R,K}} &
(_kN)_R \ar[d]^{\Eins_{(_kN)} \otimes \can_{R,K}} \\ 
M_{R_K} \ar[r]_{f_{R_K}} & (_kN)_{R_K} }
\]
commutes. Hence \eqref{GER},\eqref{OMEINS} imply $g_{R_K} \circ
(\Eins_M \otimes i_{R,K}) = \omega_{R_K} \circ f_{R_K} \circ
(\Eins_M \otimes i_{R,K}) = \omega_{R_K} \circ(\Eins_{(_kN)} \otimes
i_{R,K}) \circ f_{R_K} = f_{R_K}$, as desired.
	
\emph{Uniqueness.} Suppose $g,h\:M_K \to N$ are two polynomial laws
over $K$ making the diagram \eqref{FIEGE} commutative. We must show
$g = h$. Thanks to Exercise~\ref{pr.RESCAIN}, it actually suffices to show $_kg
=\,_kh$, i.e., $g_{R_K} = h_{R_K}$ for all $R \in \kalg$. In order
to do so, let $x_1,\ldots,x_n \in M$, $r_1,\ldots,r_n \in R$,
$a_1,\ldots,a_n \in K$ and $\bfT = (\bft_1,\ldots,\bft_n)$ a finite
chain of independent indeterminates. Write $\varphi\:R_K[\bfT] \to
R_K$ for the unique homomorphism of $R_K$-algebras sending
$\bft_\nu$ to $r_\nu \otimes a_\nu$ ($1 \leq \nu \leq n$). Since
$g,h$ are polynomial laws over $K$, the squares of the diagram
\[
\xymatrix{M_{R_K[\bfT]} \ar[r]_(.45){\Eins} \ar[d]_{\Eins_M \otimes
\varphi} & (M_K)_{R_K[\bfT]} \ar[rr]_{g_{R_K[\bfT]},h_{R_K[\bfT]}}
\ar[d]^{\Eins_{M_K} \otimes_K \varphi} && N_{R_K[\bfT]}
\ar[d]^{\Eins_N \otimes_K \varphi} \\ 
M_{R_K} \ar[r]_(.45){\Eins} & (M_K)_{R_K} \ar[rr]_{g_{R_K},h_{R_K}} && N_{R_K}}
\]
commute, and using (\ref{ss.ITSCA}.\ref{ITBA}) we conclude, for $p
\in \{g,h\}$,
\begin{align*}
p_{R_K}\big(\sum x_\nu \otimes (r_\nu \otimes a_\nu)\big)
=\,\,&p_{R_K}\big(\sum (x_\nu \otimes 1_K) \otimes_K (r_\nu \otimes
a_\nu)\big) \\
=\,\,&\big(p_{R_K} \circ (\Eins_{M_K} \otimes_K
\varphi)\big)\big(\sum (x_\nu \otimes 1_K) \otimes_K \bft_\nu\big)
\\
=\,\,&\big(p_{R_K} \circ (\Eins_M \otimes \varphi)\big)\big(\sum
x_\nu \otimes \bft_\nu\big) \\
=\,\,&\big((\Eins_N \otimes_K \varphi) \circ p_{R_K[\bfT]}\big(\sum
x_\nu \otimes \bft_\nu\big),
\end{align*}
which shows that it will be enough to prove	
\[
g_{R_K[\bfT]}(\sum x_\nu \otimes t_\nu) = h_{R_K[\bfT]}(\sum x_\nu \otimes \bft_\nu).
\]
Combining our previous identifications with the fact that \eqref{FIEGE} holds for $p$ and $R[\bfT]_K = R_K[\bfT]$ under the identification $\bft_\nu \otimes 1_K = \bft_\nu$ for all $\nu = 1,\dots,n$ , this follows from
\begin{align*}
f_{R[\bfT]}(\sum x_\nu \otimes \bft_\nu) =\,\,&(_kp)_{R[\bfT]} \circ
\big(_k(\can_{M,K})\big)_{R[\bfT]} (\sum x_\nu \otimes \bft_\nu) \\
=\,\,&p_{R[\bfT]_K} \circ (\Eins_M \otimes \can_{R[\bfT],K})(\sum x_\nu
\otimes \bft_\nu) \\
=\,\,&p_{R[\bfT]_K}\big(\sum x_\nu \otimes (\bft_\nu \otimes
1_K)\big) \\
=\,\,&p_{R_K[\bfT]}(\sum x_\nu \otimes \bft_\nu).
\end{align*}
\end{sol}

\begin{sol}{pr.CUMASE} \label{sol.CUMASE} (a) Put $H := f \ast g\:M \to N$ as a polynomial law over $K$ and $\pH := (\,_kf,\,_kg)\:\,_kM \to\,_kN$ as a polynomial law over $k$. We must show $_kH = \pH$. For any polynomial law $h$ over $K$, we apply (\ref{ss.RESC}.\ref{KAFAR}) and obtain $(_kh)_k = h_{k\otimes K} = h_K$. Hence
\[
(_kH)_k = H_K (f \ast g)_K = f_K = f
\]
as a set map $M \to N$, which is the same as
\[
_kf = (_kf)_k \big((_kf) = \ast (_kg)\big)_k = \pH_k
\]
as a set map $_kM \to\,_kN$, and we have shown $(_kH)_k = \pH_k$. Similarly, by Exercise~\ref{pr.RESCACAL}~(e),
\[
\big(D(_kH)\big)_k = \big(_k(DH)\big)_k = (DH)_K = g_K = g
\]
as a set map $M \times M \to N$, which is the same as
\[
_kg = (_kg)_k = \Big(D\big((_kf) \ast (_kg)\big)\Big)_k = (D\pH)_k
\]
as a set map $_kM \times\,_kM \to\,_kN$. Thus $(D(_kH))_k = \pH_k$, and from Exercise~\ref{pr.CUMA}~\ref{CUMA.e} we conclude $_kH = \pH$, as claimed. 

(b) By definition, \eqref{SICULA} commutes as a $\sigma$-semi-linear polynomial square if and only if
\begin{align}
\vcenter{\label{KAMKAF} \xymatrix{_kM \ar[r]_{_k\vph} \ar[d]_{_k(f \ast g)} & _k\pM \ar[d]^{_k(\pf \ast \pg)} \\_kN \ar[r]_{_k\psi} & _k\pN}}
\end{align}
does as a diagram of polynomial laws over $k$. But \eqref{KAMKAF}, thanks to (a), is the same as
\[
\xymatrix{_kM \ar[r]_{_k\vph} \ar[d]_{(_kf) \ast (_kg)} & _k\pM \ar[d]^{(_k\pf) \ast (_k\pg)} \\_kN \ar[r]_{_k\psi} & _k\pN}
\]
By the solution to Exercise~\ref{pr.CUMA}~\ref{CUMA.a}, therefore, 
\[
(_k\vph,\,_k\psi)\:(_kf) \ast (_kg) \longrightarrow (_k\pf) \ast (_k\pg)
\]
is a morphism in the category $\kholaw$. Thus Exercise~\ref{pr.CUMA}~\ref{CUMA.e} and \eqref{KFASTG} imply that
\[
(_k\vph,\,_k\psi)\:(_kf,\,_kg) \longrightarrow (_k\pf,\,_k\pg)
\]
is a morphism in the category $\kcumap$, which in turn is equivalent to the diagrams \eqref{SICUMA} of set maps being commutative.
\end{sol}

\begin{sol}{pr.NONSC} \label{sol.NONSC} Fixing $R \in \kalg$, there is clearly a unique $k$-linear map 
\[
\Psi_R\:N \otimes \Pol(M,k) \longrightarrow \Map(M_R,N_R)
\]
given by
\[
\Psi_R(v \otimes f)(x) = v \otimes f_R(x)
\]
for $v \in N$, $f \in \Pol(M,k)$, $x \in M_R$. For the first part of the problem, it will be enough to show that
$\Psi_R$ varies functorially with $R$, so let $\vph\:R \to S$ be a morphism in $\kalg$. Since $f$ is a scalar polynomial law, hence the diagram 
\[
\xymatrix{
M_R \ar[r]_{f_R} \ar[d]_{\Eins_M \otimes \vph} & R \ar[d]^{\vph} \\
M_S \ar[r]_{f_S} & S}
\]
commutes, so does
\[
\xymatrix{
M_R \ar[r]_{\Psi_R(v \otimes f)} \ar[d]_{\Eins_M \otimes \vph} & N_R \ar[d]^{\Eins_N \otimes \vph} \\
M_S \ar[r]_{\Psi_S(v \otimes f)} & N_S.}
\]
Hence we have shown that a $k$-linear map $\Psi\:N \otimes \Pol(M,k) \to \Pol(M,N)$ satisfying \eqref{PSIVE} does indeed exist and is obviously unique. It remains to prove that $\Psi$ is an isomorphism if $N$ is finitely generated projective.

To this end, $N$ still being arbitrary, we write $\Psi^N := \Psi$ to indicate dependence on $N$ and consider another $k$-module $\pN$. Then the linear projections $\pi\:N \oplus \pN \to N$, $\ppi\:N \oplus \pN \to \pN$, canonically regarded as homogeneous polynomial laws of degree $1$ via Exercise~\ref{pr.MULI}, induce $k$-linear maps
\[
\pi_\ast\:\Pol(M,N \oplus \pN) \longrightarrow \Pol(M,N), \quad \ppi_\ast\:\Pol(M,N \oplus \pN) \longrightarrow \Pol(M,\pN)
\]
given by $f \mapsto \pi \circ f$, $f \mapsto \ppi \circ f$, respectively, which in turn yield an isomorphism
\[
\pi_\ast \oplus \ppi_\ast\:\Pol(M,N \oplus \pN) \overset{\sim} \longrightarrow \Pol(M,N) \oplus \Pol(M,\pN).
\] 
After the natural identification $M \otimes (N \oplus \pN) = M \otimes N \oplus M \otimes \pN$, it is straightforward to check that the diagram
\begin{align}
\vcenter{\label{PSINPR} \xymatrix{
(N \oplus \pN) \otimes \Pol(M,k) \ar[rr]_{\Psi^{N \oplus \pN}} \ar[d]_{\can}^{\cong} && \Pol(M,N \oplus \pN) \ar[d]_{\cong}^{\pi_\ast \oplus \ppi_\ast} \\
(N \otimes \Pol(M,k)) \oplus (\pN \otimes \Pol(M,k)) \ar[rr]_(0.55){\Psi^N \oplus \Psi^{\pN}} && \Pol(M,N) \oplus \Pol(M, \pN)}}
\end{align}
commutes. Now suppose $N$ is finitely generated projective and choose a $k$-module $\pN$ making $N \oplus \pN$ free of finite rank. By a repeated application of \eqref{PSINPR}, we are reduced to the case  $N = k$. But then $\Psi = \Psi^k = \Eins_{\Pol(M,k)}$, and the problem is solved.

\begin{ermk}
In \cite[Lemma 3.3]{GPR}, the result in this exercise is obtained as a consequence of \ref{skip.other.fflat}~\ref{skip.other.fflat.3}.
\end{ermk}
\end{sol}



\solnchap{Solutions for Chapter~\ref{c.ALT}}

\solnsec{Section~\ref{s.BASIDSTRASS}}

\noexsec

\solnsec{Section~\ref{s.STRAS}}

\begin{sol}{pr.COMALT} \label{sol.COMALT}
Since the associator of $A$ is
alternating, the first part follows immediately from
(\ref{ss.ACA}.\ref{ss.ACA.1}). The second part is now obvious.
\end{sol}

\begin{sol}{pr.KLEINFELD} \label{sol.KLEINFELD} A multilinear map is alternating if and only if it vanishes provided two adjacent components of the argument are equal; hence, after linearizing, it changes signs when two adjacent components are switched. For an alternative $k$-algebra $A$, we therefore have to prove
\[
f(w,x,y,y) = 0 = f(w,x,x,z) = f(x,x,y,z)
\]
for all $w,x,y,z \in A$. The first equation follows from the fact that the associator of $A$ is alternating by definition. Combined with the left Moufang identity and the linearized left alternative law, this also gives
\begin{align*}
f(w,x,x,z) =\,\,&[wx,x,z] - x[w,x,z] = \big((wx)x\big)z - (wx)(xz) - x\big((wx)z\big) + x\big(w(xz)\big) \\
=\,\,&\big((wx)x + x(wx)\big)z - \Big(\big((wx)(xz) + x\big((wx)z\big)\Big) = 0,
\end{align*}
hence the second equation. As to the third, we write $\opgrp{f}$ for the Kleinfeld function of $\opgrp{A}$, and since $A$ and $\opgrp{A}$ have the same associator, we conclude
\[
\opgrp{f}(w,x,y,z) = f(x,w,y,z),
\]
hence $f(x,x,y,z) = -f(x,y,x,z) = -\opgrp{f}(y,x,x,z) = 0$, as desired. 
\end{sol}

\begin{sol}{pr.UNALT} \label{sol.UNALT}
The implications
(i)$\Rightarrow$(ii)$\Rightarrow$(iii)$\Rightarrow$(iv) and
(i)$\Rightarrow$(v) are obvious, while (iv)$\Rightarrow$(iii)
follows from the fact that $L_x$, $R_x$, $U_x = L_xR_x$ pairwise commute.
We are thus left with the following implications.

(iv)$\Rightarrow$(i). Given $y \in A$, we find an element $z \in A$
satisfying $xzx = y$. Moreover, since (iii) holds by what we have
just seen, there are $e,f \in A$ such that $ex = x = xf$. Now the
Moufang identities yield $ey = e(xzx) = ((ex)z)x = xzx = y$, $yf =
(xzx)f = x(z(xf)) = xzx = y$, so $e$ (resp. $f$) is a left (resp.
right) unit for $A$. Hence $A$ is unital with $1_A = e = f$.

(v)$\Rightarrow$(i). Since $L_x$ is surjective, some $e \in A$ has
$xe = x$ and any $z \in A$ can be written as $z = xw$ for some $w
\in A$. Hence $x(ez) = x(e(xw)) = (xex)w = x^2w = x(xw) = xz$,
forcing $e$ to be a left unit for $A$ since $L_x$ is also injective.
Passing to $A^{\op}$, we find a right unit for $A$ as well, and (i)
holds.
\end{sol}

\begin{sol}{pr.ONESINVALT} \label{sol.ONESINVALT}
(a) The implication (ii)$\Rightarrow$(i) is
obvious, while (i)$\Leftrightarrow$(iii) becomes
(i)$\Leftrightarrow$(ii) in $A^{\op}$. Thus we only have to prove
(i)$\Rightarrow$(ii). Using the Moufang identities in operator form
(\ref{ss.MOUF}.\ref{LLMOUF}), (\ref{ss.MOUF}.\ref{RRMOUF}), we
obtain $E^2 = L_xL_yL_xL_y = L_{xyx}L_y = L_xL_y = E$ and,
similarly, $F^2 = F$. Furthermore, by (\ref{ss.UOPALT}.\ref{ULR}),
$EF = L_xL_yR_yR_x = L_xU_yR_x = U_{xy} = U_{1_A} = \Eins_A$, hence
$\Eins_A = E^2F = EEF = E$, and this is (ii).

(b) Assume $x \in A$ is right invertible, so some $y \in A$ has $xy = 1_A$. Then $L_xL_y = \Eins_A$ by (a)(ii), so $L_x\:A \to A$ is surjective. But $A$ is a finitely generated $k$-module, forcing $L_x$ to be bijective (Prop.~\ref{p.SURBI}), hence $x$ to be invertible (Prop.~\ref{p.INVALT}) and $y = x^{-1}$ to be its two-sided inverse. The case of a left invertible element is treated analogously.

(c) Take $V$ to be the $k$-vector space of all sequences $a_0, a_1, \ldots$ of elements of $k$ with coordinatewise addition, $A := \End_k(V)$, $y$ to be the map that shifts the sequence right by one place and puts a zero in the first coordinate, and $x$ to be the map that shifts the sequence left by one place and drops the first coordinate.  Then $xy = 1_A$, but neither $x$ nor $y$ is invertible because $x$ is not injective and $y$ is not surjective. 
\end{sol}

\begin{sol}{pr.INVPRODALT} \label{sol.INVPRODALT}
Since $xy \in A$ is invertible, so is
the linear operator $U_{xy}$, by Prop.~\ref{p.INVALT}. But $U_{xy} =
R_yU_xL_y$ by (\ref{ss.UOPALT}.\ref{ULR}), so $L_y$ is injective. On
the other hand, the left Moufang identity and
(\ref{p.INVALT}.\ref{OPINVALT}) imply $L_y(x(y(xy)^{-1})) =
y[x(y(xy)^{-1})] = [y(xy)](xy)^{-1} = y = L_y1_A$, and since $L_y$
is injective, we conclude $x(y(xy)^{-1}) = 1_A$, so $x$ is right
invertible with $y(xy)^{-1}$ as a right inverse. Reading this in
$A^{\op}$ shows that $y$ is left invertible with $(xy)^{-1}x$ as a
left inverse.
\end{sol}

\begin{sol}{pr.PEIRCEALT} \label{sol.PEIRCEALT}
For $i,j = 1,2$ and $x \in A$, the
expression $c_ixc_j$ is unambiguous since it takes place in the
unital subalgebra of $A$ generated by $c$ and $x$, which is
associative by Cor.~\ref{c.ARTIN}. We now claim that the following
relations hold.
\begin{align}
\label{PEMUL} L_{c_i}L_{c_j} = \delta_{ij}L_{c_i}, \quad
R_{c_i}R_{c_j} = \delta_{ij}R_{c_i}, \quad L_{c_i}R_{c_j} =
R_{c_j}L_{c_i}.
\end{align}
For $i = j$, this follows immediately from
(\ref{ss.CONALAL}.\ref{LLALT})--(\ref{ss.CONALAL}.\ref{RLEX}),
while for $i \neq j$ it suffices to observe $L_{c_j} = \Eins_A -
L_{c_i}$, $R_{c_j} = \Eins_A - R_{c_i}$ to arrive at the desired
conclusion. Now let $i,j,l,m = 1,2$ and put $E_{ij} :=
L_{c_i}R_{c_j}$. Then \eqref{PEMUL} yields $E_{ij}E_{lm} =
L_{c_i}R_{c_j}L_{c_l}R_{c_m} = L_{c_i}L_{c_l}R_{c_j}R_{c_m} =
\delta_{il}\delta_{jm}L_{c_i}R_{c_j} =
\delta_{il}\delta_{jm}E_{ij}$, so the $E_{ij}$ are indeed orthogonal
projections that add up to $\sum_{i,j=1,2} E_{ij} = \sum_{i,j=1,2}
L_{c_i}R_{c_j} = U_{c_1 + c_2} = U_{1_A} = \Eins_A$. This proves
\eqref{PALTDEC}, where the $A_{ij} = A_{ij}(c)$ are defined by the
first among the relations in \eqref{PALTC}. The remaining ones are
now obvious. Next we establish \eqref{PALTU}--\eqref{PALTT}. In
order to do so, we let $x_{ij} \in A_{ij}$, $y_{lm} \in A_{lm}$ etc.
and first prove
\begin{align}
\label{IDPE} c_ix_{jl} = \delta_{ij}x_{jl}, \quad x_{ij}c_l =
\delta_{jl}x_{ij}.
\end{align}
Indeed, by \eqref{PALTC} the first relation is clear for $i = j$,
while for $i \neq j$ we obtain $c_ix_{jl} = x_{jl} - c_jx_{jl} =
x_{jl} - x_{jl} = 0$, as claimed. The second relation is proved
analogously. Combining now the linearized right alternative law
(\ref{ss.LIN}.\ref{LINRALT}) with \eqref{IDPE}, we obtain
$(x_{ij}y_{jl})c_l = x_{ij}(y_{jl}c_l + c_ly_{jl}) -
(x_{ij}c_l)y_{jl} = x_{ij}y_{jl} + \delta_{jl}x_{ij}y_{jl} -
\delta_{jl}x_{ij}y_{jl} = x_{ij}y_{jl}$ and, similarly,
$c_i(x_{ij}y_{jl}) = x_{ij}y_{jl}$. This proves \eqref{PALTU}.
Furthermore, for $i \neq j$, the left Moufang identity and
\eqref{IDPE} yield $x_{ii}y_{jl} = (c_ix_{ii}c_i)y_{jl} =
c_i(x_{ii}(c_iy_{jl})) = 0$ and, similarly, $y_{lj}x_{ii} = 0$. Thus
\eqref{PALTD} holds. Finally, since $c_j(x_{ij}y_{ij}) = (c_jx_{ij}
+ x_{ij}c_j)y_{ij} - x_{ij}(c_jy_{ij}) = x_{ij}y_{ij} =
(x_{ij}y_{ij})c_i$, we also have \eqref{PALTT}. In particular,
$x_{ij}^2 = c_jx_{ij}^2 = (c_jx_{ij})x_{ij} = 0$, while the
remaining assertion $A_{ij}^2 = \{0\}$ for $A$ associative is
obvious.
\end{sol}

\solnsec{Section~\ref{s.HOMOTALT}}

\begin{sol}{pr.INVISOTALT} \label{sol.INVISOTALT}
By Prop.~\ref{p.INVALT}, $x \in A$ is
invertible in $A^{(p,q)}$ iff $U_x^{(p,q)} = U_xU_{pq}$ (by
(\ref{ss.UNHOMOTALT}.\ref{UOPHOMOTALT})) is bijective iff $U_x$ is
bijective iff $x$ is invertible in $A$, in which case its inverse in
$A^{(p,q)}$ is $x^{(-1,p,q)} = (U_x^{(p,q)})^{-1}x =
U_{pq}^{-1}U_x^{-1}x = U_{pq}^{-1}x^{-1}$.
\end{sol}

\begin{sol}{pr.ALBISOT} \label{sol.ALBISOT}
By definition, $f,g,h$ are linear
bijections from $A$ to $B$ satisfying $f(xy) = g(x)h(y)$ for all
$x,y \in A$. Setting $u := g^{-1}(1_B)$, $v := h^{-1}(1_B)$, we
conclude $f(uy) = h(y)$, $f(xv) = g(x)$, equivalently, $f \circ L_u
= h$, $f \circ R_v = g$. Thus $L_u,R_v\:A \to A$ are linear
bijections, forcing $A$ to be unital (Exc.~\ref{pr.UNALT}) and $u,v
\in A$ to be invertible (Prop.~\ref{p.INVALT}), with inverses $p :=
v^{-1}$, $q := u^{-1}$ This implies $f((xp)(qy)) = g(xp)h(qy) =
f((xv^{-1})v)f(u(u^{-1}y)) = f(x)f(y)$. Thus $f\:A^{(p,q)} \to B$ is
an isomorphism; in particular, $B$ is alternative.
\end{sol}

\begin{sol}{pr.STRGRALT} \label{sol.STRGRALT}
See \cite[2.2-2.8]{MR2003j:17004}. In
particular $1_{\Str(A)} = (1_A.1_A,\Eins_A)$ is the unit element of
$\Str(A)$ and $(p,q,g)^{-1} =
(g^{-1}(q^{-1}p^{-2}),g^{-1}(q^{-2}p^{-1}),g^{-1})$ is the inverse
of $(p,q,g) \in \Str(A)$.
\end{sol}

\begin{sol}{pr.EXTLRMULT} \label{sol.EXTLRMULT}
The equation $(L_ux)u\,\cdot\,u^{-2}(L_uy)
= uxu\,\cdot\,u^{-2}(uy) = u(x(u(u^{-1}y))) = u(xy) = L_u(xy)$ for
all $x,y \in A$ shows that $L_u\:A \to A^{(u,u^{-2})}$ is an
isomorphism. Hence $\tilde{L_u} \in \Str(A)$. Reading this in
$A^{op}$ also yields $\tilde{R_u} \in \Str(A)$. Invoking
(\ref{pr.STRGRALT}.\ref{STRGRALT}) and Cor.~\ref{c.ARTIN}, we now
compute
\begin{align*}
\tilde{L_u}\tilde{L_v}\tilde{L_u}
=\,\,&(u,u^{-2},L_u)(v,v^{-2},L_v)\tilde{L_u} =
\big(uu^{-2}uvu,u^{-2}uv^{-2}uu^{-2},L_uL_v\big)\tilde{L_u} \\
=\,\,&(vu,u^{-1}v^{-2}u^{-1},L_uL_v)(u,u^{-2}L_u) = (w,z,L_uL_vL_u)
= (w,z,L_{uvu}),
\end{align*}
where
\begin{align*}
w =\,\,&vuu^{-1}v^{-2}u^{-1}uvuvu = uvu, \\
z =\,\,&u^{-1}v^{-2}u^{-1}uvu^{-2}vuu^{-1}v^{-2}u^{-1} =
u^{-1}v^{-1}u^{-2}v^{-1}u^{-1} = (uvu)^{-2}.
\end{align*}
Thus $\tilde{L_u}\tilde{L_v}\tilde{L_u} = (uvu,(uvu)^{-2},L_{uvu}) =
\tilde{L}_{uvu}$. Reading this in $A^{\op}$, we also obtain
$\tilde{R_u}\tilde{R_v}\tilde{R_u} = \tilde{R}_{uvu}$. Since
$\tilde{L}_{1_A} = 1_{\Str(A)} = \tilde{R}_{1_A}$, the relations
$(\tilde{L}_u)^{-1} = \tilde{L}_{u^{-1}}$, $(\tilde{R}_u)^{-1} =
\tilde{R}_{u^{-1}}$ are now clear. Next we compute
\begin{align*}
\tilde{L}_u\tilde{R}_v &= (u,u^{-2},L_u)(v^{-2},v,R_v) =
(uu^{-2}uv^{-2}u,u^{-2}uvuu^{-2},L_uR_v) \\
&= (v^{-2}u,u^{-1}vu^{-1},L_uR_v),
\end{align*}
as claimed. A similar computation or passing to $A^{\op}$ yields the
analogous equation for $\tilde{R}_v\tilde{L_u}$. The rest is clear.
\end{sol}

\begin{sol}{pr.UNITSTRGRALT} \label{sol.UNITSTRGRALT}
 See \cite[3.2, 3.3]{MR2003j:17004}.
\end{sol}

\begin{sol}{pr.UNITARBISOTALT} \label{sol.UNITARBISOTALT}
(a) The extended right multiplication by
$pq$ as defined in Exc.~\ref{pr.EXTLRMULT} yields an isomorphism
\[
R_{pq}\:A \overset{\sim} \longrightarrow A^{((pq)^{-2},pq)},
\]
which by functoriality \ref{ss.FULT} and Prop.~\ref{p.ITHOMOTALT}
may also be regarded as an isomorphism
\[
R_{pq}\:A^{(p,q)} \overset{\sim} \longrightarrow
(A^{((pq)^{-2},pq)})^{(p^2q,qpq)} = A^{(p^\prime,q^\prime)}
\]
with
\begin{align*}
p^\prime =\,\,&(pq)^{-2}(pq)p^2q(pq)^{-2} =
q^{-1}p^{-1}pq^{-1}p^{-1} = q^{-2}p^{-1}, \\
q^\prime =\,\,&(pq)(qpq)(pq)^{-2}(pq) = pq^2.
\end{align*}
Thus $R_{pq}\:A^{(p,q)} \to A^{pq^2}$ is an isomorphism, which may
also be verified directly, using (\ref{ss.SINVALT}.\ref{MINVALT})
and the Moufang identities:
\begin{align*}
(R_{pq}x)(pq^2)^{-1}\,\cdot\,(pq^2)R_{pq}y)
=\,\,&\big(x(pq)\big)(q^{-2}p^{-1})\,\cdot\,\big((pq)q\big)\big(y(pq)\big)
\\
=\,\,&\big(x(pq)\big)(q^{-2}p^{-1})\,\cdot\,(pq)(qy)(pq) \\
=\,\,&\Big[\Big(\big[\big(x(pq)\big)(q^{-2}p^{-1})\big](pq)\Big)(qy)\Big](pq)
\\
=\,\,&\Big[\Big(x\big((pq)(q^{-2}p^{-1})(pq)\big)\Big)(qy)\Big](pq)
\\
=\,\,&\big((xp)(qy)\big)(pq) = R_{pq}\big((xp)(qy)\big).
\end{align*}
(b) The argument uses (\ref{ss.UNTISOTALT}.\ref{ITUNTISOTALT}) to
conclude $A^{(pq)r} = (A^{pq})^r = ((A^p)^q)^r = (A^p)^{qr} =
A^{p(qr)}$, which by (\ref{ss.UNTISOTALT}.\ref{EQIS}) yields $(pq)r
= up(qr)$ for some $u \in \Nuc(A)^\times$. But the argument is
faulty since with $B := A^p$ we have $((A^p)^q)^r = (B^q)^r = B^s$,
where $s$ is the product of $q$ and $r$ \emph{not in $A$ but in $B =
A^p$}: $s = (qp^{-1})(pr)$, forcing $B^s = (A^p)^{(qp^{-1})(pr)} =
A^\trans$, $t = p[(qp^{-1})(pr)] = [p(qp^{-1})p]r = (pq)r$, in complete
agreement with the previous computation. Moreover, the conclusion
does not hold, which follows from looking at the real octonions
$\dots$.
\end{sol}

\begin{sol}{pr.NUCISOTALT} \label{sol.NUCISOTALT}
 The solution to
Exc.~\ref{pr.UNITARBISOTALT} shows that right multiplication by $pq$
yields an isomorphism from $A^{(p,q)}$ to an appropriate unital
isotope of $A$. It therefore suffices to show $\Nuc(A^p) = \Nuc(A)$
and $\Cent(A^p) = \Cent(A)$ for all $p \in A^\times$; actually,
since $A = (A^p)^{p^{-1}}$, it will be enough to verify the
inclusions $\Nuc(A) \subseteq \Nuc(A^p)$, $\Cent(A) \subseteq
\Cent(A^p)$. Writing $[\emptyslot,\emptyslot,\emptyslot]^p$ (resp. $[\emptyslot,\emptyslot]^p$) for the
associator (resp. the commutator) of $A^p$, this will follow once we
have shown
\begin{gather}
\label{ASSP} [x,y,z]^p = [x,yp^{-1},pz] - [x,p^{-1},(py)z] \quad \text{and}  \\
[x,y]^p = [x,y] + [x,p^{-1},py] - [y,p^{-1},px] \notag
\end{gather}
for all $x,y \in A$. In order to derive the first relation of
\eqref{ASSP}, we compute
\begin{align*}
[x,y,z]^p =\,\,&\Big(\big((xp^{-1})(py)\big)p^{-1}\Big)(pz) -
(xp^{-1})\Big(p\big((yp^{-1})(pz)\big)\Big) \\
=\,\,&\big(x(yp^{-1}\big)(pz) - (xp^{-1})\big((py)z\big) \\
=\,\,&[x,yp^{-1},pz] + x\big((yp^{-1})(pz)\big) -
(xp^{-1})\big((py)z\big),
\end{align*}
where the second summand on the right agrees with
\[
x\Big[p^{-1}\Big(p\big((yp^{-1})(pz)\big)\Big)\Big] =
x\Big(p^{-1}\big((py)z\big)\Big) = -[x,p^{-1},(py)z] +
(xp^{-1})\big((py)z\big).
\]
Inserting this into the previous equation yields the first relation
of \eqref{ASSP}. For the second, we obtain $[x,y]^p = (xp^{-1})(py)
- (yp^{-1})(px) = [x,p^{-1},py] + xy - [y,p^{-1},px] - yx =[x,y] +
[x,p^{-1},py] - [y,p^{-1},px]$, as claimed.
\end{sol}


\solnchap{Solutions for Chapter~\ref{c.COMAL}}

\solnsec{Section~\ref{s.CONAL}}

\begin{sol}{pr.NILRADCON} \label{sol.NILRADCON}
(a) We begin by deriving the following identities:
\begin{align}
\label{NOSQ} n_C(x^2) = \,\,&n_C(x)^2, \\
\label{NORASQ} n_C(x,x^2) =\,\,&t_C(x)n_C(x), \\
\label{TRSQ} t_C(x^2) =\,\,&t_C(x)^2 - 2n_C(x).
\end{align}
Indeed, the identities of \ref{ss.BASID} yield $n_C(x^2) = n_C(t_C(x)x - n_C(x)1_C) = t_C(x)^2n_C(x) - t_C(x)^2n_C(x) + n_C(x)^2$, hence \eqref{NOSQ}, while $n_C(x,x^2) = n_C(x,t_C(x)x - n_C(x)1_C = 2t_C(x)n_C(x) - t_C(x)n_C(x)$ yields \eqref{NORASQ} and $t_C(x^2) = t_C(t_C(x)x - n_C(x)1_C) = t_C(x)^2 - 2n_C(x)$ is \eqref{TRSQ}. Combining these identities with those of \eqref{ss.BASID} and writing $y,z \in k[x]$ as $y = \alpha_01_C + \alpha_1x$, $y = \beta_01_C + \beta_1x$ with $\alpha_0,\alpha_1,\beta_0,\beta_1 \in k$, we now compute
\begin{align*}
n_C(yz) &- n_C(y)n_C(z) =\, n_C\big((\alpha_01_C + \alpha_1x)(\beta_01_C + \beta_1x)\big) \\ 
&\quad - n_C(\alpha_01_C + \alpha_1x)n_C(\beta_01_C + \beta_1x) \\
&=\,\,n_C\big(\alpha_0\beta_01_C + (\alpha_0\beta_1 + \alpha_1\beta_0)x + \alpha_1\beta_1x^2\big) \\
&\quad - \big(\alpha_0^2 + \alpha_0\alpha_1t_C(x) + \alpha_1^2n_C(x)\big)\big(\beta_0^2 + \beta_0\beta_1t_C(x) + \beta_1^2n_C(x)\big) \\
&=\,\,\alpha_0^2\beta_0^2 + \alpha_0\beta_0(\alpha_0\beta_1 + \alpha_1\beta_0)t_C(x)  + \alpha_0\beta_0\alpha_1\beta_1t_C(x^2) \\
&\quad+ (\alpha_0\beta_1 + \alpha_1\beta_0)^2n_C(x) + (\alpha_0\beta_1 + \alpha_1\beta_0)\alpha_1\beta_1n_C(x,x^2) + \alpha_1^2\beta_1^2n_C(x^2) \\
&\quad- \big(\alpha_0^2 + \alpha_0\alpha_1t_C(x) + \alpha_1^2n_C(x)\big)\big(\beta_0^2 + \beta_0\beta_1t_C(x) + \beta_1^2n_C(x)\big) \\
&=\,\,\alpha_0^2\beta_0^2 + \alpha_0\beta_0(\alpha_0\beta_1 + \alpha_1\beta_0)t_C(x)  + \alpha_0\alpha_1\beta_0\beta_1t_C(x)^2 - 2\alpha_0\alpha_1\beta_0\beta_1n_C(x) \\
&\quad+ (\alpha_0\beta_1 + \alpha_1\beta_0)^2n_C(x) + (\alpha_0\beta_1 + \alpha_1\beta_0)\alpha_1\beta_1t_C(x)n_C(x) + \alpha_1^2\beta_1^2n_C(x)^2 \\
&\quad- \big(\alpha_0^2 + \alpha_0\alpha_1t_C(x) + \alpha_1^2n_C(x)\big)\big(\beta_0^2 + \beta_0\beta_1t_C(x) + \beta_1^2n_C(x)\big) \\
&=\,\,\alpha_0^2\beta_0^2 + \alpha_0^2\beta_0\beta_1t_C(x) + \alpha_0\alpha_1\beta_0^2t_C(x) + \alpha_0\alpha_1\beta_0\beta_1t_C(x)^2 + \alpha_0^2\beta_1^2n_C(x) \\
&\quad+ \alpha_1^2\beta_0^2n_C(x) + \alpha_0\alpha_1\beta_1^2t_C(x)n_C(x) + \alpha_1^2\beta_0\beta_1t_C(x)n_C(x) + \alpha_1^2\beta_1^2n_C(x)^2 \\
&\quad- \big(\alpha_0^2 + \alpha_0\alpha_1t_C(x) + \alpha_1^2n_C(x)\big)\big(\beta_0^2 + \beta_0\beta_1t_C(x) + \beta_1^2n_C(x)\big) \\
&=\,\,0.
\end{align*}
This completes the proof of (a).

(b) If $x$ is invertible in $k[x]$ with inverse $x^{-1} \in k[x]$, then (a) yields 
\[
n_C(x)n_C(x^{-1}) = n_C(xx^{-1}) = n_C(1_C) = 1, 
\]
hence $n_C(x) \in k^\times$ and $n_C(x^{-1}) = n_C(x)^{-1}$. Conversely, if $n_C(x) \in k^\times$, then $y := n_C(x)^{-1}\bar x$ by (\ref{ss.BASID}.\ref{CACONJC}) satisfies $xy = 1_C = yx$. Thus $x$ is invertible in $k[x]$ with inverse $y$. 

(c) If $t_C(x),n_C(x) \in \Nil(k)$, the elements $t_C(x)x,n_C(x)1_C \in k[x]$ are nilpotent and commute. But then $x^2 = t_C(x)x - n_C(x)1_C$ is nilpotent. Hence so is $x$. Conversely, assume $x$ is nilpotent. Then $x^n = 0$ for some positive integer $n$, and we show by induction on $n$ that $t_C(x)$ and $n_C(x)$ are nilpotent. If $n = 1$, there is nothing to prove. If $n > 1$, then we put $m := 1$ for $n = 2$ and $m := \lfloor\frac{n}{2}\rfloor +1 > \frac{n}{2}$ for $n > 2$. This implies $2m \geq n$, hence $(x^2)^m = 0$, and $m < n$. Thanks to the induction hypothesis, therefore, the scalars $t_C(x^2) = t_C(x)^2 - 2n_C(x)$ (by \eqref{TRSQ}) and $n_C(x^2) = n_C(x)^2$ (by \eqref{NOSQ}) are both nilpotent. Hence $n_C(x)$ is nilpotent and so is $t_C(x)^2 =t_C(x^2) + 2n_C(x)$. But then $t_C(x)$ itself must be nilpotent, which completes the induction.
\end{sol}

\begin{sol}{pr.CONALGHOMNOR} \label{sol.CONALGHOMNOR}
$n := n_{C^\prime} \circ \vph\:C \to k$ is a quadratic form such that $Dn = (Dn_{C^\prime}) \circ (\vph \times \vph)$. Since $\vph$ preserves units, we conclude $n(1_C) = 1$ and $t := (Dn)(1_C,\emptyslot) = t_{C^\prime} \circ \vph$. Now let $x \in C$. Then
\begin{align*}
\vph\big(t(x)x - n(x)1_C\big) =\,\,&t_{C^\prime}\big(\vph(x)\big)\vph(x) - n_{C^\prime}\big(\vph(x)\big)1_{C^\prime} = \vph(x)^2 = \vph(x^2),
\end{align*}
and since $\vph$ is injective, $x^2 = t(x)x - n(x)1_C$. Thus not only $n_C$ but also $n$ is a norm for the conic algebra $C$. But by the hypotheses on $C$, its norm is unique (Prop.~\ref{p.UNINOR}), which implies $n_C = n = n_{C^\prime} \circ \vph$, as desired.

If we drop the assumption that $\vph$ be injective, the asserted conclusion is false: let $C := k \times k$, $C^\prime := k$ and $\vph\:k \times k \to k$ be the projection onto the first factor. Then \ref{ss.MOTEX}~(c),(d) yields $n_k(\vph((\alpha,\beta))) = n_k(\alpha) = \alpha^2$ for all $\alpha,\beta \in k$, while $n_{k \times k}((\alpha,\beta)) = \alpha\beta$. Thus $\vph$ does not preserve norms, hence, though being a unital algebra homomorphism, is \emph{not} one of conic algebras.
\end{sol}

\begin{sol}{pr.CONMO} \label{sol.CONMO}
(a) One verifies immediately that $C$ is unital with identity element $e$. Defining $n_C\:C \to k$ by
\[
n_C(\alpha e + x) := \alpha^2 + \alpha T(x) + B(x,x)
\]
for $\alpha \in k$, $x \in M_\lambda$, we obviously obtain a quadratic form with bilinearization
\[
n_C(\alpha e + x,\beta e + y) = 2\alpha\beta + \alpha T(y) + \beta T(x) + B(x,y) + B(y,x)
\]
for $\alpha,\beta \in k$, $x,y \in M_\lambda$. We have $n_C(e) = 1$, while $t_C := (Dn_C)(\emptyslot,e)$ satisfies
\[
t_C(\alpha e + x) = 2\alpha + T(x)
\] 
for $\alpha \in k$, $x \in M_\lambda$. Hence
\begin{align*}
t_C(\alpha e + x)(\alpha e + x) - n_C(\alpha e + x)e =\,\,&\big(2\alpha + T(x)\big)(\alpha e + x) \\
&- \big(\alpha^2 + \alpha T(x) + B(x,x)\big)e \\
=\,\,&\big(\alpha^2 - B(x,x)\big)e + \Big(\big(\alpha + T(x)\big)x + \alpha x\Big) \\
=\,\,&(\alpha e + x)^2, 
\end{align*}
which shows that $C$ is a conic $k$-algebra with norm $n_C$ and trace $t_C$.

(b) We first show that $(T,B,K)$ is a conic co-ordinate system. The defining conditions of such a system are obviously fulfilled, with the following two exceptions. (i) the bilinear map $K$ takes values in $M_\lambda$, and (ii) $K$ is alternating. In order to establish (i), we let $x,y \in M_\lambda = \Ker(\lambda)$ and obtain
\[
\lambda(x \times y) = t_C(x)\lambda(y) - \lambda(xy) + \lambda(xy)\lambda(e) = -\lambda(xy) + \lambda(xy) = 0,
\]
hence $x \times y \in M_\lambda$. In order to establish (ii), we compute
\[
x \times x = t_C(x)x - x^2+ \lambda(x^2)e = n_C(x)e + \lambda\big(t_C(x)x - n_C(x)e\big)e = n_C(x)e - n_C(x)e = 0,
\]
so $K$ is indeed alternating. We now denote by ``$\cdot$'' the product in the conic algebra $C^ \prime := \Con(T,B,K)$ and obtain, for $\alpha,\beta \in k$, $x,y \in M_\lambda$, 
\begin{align*}
(\alpha e + x)\cdot(\beta e + y) =\,\,&\big(\alpha\beta + \lambda(xy)\big)e + \Big(\big(\alpha + t_C(x)\big)y \\ 
&+ \beta x - t_C(x)y + xy - \lambda(xy)e\Big) \\
=\,\,&\alpha\beta e + \alpha y + \beta x + xy = (\alpha e + x)(\beta e + y), \\
n_{C^ \prime}(\alpha e + x) =\,\,&\alpha^ 2 + \alpha t_C(x) - \lambda (x^2) = \alpha^ 2 \\ 
&+ \alpha t_C(x) - t_C(x)\lambda(x) + n_C(x)\lambda(e) \\
=\,\,&\alpha^ 2 + \alpha t_C(x) + n_C(x) = n_C(\alpha e + x).
\end{align*}
Summing up, we have proved $C = \Con(T_C,B_C,K_C)$. Finally, let $(T,B,K)$ be a conic co-ordinate system and $C : = \Con(T,B,K)$. Consulting the multiplication formula for $C$ (resp. the formula for the trace of $C$), we conclude
\begin{align*}
T_C(x) =\,\,&t_C(x) = T(x), \\
B_C(x,y) =\,\,&-\lambda(xy) = B(x,y) \\
K_C(x,y) =\,\,&t_C(x)y -xy + \lambda(xy)e = T(x)y - xy + \lambda(xy)e \\
=\,\,&\big(B(x,y + \lambda(xy)\big)e + \big(T(x)y - T(x)y + K(x,y)\big) = K(x,y).
\end{align*}
Thus $(T_C,B_C,K_C) = (T,B,K)$ and we have proved that the two constructions in (a), (b) are inverse to each other.
\end{sol}

\begin{sol}{pr.CONFIELD} \label{sol.CONFIELD}
By assumption we have $1_C \wedge x \wedge x^2 = 0$ for all $x \in C$. Replacing $x$ by $x + \alpha y$, $x,y \in C$, $\alpha \in k$, we obtain (since $k$ contains
more than two elements)
\begin{align}
\label{LIQA} 1_C \wedge y \wedge x^2 + 1_C \wedge x \wedge (x \circ y) = 0 &&(x,y \in C).
\end{align}
Now put 
\begin{align}
\label{TRAZE} U := \{0\} \cup \{u \in C \setminus k1_C \mid u^2 \in k1_C\}
\end{align}
We claim that \emph{$U$ is a vector subspace of $C$}. Clearly being closed under scalar multiplication, we must show $u + v \in U$ for all $u,v \in U$ such that $u,v$ are linearly independent. Then $1_C,u$ are linearly independent, and assuming $v = \alpha 1_C + \beta u$ for some $\alpha,\beta \in k$, we obtain $\alpha^21_C + 2\alpha\beta u + \beta^2u^2 = v^2 \in k1_C$, hence $\alpha = 0$ or $\beta = 0$, a contradiction. Thus $1_C,u,v$ are linearly independent. On the other hand, setting $x = u$, $y = v$ (resp. $x = v$, $y = u$) in \eqref{LIQA}, and observing \eqref{TRAZE}, we obtain $u \circ v = \alpha 1_C + \beta u = \gamma 1_C + \delta v$ for some $\alpha,\beta,\gamma,\delta \in k$, hence $u \circ v \in k1_C$. But this amounts to $(u + v)^2 \in k1_C$, so we have proved that $U$ is  closed under addition and thus, altogether, a vector subspace of $C$.

Now let $x \in C \setminus k1_C$. Then $x^2 = \alpha 1_C + \beta x$ for some $\alpha,\beta \in k$, which implies $(x - \frac{\beta}{2}1_C)^2 \in k1_C$, hence $x - \frac{\beta}{2}1_C \in U$, and we have shown $C = k1_C \oplus U$ as a direct sum of subspaces. Now define a quadratic form $n_C\:C \to k$ by $n_C(\alpha 1_C + u) := \alpha^2 + n_C(u)$ for $\alpha \in k$, $u \in U$, where $n_C(u) \in k$ is determined by $u^2 = -n_C(u)1_C$. This implies
\begin{align*}
n_C(\alpha 1_C + u,\beta 1_C + v) &= 2\alpha\beta + n_C(u,v), \quad u \circ v = -n_C(u,v)1_C, \\
&(\alpha,\beta \in k,\;u,v \in U),
\end{align*}
hence $t_C(\alpha 1_C + u) := n_C(1_C,\alpha 1_C + u) = 2\alpha$. Summing up, we obtain
\begin{align*}
(\alpha 1_C + u)^2 =\,\,&\alpha^21_C + 2\alpha u + u^2 = 2\alpha(\alpha 1_C + u) - \big(\alpha^ 2 + n_C(u)\big)1_C \\
=\,\,&t_C(\alpha 1_C + u)(\alpha 1_C + u) - n_C(\alpha 1_C + u)1_C.
\end{align*}
Thus $C$ is a conic algebra over $k$ with norm $n_C$.
\end{sol}

\begin{sol}{pr.DICON} \label{sol.DICON}
Since conic algebras are stable under base change, they clearly satisfy the Dickson condition. Conversely, suppose $C$ satisfies the Dickson condition. By hypothesis, the vector $1_C \in C$ is unimodular, so we find a submodule $M \subseteq C$
such that
\begin{align}
\label{UNIM} C = k1_C \oplus M
\end{align}
as a direct sum of $k$-modules. Since the assignment $x \mapsto x^2$ defines a quadratic map $M \to C$, the projections to the direct summands $k1_C \cong k$ and $M$ of the decomposition \eqref{UNIM} give rise to quadratic maps $n\:M \to k$ and $s \:M \to M$ such that
\begin{align}
\label{SQUADEC} x^2 = s(x) - n(x)1_C &&(x \in M).
\end{align}
In particular, $n$ is a quadratic form over $k$ that will eventually become the norm (restricted to $M$) of the prospective conic $k$-algebra $C$. On the other hand,
by the Dickson condition, we also have $s(x) \in kx$ for $x \in M$, and we would like to think of $s(x)$ as $t(x)x$ with $t(x) \in k$
becoming the trace of $x$ in the prospective conic $k$-algebra $C$. But since the annihilator of $x$ may not be zero, we don't even know at this stage how to make $t(x)$ a quantity that is well defined, let alone a linear form. For this reason, we bring the hypotheses on the module
structure of $C$ and the \emph{strict} Dickson condition into play.

Let us first reduce to the case that $M$ is a free $k$-module. Otherwise, $M$ is finitely generated projective by hypothesis, so there exists a finite family of elements $f \in k$
that generate the unit ideal in $k$ and make $M_f$ a free $k_f$-module of finite rank, for each such $f$. Assuming the free case has been settled, it follows that all $C_f$ are conic, and Prop.~\ref{p.UNINOR} ensures that the norms $n_{C_f}$ glue to give a quadratic form $n\:C \to k$. It is then clear that $C$ is a conic algebra with norm $n$ since this is so after changing scalars from $k$ to each $k_f$.

For the rest of the proof, we may therefore assume that $M$ is a free $k$-module, possibly of infinite rank, with basis  $(e_i)_{i \in I}$.  We let $\bfT = (\bft_i)_{i \in I}$ be a family
of independent variables and write $k[\bfT]$ for the corresponding polynomial ring. For a \emph{non-empty finite} subset
$E \subseteq I$, we consider the submodule
\begin{align}
\label{EME} M^E := \sum_{i \in E}\,ke_i \subseteq M,
\end{align}
which is a direct summand, and the element
\begin{align}
\label{EXE} \bfx^E := \sum_{i \in E}\,e_i \otimes \bft_i \in M^E_{k[\bfT]} \subseteq M_{k[\bfT]}.
\end{align}
We claim the annihilator of $\bfx^E$ in $k[\bfT]$ is zero; indeed, for $f(\bfT) \in k[\bfT]$ the relation $f(\bfT)\bfx^E = 0$ implies $\sum e_i \otimes(\bft_if(\bfT)) = 0$, hence $\bft_if(\bfT) = 0$ for all $i \in E$
and then $f(\bfT) = 0$ since $\bft_i$ is not a zero divisor in $k[\bfT]$. Invoking the strict Dickson condition, we therefore find a unique polynomial $g(\bfT) \in k[\bfT]$ such that
\[
s(\bfx^E) = g(\bfT)\bfx^E.
\]
Specializing this relation with an additional variable $\bfu$ to $\bfu\bfT$, we obtain $g(\bfu\bfT)\bfu\bfx^E = s(\bfu\bfx^E) = \bfu^2s(\bfx^E) = (\bfu g(\bfT))\bfu\bfx^E$, so the polynomial $g(\bfT) \in k[\bfT]$
is homogeneous of degree $1$. Thus there exists a unique linear form $t^E\:M^E \to k$ satisfying the relation $s(x) = t^E(x)x$ for all $x \in M^E$ because any such $t^E$ will be
converted into the linear homogeneous polynomial $g(\bfT)$ after extending scalars from $k$ to $k[\bfT]$. Here the uniqueness condition implies that the linear forms
$t^E$ as $E$ varies over the non-empty finite subsets of $I$ glue to give a linear form $t\:M \to k$ such that $s(x) = t(x)x$ for all $x \in M$. Combining with \eqref{SQUADEC}, we obtain the relation
\begin{align}
\label{SQUAM.SOL} x^2 - t(x)x + n(x)1_C = 0
\end{align}
for all $x \in M$. We now extend $t,n$ as given on $M$ to all of $C$ by
\begin{align*}
t(\alpha 1_C + x) = 2\alpha + t(x), \quad n(\alpha 1_C + x) = \alpha^2 + \alpha t(x) + n(x) &&(\alpha \in k,\;x \in M)
\end{align*}
and then conclude from a straightforward computation that \eqref{SQUAM.SOL} holds for all $x \in C$. Since we also have $t = n(1_C,\emptyslot)$, it follows that $C$ is a conic algebra. 
\end{sol}

\begin{sol}{pr.IDEMP} \label{sol.IDEMP}
For the first part, assume $k \neq \{0\}$ is connected and let $c \in C$. If $n_C(c) = 0$ and $t_C(c) = 1$, then $c$ is an idempotent by (\ref{ss.BASID}.\ref{CAQUAD}) which cannot be zero since $t_C(c) = 1$, and cannot be $1$ since $n_C(c) = 0$. Conversely, let $c \in C$ be an idempotent $\neq 0,1_C$. From Exc.~\ref{pr.NILRADCON}~(a) we deduce that $n_C(c) \in k$ is an idempotent.  If $n_C(c) = 1$, then $c$ is invertible in $k[c]$ (Exc.~\ref{pr.NILRADCON}~(b)), whence the relation $c(1_C - c) = 0$ implies $c = 1_C$, a contradiction. Hence $n_C(c) = 0$, and (\ref{ss.BASID}.\ref{CAQUAD}) reduces to $c = c^2 = t_C(c)c$. Taking traces, we conclude that $t_C(c) \in C$ is an idempotent, which cannot be zero since this would imply $c = 0$. Thus $t_C(c) = 1$. We now proceed to establish the equivalence of (i)--(iv).

(i) $\Rightarrow$ (ii). This is clear.

(ii) $\Rightarrow$ (iii). For any prime ideal $\mfp \subseteq k$, the ring $k_\mfp$ is local, hence connected, and the first part of the problem yields $n_C(c)_\mfp = n_{C_\mfp}(c_\mfp) = 0$, $t_C(c)_\mfp = t_{C_\mfp}(c_\mfp) = 1_\mfp$, and since this holds true for all $\mfp \in \Spec(k)$, condition (iii) follows.

(iii) $\Rightarrow$ (iv). $c$ is clearly an idempotent. Moreover, the second condition of (iii) shows that $c$ is unimodular, and since $t_C(1_C - c) = 2 - 1 = 1$, so is $1_C - c$.

(iv) $\Rightarrow$ (i). Unimodularity is stable under base change, so for $R \in \kalg$, $R \neq \{0\}$, the elements $c_R$ and $1_{C_R} - c_R$ are both different from zero. Hence (i) holds.
\end{sol}

\begin{sol}{pr.CONID}\label{sol.CONID} (a) Let $\mfa$ be an ideal in $k$. An element of the ideal $\mfa C \subseteq C$ has the form $x = \sum \alpha_ix_i$, $\alpha_i \in \mfa$, $x_i \in C$, and expanding $n_C(x)$ accordingly, we conclude $n_C(x) \in \mfa$; the relation $n_C(x,y) \in \mfa$ for $y \in C$ is even more obvious. Thus $(\mfa,\mfa C)$ is a conic ideal in $C$. By definition, the ideal $\mfa C \subseteq C$ is the smallest one with this property. Now let $I$ be an ideal in $C$ and write $\mfa$ for the ideal in $k$ generated by the expressions $n_C(x)$, $n_C(x,y)$ for $x \in I$, $y \in C$. By \eqref{ELCONID}, the solution of (a) will be complete once we have shown that $(\mfa,I)$ is a conic ideal in $C$. Since $n_C(x)1_C = x\bar x$ and $n_C(x,y)1_C = x\bar y + y\bar x$ by (\ref{ss.BASID}.\ref{CACONJC}) and its bilinearization, we clearly have $\mfa C \subseteq I$, while \eqref{ELCONID} is trivially fulfilled.
	
(b) follows from a straightforward computation using (\ref{ss.SELICON}.\ref{ENCEL}).

(c) Since $\mfa C \subseteq I$, the $k$-algebra $C_0 := C/I$ carries a unique $k_0$-algebra structure making $\pi$ a $\sigma$-semi-linear homomorphism of unital algebras. From \eqref{ELCONID} we deduce that there is a unique set map $n_{C_0}\:C_0 \to k_0$ making a commutative diagram
\begin{equation}
\label{ENCEPI} \xymatrix{C \ar[r]_{\pi} \ar[d]_{n_C} & C_0 \ar@{.>}[d]^{\exists!\, n_{C_0}} \\
k \ar[r]_{\sigma} &k_0.}
\end{equation}
Since $\sigma$ and $\pi$ are both surjective, one checks that $n_{C_0}$ is, in fact, a quadratic form over $k_0$ making $C_0$ a conic $k_0$-algebra. Now \eqref{ENCEPI} , being a commutative diagram of set maps, shows that $\pi\:C \to C_0$ is a $\sigma$-semi-linear homomorphism of conic algebras having $\Ker(\sigma,\pi) = (\mfa,I)$, as claimed. 
\end{sol}

\begin{sol}{pr.CONILID} \label{sol.CONILID} 
(a) Suppose first that $\mfa$ is a nil ideal in $k$. For $x \in I$, we deduce from \eqref{ELCONID} in Exc.~\ref{pr.CONID} that $n_C(x)$ and $t_C(x) = n_C(x,1_C)$ both belong to $\mfa $, so $x$ by Exc.~\ref{pr.NILRADCON}~(c) is nilpotent. Thus $I$ is a nil ideal in $C$. Conversely, let this be so. Then $\mfa C \subseteq I$ is a nil ideal in $C$, forcing in particular $\alpha 1_C$ to be nilpotent for all $\alpha \in \mfa$. By Exc.~\ref{pr.NILRADCON}~(c) again, this implies that $\alpha^2 = n_C(\alpha 1_C)$ is nilpotent, so $\mfa$ is a nil ideal in $k$. 
	
(b) We have $\Nil(k)C \subseteq \Nil(C)$ by Exc.~\ref{pr.NILRAD}, so we need only show $n_C(x)$, $n_C(x,y) \in \mfa := \Nil(k)$ for all $x \in I := \Nil(C)$, $y \in C$. The first inclusion being obvious, we turn to the second by combining the ideal property of $I$ with (\ref{ss.BASID}.\ref{CAQUADL}):
\[
I \ni xy + yx = t_C(x)y + t_C(y)x - n_C(x,y)1_C,
\]
where $t_C(x)y \in \mfa C \subseteq I$ by hypothesis and $t_C(y)x \in I$ for trivial reasons. Thus $n_C(x,y)1_C \in I$, and we conclude $n_C(x,y)^2 = n_C(n_C(x,y1)_C) \in \mfa$, hence $n_C(x,y) \in \mfa$.

(c) Let $c \in \pi^{-1}(c_0)$ be an idempotent in $C$. Then $\pi(c) = c_0$, and (\ref{ss.SELICON}.\ref{ENCEL}) shows $\sigma(n_C(c)) = n_{C_0}(c_0) = 0$, so $n_C(c) \in \mfa$ is nilpotent. But by Exc.~\ref{pr.NILRADCON}~(a), it is also an idempotent, and we conclude $n_C(c) = 0$. Similarly, $\sigma(t_C(c)) = t_{C_0}(c_0) = 1$, forcing $t_C(c) \in 1 + \mfa$ to be invertible. On the other hand, applying $t_C$ to the equation $c = c^2 = t_C(c)c$ (by (\ref{ss.BASID}.\ref{CAQUAD})), it follows that $t_C(c)$ is an idempotent. Thus $t_C(c) = 1$, and we have shown that the idempotent $c \in C$ is elementary.
\end{sol}

\begin{sol}{pr.DECID} \label{sol.DECID}
 Recall from (\ref{ss.PROPID}.\ref{PIPM}), (\ref{ss.PROPID}.\ref{EMPM}) that there is a natural identification $M_R = \vep M$ as $R$-modules such that $x_R = x \otimes \vep = \vep x$ for all $x \in M$. 

(ii) $\Leftarrow$ (i). Obvious.

(i) $\Rightarrow$ (ii). If $c \in C$ is an idempotent, then (\ref{ss.BASID}.\ref{CATQUAD}) implies $t_C(c) = t_C(c^2) = t_C(c)^2 - 2n_C(c)$, hence
\begin{align}
\label{IDNOTR} t_C(c)\big(1 - t_C(c)\big) = -2n_C(c).
\end{align}
The quantities $\vep^{(i)}$, $i = 0,1,2$, as defined in \eqref{COAL.VEPI} obviously satisfy $\sum\vep^{(i)} = 1$, while the fact that $n_C(c) \in C$ is an idempotent by Exc.~\ref{pr.NILRADCON}~(a) combined with \eqref{IDNOTR} implies $\vep^{(i)}\vep^{(j)} = 0$ for $i,j = 0,1,2$ distinct. Thus they form a complete orthogonal system of idempotents in $k$. Our preliminary remark shows $C^{(i)} = C_{k^{(i)}}$ and $c = \sum c^{(i)}$, where $c^{(i)} = c_{k^{(i)}} = \vep^{(i)}c$ for $i = 0,1,2$. We now compute $c^{(0)} = \vep^{(0)}c = (1 - n_C(c))(1 - t_C(c))c$, where $(1 - t_C(c))c = c^ 2 - t_C(c)c = -n_C(c)1_C$, which implies $c^{(0)} = 0$. Next we turn to
\[
c^{(1)} = \vep^{(1)}c = \big(1 - n_C(c)\big)t_C(c)c = \big(1 - n_C(c)\big)\big(c + n_C(c)1_C\big) = \big(1 - n_C(c)\big)c,
\]
which together with \eqref{IDNOTR} implies
\begin{align*}
t_C(c^{(1)}) =\,\,&\vep^{(1)}t_C(c) = \big(1 - n_C(c)\big)t_C(c)^2 = \big(1 - n_C(c)\big)\big(t_C(c) + 2n_C(c)\big) \\
=\,\,&\big(1 - n_C(c)\big)t_C(c) = \vep^{(1)} = 1_{k^{(1)}}, \\
n_C(c^{(1)}) =\,\,&\vep^{(1)}n_C(c) = \big(1 - n_C(c)\big)t_C(c)n_C(c) = 0.
\end{align*}
Thus $c^{(1)} \in C^{(1)}$ is an elementary idempotent. Finally,
\[
c^{(2)} = \vep^{(2)}c = n_C(c)c,
\]
which implies $n_C(c^{(2)}) = \vep^{(2)}n_C(c) = \vep^{(2)} = 1_{C^{(2)}}$. Thus $c^{(2)} \in C^{(2)}$ is an invertible idempotent in the sense of Exc.~\ref{pr.NILRADCON}~(b), forcing $c^{(2)} = 1_{C^{(2)}}$, as claimed. The statement just proved implies $t_C(c^{(2)}) = 2$ in $k^{(2)}$. This can also be proved directly by noting
\[
\big(2 - t_C(c)\big)\vep{^{(2)}} = 2n_C(c) - t_C(c)n_C(c) = n_C(c,c^2) - t_C(c)n_C(c)
\]
which is zero by \eqref{NORASQ} of the solution to Exc.~\ref{pr.NILRADCON}.

It remains to prove uniqueness of the $\vep^{(i)}$, so let $\eta^{(i)}$, $i = 0,1,2$, be any complete orthogonal system of idempotents in $C$ satisfying mutatis mutandis the conditions of (ii), and define the $\vep^{(i)}$, $i = 0,1,2$, as in \eqref{COAL.VEPI}. Then $n_C(c^{(i)}) = n_C(\eta^{(i)}c) = \eta^{(i)}n_C(c) = n_C(c)_{k^{(i)}} = n_{C^{(i)}}(c^{(i)})$ and, similarly, $t_C(c^{(i)}) = t_{C^{(i)}}(c^{(i)})$ for $i = 0,1,2$. Since $c^{(1)} \in C^{(1)}$ is an elementary idempotent and $c^{(2)} = 1_{C^{(2)}}$ by hypothesis, this implies
\begin{align*}
n_C(c) =\,\,&\big(0,n_{C^{(1)}}(c^{(1)}),n_{C^{(2)}}(c^{(2)})\big) = (0,0,1_{k^{(2)}}) = \eta^{(2)}, \\
t_C(c) =\,\,&\big(0,t_{C^{(1)}}(c^{(1)}),t_{C^{(2)}}(c^{(2)})\big) = (0,1_{k^{(1)}},2\cdot 1_{k^{(2)}}) = \eta^{(1)} + 2\eta^{(2)},
\end{align*}
hence $\eta^{(2)} = n_C(c) = \vep^{(2)}$ and $\eta^{(1)} =t_C(c) - 2n_C(c)$. But the idempotents $\eta^{(1)}, \eta^{(2)}$ are orthogonal, which implies $t_C(c)n_C(c) = 2n_C(c)$, and we conclude $\eta^{(1)} = t_C(c)(1 - n_C(c)) = \vep^{(1)}$ by \eqref{COAL.VEPI}. Since the orthogonal systems formed by the $\eta$'s and by the $\vep$'s are both complete, we also have $\eta^{(0)} = \vep^{(0)}$, as claimed.
\end{sol}

\begin{sol}{pr.TRICO} \label{sol.TRICO}
If $C$ is a conic $k$-algebra with trivial
conjugation, then $2x = t_C(x)1_C$ for all $x \in C$ by
(\ref{ss.BASID}.\ref{CACONJ}). Since the trace form of $C$ is surjective, we find an element
$u \in C$ satisfying $t_C(u) = 1$. This implies $2u = 1_C$, hence
$x = 2ux = t_C(ux)1_C$. But $C$ is a faithful $k$-module. Hence $C \cong k$ as conic algebras. The converse is obvious.
\end{sol}

\begin{sol}{pr.NORCOM} \label{sol.NORCOM} For the first part we have to show $n_C(x,[y,y]) = 0 = n_C(x,[x,y])$ for all $x,y \in C$. The first assertion is obvious, while the second one follows from $n_C(x,[x,y]) = n_C(x,xy) - n_C(x,yx)$ and (\ref{p.NORAS}.\ref{WLBILCOMP}), (\ref{p.NORAS}.\ref{WRBILCOMP}). In order to establish \eqref{PERWE}, \eqref{NORCOM}, we expand $n_C(x \circ y)$ by using (\ref{ss.BASID}.\ref{CAQUADL}), (\ref{ss.BASID}.\ref{CATN}) 
to obtain
\begin{align*}
n_C(x \circ y) =\,\,&n_C\big(t_C(x)y + t_C(y)x - n_C(x,y)1_C\big) \\
=\,\,&t_C(x)^2n_C(y) + t_C(x)t_C(y)n_C(x,y) - t_C(x)t_C(y)n_C(x,y) \\
\,\,&+ t_C(y)^2n_C(x) - t_C(x)t_C(y)n_C(x,y) + n_C(x,y)^2 \\
=\,\,&t_C(x)^2n_C(y) + t_C(y)^2n_C(x) + n_C(x,y)^2 - t_C(x)t_C(y)n_C(x,y),
\end{align*}
and (\ref{ss.BASID}.\ref{CANOCO}) yields
\begin{align}
\label{NORCIR} n_C(x \circ y) = t_C(x)^2n_C(y) + t_C(y)^2n_C(x) - n_C(x,y)n_C(x,\bar y).    
\end{align}
Now let $\vep = \pm 1$. From (\ref{ss.BASID}.\ref{CAQUADL}), (\ref{p.NORAS}.\ref{WLBILCOMP}), (\ref{p.NORAS}.\ref{WRBILCOMP}), (\ref{p.NORAS}.\ref{ASSBILT}) we deduce
\begin{align*}
n_C(xy + \vep yx) =\,\,&n_C(xy) + \vep n_C(xy,yx) + n_C(yx) \\
=\,\,&n_C(xy) + (1 - 2\vep)n_C(yx) + \vep n_C(x \circ y,yx) \\
=\,\,&n_C(xy) + (1 - 2\vep)n_C(yx) \\
\,\,&+ \vep\big(t_C(x)^2n_C(y) + t_C(y)^2n_C(x) - n_C(x,y)n_C(x,\bar y)\big),
\end{align*}
and \eqref{NORCIR} implies
\[
n_C(xy + \vep yx) = n_C(xy) + (1 - 2\vep)n_C(yx) + \vep n_C(x \circ y).
\]
Setting $\vep = 1$ in the last equation yields \eqref{PERWE}, while setting $\vep = -1$ and using \eqref{PERWE}, \eqref{NORCIR} yields \eqref{NORCOM}.
\end{sol}

\solnsec{Section~\ref{s.COALT}}

\begin{sol}{pr.NONUNNOR} \label{sol.NONUNNOR}
For the subsequent computations it is important to note
\begin{align}
\label{DRODOT} (\alpha.\beta_0)\vep = (\alpha\beta_0)\vep &&(\alpha \in k,\;\beta_0 \in k_0) 
\end{align}
since, for $\alpha = \alpha_0 + \alpha_0^\prime\vep$, $\alpha_0,\alpha_0^\prime \in k_0$, the right-hand side becomes $(\alpha_0\beta_0 + (\alpha_0^\prime\beta_0)\vep)\vep = (\alpha_0\beta_0)\vep = (\alpha.\beta_0)\vep$, as claimed.

Turning to the exercise itself, we begin by showing associativity. Writing $(e_i)_{1\leq i\leq 3}$ for the canonical basis of $k_0^3$ over $k_0$, we deduce from \eqref{ETHR} that $[k_0^3,k_0^3] = k_0e_3$ and $[e_3,k_0^3] = \{0\}$, hence $[[k_0^3,k_0^3],k_0^3] = \{0\} = [k_0^3,[k_0^3,k_0^3]]$. Now let $\alpha,\beta,\gamma \in k$ and $u,v,w \in k_0^3$. Then
\begin{align*}
\big((\alpha,u)(\beta,v)\big)(\gamma,w) =\,\,&\big(\alpha\beta,\alpha.v + \beta.u + [u,v]\big)(\gamma,w) \\
=\,\,&\big(\alpha\beta\gamma,(\alpha\beta).w + (\gamma\alpha).v + (\beta\gamma).u + \gamma.[u,v] \\
\,\,&+ \alpha.[v,w] + \beta.[u,w] + [[u,v],w]\big), \\
(\alpha,u)\big((\beta,v)(\gamma,w)\big) =\,\,&(\alpha,u)\big(\beta\gamma,\beta.w + \gamma.v + [v,w]\big) \\
=\,\,&\big(\alpha\beta\gamma,(\alpha\beta).w + (\gamma\alpha).v + \alpha.[v,w] \\
\,\,&+ (\beta\gamma).u + \beta.[u,w] + \gamma.[u,v] + [u,[v,w]]\big),
\end{align*}
and the preceding observation shows $[[u,v],w] = [u,[v,w]] = 0$. Hence $C$ is associative.

The algebra $C$ is obviously unital with identity element $1_C = (1,0)$, and the quadratic form $n_z$ obviously satisfies $n_z(1_C) = 1$. Moreover, \eqref{DRODOT} yields
\begin{align*}
n_z\big((\alpha,u)\big) = \alpha^2 + (\alpha z^\trans u)\vep &&(\alpha \in k,\; u \in k_0^3),
\end{align*}
and we have
\begin{align*}
n_z\big((\alpha,u),(\beta,v)\big) = 2\alpha\beta + (\alpha z^\trans v + \beta z^\trans u)\vep &&(\alpha,\beta \in k,\;u,v, \in k_0^3),
\end{align*}
hence
\begin{align*}
t_z\big((\alpha,u)\big) := n_z\big(1_C,(\alpha,u)\big) = 2\alpha + (z^\trans u)\vep.
\end{align*}
Comparing
\[
(\alpha,u)^2 = (\alpha^2,2\alpha.u)
\] 
with
\begin{align*}
t_z\big((\alpha,u)\big)(\alpha,u) - n_z\big((\alpha,u)\big)1_C =\,\,&\big(2\alpha + (z^\trans u)\vep\big)(\alpha,u) - \big(\alpha^2 + (\alpha z^\trans u)\vep\big)(1,0) \\
=\,\,&\big(2\alpha^2 + (\alpha z^\trans u)\vep - \alpha^2 -(\alpha z^\trans u)\vep, \\
\,\,&(2\alpha + (z^\trans u)\vep).u\big) \\
=\,\,&(\alpha^2,2\alpha.u),
\end{align*}
we see that $C$ is indeed a conic $k$-algebra with norm $n_z$ and trace $t_z$. We now turn to the equivalence of (i)--(iv). By Props.~\ref{p.MULAD} and \ref{p.NAFL}, we have (i) $\Rightarrow$ (ii) $\Rightarrow$ (iii). Thus it remains to show (iii) $\Rightarrow$ (iv) $\Rightarrow$ (i).

(iii) $\Rightarrow$ (iv). Since the conjugation of $C$ is an involution and $\Ann(C) = \{0\}$, we conclude from Prop.~\ref{p.CHAFL} that $t_z(xy) = n_z(x,\bar y)$ for all $x,y \in C$. Writing $x = \alpha \oplus u$, $y = \beta \oplus v$, $\alpha,\beta \in k$, $u,v \in k_0^3$, we apply \eqref{DRODOT} and obtain
\begin{align*}
t_z(xy) =\,\,&t_z\big((\alpha\beta,\alpha.v + \beta.u + [u,v])\big) = 2\alpha\beta + \big(z^\trans (\alpha.v + \beta.u + [u,v])\big)\vep \\
=\,\,&2\alpha\beta + \big(\alpha.(z^\trans v) + \beta.(z^\trans u) + z^\trans [u,v]\big)\vep \\ 
=\,\,& 2\alpha\beta + (\alpha z^\trans v + \beta z^\trans u + z^\trans [u,v])\vep, \\
\bar y =\,\,&t_z\big((\beta,v)\big)1_C - (\beta,v) = \big(2\beta + (z^ tv)\vep\big)(1,0) - (\beta,v) = \big(\beta + (z^ tv)\vep,-v\big), \\
n_z(x,\bar y) =\,\,&n_z\Big((\alpha,u),\big(\beta + (z^\trans v)\vep,,-v\big)\Big) \\ 
=\,\,&2\alpha\beta + 2\alpha z^\trans v\vep + \big(-\alpha z^ tv + \beta z^\trans u + (z^ tv)\vep u\big)\vep \\
=\,\,&2\alpha\beta +(\alpha z^\trans v + \beta z^\trans u)\vep.
\end{align*}
Comparing, we conclude $z^\trans [u,v] = 0$ for all $u,v \in k_0^3$, and specializing $u = e_1$, $v = e_2$, $(e_i)_{1\leq i\leq 3}$ being the $k_0$-basis of unit vectors in $k_0^3$, yields $\delta_3 = 0$, hence (iv).

(iv) $\Rightarrow$ (i). For $\alpha,\beta \in k$, $u,v \in k_0^3$, we compute
\begin{align*}
n_z\big((\alpha,u)(\beta,v)\big) =\,\,&n_z\big((\alpha\beta,\alpha.v + \beta.u + [u,v])\big) \\
=\,\,&(\alpha\beta)^2 + \big((\alpha^2\beta)z^\trans v + (\alpha\beta^2)z^\trans u + (\alpha\beta)z^\trans [u,v]\big)\vep, \\
n_z\big((\alpha,u)\big)n_z\big((\beta,v)\big) =\,\,&\big(\alpha^2 + (\alpha z^\trans u)\vep\big)\big(\beta^2 + (\beta z^\trans v)\vep\big) \\
=\,\,&\alpha^2\beta^2 + (\alpha^2\beta z^\trans v + \alpha\beta^2z^\trans u)\vep. 
\end{align*} 
By definition of the algebra $A$, the $e_1$- and $e_2$-components of $[u,v]$ are both zero, as is the $e_3$-component of $z$, by (iv). Hence $z^\trans [u,v] = 0$, and we conclude that $n_z$ permits composition, i.e., the conic algebra $C$ with norm $n_z$ is multiplicative.
\end{sol}

\begin{sol}{pr.NIRA} \label{sol.NIRA}
 By Prop.~\ref{p.MULAD}~(a), $C$ is norm-associative. We put
\[
N := \{x \in C \mid n_C(x), n_C(x,y) \in \Nil(k)\; \text{for all $y \in C$}\}
\]
It is straightforward to check that $N \subseteq C$ is a $k$-submodule. We claim that it is, in fact, an ideal. In order to see this, let $x \in N$ and $y,z \in C$. By multiplicativity, $n_C(xy) = n_C(x)n_C(y) = n_C(yx)$ is nilpotent, while norm-associativity yields the same for $n_C(xz,y) = n_C(x,y\bar z)$ and $n_C(zx,y) = n_C(x,\bar zy)$. Thus $xz$ and $zx$ both belong to $N$, forcing $N$ to be an ideal in $C$. Since $n_C(x)$ and $t_C(x) = n_C(x,1_C)$ are nilpotent, so is $x$ by Exc.~\ref{pr.NILRADCON}~(c). Hence $N \subseteq C$ is a nil ideal and as such contained in the nil radical of $C$. Conversely, let $x \in \Nil(C)$. Then $x$ is nilpotent, as is $x\bar y \in \Nil(C)$ for all $y \in C$. Consulting Exc.~\ref{pr.NILRADCON}~(c) again, we see that $n_C(x)$ and $n_C(x,y) = t_C(x\bar y)$ belong to $\Nil(k)$, which implies $x \in N$ and completes the proof.
\end{sol}

\begin{sol}{pr.ARTCONALG} \label{sol.ARTCONALG}
 For the first part we must show that the submodule $M$ of $C$ spanned by the elements $1_C,x,y,xy$ is, in fact, a subalgebra, i.e., it is closed under multiplication. In order to see this, we apply the alternative laws, the identities of \ref{ss.BASID}, and (\ref{ss.IDCO}.\ref{QAUOP}) to obtain
\begin{align*}
z^2 =\,\,&t_C(z)z - n_C(z)1_C \in M \qquad\qquad\qquad\qquad\qquad\qquad (z \in \{x,y,xy\}), \\
x(xy) =\,\,&x^2y = t_C(x)xy - n_C(x)y \in M, \\
(xy)x =\,\,&n_C(x,\bar y)x - n_C(x)\bar y = n_C(x,\bar y)x - n_C(x)t_C(y)1_C + n_C(x)y, \\
yx =\,\,&x \circ y - xy = t_C(x)y + t_C(y)x - n_C(x,y)1_C - xy \in M. 
\end{align*}
By symmetry, therefore, $M$ is closed under multiplication, giving the first part of the problem. As to the second, it suffices to verify the associative law on the generators of $M$, which is straightforward.
\end{sol}

\begin{sol}{pr.NORMA} \label{sol.NORMA}
Abbreviating $1 := 1_C$, $n := n_C$, $t := t_C$ and expanding the norm of an associator, we obtain
\[
n\big([x_1,x_2,x_3]\big) = n\big((x_1x_2)x_3\big) - n\big((x_1x_2)x_3,x_1(x_2x_3)\big) +
n\big(x_1(x_2x_3)\big),
\]
where multiplicativity of the norm yields
\begin{align}
\label{NASL} n\big([x_1,x_2,x_3]\big) = 2n(x_1)n(x_2)n(x_3) -
n\big((x_1x_2)x_3,x_1(x_2x_3)\big).
\end{align}
Turning to the second summand on the right of \eqref{NASL}, we
obtain, by (\ref{ss.MULCO}.\ref{MUPLT}),
\[
n\big((x_1x_2)x_3,x_1(x_2x_3)\big) = n(x_1x_2,x_1)n(x_3,x_2x_3) -
n\big((x_1x_2)(x_2x_3),x_1x_3\big).
\]
Since $C$ is norm associative by Prop.~\ref{p.MULAD}, we may apply apply (\ref{p.NORAS}.\ref{WLBILCOMP}) to the first summand on the right, which yields
\begin{align}
\label{NAST} n\big((x_1x_2)x_3,x_1(x_2x_3)\big) =
t(x_2)^2n(x_3)n(x_1) - n\big((x_1x_2)(x_2x_3),x_1x_3\big).
\end{align}
Manipulating the expression $(x_1x_2)(x_2x_3)$ by means of
(\ref{ss.BASID}.\ref{CAQUADL}) and the middle Moufang identity
(\ref{ss.MOUF}.\ref{MMOUF}), we obtain
\begin{align*}
(x_1x_2)(x_2x_3) =\,\,&(x_1x_2) \circ (x_2x_3) - (x_2x_3)(x_1x_2)
\\
=\,\,&t(x_1x_2)x_2x_3 + t(x_2x_3)x_1x_2 - n(x_1x_2,x_2x_3)1 -
x_2(x_3x_1)x_2, 
\end{align*}
where associativity of the trace (\ref{ss.NORAS}.\ref{ADMASSBILT}) and (\ref{p.NORAS}.\ref{ASSBILT}),(\ref{ss.IDCO}.\ref{QAUOP})
yield
\begin{align}
\label{FOURP} (x_1x_2)(x_2x_3) =\,\,&t(x_1x_2)x_2x_3 +
t(x_2x_3)x_1x_2
- n(x_1x_2,x_2x_3)1 - \\
\,\,&t(x_1x_2x_3)x_2 + n(x_2)\overline{x_3x_1}. \notag
\end{align}
Here we use (\ref{p.NORAS}.\ref{RASSNOR}),(\ref{ss.BASID}.\ref{CACONJ})
to compute
\begin{align*}
n(x_1x_2,x_2x_3) =\,\,&n(x_1,x_2x_3\bar x_2) = t(x_2)n(x_1,x_2x_3)
- n(x_1,x_2x_3x_2), 
\end{align*}
and (\ref{p.NORAS}.\ref{ASSBILT}),(\ref{ss.IDCO}.\ref{QAUOP}) give
\begin{align*}
n(x_1x_2,x_2x_3) =\,\,&t(x_2)n(x_1,x_2x_3) - t(x_2x_3)n(x_1,x_2) +
t(x_3x_1)n(x_2) \\
=\,\,&t(x_1)t(x_2)t(x_2x_3) - t(x_2)t(x_1x_2x_3) -
t(x_1)t(x_2)t(x_2x_3) + \\
\,\,&t(x_1x_2)t(x_2x_3) + t(x_3x_1)n(x_2),
\end{align*}
hence
\begin{align*}
n(x_1x_2,x_2x_3) = t(x_1x_2)t(x_2x_3) - t(x_2)t(x_1x_2x_3) +
t(x_3x_1)n(x_2).
\end{align*}
Inserting this into \eqref{FOURP}, and \eqref{FOURP} into the
second term on the right of \eqref{NAST}, we conclude
\begin{align*}
n\big((x_1x_2)(x_2x_3),x_1x_3\big) =\,\,&t(x_1x_2)n(x_2x_3,x_1x_3)
+ t(x_2x_3)n(x_1x_2,x_1x_3) - \\
\,\,&t(x_1x_2)t(x_2x_3)t(x_3x_1) + t(x_2)t(x_3x_1)t(x_1x_2x_3) -
\\
\,\,&t(x_3x_1)^2n(x_2) - t(x_1x_2x_3)n(x_2,x_1x_3) + \\
\,\,&n\big((x_3x_1)^\ast,x_1x_3\big)n(x_2) \\
=\,\,&t(x_1x_2)t(x_1x_2^\ast)n(x_3) + t(x_2x_3)t(x_2x_3^\ast)n(x_1) - \\
\,\,&t(x_1x_2)t(x_2x_3)t(x_3x_1) + t(x_2)t(x_3x_1)t(x_1x_2x_3) -
\\
\,\,&t(x_3x_1)^2n(x_2) - t(x_2)t(x_3x_1)t(x_1x_2x_3) + \\
\,\,&t(x_1x_2x_3)t(x_2x_1x_3) + t(x_3x_1^2x_3)n(x_2),
\end{align*}
where we may use
(\ref{ss.NORAS}.\ref{ADMASSBILT}), (\ref{ss.BASID}.\ref{CAQUAD}), (\ref{p.NORAS}.\ref{ASSBILT})
to expand
\begin{align*}
t(x_3x_1^2x_3)n(x_2) - t(x_3x_1)^2n(x_2) =\,\,&t(x_3^2x_1^2)n(x_2)
-t(x_3x_1)^2n(x_2) \\
=\,\,&t\big([t(x_3)x_3 - n(x_3)1][t(x_1)x_1 - n(x_1)1]\big)n(x_2)
- \\
\,\,&t(x_3x_1)^2n(x_2) \\
=\,\,&t(x_3)t(x_1)t(x_3x_1)n(x_2) - t(x_1)^2n(x_2)n(x_3) - \\
\,\,&t(x_3)^2n(x_1)n(x_2) + 2n(x_1)n(x_2)n(x_3) - \\
\,\,&t(x_3x_1)^2n(x_2) \\
=\,\,&t(x_3x_1)t(x_3\bar x_1)n(x_2) - t(x_1)^2n(x_2)n(x_3) - \\
\,\,&t(x_3)^2n(x_1)n(x_2) + 2n(x_1)n(x_2)n(x_3).
\end{align*}
Inserting the resulting expression
\begin{align*}
n\big((x_1x_2)(x_2x_3),x_1x_3\big)
=\,\,&\sum\,t(x_ix_j)t(x_ix_j^\ast)n(x_l) -
t(x_1x_2)t(x_2x_3)t(x_3x_1) + \\
\,\,&t(x_1x_2x_3)t(x_2x_1x_3) - t(x_1)^2n(x_2)n(x_3) - \\
\,\,&t(x_3)^2n(x_1)n(x_2) + 2n(x_1)n(x_2)n(x_3)
\end{align*}
into \eqref{NAST} and \eqref{NAST} into \eqref{NASL}, the assertion
follows. 
\end{sol}

\begin{sol}{pr.ISTQALT} \label{sol.ISTQALT}
From Prop.~\ref{p.ITHOMOTALT} combined with \ref{ss.ISOTALT} we deduce that $C^{(p,q)}$ is a unital alternative $k$-algebra with identity element $1_C^{(p,q)} = (pq)^{-1}$. Defining $n_{C^{(p,q)}} := n_C(pq)n_C$, Prop.~\ref{p.CHIN} yields $n_{C^{(p,q)}}(1_C^{(p,q)}) = 1$, and we have $n_{C^{(p,q)}}(x,y) =\break n_C(pq)n_C(x,y)$ for all $x,y \in C$, which implies
\[
t_{C^{(p,q)}}(x) := n_{C^{(p,q)}}(1_C^{(p,q)},x) = n_C(pq)n_C((pq)^{-1},x\big) = n_C(\overline{pq},x).
\]
Writing $x^{(m,p,q)}$ for the $m$-th power of $x \in C$ in $C^{(p,q)}$, we apply (\ref{ss.CONHOMOTALT}.\ref{HOMOTALT}), the middle Moufang identity (\ref{ss.MOUF}.\ref{MMOUF}) and (\ref{ss.IDCO}.\ref{QAUOP}) to conclude
\begin{align*}
x^{(2,p,q)} = \,\,&(xp)(qx) = x(pq)x = n_C(x,\overline{pq})x - n_C(x)\overline{pq} \\ 
=\,\,&t_{C^{(p,q)}}(x)x - n_{C^{(p,q)}}(x)(pq)^{-1} = t_{C^{(p,q)}}(x)x - n_{C^{(p,q)}}(x)1_C^{(p,q)}.
\end{align*}  
Thus $C^{(p,q)}$ is a conic alternative $k$-algebra with norm and trace as indicated. Moreover,
\begin{align*}
\iota_{C^{(p,q)}}(x) =\,\,&t_{C^{(p,q)}}(x)1_C^{(p,q)} - x = n_C(\overline{pq},x)(pq)^{-1} - x \\ 
=\,\,&n_C(pq)^{-1}\big(n_C(\overline{pq},x)\overline{pq} - n_C(\overline{pq})x\big) = n_C(pq)^{-1}\overline{pq}\bar x\,\overline{pq}, 
\end{align*}
giving the desired formula for the conjugation of $C^{(p,q)}$. Finally, if $C$ is multiplicative, then
\begin{align*}
n_{C^{(p,q)}}(x\,._{p,q}\,y) =\,\,&n_C(pq)n_C\big((xp)(qy)\big) = n_C(pq)n_C(x)n_C(p)n_C(q)n_C(y) \\ =\,\,&n_C(pq)n_C(x)n_C(pq)n_C(y) = n_{C^{(p,q)}}(x)n_{C^{(p,q)}}(y)
\end{align*}
for all $x,y \in C$, and $C^{(p,q)}$ is multiplicative as well.
\end{sol}

\solnsec{Section~\ref{s.CADI}}

\begin{sol}{pr.CDMUL} \label{sol.CDMUL}
We put $C := \Cay(B,\mu)$ and let $u_1,u_2,v_1,v_2 \in B$. Using the fact that $B$ is norm associative by Prop.~\ref{p.MULAD}~(a) and applying (\ref{p.CDCON}.\ref{CDMULT}), (\ref{p.CDCON}.\ref{CDNOR}) we obtain
\begin{align*}
n_C\big((u_1 + v_1j)(u_2 + v_2j)\big) =\,\,&n_C\big((u_1u_2 + \mu\bar v_2v_1) + (v_1\bar u_2 + v_2u_1)\big) \\
=\,\,&n_B(u_1u_2 + \mu\bar v_2v_1) - \mu n_C(v_1\bar u_2 + v_2 u_1) \\
=\,\,&n_B(u_1)n_B(u_2) + \mu n_B(u_1u_2,\bar v_2v_1) + \mu^2n_B(v_2)n_B(v_1) \\
\,\,&-\mu n_B(v_1)n_B(u_2) - \mu n_B(v_1\bar u_2,v_2u_1) - \mu n_B(v_2)n_B(u_1) \\
=\,\,&\big(n_B(u_1) - \mu n_B(v_1)\big)\big(n_B(u_2) - \mu n_B(v_2)\big) \\
\,\,&+ \mu n_B\big(v_2(u_1u_2),v_1\big) - \mu n_B\big(v_1,(v_2u_1)u_2)\big) \\
=\,\,&n_B(u_1 + v_1j)n_B(u_2 + v_2j) - \mu n_B\big([v_2,u_1,u_2],v_1\big),
\end{align*}
from which the assertion can be read off immediately.
\end{sol}

\begin{sol}{pr.CDSCAPAR} \label{sol.CDSCAPAR}
 For $u_1,u_2,v_1,v_2 \in B$, we combine the hypothesis that $a$ belongs to the nucleus of $B$ with the fact that the conjugation of $B$ is an involution\break (Prop.~\ref{p.MULAD})~(a) to compute
\begin{align*}
\vph(u_1 + v_1j)\vph(u_2 + v_2j) =\,\,&\big(u_1 + (av_1)j\big)\big(u_2 + (av_2)j\big) \\
=\,\,&(u_1u_2 + \mu \overline{av_2}av_1) + \big((av_2)u_1 + (av_1)\bar u_2\big)j \\
=\,\,&\big(u_1u_2 + \mu \bar v_2(\bar aa)v_1\big) + \big(a(v_2u_1 + v_1\bar u_2)\big)j \\
=\,\,&\big(u_1u_2 + n_B(a)\mu\bar v_2v_1\big) + \big(a(v_2u_1 + v_1\bar u_2)\big)j \\
=\,\,&\vph\Big(\big(u_1u_2 + n_B(a)\mu\bar v_2v_1\big) + (v_2u_1 + v_1\bar u_2)j\Big) \\
=\,\,&\vph\big((u_1 + v_1j)(u_2 + v_2j)\big).
\end{align*}
Thus $\vph$ is an obviously unital algebra homomorphism. Moreover, setting $C := \Cay(B,n_B(a)\mu)$, $C^\prime := \Cay(B,\mu)$, we have, for all $u,v \in B$,
\begin{align*}
n_{C^\prime} \circ \vph(u + vj) =\,\,&n_{C^\prime}\big(u + (av)j\big) = n_B(u) - \mu n_B(av) = n_B(u) - n_B(a)\mu n_B(v) \\
=\,\,&n_C(u + vj)
\end{align*}
since $B$ is multiplicative. Thus $\vph$ is a homomorphism of conic algebras. The final statement is obvious.
\end{sol}

\begin{sol}{pr.ZERDIV} \label{sol.ZERDIV}
 (iii) means the same in $C$ as in $C^{\op}$, while (i) in $C$ is equivalent to (ii) in $C^{\op}$. Hence it suffices to show (i) $\Leftrightarrow$ (iii).

(i) $\Rightarrow$ (iii). Assume $f := n_C(x)$ is not a zero divisor in $k$. We must show that $x$ is not a right zero divisor in $C$. By \ref{ss.ZEDIMO}, $f$ is not a zero divisor in $C$, whence it follows from Bourbaki \cite[II.2, Prop.~4]{MR0360549} that the natural map $z \mapsto z_f = z/1$ from $C$ to $C_f$ is injective. On the other hand, $n_{C_f}(x_f) = n_C(x)_f = f_f \in k_f^\times$, forcing $x_f$ to be invertible in $C_f$ (Prop.~\ref{p.CHIN}). Now suppose $y \in C$ satisfies $xy = 0$. Then $x_fy_f = (xy)_f = 0$ in $C_f$, hence $y_f = 0$, which by what we have just seen implies $y = 0$. Thus $x$ is not a right zero divisor of $C$. 

(iii) $\Rightarrow$ (i). Assume $x$ is not a right zero divisor of $C$. We must show that $n_C(x)$ is not a zero divisor in $k$, so let $\alpha \in k$ satisfy $\alpha n_C(x) = 0$. Then $0 = \alpha x\bar x = x(\alpha\bar x)$ implies $\alpha\bar x = 0$, hence $\alpha x = \overline{\alpha\bar x} = 0$. But this means $x(\alpha 1_C) = 0$, hence $\alpha 1_C = 0$, and since $1_C \in C$ is unimodular, we conclude $\alpha = 0$, as desired. 
\end{sol}

\begin{sol}{pr.NZKT}
We prove it for $R = k[\bft]$. A nonzero element $f$ of $A_R$ can be written uniquely as $f = \sum_{i=0}^m f_i \bft^i$ where $f_i \in A$ and $f_m \ne 0$.  The element $f_m$ is called the leading coefficient of $f$.  Let $g$ also be a nonzero element of $A_R$ and write $g = \sum_{i=0}^n g_i \bft^i$ where $g_n \ne 0$.  The product $fg$ has leading coefficient, the coefficient of $t^{m+n}$, $f_m g_n$.

Therefore, if $A_R$ has zero divisors $f$, $g$ such that $fg = 0$, then $f_m g_n = 0$ and $f_m$ is a zero divisor in $A$.

The converse is trivial.  The $k[[\bft]]$-case is similar
\end{sol}

\begin{sol}{pr.VARCD} \label{sol.VARCD}
 We put $C := \Cay(B^{\op},\mu)$, $C^\prime := \Cay^\prime(B,\mu)$ and 
define $\vph\:C \to C^{\prime\op}$ by
\begin{align*}
\vph(u + vj) := u + j^\prime v &&(u,v \in B).
\end{align*} 
Then $\vph$ is a $k$-linear bijection, and indicating the products in $B^{\op}$ as well as in $C^{\prime\op}$ by ``$\cdot$'', but the ones in $B,C$ and $C^\prime$ by juxtaposition, we obtain, for all $u_1,u_2,v_1,v_2 \in B$,
\begin{align*}
\vph(u_1 + v_1j)\cdot\vph(u_2 + v_2j) =\,\,&(u_2 + j^\prime v_2)(u_1 + j^\prime v_1) \\ 
=\,\,&(u_2u_1 + \mu v_1\bar v_2) + j^\prime(\bar u_2v_1 + u_1v_2) \\
=\,\,&(u_1\cdot u_2 + \mu \bar v_2\cdot v_1) + j^\prime(v_1\cdot\bar u_2 + v_2\cdot u_1) \\
=\,\,&\vph\big((u_1\cdot u_2 + \mu\bar v_2\cdot v_1) + (v_1\cdot\bar u_2 + v_2\cdot u_1)j\big) \\
=\,\,&\vph\big((u_1 + v_1j)(u_2 + v_2j)\big).
\end{align*}
Hence $\vph$ is an algebra isomorphism. Transporting the norm of $C$ to $C^\prime$ by means of $\vph$, therefore, gives a quadratic form $n_{C^\prime}$ on $C^\prime$ making it a conic algebra isomorphic to $C$ under $\vph$. Since the norms of $B$ and $B^{\op}$ coincide, we have $n_C = n_B \perp (-\mu)n_B$, which obviously implies
\begin{align*}
n_{C^\prime}(u + j^\prime v) =\,\,&n_B(u) - \mu n_B(v), \\
t_{C^\prime}(u + j^\prime v) =\,\,&t_B(u), \\
\overline{u + j^\prime v} =\,\,&\bar u - j^\prime v 
\end{align*}
for all $u,v \in B$.

Finally, let us assume that the conjugation of $B$ is an involution. Then we claim that
\begin{align*}
\psi\:\Cay(B,\mu)^{\op} \overset{\sim} \longrightarrow \Cay(B^{\op},\mu), \quad u + vj \longmapsto u + \bar vj,
\end{align*}
is an isomorphism of conic algebras. Since $\psi$ is obviously bijective and preserves norms, it will be enough to show that it is an algebra homomorphism. Indeed, extending the previous conventions to $\Cay(B,\mu)^{\op}$, we obtain, for all $u_1,u_2,v_1,v_2 \in B$,
\begin{align*}
\psi(u_1 + v_1j)\psi(u_2 + v_2j) =\,\,&(u_1 + \bar v_1j)(u_2 + \bar v_2j) \\ 
=\,\,&(u_1\cdot u_2 + \mu v_2\cdot\bar v_1) + (\bar v_2\cdot u_1 + \bar v_1\cdot\bar u_2)j \\
=\,\,&(u_2u_1 + \mu \bar v_1v_2) + (u_1\bar v_2 + \bar u_2\bar v_1)j \\
=\,\,&(u_2u_1 + \mu\bar v_1v_2) + \overline{v_2\bar u_1 + v_1u_2}j \\
=\,\,&\psi\big((u_2u_1 + \mu\bar v_1v_2) + (v_2\bar u_1 + v_1u_2)j\big) \\
=\,\,&\psi\big((u_2 + v_2j)(u_1 + v_1j)\big) = \psi\big((u_1 + v_1j)\cdot(u_2 + v_2j)\big),
\end{align*}
as claimed. Combining we obtain an isomorphism $\vph \circ \psi\:\Cay(B,\mu)^{\op} \to \Cay^\prime(B,\mu)^{\op}$, hence an isomorphism 
\[
\vph \circ \psi\:\Cay(B,\mu) \overset{\sim} \longrightarrow \Cay^\prime(B,\mu), \quad u + vj \longmapsto u + j^\prime\bar v,
\]
of conic algebras.
\end{sol}

\solnsec{Section~\ref{s.BAPRO}}

\begin{sol}{pr.COMPUN} \label{sol.COMPUN}
 (i) $\Rightarrow$ (ii). Since $(C,n,1_C)$ by (\ref{ss.PRECO}.\ref{PERUN}) is a pointed quadratic module, the implication follows from Lemma~\ref{l.UNBA}.

(ii) $\Rightarrow$ (iii). Obvious.

(iii) $\Rightarrow$ (i). By (\ref{ss.PRECO}.\ref{PERCOM}), (\ref{ss.MULCO}.\ref{MUPLR}), the quantity $\varepsilon := n(1_C) \in k$ is an idempotent with $\varepsilon n(x) = n(x)$, $\varepsilon n(x,y) = n(x,y)$ for all $x,y \in C$. Therefore the idempotent $\varepsilon^\prime := 1 - \varepsilon \in k$ satisfies $n(\varepsilon^\prime 1_C) = \varepsilon^\prime \varepsilon = 0$, $n(\varepsilon^\prime 1_C,y) = \varepsilon^\prime \varepsilon
n(1_C,y) = 0$ for all $y \in C$, which implies $\varepsilon^\prime 1_C = 0$ by non-degeneracy of $n$, hence $\varepsilon^\prime = 0$ by (iii). Thus $n(1_C) = 1$.
\end{sol}

\begin{sol}{pr.CONIL} \label{sol.CONIL}
 Let $K$ be a unital $F$-algebra of dimension $2$. Then $K$ is quadratic, hence commutative associative (\ref{ss.QUALG}), and we have $K = F[u]$ for any $u \in K \setminus F1_K$. Writing $\mu_u \in F[\bft]$ for the minimum polynomial of $u$ over $F$, which is monic of degree $2$, we have the following possibilities. 
 
 \step{1}
 $\mu_u$ is irreducible and separable. The $K/F$ is a separable quadratic field extension, and we are in Case~(i). 
 
 \step{2}
 $\mu_u$ is reducible and separable. Then $\mu_u = (\bft - \alpha)(\bft - \beta)$ for some $\alpha,\beta \in F$, $\alpha \neq \beta$. Here the Chinese Remainder Theorem implies
\[
K \cong F[\bft]/(\mu_u) \cong F[\bft]/(\bft - \alpha) \times F[\bft]/(\bft - \beta) \cong F \oplus F 
\]
as a direct product of ideals, and we are in Case~(ii). 

\step{3}
$\mu_u$ is irreducible and inseparable. Then $\ch(F) = 2$ and $K/F$ is an inseparable quadratic field extension. Hence we are in Case~(iii). 

\step{4}
$\mu_u$ is reducible and inseparable. Then $\mu_u = (\bft - \alpha)^2$ for some $\alpha \in F$, which implies $K = F[v]$, where $v := u - \alpha 1_K$ satisfies $v^2 = 0 \neq v$. Hence the map $F[\vep] \to K$ sending $\vep$ to  $v$ is an isomorphism of $F$-algebras, and we are in Case~(iv).

\smallskip

For the second part of the problem, assume that $F$ is perfect of characteristic $2$ and let $C$ be a conic $F$-algebra without nilpotent elements $\neq 0$. It will be enough to prove $\Ker(t_C) = F1_C$, where the inclusion ``$\supseteq$'' is obvious. Conversely, let $u \in \Ker(t_C)$ and assume $u \notin F1_C$. Then $K := F[u] \subseteq C$ is a two-dimensional subalgebra, so by what we have just proved, it satisfies one of the conditions (i)--(iv) above. If (i) or (ii) holds, then the trace of $K$ kills $1_K = 1_C$ but is not identically zero. Thus $\Ker(t_K) = F1_C$, forcing $t_C(u) = t_K(u) \neq 0$, a contradiction. Case~(iii) cannot hold since $F$ is perfect, but neither can Case~(iv) since $C$ does not contain nilpotent elements different from zero. This completes the proof.
\end{sol}

\begin{sol}{pr.COSIM} \label{sol.COSIM}
 Let $I \neq C$ be a two-sided ideal of $(C,\iota_C)$. For $x \in C$ we obtain $n_C(x)1_C = x\bar x \in I$, which implies $n_C(x) = 0$. Moreover, for $y \in C$ we obtain $n_C(x,y)1_C = x\bar y + y\bar x \in I$, which implies $n_C(x,y) = 0$. But $n_C$ is non-degenerate. Thus $x = 0$, and we have shown $I = \{0\}$, so $(C,\iota_C)$ is simple as an algebra with involution. Here Prop.~\ref{p.SIMINV} shows that either $C$ is simple, or $(C,\iota_C) \cong (A \times A^{\op},\vep)$ for some simple unital $F$-algebra $A$, where $\vep$ stands for the exchange involution. Identifying $C = A \times A^{\op}$ by means of this isomorphism and $A,A^{\op} \subseteq C$ canonically, we put $c := 1_A$, $c^\prime := 1_{A^{\op}}$ and conclude from Exc.~\ref{pr.IDEMP} that $c,c^\prime$ are elementary idempotents in $C$ such that $1_C = c + c^\prime$. For $x \in A$, we now obtain $2x = c \circ x = t_C(c)x + t_C(x)c - n_C(c,x)1_C = x + t_C(x)c - n_C(c,x)1_C$, hence
\[
x = \big(t_C(x) - n_C(c,x)\big)c - n_C(c,x)c^\prime.
\]
This shows $n_C(c,x) = 0$ and $x = t_C(x)c \in Fc$, and we have proved $A \cong F$, hence $A^{\op} \cong F$ as well. Thus $C \cong F \times F$ is split quadratic \'etale.
\end{sol}

\begin{sol}{pr.CDDIVALG} \label{sol.CDDIVALG}
We first show that \emph{the norm of a conic divsion algebra over a field $F$ is anisotropic.} Indeed, let $C$ be a conic division algebra over $F$ and $0 \neq x \in C$. Then $F[x] \subseteq C$ is a unital subalgebra of dimension at most $2$, and left multiplication by $x$ gives a linear injection from $F[x]$ to itself, which therefore is a bijection. Hence $x$ is invertible in $F[x]$, which implies $n_C(x) \neq 0$ by Exc.~\ref{pr.NILRADCON}~(b) and proves our claim. For the rest of the proof, we assume that $F$ is a field of arbitrary characteristic, $\mu \in F$ is a non-zero scalar, $B$ is an octonion algebra over $F$ and $C = \Cay(B,\mu) = B \oplus Bj$.

(a) Isotropic regular quadratic forms over $F$ are known to be universal. The assumption $\mu \notin n_B(B^\times)$ therefore implies that $n_B$ is anisotropic, forcing $B$ to be an octonion division algebra by Thm.~\ref{t.COMALT} combined with Prop.~\ref{p.CHIN}. Now (\ref{p.CDCON}.\ref{CDMULT}) implies
\[
x_1x_2 = (u_1u_2 + \mu\bar v_2v_1) + (v_2u_1 + v_1\bar u_2)j,
\]
and we conclude
\begin{align}
\label{ZECRI} x_1x_2 = 0 \quad \Longleftrightarrow \quad u_1u_2 = -\mu\bar v_2v_1,\;\;v_2u_1 = -v_1\bar u_2.
\end{align}
Suppose first $x_1x_2 = 0$. If $u_1 = 0$, then $v_1 \neq 0$ and \eqref{ZECRI} implies $v_2 = u_2 = 0$, hence $x_2 = 0$, a contradiction. Thus $u_1 \neq 0$. If $u_2 = 0$, then \eqref{ZECRI} again implies $v_2 = 0$. This contradiction shows $u_2 \neq 0$. Since $B$ is a division algebra, the first equation of \eqref{ZECRI} now implies $v_1 \neq 0 \neq v_2$. We have thus shown $u_i \neq 0 \neq v_i$ for $i = 1,2$. Next, since the norm of $B$ is multiplicative, \eqref{ZECRI} yields first $n_B(u_1)n_B(u_2) = \mu^2n_B(v_1)n_B(v_2)$ and then
\[
n_B(u_1)^2n_B(u_2) = \mu^2n_B(v_1)n_B(v_2)n_B(u_1) = \mu^2n_B(v_1)^2n_B(u_2).
\]
Thus $n_B(u_1) = \pm\mu n_B(v_1)$. But $n_B(u_1) = \mu n_B(v_1)$ is impossible since this and multiplicativity of $n_B$ would yield the contradiction $\mu \in n_B(B^\times)$. Hence (i) holds. Turning to (ii), we use \eqref{ZECRI}, (i) and Kirmse's identities (\ref{ss.IDCO}.\ref{KID}) to compute
\begin{align*}
(u_1u_2)\bar v_1 =\,\,&-\mu(\bar v_2v_1)\bar v_1 = -\mu n_B(v_1)\bar v_2 = n_B(u_1)\bar v_2 = u_1(\bar u_1\bar v_2) \\ 
=\,\,& -u_1\overline{v_1\bar u_2} = -u_1(u_2\bar v_1),
\end{align*} 
and this is (ii). Finally, (iii) is equivalent to the second equation of \eqref{ZECRI}. Conversely, suppose $u_i \neq 0 \neq v_i$ for $i = 1,2$ and (i)--(iii) hold. Then so does the second equation of \eqref{ZECRI}, and (i), (ii) imply
\begin{align*}
(u_1u_2)\bar v_1 =\,\,&-u_1(u_2\bar v_1) = -u_1\overline{v_1\bar u_2} = u_1(\bar u_1\bar v_2) = n_B(u_1)\bar v_2 \\ 
=\,\,&-\mu n_B(v_1)\bar v_2 = -\mu(\bar v_2v_1)\bar v_1.
\end{align*} 
But $\bar v_1$, being different from zero, is invertible in $B$. Hence the first equation of \eqref{ZECRI} holds as well. This completes the proof of (a).

(b) We know already that if $C$ is a division algebra, then $n_C$ is anisotropic. Conversely, suppose $n_C = n_B \perp (-\mu)n_B$ (cf. Remark~\ref{r.NOCD}) is anisotropic. Assuming $\mu = n_B(u)$ for some $u \in B^\times$ would lead to the contradiction $n_C(u + j) = 0$. Thus $\mu \notin n_B(B^\times)$. If $C$ were not a division algebra, some $x_i = u_i + v_ij$, $u_i,v_i \in B$, $i = 1,2$, would satisfy the conditions of (a). In particular, since we are in characteristic $2$, condition (i) would imply $n_C(x_1) = n_B(u_1) + \mu n_B(v_1) = 0$, in contradiction to $n_C$ being anisotropic. Hence $C$ is a division algebra.

(c) Suppose first $x,y \in A$, $z \in A^\perp$ for some quaternion subalgebra $A \subseteq B$ and $x \circ y = 0$, equivalently, $xy = -yx$. Then we may identify $B = \Cay(A,\nu) = A \perp Ai$ for some non-zero scalar $\nu \in F$ (Thm.~\ref{t.COMSUBALG}~(a)) and have $A^\perp = Ai$, hence $z = ui$ for some $u \in A$. Now we use (\ref{p.CDCON}.\ref{CDMULT}) and compute
\begin{align*}
(xy)z =\,\,&(xy)(ui) = (uxy)i, \\
x(yz) =\,\,&x\big(y(ui)\big) = x\big((uy)i\big) = (uyx)i = -(uxy)i.
\end{align*}
Hence $(xy)z = -x(yz)$, as desired. Conversely, let this be so. Since $\ch(F) \neq 2$ and $B$ is alternative, we have $x \notin F1_B$, $y \notin F[x]$. Let $A$ be the (unital) subalgebra of $B$ generated by $x,y$. Since $B$ is a division algebra, so is $A$, forcing its norm to be anisotropic, hence regular. Since, therefore, $A$ is a composition algebra of dimension at least $3$ and at most $4$ (Exc.~\ref{pr.ARTCONALG}), it must, in fact, be a quaternion algebra. As before, we may assume $B = \Cay(A,\nu) = A \perp Ai$ for some $\nu \in F^\times$. Then $z = u + vi$, $u,v \in A$, and we conclude
\begin{align*}
(xy)z =\,\,&xyu + (vxy)i, \\
-x(yz) =\,\,&-x\big(yu + (vy)i\big) = -xyu - (vyx)i.
\end{align*}
Comparing coefficients, this implies $xyu = v(x \circ y) = 0$, hence $u = 0$. But then $v \neq 0$, forcing $x \circ y = 0$ and $z = vi \in A^\perp$.

(d) Suppose first that $C$ is a division algebra. Then its norm is anisotropic, and we have seen in (b) that this implies $\mu \notin n_B(B^\times)$. Assume now $-\mu = n_B(z)$ for some element $z \in B$ of trace zero. Pick any non-zero element $y \in F[z]^\perp$ and write $A^\prime$ for the subalgebra of $B$ generated by $y,z$. The argument produced in the proof of (c) shows that $A^\prime$ is a quaternion subalgebra of $B$. Pick any non-zero element $x \in A^{\prime\perp} \subseteq F[y]^\perp$. Then the subalgebra $A$ of $B$ generated by $x,y$ is again a quaternion algebra, and we have not only $n_B(z,1_B) = n_B(z,x) = n_B(z,y) = 0$ but also (by norm associativity) $n_B(z,xy) = n_B(z\bar y,x) \in n_B(A^\prime,x) = \{0\}$. Thus $z \in A^\perp$, while (\ref{ss.BASID}.\ref{CAQUADL}) yields $x \circ y = 0$. Setting $u_1 := x$, $u_2 := y$, $v_1 := u_1z^{-1} = - n_B(z)^{-1}u_1z \in A^\perp$, $v_2 := -(v_1\bar u_2)u_1^{-1}$, we have $\bar v_1 = -v_1 \in A^\perp$, whence (c) implies that conditions (i)--(iii) of (a) hold. Hence (a) produces zero divisors in $C$, a contradiction. Conversely, suppose $\mu \notin n_B(B^\times)$, $-\mu \notin n_B(B^0 \cap B^\times)$. If $C$ were not a division algebra, we would find elements $u_i,v_i \in B$, $i = 1,2$ satisfying the conditions of (a). Then (c) would lead to a quaternion subalgebra $A \subseteq B$ such that $u_1,u_2 \in A$, $\bar v_1 \in A^\perp$ and $u_1 \circ u_2 = 0$. This would imply $\bar v_1 = -v_1 \in A^\perp$ and $-\mu = n_B(z)$, $z := u_1v_1^{-1} \in A^\perp$ and, in particular, $t_B(z) = 0$. This contradiction completes the proof. 

(e) The norm of $\IO$ is positive definite, and so is its restriction to the trace zero elements. Hence $1 = -(-1)$ cannot avoid being the norm of a trace zero element of $\IO$.
\end{sol} 

\begin{sol}{pr.FROBTHE} \label{sol.FROBTHE}
Let $C$ be a finite-dimensional alternative division algebra over $\IR$. Then Exc.~\ref{pr.UNALT} shows that $C$ is unital. Moreover, for $0 \neq x \in C$, the unital subalgebra $\IR[x] \subseteq C$ is finite-dimensional, commutative associative and without zero divisors. Thus $\IR[x]/\IR$ is a finite algebraic field extension, which implies $\IR[x] = \IR1_C$ or $\IR[x] \cong \IC$ as $\IR$-algebras. In particular, the elements $1_C,x,x^2$ are linearly dependent over $\IR$, and we deduce from Exc.~\ref{pr.CONFIELD} that $C$ is a conic algebra. Since the non-zero elements of $C$ are invertible, Prop.~\ref{p.CHIN} shows that the norm of $C$ is anisotropic. On the other hand, it represents $1$ and hence (by the intermediate value theorem) must be positive definite. Now Thm.~\ref{t.COMALT} shows that $C$ is a composition division algebra over $\IR$. By Cor.~\ref{c.ENUM}, therefore, $C \cong \Cay(\IR,\alpha_1,\dots,\alpha_n)$ for some $n = 0,1,2,3$ and $\alpha_1,\dots,\alpha_n \in \IR^\times$. Here the assumption $\alpha_i > 0$ for some $i = 0,\dots,n$ would contradict Remark~\ref{r.NOCD} in conjunction with the property of $n_C$ to be positive definite. Thus $\alpha_i < 0$ for $0 \leq i \leq n$, and invoking Exc.~\ref{pr.CDSCAPAR}, we are actually reduced to the case $\alpha_i = -1$, $0 \leq i \leq n$, which means $C \cong \IR,\IC,\IH,\IO$ for $n = 0,1,2,3$, respectively.
\end{sol}

\begin{sol}{pr.CDISOT} \label{sol.ISOT}
For the first part of the problem, put $C := \Cay (B,\mu)$ and write $\vph\:C \to C^p$ for the map defined by 
\begin{align*}
\vph (u + vj) := (p^{-1}up) + vj &&(u,v \in B),
\end{align*} which is obviously a linear bijection preserving units but also norms since $B$ was assumed to be multiplicative. It therefore suffices to show that $\vph$ is an algebra homomorphism. In order to do so, write ``$\cdot$'' for the product in $C^p$ and let $u_1,u_2,v_1,v_2 \in B$. Combining (\ref{p.CDCON}.\ref{CDMULT}) and associativity of $B$ with the fact that its conjugation is an involution (Prop.~\ref{p.MULAD}~(a)), we obtain
\begin{align*}
\vph(u_1 + v_1j)\cdot\vph(u_2 + v_j) =\,\,&\big((p^{-1}u_1p) + v_1j\big)\cdot\big((p^{-1}u_2p) + v_2j\big) \\
=\,\,&\Big(\big((p^{-1}u_1p) + v_1j\big)p^{-1}\Big)\Big(p\big((p^{-1}u_2p) + v_2j\big)\Big) \\
=\,\,&\big(p^{-1}u_1 + (v_1\bar p^{-1})j\big)\big(u_2p + (v_2p)j\big) \\
=\,\,&(p^{-1}u_1u_2p + \mu\bar p\bar v_2v_1\bar p^{-1}) + (v_2pp^{-1}u_1 + v_1\bar p^{-1}\bar p \bar u_2)j. 
\end{align*}
But $\bar p = n_B(p)p^{-1}$, $\bar p^{-1} = n_B(p)^{-1}p$. Hence
\begin{align*}
\vph(u_1 + v_1j)\cdot\vph(u_2 + v_j) =\,\,&\big(p^{-1}(u_1u_2 + \mu\bar v_2v_1)p\big) + (v_2u_1 + v_1\bar u_2)j \\
=\,\,&\vph\big((u_1u_2 + \mu\bar v_2v_1) + (v_2u_1 + v_1\bar u_2)j\big) \\
=\,\,&\vph\big((u_1 + v_1j)(u_2 + v_2j)\big),
\end{align*}
as desired.

Turning to the second part of the problem, let $C$ be an octonion algebra over $k$ and suppose $p,q \in C^\times$ have the property that $pq^2$ belongs to some quaternion subalgebra $B \subseteq C$. By Exc.~\ref{pr.UNITARBISOTALT}~(a) and its solution, it will be enough to show for $p \in B^\times$ and $B$ as above that $C$ and $C^p$ are isomorphic. In order to do so, we note $C = B \oplus B^\perp$ as a direct sum of submodules (Lemma~\ref{l.NOSU}) and define $\vph\:C \to C^p$ by $\vph(u + v) := (p^{-1}up) + v$ for $u \in B$, $v \in B^\perp$. We claim that $\vph$ is an isomorphism of conic algebras. This assertion is local on $k$, so we may assume that $k$ is a local ring. By Thm.~\ref{t.COMSUBALG}~(a), we may identify $C = \Cay(B,\mu) = B \oplus Bj$, where the relation $B^\perp = Bj$ can be read off from (\ref{p.CDCON}.\ref{CDNORL}) and the regularity of $n_B$. Thus our new map $\vph$ agrees with the one defined in the first part and hence must be an isomorphism of conic algebras.
\end{sol}

\begin{sol}{pr.AZUQUAD} \label{sol.AZUQUAD}
(a) We begin by showing that \emph{for a quadratic $k$-algebra $R$ and $u \in R$,}
\begin{align}
\label{UMIBA} n_R(u - \bar u) = 4n_R(u) - t_R(u)^2,
\end{align}
which follows from
\[
n_R(u - \bar u) = n_R(u) - n_R(u,\bar u) + n_R(\bar u) = 2n_R(u) - t_R(u)^2 + n_R(u,u) = 4n_R(u) - t_R(u)^2.
\]
Assume first that $R$ is \'etale and $u$ generates $R$ as a $k$-algebra. Combining Prop.~\ref{p.SMASU} with \eqref{UMIBA}, we conclude that $u - \bar u$ is invertible. Conversely, suppose that $u - \bar u$ is invertible in $R$. Again by Prop.~\ref{p.SMASU} and \eqref{UMIBA}, $D :=k[u] \subseteq R$ is a quadratic \'etale subalgebra. By Lemma~\ref{l.NOSU}, therefore, we obtain the decomposition $R = D \oplus D^\perp$ as $k$-modules and conclude, by comparing ranks, that $R = D$ is quadratic \'etale. 

Next suppose $k$ is an LG ring and $D$ is a quadratic \'etale $k$-algebra. By Prop.~\ref{p.PROSEM}, $D$ is a free $k$-module of rank $2$ containing $1_D$ as a unimodular element (\ref{ss.QUALG}), which therefore can be extended to a basis $(1_D,u)$ of $D$. Thus $D = k[u]$, and Prop.~\ref{p.SMASU} combined with \eqref{UMIBA} shows $n_D(u - \bar u) \in k^\times$, hence $u - \bar u \in D^\times$.

Finally, returning to an arbitrary base ring $k$ and a quadratic \'etale $k$-algebra $D$, it remains to show $H(D,\iota_D) = k1_D$. In order to do so, we may assume $D = k[u]$ for some $u \in D$ such that $u - \bar u \in D^\times$. In particular, by Prop.~\ref{p.SMASU}, $D$ is free (of rank $2$) as a $k$-module with basis $1_D,u$. We clearly have $k1_D \subseteq H(D,\iota_D)$. Conversely, write $x \in H(D,\iota_D)$ as $x = \alpha 1_D + \beta u$ for some $\alpha,\beta \in k$. Then $x = \bar x = \alpha 1_D + \beta\bar u$ implies $\beta(u - \bar u) = 0$, hence $\beta = 0$ (since invertible elements, being the image of $1_D$ under a linear bijection, are unimodular), and we end up with $x = \alpha 1_D \in k1_D$.

(b) Let $B$ be a quaternion algebra over $k$. We must show that an element $x \in B$ satisfying $[x,y] = 0$ for all $y \in B$ is a scalar multiple of $1_C$. This assertion being local on $k$, we may assume that $k$ is a local ring. By Cor.~\ref{c.ENUM}, we may further assume $B = \Cay(D,\mu) = B \oplus Bj$ for some quadratic \'etale $k$-algebra $D$ and some $\mu \in k$. From Thm.~\ref{t.COALCD} we deduce $\mu \in k^\times$, while part (a) of our problem yields an element $u \in D$ making $u - \bar u$ invertible. Now write $x = u_1 + v_1j$, $u_1,v_1 \in D$. By hypothesis we have $[x,u_2 + v_2j] = 0$ for all $u_2,v_2 \in D$, which by (\ref{ss.CDCOMASS}.\ref{CDCOMONE}) implies
\[
0 = [u_1,u_2] + \mu(\bar v_2v_1 - \bar v_1v_2) = \mu(\bar v_2v_1 - \bar v_1v_2).
\] 
Since $\mu$ is invertible, we conclude $\bar v_2v_1 = \bar v_1v_2$ for all $v_2 \in D$. Setting $v_2 = 1_B$ gives $\bar v_1 = v_1$, while setting $v_2 = u$ now gives $(u - \bar u)v_1 = 0$, hence $v_1 = 0$ since $u - \bar u$ is invertible. Next applying (\ref{ss.CDCOMASS}.\ref{CDCOMTWO}), we obtain $0 = v_2(u_1 - \bar u_1) - v_1(u_2 - \bar u_2) = v_2(u_1 - \bar u_1)$ for all $v_2 \in B$, which for $v_2 = 1_B$ amounts to $u_1 = \bar u_1$, hence by (a) to $u_1 \in k1_B$. Thus we have $x \in k1_C$, as desired.

Finally, let $C$ be an octonion algebra over $k$. We have
\[
k1_C \subseteq \{x \in C\mid \forall y \in C:\;xy = yx\} \cap \Nuc(C),
\]
so it remains to show that both sets taking part in the intersection on the right-hand side belong to $k1_C$. Since these are local questions, we may assume that $k$ is a local ring and then obtain, by Cor.~\ref{c.ENUM}, $C = \Cay(B,\mu) = B \oplus Bj$ up to isomorphism, for some quaternion subalgebra $B \subseteq C$ and some scalar $\mu \in k^\times$ (Thm.~\ref{t.COALCD}). Now let $x = u_1 + v_1j \in C$ with $u_1,v_1 \in B$ and first assume $[x,u_2 + v_2j] = 0$ for all $u_2,v_2 \in B$. From (\ref{ss.CDCOMASS}.\ref{CDCOMONE}), (\ref{ss.CDCOMASS}.\ref{CDCOMTWO}), we conclude
\begin{align}
\label{OCOM} [u_1,u_2] + \mu(\bar v_2v_1 - \bar v_1v_2) = 0 = v_2(u_1 - \bar u_1) - v_1(\bar u_2 - u_2)
\end{align} 
for all $u_2,v_2 \in B$. By what we have seen before, the first of these equations for $v_2 = 0$ shows $u_1 \in \Cent(B) = k1_B$. Now the second one implies $v_1(u_2 - \bar u_2) = 0$ for all $u_2 \in B$. By Thm.~\ref{t.COMSUBALG}~(b), we find a quadratic \'etale subalgebra $D \subseteq B$, which in turn, by part (a) of  this exercise, contains an element $u_2$ making $u_2 - \bar u_2$ invertible. Hence $v_1 = 0$ and thus $x \in k1_C$. Next assume $x \in \Nuc(C)$. Applying (\ref{ss.CDCOMASS}.\ref{CDASSONE}),  for $u_2 = v_3 = 0$, we obtain $(\bar v_2v_1)u_3 = (u_3\bar v_2)v_1$ for all $v_2,u_3 \in B$. Setting $v_2 = 1_B$, this implies $v_1u_3 = u_3v_1$ for all $u_3 \in B$, hence $v_1 \in \Cent(B) = k1_B$. Thus $v_1 = \beta 1_B$ for some $\beta \in k$, and our original equation reduces to $(\beta\bar v_2)u_3 = u_3(\beta\bar v_2)$ for all $v_2,u_3 \in B$ This means $\beta v_2 \in \Cent(B) = k 1_B$ for all $v_2 \in B$. But $1_B$, being unimodular, may be extended to basis of $B$ as a $k$-module, and picking for $v_2$ one of the basis vectors distinct from $1_B$, we conclude $\beta = 0$, hence $v_1 = 0$. Now we consider (\ref{ss.CDCOMASS}.\ref{CDASSTWO}) for $u_2 = u_3 = v_2 = 0$ and $v_3 = 1_B$. Then $[u_1,u_2] = 0$ for all $u_2 \in B$, which yields $u_1 \in \Cent(B) = k1_B$, hence $x \in k1_C$, and completes the proof.
\end{sol}

\begin{sol}{pr.ETAUT} \label{sol.ETAUT}
(a) Assume first $2 \in k^\times$. Then $k1_D \subseteq D$ is a regular composition subalgebra, and Thm.~\ref{t.COMSUBALG} yields an identification $D = \Cay(k,\mu) = k \oplus kj$, for some $\mu \in k^\times$. By (\ref{p.CDCON}.\ref{CDNOR}), therefore,  $u := \frac{1}{2}\cdot 1_D + j$ does the job. We are left with the case $2 \in \mfm$, where $\mfm$ stands for the maximal ideal of $k$. Since all finitely generated projectives over $k$ are free and $1_D \in D$ is unimodular, some $v \in D$ makes $(1_D,v)$ a basis of $D$ over $k$. Now Prop.~\ref{p.SMASU} implies $\mu := t_D(v)^2 - 4n_D(v) \in k^\times$. Since $2 \in \mfm$, we conclude $t_D(v) \in k^\times$ and $u := t_D(v)^{-1}v$ does the job.

(b) If $k$ allows a decomposition of the desired kind, then $\vph$ is clearly an automorphism of $D$. Conversely, let this be so. We first treat the case that $k$ is a local ring. By (a), $D = k[u]$ for some element $u \in D$ of trace $1$. By Prop.~\ref{p.SMASU}, therefore, $(1_D,u)$ is a basis of $D$, and there are $\alpha,\beta \in k$ such that $\vph(u) = \alpha 1_D + \beta u$. But $\vph$, being an automorphism of $D$, preserves norms, traces and unit. Hence
\begin{align*}
2\alpha + \beta =\,\,& t_D\big(\vph(u)\big) = t_D(u) = 1, \\
\alpha^2 + \alpha\beta + \beta^2n_D(u) =\,\,&n_D\big(\vph (u)\big) = n_D(u),
\end{align*}
and we conclude $\beta = 1 - 2\alpha$,
\[
\alpha^2 + \alpha - 2\alpha^2 + (1 - 4\alpha + 4\alpha^2)n_D(u) = n_D(u).
\]
Thus $(\alpha - \alpha^2)(1 - 4n_D(u)) = 0$. On the other hand, since $D$ is quadratic \'etale, Prop.~\ref{p.SMASU} shows $1 - 4n_D(u) \in k^\times$. Hence $\alpha \in k$ is an idempotent. But as a local ring, $k$ is connected. We therefore conclude $\alpha = 0$, $\beta = 1$ or $\alpha = 1$, $\beta = -1$, which implies $\vph = \Eins_D$ in the first case, $\vph = \iota_D$ in the second.

Now let $k$ be arbitrary  and put $X := \Spec(k)$. Since $D$ is finitely generated as a $k$-module,
\[
X_+ := \{\mfp \in X \mid \vph_\mfp = \Eins_{D_\mfp}\}, \quad X_- := \{\mfp \in X \mid \vph_\mfp = \iota_{D_\mfp}\}
\]
by Exc.~\ref{pr.LOCHOM}, are Zariski-open subsets of $X$. Moreover, they are disjoint by Exc.~\ref{pr.AZUQUAD}~(a) and cover $X$ by the special case just treated. Hence Exc.~\ref{pr.IDEPAR}, yields a complete orthogonal system $(\vep_+,\vep_-)$ of idempotents in $k$ such that $X_\pm = D(\vep_\pm)$. Put $k_\pm := \vep_\pm k$. Then $k = k_+ \times k_-$ as a direct product of ideals, and with the canonical projection $\pi_+\:k \to k_+$, we consult (\ref{ss.PROPID}.\ref{SPEPI}) to conclude $\mfp := \Spec(\pi_+)(\mfp_+) = \mfp_+ \times k_- \in D(\vep_+)$ for any prime ideal $\mfp_+ \subseteq k_+$. Now \ref{ss.PROPSP}, \ref{ss.PROPID} show that $\vph_{+\mfp_+} = \vph_\mfp \otimes_{k_\mfp} k_{+\mfp_+} = \Eins_{D_{+\mfp_+}}$. Thus $\vph_+ = \Eins_{D_+}$ and, similarly, $\vph_- = \iota_{D_-}$.
\end{sol} 

\begin{sol}{pr.IDEALSCOAL} \label{sol.IDEALSCOAL} See \cite[pp. 290-291]{PeMS} 
\end{sol}

\begin{sol}{pr.PEIRCEELID} \label{sol.PEIRCEELID}
 We begin by proving \eqref{ONEONEPEIR} and first note that not only $c_1 = c$ but also $c_2 = 1_C - c$ is an elementary idempotent (cf. Exc.~\ref{pr.IDEMP}, particularly (iv)). Let $ i = 1,2$ and $x \in C_{ii}$. Then $c_ix = x = xc_i$, and (\ref{ss.IDCO}.\ref{QAUOP}) yields $x = c_ixc_i = n_C(c_i,\bar x)c_i - n_C(c_i)\bar x = n_C(c_i,\bar x)c_i \in kc_i$, and \eqref{ONEONEPEIR} holds. In view of this, arbitrary elements $x,y \in C$ can be written as in \eqref{PEXY}. Expanding the norm of $x$ in the obvious manner, we obtain
\begin{align*}
n_C(x) =\,\,&n_C(\alpha_1c_1 + x_{12} + x_{21} + \alpha_2c_2) \\
=\,\,&\alpha_1^2n_C(c_1) + \alpha_1n_C(c_1,x_{12}) + \alpha_1n_C(c_1,x_{21}) + \alpha_1\alpha_2n_C(c_1,c_2) + n_C(x_{12}) \\
\,\,&+ n_C(x_{12},x_{21}) + \alpha_2n_C(x_{12},c_2) + n_C(x_{21}) + \alpha_2n_C(x_{21},c_2) + \alpha_2^2n_C(c_2).
\end{align*} 
Applying (\ref{ss.MULCO}.\ref{MUPLR}), (\ref{ss.MULCO}.\ref{MUPLL}) and observing $n_C(c_1) = n_C(c_2) = 0$, we obtain, for $\{i,j\} = \{1,2\}$,
\begin{align*}
n_C(c_i) =\,\,&0, \\
n_C(c_i,x_{ij}) =\,\,&n_C(c_ic_i,c_ix_{ij}) = n_C(c_i)n_C(c_i,x_{ij}) = 0, \\
n_C(c_i,x_{ji}) =\,\,&n_C(c_ic_i)n_C(x_{ji}c_i) = n_C(c_i,x_{ji})n_C(c_i) = 0, \\
n_C(c_1,c_2) =\,\,&n_C(c,1_C - c) = t_C(c) - 2n_C(c) = 1, \\
n_C(x_{ij}) =\,\,&n_C(c_ix_{ij}) = n_C(c_i)n_C(x_{ij}) = 0. 
\end{align*}
Inserting these relations into the preceding equation, we end up with \eqref{NORPEIR}, while linearizing \eqref{NORPEIR} gives \eqref{NORPEIRLIN}, which for $y = 1_C$ yields \eqref{TRPEIR}. Finally, \eqref{CONJPEIR} follows immediately from the definition. 

Now let $C$ be a composition algebra. We must show for $\{i,j\} = \{1,2\}$ that $Dn_C$ determines a duality between $C_{ij}$ and $C_{ji}^\ast$. Let $x_{ij} \in C_{ij}$ and suppose $n_C(x_{ij},C_{ji}) = \{0\}$. Then \eqref{NORPEIRLIN} yields
\[
n_C(x_{ij},C) = n_C(x_{ij},kc_i + C_{ij} + C_{ji} + kc_j) = n_C(x_{ij},C_{ji}) = \{0\}
\]	
hence $x_{ij} = 0$ since $n_C$ is non-singular. Next let $\lambda_{ji}\:C_{ji} \to k$ be a linear form and $\lambda\:C \to k$ the linear extension of $\lambda_{ji}$ that kills $kc_1 \oplus C_{ij} \oplus kc_2$. By non-singularity of $n_C$, there exists $x \in C$ such that $\lambda = n_C(x,\emptyslot)$. Writing $x$ in the form \eqref{PEXY} and choosing $y_{ji} \in C_{ji}$, we apply \eqref{NORPEIRLIN} and obtain $\lambda_{ji}(y_{ji}) = \lambda(y_{ji}) = n_C(x,y_{ji}) = n_C(x_{ij},y_{ji})$. Summing up, $Dn_C$ does indeed determine a duality between $C_{ij}$ and $C_{ji}$. In order to prove the final statement, we must show $t_C(x^2z) = t_C(xy^2) = 0$ for all $x,y,z \in C_{ij}$. This follows from the final statement in 
Exercise~\ref{pr.PEIRCEALT}.
\end{sol}

\begin{sol}{pr.MISPLET} \label{sol.MISPLET}
(a) Since scalar multiplication by $D$ acts on the second factor of $D \otimes D$, the map $\vph$ is a unital homomorphism of $D$-algebras. Moreover, for $u \in D$ we have
\begin{align*}
\vph(u \otimes 1_D - 1_D \otimes \bar u) =\,\,&(u,\bar u) - (\bar u,\bar u) = (u - \bar u,0), \\
\vph(1_D \otimes u - u \otimes 1_D) =\,\,&(u,u) - (u,\bar u) = (0,u - \bar u).
\end{align*}
In order to prove that $\vph$ is an isomorphism, we may assume that $k$ is a local ring. By Exc.~\ref{pr.AZUQUAD}~(a), the element $u$ may be so chosen that $u - \bar u \in D^\times$. The preceding equations therefore show that $(1_D,0)$ and $(0,1_D)$ both belong to $\Img(\vph)$. But $D \times D$ is generated by these elements as a $D$-algebra. Hence $\vph$ is surjective. On the other hand, $D \otimes D$ and $D \times D$ are finitely generated projective $D$-modules of rank $2$. Hence $\vph$ is an isomorphism by Exc.~\ref{pr.SURPRO}.

(b) By definition, we have $\sigma^\prime \circ \vph = \vph \circ \sigma^\prime$. Hence, for $x,d \in D$,
\begin{align*}
\sigma^\prime\big((xd,\bar xd)\big) =\,\,&\sigma^\prime \big(\vph(x \otimes d)\big) = \vph\big(\sigma(x \otimes d)\big) = \vph(x \otimes \bar d) =(x\bar d,\bar x\bar d) = (\overline{\bar xd},\overline{xd}),
\end{align*}  
which proves the assertion since, by (a), any element of $D \times D$ can be written as a finite sum of expressions of the form $(xd,\bar xd)$ with $x,d \in D$.
\end{sol}

\begin{sol}{pr.QUACOMQUA} \label{sol.QUACOMQUA}
(a) Let $x,y \in R$. Since $R$ by Prop.~\ref{p.MULAD}~(b) is a multiplicative conic algebra, we may apply (\ref{p.NORAS}.\ref{ASSBILT}), (\ref{ss.BASID}.\ref{CANOCO}) and (\ref{ss.MULCO}.\ref{MUPLR}) to obtain
\begin{align*}
t_R\big(e(xy)\big) =\,\,&t_R\big((ex)(ey)\big) = t_R(ex)t_R(ey) - n_R(ex,ey) \\ 
=\,\,&t_R(ex)t_R(ey)- n_R(e)n_R(x,y) = t_R(ex)t_R(ey).
\end{align*}
Hence $t_R \circ L_e\:R \to k$ is an algebra homomorphism, which by Exc.~\ref{pr.IDEMP} is unital if and only if $e$ is an elementary idempotent.

(b) For $x,y \in R$, we combine (a) with the multiplicativity of $n_R$ and obtain
\begin{align*}
q(x)q(y) =\,\,&\big(t_R(ex^2) + \vep n_R(x)\big)\big(t_R(ey^2) + \vep n_R(y)\big) \\
=\,\,&t_R(ex^2)t_R(ey^2) + \vep t_R(ex^2)n_R(y) + \vep n_R(x)t_R(ey^2) + \vep^2n_R(xn_R(y) \\
=\,\,&t_R\big(e(xy)^2\big) + t_R(\vep ex^2)n_R(y) + n_R(x)t_R(\vep ey^2) + \vep n_R(xy) \\
=\,\,&t_R\big(e(xy)^2\big) + \vep n_R(xy) = q(xy).
\end{align*}
Hence $q$ permits composition.

(c) The condition is clearly sufficient. Conversely, suppose $q$ permits composition. We put $e_1 := (1,0)$, $e_2 := (0,1)$ and
\[
\vep_1 := q(e_1), \quad \vep_2 := q(e_1,e_2), \quad \vep_3 := q(e_2).
\]
Then \eqref{MAQU} holds. Moreover, since $q$ permits composition, $(\vep_1,\vep_3)$ is an orthogonal system of idempotents such that $\vep_1\vep_2 = q(e_1)q(e_1,e_2) = q(e_1^2,e_1e_2) = 0$ and, similarly, $\vep_3\vep_2 = 0$. On the other hand, $\vep := q(1_D) = \sum\vep_i$ is an idempotent in $k$ satisfying $q(x) = \vep q(x)$, $q(x,y) = \vep q(x,y)$ for all $x,y \in D$. In particular, we have $\vep\vep_i = \vep_i$ for $i = 1,2,3$, whence $\vep_2 = \vep - \vep_1 - \vep_3$ must be an idempotent as well.

(d) Let $q\:D \to k$ be a quadratic form permitting composition. Since, therefore, its scalar extension $q_D\:D \otimes D \to D$ permits composition, so does $q^\prime := q_D \circ \vph^{-1}\:D \times D \to D$, where $\vph$ is the isomorphism of Exc.~\ref{pr.MISPLET}~(a). Hence (c) yields an orthogonal system $(d_1,d_2,d_3)$ of idempotents in $D$ such that
\[
q\big((a,b)\big) = d_1a^2 + d_2ab + d_3b^2
\]
for all $a,b \in D$. With $\sigma$ as in Exc.~\ref{pr.MISPLET}~(b), we now claim
\begin{align}
\label{QUDE} q_D \circ \sigma = \iota_D \circ q_D. 
\end{align}
In order to prove this, we let $x,x^\prime,d,d^\prime \in D$ and compute
\begin{align*}
(q_D \circ \sigma)(x \otimes d) =\,\,&q_D(x \otimes \bar d) = q(x)\bar d^2, \\
(\iota_D \circ q_D)(x \otimes d) =\,\,&\overline{q_D(x \otimes d)} =\overline{q(x)d^2} = q(x)\bar d^2, \\
\big(q_D \circ (\sigma \times \sigma)\big)(x \otimes d,x^\prime \otimes d^\prime) =\,\,&q_D(x \otimes \bar d,x^\prime \otimes \bar d^\prime) = q(x,x^\prime)\bar d\bar d^\prime, \\
(\iota_D \circ q_D)(x \otimes d,x^\prime \otimes d^\prime) =\,\,&\overline{q(x,x^\prime)dd^\prime} = q(x,x^\prime)\bar dd^\prime.
\end{align*}
Hence \eqref{QUDE} holds. Using this, and keeping the notation of Exc.~\ref{pr.MISPLET}~(b), we conclude
\[
q^\prime \circ \sigma^\prime = q_D \circ \vph^{-1} \circ \vph \circ \sigma \circ \vph^{-1} = q_D \circ \sigma \circ \vph^{-1} = \iota_D \circ q_D \circ \vph^{-1} = \iota_D \circ q^\prime.
\]
deviating slightly from the notation used in (c), we now put $e_1 := (1_D,0)$, $e_2 := (0,1_D$ and obtain
\[
\bar d_1 = \overline{q^\prime(e_1)} = (\iota_D \circ q^\prime)(e_1) = 
\]
hence $\bar d_3 = d_1$. Similarly,
\[
\bar d_2 = (\iota_D \circ q^\prime)(e_1,e_2) = q^\prime\big(\sigma^\prime(e_1),\sigma^\prime(e_2)\big) = q^\prime (e_2,e_1) = d_2. 
\]
Thus $d_2 \in H(D,\iota_D) = k1_D$ (by Exc.~\ref{pr.AZUQUAD}~(a)). Since $1_D$ is unimodular, we find a unique idempotent $\vep \in k$ such that $d_2 = \vep 1_D$. Then $e := d_1$ is an idempotent in $D$ such that $\vep e = d_2d_1 = 0$, while the relation $n_D(e)1_D = e\bar e = d_1d_3 = 0$ implies $n_D(e) = 0$. Finally, for $x \in D$, we compute
\begin{align*}
q(x)1_D =\,\,&q_D(x \otimes 1_D) = (q^\prime \circ \vph)(x \otimes 1_D) = q^\prime(x \oplus \bar x) = ex^2 + \vep x\bar x + \bar e\bar x^2 \\
=\,\,&\big(t_D(ex^2) + \vep n_D(x)\big)1_D,
\end{align*}
and the assertion follows.
\end{sol}

\begin{sol}{pr.COVQUATALG} \label{sol.COVQUATALG}
We begin by reducing to the case that $k$ is a field. Since all finitely generated projectives over $k$ are free and $1_C \in R$ is unimodular, it can be extended to a basis $(1_C,v)$ of $R$ as a $k$-module. Moreover, the $k$-module $k \times k \times C$ is free, so the polynomial law $f\:k \times k \times C \to k$ given by the set maps
\[
f_S\:S \times S \times C_S \longrightarrow S, (s_1,s_2,u) \longmapsto s_2\det\Big(\big(n_{C_S}(u_i,u_j)\big)_{0\le i,j\le 3}\Big)
\]
for $S \in \kalg$, where $u_0 := 1_{C_S}$, $u_1 := s_11_{C_S} + s_2v_S$, $u_2 := u$, $u_3 := u_1u_2$, may be regarded as a polynomial $f \in k[\bft_1,\dots,\bft_{10}]$ (Cor.~\ref{c.FREPOLA}). This polynomial represents a unit over $k$ if and only if, for some $\xi_1,\xi_2 \in k$, $u \in C$, both $\xi_2$ and
$\det(((n_C(u_i,u_j))_{0\le i,j\le 3})$ with $u_0 := 1_C$, $u_1 = \xi_11_C + \xi_2v$, $u_2 := u$, $u_3 = u_1u_2$ belong to $k^\times$, equivalently, $R = k[u_1]$ and $B := \sum_{i=0}^3 ku_i$ is a quaternion subalgebra of of $C$ (Exc.~\ref{pr.ARTCONALG}) (which, in addition, contains $R$). Hence the proof will be complete once we have shown that $f$ represents a unit over $k$. But $k$ is LG, so it suffices to show that $F$ represents a unit over every field in $\kalg$. This completes the reduction to the field case.

From now on, therefore, let $k = F$ be a field. Then $R$ is a two-dimensional $F$-algebra. If $R$ is quadratic \'etale, then Thm.~\ref{t.COMSUBALG}~(a) yields a scalar $\mu \in k^\times$ such that the inclusion $R \hookrightarrow C$ extends to an embedding $\vph\:\Cay(R,\mu) = R \oplus Rj \to C$. Its image, therefore, is a quaternion subalgebra of $C$ generated by $R$ and $\vph(j)$ in $C$, and we conclude $f(0,1,\vph(j)) \ne 0$. Hence we may assume that $R$ is not quadratic \'etale. Then Exc.~\ref{pr.CONIL} implies that $R/F$ is either an inseparable quadratic field extension of characteristic $2$ or the $F$-algebra of dual numbers. In any event, setting $u_0 := 1_C$, there is an element $u_1 \in C$ such that $R = Fu_0 \oplus Fu_1$ and $n_C(u_0,u_1) = t_C(u_1) = 2\alpha = 0$, where $\alpha := n_C(u_1)$. Since $u_0,u_1$ are linearly independent, any $\beta \in k$ yields a linear form $\lambda$ on $C$ having $\lambda(u_0) = 1$, $\lambda(u_1) = \beta$, and by regularity of $n_C$ we obtain $\lambda = n_C(\emptyslot,u_2)$ for some $u_2 \in C$. Thus $n_C(u_0,u_2) = t_C(u_2) = 1$ and $n_C(u_1,u_2) = \beta$. Setting $\gamma := n_C(u_2)$ and $u_3 := u_1u_2$, we conclude
\begin{align*}
n_C(u_0,u_3) =\,\,&n_C(1_C,u_1u_2) = n_C(\bar u_1,u_2) = -n_C(u_1,u_2) = -\beta, \\
n_C(u_1,u_1) =\,\,&2\alpha = 0, \\
n_C(u_1,u_3) =\,\,&n_C(u_1,u_1u_2) = n_C(u_1)t_C(u_2) = \alpha, \\
n_C(u_2,u_2) =\,\,&2n_C(u_2) = 2\gamma, \\
n_C(u_2,u_3) =\,\,&n_C(u_2,u_1u_2) = n_C(u_2)t_C(u_1) = 0, \\
n_C(u_3,u_3) =\,\,&2n_C(u_1u_2) = 2n_C(u_1)n_C(u_2) = 2\alpha n_C(u_2) = 0. 
\end{align*}
This implies
\begin{align*}
\det\Big(\big(n_C(u_i,u_j)\big)_{0\leq i,j\leq 3}\Big) =\,\,&\det(\left(
\begin{matrix}
2 & 0 & 1 & -\beta \\
0 & 0 & \beta & \alpha \\
1 & \beta & 2\gamma & 0 \\
-\beta & \alpha & 0 & 0
\end{matrix}
\right)) \\ 
=\,\,&\det(\left(
\begin{matrix}
0 & -2\beta & 1 - 4\gamma & -\beta \\
0 & 0 & \beta & \alpha \\
1 & \beta & 2\gamma & 0 \\
0 & \alpha + \beta^2 & 2\beta\gamma & 0
\end{matrix}
\right)) \\
=\,\,&\det(\left(
\begin{matrix}
-2\beta & 1 - 4\gamma & -\beta \\0 & \beta & \alpha \\
\alpha + \beta^2 & 2\beta\gamma & 0
\end{matrix}
\right)) \\
=\,\,&\alpha(\alpha + \beta^2)(1 - 4\gamma) + (\alpha + \beta^2)\beta^2 + 4\alpha\beta^2\gamma \\
=\,\,&(\alpha + \beta^2)(\alpha + \beta^2) = (\alpha + \beta^2)^2
\end{align*}
since $2\alpha = 0$. Setting $\beta = 1$ (resp. $\beta = 0$) for $\alpha = 0$ (resp. $\alpha \neq 0$), the above determinant will be different from zero. Writing $u_1 = \xi_11_C + \xi_2v$ for some $\xi_1,\xi_2 \in F$, we conclude $\xi_2 \ne 0$ (since $u_0,u_1$ are linearly independent), hence $f(\xi_1,\xi_2,u_2) \ne 0$, as desired. 
\end{sol}

\begin{sol}{pr.RINV} \label{sol.RINV}
(a) Let $\sigma$ be a reflection of $C$ and write $B^\pm := \{x \in C\mid \sigma(x) = \pm x\}$. Then $B^+ = \Fix(\sigma) \subseteq C$ is a unital subalgebra and since $\frac{1}{2} \in k$, we have $C = B^+ \oplus B^-$ as a direct sum of submodules. In fact, this sum is even an orthogonal one since $\sigma$ preserves $n_C$ (Exc.~\ref{pr.CONALGHOMNOR}), so for $x^\pm \in B^\pm$ we obtain $n_C(x^+,x^-) = n_C(\sigma(x^+),\sigma(x^-)) = -n_C(x^+,x^-)$, forcing $n_C(x^+,x^-) = 0$. It follows that $\Fix(\sigma) = B^+ \subseteq C$ is a composition subalgebra. In order to prove that its rank is $\frac{r}{2}$, we may assume that $k$ is a local ring. Then Thm.~\ref{t.COMSUBALG}~(a) yields an element $\mu \in k^\times$ such that $\Cay(B^+,\mu) = B^+ \oplus B^+j$ may be viewed as a subalgebra of $C$. In particular, $j \in C^\times$, and since $\sigma$ is an automorphism of $C$, we conclude $jB^\pm \subseteq B^\mp$. This shows that $B^+$ and $B^-$ have the same rank. But the sum of the two ranks is $r$. Thus $\rk(B^\pm) = \frac{r}{2}$.

Conversely, suppose $B \subseteq C$ is a composition subalgebra of rank $\frac{r}{2}$, automatically regular since $\frac{1}{2} \in k$. Hence $C = B \oplus B^\perp$ as direct sum of submodules, and we have
\begin{align}
\label{BEBEP} BB \subseteq C, \quad BB^\perp \subseteq B^\perp, \quad B^\perp B \subseteq B^\perp, \quad B^\perp B^\perp \subseteq B.
\end{align}
Here the first inclusion is obvious, while the second and third one follow from the fact that $C$ is norm associative, particularly from (\ref{p.NORAS}.\ref{RASSNOR}), (\ref{p.NORAS}.\ref{LASSNOR}). Finally, in order to establish the fourth inclusion, we may assume that $k$ is a local ring, in which case Thm.~\ref{t.COMSUBALG}~(a) yields an identification $C = \Cay(B,\mu) = B \oplus Bj$ for some $\mu \in k^\times$ and then $B^\perp = Bj$ by (\ref{p.CDCON}.\ref{CDNORL}). But now the assertion follows from (\ref{p.CDCON}.\ref{CDMULT}). By \eqref{BEBEP}, the map $\sigma_B := \Eins_B \oplus (-\Eins_{B^\perp})$ is a reflection of $C$, and it is straightforward to check that the assignments $\sigma \mapsto \Fix(\sigma)$, $B \mapsto \sigma_B$ define inverse bijections between the set of reflections of $C$  and the set of composition subalgebras of $C$ having rank $\frac{r}{2}$.

Let $\rho \in \Aut(C)$. If $\sigma$ is a reflection of $C$, then so is $\sigma^\rho := \rho^{-1} \circ \sigma \circ \rho$ and $\Fix(\sigma^\rho) = \rho^{-1}(\Fix(\sigma))$. On the other hand, if $B \subseteq C$ is a composition subalgebra of rank $\frac{r}{2}$, then $\sigma_{\rho^{-1}(B)} = (\sigma_B)^\rho$. This shows that two reflections of $C$ are conjugate under $\Aut(C)$ if and only if their fixed algebras are, and the proof of (a) is complete.

(b) Let $\tau$ be an involution of $C$. Since $C$ and $C^{\op}$ have the same unit, norm and trace, Exc.~\ref{pr.CONALGHOMNOR} shows that they are preserved by $\tau$. For $x \in C$, this implies
\[
\tau(\bar x) = \tau\big(t_C(x)1_C - x\big) =t_C\big(\tau(x)\big)1_C - \tau(x) = \overline{\tau(x)}.
\]
Thus $\tau$ commutes with $\iota_C$.

Now let $\tau$ be an involution of $C$ distinct from $\iota_C$. Since $\tau$, as we have just seen, commutes with $\iota_C$, the map $\sigma_\tau := \tau \circ \iota_C = \iota_C \circ \tau$ is a reflection of $C$. Conversely, let $\sigma$ be a reflection of $C$. Being an automorphism of $C$, $\sigma$ preserves conjugation, forcing $\tau_\sigma := \sigma \circ \iota_C = \iota_C \circ \sigma$ to be an involution other than $\iota_C$. It is clear that the assignments $\tau \mapsto \sigma_\tau$ and $\sigma \mapsto \tau_\sigma$ define inverse bijections between the set of involutions of $C$ distinct from $\iota_C$ and the set of reflections of $C$. Moreover, the identities $\tau_{\sigma^\rho} = \rho^{-1} \circ \tau_\sigma \circ \rho$ and $\sigma_{\rho^{-1}\circ\tau\circ\rho} = (\sigma_\tau)^\rho$ are immediately verified. Invoking (a), the aforementioned bijections, therefore, induce canonically inverse bijections between the isomorphism classes of involutions of $C$ other than $\iota_C$ and the conjugacy classes under $\Aut(C)$ of the composition subalgebras of $C$ having rank $\frac{r}{2}$.

Finally, let $\tau \neq \iota_C$ be an involution of $C$ and $B \subseteq C$ the corresponding composition subalgebra of rank $\frac{r}{2}$. Then $\tau = \sigma_B \circ \iota_C$, and for $x \in B$, $y \in B^\perp$, we obtain $\tau(x + y) = \sigma_B(\bar x - y) = \bar x + y$, hence
\[
H(C,\tau) = H(B,\iota_B) \oplus B^\perp = k1_B \oplus B^\perp \cong k \oplus B^\perp,
\]
and since $B^\perp$ is finitely generated projective of rank $\frac{r}{2}$, the final assertion follows.  
\end{sol}


\solnsec{Section ~\ref{s.HERFOR}}

\noexsec


\solnsec{Section~\ref{s.TERHER}}

\begin{sol}{pr.BACSES} \label{sol.BACSES} Let $x_1,\dots,x_n \in M$, $y_1,\dots,y_n \in N$, $r_1,\dots,r_n,s_1,\dots,s_n \in R$. We put $r := r_1\cdots r_n$, $s := s_1\cdots s_n \in R$, $A := \diag(r_1,\dots,r_n)$, $B := \diag(s_1,\dots,s_n) \in \Mat_n(R)$ and obtain
\begin{align*}
 (\bigwedge^nh)_R\big((x_1 \otimes r_1)& \wedge \cdots \wedge (x_n \otimes r_n),(y_1 \otimes s_1) \wedge \cdots \wedge (y_n \otimes s_n)\big) \\
 \,\,&=(\bigwedge^nh)_R\big((x_1 \wedge \cdots \wedge x_n) \otimes r,(y_1 \wedge \cdots \wedge y_n) \otimes s\big) \\
 \,\,&= \big((\bigwedge^nh)(x_1 \wedge \cdots \wedge x_n,y_1 \wedge \cdots \wedge y_n)\big) \otimes rs \\
 \,\,&= \det\Big((h(x_i,y_j)\big)_{1\le i,j\le n}\Big) \otimes rs \\
 \,\,&=rs\det\Big(\big((h(x_i,y_j)\big)_{1\le i,j\le n}\Big)_R \\
 \,\,&= rs\det\Big(\big(h(x_i,y_j)_R\big)_{1\le i,j\le n}\Big) \\
 \,\,&= rs\det\Big(\big(h_R(x_{iR},y_{jR})\big)_{1\le i,j\le n}\Big) \\
 \,\,&= \det(A)\det\Big(\big(h_R(x_{iR},y_{jR})\big)_{1\le i,j\le n}\Big)\det(B) \\
 \,\,&= \det\Big(\big(r_is_jh_R(x_{iR},y_{jR})\big)_{1\le i,j\le n}\Big) \\
 \,\,&= \det\Big(\big(h_R(x_i \otimes r_i,y_j \otimes s_j)\big)_{1\le i,j\le n}\Big) \\
 \,\,&= \bigwedge^n(h_R)\big((x_1 \otimes r_1) \wedge \cdots \wedge (x_n \otimes r_n),(y_1 \otimes s_1) \wedge \cdots \wedge (y_n \otimes s_n)\big),
\end{align*} 
which completes the proof.
\end{sol}

\begin{sol}{pr.INQUA} \label{sol.INQUA}
(a) We begin with the following:
\begin{claim*}
If $k = F$ is a field and $M \neq \{0\}$, then the map $M \to F$, $x \mapsto h(x,x)$, is not identically zero.
\end{claim*}
\begin{proof}
We argue indirectly and assume $h(x,x) = 0$ for all $x \in M$. Linearizing this relation yields
\[
t_D\big(h(x,y)\big) = h(x,y) + \overline{h(x,y)} = h(x,y) + h(y,x) = 0
\]
for all $x,y \in M$, where we may replace $x$ by $xa$, $a \in D$, to obtain $n_D(a,h(x,y)) = t_D(\bar ah(x,y)) = t_D(h(xa,y)) = 0$ for all $a \in D$, $x,y \in M$. But since the norm of $D$ is non-singular, this leads to the contradiction $h = 0$ and proves our claim.
\end{proof}

Returning to our LG ring $k$, 
we argue by induction on $n$ and have nothing to prove for $n = 0$. If $n > 0$, our claim shows that the scalar polynomial law $x \mapsto h(x,x)$ over k represents a unit over any field in $\kalg$. Hence it represents a unit over our LG ring $k$: we deduce $\alpha_1 := h(e_1,e_1) \in k^\times$ for some $e_1 \in M$. Now we put
\[
N := (e_1D)^\perp = \{y \in \mid h(e_1,y) = 0\} = \{x \in M \mid h(x,e_1) = 0\}
\]
and see by a straightforward verification that $M = (e_1D) \oplus N = (e_1D) \perp N$ relative to $h$. It follows that $(N,h\vert_{N\times N })$ is a hermitian space of rank $n - 1$ over $D$, and applying the induction hypothesis completes the proof.

(b) We may clearly assume that $k$ is a local ring. Let the $e_i \in M$, $\alpha_i \in k$, $1 \leq i \leq n$, be as in (a). Then
\[
M = (e_1D) \oplus \cdots \oplus (e_nD)
\]
as a direct sum of free $D$-submodules of rank $1$. Since $D$ is a free $k$-module of rank $2$, this shows that $M$ is a free $k$-module of rank $2n$. On the other hand, the relations $h(e_id,e_jd^\prime) = \delta_{ij}\alpha_i\bar dd^\prime$ for $d,d^\prime \in D$, $1 \leq i,j \leq n$ show that the preceding decomposition is an orthogonal one relative to $h$, hence relative to $q$ as well, and $q(e_id) = h(e_id,e_id) = \alpha_i\bar dd = \alpha_in_D(d)$ for $d \in D$, $1 \leq i \leq n$. Thus
\[
q \cong \alpha_1n_D \perp \cdots \perp \alpha_nn_D \cong \la\alpha_1,\dots,\alpha_n\ra \otimes n_D
\]
is the finite orthogonal sum of of regular quadratic forms over $k$, hence must be regular itself. 
\end{sol}

\begin{sol}{pr.HERISOT} \label{sol.HERISOT}
We put $C := \Ter(D;M,h,\Delta)$ and have $C = D \times M$ as $k$-modules. After the natural identification $D \subseteq C$, we obtain $p^{\pm 1} = (p^{\pm 1},0)$, and writing ``$\cdot$'' for the product in $C^p$, we conclude, for $a,b \in D$, $x,y \in M$,
\begin{align*}
(a,x)\cdot(b,y) =\,\,&\big((a,x)(p^{-1},0)\big)\big((p,0)(b,y)\big) \\
=\,\,&(ap^{-1},xp^{-1})(pb,y\bar p) \\
=\,\,&\big(ap^{-1}pb - h(xp^{-1},y\bar p),y\bar p\overline{ap^{-1}} + xp^{-1}pb + (xp^{-1}) \times_{h,\Delta}(y\bar p)\big) \\
=\,\,&\big(ab - \bar p^{-1}h(x,y)\bar p,y\bar a + xb + (x \times_{h,\Delta} y)\bar p^{-1}p\big)
\end{align*}
Hence
\begin{align}
\label{PROP.SOL} (a,x)\cdot(b,y) = \big(ab - h(x,y),y\bar a + xb + (x \times_{h,\Delta} y)\bar p^{-1}p\big)  &&(a,b \in D,\;x,y \in M).
\end{align}  
We will now be able to compute the constituents of $C^p = \Ter(D;M^p,h^p,\Delta^p)$ by appealing to Thm.~\ref{t.TERHER} as follows.

As a $k$-module, $M^p$ is the orthogonal complement of $D = D^p$ in $C^p$ relative to the bilinearized norm of $C^p$, which by Exc.~\ref{pr.ISTQALT} agrees with the bilinearized norm of $C$. Thus $M^p = M$ as $k$-modules. Consulting \eqref{PROP.SOL}, the right action of $D^p = D$ on $M^p = M$ may be read off from
\[
(x,a) \longmapsto (0,x)\cdot(a,0) = (0,xa),
\]
which amounts to $M^p = M$ as right $D$-modules.

For $x,y \in M$, we deduce from \eqref{PROP.SOL} that $h^p(x,y)$ is the negative of the $D$-component of 
\begin{align}
\label{PROPER} (0,x)\cdot(0,y) = \big(-h(x,y),(x \times_{h,\Delta} y)\bar p^{-1}p\big),
\end{align}
hence agrees with $h(x,y)$. Thus $h^p = h$.

Finally, let $a \in D^\times$. Then $a\Delta\:\bigwedge^3M \overset{\sim}\to D$ is an orientation of $M$, and \ref{ss.HERVE} shows 
\[
h(x \times_{h,a\Delta} y,z) = a\Delta(x \wedge y \wedge z) = ah(x \times_{h,\Delta} y,z) = h\big((x \times_{h,\Delta} y)\bar a,z\big)
\]
for all $x,y,z \in M$, which is equivalent to
\begin{align}
\label{ADEL} x \times_{h,a\Delta} y = (x \times_{h,\Delta} y)\bar a &&(x,y \in M,\;a \in D^\times).
\end{align}
Now by Thm.~\ref{t.TERHER}, $x \times_{h,\Delta^p} y$ is the $M$-component of $(0,x)\cdot(0,y)$, hence by \eqref{PROPER} agrees with $(x \times_{h,\Delta} y)\bar p^{-1}p$. Since, fixing $h$, the hermitian vector product uniquely determines the orientation it corresponds to,
 we conclude from this and \eqref{ADEL} that $\Delta^p = \bar pp^{-1}\Delta$.
 
Summing up, we have proved
\begin{align}
\label{HERP} C^p = \Ter(D;M,h,\bar pp^{-1}\Delta).
\end{align} 
 
Finally, we wish to understand what all this means for the algebra $C = \Zor(k)$ of Zorn vector matrices over $k$. By \ref{ss.ZOVE}, we have
\begin{align*}
&C = \Ter(D;D^3,h,\Delta), \quad D = k \times k, \quad D^3 = k^3 \times k^3, \quad h = \langle\Eins_3\ras, \\ 
&\Delta(\bfe_1 \wedge \bfe_2 \wedge \bfe_3) = 1, \quad \bfe_i = (e_i,e_i) \;\;(1 \leq i \leq 3), 
\end{align*}
where $(e_i)_{1\leq i\leq 3}$ is the canonical basis of unit vectors in $k^3$. The quantity $\bar p^{-1}p$ is essentially an arbitrary element of $D$ having norm $1$, hence can be written in the form $\gamma \oplus \gamma^{-1}$ for some $\gamma \in k^\times$. Writing $x,y \in M$ as $x = u_1 \oplus u_2$, $y = v_1 \oplus v_2$, we may apply (\ref{ss.ZOVE}.\ref{VERE}) and obtain
\begin{align*}
(x \times_{h,\Delta^p} y) =\,\,&(x \times_{h,\Delta} y)\bar p^{-1}p =\big((u_2 \times v_2) \oplus (u_1 \times v_1)\big)(\gamma \oplus \gamma^{-1}) \\
=\,\,&\big(\gamma(u_2 \times v_2)\big) \oplus \big(\gamma^{-1}(u_1 \times v_1)\big).
\end{align*}
Repeating the computations of \ref{ss.ZOVE}, we therefore conclude that $\Zor(k)^p$ is the $k$-module 
\[
\left(\begin{matrix}
k & k^3 \\
k^3 & k
\end{matrix}\right)
\] 
under the multiplication
\[
\left(\begin{matrix}
\alpha_1 & u_2 \\
u_1 & \alpha_2
\end{matrix}\right)\left(\begin{matrix}
\beta_1 & v_2 \\
v_1 & \beta_2
\end{matrix}\right) = \left(\begin{matrix}
\alpha_1\beta_1 - u_2^\trans v_1 & \alpha_1v_2 + \beta_2u_2 + \gamma^{-1}(u_1 \times v_1) \\
\beta_1u_1 + \alpha_2v_1 + \gamma(u_2 \times v_2) & \alpha_2\beta_2 - u_1^\trans v_2
\end{matrix}\right)
\]
for $\alpha_i,\beta_i \in k$, $u_i,v_i \in k^3$, $i = 1,2$. Since $p$ belongs to the split quaternion subalgebra
\[
\left(\begin{matrix}
k & ke_1 \\
ke_1 & k
\end{matrix}\right) \subseteq \Zor(k),
\]
it follows from Exc.~\ref{pr.CDISOT} that $\Zor(k)$ and $\Zor(k)^p$ are isomorphic.
\end{sol}

\begin{sol}{pr.TERFU} \label{sol.TERFU} (i) $\Leftrightarrow$ (ii). Let $\chi\:(M_1,h_1) \to (M_2,h_2)$ be an isometry and $x,y,z \in M_1$. Then
\begin{align*}
h_2\big(\chi(x) \times_{h_2,\Delta_2} \chi(y),\chi(z)\big) =\,\,&\Delta_2\big(\chi(x) \wedge \chi(y) \wedge \chi(z)\big) = (\Delta_2 \circ \bigwedge^3\chi)(x \wedge y \wedge z), \\
h_2\big(\chi(x \times_{h_1,\Delta_1} y),\chi(z)\big) =\,\,&h_1(x \times_{h_1,\Delta_1} y,z) = \Delta_1(x \wedge y \wedge z),
\end{align*}
and comparing yields the asserted equivalence.

(ii) $\Leftrightarrow$ (iii). Apply (\ref{t.TERHER}.\ref{TEMU}).
\end{sol}

\begin{sol}{pr.PRESPLI} \label{sol.PRESPLI}
We proceed in several steps.

We write $(e_i)_{1\leq i\leq 3}$ for the canonical basis of unit vectors in $k^3$ and put
\[
E := \left(\begin{matrix}
1 & 0 \\
0 & 0
\end{matrix}\right), \quad X_1 := \left(\begin{matrix}
0 & -e_1 \\
e_1 & 0
\end{matrix}\right), \quad X_2 := \left(\begin{matrix}
0 & -e_2 \\
e_2 & 0
\end{matrix}\right).
\]
Then (\ref{ss.ZOVE}.\ref{ZOMU}) implies
\begin{align*}
X_3 :=\,\,&X_1X_2 = \left(\begin{matrix}
0 & e_3 \\
e_3 & 0
\end{matrix}\right), \quad
X_1^2 = \left(\begin{matrix}
1 & 0, \\
0 & 1
\end{matrix}\right)= 1_C = X_2^2, \\
X_1X_2X_1 =\,\,&X_3X_1 = \left(\begin{matrix}
0 & e_3 \\
e_3 & 0
\end{matrix}\right)\left(\begin{matrix}
0 & -e_1 \\
e_1 & 0 
\end{matrix}\right) = \left(\begin{matrix}
0 & e_2 \\
-e_2 & 0
\end{matrix}\right) = -X_2.
\end{align*}
Moreover, since $E$ is an elementary idempotent, we have
\[
\bar E = 1_C - E = \left(\begin{matrix}
0 & 0 \\
0 & 1
\end{matrix}\right),
\]
and using (\ref{ss.IDCO}.\ref{QAUOP}), we see $X_iEX_i = n_C(X_i,\bar E)X_i - n_C(X_i)\bar E$ which is $\bar E$ for $i = 1,2$ and $-\bar E$ for $i = 3$. We have thus established all the relations of \eqref{PRESPLI}. Now observe that an easy computation yields
\begin{align*}
EX_i = \left(\begin{matrix}
0 & -e_i \\
0 & 0
\end{matrix}\right), \quad EX_3 = \left(\begin{matrix}
0 & e_3 \\
0 & 0
\end{matrix}\right) &&(i = 1,2),
\end{align*}
which shows that the quantities
\[
E,\bar E, X_1,X_2,X_3,EX_1,EX_2,EX_3
\]
form a basis of $C$ as a $k$-module. In particular, $C$ is generated by $E,X_1,X_2$ as a $k$-algebra. 

\step{2}
Now let $A$ be any unital $k$-algebra and suppose that $e,x_1,x_2 \in A$ satisfy the relations 
\begin{align}
\label{SMAPR} &e^2 = e, \quad x_1^2 = x_2^2 = 1_A, \quad x_1x_2x_1 = -x_2, \\
&x_1ex_1 = x_2ex_2 = \bar e := 1_A - e, \quad (x_1x_2)e(x_1x_2) = -\bar e. \notag
\end{align}
Then we put
\begin{align}
\label{SIGN} \vep_{12} := 1, \quad \vep_{21} := -1, \quad x_3 := x_1x_1.
\end{align}
Since \eqref{SMAPR}, \eqref{SIGN} yield $x_2x_1x_2 = x_2x_1x_2x_1^2 = -x_2x_2x_1 = -x$ and $x_2x_1 = x_1^2x_2x_1 = -x_1x_2 = -x_3$, we obtain
\begin{align}
\label{UNIF} &\vep_{ji} = -\vep_{ij}, \quad e^2 = e, \quad x_i^2 = 1_A, \quad x_ix_jx_i = -x_j, \\
&x_iex_i = \bar e, \quad(x_ix_j)e(x_ix_j) = -\bar e, \quad x_ix_j = \vep_{ij}x_3 &&(\{i,j\} = \{1,2\}). \notag  
\end{align}
Next we put
\begin{align}
\label{BEX} B := ke + k\bar e + \sum_{i=1}^3 kx_i + \sum_{i=1}^3 ky_i \subseteq A &&(y_i := ex_i,\;1 \leq i \leq 3)
\end{align}
and claim that \emph{the spanning elements of $B$ as a $k$-module displayed in \eqref{BEX} satisfy the following multiplication rules.}
\begin{align}
\label{SQR} e^2=\,\,&e, \quad \bar e^2 = \bar e, \quad x_i^2 = 1_C, \quad x_3^2 = -1_C, \quad y_i^2 = y_3^2 = 0 &&(i = 1,2), \\
\label{XYE} x_iy_i =\,\,&\bar e, \quad y_ix_i = e, \quad x_3y_3 = -\bar e, \quad y_3x_3 = -e &&(i = 1,2), \\
\label{ELE} e\bar e =\,\,&0, \quad ex_i = y_i, \quad ey_i = y_i &&(1 \leq i \leq 3), \\
\label{ERI} \bar ee =\,\,& 0, \quad x_ie = x_i - y_i, \quad y_ie = 0 &&(1 \leq i \leq 3), \\
\label{BELE} \bar ex_i =\,\,&x_i - y_i, \quad \bar ey_i = 0 &&(1 \leq i \leq 3), \\
\label{BERI} x_i\bar e =\,\,&y_i, \quad y_i\bar e = y_i &&(1 \leq i \leq 3), \\
\label{XLE} x_ix_j =\,\,&\vep_{ij}x_3, \quad x_ix_3 = \vep_{ij}x_j, \quad x_iy_j = \vep_{ij}(x_3 - y_3), \quad x_iy_3 = \vep_{ij}(x_j - y_j) &&(\{i,j\} = \{1,2\}), \\
\label{XRI} x_3x_i =\,\,&-\vep_{ij}x_j, \quad y_jx_i = -\vep_{ij}(x_3 - y_3), \quad  y_3x_i = -\vep_{ij}(x_j - y_j) &&(\{i,j\} = \{1,2\}), \\
\label{XTHR} x_3y_i =\,\,&-\vep_{ij}(x_j - y_j), \quad y_ix_3 = \vep_{ij}(x_j - y_j)  &&(\{i,j\} = \{1,2\}), \\
\label{YLE} y_iy_j =\,\,&\vep_{ij}(x_3 - y_3), \quad y_iy_3 = \vep_{ij}(x_j - y_j), \quad y_3y_i  = -\vep_{ij}(x_j - y_j) &&(\{i,j\} = \{1,2\}).
\end{align}
Suppose these relations have been proved. Then they hold, mutatis mutandis, for $E,X_1,X_2 \in C$ as well, and it follows that the linear map $C \to A$ acting on the basis vectors of $C$ exhibited in $1^0$ according to
\begin{align*}
E \mapsto e, \quad \bar E \mapsto \bar e, \quad X_i \mapsto x_i, \quad EX_i \mapsto y_i &&(1 \leq i \leq 3) 
\end{align*} 
is the unique unital homomorphism satisfying the conditions of the problem. The kernel of this homomorphism is an ideal in $C$, which, by Exc.~\ref{pr.IDEALSCOAL}, has the form $\mfa C$ for some ideal $\mfa \subseteq k$. Hence $B \cong C/\mfa C \cong C \otimes (k/\mfa) \cong \Zor(k/\mfa)$ as $k$-algebras. 

\step{3}
It remains to establish \eqref{SQR}--\eqref{YLE}, which we now proceed to do by making use, among other things, of Artin's theorem (Cor.~\ref{c.ARTIN}) allowing us to drop parentheses in products with at most two distinct factors.

\emph{Proof of} \eqref{SQR}. The first four relations are clear in view of \eqref{UNIF}. Moreover, by \eqref{SMAPR}, \eqref{SIGN}, $x_3^2 = x_1x_2x_1x_2 = -x_2^2 = -1_C$, $y_i^2 = ex_iex_i = e\bar e = 0$ for $i = 1,2$ and $y_3^2 = ex_3ex_3 = -e\bar e = 0$, which completes the proof.

\emph{Proof of} \eqref{XYE}. For $i = 1,2$, the definitions and \eqref{UNIF} yield $x_iy_i = x_iex_i = \bar e$, $y_ix_i = ex_i^2 = e$, $x_3y_3 = x_3ey_3 = -\bar e$, $y_3x_3 = ex_3^2 = -e$, as desired. 

\emph{Proof of} \eqref{ELE}. The first relation is obvious, the next three hold by definition (cf. \eqref{BEX}), while the remaining ones follow from $ey_i = eex_i = ex_i = y_i$ for $i = 1,2,3$.

\emph{Proof of} \eqref{ERI}. The first equation is again obvious. As to the next three, we apply \eqref{UNIF}, \eqref{SQR} to obtain $x_ie = x_iex_i^2 = \bar ex_i = x_i - ex_i = x_i - y_i$ for $i = 1,2$ and $x_3e = -x_3ex_3^2 = \bar ex_3 = x_3 - ex_3 = x_3 - y_3$, as claimed.

\emph{Proof of} \eqref{BELE} We have $\bar ex_i = x_i - ex_i$ for $1 \leq i \leq 3$, hence the first three equations of \eqref{BELE}, while \eqref{ERI} yields $y_i\bar e = y_i - y_ie = y_i$ for $1 \leq i \leq 3$, which completes the proof.

\emph{Proof of} \eqref{BERI}. From \eqref{ERI} we deduce $x_i\bar e = x_i - x_ie = y_i$, $y_i\bar e = y_i - y_ie = y_i$ for $1 \leq i \leq 3$, and the assertion follows.

\emph{Proof of} \eqref{XLE}. Let $\{i,j\} = \{1,2\}$. The first equation is already in \eqref{UNIF}. By the same token, $x_ix_3 = \vep_{ij}x_ix_ix_j = \vep_{ij}x_j$ yields the second equation. Applying \eqref{UNIF}, \eqref{BELE} and the left Moufang identity, we obtain $x_iy_j = x_i(ex_j) = \vep_{ij}x_i(e(x_ix_3)) = \vep_{ij}(x_iex_i)x_3 = \vep_{ij}\bar ex_3 = \vep_{ij}(x_3 - y_3)$ and $x_iy_3 = x_i(ex_3) = \vep_{ij}x_i(e(x_ix_j)) = \vep_{ij}(x_iex_i)x_j = \vep_{ij}\bar ex_j = \vep_{ij}(x_j - y_j)$, which completes the proof of \eqref{XLE}.

\emph{Proof of} \eqref{XRI}. Let ${i,j} = \{1,2\}$. Then \eqref{XLE} yields $x_3x_i = \vep_{ij}x_ix_jx_i = -\vep_{ij}x_j$, while this, \eqref{ERI}, \eqref{UNIF}, \eqref{XLE}, \eqref{BERI} imply $y_jx_i = x_jx_i - (x_je)x_i = \vep_{ji}x_3 + \vep_{ij}((x_3x_i)e)x_i = \vep_{ij}(x_3(x_iex_i) -x_3) = \vep_{ij}(x_3\bar e - x_3) = -\vep_{ij}(x_3 - y_3)$. Finally, applying \eqref{ERI}, \eqref{XRI}, \eqref{XLE}, \eqref{BERI}, we obtain $y_3x_i = x_3x_i - (x_3e)x_i = -\vep_{ij}x_j -\vep_{ji}((x_jx_i)e)x_i = \vep_{ij}(x_j(x_iex_i) - x_j) = \vep_{ij}(x_j\bar e - x_j) = -\vep_{ij}(x_j - y_j)$. 

\emph{Proof of} \eqref{XTHR}. Let $\{i,j\} = \{1,2\}$. Then \eqref{SQR}, \eqref{XLE}, \eqref{ERI}, \eqref{XRI} yield
\begin{align*}
x_3y_i =\,\,&-x_3(ex_i)x_3^2 = -[x_3(ex_i)x_3]x_3 = -[(x_3e)(x_ix_3)]x_3 = -\vep_{ij}[(x_3e)x_j]x_3 \\
=\,\,&-\vep_{ij}x_3x_jx_3 + \vep_{ij}[(ex_3)x_j]x_3 = \vep_{ij}\vep_{ji}x_ix_3 + \vep_{ij}e(x_3x_jx_3) \\ 
=\,\,&= -\vep_{ij}x_j - \vep_{ij}\vep_{ji}e(x_ix_3) -\vep_{ij}x_j + \vep_{ij}ex_j = -\vep_{ij}(x_j - y_j), \\
y_ix_3 =\,\,&[x_i(x_ie)x_i]x_3 = x_i[(x_ie)(x_ix_3)] = \vep_{ij}x_i[(x_ie)x_j] = \vep_{ij}x_i[(x_i - y_i)x_j] \\
=\,\,&\vep_{ij}\big(x_i^2x_j - x_i(y_ix_j)\big) = \vep_{ij}\big(x_j + \vep_{ji}x_i(x_3 - y_3)\big) = \vep_{ij}x_j - x_ix_3 + x_iy_3 \\
=\,\,&\vep_{ij}x_j - \vep_{ij}x_j + \vep_{ij}(x_j - y_j) = \vep_{ij}(x_j - y_j),
\end{align*}
hence \eqref{XTHR}.

\emph{Proof of} \eqref{YLE}. Again, let $\{i,j\} = \{1,2\}$. Then \eqref{ERI}, \eqref{XRI}, \eqref{ELE} yield $y_iy_j =  y_ix_j - (ex_i)(x_je) = \vep_{ij}(x_3 - y_3) -e(x_ix_j)e = \vep_{ij}(x_3 - y_3 -ex_3e) = \vep_{ij}(x_3 - y_3)$, giving the first equation of \eqref{YLE}. Similarly, $y_iy_3 = y_ix_3 - (ex_i)(x_3e) = \vep_{ij}(x_j - y_j) - e(x_ix_3)e = \vep_{ij}(x_j - y_j - ex_je) = \vep_{ij}(x_j - y_j)$. And finally, $y_3y_i = y_3(x_i - x_ie) = y_3x_i - (ex_3)(x_ie) = -\vep_{ij}(x_j - y_j) -e(x_3x_i)e = -\vep_{ij}(x_j - y_j) + \vep_{ij}ex_je = -\vep_{ij}(x_j - y_j)$, which completes the entire proof.
\end{sol}

\solnsec{Section~\ref{s.REDCO}}

\begin{sol}{pr.DIPRAB} \label{sol.DIPRAB} (i) $\Rightarrow$ (ii). $c \ne 0$ implies $k \ne \{0\}$, and if $k^\pm \ne \{0\}$, then $c^\pm = c_{k^\pm} \in A_{k^\pm} = A^\pm$ is an absolutely primitive idempotent by definition.

(ii) $\Rightarrow$ (i). For each sign $\pm$, $c^\pm \in A^\pm$ is an idempotent, and for some sign $\pm$, we have $c^\pm \ne 0$. Hence $c \ne 0$ is an idempotent in $A$. Now assume we are given orthogonal 
idempotents $d = (d^+,d^-), e = (e^+,e^-) \in A$ such that $c = d + e$. Then $d^\pm$, $e^\pm$ are orthogonal idempotents in $A^\pm$, $c^\pm = d^\pm + e^\pm$, and since $k^\pm = \{0\}$ or $c^\pm$ is primitive, we find complete orthogonal systems $(\vep_1^\pm,\vep_2^\pm)$ of idempotents in $k^\pm$ such that $d^\pm =\vep_1^\pm c^\pm$, $e^\pm = \vep_2^\pm c^\pm$. Hence $(\vep_1,\vep_2)$, $\vep_1 = (\vep_1^+,\vep_1^-)$, $\vep_2 = (\vep_2^+,\vep_2^-)$, is a complete orthogonal system of idempotents in $k$ satisfying $d = \vep_1c$, $e = \vep_2c$. Thus $c$ is a primitive idempotent in $A$.

It remains to show that the same conclusion holds under every non-zero base change of $A$, equivalently, that condition (ii) is stable under such a non-zero base change, so let $R \in \kalg$, $R \neq \{0\}$. Setting $\vep^+ := (1,0)$, $\vep^- := (0,1) \in k$ and applying \ref{ss.PROPID}, we conclude $R^+ \times R^-$, $R^\pm = \vep^\pm R = R_{k^\pm}$, $A_R = A_R^+ \times A_R^-$, $A_R^\pm = (A_R)_{R^\pm} = A_{R^\pm} = (A_{k^\pm})_{R^\pm} = A^\pm_{R^\pm}$. We have $c_R = (c_{R^+}^+, c_{R^-}^-)$, and $R^\pm \ne \{0\}$ for some sign $\pm$ implies $k^\pm \ne \{0\}$, hence by (ii) that $c^\pm$ is an absolutely primitive idempotent of $A^\pm$. But then $c_{R^\pm}^\pm$ must be an absolutely primitive idempotent of $A_{R^\pm}^\pm$, and we have shownn that (ii) is indeed stable under non-zero base change.

\end{sol}

\begin{sol}{pr.DECAPID} \label{sol.DECAPID} (i) $\Rightarrow$ (ii). We clearly have $k \ne \{0\}$. Assuming $k^{(0)} \ne \{0\}$ in the decomposition (ii) of Exc.~\ref{pr.DECID}, the base change $0 = c_{k^{(0)}} \in C^{(0)} = C_{k^{(0)}}$ would be primitive, a contradiction. Hence $k^{(0)} = \{0\}$. Similarly, if $k^{(2)} \ne \{0\}$, then $1_{C^{(2)}} = c_{k^{(2)}}$ is absolutely primitive in $C^{(2)} = C_{k^{(2)}}$. Hence every composition algebra appearing in the rank decomposition of $C_{k^{(2)}}$ (Exc.~\ref{pr.RANDEC}) contains its identity element as an absolutely primitive idempotent and thus, by Prop.~\ref{p.ABPRIELID}, has rank $1$. We conclude $C^{(2)} \cong k^{(2)}$ and have established decompositions \eqref{DECAPBAS}, \eqref{DECAPAL}.

(ii) $\Rightarrow$ (i). If $k^{(1)} \ne \{0\}$, then $c^{(1)}$ is absolutely primitive in $C^{(1)}$ by  Prop.~\ref{p.ABPRIELID}, and if $k^{(2)} \ne \{0\}$, then $1_{k^{(2)}}$ is an absolutely primitive idempotent of $k^{(2)}$. Thus (i) follows from Exc.~\ref{pr.DIPRAB}.

The final statement of the exercise is now obvious.
\end{sol}

\begin{sol}{pr.ELIDREDQUAT} \label{sol.ELIDREDQUAT}
(a)\;(i) $\Rightarrow$ (ii). We view $c$ canonically as an idempotent linear map $c\:k \oplus L_0 \to k \oplus L_0$. With $L_c := \Img(c)$, $\bar L_c := \Ker(c) = L_{\bar c}$, we then have
\begin{align}
\label{IDEC} k \oplus L_0 = L_c \oplus \bar L_c = L_c \oplus L_{\bar c}
\end{align}
as a directi sum of submodules. Hence $L_c$ and $L_{\bar c}$ are both finitely generated projective $k$-modules. For $\mfp \in \Spec(k)$, we have $c_\mfp \neq 0 \neq \bar c_\mfp$ since $c$ is elementary, so $(L_c)_\mfp = L_{c_\mfp}$ and $(L_{\bar c})_\mfp = L_{\bar c_\mfp}$ are both free $k_\mfp$-modules of positive rank. On the other hand, \eqref{IDEC} yields $\rk_\mfp(L_c) + \rk_\mfp(L_{\bar c}) = 2$, which altogether implies $\rk_\mfp(L_c) = \rk_\mfp(L_{\bar c}) = 1$. Hence $L_c$ is a line bundle, as claimed.

(ii) $\Rightarrow$ (iii). Regardless of whether the idempotent $c$ is elementary or not, \eqref{IDEC} holds, forcing $L_c$ and $L_{\bar c}$ to be finitely generated projective $k$-modules. Now suppose as in (ii) that $L_c$ is a line bundle. Then \cite[III.7, Cor. of Prop.~10]{MR0354207} implies
\[
L_0 \cong k \otimes L_0 \cong \biw^2(k \oplus L_0) \cong \biw^2(L_c \oplus L_{\bar c}) \cong L_c \otimes L_{\bar c}.
\]
Thus $L_{\bar c} \cong L_0 \otimes L_c^\ast$ is a line bundle as well, and we have established \eqref{ELCE}. At this stage, we require the following well-known fact:
\begin{fact*} 
Let $M,N$ be finitely generated projective $k$-modules. Then there is a natural identification $\Hom_k(M,N) = N \otimes M^\ast$ such that $(y \otimes x^\ast)(x^\prime) = \la x^\ast,x^\prime\ra y$ for all $x^\prime \in M$, $y \in N$, $x^\ast \in M^\ast$.
\end{fact*}
Using this and \eqref{ELCE}, we now obtain a chain of natural isomorphism as follows:
\begin{align*}
B \cong\,\,&\End_k(k \oplus L_0) \cong \End_k(L_c \oplus L_{\bar c}) \cong \End_k\left(\begin{matrix}
L_c \\
L_{\bar c}
\end{matrix}\right) \\
\cong\,\,&\left(\begin{matrix}
\Hom_k(L_c,L_c) & \Hom_k(L_{\bar c}, L_c) \\
\Hom_k(L_c,L_{\bar c}) & \Hom_k(L_{\bar c},L_{\bar c})
\end{matrix}\right) \cong \left(\begin{matrix}
L_c \otimes L_c^\ast & L_c \otimes L_{\bar c}^\ast \\
L_{\bar c} \otimes L_c^\ast & L_{\bar c} \otimes L_{\bar c}^\ast
\end{matrix}\right) \\
\cong\,\,&\left(\begin{matrix}
k & L_0^\ast \otimes L_c^{\otimes 2} \\
L_0 \otimes L_c^{\ast\otimes 2} & k
\end{matrix}\right) \cong B^\prime := \End_k(k \oplus L)
\end{align*}
with $L := L_0 \otimes L_c^{\ast\otimes 2}$ an appropriate line bundle over $k$. Writing $\Phi\:B \to B^\prime$ for the composite of these isomorphisms, and following step-by-step their effect on the idempotent $c$, we conclude
\begin{align}
\label{PHICE} \Phi(c) = \left(\begin{matrix}
1 & 0 \\
0 & 0
\end{matrix}\right),
\end{align}
which completes the proof of (iii). Now suppose $L$ is any line bundle over $k$ such that there exists an isomorphism $\Phi\:B \to B^\prime := \End_k(k \oplus L)$ satisfying \eqref{PHICE}. Then $L \cong B^\prime_{21}(\Phi(c)) \cong B_{21}(c)$, $L^\ast \cong B^\prime_{12}(\Phi(c)) \cong B_{12}(c)$, which not only shows that $L$ is unique up to isomorphism but also completes the proof of \eqref{BETCE}.

(iii) $\Rightarrow$ (i). Since $\Phi(c)$ is obviously elementary, so is $c$.

(b) Let $g \in \GL(k \oplus L_0) = B^\times$ such that $d = gcg^{-1}$. Then
\[
L_d =d\left(\begin{matrix}
k \\
L_0
\end{matrix}\right) = gcg^{-1}\left(\begin{matrix}
k \\
L_0
\end{matrix}\right) = gc\left(\begin{matrix}
k \\
L_0
\end{matrix}\right) = g(L_c).
\]
Hence $g$ determines via restriction an isomorphism $L_c \overset{\sim} \to L_d$. Conversely, suppose $L_c$ and $L_d$ are isomorphic. Then, by \eqref{ELCE}, so are $L_{\bar c}$ and $L_{\bar d}$. Let
\[
\rho\:L_c \overset{\sim} \longrightarrow L_d, \quad \bar\rho\:L_{\bar c} \overset{\sim} \longrightarrow L_{\bar d}
\]
be isomorphisms. Then so is
\[
g := \rho \oplus \bar\rho\:k \oplus L_0 = L_c \oplus L_{\bar c} \overset{\sim} \longrightarrow L_d \oplus L_{\bar d} = k \oplus L_0,
\]
which amounts to $g \in \GL(k \oplus L_0) = B^\times$. Moreover, for $x \in L_d$, $y \in L_{\bar d}$ we have $\rho^{-1}(x) \in L_c$, $\bar\rho^{-1}(y) \in L_{\bar c}$, hence
\[
gcg^{-1}(x + y) =gc\big(\rho^{-1}(x) + \bar\rho^{-1}(y)\big) = (\rho \oplus \bar\rho)\big(\rho^{-1}(x)\big) = x.
\]  
But this means $gcg^{-1} = d$, and $c,d$ are conjugate under inner automorphisms of $B$. Now suppose $\vph(c) = d$ for some automorphism $\vph$ of $B$. Then $B_{12}(c) \cong B_{12}(d)$, which by \eqref{BETCE} implies $L_c^{\otimes 2} \cong L_d^{\otimes 2}$. Conversely, let this be so and put $L := L_0 \otimes L_c^{\ast\otimes 2}$. Then (iii) yields isomorphisms $\Phi\:B \overset{\sim} \to B^\prime := \End_k(k \oplus L)$, $\Psi\:B \overset{\sim} \to B^\prime$ such that
\[
\Phi(c) = \left(\begin{matrix}
1 & 0 \\
0 & 0
\end{matrix}\right) = \Psi(d).
\]
Hence $\vph := \Psi^{-1} \circ \Phi$ is an automorphism of $B$ sending $c$ to $d$.

(c) By \eqref{ELCE}, the line bundle $L_c$ is a direct summand of $k \oplus L_0$. Conversely, let $L$ be a direct summand of $k \oplus L_0$, so $k \oplus L_0 = L \oplus L^\prime$ for some submodule $L^\prime \subseteq k \oplus L_0$. Then the projection from $k \oplus L_0$ onto $L$ alongside $L^\prime$ yields an idempotent $c \in B$ such that $L = L_c$, and by (a), $c$ is elementary.
\end{sol}

\begin{sol}{pr.LITWO} \label{sol.LITWO}
(a) Let $x_i/1$ be a basis of $L_{f_i}$ over $k_{f_i}$, for
$i = 1,\dots,n$. Given $x \in L$, there exist $m \in \IN$ and
$\alpha_1,\dots,\alpha_n \in k$ such that $x/1 =
(\alpha_i/f_i^m)(x_i/1) = (\alpha_ix_i)/f_i^m$ in $L_{f_i}$
for all $i = 1,\dots,n$. Hence for some integer $p \geq m$ and all
$i = 1,\dots,n$, $f_i^px = \beta_ix_i$, $\beta_i :=
f_i^{p-m}\alpha_i$. Since the $f_1^p,\dots,f_n^p$ continue to
generate $k$ as an ideal, we find $g_1,\dots,g_n \in k$ such that
$\sum f_i^pg_i = 1$. But this implies $x = \sum \beta_ig_ix_i$, so
$L = \sum kx_i$ is generated by $x_1,\dots,x_n$. 

(b) (i) $\Rightarrow$ (ii). Since $L$ is generated by two elements, we obtain a short exact sequence
\[
\xymatrix{0 \ar[r] & L^\prime \ar[r] & k \oplus k \ar[r] & L \ar[r] & 0}
\]
of $k$-modules, which splits since $L$ is projective. Thus $L \oplus L^\prime \cong k \oplus k$ is free of rank $2$, and taking determinants ($=$ second exterior powers), we obtain $L \otimes L^\prime \cong k$, i.e., $L^\prime \cong L^\ast$.

(ii) $\Rightarrow$ (iii). This is Exc.~\ref{pr.ELIDREDQUAT} for $L_0  := k$.

(iii) $\Rightarrow$ (iv). Write 
\[
c = \left(\begin{matrix}
\alpha & \beta \\
\gamma & \delta
\end{matrix}\right)
\]
with $\alpha,\beta,\gamma,\delta \in k$. Since $c$ has trace $1$ and determinant $0$, we have
\begin{align}
\label{EQEL} \alpha + \delta = 1, \quad \alpha\delta = \beta\gamma.
\end{align}
From Exc.~\ref{pr.ELIDREDQUAT} we deduce
\[
L \cong L_c = \Img(c) = k\left(\begin{matrix}
\alpha \\
\gamma 
\end{matrix} \right) + k\left(\begin{matrix}
\beta \\
\delta
\end{matrix} \right).
\]
Observing
\begin{align}
\label{COLELIDMAT} \beta\left(
\begin{matrix}
\alpha \\
\gamma
\end{matrix}
\right) = \alpha\left(
\begin{matrix}
\beta             \\
\delta
\end{matrix}
\right), \quad \delta\left(
\begin{matrix}
\alpha \\
\gamma
\end{matrix}
\right) = \gamma\left(
\begin{matrix}
\beta             \\
\delta
\end{matrix}
\right),
\end{align}
we now put $f_1 = \alpha$, $f_2 = \delta$. Then $kf_1 + kf_2 = k$ by \eqref{EQEL}, and \eqref{COLELIDMAT} shows that $L_{f_1}$ is the free $k_{f_1}$-module of rank $1$ with basis $\left(\begin{smallmatrix}
\alpha/1 \\
\beta/1
\end{smallmatrix}\right)$, while $L_{f_2}$ is the free $k_{f_2}$-module of rank $1$ with basis $\left(\begin{smallmatrix}
\beta/1 \\
\delta/1
\end{smallmatrix}\right)$. Thus (iv) holds.

(iv) $\Rightarrow$ (i). This is just a special cse of (a).

Now suppose $L$ satisfies one (hence all) of the preceding four conditions. Then, by (ii), so does $L^\ast$, allowing us to assume $n > 0$. But then (iv) holds for $L^{\otimes n}$. In partiuclar, by (ii), an elementary idempotent $c^{(n)} \in \Mat_2(k)$ satisfying $L^{\otimes n} \cong L_{c^{(n)}}$ exists for any $n \in \IZ$ and is unique up to conjugation by inner automorphisms (Exc.~\ref{pr.ELIDREDQUAT}~(b)). Since $L^{\otimes 0} \cong k$ is free of rank $1$, we may put $c^{(0)} = \left(\begin{smallmatrix}
1 & 0 \\
0 & 0
\end{smallmatrix}\right)$. Moreover, combining (ii) with \eqref{ELCE} in Exc.~\ref{pr.ELIDREDQUAT}, we may also put $c^{(-1)} := \bar c = \Eins_2 - c$, and it will be enough to treat the case $n > 0$. Writing $c$ as in the proof of the implication (iii)$\Rightarrow$(iv), we claim that
\begin{align}
\label{TREN} \alpha^n\alpha_n + \delta^n\delta_n = 1,
\end{align}
\emph{where
\[
\alpha_n := \sum_{i=0}^{n-1}
\binom{2n-1}{i}\alpha^{n-1-i}\delta^i, \quad \delta_n
:= \sum_{i=0}^{n-1}
\binom{2n-1}{n+i}\alpha^{n-1-i}\delta^i,
\]
and
\begin{align}
\label{CEEN} c^{\langle n\rangle} := \left(
\begin{matrix}
\alpha^n\alpha_n & \beta^n\delta_n             \\
\gamma^n\alpha_n & \delta^n\delta_n
\end{matrix}
\right) \in \Mat_2(k)
\end{align}
is an elementary idempotent satisfying $L_{c^{\langle n \rangle}}
\cong L_c^{\otimes n}$.}

We begin by proving \eqref{TREN}, which follows from the computation
\begin{align*}
\alpha^n\alpha_n + \delta^n\delta_n =\,\,&\sum_{i=0}^{n-1} \binom{2n-1}{i}\alpha^{2n-1-i}\delta^i + \sum_{i=0}^{n-1} \binom{2n-1}{n+i}\alpha^{n-1-i}\delta^{n+i} \\
=\,\,&\sum_{i=0}^{n-1} \binom{2n-1}{i}\alpha^{2n-1-i}\delta^i + \sum_{j=n}^{2n-1} \binom{2n-1}{j}\alpha^{2n-1-j}\delta^j \\
=\,\,&\sum_{i=0}^{2n-1} \binom{2n-1}{i}\alpha^{2n-1-i}\delta^i = (\alpha + \delta)^{2n-1} = 1.
\end{align*}
Combining \eqref{TREN} with \eqref{COLELIDMAT}, we conclude that $c^{(n)}$ as defined in \eqref{CEEN} has trace $1$ and determinant $0$, hence is an elementary idempotent. Moreover, 
\begin{align}
\label{COLEN} \beta^n\left(\begin{matrix}
\alpha^n \\
\gamma^n
\end{matrix}\right) = \alpha^n\left(\begin{matrix}
\beta^n \\
\delta^n
\end{matrix}\right), \quad \delta^n\left(\begin{matrix}
\alpha^n \\
\gamma^n
\end{matrix} \right) = \gamma^n\left(\begin{matrix}
\beta^n \\
\delta^n
\end{matrix}\right).
\end{align}
We now put
\begin{align}
\label{EFEN} f^{(n)} := \alpha^n\alpha_n, \quad g^{(n)} := \delta^n\delta_n
\end{align}
and have $f^{(n)} + g^{(n)} = 1$ by \eqref{TREN}. Setting $L := L_c$, $L_{c^{(n)}}$, we must show $L^{\otimes n} \cong L^{(n)}$. To this end, we examine the situation of $k_{f^{(n)}}$ and $k_{g^{(n)}}$.

Over $k_{f^{(n)}}$, the quantities $\alpha,\alpha_n$ are both invertible. By \eqref{COLELIDMAT}, therefore, $L_{f^{(n)}}$ is the free $k_{f^{(n)}}$- module with basis $\left(\begin{smallmatrix}
\alpha \\
\gamma
\end{smallmatrix}\right)$. Consequently, $(L^{\otimes n})_{f^{(n)}} = (L_{f^{(n)}})^{\otimes n}$ is the free $k_{f^{(n)}}$-module with basis $\left(\begin{smallmatrix}
\alpha \\
\gamma
\end{smallmatrix}\right)^{\otimes n}$. By the same token, $(L^{(n)})_{f^{(n)}}$ is the free $k_{f^{(n)}}$-module with basis $\left(\begin{smallmatrix}
\alpha^n\alpha_n \\
\gamma^n\alpha_n 
\end{smallmatrix}\right) = \left(\begin{smallmatrix}
\alpha^n \\
\gamma^n  
\end{smallmatrix}\right)\alpha_n$, which is the same as the free $k_{f^{(n)}}$-module with basis $\left(\begin{smallmatrix}
\alpha^n \\
\gamma^n 
\end{smallmatrix}\right)$. Summing up, we find a unique isomorphism
\[
\Phi^{(n)}\:(L^{\otimes n})_{f^{(n)}} \overset{\sim}\longrightarrow (L^{(n)})_{f^{(n)}}\;\;\text{such that}\;\;\Phi^{(n)}\Big(\left(\begin{matrix}
\alpha \\
\gamma 
\end{matrix}\right)^{\otimes n}\Big) = \left(\begin{matrix}
\alpha^n \\
\gamma^n 
\end{matrix}\right).
\]
Similarly, by \eqref{EFEN}, the quantities $\delta,\delta^{(n)}$ are invertible over $k_{g^{(n)}}$, and there is a unique isomorphism
\[
\Psi^{(n)}\:(L^{\otimes n})_{g^{(n)}} \overset{\sim}\longrightarrow (L^{(n)})_{g^{(n)}}\;\;\text{such that} \;\;\Psi^{(n)}\Big(\left(\begin{matrix}
\beta \\
\delta 
\end{matrix}\right)^{\otimes n}\Big) = \left(\begin{matrix}
\beta^n \\
\delta^n 
\end{matrix}\right).
\]
Over $k_{f^{(n)}g^{(n)}}$, all the quantities $\alpha,\alpha_n,\delta,\delta_n,\beta,\gamma$ become invertible, and we have $\left(\begin{smallmatrix}
\beta \\
\delta 
\end{smallmatrix} \right) = \gamma^{-1}\delta\left(\begin{smallmatrix}
\alpha \\
\gamma 
\end{smallmatrix}\right)$, hence $\left(\begin{smallmatrix}
\beta \\
\gamma 
\end{smallmatrix}\right)^{\otimes n} = \gamma^{-n}\delta^n\left(\begin{smallmatrix}
\alpha \\
\gamma 
\end{smallmatrix}\right)^{\otimes n}$, but also $\left(\begin{smallmatrix}
\beta^n \\
\delta^n 
\end{smallmatrix}\right)= \gamma^{-n}\delta^n\left(\begin{smallmatrix}
\alpha^n \\
\gamma^n 
\end{smallmatrix}\right)$. Hence the isomorphisms $\Phi^{(n)}$, $\Psi^{(n)}$ agree over $k_{f^{(n)}g^{(n)}}$ and thus glue to an isomorphism $L^{\otimes n} \overset{\sim} \to L^{(n)}$. This completes the proof.

\begin{rmk*}
The final proof of this exercise would have become much more natural if we had used the identification of $\Pic(k)$ with $H_{\Zar}^1(k,\GL_1)$, compare Example \ref{e.Pic.H1}. 
\end{rmk*}
\end{sol}

\begin{sol}{pr.ISOREDQUAT} \label{sol.ISOREDQUAT}
(i) $\Rightarrow$ (ii). Let $L$ be any line bundle over $k$ that is a direct summand of $k \oplus L_0$. By Exc.~\ref{pr.ELIDREDQUAT}~(c), there exists an elementary idempotent $c \in B$ such that $L \cong L_c$. Using (i) to select an isomorphism $\Phi\:B \overset{\sim} \to B^\prime$, we obtain an  elementary idempotent $c^\prime := \Phi(c) \in B^\prime$ and, again by Exc.~\ref{pr.ELIDREDQUAT}, we find in $L^\prime := L_{c^\prime}$ a line bundle over $k$ that is a direct summand of $k \oplus L_0^\prime$. Now \eqref{BETCE} of Exc.~\ref{pr.ELIDREDQUAT}
\[
L_0^\ast \otimes L^{\otimes 2} \cong L_0^\ast \otimes L_c^{\otimes 2} \cong B_{12}(c) \cong B^\prime_{12}(c^\prime) \cong L_0^{\prime\ast} \otimes L^{\prime\otimes 2}, 
\] 
and \eqref{ISOCOND} holds.

(ii) $\Rightarrow$ (iii). This is clear since line bundles $L$ over $k$ satisfying the requirements of (ii) exist, e.g., $L := k \oplus \{0\} \subseteq k \oplus L_0$.

(iii) $\Rightarrow$ (i). Let $L$, $L^\prime$ be line bundles over $k$ that are direct summands of $k \oplus L_0$, $k \oplus L_0^\prime$, respectively, and satisfy \eqref{ISOCOND}. By Exc.~\ref{pr.ELIDREDQUAT}~(c), there exist elementary idempotents $c \in B$, $c^\prime \in B^\prime$ such that $L \cong L_c$, $L^\prime \cong L_{c^\prime}$. Now \eqref{ISOCOND} yields
\[
M :=L_0^\ast \otimes L_c^{\otimes 2} \cong L_0^{\prime\ast} \otimes L_{c^\prime}^{\otimes 2},
\]
and Exc.~\ref{pr.ELIDREDQUAT}~(a) produces isomorphisms $\Phi\:B \overset{\sim} \to C := \End_k(k \oplus M)$, $\Phi^\prime\:B^\prime \overset{\sim} \to C$ sending $c,c^\prime$ to $\left(\begin{smallmatrix}
1 & 0 \\
0 & 0
\end{smallmatrix}\right)$, respectively. In particular, $\Phi^{\prime -1} \circ \Phi$ is an isomorphism from $B$ to $B^\prime$, and (i) holds.
\end{sol}

\begin{sol}{pr.FQUATSPL} \label{sol.FQUATSPL}
We apply Exc.~\ref{pr.ISOREDQUAT} to $B := \End_k(k \oplus L_0)$, $L_0 := L$, and $B^\prime := \Mat_2(k) = \End_k(k \oplus L_0^\prime)$, $L_0^\prime := k$. Then $M := k \oplus \{0\} \subseteq k \oplus L_0$ is a free submodule of rank $1$ and a direct summand at the same time. If $B$ is split, then $B \cong B^\prime$ and condition (ii) of Exc.~\ref{pr.ISOREDQUAT} yields a line bundle $M^\prime \subseteq k \oplus L_0^\prime = k \oplus k$ which is a direct summand of $k \oplus k$ and satisfies $L_0 \otimes M^{\prime\otimes 2} \cong L_0^\prime \otimes M^{\otimes 2}$. This means $L_0 \cong M^{\prime\ast\otimes 2}$. Moreover, $k \oplus k \cong M^\prime \oplus M^{\prime\prime}$ for some line bundle $M^{\prime\prime}$, and taking determinants ($=$ second exterior powers), we obtain $M^\prime \otimes M^{\prime\prime} \cong k$. Hence $M^{\prime\ast} \cong M^{\prime\prime}$, being a homomorphic image of $k \oplus k$, is generated by two elements.

Conversely, suppose $L_0 = L \cong M^{\prime\otimes 2}$ for some line bundle $M^\prime$ on two generators over $k$. Then Exc.~\ref{pr.LITWO}~(b) implies $M^\prime \oplus M^{\prime\ast} \cong k \oplus k$. Since $L_0 \otimes M^{\prime\ast\otimes 2} \cong L \otimes M^{\prime\ast\otimes 2} \cong k \cong L_0^\prime \otimes M^{\otimes 2}$, we deduce from Exc.~\ref{pr.ISOREDQUAT} that $B$ and $B^\prime = \Mat_2(k)$ are isomorphic, i.e., $B$ is split.

Finally, suppose we are given a line bundle $L$ on two generators over $k$ that is not a square in $\Pic(k)$. (For example, take $k$ to be the ring of algebraic integers in $\IQ(\sqrt{-14})$.  Then $\Pic(k) \cong \IZ/4$ and we take $L$ to be a generator.  It is generated by 2 elements because $k$ is Dedekind.)  By what we have just shown, the quaternion algebra
\[
B := \End_k(k \oplus L) = \left(\begin{matrix}
k & L^\ast \\
L & k
\end{matrix}\right)
\]
is reduced but not split. It will be free as a $k$-module once we have shown that $L \oplus L^\ast$ is a free $k$-module. But this follows immediately from Exc.~\ref{pr.LITWO}~(b).
\end{sol}

\begin{sol}{pr.FAIWITCAN} \label{sol.FAIWITCAN}
By Exc.~\ref{pr.FQUATSPL}, the quaternion algebra
\[
B := \End_k(k \oplus L^\prime) = \left(\begin{matrix}
k & L^{\prime\ast} \\
L^\prime & k 
\end{matrix}\right) 
\] 
is split since $L^\prime$ is the square (in $\Pic(k)$) of some line bundle on two generators over $k$. Computing the norm of $B$ in two ways by means of (\ref{ss.TWIMA}.\ref{TWINO}), we deduce $\bfh \perp \bfh \cong n_B \cong \bfh \perp (-\bfh_{L^\prime}) \cong \bfh \perp \bfh_{L^\prime}$ and hence obtain the first assertion. As to the second, it suffices to consult Exc.~\ref{pr.SPLHYPLA}, which shows that the hyperbolic plane $\bfh_{L^\prime}$ is not split.
\end{sol}

\begin{sol}{pr.GRASSSPL} \label{sol.GRASSSPL}
Localizing if necessary, we may assume that $M$ is free (of rank $3$). Let $(e_i)_{1\leq i\leq 3}$ be a basis of $M$ and $(e_i^\ast)_{1\leq i\leq 3}$ the corresponding dual basis of $M^\ast$. Setting $\alpha := \theta(e_1 \wedge e_2 \wedge e_3)$, $\beta := \theta^{\ast -1}(e_1^\ast \wedge e_2^\ast \wedge e_3^\ast)$, we apply (\ref{ss.TWIZO}.\ref{DUVO}) to conclude $\alpha\beta = 1$. In particular, $\alpha$ and $\beta$ are both invertible.  

Letting $(ijl)$ vary over the cyclic permutations of $(123)$ and $s = 1,2,3$, we now apply (\ref{ss.TWIZO}.\ref{TRIVE}) to obtain $\la e_i \times_\theta e_j,e_s \ra =  \theta(e_i \wedge e_j \wedge e_s) = \alpha\delta_{ls} = \la \alpha e_l^\ast,e_s \ra$ and $\la e_s^\ast,e_i^\ast \times_\theta e_j^\ast \ra = \theta^{\ast -1}(e_i^\ast \wedge e_j^\ast \wedge e_s^\ast) = \beta\delta_{ls} = \la e_s^\ast,\beta e_l\ra$, which amounts to
\begin{align}
\label{TWIBA} e_i \times_\theta e_j = \alpha e_l^\ast, \quad e_i^\ast \times_\theta e_j^\ast = \beta e_l.
\end{align}
Now observe that our asserted equations are alternating in $u,v$ (resp. $u^\ast,v^\ast$).
For the first equation, we may therefore assume $u = e_i$, $v = e_j$, $w^\ast = e_s^\ast$. Then \eqref{TWIBA} yields
\begin{align*}
(e_i \times_\theta e_j) \times_\theta e_s^\ast = \alpha e_l^\ast \times_\theta e_s^\ast =\,\,&
\begin{cases}
\alpha\beta e_j & \text{for}\;\;s = i, \\
-\alpha\beta e_i & \text{for}\;\;s = j, \\
0 & \text{for}\,\,s = l  
\end{cases} \\
=\,\,&\begin{cases}
e_j & \text{for}\;\;s = i, \\
-e_i & \text{for}\;\;s = j, \\
0 & \text{for}\,\,s = l  
\end{cases} \\
=\,\,&\la e_s^\ast,e_i\ra e_j - \la e_s^\ast,e_j\ra e_i,
\end{align*} 
as claimed. Similarly, for the second equation, we may assume $u^\ast = e_i^\ast$, $v^\ast = e_j^\ast$, $w = e_s$ and obtain
\begin{align*}
(e_i^\ast \times_\theta e_j^\ast) \times_\theta e_s = \beta e_l \times_\theta e_s =\,\,&
\begin{cases}
\alpha\beta e_j^\ast & \text{for}\,\,s = i, \\
-\alpha\beta e_i^\ast & \text{for}\;\;s = j, \\
0 & \text{for}\;\;s = l
\end{cases} \\
=\,\,&\begin{cases}
e_j^\ast & \text{for}\,\,s = i, \\
-e_i^\ast & \text{for}\;\;s = j, \\
0 & \text{for}\;\;s = l
\end{cases} \\
=\,\,&\la e_i^\ast,e_s\ra e_j^\ast - \la e_j^\ast,e_s\ra e_i^\ast.
\end{align*}
This completes the proof.
\end{sol}

\begin{sol}{pr.OCTDEDDOM} \label{sol.OCTDEDDOM}
Let $C$ be a reduced octonion algebra over the Dedekind domain $k$. By Thm.~\ref{t.REDCOE}, there exist a finitely generated projective $k$-module $M$ of rank $3$ and an orientation $\theta$ of $M$ such that $C \cong \Zor(M,\theta)$ in the sense of \ref{ss.TWIZO}. By \cite[\S{VII.4.10}, Prop.~24]{MR0360549}, we can find an ideal $\mfa \subseteq k$ satisfying $M \cong k^2 \oplus \mfa$. But $\bigwedge^3(M) \cong k$ by means of $\theta$. By \cite[III, \S7.7, Cor. of Prop.~10]{MR0354207}, therefore, $\bigwedge^2(k^2) \otimes \mfa \cong k \otimes \mfa \cong \mfa$ is free of rank $1$, forcing $M$ to be free of rank $3$. Now Thm.~\ref{t.REDCOE} shows that $C \cong \Zor(M,\theta) \cong \Zor(k)$ is split.
\end{sol}


\solnsec{Section~\ref{s.NOREQ}}

\begin{sol}{pr.SKONOCO} \label{sol.SKONOCO}
The cases $s = 1$ and $s = r$ are obvious. Hence we may assume $1 < s < r$, which by Cor.~\ref{c.COMALGLOC} implies that $C,B,B^\prime$ are all regular. Arguing by induction on $r - s$, we let $\vph\:B \overset{\sim} \to B^\prime$ be an isomorphism. Then $C = B \perp B^\perp = B^\prime \perp B^{\prime\perp}$, hence
\[
n_B \perp n_C\vert_{B^\perp} \cong n_C \cong n_{B^\prime} \perp n_C\vert_{B^{\prime\perp}}
\]
On the other hand, $\vph$ is an isometry from $n_B$ onto $n_{B^ \prime}$. Thus, by Witt cancellation (Thm.~\ref{c.WITCAN}), $n_C\vert_{B^\perp}$ and $n_C\vert_{B^{\prime\perp}}$ are isometric. Since $s < r$, applying Lemma~\ref{l.QUAFORUNT} yields an element $l \in B^\perp$ such that $n_C(l) \in k^\times$, which in turn leads to an element $l^\prime \in B^{\prime\perp}$ satisfying $n_C(l^\prime) = n_C(l) = :-\mu$. Write $B_1$ (resp. $B_1^\prime$) for the subalgebra of $C$ generated by $B$ (resp. $B^\prime$) and $l$ (resp. $l^\prime$). Then Cor.~\ref{c.CDEMB} yields unique isomorphisms $h\:\Cay(B,\mu) = B \oplus Bj \overset{\sim} \to B_1$ (resp. $h^\prime\:\Cay(B^\prime,\mu) = B^\prime \oplus B^\prime j^\prime \overset{\sim} \to B_1^\prime$) extending the identity of $B$ (resp. $B^\prime$) and sending $j$ (resp. $j^\prime$) to $l$ (resp. $l^\prime$). The isomorphisms thus obtained fit into the diagram
\[
\xymatrix{\Cay(B,\mu) \ar[r]^(.6){h}_(.6){\cong} \ar[d]_{\Cay(B,\vph)}^{\cong} & B_1 \ar@{.>}[d]^{\vph_1}_{\cong} \\
\Cay(B^\prime,\mu) \ar[r]_(.6){h^\prime}^(.6){\cong} & B_1^\prime,}
\]
which can be completed uniquely to a commutative square by the dotted isomorphism $\vph_1\:B_1 \overset{\sim} \to B_1^\prime$ as indicated. Since $B_1$ and $B_1^\prime$ both have rank $2s$, the induction hypothesis applies to $\vph_1$ and completes the proof.  To conclude the solution to the problem, let $B \subseteq \IO$ be a non-zero subalgebra and let $0 \neq x \in B$. Then $n_\IO(x) \neq 0$ and $B$ contains the element $x^2 = t_\IO(x)x - n_\IO(x)1_\IO$, hence $1_\IO$. Thus $B$ is a unital subalgebra of $\IO$ on which $n_\IO$ continues to permit composition and to be anisotropic, hence regular as well since the characteristic is not $2$. Summing up, therefore, $B \subseteq \IO$ is a composition division subalgebra. From \ref{ss.REAL} we now conclude that $B$ is isomorphic to one of the algebras $\IR, \IC, \IH, \IO$, and the Skolem-Noether theorem completes the proof.  
\end{sol}

\begin{sol}{pr.UNINOR} \label{sol.UNINOR}
By Prop.~\ref{p.NORSIM}~(b), $f$ preserves not only norms and units, but also traces and conjugations.

(a) Applying (\ref{ss.IDCO}.\ref{QAUOP}), we obtain
\begin{align*}
f(U_xy) =\,\,&n_C(x,\bar y)f(x) - n_C(x)f(\bar y) = n_C\big(f(x),\overline{f(y)}\big)f(x) - n_C\big(f(x)\big)\overline{f(y)} \\  
=\,\,&U_{f(x)}f(y),
\end{align*}
hence \eqref{COMUOP}. In order to prove \eqref{ISAN}, we first note
\[
x^2 = U_x1_C, \quad \big(U_{x+y} - U_x - U_y\big)z = x(yz) + z(yx),
\]
so by \eqref{COMUOP}, $f$ preserves squares and the expression $x(yz) + z(yx)$. Abbreviating the left-hand side of \eqref{ISAN} by $A$ and using the fact that the associator of $C$ is alternating, we can now compute
\begin{align*}
A =\,\,&f(xy)^2 - f(xy)\big(f(y)f(x)\big) - \big(f(x)f(y)\big)f(xy) + f(x)f(y)^2f(x) \\
=\,\,&f\big((xy)^2\big) - f(xy)\big(f(y)f(x)\big) - f(x)\big(f(y)f(xy)\big) - [f(x),f(y),f(xy)] \\ 
\ \,\,&+ f(xy^2x) \\
=\,\,&[f(x),f(xy),f(y)] + f\Big((xy)^2 - (xy)(yx) - x\big(y(xy)\big) + xy^2x\Big) \\ 
=\,\,&[f(x),f(xy),f(y)]
\end{align*}
since, thanks to Artin's Theorem (Cor.~\ref{c.ARTIN}), the subalgebra of $C$ generated by two elements is associative. This completes the proof of \eqref{ISAN}. In order to complete the proof of (a), it will therefore suffice to show that the \emph{right}-hand side of \eqref{ISAN} is symmetric in $x,y$. To see this, we note that, since squares are preserved by $f$, so is the circle product. Hence $f(x \circ y) = f(x) \circ f(y)$, which implies 
\begin{align*}
[f(y),f(yx),f(x)] =\,\,&[f(y),f(x \circ y),f(x)] - [f(y),f(xy),f(x)] \\
=\,\,&[f(y),f(x) \circ f(y),f(y)] + [f(x),f(xy),f(y)] \\ 
=\,\,&[f(x),f(xy),f(y)],
\end{align*}
again by Artin's Theorem, and the assertion follows.

(b) Since $f$ preserves norms and traces, we have $t_{C^\prime}(c^\prime) = t_{C^\prime}(f(c)) = t_C(c) = 1$, $n_{C^\prime}(c^\prime) = n_{C^\prime}(f(c)) = n_C(c) = 0$, so $c^\prime \in C^\prime$ is an elementary idempotent. Now, as we have seen in (a), $f$ preserves the circle product: $f(x \circ y) = f(x) \circ f(y)$ for all $x,y \in C$. In particular, $f(c \circ y) = c^\prime \circ f(y)$. Hence, writing $C_{ij}$, $i,j = 1,2$, for the Peirce components of $C$ relative to $c$, it will be enough to show
\begin{align}
\label{CIROFF} C_{12} + C_{21} = \{x \in C\mid c \circ x = x\}.
\end{align}
For $x \in C$, we let $x = x_{11} + x_{12} + x_{21} + x_{22}$ be its Peirce decomposition relative to $c$. Then
\[
c \circ x = cx + xc =x_{11} + x_{12} +x_{11} +x_{21} = 2x_{11} + x_{12} + x_{21},
\]
and comparing Peirce components we see that $c \circ x = x$ if and only if $x_{11} = x_{22} = 0$. Hence \eqref{CIROFF} is proved.
\end{sol}

\begin{sol}{pr.UNORQUA} \label{sol.UNORQUA}
(a) We begin by assuming that the case of a local ring has been settled and then let $k$ be arbitrary. We put $X := \Spec(k)$ and deduce from 
Exc.~\ref{pr.LOCHOM}~(b) that
\begin{align*}
X_+ :=\,\,&\{\mfp \in X \mid f_\mfp\:B_\mfp \to \pB_\mfp\;\text{is an isomorphism}\}, \\
X_- :=\,\,&\{\mfp \in X \mid f_\mfp\:B_\mfp \to \pB_\mfp\;\text{is an anti-isomorphism}\} 
\end{align*}
are Zariski-open subsets of $X$. Since a quaternion algebra is not commutative, they are also disjoint, and our assumption implies that they cover $X$. Hence Exc.~\ref{pr.IDEPAR}, yields a complete orthogonal system $(\vep_+,\vep_-)$ of idempotents in $k$ satisfying $X_\pm = D(\vep_\pm)$. Now it suffices to put $k_\pm = k\vep_\pm$ and to invoke (\ref{ss.PROPID}.\ref{SPEPI}), which implies for any $\mfp_+ \in \Spec(k_+)$ that $\mfp := \mfp_+ \times k_- \in D(\vep_+)$ makes $f_{+\mfp_+} = (f_\mfp)_{k_{+\mfp_+}}\:B_{+\mfp_+} \to B^\prime_{+\mfp_+}$ (by (\ref{ss.LOPROM}.\ref{MPRLO})) an isomorphism. Similarly, for $\mfp_- \in \Spec(k_-)$, $f_{-\mfp_-}\:B_{-\mfp_-} \to B^\prime_{-\mfp_-}$ turns out to be an anti-isomorphism. Hence $f_+$ is an isomorphism and $f_-$ is an anti-isomorphism.

We are thus left with the case that $k$ is a local ring and must show that $f$ is either an isomorphism or an anti-isomorphism. Thm.~\ref{t.COMSUBALG}~(b) and Exc.~\ref{pr.ETAUT}~(a) yield a quadratic \'etale subalgebra $D = k[u] \subseteq B$, for some $u \in B$ having trace $1$. Since $f$ preserves units, norms, and traces by Prop.~\ref{p.NORSIM}, we conclude that $\pD := f(D) = k[\pu]$, $\pu := f(u)$, is a quadratic \'etale subalgebra of $\pB$ (Prop.~\ref{p.SMASU}), and $f\vert_D\:D \to \pD$ is an isomorphism (Prop.~\ref{p.ISNEQ}). Now apply Thm.~\ref{t.COMSUBALG}~(a) to reach $B$ from $D$ by means of the Cayley-Dickson construction: there exists a unit $\mu \in k^\times$ such that the inclusion $D \hookrightarrow B$ extends to an identification $B = \Cay(D,\mu) = D \oplus Dj$, $j \in D^\perp$, $n_B(j) = -\mu$. Setting $\pj := f(j) \in D^{\prime\perp} \subseteq \pB$, we obtain $n_{\pB}(\pj) = -\mu$, and from Prop.~\ref{p.CDUNV} we conclude that $f\vert_D$ extends to a homomorphism $g\:B \to \pB$ of conic algebras satisfying $g(j) = \pj$. By Cor.~\ref{c.CDEMB}, therefore, $g$ is an isomorphism from $B$ onto the subalgebra of $\pB$ generated by $\pD$ and $\pj$. Counting ranks we conclude that $g$ is in fact an isomorphism from $B$ onto $\pB$. Hence $f_1 := g^{-1} \circ f\:B \to B$ is a unital norm equivalence inducing the identity on $D$ and satisfying $f_1(j) = j$. Since $f_1$ stabilizes $D^\perp = Dj$, we find a $k$-linear bijection $\vph\:D \to D$ such that $f_1(vj) = \vph(v)j$ for all $v \in D$. Then $\vph(1_D) = 1_D$, and since $n_B$ permits composition, $\vph$ leaves $n_D$ invariant and is thus a unital norm equivalence of $D$, hence an automorphism (Prop.~\ref{p.ISNEQ}~(b)). By Exc.~\ref{pr.ETAUT}, therefore, we are left with two cases depending on $\vph$.

Suppose first that $\vph = \Eins_D$. Then $f_1 = \Eins_B$, and $f = g\:B \to \pB$ is an isomorphism. 

Suppose next that $\vph = \iota_D$. By 
Exc.~\ref{pr.CDSCAPAR}, the map $\psi\:B \to B$ defined by $\psi(v + wj) := v - wj$ for all $v,w \in D$ is an automorphism, and one checks that $f_1 = \psi \circ \Int(j) \circ \iota_B$, where $\Int(j)$ stands for the inner automorphism $x \mapsto jxj^{-1}$ of $B$ affected by $j$. Hence $f = g \circ \psi \circ \Int(j) \circ \iota_B\:B \to \pB$ is an anti-isomorphism.

(b) The implications (i) $\Rightarrow$ (ii) $\Rightarrow$ (iii) are obvious, so let us prove (i) under the assumption (iii). By Cor.~\ref{c.NORSIM}, there exists a unital norm equivalence $f\:B \to \pB$. With the notation of \eqref{QUATPM}, we apply (a) to conclude that $g := f_+ \times (f_- \circ \iota_{B_-})\:B \to \pB$ is an isomorphism.
\end{sol}

\begin{sol}{pr.ELIDNO} \label{sol.ELIDNO}
By Prop.~\ref{p.ISNEQ}~(a), every automorphism of $B$ is a unital norm equivalence. Hence it suffices to show that, conversely, if $f\:B \to B$ is a unital norm equivalence sending $e_{11}$ to $c$, then an automorphism of $B$ exists having the same property. From Exc.~\ref{pr.UNINOR}~(b) we deduce that $c$ is an elementary idempotent in $B$. Hence Exc.~\ref{pr.ELIDREDQUAT}~(a) yields a line bundle $L$ over $k$ and an isomorphism
\begin{align*}
\Phi\:B \overset{\sim} \longrightarrow B^\prime := \End_k(k \oplus L) = \left(\begin{matrix}
k & L^\ast \\
L & k 
\end{matrix}\right) \quad \text{such that} \quad \Phi(c) = \left(\begin{matrix}
1 & 0 \\
0 & 0 
\end{matrix}\right). 
\end{align*}
Viewing $B_{12}(c) \oplus B_{21}(c)$ (resp. $L^\ast \oplus L$) canonically as quadratic submodules of $(B,n_B)$ (resp. $(B^\prime,n_{B^\prime})$,
\[
\Phi\:B_{12}(c) \oplus B_{21}(c) \overset{\sim} \longrightarrow L^\ast \oplus L
\]
is an isometry and by (\ref{ss.TWIMA}.\ref{TWINO}), (\ref{ss.HYSP}.\ref{HYPA}), there is a natural identification of $L^\ast \oplus L$ with the hyperbolic plane $\bfh_L$. On the other hand, $f$, being a unital norm equivalence, by Exc.~\ref{pr.UNINOR}~(b) induces an isometry from the split hyperbolic plane $k \oplus k \cong B_{12}(e_{11}) \oplus B_{21}(e_{11})$ onto $B_{12}(c) \oplus B_{21}(c)$. Hence the hyperbolic plane $\bfh_L$ is split, which by Exc.~\ref{pr.SPLHYPLA} forces the line bundle $L$ to be free of rank $1$. Thus we may identify $L =k$ and in this way view $\Phi$ as an automorphism of $B$ sending $e_{11}$ to $c$.
\end{sol}

\begin{sol}{pr.OCTNOREQU} \label{sol.OCTNOREQU}
By Exc.~\ref{pr.UNITSTRGRALT}, $\vph := L_pR_{p^{-1}}$ is an isomorphism from $C$ to $C^q$ with $q := p^{-3}$. 

(a) By Prop.~\ref{p.NORSIM}~(a), $\vph$ is a norm similarity fixing $1_C$ and hence a unital norm equivalence.

(b) $\vph$ is an automorphism of $C$ if and only if $C = C^q$, which by (\ref{ss.UNTISOTALT}.\ref{EQIS}) is equivalent to $q \in \Nuc(C) = k1_C$ (Exc.~\ref{pr.AZUQUAD}~(b)), hence to $p^3 = q^{-1} \in k1_C$.

(c) Assume $\vph$ is an anti-automorphism of $C$. Then $C^{\op} = C^q$, so we have $yx = (xq^{-1})(qy)$, equivalently, $y(xq) = x(qy)$, for all $x,y \in C$. Setting $x = 1_C$ gives $yq = qy$ for all $y \in C$. By Exc.~\ref{pr.AZUQUAD}~(b), this implies $q \in k1_C$ and then $xy = yx$ for all $x,y \in C$, a contradiction.
\end{sol}

\begin{sol}{pr.ISTIS} \label{sol.ISTIS}
Let $p,q \in C^\times$. By Exc.~\ref{pr.ISTQALT}, $C^{(p,q)}$ is a multiplicative conic alternative $k$-algebra with norm $n_{C^{(p,q)}} = n_C(pq)n_C$. By multiplicativity, $L_{pq}$ is an isometry from $n_{C^{(p,q)}}$ to $n_C$. Hence $C^{(p,q)}$ is a composition algebra which is regular if $C$ is. Moreover, if $C$ is associative, then $L_{pq}\:C^{(p,q)} \to C$ is an isomorphism. And finally, if $k$ is LG, the norm equivalence theorem \ref{t.NOREQ} shows that $C$ and $C^{(p,q)}$ are isomorphic.
\end{sol}

\begin{sol}{pr.ZEPA} \label{sol.ZEPA}
Let $x_i = u_i + v_ij$, $u_i,v_i \in \IO$, $i = 1,2$ be non-zero elements of $\IS$, Then Exc.~\ref{pr.CDDIVALG}~(a) for $\mu = -1$ shows that $x_1x_2 = 0$ if and only if $u_i \neq 0 \neq v_i$ for $i = 1,2$ and
\begin{align}
\label{ZEROC} n_\IO(u_1) = n_\IO(v_1), \quad (u_1u_2)\bar v_1 = -u_1(u_2\bar v_1), \quad v_2 = -(v_1\bar u_2)u_1^{-1}.
\end{align}
Now let $(w_i)_{1\leq i\leq 7}$ be a Cartan-Shouten basis of $\IO$ in the sense of \ref{ss.DECS}. Combining (\ref{p.CDCON}.\ref{CDMULT}) with Fig.~\vref{fig.fano}, we compute
\begin{align*}
(w_1 + w_3j)(w_2 - w_6j) =\,\,&(w_1w_2 - (-\bar w_6)w_3) + (w_3\bar w_2 - w_6w_1)j \\
=\,\,&(w_1w_2 - w_6w_3) + (-w_3w_2 - w_6w_1)j,
\end{align*}
where $w_6w_3 = w_4 = w_1w_2$ and $w_6w_1 = w_5 = -w_3w_2$. Thus the quantities $a := w_1 + w_3j$, $b := w_2 - w_6j$ have $\IS$-norm $2$ and satisfy $ab = 0$, so we conclude $(a,b) \in \Zer(\IS)$, while \eqref{ZEROC} (or the property of $(w_i)$ being a Cartan-Shouten basis) yields
\begin{align}
\label{DUSIX} w_6 = (w_3\bar w_2)w_1^{-1} = (w_3w_2)w_1.
\end{align}
We must show that the action of $G$ on $\Zer(\IS)$ is simply transitive. Let $(x_1,x_2) \in \Zer(\IS)$ and write $x_i = u_i + v_ij$ with $u_i,v_i \in \IO$ for $i = 1,2$ as before. Then \eqref{ZEROC} shows $n_\IO(u_i) = n_\IO(v_i)$ for $i = 1,2$, whence the relations $n_\IS(x_i) = 2$ imply $n_\IO(u_i) = n_\IO(v_i) = 1$. From Exc.~\ref{pr.CDDIVALG} and (\ref{ss.BASID}.\ref{CAQUADL}) we deduce $0 = x := u_1 \circ u_2 = \alpha 1_\IO + \alpha_1u_1 + \alpha_2u_2$, $\alpha := -n_\IO(u_1,u_2)$,  $\alpha_i := t_\IO(u_{3-i})$, $i = 1,2$, hence $0 = u_1 \circ x =  2\alpha u_1 + 2\alpha_1u_1^2 = 2u_1(\alpha 1_\IO  + \alpha_1u_1)$. Thus $0 = y := \alpha 1_\IO + \alpha_1u_1$, which implies $0 = y \circ u_2 = 2\alpha 1_\IO$, and we conclude $0 = \alpha = \alpha_1 = \alpha_2$. Summing up, we have shown that $(u_1,u_2)$ is an orthonormal system in the euclidean space $\IO^0$. Consulting Exc.~\ref{pr.CDDIVALG}~(b) again, we also get $n_\IO(u_1u_2,v_1) = 0$. By Exc.~\ref{pr.CHACS}, therefore, the quantities $u_1,u_2,v_1$ can be extended to a Cartan-Shouten basis of $\IO$. Hence there exists an automorphism $\sigma$ of $\IO$ sending $u_1,u_2,v_1$ respectively to $w_1,w_2,w_3$. Consulting \eqref{ZEROC}, \eqref{DUSIX}, we conclude that $\sigma$ also sends $v_2$ to $-w_6$. But this means $\sigma((x_1,x_2)) = (a,b)$, and we have proved that the action is transitive. It remains to show that $\sigma(a) = a$, $\sigma(b) = b$ implies $\sigma = \Eins_\IO$. But this is clear since we then have $\sigma(w_i) = w_i$ for $i = 1,2,3$, and $w_1,w_2,w_3$ by Exc.~\ref{pr.CHACS} generate the octonion algebra $\IO$. 
\end{sol}

\begin{sol}{pr.ZOFI} \label{sol.ZOFI}
We begin by proving the following general statement. 
\begin{claim*}
Let $(V,Q)$ be a hyperbolic quadratic space of dimension $2n$ over $\IF_q$. Then the number of anisotropic vectors in $(V,Q)$ is
\[
q^{n-1}(q - 1)(q^n - 1).
\]
\end{claim*}

\begin{proof}
We identify $(V,Q) = (\IF_q^n \times \IF_q^n,(x,y) \mapsto x^\trans y)$ and note for $x \in \IF_q^n$ that
\[
\vert \{y \in \IF_q^n\mid x^\trans y = 0\}\vert =
\begin{cases}
q^n & \text{for $x = 0$,} \\
q^{n-1} & \text{for $x \neq 0$}.
\end{cases}
\]
Thus $\vert\{v \in V \mid Q(v) = 0\}\vert = q^n + (q^n - 1)q^{n-1}$, and we conclude that $(V,Q)$ contains precisely 
\begin{align*}
q^{2n} - q^n - (q^n - 1)q^{n-1} &= q^n(q^n - 1) - (q^n - 1)q^{n-1} = (q^n - q^{n-1})(q^n - 1)\\
&= q^{n-1}(q - 1)(q^n - 1)
\end{align*}
anisotropic vectors.
\end{proof}

We now turn to the assertions of our problem. Since the anisotropic vectors of $C$ relative to $n_C$ are just its invertible elements (Prop.~\ref{p.CHIN}), and the norm of the split octonions is hyperbolic (Cor.~\ref{c.SPLIDI}), the preceding claim reduces to the first equation for $n = 3$. 

To prove the second equation, we note that $n\:C^\times \to \IF_q^\times$, $x \mapsto n(x) := n_C(x)$, is a surjective multiplicative map. Writing $N :=\{x \in C\mid n_C(x) = 1\}$ for its ``kernel'', we let $x,y \in C^\times$, use Prop.~\ref{p.INVALT} and obtain
\[
y \in n^{-1}\big(n(x)\big) \Leftrightarrow n(y) = n(x) \Leftrightarrow n(x^{-1}y) = 1 \Leftrightarrow x^{-1}y \in N \Leftrightarrow y \in xN,
\] 
so the fibers of $n$ all have the same cardinality. This shows $\vert N\vert = \vert C^\times\vert/\vert\IF_q^\times\vert = q^3(q^4 - 1)$. 
\end{sol}

\begin{sol}{pr.MIVE} \label{sol.MIVE}
(a) is trivial: if $x \in M$ is not minimal, then $x = y + z$ for some $y,z \in M$ such that $Q(y) < Q(x)$ and $Q(z) < Q(x)$. Proceeding with $y,z$ in the same manner, we eventually arrive at a decomposition of $x$ as a finite sum of minimal elements.

(b) Again this is trivial, demanding not even the \emph{sketch} of a proof.

(c) We begin with a simple lemma. 

\begin{lem*}
(i) {If $N \subseteq M$ is a submodule, then $\Min(M,Q) \cap N \subseteq \Min(N,Q\vert_N)$.} 
\begin{holgerenum2}[(i),start=2]
\item {If $M = N \perp N^\prime$ is a decomposition of $M$ into the orthogonal sum (relative to $Q$) of two submodules $N,N^\prime \subseteq M$, then $\Min(M,Q) = \Min(N,Q\vert_N) \cup \Min(N^\prime,Q\vert_{N^\prime})$}.
\end{holgerenum2}
\end{lem*}

\begin{proof}
In order to prove (i), let $x \in \Min(M,Q) \cap N$ and suppose we have a decomposition $x = y + z$ with $y,z \in N$, $Q(y) < Q(x)$, $Q(z) < Q(x)$. Then this decomposition takes place also in $M$, contradicting minimality of $x$ in $(M,Q)$. Thus $x \in \Min(N,Q\vert_N)$. But note for $x \in \Min(N,Q\vert_N)$, that it is conceivable to have a decomposition $x = y + z$ with $y,z \in M \setminus N$ and $Q(y) < Q(x)$, $Q(z) < Q(x)$, so in general we won't have equality in (i).

In order to prove (ii), let $x \in \Min(M,Q)$ and write $x = y + y^\prime$ with $y \in N$, $y^\prime \in N^\prime$. Then $Q(x) = Q(y) + Q(y^\prime)$, and since $x$ is minimal in $M$, we conclude $Q(y) = 0$ or $Q(y^\prime) = 0$, hence $y = 0$ or $y^\prime = 0$. Thus $\Min(M,Q) \subseteq N \cup N^\prime$, which implies
\begin{align*}
\Min(M,Q) &= \Min(M,Q) \cap (N \cup N^\prime) \\
&= \left(\Min(M,Q) \cap N\right) \cup \left(\Min(M,Q) \cap N^\prime\right).
\end{align*}
By symmetry and (i), it therefore suffices to show $\Min(N,Q\vert_N) \subseteq \Min(M,Q)$. Let $x \in \Min(N,Q\vert_N )$ and suppose $x = y + z$ for some $y,z \in M$, $Q(y) < Q(x)$, $Q(z) < Q(x)$. Write $y = u + u^\prime$, $z = v + v^\prime$ with $u,v \in N$, $u^\prime,v^\prime \in N^\prime$. Since $x$ belongs to $N$, this implies $x = u + v$ and $Q(u) \leq Q(y) < Q(x)$, $Q(v) \leq Q(z) < Q(x)$, in contradiction to $x$ being a minimal vector of $(N,Q\vert_N)$. Thus $x \in \Min(M,Q)$. 
\end{proof}

We can now prove (c). Let $x,y \in \Min(M,Q)$ be inequivalent. For $u \in [x]$, $v \in [y]$, the assumption $Q(u,v) \neq 0$ would imply that $u,v \in \Min(M,Q)$ are equivalent. Since $M_{[x]}$ is spanned as a $\IZ$-module by the minimal elements of $(M,Q)$ belonging to $[x]$, ditto for $M_{[y]}$,  this contradiction shows that $M_{[x]}$ and $M_{[y]}$ are orthogonal. But every $x \in \Min(M,Q)$ belongs to $M_{[x]}$, and $M$ is spanned by $\Min(M,Q)$ as a $\IZ$-module. Thus $M$ is the orthogonal sum of the distinct $M_{[x]}$'s, $x \in \Min(M,Q)$, and it remains to show that each $M_{[x]}$ is indecomposable. Thanks to the obvious extension of (b) in the preceding lemma to more than two orthogonal summands, it will actually be enough to show that every positive definite integral quadratic module $(M,Q)$ all of whose minimal elements are equivalent is indecomposable. Indeed, suppose we have an orthogonal decomposition $M = N \perp N^\prime$ relative to $Q$, with submodules $N,N^\prime \subseteq M$, and let $x,y \in \Min(M,Q)$. By part (b) of the lemma, we may assume $x \in \Min(N,Q\vert_N)$, and by definition, there is a finite sequence $x = x_0,x_1,\dots,x_{k-1},x_k = y$ of minimal vectors in $(M,Q)$ such that $Q(x_{i-1},x_i) \neq 0$ for $1 \leq i \leq k$. We claim $x_i \in \Min(N,Q\vert_N)$ for all $i = 0,\dots,k$ and argue by induction. For $i = 0$, there is nothing to prove. If $i > 0$ and the assertion holds for $i - 1$, then $x_i \notin N^\prime$ since $x_{i-1} \in N$ by the induction hypothesis and $Q(x_{i-1},x_i) \neq 0$, whence aprt (b) of the lemma implies $x_i \in \Min(N,Q\vert_N)$, and the induction is complete. In particular, $y \in \Min(N,Q\vert_N)$, and we have shown $\Min(M,Q) = \Min(N,Q\vert_N)$, which obviously implies $M = N$, $N^\prime = \{0\}$. 
\end{sol}

\begin{sol}{pr.ODDI} \label{sol.ODDI}
We put $C^\dag = C/2C = C \otimes_\IZ \IF_2$ as a conic algebra over $\IF_2$ and write $x \mapsto x^\dag$ for the natural epimorphism from $C$ to $C^\dag$. This epimorphism canonically induces a map from $C^\times$ to  $C^{\dag\times}$, which may or may not be surjective. Note that, since the discriminant of $C$ is odd, $C^\dag$ is a composition algebra of dimension $8$ over $\IF_2 $, hence an octonion algebra and thus isomorphic to $\Zor(\IF_2)$ (Cor.~\ref{c.SPLIDI}).  

\begin{lem*} 
For $x,y \in C^\times$ we have $x^\dag = y^\dag$ if and only if $y = \pm x$.
\end{lem*}

\begin{proof}
$y = \pm x$ clearly implies $x^\dag = y^\dag$. Conversely, let this be so. Then $y = x + 2z$ for some $z \in C$. If $z = 0$, we are done, so we may assume $z \neq 0$. Since $x,y \in C$ are both units, we obtain $1 = n_C(x + 2z) = 1 + 2n_C(x,z) + 4n_C(z)$, hence $n_C(x,z) = -2n_C(z)$, and applying the Cauchy-Schwarz inequality yields
\[
4n_C(z)^2 = n_C(x,z)^2 \leq n_C(x,x)n_C(z,z) = 4n_C(z).
\]
But $n_C(z)$ is a positive integer, forcing $n_C(z) = 1$. Thus the above inequality is, in fact, an equality, which implies $z = \alpha x$ for some $\alpha \in \IQ$. Taking norms, we deduce $\alpha^2 = 1$, hence $\alpha = \pm 1$, where the assumption $\alpha = 1$ would lead to the contradiction that $y = 3x$ is not invertible in $C$. Hence $\alpha = -1$ and $y = -x$.
\end{proof}

We can now establish the first part of the problem. By Exc.~\ref{pr.ZOFI} for $q = 2$, we have $\vert C^{\dag\times}\vert = 120$, while by the lemma, the fiber of an arbitrary element $u \in C^{\dag\times}$ under the natural map $C^\times \to C^{\dag\times}$ is either empty or consists of two elements. Hence $\vert C^\times\vert \leq 240$.

For the rest of the proof, we may assume that $\vert C^\times\vert = 240$. Writing $C^{\times\dagger}$ for the image of $C^\times$ under the natural map $C^\times \to C^{\dag\times}$, and putting $r := \vert C^{\times\dagger}\vert$, we note that the fiber of any point in $C^{\times\dag}$ consists of two elements. Thus $2r = 240$, hence $r = 120 = \vert C^{\dag\times}\vert$, and we conclude $C^{\times\dagger} = C^{\dag\times}$. In other words, \emph{the natural map $C^\times \to C^{\dag\times}$ is surjective.} 

We can now prove that the positive definite integral quadratic lattice $(C,n_C)$ is indecomposable. We have $C^\times = \{x \in C \mid n_C(x) = 1\} \subseteq \Min(C,n_C)$, and for $x,y \in C^\times$, Exc.~\ref{pr.LIQUA}~(b) yields a sequence $x^\dag = x_0^\prime,x_1^\prime,\dots,x_k^\prime = y^\dag$ of elements in $C^{\dag\times}$ such that $k \leq 2$ and $n_{C^\dag}(x_{i-1}^\prime,x_i^\prime) \neq 0$ for $1 \leq i \leq k$. Since the natural map $C^\times \to C^{\dag\times}$ is surjective, each $x_i^\prime$, $1 \leq i < k$, lifts to an invertible element $x_i \in C$, so we have got a sequence $x = :x_0,x_1,\dots,x_k := y$ of elements in $C^\times$, hence of minimal vectors in $(C,n_C)$, such that $n_C(x_{i-1},x_i)$ is an odd integer and, in particular, different from $0$, for $1 \leq i \leq k$. This shows that the minimal vectors $x$ and $y$ of $(C,n_C)$ are equivalent, in other words, the invertible elements of $C$ all belong to a single equivalence class of minimal vectors, which may also be expressed by saying that $C^\times \subseteq C_{[1_C]}$. Now suppose that $x \in \Min(C,n_C)$ is \emph{not} equivalent  to $1_C$.  From Exc.~\ref{pr.MIVE}~(c) we deduce $n_C(C^\times,x) = \{0\}$, which  implies $n_{C^\dag}(C^{\dag\times},x^\dag) = n_{C^\dag}(C^{\times\dagger},x^\dag) = \{0\}$. But $C^\prime$ is an octonion algebra over $\IF_2$ that is spanned as a vector space by $C^{\dag\times}$ (Exc.~\ref{pr.LIQUA}~(a) and Prop.~\ref{p.CHIN}). Thus $x^\dag = 0$, hence $x \in 2C$, and we conclude that the indecomposable submodules $C_{[x]} \subseteq C$, with $x \in \Min(C,n_C)$ not equivalent to $1_C$, all belong to $2C$; in particular, they have even discriminant. But the discriminant of $(C,n_C)$ is odd by hypothesis, which, thanks to Exc.~\ref{pr.MIVE}~(c), implies that $(C,n_C) = C_{[1_C]}$ is indeed indecomposable. 
\end{sol}

\begin{sol}{pr.QEZE} \label{sol.EZE}
Let $R$ be a quadratic \'etale $\IZ$-algebra. Then $R$ is a free $\IZ$-module of rank $2$, and since $1_R \in R$ is unimodular, it can be extended to a basis $(1_R,u)$ of $R$. By Lemma~\ref{l.COMTRI}, the trace form of $R$ is surjective. Since $t_R(1_R) = 2$, we conclude that $t_C(u) = 2n + 1$ (for some integer $n$) is odd. In fact, replacing $u$ by $u - n\cdot 1_R$, we may assume $t_R(u) = 1$. Now Prop.~\ref{p.SMASU} shows $1 - 4n_C(u) = \pm 1$, which implies $n_C(u) = 0$ or $n_C(u) = \frac{1}{2}$. The latter option obviously being impossible, we conclude that $R$ has zero divisors (Exc.~\ref{pr.ZERDIV}), hence is split (Cor.~\ref{ss.OCTINT}). This solves the first part of the problem.

As to the second, assume $\Hur(\IH) \cong \Cay(R,\mu)$ for some quadratic $\IZ$-algebra $R$ and some $\mu \in \IZ$. Combining the definition of the discriminant (\ref{ss.DISC}) with (\ref{t.HURO}.\ref{DISCH}) and Remark~\ref{r.NOCD}, we conclude with $d := \disc(R)$ that $4 = d^2\mu^2$. On the other hand, since the norm of $R \subseteq \Hur(\IH)$ is positive definite, and $R$ itself by the first part of this exercise is not \'etale, we have $d \geq 2$ and hence $d = 2$, $\mu = -1$. Extending $1_R$ to a basis $(1_R,u)$ of $R$ over $\IZ$, we thus obtain $2 = t_C(u)^2 - 4n_R(u)$ from Prop.~\ref{p.SMASU}, where the left-hand side is $\equiv 2 \bmod 4$, while the right-hand is $\equiv 0$ or $\equiv 1 \bmod 4$, a contradiction.
\end{sol}

\begin{sol}{pr.OCZE} \label{sol.OCZE}
(a) Since $C$ is not split, it has no zero divisors (Cor.~\ref{ss.OCTINT}). Hence neither has $C_\IQ$, whence $C_\IQ$ is an octonion division algebra over the rationals, up to isomorphism actually the only one (Cor.~\ref{c.NUMOC}). This shows $C_\IQ\cong \Cay(\IQ;-1,-1,-1)$. In particular, the norm of $C_\IQ$ is positive definite. Hence so is the norm of $C$. Summing up we have shown that $(C,n_C)$ is a positive definite inner product space of rank $8$ over $\IZ$. But by \ref{ss.VIPO}~(c), up to isometry there is only one. Hence, by Cor.~\ref{c.NORSIM}, there exists a unital norm equivalence $f\:C \to \Cox(\IO)$. Transferring the unit element, norm  and multiplication of $C$ to $\Cox(\IO)$ by means of $f$, we may actually assume $C = \Cox(\IO)$ as additive groups with $1_C = 1_\IO$ and $n_C = n_{\Cox(\IO)}$. In particular $\vert C^\times\vert = 240$ (Exc.~\ref{pr.ROOC}~(d)).

(b) Put $D := \Cox(\IO)$. Then $C^\dag := C/2C = C \otimes_\IZ \IF_2$ and $D^\dag = D/2D = D \otimes_\IZ \IF_2$ are both octonion algebras over $\IF_2$, hence split (\ref{ss.FINFIE}) and, in particular, isomorphic. We therefore find an isomorphism $\psi\:C^\dag \overset{\sim} \to D^\dag$ of $\IF_2$-algebras. Note by the normalization provided in (a) that $(C^\dag,n_{C^\dag},1_{C^\dag}) = (D^\dag,n_{D^\dag},1_{D^\dag})$ as ``pointed'' quadratic spaces over $\IF_2$. In particular, $\psi$ is an orthogonal transformation of the quadratic space $(C^\dag,n_{C^\dag})$. Applying \cite[Prop.~14]{MR0024439}, we see that the orthogonal group of this space is generated by the orthogonal transvections
\begin{align*}
\tau_a\: \quad x \longmapsto x - \frac{n_{C^\dag}(a,x)}{n_{C^\dag}(a)}a &&(a \in  C^{\dag\times}),
\end{align*}
each of which can be lifted to an orthogonal transformation from $(C,n_C)$ onto itself since, as we have seen in the solution to Exc.~\ref{pr.ODDI}, the natural map from $C^\times$ to $C^{\dag\times}$ is surjective. Hence we can also lift $\psi$ to an orthogonal transformation $\vph$ of $(C,n_C)$. In other words, $\vph\:C \to \Cox(\IO)$ is a linear bijection preserving norms such that the diagram
\[
\xymatrix{C \ar[r]^{\vph} \ar[d]_{\mathrm{can}} & \Cox(\IO) \ar[d]^{\mathrm{can}} \\
C/2C \ar[r]_(.33){\psi} & \Cox(\IO)/2\Cox(\IO),}
\]
commutes. Writing $x^\dag := \can(x)$ for $x \in C^\times$ and using the fact that $\psi$ is an algebra isomorphism, we therefore obtain, for all $x,y \in C^\times$,
\begin{align*}
\big(\vph(x\cdot y)\big)^ \dag =\,\,&\psi\big((x\cdot y)^\dag\big) = \psi(x^\dag\cdot y^\dag) = \psi(x^\dag)\psi(y^\dag) = \big(\vph(x)\big)^\dag\big(\vph(y)\big)^\dag = \big(\vph(x)\vph(y)\big)^\dag.
\end{align*}
But for $x \in C^\times$, the fiber of $x^\dag \in C^{\dag\times}$ under the vertical arrows in the above diagram consists of the elements $\pm x$. Thus $\vph(x\cdot y) = \pm\vph(x)\vph(y)$.

(c) According to (b), there exists a map $\vep\:C^\times \times C^\times \to \{\pm 1\}$ such that $\vph(x\cdot y) = \vep(x,y)\vph(x)\vph(y)$ for all $x,y \in C^\times$. We claim that \emph{$\vep$ is constant}. Indeed, let $x,y,z \in C^\times$. Since $\vph$ preserves norms, and these are the same for $C$ and $\Cox(\IO)$, an application of (\ref{ss.MULCO}.\ref{MUPLR}) yields
\begin{align*}
\vep(x,y)\vep(x,z)n_C(x)n_C(y,z) =\,\,&\vep(x,y)\vep(x,z)n_C\big(\vph(x)\big)n_C\big(\vph(y),\vph(z)\big) \\
=\,\,&\vep(x,y)\vep(x,z)n_C\big(\vph(x)\vph(y),\vph(x)\vph(z)\big) \\
=\,\,& n_C\big(\vph(x\cdot y),\vph(x\cdot z)\big) \\
=\,\,&n_C(x\cdot y,x\cdot z) = n_C(x)n_C(y,z).
\end{align*}
Hence $\vep(x,y) = \vep(x,z)$ provided $n_C(y,z) \neq 0$. On the other hand, if $n_C(y,z) = 0$, then Exc.~\ref{pr.LIQUA}~(b) combined with the surjectivity of the natural map $C^\times \to C^{\dag\times}$ yields a $w \in C^\times$ satisfying $n_C(y,w) \neq 0 \neq n_C(w,z)$, which in turn implies $\vep(x,y) = \vep(x,w) = \vep (x,z)$. Thus $\vep(x,y)$ does not depend on $y \in C^\times$. Interchanging the role of $x$ and $y$, and replacing (\ref{ss.MULCO}.\ref{MUPLR}) by (\ref{ss.MULCO}.\ref{MUPLL}) in the process, the same argument shows that $\vep(x,y)$ is also independent of $x \in C^\times$, showing that the map $\vep$ is indeed constant. Finally, replacing $\vph$ by $-\vph$ if necessary,  we may assume $\vph(x\cdot y) = \vph(x)\vph(y)$ for all $x,y \in C^\times = \Cox(\IO)^\times$. But since $C$ is generated by $C^\times$ as an additive group (Exc.~\ref{pr.ROOC})~(b),(d)), this means that $\vph\:C \overset{\sim} \to \Cox(\IO)$ is an isomorphism of octonion algebras.
\end{sol}

\begin{sol}{pr.SPLIFICO} \label{sol.SPLIFICO}
(a) We begin by proving sufficiency The assertion is clear if (i) or (ii) holds. As to (iii), assuming there exists an $F$-embedding $K \hookrightarrow C$, we obtain an induced $K$-embedding $K_K \hookrightarrow C_K$. If we can show that $K_K$ is not a division algebra, then neither is $C_K$, forcing $C_K$ to be split by Cor.~\ref{c.SPLIDI}. Hence it suffices to show that $K_K = K \otimes K$ contains zero divisors. If $K/F$ is separable, this follows from Exc.~\ref{pr.MISPLET}~(a). On the other hand, if $K/F$ is inseparable, then $\ch(F) = 2$ and $\xi^2 \in F$ for all $\xi \in K \setminus F$. This implies $0 \neq u := 1_K \otimes \xi + \xi \otimes 1_K \in K \otimes K$ and $u^2 = 1_K \otimes \xi^2 + \xi^2 \otimes 1_K = 2\xi^21_{K \otimes K} = 0$. Thus $K_K$ contains non-zero nilpotent elements, and we are done.

In order to prove necessity, let us assume that $K$ is a splitting field of $C$ and (i), (ii) do not hold. If $C$ is split, then it has dimension $4$ or $8$, and $C_1 := \Cay(K,1)$ in the first case, $C_2 := \Cay(K;1,1)$ in the second, is a split composition algebra over $F$ of the same dimension as $C$. Thus $C \cong C_i$ for some $i = 1,2$, and since $C_i$ contains an isomorphic copy of $K$, so does $C$. We are left with the case that $C$ is a division algebra, so by Cor.~\ref{c.SPLIDI} its norm is anisotropic. Pick any $\xi \in K \setminus F$. Then $\xi^2 = \alpha\xi - \beta$ with $\alpha := t_K(\xi)$, $\beta := n_K(\xi)$. Since $C_K = C \oplus (C \otimes \xi)$ is split over $K$, we find elements $u,v \in C$ not both zero such that
\begin{align*}
0 =\,\,&n_C(u - v \otimes \xi) = n_C(u) - \xi n_C(u,v) + \xi^2n_C(v) \\
=\,\,& \big(n_C(u) - \beta n_C(v)\big) + \big(\alpha n_C(v) - n_C(u,v)\big)\xi,
\end{align*}
hence $n_C(u) = \beta n_C(v)$ and $n_C(u,v) = \alpha n_C(v)$. This implies $v \in C^\times$ and with $x := uv^{-1}$ we deduce $\beta = n_C(x)$, $\alpha = n_C(u,n_C(v)^{-1}v) = n_C(u,\bar v^{-1}) = t_C(x)$. Summing up, $F[x] \subseteq C$ is a subalgebra isomorphic to $K$, and we have found an embedding $K \hookrightarrow C$.

(b) We claim that the given quadratic form $q$ is isotropic iff there exists $x \in C^0$ such that $n_C(x) = -\alpha$.  The ``if'' direction is obvious.  For ``only if'', suppose $c = \beta 1_C + y$ for $\beta \in k$ and $y \in C^0$ such that $c \ne 0$ but $q(c) = \alpha \beta^2 + n_C(y) = 0$.  If $\beta = 0$, then $n_C\vert_{C^0}$ is isotropic, so it contains a hyperbolic plane and such an $x$ exists.  If $\beta \ne 0$, then $n_C(y/\beta) = n_C(y)/\beta^2 = -\alpha$.  This verifies the claim.

For $x \in C^0$, we have $n_C(x) = -\alpha$ if and only if $x^2 = \alpha$.  By part (a), such an $x$ exists iff there is an $F$-embedding $K \hookrightarrow C$.
\end{sol}

\begin{sol}{pr.skip.antiSPLIFICO} \label{sol.skip.antiSPLIFICO} Write $C \cong F[\bft] / (f)$ for some irreducible polynomial $f \in F[\bft]$.  By hypothesis, $K[\bft]/(f) \cong C_K$ is not a field, so $f$ is reducible in $K[\bft]$.  Because $f$ has degree 2 or 3, one of the factors of $f$ has degree 1, i.e., $f$ has a root $\alpha \in K$.  Then $F[\alpha]$ is a subfield of $K$ isomorphic to $C$.
\end{sol}

\begin{sol}{pr.skip.detlev.biquad} \label{sol.skip.detlev.biquad} 
$C_{\IR} \cong \IH$ is division, so $C$ is division.  Since
\[
1 + 7\cdot i^2 + 15 \cdot i^2 +7 \cdot 15 \cdot (1/\sqrt{5})^2 = 0,
\]
$(n_C)_K$ is isotropic, so $C_K$ is split.

A subalgebra $E$ of $K$ properly containing $\IQ$ is a quadratic field extension of $\IQ$.  Now $K$ is Galois over $\IQ$ with Galois group $\IZ/2 \times \IZ/2$, so the three possibilities for $E$ are $\IQ(\sqrt{d})$ for $d = 5$, $-1$, or $-5$.  

By Exc.~\ref{pr.SPLIFICO}, it suffices to show that $q_d := \la d, 7, 15, 7 \cdot 15 \raq$ is anisotropic over $\IQ$.  For $d = 5$, we have $(q_5)_\IR \cong (n_C)_\IR$, which is anisotropic.  For $d = -1$ or $-5$, we verify that $(q_d)_{\IQ_p}$ is anisotropic, where $\IQ_p$ denotes the $p$-adic numbers with $p = 5$ or 3, respectively.

We have $q_{-1} \cong \la -1, 7 \raq \perp  5 \la 3, 21 \raq$.  We apply Exc.~\ref{pr.DVR.FORMS} or \cite[Lemma 19.5]{MR2427530}, to see that $q_{-1}$ is anisotropic over $\IQ_5$ if the residue forms $\la -1, 7 \raq$ and $\la 3, 21 \raq$ are anisotropic over $\IF_5$.  To check this, by Exc.~\ref{pr.HYPDISC} it suffices to note that the forms have determinant 3 and 2, respectively, which are not congruent to $-1$ mod squares.

We have $q_{-5} \cong \la -5, 7 \raq \perp 3 \la 5, 35 \raq$.  For $(q_{-5})_{\IQ_3}$, the residue forms are $\la -5, 7 \raq \cong \la 1, 1 \raq$ and $\la 2, 2 \raq$ over $\IF_3$, which both have determinant a square.  Since $-1$ is not a square in $\IF_3$, these forms are anisotropic, whence so is $q_{-5}$.
\end{sol}

\solnsec{Section~\ref{s.AFF}}

\begin{sol}{pr.quad.exercise} \label{sol.quad.exercise}
Let $R$ be a $k$-algebra and suppose $R \cong R_0 \times R_1$ for $k$-algebras $R_i$.  Then we obtain an element of $\Aut(D_R)$ that acts as the identity on $D_{R_0}$ and as conjugation on $D_{R_1}$.  Exercise \ref{pr.ETAUT} says that every element of $\Aut(D_R)$ arises in this way and conversely each expression $R \cong R_0 \times R_1$ gives an automorphism of $D_R$.  In this way, we may identify $\Aut(D_R)$ as a set with the collection of complete orthogonal system of idempotents $e_0, e_1$ in $R$, i.e., with the set of $k$-algebra homomorphisms $\Hom(k \times k, R)$.

Thus we have defined a bijection $\Aut(D_R)$ between and $(\IZ/2)(R)$ as sets, in a way that is functorial in $R$.  We check that this bijection is compatible with the group operation.  Let $(c_0, c_1)$, $(d_0, d_1)$ be complete orthogonal systems of idempotents in $R$.  Then $(c_0 d_0, c_0 d_1, c_1 d_0, c_1 d_1)$ is a complete orthogonal system of idempotents, giving rise to a direct product decomposition $R = \prod_{i, j} R_{ij}$ for $R_{ij} := c_i d_j R$ as well as $D_R = \prod_{i, j} D_{ij}$ for $D_{ij} := D_{R_{ij}}$.  The product of $(c_0, c_1)$, $(d_0, d_1)$ in $\Aut(D_R)$ is an automorphism of $D_R$ that is the identity on $D_{00} \times D_{11}$ and conjugation on $D_{01} \times D_{10}$.  That is, the product in $\Aut(D_R)$ corresponds to the complete orthogonal system of idempotents $(c_0 d_0 + c_1 d_1, c_0 d_1 + c_1 d_0)$.

We describe the product in $(\IZ/2)(R)$, which is a $k$-algebra homomorphism $k \times k \to R$, by evaluating it on an element $(x_0, x_1) \in k \times k$.  Put $\Delta$ for the co-multiplication homomorphism defined in \ref{skip.constant}, $\phi_i$ for the map $k \times k \to k$ that is projection on factor $i = 0, 1$, and $1_i$ for the identity element of the $i$-th factor of $k \times k$.  Then the product of $(c_0, c_1)$, $(d_0, d_1)$ computing in $(\IZ/2)(R)$ maps $(x_0, x_1)$ to 
\begin{align*}
\big( (c_0 \phi_0 &+ c_1 \phi_1) \otimes (d_0 \phi_0 + d_1 \phi_1) \big) \, \Delta(x_0, x_1) \\
&= (\sum_{i,j} c_i d_j \phi_i \otimes \phi_j) \big( x_0 (1_0 \otimes 1_0 + 1_1 \otimes 1_1) + x_1(1_0 \otimes 1_1 + 1_1 \otimes 1_0) \big) \\
&= x_0 (c_0 d_0 + c_1 d_1) + x_1 (c_0 d_1 + c_1 d_0),
\end{align*}
confirming that the two products agree.
\end{sol}

\begin{sol}{pr.skip.GLstab} \label{sol.skip.GLstab}
For each $i$ and $\chi \in M_i^*$, there is a polynomial function $f_{i,\chi}$ on $\bfGL(M)$ defined by
\[
f_{i,\chi}(g) = \chi(\rho_i(g)m_i - m_i).
\]
Since $M_i$ is finitely generated projective, $f_{i,\chi}(g) = 0$ for all $\chi$ if and only if $\rho_i(g)m_i - m_i = 0$.  That is,
for $\mcF := \cup_{i \in I} \cup_{\chi_i \in M_i^*} f_{i,\chi}$, we have
\[
\bfH(R) = \{ g \in \bfGL(M)(R) \mid \text{$f(g) = 0$ for all $f \in \mcF$} \}.
\]
\end{sol}

\solnsec{Section~\ref{s.ETFPPF}}

\begin{sol}{pr.TRAFIPR} \label{sol.TRAFIPR}
 (a) Let $\pR$ be a finitely presented $R$-algebra. Then there are presentations
\begin{eqnarray}
\label{IKS} \xymatrix{
0 \ar[r] & I \ar[r] & k[\bfS] \ar[r]_{\pi} & R \ar[r] & 0}, \\
\label{JERT} \xymatrix{
0 \ar[r] & J \ar[r] & R[\bfT] \ar[r]_{\rho} & \pR \ar[r] & 0}
\end{eqnarray}
of $R$ over $k$, $\pR$ over $R$, respectively, with chains $\bfS = (\bfs_1,\dots,\bfs_m)$, $\bfT= (\bft_1,\dots,\bft_n)$ of independent variables and finitely generated ideals $I \subseteq k[\bfS]$, $J \subseteq k[\bfT]$. Extending \eqref{IKS} from $k$ to $k[\bfT]$, we conclude that the sequence
\begin{align}
\label{IKIK} \xymatrix{
0 \ar[r] & I \otimes k[\bfT] \ar[r] & k[\bfS,\bfT] = k[\bfS] \otimes k[\bfT] \ar[r]_(.55){\pi_{k[\bfT]}} & R \otimes k[\bfT] = R[\bfT] \ar[r] & 0}
\end{align}
is exact, whence
\begin{align}
\label{ROPIK} \vph := \rho \circ \pi_{k[\bfT]}\:k[\bfS,\bfT] \longrightarrow \pR 
\end{align}
is a surjective morphism in $\kalg$ which by \eqref{JERT} satisfies
\begin{align}
\label{PINV} \Ker(\vph) = \pi_{k[\bfT]}^{-1}\big(\Ker(\rho)\big) = \pi_{k[\bfT]}^{-1}(J).
\end{align}
Now pick $g_1,\dots,g_r \in R[\bfT]$ such that $J = \sum_{j=1}^r R[\bfT]g_j$ and use \eqref{IKIK} to find a lift $h_j$ of $g_j$ in $k[\bfS,\bfT]$ under $\pi_{k[\bfT]}$: $\pi_{k[\bfT]}(h_j) = g_j$ for $1 \leq j \leq r$. By \eqref{IKIK}--\eqref{PINV}, this is easily seen to imply
\[
\Ker(\vph) = I \otimes k[\bfT] + \sum_{j=1}^r k[\bfS,\bfT]h_j.
\]
But $I \otimes k[\bfT]$ is a finite $k[\bfT]$-module, forcing $\Ker(\vph) \subseteq k[\bfS,\bfT]$ to be a finitely generated ideal. Hence $\pR$ is finitely presented as a $k$-algebra.

(b) By (a), we may assume $R = k$. Since the $k$-algebra $E := k[\bft]/(f\bft - 1)$ is clearly finitely presented, it suffices to show that $k_f$ and $E$ are canonically isomorphic. Note first that the natural projection $\pi\:k[\bft] \to E$ makes $f$ invertible: $\pi(f) \in E^\times$ with inverse $\pi(\bft)$. Hence the unit homomorphism $k \to E$ extends to a homomorphism $\vph\:k_f \to E$ such that $\vph(1/f) = \pi(\bft)$. On the other hand, the $k$-homomorphism $k[\bft] \to k_f$ sending $\bft$ to $1/f
$ kills $f\bft - 1$ and hence induces a $k$-homomorphism $\psi\:E \to k_f$. One checks that $\vph$ and $\psi$ are inverse to one another, and the assertion follows.
\end{sol}

\begin{sol}{pr.IDFIPR} \label{sol.IDFIPR}
 (a) Let $\pi\:k[\bfT] \to R$ be a surjective morphism in $\kalg$, giving rise to a presentation 
\[
\xymatrix{
0 \ar[r] & I \ar[r] & k[\bfT] \ar[r]_{\vph \circ \pi} & \pR \ar[r] & 0}
\] 
of $\pR$ as a $k$-algebra. By Prop.~\ref{p.APREFI}, therefore, $I \subseteq k[\bfT]$ is a finitely generated ideal. But $I = \pi^{-1}(\Ker(\vph))$, and we conclude that $\Ker(\vph) = \pi(I) \subseteq R$ is a finitely generated ideal.

(b) Let
\begin{equation}
\label{IKIP} \xymatrix{
0 \ar[r] & I \ar[r] & k[\bfT] \ar[r]_(.55){\pi} & R \ar[r] & 0}
\end{equation}
be a finite presentation of $R$. Setting $J := \Ker(\vph)$, this induces a presentation
\begin{equation}
\label{PINJ} \xymatrix{
0 \ar[r] & \pi^{-1}(J) \ar[r] & k[\bfT] \ar[r]_{\vph \circ \pi} & \pR \ar[r] & 0}
\end{equation}
of $\pR$. Apply \eqref{IKIP} and let $\pi(f_i)$, $f_i \in k[\bfT]$, $1 \leq i \leq m$, be finitely many generators of $J$ as an ideal in $R$. Then one checks that
\[
\pi^ {-1}(J) = I + \sum_{i=1}^ m k[\bfT]f_i,
\]
so \eqref{PINJ} makes $\pR$ a finitely presented $k$-algebra.
\end{sol}

\begin{sol}{pr.EMMA} \label{sol.EMMA}
 By \ref{e.EXMA}, $k[M_{\bfa}] \cong S(M^\ast)$ is the symmetric algebra of $M^\ast$, the dual of $M$. Since $M^\ast$ is a finitely generated projective $k$-module as well, it suffices to show: if $M$ is a finitely generated projective $k$-module, then $S(M)$ is a finitely presented $k$-algebra.

Applying Exc.~\ref{pr.RANDEC}, we find a complete orthogonal system $(\vep_i)_{i\in\IN}$ of idempotents in $k$ such that $M_i = M \otimes k_i$ for each $i \in \IN$ is a finitely generated projective module of rank $i$ over $k_i := k\vep_i$. Note that $\vep_i = 0$, hence $k_i = \{0\}$ and $M_i = \{0\}$, for almost all $i \in \IN$. Since passing to the symmetric algebra of a module is compatible with base change \cite[III.6, Prop.~7]{MR0354207}, we conclude
\[
S(M) = \prod_{i\in\IN}S(M_i)
\]
as $k$-algebras. We are thus reduced to the case that $M$ is finitely generated projective of rank $r \in \IN$ over $k$. But then \ref{skip.other.fflat}~\ref{skip.other.fflat.2} produces a faithfully flat $k$-algebra $\pk$ making $M_{\pk}$ a free $\pk$-module of rank $r$. Since $S(M)_{\pk} = S(M_{\pk})$, therefore being a polynomial ring in $r$ variables \cite[III.6]{MR0354207}, hence finitely presented over $\pk$, so is $S(M)$ over $k$, by Cor.~\ref{c.FAFLAD}.
\end{sol}

\begin{sol}{pr.SGLT} \label{sol.SGLT}
 It is straightforward to verify that $\tiw$ is a subfunctor of $M_{\bfa}$. In order to derive the second part of the exercise, let us assume that $M$ is finitely generated projective. For $u^\ast \in M^\ast$, the set maps
\[
f_{u^\ast}(R)\:M_R \longrightarrow R, \quad x \longmapsto \la u_R^\ast,x - w_R\ra, 
\]
vary functorially with $R \in \kalg$ and hence define an element $f_{u^\ast} \in k[M_\bfa]$. Since $M^\ast$ is finitely generated projective as well, and the canonical pairing $M^\ast \times M \to k$ is regular, we conclude that $x \in M_R$ agrees with $w_R$ if and only if $f_{u^\ast}(R)(x) = 0$ for all $u^\ast \in M^\ast$. By Exercises \ref{pr.IDFIPR}, \ref{pr.EMMA}~(b), therefore, $\tiw \subseteq M_\bfa$ is a finitely presented closed affine subscheme.
\end{sol}

\begin{sol}{pr.FFSEP} \label{sol.FFSEP}
 By definition of non-singularity (\ref{ss.REQUA}), we have to prove that the quadratic form $q_K:M_K \to K$ over $K$ is non-degenerate, for any field $K \in \kalg$. By Exc.~\ref{pr.COVERS} there is a field $L \in \Ralg$ containing $K$ such that the maps $k \to K \to L$ and $k \to R \to L$ agree.  Since $q_R$ is non-singular, the quadratic form $(q_K)_L \cong q_L \cong (q_R)_L$ over $L$ is non-degenerate. Hence so is $q_K$ over $K$.
\end{sol}

\begin{sol}{pr.FFPOLA} \label{sol.FFPOLA}
 (a) For $i = 0,1$, $T \in \Salg$ and all $k$-modules $P$, we have the natural identifications
\[
(P_S)_T = P_T = P_{T_i} = (P_R)_{T_i}
\]
via (\ref{ss.ITSCA}.\ref{ITBA}), so $g_{iT}$ is a set map from $(M_S)_T$ to $(N_S)_T$, as claimed. It remains to show that these set maps vary functorially with $T$. In order to see this, we regard the identity map of $T$ as an isomorphism $T \to T_i$, $t \mapsto t_i$, of $k$-algebras. Given a morphism $\vph\:T \to \pT$ in $\Salg$, we therefore obtain a morphism $\vph_i\:T_i \to \pT_i$ in $\Ralg$ by defining $\vph_i(t_i) := \vph(t)_i$ for all $t \in T$. One checks that $\Eins_{M_R} \otimes_R \vph_i = \Eins_{M_S} \otimes_S \vph$, and it follows that the diagram
\[
\xymatrix{
(M_S)_T \ar[r]_{g_{iT}} \ar[d]_{\Eins_{M_S} \otimes_S \vph} & (N_S)_T \ar[d]^{\Eins_{N_S} \otimes_S \vph} \\
(M_S)_{\pT} \ar[r]_{g_{i\pT}} & (N_S)_{\pT}}
\] 
commutes. Hence $g_i\:M_S \to N_S$ is a polynomial law over $S$. 

\smallskip

(b) If $f\:M \to N$ is a polynomial law over $k$ such that $f \otimes R = g$, then for all $T \in \Salg$ we deduce that $g_{iT} = g_{T_i} = (f \otimes R)_{T_i} = f_{T_i} = f_T$ is independent of $i = 0,1$. Thus $g_0 = g_2$. Conversely, assume $g_0 = g_1 =: h$. We first prove uniqueness of $f$ and let $\pk \in \kalg$. One checks that the outer squares in \eqref{CHARF} commute. Hence so does the inner one since $g_{R_{\pk}} = f_{R_{\pk}}$ by hypothesis. Moreover, the vertical arrows in \eqref{CHARF} are injective by Prop.~\ref{p.CHAFF}, proving uniqueness of the set maps $f_{\pk}$, hence of $f$ as a polynomial law over $k$.

In order to prove existence, let $\pk \in \kalg$. To simplify notation, we indicate the base change from $k$ to $\pk$ simply by a dash, for example $\pM = M_{\pk}$, $\pR = R \otimes \pk$, $\pS = S \otimes \pk = \pR \otimes_{\pk} \pR$. Furthermore, the notational conventions fixed in (a) on the level of $R$ and $S$ will now be employed for $\pR$ and $\pS$. In particular, the very definition of the $\pR$-algebras $\pS_i$ ($i = 0,1$) yield morphisms 
\[
\sigma_i^{\pR}\:\pR \longrightarrow \pS_i, \quad \pr \longmapsto \sigma_i^{\pR}(\pr) := d^i_{\pR}(\pr)_i
\]
in $\pRalg \subseteq \Ralg$. Now we consider the following diagram. 
\begin{equation}
\vcenter{\label{DESSET} \xymatrix{
0 \ar[d] &&& 0 \ar[d] \\
\pM \ar@{-->}[rrr]^{\exists! f_{\pk}} \ar[d]_{\can_{\pM,\pR}} &&& \pN \ar[d]^{\can_{\pN,\pR}} \\
\pM_{\pR} \ar[r]_{\Eins} \ar[d]_{\rho_i^{\pR,\pM}} & (M_R)_{\pR} \ar[r]_{g_{\pR}} \ar[d]^{\Eins_{M_R} \otimes_R \sigma_i^{\pR}} & (N_R)_{\pR} \ar[r]_{\Eins} \ar[d]^{\Eins_{N_R} \otimes_R \sigma_i^{\pR}} & \pN_{\pR} \ar[d]^{\rho_i^{\pR,\pN}} \\
\pM_{\pS} \ar[r]_{\Eins} & (M_R)_{\pS_i} \ar[r]_{h_{\pS} = g_{\pS_i}} & (N_R)_{\pS_i} \ar[r]_{\Eins} & \pN_{\pS},}}
\end{equation}
whose inner lower square commutes by the definition of polynomial laws, while the outer lower ones do thanks, e.g., to the following computation, for all $x \in M$, $r \in R$, $\palpha,\pbeta \in \pk$:
\begin{align*}
\rho_1^{\pR,\pM}\big((x \otimes \palpha) \otimes_{\pk} (r \otimes \pbeta)\big) =\,\,&(x \otimes \palpha) \otimes_{\pk} \big((r \otimes \pbeta) \otimes_{\pk} (1_R \otimes 1_{\pk})\big) \\
=\,\,&(x \otimes \palpha) \otimes_{\pk} \big((r \otimes 1_R) \otimes \pbeta\big) \\ 
=\,\,&x \otimes \big((r \otimes 1_R) \otimes \palpha\pbeta\big) \\
=\,\,&x \otimes \big((r \otimes 1_R) \otimes \palpha\pbeta\big)_1 \\ 
=\,\,&(x \otimes 1_R) \otimes_R \big((r \otimes 1_R) \otimes \palpha\pbeta\big)_1 \\
=\,\,&(x \otimes 1_R) \otimes_R \big((r \otimes \palpha\pbeta) \otimes_{\pk} (1_R \otimes 1_{\pk})\big)_1 \\
=\,\,&(x \otimes 1_R) \otimes_R \sigma_1^{\pR}(r \otimes \palpha\pbeta) \\
=\,\,&(\Eins_{M_R} \otimes_R \sigma_1^{\pR})\big((x \otimes 1_R) \otimes_R (r \otimes \palpha\pbeta)\big) \\
=\,\,&(\Eins_{M_R} \otimes_R \sigma_1^{\pR})\Big(x \otimes\big(\palpha(r \otimes \pbeta)\big)\Big) \\
=\,\,&(\Eins_{M_R} \otimes_R \sigma_1^{\pR})\big((x \otimes \palpha) \otimes_{\pk} (r \otimes \pbeta)\big).
\end{align*}
The universal property (\ref{ss.EQU}) of the equalizer $\can_{\pN,\pR}$ of $\rho_1^{\pR,\pN}$, $\rho_0^{\pR,\pN}$\break (Prop.~\ref{p.RHOIQ}), applied to the map $ u := g_{\pR} \circ \can_{\pM,\pR}\:\pM \to \pN_{\pR}$, by \eqref{DESSET} implies that there exists a unique set map $f_{\pk}\:\pM \to \pN$ making \eqref{DESSET} totally commutative. In other words, we have a commutative diagram
\begin{equation}
\vcenter{\label{EXDES} \xymatrix{
\pM \ar[r]_{f_{\pk}} \ar[d]_{\can_{\pM,\pR}} & \pN \ar[d]^{\can_{\pN,\pR}} \\
\pM_{\pR} = (M_R)_{{\pR}} \ar[r]_{g_{\pR}} & (N_R)_{\pR} = \pN_{\pR}.}}
\end{equation}
Next we show that $f_{\pk}$ depends functorially on $\pk$, so let $\vph\:\pk \to \ppk$ be a morphism in $\kalg$. Indicating the base change from $k$ to $\ppk$ by a double dash, we see that
\begin{align}
\label{PHEER} \vph_R := \Eins_R \otimes \vph\:\pR \longrightarrow \ppR
\end{align}
is a morphism in $\Ralg$, and one checks that the diagram
\begin{equation}
\vcenter{\label{EMPHER} \xymatrix{
\pM \ar[r]_{\Eins_M \otimes \vph} \ar[d]_{\can_{\pM,\pR}} & \ppM \ar[d]^{\can_{\ppM,\ppR}} \\
\pM_{\pR} = (M_R)_{\pR} \ar[r]_{\Eins_{M_R} \otimes_R \vph_R} & (M_R)_{\ppR} = \ppM_{\ppR}}}
\end{equation}
commutes. Now consider the cube
\[
\xymatrix{
\pM \ar[rr]^{f_{\pk}} \ar[dd]_{\can_{\pM,\pR}} \ar[rd]^{\Eins_M \otimes \vph} && \pN \ar@{.>}[dd]^(.4){\can_{\pN,\pR}} \ar[rd]^{\Eins_N \otimes \vph} \\
& \ppM \ar[rr]_(.4){f_{\ppk}} \ar[dd]_(.4){\can_{\ppM,\ppR}} && \ppN \ar[dd]^{\can_{\ppN,\ppR}} \\
\pM_{\pR} = (M_R)_{\pR} \ar@{.>}[rr]_(.6){g_{\pR}} \ar[rd]_{\Eins_{M_R} \otimes_R \vph_R} &&(N_R)_{\pR} = \pN_{\pR} \ar@{.>}[rd]_{\Eins_{N_R} \otimes_R \vph_R} \\
& \ppM_{\ppR} = (M_R)_{\ppR} \ar[rr]_{g_{\ppR}} && (N_R)_{\ppR} = \ppN_{\ppR},}
\]
where by \eqref{EXDES}--\eqref{EMPHER} all rectangles commute, with the possible exception of the top one. By diagram chasing, therefore, so does the top one after being composed with $\can_{\ppN,\ppR}$. But $\ppR$ is faithfully flat over $\ppk$, forcing $\can_{\ppN,\ppR}$ by Prop.~\ref{p.CHAFF} to be injective, and we have proved that the top rectangle commutes as well. Summing up, $f\:M \to N$ is a polynomial law over $k$. Hence the solution of the problem will be complete once we have shown $f \otimes R = g$. In order to see this, let $\pk \in \Ralg$ and $\rho\:R \to \pk$ the corresponding unit morphism in $\kalg$. With the notational simplifications introduced before, we obtain an induced morphism $\prho\:\pR\to \pk$ given by
\[
\prho(r \otimes \palpha) = \rho(r)\palpha = r\palpha 
\] 
for $r \in R$, $\palpha \in \pk$. In particular, $\prho$ is a morphism in $\Ralg$, and we obtain a diagram
\[
\xymatrix{
\pM \ar[r]_{f_{\pk}} \ar[d]^{\can_{\pM,\pR}} \ar@/_3pc/[dd]_{\Eins} & \pN \ar[d]_{\can_{\pN,\pR}} \ar@/^3pc/[dd]^{\Eins} \\
\pM_{\pR} = (M_R)_{\pR} \ar[r]_{g_{\pR}} \ar[d]^{\Eins_{M_R} \otimes_R\, \prho} & (N_R)_{\pR} = \pN_{\pR} \ar[d]_{\Eins_{N_R} \otimes_R\, \prho} \\
\pM = (M_R)_{\pk} \ar[r]_{g_{\pk}} & (N_R)_{\pk} = \pN,}
\]
where the outer triangles are easily seen to be commutative. Hence $f_{\pk} = g_{\pk}$, as desired.
\end{sol}

\begin{sol}{pr.SYHOPO} \label{sol.SYHOPO} In the special case where $M$ and $N$ are free, this holds by Cor.~\ref{c.FREPOLA}.  

If $M$ and $N$ have constant rank, then there is a faithfully flat $R \in \kalg$ such that $M_R$ and $N_R$ are free $R$-modules.  Then we have a commutative diagram
\[
\xymatrix{S^d(M^*) \otimes N \otimes R \ar[r] \ar[d] & \Pol^d(M,N) \otimes R \ar[d]\\
S^d(M_R^*) \otimes N_R \ar[r] & \Pol^d(M_R, N_R)}
\]
where the left arrow is trivially an isomorphism, the bottom arrow is an isomorphism by the free case, and the right arrow is an isomorphism by \ref{skip.other.fflat}~\ref{skip.other.fflat.3}, so the top arrow is also an isomorphism.  Since $\Phi$ becomes an isomorphism after tensoring with the faithfully flat $R$, it is itself an isomorphism.

If $M$ and $N$ do not have constant rank, then there is an $n$ such that $k = \prod_{i=1}^n k_i$, $M = \prod_{i=1}^n M_i$, and $N = \prod_{i=1}^n N_i$ such that each $M_i$ and $N_i$ is a finitely generated projective $k_i$-module of constant rank.  Then $S^d(M^*) \otimes N$ is naturally isomorphic to $\prod_{i=1}^n (S^d(M_i^*) \otimes N_i)$ and similarly for $\Pol^d(M, N)$, and the conclusion follows from the previous case.
    
\end{sol}

\begin{sol}{pr.FFCON} \label{sol.FFCON}
 We systematically adhere to the notation and terminology of Exc.~\ref{pr.FFPOLA}.

(a) Let $e$ be the unit element of $C_R$ and write $g := \tie\:C_R \to C_R$ for the constant polynomial law over $R$ determined by $e$, as defined in 
Exc.~\ref{pr.COPOLA}.  For $T \in\Salg$ and $i = 1,2$, $e_{T_i} \in (C_R)_{T_i} = (C_S)_T$ is the unit element of $(C_S)_T$, and we conclude $e_{T_1} = e_{T_2}$. Thus $g_1 = g_2$ as polynomial laws $C_S \to C_S$ over $S$. Applying Exc.~\ref{pr.FFPOLA}, we find a unique polynomial law $f\:C \to C $ over $k$ such that $f \otimes R = g$. Put $e_0 := f_k(0) \in C$. Since the vertical arrows in \eqref{CHARF} are injective, we conclude not only $e_{0R} = e$ but also that $C$ is unital with $1_C = e_0$.

(b) Let $n\:C_R \to R$ be a quadratic form over $R$ making $C_R$ a conic $R$-algebra. For $T \in \Salg$ and $i = 1,2$, there are two quadratic forms over $T$ making $(C_S)_T = (C_R)_{T_i}$ a conic $T$-algebra, namely, $n_{T_i}$. Hence Prop.~\ref{p.UNINOR} implies $n_{T_1} = n_{T_2}$. Viewing $n$ as a homogeneous polynomial law $g\:C_R \to R$ of degree $2$ over $R$, we have $g_1 = n_{T_1} = n_{T_2} = g_2$. Thus Exc.~\ref{pr.FFPOLA} yields a unique polynomial law $f\:C \to k$ over $k$ such that $f \otimes R = g$. Hence the set map $n_0 := f_k\:C \to k$ is a quadratic form over $k$, and from \eqref{CHARF} we deduce as before that $C$ is a conic $k$-algebra with norm $n_0$ whose base change to $R$ is $C_R$ as a conic $R$-algebra with norm $n$.
\end{sol}

\begin{sol}{skip.ex.charpoly} \label{sol.skip.ex.charpoly}(a): If $M$ is a finitely generated free module, say of rank $d$, then the characteristic polynomial as in \ref{ss.charpoly} does not depend on the choice of isomorphism $M \xrightarrow{\sim} k^d$, equivalently $\End_k(M) \xrightarrow{\sim} \Mat_d(k)$, and so defines a polynomial law on $\End_k(M)$ in that case determined by $d$.  If $M$ is merely assumed to have constant rank $d$, then there is a faithfully flat $R \in \kalg$ so that $M \otimes R$ is free of rank $d$, and we apply faithfully flat descent as in Exercises \ref{pr.FFPOLA} and \ref{pr.FFCON} to obtain a characteristic polynomial over $k$.  Finally, if $M$ does not have constant rank, we may write $k = \prod_{i = 1}^d k_i$ where $M_i := M \otimes k_i$ is projective of constant rank for each $i$, and we define the characteristic polynomial of $x \in \End_k(M)$ in $k[\bft]$ to be $(f_1, \ldots, f_d)$, where $f_i \in k_i[\bft]$ is the characteristic polynomial of $x \otimes k_i$.

(b): For (i), note that $(x, y) \mapsto \det(xy) - (\det x) (\det y)$ is a polynomial law $\End_k(M) \times \End_k(M) \to k$ that is zero if $M$ is free.  Note that this implies that $x$ invertible implies $\det x$ invertible in $k$. 

For (ii), evaluating the characteristic polynomial of $x$ at $x$ defines a polynomial law $\End_k(M) \to \End_k(M)$ that is zero whenever $M$ is free.

For (iii), we reduce to the case where $M$ has constant rank $d$.  The characteristic polynomial of $x$ is $\bft^d - \alpha_1 \bft^{d-1} + \cdots + (-1)^d \alpha_d$, with $\alpha_d := \det x$.  Define $x^\sharp := (-1)^{d+1} (x^{d-1} - \alpha_1 x^{d-2} + \cdots + (-1)^{d-1} \alpha_{d-1})$.  (The element $x^\sharp$ is known as the \emph{classical adjoint} of $x$.)  Then $x x^\sharp = (-1)^{d+1}(x^d - \alpha_1 x^{d-1} + \cdots (-1)^{d-1} \alpha_{d-1} x)$, which is almost the characteristic polynomial of $x$ evaluated at $x$.  Applying (ii), we find:
\[
xx^\sharp = (-1)^{d+1} (0 - (-1)^d \alpha_d) = \alpha_d 1_{\End_k(M)}.
\]
Therefore, if $\alpha_d$ is invertible, we have $x^{-1} = x^\sharp \alpha_d^{-1}$, proving (iii). 
    \end{sol}

\begin{sol}{pr.SMODI} \label{sol.SMODI}  First of all, assume $\bfX_1$, $\bfX_2$ are both smooth over $k$. Then $k[\bfX_1]$ is finitely presented over $k$. But the property of a $k$-algebra in $\kalg$ to be finitely presented is stable under base change (\ref{ss.FIPRAL}), forcing $k[\bfX_1 \times \bfX_2] \cong k[\bfX_1] \otimes k[\bfX_2]$ (\ref{ss.PRAFF}) to be finitely presented over $k[\bf X_2]$. On the other hand, $k[\bfX_2]$ is finitely presented over $k$ as well, and we deduce from Exc.~\ref{pr.TRAFIPR}~(a) that $k[\bfX_1 \times \bf X_2]$ is finitely presented over $k$. Now  let $R \in \kalg$ and suppose $I \subseteq R$ is an ideal that squares to zero. Since the set maps $\bfX_i(R) \to \bfX_i(R/I)$ are surjective for $i = 1,2$ by definition, so is their product
\begin{align}
\label{EXIR} (\bfX_1 \times \bfX_2)(R) = \bfX_1(R) \times \bfX_2(R) \longrightarrow \bfX_1(R/I) \times \bfX_2(R/I) = (\bfX_1 \times \bfX_2)(R/I).	
\end{align}
Thus $\bfX_1 \times \bfX_2$ is smooth. Conversely, let this be so and assume that $\bfX_1$ is finitely presented over $k$ as well as $\bfX_2(k) \neq \emptyset$. Pick $x_2 \in \bfX_2(k)$ and write $\px_2$ for its canonical image in $\bfX_2(R/I)$. Given any $\px_1 \in \bfX_1(R/I)$, the surjectivity of \eqref{EXIR} yields a pre-image $(x_1,x_2) \in (\bfX_1 \times \bfX_2)(R)$ of $(\px_1,\px_2) \in (\bfX_1 \times \bfX_2)(R/I)$, forcing $x_1$ to be a pre-image of $\px_1$ in $\bfX_1(R)$. Thus the natural map $\bfX_1(R) \to \bfX_1(R/I)$ is surjective, and $\bfX_1$ is smooth.
\end{sol}

\begin{sol}{pr.roots.etale}
$R$ is evidently finitely presented. Since it is free with basis $1_k, \bft, \ldots, \bft^{d-1}$, it is flat.

In case $k$ is a field, $R$ is \'etale if and only if the polynomial $\bft^d - x$ has $d$ distinct roots in an algebraic closure $\aclosk$ of $k$.  The roots of this polynomial are in any case $d$-th roots of $x$.   The derivative of this polynomial is $dx^{d-1}$.  If $d$ and $x$ are both nonzero in $k$, then the derivative only has zero as a root, which is not a root of the polynomial, so it is separable.  Conversely, if $d$ is divisible by the characteristic $p$ of $k$, then $\bft^d - x = (\bft^{d/p} - y)^p$ for a $p$-th root $y$ of $x$ in $\aclosk$.  And if $x = 0$ in $k$, then $\bft^d - x = \bft^d$.  In either case, $R$ is not a direct product of finite separable extensions of $k$.

For general $k$, $d$ and $x$ are invertible in $k$ if and only if they are not contained in any maximal ideal $\mfm$ of $k$, if and only if they have non-zero image in $k(\mfm)$.  This holds if and only if --- by the case for fields from the preceding paragraph --- $R \otimes k(\mfm)$ is an \'etale $k(\mfm)$-algebra.
\end{sol}

\begin{sol}{pr.ROOUNT} \label{sol.ROOUNT}  
(a) By \ref{ss.CLOSU}, $\bfmu_n$ is a closed subscheme of $\IA_k^1$, with co-ordinate algebra $k[\bft]/k[\bft](\bft^n - 1)$. But in actual fact, $\bfmu_n$ is also a closed subgroup functor, and hence a closed $k$-group subscheme, of $\bfGm$, with co-ordinate algebra $k[\bft,\bft^{-1}]/k[\bft,\bft^{-1}](\bft^n - 1)$. This proves the isomorphism in the second displayed equation (which can also be derived quite easily directly) and completes the proof of part (a). 

\smallskip

(b) There are $i,j \in \IZ$ such that $il + jm = 1$. For $R \in \kalg$, the assignments
\[
x \longmapsto (x^{jm},x^{il}), \quad (\text{resp.} (x,y) \longmapsto xy)
\]	
define group homomorphisms $\bfmu_n(R) \to \bfmu_l(R) \times \bfmu_m(R)$ (resp. $\bfmu_l(R) \times \bfmu_m(R) \to \bfmu_n(R)$) that are inverse to one another and compatible with base change. Hence they yield an isomorphism $\bfmu_n \cong \bfmu_l \times \bfmu_m$. 

\smallskip

(c) (i) and (iii) are equivalent by part (a) and Exc.~\ref{pr.roots.etale}.  (i) implies (ii) by \ref{ss.PROSMO}(d).  (ii) implies (i) as recalled in \ref{skip.et.egs}\ref{et.egs.smooth}.
\end{sol}

\begin{sol}{skip.pr.etale.gen} \label{sol.skip.pr.etale.gen} (a): The function
\[
e \mapsto 1_E \wedge e \wedge e^2 \wedge \cdots \wedge e^{d-1} \in \wedge^d E \cong k
\]
is a polynomial function $E \to k$.  It is non-zero for some $e$ if and only if $k[e] = E$ \cite[III.7, Thm.~2]{MR0354207}, proving that the set is Zariski-open.  The fact that it is not empty when $k$ is infinite is standard, see for example \cite[V.7, Prop.~7]{MR643362}.

(b): Write an element $e$ as $(e_1, \ldots, e_d)$ for $e_i \in k$.  Then $k[e] = E$ if and only if $1, e, \ldots, e^{d-1}$ are linearly independent, i.e., if and only if the matrix
\[
\left( \begin{matrix}
1 & e_1 & e_1^2 & \cdots & e_1^{d-1} \\
1 & e_2 & e_2^2 & \cdots & e_2^{d-1} \\
\vdots & \vdots & \vdots & \ddots & \vdots\\
1 & e_d & e_d^2 & \cdots & e_d^{d-1} 
\end{matrix} \right)
\]
is invertible.  The determinant of this matrix is the Vandermonde determinant, $\prod_{i < j} (e_i - e_j)$, and therefore $k[e] = E$ if and only if the $e_i$ are all distinct.  Such an $e$ exists if and only if $|k| \ge d$.
    
\end{sol}

\begin{sol}{skip.pr.simple} \label{sol.skip.pr.simple} For elements $a, a' \in A \otimes R$, we write $a$ as a finite sum $a = \sum_i a_i \otimes r_i$ with $a_i \in A$, $r_i \in R$ and similarly for $a'$.  Then
\[
aa' = \sum_i \sum_j a_i a_j \otimes r_i r_j,
\]
because $R \subseteq \Cent(A \otimes R)$.  In particular, since there exist $a, a' \in A \otimes R$ with $aa' \ne 0$, there are $a_i, a'_j \in A$ for some $i$, $j$ such that $a_i a'_j \ne 0$ in $A$.

Next suppose that $I$ is an ideal in $A$.  The displayed equation in the preceding paragraph shows that $I \otimes R$ is an ideal in $A \otimes R$, and therefore $I \otimes R = \{ 0 \}$ or $A \otimes R$.  Faithful flatness of $R$ implies that $I = \{ 0 \}$ or $A$, as claimed.
    \end{sol}

\begin{sol}{skip.pr.homdesc} \label{sol.skip.pr.homdesc} See \cite[p.~33, Prop.~2.5]{MR0417149}.
\end{sol}

\begin{sol}{skip.ex.trivtorsor} \label{sol.skip.ex.trivtorsor} For each $R \in \kalg$, define a map $g \mapsto \gamma_g$ from $\bfG(R)$ to the group of automorphisms of $\bfG_R$ viewed as a $\bfG_R$ torsor via 
$\gamma_g(h) = gh$ for $h \in \bfG(R)$.  (Note that $\gamma_g$ is indeed an automorphism of $\bfG_R$ because it acts on the left, whereas $\bfG_R$ acts on itself on the right.)  The map is a homomorphism because 
\[
\gamma_{gg'}(h) = gg'h = \gamma_g \gamma_{g'}(h)
\]
for all $g, g' \in \bfG(R)$.  It is injective because $\gamma_g = \gamma_{g'}$ implies that $g = \gamma_g(1_G) = \gamma_{g'}(1_G) = g'$.  

To see that it is surjective, let $\eta$ be an automorphism of $\bfG_R$ as a $\bfG_R$-torsor.  Then for $h \in \bfG(R)$, we have
\[
\eta(h) = \eta(1_G h) = \eta(1_G) h = \gamma_{\eta(1_G)}(h),
\]
i.e., $\eta = \gamma_{\eta(1_G)}$.
\end{sol}

\begin{sol}{etale.aut} \label{sol.etale.aut} \ref{etale.aut.1}:
Certainly, permuting the coordinates of elements of $E_\Gamma$ according to a permutation of $\Gamma$ does provide a $k$-algebra automorphism of $\Gamma$, so the group of permutations is a subgroup of $\mathrm{Aut}_{\kalg}(E_\Gamma)$. 

Conversely, since $k$ is connected, every idempotent in $E_\Gamma$ has a 0 or 1 in each coordinate.  It follows that the collection of elements $\{ 1_\gamma \mid \gamma \in \Gamma \}$ from Example \ref{skip.constant} is the unique maximal complete system of orthogonal idempotents.  That is, every $k$-algebra automorphism of $E_\Gamma$ permutes this set.  Given a $k$-algebra automorphism $f$ of $E_\Gamma$, we may modify it by a permutation of $\Gamma$ as in the preceding paragraph and so assume that $f$ fixes each $1_\gamma$.   As $f$ is $k$-linear, it follows that $f$ is the identity on each copy of $k$ in the product for $E_\Gamma$, verifying the claim.

\ref{etale.aut.2}: 
We abbreviate $E_\Gamma$ to $E$.  Note that $\bfAut(E)$ is a $k$-group scheme by Example \ref{e.AUTNAS}.
Put $\bfG$ for the constant group scheme corresponding to the permutation group on the set $\Gamma$.  There is a morphism $\bfG \to \bfAut(E)$ such that $\bfG(R) \to \bfAut(E)(R)$ is injective for every $R \in \kalg$, constructed in the same way as the injection in the case $|\Gamma| = 2$ in Exc.~\ref{pr.quad.exercise}.  Our aim is to show that this injection is an isomorphism of group schemes, which we do following Exc.~\ref{pr.ETAUT}.

Identify $\Gamma$ with $\{ 1, \ldots, n \}$ and write $e_i$ for the $i$-th standard basis vector of $k^n$.  An element $A \in \Aut(E) = \bfAut(E)(k)$ is an element of $\GL_n(k)$ that preserves the coordinate-wise multiplication on $k^n$, which we denote by $*$.  For each $i$, we have $(Ae_i) * (Ae_i) = A(e_i * e_i) = Ae_i$, i.e., every entry of $A$ is an idempotent in $k$.  For $i \ne j$, we have $0 = e_i * e_j = (Ae_i) * (Ae_j)$, so $A_{ri} A_{rj} = 0$ for all $r$.  That is, for each row of $A$, the entries form a complete orthogonal system of idempotents in $k$.

Given any two complete orthogonal systems of idempotents, the collection of products (one from each set) forms a new complete orthogonal system of idempotents.  Therefore, the collection of elements $\{ A_{1i_1} A_{2i_2} \cdots A_{ni_n} \mid i_1, \ldots, i_n \in \{ 1, \ldots n \} \}$ is a complete orthogonal system of idempotents.  Let $f_1, \ldots, f_N$ be the nonzero elements.  Then $k \cong \prod_{i=1}^N kf_i$ and we can view $A$ as an $N$-tuple of matrices with entries in $kf_1, \ldots, kf_N$.

Focus on one of these $N$ matrices.  Its entries are 0 or 1.  Since each row is a complete orthogonal system of idempotents, it has at most one 1.  Since $A$ is invertible, each row has at least one 1.  In summary, each row has exactly one 1 and since $A$ is invertible it must be a monomial matrix.

In this way, we have recovered the $n$ arbitrary version of the $n = 2$ result contained in Exc.~\ref{pr.ETAUT}(b), which shows that surjectivity of $\bfG(R) \to \bfAut(E)(R)$ as desired.

%
%
%
%

\end{sol}

\begin{sol}{pr.etale.descent} \label{sol.pr.etale.descent} \ref{pr.etale.descent.1}: Since $E$ is finitely presented and flat over $k$, so is $E_R$ over $R$.  For each $\mfq \in \Spec(R)$, put $\mfp := \mfq \cap k$, so $R(\mfq)$ is a field containing $k(\mfp)$.  Now
\[
E_R \otimes_R R(\mfq) \cong E \otimes_k R(\mfq)  \cong (E \otimes k(\mfp)) \otimes_{k(\mfp)} R(\mfq),
\]
where $E \otimes k(\mfp)$ is an \'etale algebra over the field $k(\mfp)$.  The property of being an \'etale algebra over a field is stable under enlarging the field \cite[V.6, Cor.~2]{MR643362}, i.e., the tensor product on the right is an \'etale $R(\mfq)$-algebra, completing the proof of the claim.

\smallskip
\ref{pr.etale.descent.2}: The ``only if'' direction is trivial by taking $R = k$, so we prove ``if''.  Suppose $R \in \kalg$ is faithfully flat and $E_R$ is \'etale.  Because $E_R$ is finitely presented and projective over $R$, the same is true for $E$ over $k$ by Cor.~\ref{c.FAFLAD} and \ref{skip.other.fflat}.

For each $\mfp \in \Spec(k)$, Exc. \ref{pr.COVERS} provides a field $L \in \Ralg$ containing $k(\mfp)$.  By hypothesis, $(E_{k(\mfp)})_L = E_L$ is a finite direct product of separable field extensions of $L$, so it is a separable algebra.  It follows by \cite[V.15, Prop.~3d]{MR643362}  that $E_{k(\mfp)}$ is a separable $k(\mfp)$-algebra of finite rank, and so is a finite product of finite separable extensions of $k(\mfp)$ by ibid., V.6, Thm.~4, proving that $E$ is \'etale.
\end{sol}

\solnsec{Section~\ref{s.FPPFCA}}

\begin{sol}{pr.FFSPLIQ} \label{sol.FFSPLIQ}
 By Exc.~\ref{pr.EMMA}, the $k$-functor $M_\bfa^{2n} = (M^{2n})_\bfa$ of Example~\ref{e.EXMA} is a finitely presented $k$-scheme. Hence so is $\bfX := \bfHyp(Q)$, defined as a closed subfunctor of $M_\bfa^{2n}$ by finitely many equations (cf. \eqref{HYPB} and Exc.~\ref{pr.IDFIPR}~(b)). After an appropriate base change, it therefore remains to show that the natural set map $\bfX(k) \to \bfX(\bar k)$ is surjective, where $\bar k := k/I$ and $I \subseteq k$ is an ideal satisfying $I^2 = \{0\}$. Note that \ref{ss.REDID} and Cor.~\ref{c.BAQUA} imply $\bar M := M_{\bar k} = M/IM$, and the quadratic form $\bar q := q_{\bar k}\:\bar M \to \bar k$ over $\bar k$ is given by $\bar q(\bar x) = \overline{q(x)}$ for all $x \in M$, where $\alpha \mapsto \bar\alpha$ (resp. $x \mapsto \bar x$) stand for the natural maps $k \to \bar k$ (resp. $M \to \bar M$). We first note that a hyperbolic basis is indeed a basis of the underlying module and must show that any hyperbolic basis of $\bar Q := (\bar M,\bar q)$ can be lifted to one of $Q$, so let $(\pw_i)_{1\leq i\leq 2n}$ be a hyperbolic basis of $\bar Q$. We argue by induction on $n$. While the case $n = 0$ is trivial, the case $n = 1$ amounts to showing that a hyperbolic pair $(\pu,\pv)$ of $\bar Q$ in the sense of \ref{ss.ISVE} can be lifted to a hyperbolic pair of $Q$. To this end, we pick any elements $a,b \in M$ such that $\bar a = \pu$, $\bar b = \pv$. Then $q(a),q(b) \in I$, $q(a,b) \in 1 + I \subseteq k^\times$, and since $I$ squares to zero, one checks that
\[
u := a - q(a,b)^{-1}q(a)b, \quad v := -q(a,b)^{-2}q(b)a + q(a,b)^{-1}b
\]
form a hyperbolic pair of $Q$ lifting $(\pu,\pv)$. Now assume $n > 1$. By what we have just seen, the hyperbolic pair $(\pw_n,\pw_{2n})$ of $\bar Q$ can be lifted to a hyperbolic pair $(w_n,w_{2n})$ of $Q$. Put $V := kw_n \oplus kw_{2n}$ and $M_0 := V^\perp$ (relative to $Dq$). Then Lemma~\ref{l.NOSU} implies $M = M_0 \perp V$, and $Q_0 := (M_0,q\vert_{M_0})$ is a quadratic space of rank $2(n - 1)$ over $k$. Moreover, $\bar Q_0 = (\bar M_0,\bar q\vert_{\bar M_0})$ contains $(\pw_i)_{1\leq i\leq 2(n-1)}$ as a hyperbolic basis, which, by the induction hypothesis, can be lifted to a hyperbolic basis $(w_i)_{1\leq i\leq 2(n-1)}$ of $Q_0$. Hence $(w_i)_{1\leq i\leq 2n}$ is the desired lift of $(\pw_i)_{1\leq i\leq 2n}$.

(b) We proceed in several steps. 

\step{1}
For the time being, we dispense ourselves from the quadratic space $Q$ of (a) and consider instead arbitrary quadratic modules $Q = (M,q)$, $\pQ= (\pM,\pq)$ over $k$. By an \emph{isometry} from $Q$ to $\pQ$ we mean a $k$-linear bijection $\eta\:M \to \pM$ satisfying $\pq \circ \eta = q$. We define the set
\begin{align*}
\Isom(Q,\pQ) := \{\eta \mid \eta\:Q \to \pQ\;\text{is an isometry}\},
\end{align*}
which in general will be empty but gives rise to a $k$-functor
\[
\bfIsom_k(Q,\pQ) := \bfIsom(Q,\pQ)\:\kalg \rightarrow \set
\]
by defining
\begin{align}
\label{ISQU} \bfIsom(Q,\pQ)(R) := \Isom_R(Q_R,\pQ_R)
\end{align}
for all $R \in \kalg$ and by
\begin{align}
\label{ISOMET} \bfIsom(Q,\pQ)(\vph)\:&\bfIsom(Q,\pQ)(R) \longrightarrow \bfIsom(Q,\pQ)(S), \\
\Isom_R(Q_R,\pQ_R&) \ni \eta \longmapsto \eta_S \in \Isom_S(Q_S,\pQ_S) \notag
\end{align}
for all morphisms $\vph\:R \to S$ in $\kalg$. The $k$-group functor $\bfOr(Q)$ of Example~\ref{e.ORGR} acts canonically in a simply transitive manner on $\bfIsom(Q,\pQ)$ from the right via
\[
\Isom_R(Q_R,\pQ_R) \times \Or(Q_R) \longrightarrow \Isom_R(Q_R,\pQ_R), \quad (\eta,\zeta) \longmapsto \eta \circ \zeta
\]
for all $R \in \kalg$. 

\step{2}
After this digression, we return to our quadratic space $Q = (M,q)$ of rank $2n$ over $k$ as considered in (a).We denote by $\bfh_k^{2n}$ the split hyperbolic quadratic space of rank $2n$ over $k$ defined on $k^{2n} = \left(\begin{smallmatrix}
k^n \\
k^n
\end{smallmatrix}\right)$ by the quadratic form
\begin{align}
\label{SPLIHYP} \la S\raq, \quad S := \left(\begin{matrix}
0 & \Eins_n \\
0 & 0
\end{matrix}\right)
\end{align}
in the sense of (\ref{ss.MAFO}.\ref{QUAS}). The canonical basis of $k^{2n}$ is a hyperbolic basis for $\bfh_k^{2n}$ and obviously compatible with base change. It follows immediately from the definitions that the assignment
\begin{align}
\label{ETEI} \eta \longmapsto \big(\eta(e_i)\big)_{1\leq i\leq 2n}
\end{align}
defines a bijection
\begin{align}
\label{ISHY} \Phi = \Phi(k)\:\Isom(\bfh_k^{2n},Q) \overset{\sim} \longrightarrow \Hyp(Q).
\end{align}

\step{3}
Putting $\bfX:= \bfIsom(\bfh_k^{2n},Q)$, the set maps
\[
\Phi(R)\:\bfX(R) \longrightarrow \bfHyp(Q)(R)
\]
given by \eqref{ETEI}, \eqref{ISHY} for all $R \in \kalg$ are bijective and obviously compatible with base change, hence give rise to an isomorphism
\[
\Phi\:\bfX \overset{\sim} \longrightarrow \bfHyp(Q)
\]
of $k$-functors. By (a), therefore, $\bfX$ is a smooth affine $k$-scheme. Moreover, $\bfX$ is faithful in the sense of \ref{ss.FAIAFF} since regular even-dimensional quadratic forms over an algebraically closed field are (split) hyperbolic \cite[Exc.~7.34]{MR2427530}. Summing up, \ref{ss.THREEFF}~(ii) implies that there exists a faithfully flat \'etale $k$-algebra $R$ satisfying $\bfX(R) \neq \emptyset$. But this means that $Q_R$ is split hyperbolic. Furthermore, $\bfX$ is a torsor with structure group $\bfOr(\bfh_k^{2n})$ in the \'etale topology, forcing
\[
\bfOr(Q)_R \cong \bfOr(Q_R) \cong \bfOr(\bfh_k^{2n}) \cong \bfX_R
\]
to be smooth over $R$. By \ref{ss.THREEFF}~(iii), therefore $\bfOr(Q)$ is smooth over $k$.
\end{sol}



\solnchap{Solutions for Chapter~\ref{c.JORAL}}

\solnsec{Section~\ref{s.LIJO}}

\begin{sol}{pr.jord.assoc}
(a): $A^+$ is associative if and only if the following is zero for all $x, y, z \in A$:
\begin{align*}
4(x\bu(y\bu z) - (x \bu y) \bu z) &= xzy + yzx - yxz - zxy = [x,z]y - y[x,z] \\
&= [[x,z],y].
\end{align*}
Since $A$ is associative, $[x,z]$ is central if and only if $[[x,z],y] = 0$ for all $y \in A$, proving the claim.

(b): Consider first the case $n = 2$ and take 
\[
x = \stbtmat{0}{1}{0}{0} \quad \text{and} \quad y = \stbtmat{0}{0}{1}{0}.
\]
Then
\[
[x,y] = \stbtmat{1}{0}{0}{-1}
\]
and $[[x,y],x] = 2x \ne 0$, so $[x,y]$ is not central.  In view of part (a), we conclude that $\Mat_2(k)^+$ is not associative.

For $n > 2$, we can perform the same argument, where the role of $x$ is played by a matrix with the 2-by-2 $x$ in the upper left corner and zeros elsewhere, and similarly for $y$.  The same calculations show that $[x,y]$ is not central and therefore that $\Mat_n(k)^+$ is not associative.
\end{sol}


\solnsec{Section~\ref{s.PAQUA}}

\begin{sol}{pr.PARNIL} \label{sol.PARNIL}
(a) By induction on $m$. For $m = 0,1$ there is nothing to prove. Now let $m > 1$ and suppose the assertion has been established for all natural numbers $< m$. By the induction hypothesis we have $x^m = U_xx^{m-2} \in \Mon_m(\{x\})$. Conversely, let $x \in \Mon_m(\{x\})$. Then $x = U_yz$, $y \in \Mon_n(\{x\})$, $z \in \Mon_p(\{x\})$, $n,p \in \IN$, $n > 0$, $2n + p = m$. By the induction hypothesis, $y = x^n$, $z = x^p$, and power associativity yields $x = U_{x^n}x^p = x^{2n+p} = x^m$, as claimed. 

(b) For any subset $X \subseteq J$ and any homomorphism $\varphi\:J \to J^\prime$ of para-quadratic $k$-algebras, one makes the following observations by straightforward induction:
\begin{enumerate}[(i)]
\item $\varphi(\Mon_m(X)) = \Mon_m(\varphi(X))$ for all $m \in \IN$.

\item $\Mon_m(Y) \subseteq \Mon_{mn}(X)$ for all $Y
\subseteq \Mon_n(X)$ and all $m,n \in \IN$.
\end{enumerate}
If $x \in J$ is nilpotent, then (i) implies that $\varphi(x) \in J^\prime$
is nilpotent. Hence if $I$ is a nil ideal in $J$, then $I^\prime$ and $I/I^\prime$ are nil ideals in $J$ and $J/I^\prime$, respectively. Conversely, let this be so. Writing $\pi\:J \to J/I^\prime$ for the canonical epimorphism, any $x \in I$ makes $\pi(x) \in I/I^\prime$ nilpotent, so by (i), $\Mon_n(\{x\})$ meets $I^\prime$ for
some $n \in \IN$. But since $I^\prime$ is nil, we conclude that some $y
\in \Mon_n(\{x\})$ is nilpotent, leading to a positive integer $m$
such that $0 \in \Mon_{mn}(\{x\})$ by (ii). Thus $x$ is nilpotent,
forcing $I$ to be a nil ideal. Now write $\mfn := \Nil(J)$ for the
sum of all nil ideals in $J$. Given $x \in \mfn$, there are finitely
many nil ideals $\mfn_1,\dots,\mfn_r \subseteq J$ such that $x \in
\mfn_1 + \cdots + \mfn_r$. Thus we only have to show that the sum of
finitely many nil ideals in $J$ is nil. Arguing by induction, we are
actually reduced to the case $r = 2$. But then the isomorphism
$(\mfn_1 + \mfn_2)/\mfn_1 \cong \mfn_2/(\mfn_1 \cap \mfn_2)$
combines with what we have proved earlier to yield the assertion. 

(c) It suffices to prove that $\Nil(k)J \subseteq J$ is a nil ideal. The ideal property being obvious, we are reduced to showing that $\Nil(k)J$ consists entirely of nilpotent elements. An arbitrary element $x \in \Nil(k)J$ actually belongs to $IJ$, where $I \subseteq \Nil(k)$ is some finitely generated ideal in $k$. A straightforward induction now shows $x^n \in I^nJ$ for all $n \in \IN$, and the assertion follows from the lemma in solution to Exc.~\ref{pr.NILRAD}. 
\end{sol}

\begin{sol}{pr.EVALPOL} \label{sol.EVALPOL}
(a) For the first assertion, since $\vep(1) = 1_J$, we need only show $(f^2g)(x) = U_{f(x)}g(x)$ for all $f,g \in k[\bft]$. By linearity, we may assume $g = \bft^n$ for some $n \in \IN$ and, writing $f = \sum\alpha_i\bft^i$ with coefficients $\alpha_i \in k$, we obtain
\begin{align*}
f^2g =\,\,&(\sum\alpha_i^2\bft^{2i} + 2\sum_{i<j}\alpha_i\alpha_j\bft^{i+j})\bft^n = \sum\alpha_i^2\bft^{2i+n} + 2\sum_{i<j}\alpha_i\alpha_j\bft^{i+j+n}, 
\end{align*} 
hence, as $J$ is power-associative,
\begin{align*}
(f^2g)(x) =\,\,&\sum\alpha_i^2x^{2i + n} + \sum_{i<j}\alpha_i\alpha_j2x^{i+j+n} = \sum\alpha_i^2U_{x^i}x^n + \sum_{i<j}\alpha_i\alpha_j\{x^ix^nx^j\} \\
=\,\,&\big(\sum\alpha_i^2U_{x^i} + \sum_{i<j}\alpha_i\alpha_j U_{x^i,x^j}\big)x^n = U_{\sum\alpha_ix^i}x^n = U_{f(x)}g(x).
\end{align*}
Thus $\vep_x$ is a homomorphism of para-quadratic algebras, and we conclude that $I \subseteq k[\bft]^{(+)}$ is an ideal. Our next aim will be to show that $I^0$ is an ideal in $k[\bft]$. By linearity, it suffices to show for $f \in I^0$ that $\bft^nf$ belongs to $I^0$ for all $n \in \IN$. We do so by induction on $n$. For $n = 0$ there is nothing to prove. Now assume $n > 0$ and that the assertion holds for $n - 1$. Since $\vep_x$ is a homomorphism of para-quadratic algebras, the induction hypothesis yields not only $(\bft^nf)(x) = 0$, but also $(\bft^{n+1}f)(x) = (\bft^2\bft^{n-1}f)(x) = U_x(\bft^{n-1}f)(x) = 0$, and we have shown $\bft^nf \in I^0$. Next let $I^{0\prime}$ be any ideal of $k[\bft]$ contained in $I$. Given $f \in I^{0\prime} \subseteq I$, we also have $\bft f \in I^{0\prime} \subseteq I$, hence $f(x) = (\bft f)(x) = 0$, which show $I^{0\prime} \subseteq I^0$, so $I^0$ is indeed the largest ideal of $k[\bft]$ contained in $I$. Finally, let $f \in I$. Then $f^2(x) = U_{f(x)}1_J = 0$ and $(\bft f^2)(x) = U_{f(x)}x = 0$, hence $f^2 \in I^0$. Similarly, not only $2f \in I$ but also $(2\bft f)(x) = (\bft \circ f)(x) = x \circ f(x) = 0$, hence $2f \in I^0$. The assertions about $R$ and $\pi$ making a commutative diagram as shown in the exercise are now obvious, as is the statement that $k[x] \subseteq J$ is a para-quadratic subalgebra. The assertion about $\Ker(\pi) = I/I^{(0)}$ follows immediately from the fact that $f^2$ and $2f$ belong to $I^{(0)}$ for all $f \in I$. This completes the solution to (a).

(b) Write $f = \alpha_n\bft^n + h \in I_n$ with $\alpha_n \in k$ and $h \in \bft^{n+2}k[\bft]$. Since $2 = 0$ in $k$ and $n \geq 2$, we have $f^2 = \alpha_n^2\bft^{2n} + h^2 \in \bft^{n+2}k[\bft]$, which is an ideal in $k[\bft]$. Thus $U_fg = f^2g \in \bft^{n+2}k[\bft] \subseteq I_n$ for all $g \in k[\bft]$. Conversely, let $g = \sum\beta_i\bft^i \in k[\bft]$ with coefficients $\beta_j \in k$ for $j \in \IN$. Then $U_gf =g^2f \equiv \alpha_n\beta_0^2\bft^n \bmod \bft^{n+2}k[\bft]$ and hence $U_gf \in I_n$. Since $\{g_1g_2f\} = 2fg_1g_2 = 0 \in I_n$, it follows that $I_n$ is an ideal in $k[\bft]^{(+)}$. But it is not an ideal in $k[\bft]$ since $\bft^n$ belongs to $I_n$ but $\bft\bft^n = \bft^{n+1}$ does not. And finally, if $J_n \cong A^{(+)}$ for some unital flexible $k$-algebra $A$, then $x^n = 0$ would imply the contradiction $x^{n+1} = 0$ since powers in $A$ and $A^{(+)}$ coincide.
\end{sol}

\begin{sol}{pr.PARALIFT} \label{sol.PARALIFT}
(a) In the para-quadratic $k$-algebra $k[\bft]^{(+)}$ we have $U_{fg}h = (fg)^2h
= f^2g^2h = U_fU_gh$ for all $f,g,h \in k[\bft]$, hence 
\begin{align}
\label{UFG} U_{fg} =U_fU_g &&(f,g \in k[\bft]).
\end{align}
Since $J$ is power-associative, the evaluation at $x$ (Exc.~\ref{pr.EVALPOL}) is a homomorphism $k[\bft]^{(+)} \to J$ of para-quadratic algebras with image $k[x]$. Applying this homomorphism to \eqref{UFG} therefore gives \eqref{MULU} of the exercise, hence not only the first equation of \eqref{COMU} of the exercise but also
\begin{align*}
U_{U_{f(x)}^\prime g(x)}^\prime =\,\,&U_{(f^2g)(x)}^\prime = U_{f^2(x)}^\prime U_{g(x)}^\prime = U_{f(x)}^{\prime 2}U_{g(x)}^\prime, \\
U_{x^n}^\prime =\,\,&U_{\bft^n(x)}^\prime = U_{\bft(x)}^{\prime n} = U_x^{\prime n},
\end{align*}
and the proof of (a) is complete.

(b) The polynomial
\[
f := \alpha_d\bft^d + \cdots + \alpha_{n-1}\bft^{n-1} + \bft^n \in k[\bft] 
\] kills $x$ and hence belongs to $I := \Ker(\vep_x)$. From Exc.~\ref{pr.EVALPOL}~(a) we therefore conclude that
\[
f^2 = \alpha_d^2\bft^{2d} + \cdots + \bft^{2n} \in k[\bft]
\]
belongs to $I^0$, so by Exc.~\ref{pr.EVALPOL}~(a) we have $(\bft^if^2)(x) = 0$, i.e., 
\begin{align}
\label{EFSQ} \alpha_d^2x^{2d+i} + \cdots + x^{2n+i} = 0 &&(i \in \IN).
\end{align}
Now, setting $k_r[x] := \sum_{s\geq r}kx^s$ for $r \in \IN$, we find in $(k_r[x])_{r\in \IN}$ a descending chain of submodules of $k[x]$. Invoking \eqref{EFSQ}, we therefore conclude that $k_{2d}[x] = k_{2d+1}[x] = k_{2(d+1)}[x]$ is a finitely generated $k$-module. Thus the linear map $U_x^\prime\:k_{2d}[x] \to k_{2d}[x]$, being surjective, is in fact bijective (\ref{p.SURBI}). This property carries over to $(U_x^\prime)^d = U_{x^d}^\prime\:k_{2d}[x] \to k_{2d}[x]$, whence there exists a unique element $v \in k_{2d}[x]$ such that $U_{x^d}^\prime v = x^{2d}$. Here (a) implies $U_{x^{2d}}^\prime U_v^\prime = U_{x^{2d}}^\prime$, and we conclude that $U_v$ is the identity on $k_{2d}[x]$. With $c := v^2$, this means $U_c = (U_v)^2 = \Eins$ on $k_{2d}[x]$ and $c^2 = v^4 = U_vv^2 = c$, $c^3 = U_cc = c$, so $c \in k_{2d}[x]$ is an idempotent of the desired kind. 

(c) Let $x \in J$ be a pre-image of $c^\prime$ under $\vph$. Then
$x^2 - x$ is nilpotent, so some integer $d > 0$ has $0 = (x^2 - x)^d
= x^{2d} - \cdots + (-1)^dx^d$. Pick $v \in k_{2d}[x]$ as in (b) and put $c := v^2$. Writing $\bar u := \vph(u)$ for $u \in J$, we conclude $U_{c^\prime}\bar v = U_{\bar x}\bar v = U_{\bar x^d}\bar v = \bar x^{2d} = c^\prime$. On the other hand, $v =
\sum_{i\geq 2d}\alpha_ix^i$ for some scalars $\alpha_i \in k$, $i
\geq 2d$, so with $\alpha := \sum_{i\geq 2d}\alpha_i$ we obtain $\bar
v = \bar\alpha c^\prime$, hence $c^\prime = U_{c^\prime}\bar v = \bar\alpha U_{c^\prime}c^\prime = \bar\alpha c^\prime = \bar v$. But this implies $\bar c = \bar v^2 = c^{\prime 2} = c^\prime$, and the problem is solved. 
\end{sol}

\begin{sol}{pr.PAPONI} \label{sol.PAPONI} (i) $\Rightarrow$ (ii). We may assume $J \ne \{0\}$. Since $x$ is nilpotent and $J$ is power-associative at $x$, some positive integer $m$ has $0 \in \Mon_m(\{x\}) = \{x^m\}$ (by Exc.~\ref{pr.PARNIL}~(a)). Hence $x^m = 0$.

(ii) $\Rightarrow$ (iii). Let $n \in \IZ$, $n \ge 2m$. Then $x^n = U_{x^m}x^{n-2m} = 0$. 

(iii) $\Rightarrow$ (i). $0 = x^m \in \Mon_m(\{x\}) \subseteq \Mon(\{x\})$. Hence $x$ is nilpotent.

It remains to show that
\[
I := \{x \in k[x] \mid x\,\text{is nilpotent}\}
\]
is an ideal in $k[x]$. To this end, let $v,w \in k[x]$. We claim
\begin{align}
\label{UPOW} (U_vw)^n = U_{v^n}w^n    
\end{align}
for all $n \in \IN$. This is clear for $n = 0,1$. Assuming $n \ge 2$ and arguing inductively, we use notation and results of Exc.~\ref{pr.PARALIFT} to obtain
\begin{align*}
 (U_vw)^n =\,\,&(U_v^\prime w)^n = U^\prime_{U^\prime_vw}(U_vw)^{n-2} = U_v^{\prime 2}U_w^\prime U_v^{\prime n-2}w^{n-2} = U_{v^n}^\prime w^n = U_{v^n}w^n,  
\end{align*}
and the proof of \eqref{UPOW} is complete. From \eqref{UPOW} we conclude that $U_vw$ belongs to $I$ provided $v$ or $w$ do. Next we show that that $I$ is additively closed, so let $v,w \in I$ be written as $v = f(x)$, $w = g(x)$ for some $f,g \in k[\bft]$. Applying (ii) above, we find a positive $n$ such that $v^n = w^n = 0$. Then
\[
(v + w)^{4n} = (f + g)^{4n}(x) = \Big(\sum_{i=0}^{4n} \binom{4n}{i} f^ig^{4n-i}\Big)(x).
\]
For $0 \le i \le 2n$, we obtain
\[
(f^ig^{4n-i})(x) = (g^{2n}f^ig^{2n-i})(x) = U_{g^n(x)}\big((f^ig^{2n-i})(x)\big) = U_{w^n}\big((f^ig^{2n-i})(x)\big) = 0
\]
and similarly, for $2n < i \le 4n$,
\[
(f^ig^{4n-i})(x) =  (f^{2n}f^{i-2n}g^{4n-i})(x) = U_{f^n(x)}\big((f^{i-2n}g^{4n-i})(x)\big) = U_{v^n}\big((f^{i-2n}g^{4n-i}(x)\big) = 0.
\]
Summing up, $v + w \in I$, as desired. Finally we have to show $\{vwz\} \in I$ for $v,w \in k[x]$ and $z \in I$. With $v = f(x)$, $w = g(x)$ as before and $z = h(x)$, $h \in k[\bft]$, we obtain $\{vwz\} =2(fgh)(x)$, hence
\[
\{vwz\}^n = \big(2fgh\big)^n(x) = 2^n(f^ng^nh^n)(x) = 2^{n-1}\{v^nw^nz^n\}
\]
for all integers $n > 0$, from which the assertion follows.
\end{sol}

\begin{sol}{pr.OUTCEN} \label{sol.OUTCEN}
We closely follow the arguments used in the proofs of \ref{ss.c-CENTSIM1}$-$\ref{ss.c-CENTSIM4}.

(a) Since $J$ is outer simple, $\Mult(J)$ acts irreducibly on $J$. Thus the assertion follows from Schur's lemma \cite[p.~118]{MR1009787}.

(b)  Since $\Mult(J)$ acts faithfully and irreducibly on $J$, it is a primitive Artinian $k$-algebra \cite[Def.~4.1]{MR1009787}. Moreover, by definition, its centralizer in $\End_k(J)$ is $\Cent_0(J)$. Thus the assertion follows from the double centralizer theorem \cite[Thm.~4.10]{MR1009787}.

(c) If $J$ is outer central and outer simple, the assertion follows from
(b). Conversely, suppose $J \neq \{0\}$ and
$\Mult(J) = \End_k(J)$. Then $J$ is outer simple, and
(b) implies that $\Cent_0(J)$ belongs to the centre of $\End_k(J)$. Hence $J$ is outer central.

(d) (i) $\Rightarrow$ (ii). Let $k^\prime$ be an extension
field of $k$ and put $J^\prime = J \otimes k^\prime$ as a para-quadratic
$k^\prime$-algebra. After identifying $\End_{k^\prime}(J^\prime) =
\End_k(J) \otimes k^\prime$ canonically, a moment's reflection shows
$\Mult(J^\prime) = \Mult(J) \otimes k^\prime$. Hence the assertion
follows from (c). \\
(ii) $\Rightarrow$ (iii). Obvious. \\
(iii) $\Rightarrow$ (i). If $I \subseteq J$ is a non-trivial outer ideal, then so is $I \otimes \bar k \subseteq J \otimes \bar k$, a contradiction. Hence $J$ is outer simple. Since the multiplication algebra of $J$ behaves nicely under base field extensions, so does its centralizer. Therefore $\Cent_0(J) \otimes \bar k$ agrees with the outer centroid of $J \otimes \bar k$ and hence is a finite algebraic field extension of $\bar k$. As such, it has degree 1, forcing $\Cent_0(J) = k1_J$, and $J$ is outer
central.
\end{sol}

\begin{sol}{pr.COS} \label{sol.COS}
(a) Setting $c := \sum_{i=1}^rc_i$, we have
\begin{align*}
c^2 =\,\,&\sum_{i=1}^rc_i^2 + \sum_{i<j}c_i \circ c_j = \sum_{i=1}^rc_i = c, \\
c^3 =\,\,&U_cc = \sum_{i,j=1}^rU_{c_i}c_j + \sum_{i<l}\sum_{j=1}^r\{c_ic_jc_l\} = \sum_{i=1}^rc_i^3 = \sum_{i=1}^rc_i = c, 
\end{align*}
and we have shown that $c$ is an idempotent in $J$. Hence so is $c_{r+1} := 1_J - c$ (\ref{ss.PAIM}) and it remains to show that $(c_1,\dots,c_r,c_{r+1})$ is an orthogonal system of idempotents, completeness being obvious. For $1 \leq i,l \leq r$ distinct, we compute
\begin{align*}
U_{c_i}c_{r+1} =\,\,&U_{c_i}1_J - \sum_{j=1}^rU_{c_i}c_j = c_i^2 - U_{c_i}c_i = c_i - c_i^3 = c_i - c_i = 0, \\
U_{c_{r+1}}c_i =\,\,&U_{1_J-c}c_i = c_i - c \circ c_i + U_cc_i = c_i - \sum_{m=1}^rc_m \circ c_i + \sum_{m=1}^rU_{c_m}c_i + \sum_{m<n}\{c_mc_ic_n\} \\
=\,\,&c_i - 2c_i + c_i = 0, \\
\{c_ic_ic_{r+1}\} =\,\,& c_i \circ c_i - \sum_{m=1}^r\{c_ic_ic_m\} = 2c_i^2 - 2U_{c_i}c_i = 2c_i^2 - 2c_i^3 = 0, \\
\{c_{r+1}c_{r+1}c_i\} =\,\,&\{(1_J - c)(1_J - c)c_i\} = 2c_i - c \circ c_i - c \circ c_i + \{ccc_i\} \\
=\,\,&2c_i - 2c_i - 2c_i + \sum_{m,n=1}^r\{c_mc_nc_i\} = -2c_i + 2c_i = 0, \\
c_i \circ c_{r+1} =\,\,&c_i \circ 1_J - c_i \circ c = 2c_i - 2c_i = 0, \\
\{c_ic_lc_{r+1}\} =\,\,&c_i \circ c_l - \sum_{m=1}^r\{c_ic_lc_m\} = 0.
\end{align*}
Hence $(c_1,\dots,c_r,c_{r+1})$ is a complete orthogonal system of idempotents in $J$.

(b) If $(c_1,\dots,c_r)$ is a (complete) orthogonal system  of idempotents in $A$, then it is clearly one in $J$. Conversely, suppose it is a (complete) orthogonal system of idempotents in $J$. Comparing \eqref{PACOSY} with (\ref{ss.PAIM}.\ref{PAORID}), we see that the idempotents $c_i,c_j$, $1 \leq i < j \leq r$, are orthogonal in $J$, hence, as we have seen in \ref{ss.PAIM}, in $A$. But this means that $(c_1,\dots,c_r)$ is a (complete) orthogonal system of idempotents in $A$.
\end{sol}


\solnsec{Section~\ref{s.JABAS}}

\begin{sol}{pr.CHRAJO} \label{sol.CHRAJO}  (i) $\Rightarrow$ (iii). Obvious.

(iii) $\Rightarrow$ (i). Since \eqref{RUUXY}, \eqref{RUXVY} hold in all scalar extensions of $J$, by Cor.~\ref{c.VANCRI}, $J$ is a Jordan algebra.

(i) $\Rightarrow$ (ii). This follows from \ref{ss.JABAS} since the identities \eqref{RUUXYZ}--\eqref{RUXZVY} are respectively the same as (resp. part of) (\ref{ss.JABAS.fig}.\ref{UUXYZ})--(\ref{ss.JABAS.fig}.\ref{UXZVY}).

(ii) $\Rightarrow$ (i). We will show successively that \eqref{RUUXYUZ}--\eqref{RUXZVY}, \eqref{RUUXYZ}, \eqref{RUXVY}, \eqref{RUUXY} hold in $J_R$, for all $R \in \kalg$. In the sense of Exc.~\ref{pr.MULQUA}, the left (resp. right) of \eqref{RUUXYUZ} defines a $k$-$3$-quadratic map $F$ (resp. $G$) from $J^3$ to $\End_k(J)$, and we have $F = G$ by hypothesis. Thus $F_R = G_R$ by Exc.~\ref{pr.MULQUA}, which by \ref{ss.PASCA}implies that \eqref{RUUXYUZ} holds in $J_R$. For analogous reasons, \eqref{RUUXW}, \eqref{RUXZVY} holds in $J_R$ as well.

\smallskip

\emph{Proof of} \eqref{RUUXYZ}. By linearity in $z$, we may assume $z = w_R$, $w \in J$. Assume first $x = u_R$, $u \in J$. Since both sides of \eqref{RUUXYZ} are quadratic in $y$ and agree for $y = v_R$, $v \in J$, by hypothesis, it follows that \eqref{RUUXYZ} holds for all $y \in J_R$. Now fix $y \in J_R$ and put
\[
M := \{x \in J_R \mid U_xU_yU_{x,z} + U_{x,z}U_y U_x = U_{U_xy,U_{x,z}y}\},
\]
which is obiously closed under multiplication by scalars in $R$ and, as we have just seen, contains all elements $u_R$, $u \in J$. Hence \eqref{RUUXYZ} will follow once we have shown that $M$ is closed under addition. But this is implied by \eqref{RUUXW} and a straightforward computation. 

\smallskip

\emph{Proof of} \eqref{RUXVY}. By linearity in $y$, we may assume $y = v_R$, $v \in J$. Arguing as before, it suffices to show that
\[
M := \{x \in J_R \mid U_xV_{y,x} = V_{x,y}U_x\}
\] 
is closed under addition. But this is immediate, by \eqref{RUXZVY} and a straightforward computation. 

\smallskip

\emph{Proof of} \eqref{RUUXY}.  Since \eqref{RUUXY} is quadratic in $y$, it suffices to prove it for $y = v_R$, $v \in J$. Again we will be through once we have shown that
\[
M := \{x \in J_R \mid U_{U_xy} = U_xU_yU_x\}
\]
is additively closed. This follows from \eqref{RUUXYZ}, \eqref{RUUXYUZ} by a lengthy, but still straightforward, computation.

Since a Jordan algebra satisfies \eqref{RUUXY}--\eqref{RUXZVY}, so do all its subalgebras and homomorphic images. By what we ahve just shown, therefore, they are all Jordan algebras.
\end{sol}

\begin{sol}{pr.OUTIJO} \label{sol.OUTIJO}  (a) By definition it suffices to show that $I$ is an outer ideal if $U_JI \subseteq I$. Indeed, if this inclusion holds, then, for all $x,y \in J$, we have $U_{x,y}I \subseteq I$ and given $z \in I$, equations (\ref{ss.JABAS.fig}.\ref{UXYV}), (\ref{ss.JABAS.fig}.\ref{VXUX}) imply
\[
\{xyz\} = V_{x,y}z = V_xV_yz - U_{x,y}z = U_{x,1_J}U_{y,1_J}z - U_{x,y}z \in I,
\]
which shows $\{JJI\} \subseteq I$ and completes the proof.

(b) If $c^2 = c$, then , then $U_c^2 = U_{c^2} = U_c$ by (\ref{ss.TOPOAS}.\ref{UXN}), hence $c^3 = U_cc = U_cc^2 = U_cU_c1_J = U_c^21_J = U_c1_J = c^2 = c$. Thus $c$ is an idempotent.

(c) Let $x,y \in J$ and $z \in \Rex(J)$. Then (\ref{ss.JABAS.fig}.\ref{UXYV}) implies $V_{x,z} = V_xV_z - U_{x,z} = V_xU_{z,1_J} - U_{z,x} = 0$ and, similarly, $V_{z,x} = 0$. From the fundamental formula (\ref{ss.JABAS.fig}.\ref{UUXY}) we deduce $U_{U_xz} = 0 = U_{U_zx}$. And, finally, applying (\ref{ss.JABAS.fig}.\ref{UXZVY}), we obtain
\begin{align*}
U_{U_xz,y} =\,\,&U_{y,U_xz} = U_{x,y}V_{z,x} + U_xV_{z,y}- U_{x,U_{x,y}z} = -U_{x,V_{x,z}y} = 0, \\
U_{U_zx,y} =\,\,&U_{y,U_zx} = U_{z,y}V_{x,z} + U_zV_{x,y} - U_{z,U_{z,y}x} = 0.
\end{align*}
Hence $U_J\Rex(J) + U_{\Rex(J)}J \subseteq \Rex(J)$, and we have shown by (a) that $\Rex(J) \subseteq J$ is an ideal. Now suppose that $J$ is simple. Then $J \neq \{0\}$ by definition, hence $U_{1_J} = \Eins_J \neq 0$, which implies $\Rex(J) = \{0\}$ by simplicity. The final statement now follow from Thm.~\ref{t.PACENSI}
\end{sol}

\begin{sol}{pr.CENCENT} \label{sol.CENCENT}  Let $a \in \Cent(J)$. Then $L_aL_x = L_xL_a = L_{ax}$ for all $x \in J$. Since, therefore, $L_a$ and $L_x$ commute, so do $L_a$ and $U_x = 2L_x^2 - L_{x^2}$, for all $x \in J$. Moreover, this and $a^2 \in \Cent(J)$ imply
\begin{align*}
U_{L_ax} =\,\,&2L_{ax}^2 - L_{(ax)^2} = 2L_a^2L_x^2 - L_{a^2x^2} = 2L_a^2L_x^2 - L_{a^2}L_{x^2} = L_a^2U_x, \\
U_{L_ax,y} =\,\,&2\big(L_{ax}L_y + L_yL_{ax} - L_{(ax)y}\big) = 2L_a(L_xL_y + L_yL_x - L_{xy}) =L_aU_{x,y}
\end{align*}
for all $x,y \in J$, whence (\ref{ss.PACE}.\ref{PACE}) yields $L_a \in \Cent(\Jquad)$. Conversely,  let $\vph \in \Cent(\Jquad)$. By (\ref{ss.PACE}.\ref{AUV}) we have $\vph V_x = V_x\vph$ for $x \in J$, and from \ref{ss.COLIJO} we deduce $\vph L_x = L_x\vph$, hence $\vph(xy) = x\vph(y)$ for all $x,y \in J$. Setting $y = 1_J$, we obtain $\vph = L_a$ with $a =\vph(1_J)$. In particular, $L_a$ commutes with $L_x$ for all $x \in J$, whence $a$ belongs to the centre of $J$. Summing up we have shown that the injective algebra homomorphism $L\:\Cent(J) \to \End_k(J)$ has image $\Cent(\Jquad)$, and the assertion follows.
\end{sol}

\begin{sol}{pr.CHAJOCLI} \label{sol.CHAJOCLI}  If $J = J(M,q,e)$, then \eqref{QUACUB} holds strictly by (\ref{ss.JOPOID.fig}.\ref{JOQUAD}), (\ref{ss.JOPOID.fig}.\ref{JOQUAT}). Conversely, suppose \eqref{QUACUB} holds strictly. Then the second equation of \eqref{QUACUB} may be written as $U_xx = t(x)x^2 - q(x)x$, and linearizing we may apply (\ref{ss.JABAS.fig}.\ref{UXYX}) to conclude
\begin{align}
\label{LICUB} x^2 \circ y + U_xy = U_{x,y}x + U_xy = t(x)x \circ y + t(y)x^2 - q(x,y)x - q(x)y.
\end{align}
Applying the trace (resp. $q(\emptyslot,y)$) to the first equation of \eqref{QUACUB}, we obtain
\[
t(x^2) = t(x)^2 - 2q(x), \quad q(x^2,y) = t(x)q(x,y) -q(x)t(y),
\]
while linearizing it yields
\[
x \circ y = t(x)y + t(y)x - q(x,y)e.
\]
Hence on the one hand
\begin{align*}
x^2 \circ y + U_xy =\,\,&U_xy + t(x^2)y + t(y)x^2 - q(x^2,y)e \\
=\,\,&U_xy + \big(t(x)^2 - 2q(x)\big)y + t(x)t(y)x - q(x)t(y)e - t(x)q(x,y)e + q(x)t(y)e \\
=\,\,&U_xy + \big(t(x)^2 - 2q(x)\big)y + t(x)t(y)x -t(x)q(x,y)e, 
\end{align*}
while on the other
\begin{align*}
t(x)x \circ y + t(y)x^2 - q(x,y)x -q(x)y =\,\,&t(x)^2y + t(x)t(y)x - t(x)q(x,y)e \\
\,\,&+t(x)t(y)x - q(x)t(y)e -q(x,y)x - q(x)y
\end{align*}
Comparing by means of \eqref{LICUB} and writing $x \mapsto \bar x$ for the conjugation of $(M,q,e)$, we obtain
\begin{align*}
U_xy =\,\,&q(x)y + t(x)t(y)x - q(x)t(y)e - q(x,y)x =q(x,\bar y)x - q(x)\bar y.
\end{align*}
Thus $J$ and $J(M,q,e)$ not only have the same unit element but also the same $U$-operator, which implies $J = J(M,q,e)$.
\end{sol}

\begin{sol}{pr.NIPOQUA} \label{sol.NIPOQUA}  Concerning the first part, assume that $x \in J$ is nilpotent. Then (\ref{ss.JOPOID.fig}.\ref{JOQUU}) shows that $q(x)$ is nilpotent, while Exc.~\ref{pr.PAPONI} yields a positive integer $m$ satisfying $x^n = 0$ for all integers $n \geq m$. We must show that $t(x)$ is nilpotent, and we will do so by induction on $m$. For $m = 1$, there is nothing to prove. For $m = 2$, we have $x^2 = 0$ and conclude from (\ref{ss.JOPOID.fig}.\ref{JOTQUAD}) that $t(x)^2 = 2q(x)$ is nilpotent. Hence so is $t(x)$. We may therefore assume $m \geq 3$. Since $\lfloor\frac{m}{2}\rfloor \leq \frac{m}{2} < \lfloor\frac{m}{2}\rfloor + 1$ and $m > 2$, we deduce
\begin{align*}
\lfloor\frac{m}{2}\rfloor + 1 \leq \frac{m}{2} + 1 = m + 1 - \frac{m}{2} < m + 1 - 1 = m, 
\end{align*}
and if $n \geq \lfloor\frac{m}{2}\rfloor + 1$, then $2n \geq 2(\lfloor\frac{m}{2}\rfloor + 1) > 2\frac{m}{2} = m$. On the other hand, (\ref{ss.PAOW}.\ref{PAIT}) yields $(x^2)^n = x^{2n} = 0$, and the induction hypothesis implies that $t(x^2)$ is nilpotent. But then so is $t(x)^2  = t(x^2) + 2q(x)$ by (\ref{ss.JOPOID.fig}.\ref{JOTQUAD}), which shows that $t(x)$ is nilpotent and completes the induction. Summing up we have shown that if $x$ is nilpotent, then so are $t(x)$ and $q(x)$. Conversely, assume that $t(x)$ and $q(x)$ are nilpotent. Then so are $t(x)x$ and $q(x)1_J$, which both belong to $k[x]$. Since $J$ is power-associative by Prop.~\ref{p.JOAPOA}, Exc.~\ref{pr.PAPONI} now implies that $x^2 = t(x)x - q(x)1_J$ is nilpotent. Hence so is $x$, and the first part of the problem is solved.

To establish the second part, put
\[
N := \{x \in M \mid q(x), q(x,y) \in \Nil(k)\;\text{for all $y \in M$}\}
\]
and let $x \in N$, $y,z \in M$. Then (\ref{ss.JOPOID.fig}.\ref{JOQUU}) yields $q(U_xz) = q(x)^2q(z) \in \Nil(k)$ and, similarly, $q(U_zx) \in \Nil(k)$. Moreover, since the conjugation of $J$ leaves $q$ invariant, the definition of the $U$-operator implies 
\begin{align*}
q(U_xz + U_zx,y) =\,\,&q(x,\bar z)q(x,y) + q(z,\bar x)q(z,y) - q(x)q(\bar z,y) - q(z)q(\bar x,y) \\
=\,\,&q(x,\bar z)q(x + z,y) - q(x)q(y,\bar z) - q(x,\bar y)q(z) \in \Nil(k). 
\end{align*} Thus $U_NJ + U_JN \subseteq N$, and we have shown in view of Exc.~\ref{pr.OUTIJO}~(a) that $N \subseteq J$ is an ideal which, by the first part of this problem, is nil. This proves $N \subseteq \Nil(J)$. Conversely, let $x \in \Nil(J)$. Then $q(x)$ and $t(x)$ are nilpotent and $\Nil(J)$ contains $x \circ y = t(x)y + t(y)x - q(x,y)1_J$. Hence $t(x)y - q(x,y)1_J \in \Nil(J)$, forcing
\[
q\big(t(x)y - q(x,y1_J\big) = t(x)^2q(y) - t(x)q(x,y)t(y) + q(x,y)^2
\] 
to be nilpotent. Since, therefore, $q(x,y)^2$ is nilpotent, so is $q(x,y)$ and we have shown $x \in N$. This shows $\Nil(J) \subseteq N$ and completes the proof.
\end{sol}

\begin{sol}{pr.EVPOQUA} \label{sol.EVPOQUA}  We begin with a simple lemma, which follows immediately from (\ref{ss.JOPOID.fig}.\ref{JOQUAT}). 

\begin{lem*}
Let $(M,q,e)$ be a pointed quadratic module over any commutative ring and $J := J(M,q,e)$. If $x \in J$ satisfies $x^2 = 0$, then $x^3 = -q(x)x$. \Proofend
\end{lem*}

(a) We have $I_2 = k\bft^2 \oplus k\bft^4 \oplus k\bft^5 \oplus \cdots$. Hence, writing $x$ for the imagae of $\bft$ under the canonical map $k[\bft]^{(+)} \to J_2$, we conclude that $J_2$ is free as a $k$-module, with basis $1_{J_2},x,x^3$. In particular, $x$ and $x^3$ are linearly independent. Since $x^2 = 0$, and by the preceding lemma, therefore, $J_2$ cannot be isomorphic to the Jordan algebra of a pointed quadratic form over $k$.

(b) Let $k_0$ be a commutative ring and assume $2 = 0$ in $k_0$. Furthermore, let $k = k_0[\vep]$, $\vep^2 = 0$, be the $k_0$-algebra of dual numbers. As in Exc.~\ref{pr.NONUNNOR}, we may view $k_0$ as a $k$-algebra by means of the homomorphism $k \to k_0$, $\alpha \mapsto \bar\alpha$, given by $\bar\vep = 0$. Now let $n$ be a positive integer, write $(e_1,\dots,e_n)$ for the canonical basis of $k^n$ over $k$, put $M := k_0 \times k^n$ as a $k$-module, define $q\:M \to k$ by
\begin{align*}
q\big((\bar\alpha, \sum_{i=1}^n\alpha_ie_i\big) := \alpha^2 + \vep\alpha_1^2 + \sum_{i=2}^n\alpha_i^2 &&(\alpha,\alpha_1,\dots,\alpha_n \in k),
\end{align*}
and put $e := (1_{k_0},0)$. We claim that \emph{the map $q$ is well-defined.} Indeed, if  $\alpha,\beta \in k$ satisfy $\bar\alpha = \bar\beta$, then $\alpha = \beta + \vep\gamma$, for some $\gamma \in k$, hence $\alpha^2 = \beta^2 + \vep^2\gamma^2 = \beta^2$, as desired. It now follows trivially that $(M,q,e)$ is a pointed quadratic module over $k$ whose trace is identically zero. In particular, putting $J := J(M,q,e)$ and $x := e_1$, we have $t(x) = 0$ and $q(x) = \vep$, hence $q(x)1_J = \vep(1_{k_0},0) = (\bar\vep,0) = 0 = t(x)x$, so (\ref{ss.JOPOID.fig}.\ref{JOQUAD}) implies $x^2 = 0$. Now the preceding lemma yields $x^3 = -q(x)x = -\vep e_1 \neq 0$, and we have obtained an example of the desired kind.  
\end{sol}

\begin{sol}{pr.PLUFU} \label{sol.PLUFU}  For an appropriate commutative ring $k$, we must find examples of $R,S \in \kalg$ and a homomorphism $R^{(+)} \to S^{(+)}$ of Jordan algebras that is not a morphism in $\kalg$. To this end, let $F$ be a field of characteristic two, $V$ a vector space over $F$, $a \in V$ a non-zero element, and $\sigma\:V \times V \to F$ a non-zero alternating bilinear form having $\sigma(a,V) = \{0\}$. Since $\sigma$ is alternating, hence symmetric ($F$ having characteristic two), the multiplication $uv := \sigma(u,v)a$ for $u,v \in V$ makes $V$ a commutative $F$-algebra. Since $a$ belongs to radical of $\sigma$, we also have $(uv)w = 0 = u(vw)$ for all $u,v,w \in V$. Hence in actual fact $V$ is a commutative associative algebra over $F$, as therefore is $R$, the algebra arising from $V$ by adjoining an identity element. Also, $v^2 = \sigma(v,v)a = 0$ for all $v \in V$, from which we easily deduce that any $\vph \in \GL(R)$ fixing $1_R$ and stabilizing $V$ is an automorphism of $R^{(+)}$. On the other hand, $\vph$ will never be an automorphism of $R$ unles $a$ and $\vph(a)$ are linearly dependent.
\end{sol}

\begin{sol}{skip.etale.plus.aut} \label{sol.skip.etale.plus.aut}
Given a positive integer $n$ and writing $E := k^n$, the direct product of $n$ copies of $k$ as a $k$-algebra, we will show that every $k$-automorphism of the Jordan algebra $E^{(+)}$ is a $k$-automorphism
 of the commutative associative $k$-algebra $E$.  

We first note that orthogonal idempotents $c,d$ in $E^{(+)}$ are orthogonal idempotents in $E$. Indeed, they are clearly idempotents in $E$, while orthogonality in $E$ follows from $cd = c^2d = cdc = U_cd = 0$. Thus, complete orthogonal systems of idempotents in $E$ and in $E^{(+)}$ are the same thing.

Now let $f$ be an automorphism of $E^{(+)}$ and write $(e_1, \ldots,e_n)$ for the complete orthogonal system of idempotents in $E$ given by the unit vectors of $k^n$. By the above, $(f(e_1), \ldots, f(e_n))$ is a complete orthogonal system of idempotents in $E$ such that $E = \sum kf(e_i)$.  Extending the relations $f(e_ie_j) = \delta_{ij}f(e_i) = f(e_i)f(e_j)$ for all $i,j = 1, \ldots, n$ bilinearly to all of $E$ now shows that $f$ is an automorphism of $E$.
\end{sol}

\begin{sol}{pr.COBRPL} \label{sol.COBRPL}  If $C$ is alternative, then (\ref{ss.IDCO}.\ref{QAUOP}) shows $C^{(+)} = J(C,n_C,1_C)$. Conversely, let this be so. For $x,y \in C$ we apply (\ref{ss.BASID}.\ref{CAQUADL}), (\ref{ss.BASID}.\ref{CANOCO}), (\ref{ss.BASID}.\ref{CACONJ})  and obtain
\begin{align*}
n_C(x,\bar y)x - n_C(x)\bar y =\,\,&xyx = (x \circ y)x - (yx)x = t_C(x)yx + t_C(y)x^2 - n_C(x,y)x - (yx)x \\
=\,\,&t_C(x)yx + \big(t_C(x)t_C(y) - n_C(x,y)\big)x - t_C(y)n_C(x)1_C - (yx)x \\
=\,\,&n_C(x,\bar y)x - n_C(x)\bar y - n_C(x)y + (y\bar x)x.
\end{align*}
Comparing, we arrive at one of Kirmse's identities (\ref{ss.IDCO}.\ref{KID}), i.e., at $(y\bar x)x = n_C(x)y$. But this means that $C$ is right alternative. Since it is also flexible by hypothesis, it is, in fact, alternative.
\end{sol}

\begin{sol}{pr.POJOCAT} \label{sol.POJOCAT} Homomorphism of pointed quadratic modules are clearly homomorphisms of the corresponding Jordan algebras. Hence the assignment in question defines a functor from $k$-$\mathbf{poqua}_{\mathrm{inj}}$ to $k$-$\mathbf{jord}_{\mathrm{inj}}$, which by its very definition is faithful Before we can show that it is full, we need the following lemma. 

\begin{lem*}
Let $(M,q,e)$ and $(M,q_1,e)$ be pointed quadratic modules over $k$, with the same underlying $k$-module $M$ and the same base point $e$. If $M$ is projective and $J(M,q,e) = J(M,q_1,e)$, then $q = q_1$.
\end{lem*}

\begin{proof}
Since $e \in M$ is unimodular by Lemma~\ref{l.UNBA}, the proof of Prop.~\ref{p.UNINOR} can be extended verbatim to this slightly different setting. 
\end{proof}

Now let $(M,q,e)$ and $(M^\prime,q^\prime,e^\prime)$ be pointed quadratic modules over $k$ such that $M$ and $M^\prime$ are both projective and let $\vph\:J(M,q,e) \to J(M^\prime,q^\prime,e^\prime)$ be an injective homomorphism of Jordan algebras over $k$. Setting $q_1 := q^\prime \circ \vph$, we conclude that $(M,q_1,e)$ is a pointed quadratic module over $k$ (since $\vph$ preserves identity elements) and $\vph\:(M,q_1,e) \to (M^\prime,q^\prime,e^\prime)$ is a homomorphism of pointed quadratic modules. Hence $\vph\:J(M,q_1,e) \to J(M^\prime,q^\prime,e^\prime)$ is a homomorphism of Jordan algebras. But $\vph$ was assumed to be injective. Hence $J(M,q,e) = J(M,q_1,e)$, which by the preceding lemma implies $q_1 = q$. Hence $\vph\:(M,q,e) \to (M^\prime,q^\prime,e^\prime)$ is a homomorphism of pointed quadratic modules. Thus the functor in question is full. 
\end{sol}

\begin{sol}{pr.POQUAIN} \label{sol.POQUAIN}  We put $J := J(V,q,e)$. Since $F$ has characteristic not $2$, we may extend $e = 1_J$ to an orthogonal basis of $V$ relative to $q$, allowing us to identify $q = \la 1,\delta,\delta\vep\raq$ for some $\delta,\vep \in F^\times$. Now put
\[
K := \Cay(F,-\vep) = F \oplus Fj_1, \quad B := \Cay(K,-\delta) = K \oplus Kj_2.
\]
as in \ref{ss.EXTCD}. Then $B$ is a quaternion algebra over $F$, while $K \subseteq B$ is a quadratic \'etale subalgebra. Let $\tau$ be the involution of $B$ corresponding to $K$ by Exc.~\ref{pr.RINV}~(b). Then $\tau(a \oplus bj_2) = \bar a \oplus bj_2$ for all $a,b \in K$, which implies $H(B,\tau) = W := F1_B \oplus Kj_2 = F1_B \oplus Fj_2 \oplus Fj_1j_2$. By Exc.~\ref{pr.COBRPL} we have $B^{(+)} = J(B,n_B,1_B)$ and therefore $H(B,\tau) = J(W,n_B\vert_W,1_B)$. By standard properties of the Cayley-Dickson construction, the vectors $1_B,j_2,j_1j_2$ are orthogonal relative to $n_B $ and satisfy $n_B(1_B) = 1$, $n_B(j_2) = \delta$, $n_B(j_1j_2) = n_B(j_1)n_B(j_2) = \delta\vep$. Hence the pointed quadratic modules $(V,q,e)$ and $(W,N_B\vert_W,1_B)$ are isometric, forcing $J$ and $H(B,\tau)$ as the corresponding Jordan algebras by Exc.~\ref{pr.POJOCAT} to be isomorphic.

It remains to prove that $(B,\tau)$ as a quaternion algebra with involution is unique up to isomorphism. To this end, we require a lemma. 

\begin{lem*}
Let $B$, $B^\prime$ be quaternion algebras over $F$ and suppose $W \subseteq (B,n_B)$, $W^\prime \subseteq (B^\prime,n_{B^\prime})$ are regular subspaces of co-dimension $1$. If $W$ and $W^\prime$ are isometric, then $B$ and $B^\prime$ are isomorphic.
\end{lem*}

\begin{proof}
By Lemma~\ref{l.NOSU} we have $B =  W \perp Fu$, $B^\prime = W^\prime \perp Fu^\prime$ for some non-zero $u \in B$, $u^\prime \in B^\prime$. Recall that the determinant of a quadratic space like $(B,n_B)$ (resp. $(B^\prime,n_{B^\prime})$) is unique up to non-zero square factors in $F$. As a matter of fact, we deduce from (\ref{ss.QUATALG}.\ref{QUATNOR}) that this determinant is a square. On the other hand, it agrees with $\det(W)n_B(u)$ (resp. $\det(W^\prime)n_{B^\prime}(u^\prime)$). But $W$ and $W^\prime$ being isometric, their determinants differ by a square factor in $F$. Hence so do $n_B(u)$ and $n_{B^\prime}(u^\prime)$. But this means that $B$ and $B^\prime$ are norm-equivalent. By the norm equivalence theorem~\ref{t.NOREQ}, therefore, they are isomorphic. 
\end{proof}

Now let $(B^\prime,\tau^\prime)$ be another quaternion algebra with involution over $F$ having $W^\prime := H(B^\prime,\tau^\prime) \cong J(V,q,e)$. By Exc.~\ref{pr.POJOCAT}, the quadratic subspaces $W \subseteq (B,n_B)$ and $W^\prime \subseteq (B^\prime,n_{B^\prime})$ are not only regular but also isometric. Our preceding lemma therefore implies that $B$ and $B^\prime$ are isomorphic. We may thus assume $B = B^\prime$. Let $(1_B, j_2^\prime,j_3^\prime)$  be an orthogonal basis of $W^\prime$ such that $n_B(j_2^\prime) = \delta$, $n_B(j_3^\prime) = \delta\vep$, with $\delta,\vep$ as above. Since $(B,n_B)$ has trivial determinant, we conclude $W^{\prime\perp} = Fj_1^\prime$, $n_B(j_1^\prime) = \vep$. Then $F[j_1^\prime] \cong \Cay(F,-\vep) = K$. Now (\ref{ss.BASID}.\ref{CAQUADL}) implies for $i = 2,3$:
\begin{align*}
0 =\,\,&\tau^\prime\big(t_B(j_1^\prime)j_i^\prime + t_B(j_i^\prime)j_1^\prime - n_B(j_1^\prime,j_i^\prime)1_B\big) = \tau^\prime(j_1^\prime \circ j_i^\prime) = \tau^\prime(j_1^\prime) \circ j_i^\prime \\
=\,\,&t_B(j_1^\prime)j_i^\prime + t_B(j_i^\prime)\tau^\prime(j_1^\prime) - n_B\big(\tau^\prime(j_1^\prime),j_i^\prime\big)1_B = -n_B\big(\tau^\prime(j_1^\prime),j_i^\prime\big)1_B.
\end{align*}
Thus $n_B(\tau^\prime(j_1^\prime),j_i^\prime) = 0$, and since we trivially have $n_B(\tau^\prime(j_1^\prime),1_B) = 0$, we deduce $\tau^\prime(j_1^\prime) \in W^{\prime\perp} = Fj_1^\prime$, Hence $\tau^\prime(j_1^\prime) = \pm j_1^\prime$. Here $\tau^\prime(j_1^\prime) = j_1^\prime$ would imply $\tau^\prime = \Eins_B$, a contradiction. Hence $\tau^\prime(j_1^\prime) = -j_1^\prime$. Thus  Thus $\tau^\prime$ is the involution of $B$ corresponding to $F[j_1^\prime] \cong K$, which implies $(B,\tau^\prime) \cong (B,\tau)$ by Exc.~\ref{pr.RINV}~(b).
\end{sol}

\solnsec{Section~\ref{s.POID}}

\begin{sol}{pr.ABZECLI} \label{sol.ABZECLI}  (a) $x \in \Rad(q)$ implies $U_xy = q(x,\bar y)x - q(x)\bar y = 0$ for all $y \in J$, hence $U_x = 0$, and $x$ is an absolute zero divisor of $J$. Conversely, let this be so. Then (\ref{ss.JOPOID.fig}.\ref{JOCONJC}) yields $q(x)^2e = U_x\bar x^2 = 0$, and applying $q$ we conclude $q(x)^4 = 0$, hence $q(x) = 0$ since $k$ is reduced. For any $y \in J$, we therefore obtain $0 = U_x\bar y
= q(x,y)x$ and then $q(x,y)^2 = q(q(x,y)x,y) = 0$. Thus $q(x,y) = 0$, and we have proved $x \in \Rad(q)$.

(b) Here $k$ is arbitrary. Then $\bar k := k/\Nil(k) \in \kalg$ is reduced, and we have a canonical identification $J_{\bar k} = \bar J := J/\Nil(k)J$ via \ref{ss.REDID}, matching $z_{\bar k}$ for $z \in J$ with $\bar z$, the image of $z$ under the natural map $J \to \bar J $. Since $x$ is an absolute zero divisor in $J$, $\bar x$ is one in $\bar J$, so by the special case just treated we have $\bar q(\bar x) = \bar q(\bar x,\bar y) = 0$ for all $y \in J$, where $\bar q = q_{\bar k}$ is the $\bar k$-quadratic extension of $q$. But this means that $q(x)$ and $q(x,y)$ are nilpotent elements of $k$, for all $y \in J$, which implies $x \in \Nil(J)$ by Exc.~\ref{pr.NIPOQUA}. The rest is clear, by Thm.~\ref{t.AZLOC}.
\end{sol}

\begin{sol}{pr.DICKJO} \label{sol.DICKJO} 
If $J = J(M,q,e)$ for some pointed quadratic module $(M,q,e)$ over $F$, then the Dickson condition clearly holds. Conversely, suppose $J$ satisfies the Dickson condition. Then we are forced to put $M := J$ as a vector space over $F$ and $e := 1_J$. The problem is to construct the quadratic form $q\:M \to F$. We do so by considering the following cases. 

\case{1}
$\ch(F) \neq 2$. Then $J$ becomes a linear Jordan algebra under the product $x\cdot y := \frac{1}{2}x \circ y$, and by hypothesis, the quantities $1_J,x,x^2$ are linearly dependent over $F$. Hence Exc.~\ref{pr.CONFIELD} implies that $J$ is a commutative conic algebra over $F$, with norm $q$ and trace $t$, say. Linearizing the relation $x^2 = t(x)x -q(x)e$ and invoking (\ref{ss.LIJPOQ}.\ref{LIMPOQ}), we conclude $J = J(M,q,e)$ as linear Jordan algebras, hence also as quadratic ones. 

\case{2}
$\ch(F) = 2$. We may clearly assume that $J$ has dimension at least $2$. 

\case{2.1}
$\dim_F(J) = 2$. Let $e = 1_J,u$ be a basis of $J$ over $F$. Then
\begin{align}
\label{USQ} u^2 = \gamma u + \delta e
\end{align}
for some $\gamma,\delta \in F$, and we define a quadratic form
\begin{align}
\label{QUEMF} q(\xi_0e + \xi_1u) := \xi_0^2 + \gamma\xi_0\xi_1 + \delta\xi_1^2 
\end{align}
for $\xi_0,\xi_1 \in F$. One checks that $(M,q,e)$ is a pointed quadratic module over $F$, with bilinearized norm, trace and conjugation respectively given by
\begin{align}
\label{QUXEZ} q(\xi_0e + \xi_1u,\eta_0e + \eta_1u) =\,\,&\gamma(\xi_0\eta_1 + \xi_1\eta_0), \\
\label{TEXEZ} t(\xi_0e + \xi_1u) =\,\,&\gamma\xi_1, \\
\label{CONXEZ} \overline{\xi_0e + \xi_1u} =\,\,&(\xi_0 + \gamma\xi_1)e + \xi_1u
\end{align}
for $\xi_i,\eta_i \in F$, $i = 0,1$. Moreover, since $J$ is locally linear, we have $u^3 = uu^2 =\gamma u^2 + \delta u$, hence
\begin{align}
\label{UCUB} u^3 = \gamma\delta e + (\gamma^2 + \delta)u.
\end{align}
Using \eqref{USQ}--\eqref{UCUB}, it is now straightforward to check that the $U$-operators of $J$ and $J(M,q,e)$ are the same. Hence $J = J(M,q,e)$. 

\case{2.2}
$\dim_F(J) \geq 3$. Choose a subspace $V \subseteq M$ that is complementary to $Fe$:
\begin{align}
\label{EMEF} M = Fe \oplus V.
\end{align}
Since $J$ satisfies the Dickson condition, we obtain
\begin{align}
\label{VESQ} v^2 = q(v)e + s(v), \quad q(v) \in F, \quad s(v) \in V
\end{align}
strictly for all $v \in V$ such that there eists a map $t\:V \to F$ satisfying $t(0) = 0$ and 
\begin{align}
\label{ESOVE} s(v) = t(v)v, 
\end{align}
again strictly for all $v \in V$ . Here $q\:V \to F$ is a quadratic form, $s\:V \to V$ is a quadratic map, and $t\:V \to F$ is rational as well as homogeneous of degree $1$. We wish to show that $t$ is linear. To this end, we may assume that $F$ is infinite, replace $v$ by $v + w$ in \eqref{VESQ}, \eqref{ESOVE} and collect mixed terms to obtain
\begin{align}
\label{VECIR} v \circ w =\,\,&\big(t(v + w) - t(v)\big)v + \big(t(v + w) - t(w)\big)w - q(v,w)e.
\end{align}                     
Since $F$ has characteristic $2$, this implies
\begin{align}
\label{VECIRVE} v \circ (v \circ w) =\,\,&\big(t(v + w) - t(w)\big)v \circ w, \\
\label{VESQCIR} v^2 \circ w =\,\,&t(v)v \circ w.
\end{align}
By the same token and (\ref{ss.JABAS.fig}.\ref{TWOUX}), we have $v \circ (v \circ w) = v^2 \circ w$, hence
\begin{align}
\label{TEVE} \big(t(v + w) - t(v) - t(w)\big)v \circ w = 0.
\end{align}
Now, if $v \circ w = 0$ for all $v,w \in V$, then \eqref{VECIR} implies $t(v + w) = t(v) = t(w)$ provided $v,w$ are linearly independent. By Zariski density, therefore, $t$ is constant, hence zero, hence linear. On the other hand, if $v \circ w \neq 0$ for some $v,w \in V$, Zariski density again and \eqref{TEVE} imply that $t$ is linear. Thus $t$ is linear under all circumstances, and we have
\[
v^2 = t(v)v - q(v)e
\] 
for all $v \in V$. Extending $q,t$ to all of $M$ via
\[
q(\xi e + v) = \xi^2 + \xi t(v) + q(v), \quad t(\xi e + v) = t(v)
\]
for $\xi \in k$, $v \in V$, one checks that $x^2 - t(x)x + q(x)e = 0$ for all $x \in M$. But $J$ is strictly locally linear, forcing $x^3 - t(x)x^2 + q(x)x = 0$ strictly for all $x \in M$. Now $J = J(M,q,e)$ follows from Exc.~\ref{pr.CHAJOCLI}. 
\end{sol}


\solnsec{Section~\ref{s.INIS}}

\begin{sol}{pr.INVUNI} \label{sol.INUNI} (i) $\Rightarrow$ (ii). This implication follows from the fact that $J$ contains invertible elements (e.g., the identity).

(ii) $\Rightarrow$ (iii) Let $x \in J$ be unimodular. For $\mfp \in \Spec(k)$, Lemma~\ref{l.FAITH} implies $x(\mfp) \ne 0$, hence $J(\mfp) \ne \{0\}$, and we conclude $\rk_\mfp(J) = \dim_{k(\mfp)}(J(\mfp)) > 0$.

(iii) $\Rightarrow$ (i). For $x \in J^\times$, the linear map $U_x\:J \to J$ is bijective (Prop.~\ref{p.CHAIN.JOAL}), and hence so is $U_{x(\mfp)} = U_x{(\mfp)}$. But from (iii) we deduce $J(\mfp) \ne \{0\}$, which implies $x(\mfp) \ne 0$, and Lemma~\ref{l.FAITH} show that $x$ is unimodular. 
\end{sol}

\begin{sol}{pr.EVINV} \label{sol.EVINV}  (a) By the results of \S\S\ref{s.JABAS},\ref{s.POID}, equations \eqref{NEUXN}--\eqref{NEUUXIJMN} hold for all $i,j,m,n,p \in \IN$. For $n \in \IZ$, $n > 0$, we combine \eqref{NEPO} with (\ref{p.PROIN}.\ref{FININ}) and obtain $U_{x^{-n}} = U_{(x^{-1})^n} = U_{x^{-1}}^n = (U_x^{-1})^n = U_x^{-n}$. Hence \eqref{NEUXN} holds for all $n \in \IZ$.

By \eqref{NEUXN}, the set of all integers $m$ such that \eqref{NEUXMN} holds for all integers $n$ is a subgroup of $\IZ$. Hence, in order to prove \eqref{NEUXMN}, we may assume $m = 1$. Then we must show $U_xx^{-n} = x^{2-n}$ for all $n \in \IN$ and we do so by induction on $n$. For $n = 0,1$, there is nothing to prove. For $n \geq 2$, assuming the validity of the assertion for $n - 2$ yields $U_xx^{-n} = U_x(x^{-1})^n = U_xU_{x^{-1}}(x^{-1})^{n-2} = x^{2-n}$, as claimed. Thus \eqref{NEUXMN} holds.

Turning to \eqref{NEIT}, we first assume $n \in \IN$ and then argue by induction. For $n = 0,1$ there is nothing to prove. Now let $n \geq 2$ and assume the assertion is valid for $n - 2$. Then $(x^m)^n = U_{x^m}(x^m)^{n-2} = U_{x^m}x^{mn-2m} = x^{2m+mn-2m} = x^{mn}$. as claimed. Next assume $n < 0$. The case just treated and \eqref{NEUXN}, \eqref{NEUXMN} yield $U_{(x^m)^{-n}}(x^m)^n = (x^m)^{-n} = x^{-mn} = U_{x^{-mn}}x^{mn} = U_{(x^m)^{-n}}x^{mn}$, and since $U_{(x^m)^{-n}}$ is bijective, the proof of \eqref{NEIT} is complete.

Next we establish \eqref{NEUUXMN}, where we may clearly assume $i = 1$. For $l \in \IN$ such that $2l + m, 2l + n \in \IN$, we obtain, using \eqref{NEUXN}, (\ref{ss.JABAS.fig}.\ref{UUXYZU}), \eqref{NEUXMN}, (\ref{p.LIPO}.\ref{PROU}),
\begin{align*}
U_xU_{x^m,x^n} =\,\,&U_{x^{-l}}U_xU_{x^l}U_{x^m,x^n}U_{x^l}U_{x^{-l}} = U_{x^{-l}}U_xU_{U_{x^l}x^m,U_{x^l}x^n}U_{x^{-l}} = U_{x^{-l}}U_xU_{x^{2l+m},x^{2l+n}}U_{x^{-l}} \\
=\,\,&U_{x^{-l}}U_{x^{2l+m+1},x^{2l+n+1}}U_{x^{-l}} = U_{U_{x^{-l}}x^{2l+m+1},U_{x^{-l}}x^{2l+n+1}} = U_{x^{m+1},x^{n+1}}
\end{align*}
and, similarly, $U_{x^m,x^n}U_x = U_{x^{m+1},x^{n+1}}$.Thus \eqref{NEUUXMN} holds.

Now we can prove \eqref{NETRI} by letting $l \in \IN$ satisfy $2l + m, 2l + n, 2l + p \in \IN$. Then \eqref{NEUUXMN} yields
\begin{align*}
U_{x^{3l}}\{x^mx^nx^p\} =\,\,&U_{x^{3l}}U_{x^m,x^p}x^n = U_{x^l}U_{x^m,x^p}U_{x^l}U_{x^l}x^n = U_{U_{x^l}x^m,U_{x^l}x^p}U_{x^l}x^n =\{x^{2l+m}x^{2l+n}x^{2l+p}\} \\
=\,\,&2x^{6l+m+n+p} = U_{x^{3l}}(2x^{m+n+p}).
\end{align*} 
Thus \eqref{NETRI} holds.

Turning to \eqref{NEVXMN}, we use \eqref{NEUXN}, (\ref{ss.JABAS.fig}.\ref{UXVY}), \eqref{NEUXMN}, \eqref{NEUUXMN} and obtain 
\begin{align*}
V_{x^m,x^n}U_{x^{m+n}} =\,\,&V_{x^m,x^n}U_{x^m}U_{x^n} = U_{x^m,U_{x^m}x^n}U_{x^n} = U_{x^m,x^{2m+n}}U_{x^n} = U_{x^{m+n},x^{2(m+n)}} \\
=\,\,&U_{x^{m+n},U_{x^{m+n}}1_J} = V_{x^{m+n},1_J}U_{x^{m+n}} = V_{x^{m+n}}U_{x^{m+n}}.
\end{align*}
Hence \eqref{NEVXMN} follows.

Finally,in order to establish \eqref{NEUUXIJMN}, we let $l \in \IN$ satisfy $l + i, l + j, l + m, l + n \in \IN$. Then, by \eqref{NEUXN}, \eqref{NEUUXMN}, (\ref{p.LIPO}.\ref{PROUL}), \ref{ss.JABAS.fig}.\ref{UUXYZU}),
\begin{align*}
U_{x^i,x^j}U_{x^m,x^n} =\,\,&U_{x^{-l}}U_{x^l}U_{x^i,x^j}U_{x^l}U_{x^m,x^n}U_{x^{-l}} = U_{x^{-l}}U_{x^{l+i},x^{l+j}}U_{x^{l+m},x^{l+n}}U_{x^{-l}} \\
=\,\,&U_{x^{-l}}\Big(U_{x^{2l+i+m},x^{2l+j+n}} + U_{x^{2l+i+n},x^{2l+j+m}}\Big)U_{x^{-l}} \\
=\,\,&U_{x^{i+m},x^{j+n}} + U_{x^{i+n},x^{j+m}},
\end{align*} 
as claimed.

(b) Arguing as in the solution to Exc.~\ref{pr.PARNIL}~(a), it follows that $\vep_x^\times\:k[\bft,\bft^{-1}]^{(+)} \to J$ is a homomorphism of Jordan algebras forcing $k[x,x^{-1}] \subseteq J$ to be a subalgebra containing $\vep_x^\times (\bft^{\pm 1}) = x^{\pm 1}$. Writing $J^\prime$ for the subalgebra of $J$ generated by $x$ and $x^{-1}$, we therefore have $J^\prime \subseteq k[x,x^{-1}]$. On the other hand, it follows immediately from the definitions that $J^\prime$ contains all powers of $x$ with integer exponents. Since, therefore, $k[x,x^{-1}] \subseteq J^\prime$, we have equality: $J^\prime = k[x,x^{-1}]$. 

(c) Arguing as in \ref{ss.FIRE}, replacing $\IN$ by $\IZ$ everywhere and using \eqref{NEUUXMN}, \eqref{NEUUXIJMN}, we deduce \eqref{NEUFG}. The final assertion may be proved in exactly the same manner as Thm.~\ref{t.AZLOC}.
\end{sol}

\begin{sol}{pr.GEOSER} \label{sol.GEOSER}  Put $x := 1_J - u$ and $y := \sum_{n\geq 0}u^n$. Then, since $J$ is power-associative,
\begin{align*}
U_xy =\,\,&y - u \circ y + U_uy = \sum_{n\geq 0}u^n - 2\sum_{n\geq 0}u^{n+1} + \sum_{n\geq 0}u^{n+2} \\
=\,\,&\sum_{n\geq 0}u^n - 2\sum_{n\geq 1}u^n + \sum_{n\geq 2}u^n = 1_J - u = x,
\end{align*}
while Thm.~\ref{t.UPOL} yields
\begin{align*}
U_xy^2 =\,\,&U_{1_J - u}U_{\sum_{n\geq 0}u^n}1_J = U_{\sum_{n\geq 0}u^n-\sum_{n\geq 1}u^n}1_J = U_{1_J}1_J = 1_J,
\end{align*}
solving the first part of the problem. As to the second, if $x$ is invertible in $J$, then by Prop.~\ref{p.HOIN} so is $\bar x$ in $\bar J$. Conversely, let this be so. Then Prop.~\ref{p.CHAIN.JOAL} implies that some $y \in J$ satisfies $U_xy = 1_J -u$, where $u$ belongs to $I$ and hence is nilpotent. By the first part, therefore, $U_xy$ is invertible in $J$ and thus so is $x$, by Prop.~\ref{p.PROIN}.
\end{sol}

\begin{sol}{pr.INVLIN} \label{sol.INVLIN}  (a) By (\ref{ss.ELEIDE}.\ref{ONELE}) we have
\[
[L_{x^2},L_y] = 0 \quad \Longleftrightarrow \quad [L_x,L_{xy}] = 0.
\]

(b) (i) $\Rightarrow$ (ii). Let $y := x^{-1}$. Then $2xy = x \circ y = \{x1_Jx^{-1}\} = 2\cdot 1_J$ by Exc.~\ref{pr.EVINV}, \eqref{NETRI}. Hence $xy = 1_J$. Moreover, Exc.~\ref{pr.EVINV}, \eqref{NEUUXIJMN} shows $4[L_x,L_y] = [U_{x,1_J},U_{1_j,y}] = 0$, so $x$ and $y$ operator commute.

(ii) $\Rightarrow$ (iii). We have $xy = 1_J$, and since $x$ and $y$ operator commute, we also obtain $x^2y = y(xx) = x(yx) = x(xy) = x$,

(iii) $\Rightarrow$ (i). Since $x$ and $xy = 1_J$ operator commute, so do $x^2$ and $y$, by (a). For the same reason $y^2$ and $x$ operator commute. Hence
\[
U_xy^2 =2x(y^2x) - x^2y^2 = x^2(yy) = y(x^2y) = yx = 1_J,
\]
and $x$ is invertible in $J$ by Prop.~\ref{p.CHAIN.JOAL}. 

Next suppose $y \in J$ satisfies (ii). Then $x$ is invertible and $U_xy = 2x(xy) - y(xx) = 2x - x(yx) = x$, hence $y = U_x^{-1}x = x^{-1}$. Similarly, if $y$ satisfies (iii), then $U_xy = 2x(xy) -x^2y = x$ and again $y = x^{-1}$.
\end{sol}

\begin{sol}{pr.INVPOI} \label{sol.INVPOI}  If $x$ is invertible in $J$, then $1 = q(1_J) = q(U_xx^{-2}) = q(x)^2q(x^{-2})$ by (\ref{ss.JOPOID.fig}.\ref{JOQUU}), so $q(x)$ is invertible in $k$. Conversely, let this be so and put $y := q(x)^{-1}\bar x$. Then (\ref{ss.JOPOID.fig}.\ref{JOCONJC}) yields $U_xy = q(x)^{-1}U_x\bar x = x$. On the other hand, writing $t$ for the trace of $(M,q,e)$, we apply (\ref{ss.JOPOID.fig}.\ref{JOCONJO}), (\ref{ss.JOPOID.fig}.\ref{JOQUAD}) and obtain not only $q(y) = q(x)^{-1}$ but also
\[
U_xy^2 = t(y)U_xy -q(y)U_x1_J = q(x)^{-1}t(x)x - q(x)^{-1}x^2 = 1_J.
\] 
Hence (\ref{ss.COIN}.\ref{ININ}) shows that $x$ is invertible in $J$ with inverse $y$. The rest is clear.
\end{sol}

\begin{sol}{pr.ISTPOI} \label{sol.ISTPOI}  (a) $q^{(f)}\:M \to k$ is a quadratic form and $e^{(f)} \in M$ satisfies $q^{(f)}(e^{(f)}) = q(f)q(f^{-1}) = 1$ by Exc.~\ref{pr.INVPOI}, \eqref{INVPOI}. Hence $(M,q,e)^{(f)}$ is a pointed quadratic module over $k$. For $x \in M$, the definitions imply
\begin{align*}
t^{(f)}(x) =\,\,&q^{(f)}(e^{(f)},x) = q(f)q\big(q(f)^{-1}\bar f,x\big) = q(\bar f,x), \\
\bar x^{(f)} =\,\,&t^{(f)}(x)e^{(f)} - x = q(f)^{-1}q(\bar f,x)\bar f - x,
\end{align*}
hence \eqref{TEF} and \eqref{CONEF}. In \eqref{JEF}, both sides live on the same $k$-module and have the same identity element. Hence it suffices to show that they have the same $U$-operator. Let $x,y \in M$ and write $U_x$ (resp. $U^\prime_x$) for the $U$-operator of $J(M,q,e)$ (resp. $J((M,q,e)^{(f)}$). We first note
\begin{align}
\label{UEFY} \overline{U_fy} = U_{\bar f}\bar y = q(\bar f,y)\bar f - q(f)y
\end{align} 
by (\ref{ss.JOPOID.fig}.\ref{JOCONJ}), (\ref{ss.JOPOID.fig}.\ref{JOQUADU}) and (\ref{ss.JOPOID.fig}.\ref{JOCONJO}). Hence \eqref{TEF} and \eqref{CONEF} imply
\begin{align*}
U^\prime_xy =\,\,&q^{(f)}(x,\bar y^{(f)})x - q^{(f)}(x)\bar y^{(f)} \\
=\,\,&q(f)q\big(x,q(f)^{-1}q(\bar f,y)\bar f  - y\big)x - q(f)q(x)\big(q(f)^{-1}q(\bar f,y)\bar f - y\big) \\
=\,\,&\big(q(\bar f,x)q(\bar f,y) - q(f)q(x,y)\big)x - q(x)\big(q(\bar f,y)\bar f - q(f)y\big) \\
=\,\,&q\big(x,q(\bar f,y)\bar f - q(f)y\big)x - q(x)\big(q(\bar f,y)\bar f - q(f)y\big) \\
=\,\,&q\big(x,\overline{U_fy}\big)x - q(x)\overline{U_fy} = U_xU_fy = U^{(f)}_xy,
\end{align*}
hence \eqref{JEF}.

(b) From (a) we conclude that $g \in M$ is invertible in $(M,q,e)^{(f)}$ iff $k^\times \ni q^{(f)}(g)  = q(f)q(g)$ iff $q(g) \in k^\times$ iff $g$ is invertible in $(M,q,e)$. In this case, $((e^{(f)})^{(g)} = g^{(-1,f)}$ (by (a)) $= U_{f^{-1}}g^{-1}$ (by Thm.~\ref{t.ISTJIT}~(a)) $= (U_fg)^{-1}$ (by (\ref{p.PROIN}.\ref{UIN})) $= e^{(U_fg)}$ and $(q^{(f)})^{(g)}(x) = q^{(f)}(g)q^{(f)}(x) = q(f)^2q(g)q(x) = q(U_fg)q(x)$ (by (\ref{ss.JOPOID.fig}.\ref{JOQUU})) $= q^{(U_fg)}(x)$ for all $x \in M$. This proves (b).

(c) $\eta \in \Str(J)$ is an isomorphism $\eta\:J \overset{\sim} \to J^{(p)}$, $J^{(p)} = J(M,e^{(p)},q^{(p)})$ by (a), $p := \eta(1_J) \in J^\times$ (by (\ref{ss.HOMOT}.\ref{HOMOTONE})). By Exc.~\ref{pr.POJOCAT}, therefore, $\eta\:(M,q) \overset{\sim} \to (M,q^{(p)})$ is an isometry of quadratic modules, which implies $q \circ \eta = \alpha q$, $\alpha :=q(p)^{-1}$. Conversely, let $\eta\:M \to M$ be a linear bijection such that $q \circ \eta = \alpha q$, for some $\alpha \in k^\times$. Then $\alpha = q(\eta(e))$, so $\eta(e)$ is invertible in $(M,q,e)$ and $q^{(p)} \circ \eta = q$, $p := \eta(e)^{-1}$. Thus $\eta\:(M,q,e) \to (M,q,e)^{(p)}$ is an isomorphism of pointed quadratic modules and hence, again by Exc.~\ref{pr.POJOCAT}, and isomorphism $\eta\:J \to J^{(p)}$ of Jordan algebras. Summing up, $\eta \in \Str(J)$ and we are done.

(d) With $p := e_3$ we have $q(p) = -1$ and hence $q^{(p)} = -q$. The signature of $q$ is $-1$, while the signature of $q^{(p)}$ is $1$. Since, therefore, $q$ and $q^{(p)}$ are not isometric, $J$ and $J^{(p)}$ again by (a) and Exc.~\ref{pr.POJOCAT} are not isomorphic.

(e) The linear trace, $t^{(f)}$, of $J^{(f)} = J((M,q,e)^{(f)})$ by (a) satisfies $t^{(f)}(x) = q(f,\bar x)$, hence is different from zero. It follows that $(M,q,e)$ and $(M,q,e)^{(f)}$ are not isomorphic. By Exc.~\ref{pr.POJOCAT}, therefore, neither are $J$ and $J^{(f)}$. As an example put $M := F^3$ with the canonical basis $e_1,e_2,e_3$, $e := e_1$ and $q := \la S \raq$, where
\[
S = \left(\begin{matrix}
1 & 0 & 0 \\0 & 1 & 1 \\
0 & 0 & 0
\end{matrix}\right).
\]
Then $(M,q,e)$ is a pointed quadratic module over $F$ and $Dq = \la S + S^\trans \ra$, where
\[
S + S^\trans = \left(\begin{matrix}
    0 & 0 & 0 \\0 & 0 & 1 \\0 & 1 & 0
\end{matrix}\right)
\]
Thus $t = q(e,\emptyslot) = 0$ and $f := e_2$ satisfies $q(f) = q(f,e_3) = 1$. Hence all hypotheses of (e) are fulfilled, irrespective of whether $F$ is or is not algebraically closed.
\end{sol}

\begin{sol}{pr.FOUROP} \label{sol.FOUROP}  Setting $1^{(y)} := 1_{J^{(y)}}$ as in (\ref{ss.JOISOT}.\ref{UNIP}), we apply (\ref{ss.JOISOT}.\ref{UOPP}) and obtain $U_{x^{-1} + y^{-1}}U_y = U^{(y)}_{1^{(y)} + z}$, $z := x^{-1}$. Here Thm.~\ref{t.ISTJIT}~(a) yields $z^{(-1,y)} = U_y^{-1}(x^{-1})^{-1} = = U_y^{-1}x$, hence $U_yz^{(-1,y)} = x$. Using Exc.~\ref{pr.EVINV}, \eqref{NEUFG}, we therefore conclude
\begin{align*}
U_xU_{x^{-1} + y^{-1}}U_y =\,\,&U_{U_yz^{(-1,y)}}U_{1^{(y)} + z}^{(y)} = U_y U_{z^{(-1,y)}}^{(y)}U_{1^{(y)} + z}^{(y)} = U_yU_{z^{(-1,y)} + 1^{(y)}}^ {(y)} = U_y U_{z^{(-1,y)} + y^{-1}}U_y \\
=\,\,&U_{U_yz^ {(-1,y)} + U_yy^ {-1}} = U_{x + y}.
\end{align*}
\end{sol}

\begin{sol}{pr.LININV} \label{sol.LININV}  We have $xx^{-1} = L_xL_x^{-1}1_J = 1_J$, and the Jordan identity (\ref{ss.LIJO}.\ref{LIJO}) yields $x(x^2x^{-1}) = x^2(xx^{-1}) = x^2$, hence $x^2x^{-1} = x$ since $L_x$ is bijective. We therefore conclude from Exc.~\ref{pr.INVLIN} that $x$ is invertible in $J$ with inverse $x^{-1}$. Now (\ref{p.PROIN}.\ref{FININ}) and (\ref{ss.UOPLIJ}.\ref{VELELI}) imply that $L_{x^{-1}} = \frac{1}{2}V_{x^{-1}} = \frac{1}{2}U_x^{-1}V_x = U_x^{-1}L_x$ is bijective, forcing $x^{-1}$ to be linearly invertible with linear (= ordinary) inverse $(x^{-1})^{-1} = x$.

Now let $F$ be a field of characteristic not $2$ and $(M,q,e)$ a pointed quadratic module over $F$. If $q$ is anisotropic and $J := J(M,q,e)$ has dimension at most $2$, then Exc.~\ref{pr.CONIL} shows that $J$ is a field and, in particular, a linear division algebra. Conversely, let $J$ be a linear division algebra. By what we have just seen, $J$ is a Jordan division algebra, forcing $q$ to be anisotropic by Exc.~\ref{pr.INVPOI}. Assume $J$ has dimension $n \geq 3$ and extend $1_J$ to an orthogonal basis $(e_1 = 1_J,e_2,e_3,\dots,e_n)$ of $J$ relative to $Dq$. Applying (\ref{ss.JOPOID.fig}.\ref{JOQUADL}), we obtain
\[
e_2e_3 = \frac{1}{2}\big(q(e_1,e_2)e_3 + q(e_1,e_3)e_2 - q(e_2,e_3)e_1\big) = 0.
\]
Hence $e_2$, though invertible, is not linearly invertible, and $J$ cannot be a linear division algebra. This contradiction shows $n \leq 2$ and completes the proof.
\end{sol}

\begin{sol}{pr.STRHOT} \label{sol.STRHOT}  (a) Let $\eta\:J \to J^\prime$ be a strong homotopy, so some $p \in J^\times$ makes $\eta\:J^{(p)} \to J^\prime$ a homomorphism. Hence $\eta\:J \to J^\prime$ is a homotopy by Prop.~\ref{p.HOMOTIST}~(b). Conversely, let $\eta\:J \to J^\prime$ be a homotopy such that $\eta(J^\times) = J^{\prime\times}$. Then there exists a $p^\prime \in J^{\prime\times}$ such that $\eta\:J \to J^{\prime(p^\prime)}$ is a homomorphism. By hypothesis, we find an element $p \in J^\times$ satisfying $\eta(p) = p^{\prime -2}$, Then \ref{ss.JOISOT}~(ii) and Cor.~\ref{c.ISTJYM} show that
\[
\eta\:J^{(p)} \longrightarrow J^{\prime(p^\prime)(\eta(p))} = J^{\prime(p^\prime)(p^{\prime -2})} = J^\prime
\]
is a homomorphism, forcing $\eta\:J \to J^\prime$ to be a strong homotopy.

(b) Let $F$ be a field and $\mu \in F[\bft]$ be an irreducible polynomial of degree at least $2$, making $K := F[\bft]/(\mu)$ a finite algebraic extension field of $F$. We denote by $\pi\:F[\bft] \to K$ the canonical projection and put $J := F[\bft]^{(+)}$, $J^\prime := K^{(+)}$. Pick an element $p^\prime \in K \setminus F1_K$. Then $p^\prime \in J^{\prime\times}$ and the linear map
\begin{align}
\label{ETAM} \eta\:J \longrightarrow J^\prime, \quad f \longmapsto \eta(f) :=\pi(f)p^{\prime -1}
\end{align} fits into the commutative diagram
\[
\xymatrix{J \ar[rr]^{\eta}\ar[rdd]^{\pi}  && J^{\prime(p^\prime)} \\ \\
& J^\prime \ar[ruu]_{R_{p^{\prime -1}}}^{\cong}.}
\]
where, by (\ref{e.ALIST}.\ref{PLUSP}), the horizontal arrow is a homomorphism. Thus $\eta\:J \to J^\prime$ is a homotopy. Assume that this homotopy is strong. Then there exists a $p \in J^\times$ such that $\eta\:J^{(p)} \to J^\prime$ is a homomorphism. But Ex.~\ref{e.ALIN} shows $J^\times = F[\bft]^\times = F^\times$, hence $p \in F^\times$. By \eqref{ETAM}, therefore, 
\[
1_K = 1_{J^\prime} = \eta(1_{J^{(p)}}) = \eta(p^{-1}) = p^{-1}\eta(1_J) = p^{-1}p^{\prime -1},
\]
and we arrive at the contradiction $p^\prime = p^{-1}1_K \in F1_K$. Thus $\eta$ is not a strong homotopy.
\end{sol}

\begin{sol}{pr.VALJO} \label{sol.VALJO} (a) We argue by induction if $n \geq 0$. For $n = 0$, we have to show $\lambda(1_J) = 0$, which follows from \eqref{UMULT} by setting $x = y = 1_J$: $\lambda(1_J) = \lambda(U_{1_J}1_J) = 3\lambda(1_J)$, and the assertion follows. For $n = 1$, there is nothing to prove, and if the formula is valid for all natural numbers $< n$ where $n \geq 2$, then \eqref{UMULT} implies $\lambda(x^n) = \lambda(U_xx^{n-2}) = 2\lambda(x) + \lambda(x^{n-2}) = 2\lambda(x) + (n - 2)\lambda(x) = n\lambda(x)$. If $n = -1$, we obtain $\lambda(x) = \lambda(U_xx^{-1}) = 2\lambda(x) + \lambda(x^{-1})$, hence $\lambda(x^{-1}) = -\lambda(x)$. And finally, for any integer $n > 0$, the definition (cf. Exc.~\ref{pr.EVINV}~\eqref{NEPO}) implies $\lambda(x^{-n}) = \lambda((x^{-1})^n) = n\lambda(x^{-1}) = -n\lambda(x)$. This completes the proof.

(b) We have $(vw)^2 = v^2w^2 = U_vw^2$. Hence (a) implies $2\lambda(vw) = \lambda((vw)^2) = \lambda(U_vw^2) = 2\lambda(v) + \lambda(w^2) = 2(\lambda(v) + \lambda(w))$, and the first part of (b) follows. The rest is clear.

(c) is straightforward.

(d) The statement  is obvious if $\lambda$ is trivial, i.e., if $\lambda(x) = 0$ for all $x \in J^\times$. Otherwise, (a) yields a least positive integer $m \in \Gamma := \lambda(J^\times)$, and we have $\IZ m \subseteq \Gamma$. Conversely, let $x \in J^\times$. Then $\lambda(x) \in \IZ$ can be written as $\lambda(x) = (2i + j)m + l$, with $i,j,l \in \IZ$, $j = 0,1$, $0 \leq l < m$. This implies $jm + l = -(2im - \lambda(x)) \in \Gamma$ since $\Gamma$ by \eqref{UMULT} and (a) is closed under the operation $(\alpha,\beta) \mapsto 2\alpha - \beta$. If $j = 0$, we conclude $l \in \Gamma$, hence $l = 0$ by the minimality of $m$. On the other hand, if $j = 1$, we conclude $m - l = (2m - (m + l)) \in \Gamma$, forcing again $l = 0$ by the minimality of $m$. In any event, $\lambda(x) = (2i + j)m \in \IZ m$, which completes the proof.

(e) We have to verify \eqref{VALDEF}--\eqref{NARTRI} and do so by a straightforward computation: For $x,y \in J$ we obtain
\begin{align*}
\lambda^{(p)}(x ) =\,\,&\infty\;\Longleftrightarrow\;\lambda(p) + \lambda(x) = \infty\;\Longleftrightarrow\;\lambda(x) = \infty\;\Longleftrightarrow\;x = 0, \\
\lambda^{(p)}(U^{(p)}_xy) =\,\,&\lambda(p) + \lambda(U_xU_py) = \lambda(p) + 2\lambda(x) + 2\lambda(p) + \lambda(y) \\
=\,\,&2\big(\lambda(p) + \lambda(x)\big) + \big(\lambda(p) + \lambda(y)\big) = 2\lambda^{(p)}(x) + \lambda^{(p)}(y), \\
\lambda^{(p)}(x + y) =\,\,&\lambda(p) + \lambda(x + y) \geq \lambda(p) + \min\big\{\lambda(x),\lambda(y)\big\} \\
=\,\,&\min\big\{\lambda(p) + \lambda(x),\lambda(p) + \lambda(y)\big\} = \min\big\{\lambda^{(p)}(x),\lambda^{(p)}(y)\big\},
\end{align*}
and finally,
\begin{align*}
(\lambda^{(p)})^{(q)}(x) =\,\,&\lambda^{(p)}(q) + \lambda^{(p)}(x) = \lambda(p) + \lambda(q) + \lambda(p) + \lambda(x) = 2\lambda(p) + \lambda(q) + \lambda(x) \\
=\,\,&\lambda(U_pq) + \lambda(x) = \lambda^{(U_pq)}(x).
\end{align*}
(f) There is clearly no harm in assuming $x,y,z \in J^\times$. Suppose first that the assertion holds for $y = 1_J$, so we have
\begin{align}
\label{LACI} \lambda(x \circ z) \geq \lambda(x) + \lambda(z) &&(x,z \in J),
\end{align} 
and let $y \in J^\times$ be arbitrary. Then we pass to the $y$-isotopes $J^{(y)}$ and $\lambda^{(y)}$ and combine (\ref{ss.JOISOT}.\ref{CIRP}) with \eqref{LACI} to conclude
\begin{align*}
\lambda(y) + \lambda(\{xyz\}) =\,\,&\lambda^{(y)}(x \circ^{(y)} z) \geq \lambda^{(y)}(x) +  \lambda^{(y)}(z) = 2\lambda(y) + \lambda(x) + \lambda(z).
\end{align*}
Hence \eqref{LATRI} follows. We are thus reduced to the case $y = 1_J$ and have to prove \eqref{LACI}. Fixing $x \in J^\times$ and letting $z \in J^\times$ be arbitrary, the estimate
\begin{align*}
\lambda(x \circ z) =\,\,&\lambda(U_{x,1_J}z) = \lambda(U_{x+1_J}z - U_xz - z) \\
\geq\,\,&\min\big\{2\lambda(x + 1_J) + \lambda(z),2\lambda(x) + \lambda(z),\lambda(z)\big\} \\
\geq\,\,&\min\big\{2\lambda(x),0\big\} + \lambda(z) = \min\big\{\lambda(x),-\lambda(x)\big\} + \lambda(x) + \lambda(z) \\
=\,\,&\lambda(x) + \lambda(z) - \max\big\{\lambda(x),-\lambda(x)\big\} = \lambda(x) + \lambda(z) - \vert\lambda(x)\vert
\end{align*}
shows that there exists a constant $\kappa \geq 0$, depending on $x$, such that
\begin{align}
\label{LACIKAP} \lambda(x \circ z) \geq \lambda(x) + \lambda(z) - \kappa &&(z \in J^\times).
\end{align}
We now consider (\ref{ss.JABAS.fig}.\ref{UXCY}) for $z$ in place of $y$ and apply the result to $1_J$. Then
\[
(x \circ z)^2 + z \circ U_xz = U_xz^2 + U_zx^2 + 2x \circ U_zx, 
\]
which in turn may be combined with (\ref{ss.JABAS.fig}.\ref{UXYC}) to yield
\[
(x \circ z)^2 = U_xz^2 + U_zx^2 + x \circ U_zx.
\]
This and \eqref{LACIKAP} imply
\begin{align*}
2\lambda(x \circ z) =\,\,&\lambda\big((x \circ z)^2\big) = \lambda(U_xz^2 + U_zx^2 + x \circ U_zx) \\
\geq\,\,&\min\big\{2\big(\lambda(x) + \lambda(z)\big),2\big(\lambda(z) + \lambda(x)\big),\lambda\big(x \circ U_zx\big)\big\} \\
\geq\,\,&\min\big\{2\big(\lambda(x) + \lambda(z)\big),2\big(\lambda(x) + \lambda(z)\big) - \kappa\big\} \\
=\,\,&2\big(\lambda(x) + \lambda(z)\big) - \kappa.
\end{align*}
Thus
\begin{align*}
\lambda(x \circ z) \geq \lambda(x) + \lambda(z) - \frac{\kappa}{2} &&(z \in J^\times),
\end{align*}
and iterating this procedure, we end up with
\begin{align*}
\lambda(x \circ z) \geq \lambda(x) + \lambda(z) - \frac{\kappa}{2^n} &&(z \in J^\times, n \in \IN).
\end{align*}
Hence, as $n \to \infty$, we arrive at \eqref{LACI}.
\end{sol}

\begin{sol}{pr.ALJODI} \label{sol.ALJODI} (a) Since $J$ by Thm.~\ref{t.AZLOC} is locally linear, $F[x] \subseteq J$ is a commutatuve associative $F$-algebra. Every $y \neq 0$ in $F[x]$ is invertible in $J$ and hence $U_y\:J \to J$ is bijective (Prop.~\ref{p.CHAIN.JOAL}). Thus $U_y$ via restriction induces a linear injection $F[x] \to F[x]$, and a dimension argument shows that $U_y\:F[x] \to F[x]$ is, in fact, bijective. Since, therefore, $y \in F[x]^\times$, it follows that $F[x]/F$ is indeed a finite algebraic field extension.

(b) For $x \in J$, the finite algebraic field extension $F[x]/F$ has degree $1$ since $F$ is algebraically closed. Hence $J = F\cdot 1_J$.

(c) For $x \in J$, the finite algebraic field extension $\IR[x]/\IR$ has degree at most $2$. Thus $1_J,x,x^ 2$ are linearly independent, and invoking Exc.~\ref{pr.CONFIELD}, we conclude that $J$ (viewed as a linear Jordan algebra with bilinear product $xy$) is conic. Writing $q$ for its norm and $M$ for its underlying real vector space, $(M,q,e)$, $e := 1_J$, is a pointed quadratic module over $\IR$. Now (\ref{ss.JOPOID.fig}.\ref{JOQUADL}) and (\ref{ss.LIJPOQ}.\ref{LIMPOQ}) imply $J = J(M,q,e)$. Since by Exc.~\ref{pr.INVPOI}, therefore, $q$ is anisotropic, it is either positive or negative definite. The latter alternative being excluded since $q(e) = 1 > 0$, the problem is solved.
\end{sol}


\solnsec{Section~\ref{s.PEIDEC}}

\begin{sol}{pr.LOLIC} \label{sol.LOLIC} Let $f = \sum_{i\geq 0}\alpha_i\bft^i \in k[\bft]$ and suppose $f(c) = 0$. By Prop.~\ref{p.LIKEX}, we must show $(\bft f)(c) = 0$. With $d := 1_J - c$ and $\alpha := \sum_{i\geq 1}\alpha_i \in k$ we have $0 = f(c) = \alpha_01_J + \alpha c = (\alpha_0 + \alpha)c + \alpha_0d$, and since $c,d$ are orthogonal idempotents, we conclude $(\alpha_0 + \alpha)c = \alpha_0d = 0$, hence $(\bft f)(c) = (\sum_{i\geq 0}\alpha_i)c = (\alpha_0 + \alpha)c = 0$, as desired.
\end{sol}

\begin{sol}{pr.ISTTWOZER} \label{sol.ISTTWOZER}  Write $J_i := J_i(c)$ for $i = 0,1,2$ and put $d := 1_J - c$. Applying (\ref{t.PEDESI}.\ref{UJIJJ}), we see that $U_v$ maps $J_2$ to $J_0$, hence via restriction induces a linear map $\vph\:J_2 \to J_0$. Write $w := v^{-1}$ as $w = w_2 + w_1 + w_0$, $w_i \in J_i$ for $i = 0,1,2$. Then (\ref{p.CHAIN.JOAL}.\ref{FOIN}) and (\ref{t.PEDESI}.\ref{UJIJJ}) again yield
\begin{align*}
v = U_vw = U_vw_2 + U_vw_1 + U_vw_0, \quad U_vw_i \in J_{2-i} &&(i = 0,1,2).
\end{align*}
Comparing Peirce components and using bijectivity of $U_v$ (Prop.~\ref{p.CHAIN.JOAL}), we conclude $w_2 = w_0 = 0$. Thus $w \in J_1$, and by (\ref{t.PEDESI}.\ref{UJIJJ}) combined with Cor.~\ref{c.PEIBAS}~(b) we have
\begin{align}
\label{WESQU} w^2 = p + q, \quad p := U_wd \in J_2^\times, \quad q := U_wc \in J_0^\times,
\end{align}
which implies
\begin{align}
\label{UVEP} U_vp = d, \quad U_vq = c.
\end{align}
Next,writing $p^{-1}$ (resp. $q^{-1}$) for the inverse of $p$ (resp. $q$) in $J_2$ (resp. $J_0$), we claim
\begin{align}
\label{UVED} U_vd = p^{-1}, \quad U_vc = q^{-1}. 
\end{align}
Indeed, from \eqref{WESQU}, \eqref{UVEP}, Cor.~\ref{c.PEIBAS}~(b) and Exc.~\ref{pr.EVINV} we deduce $U_qU_vc = U_qU_{v^2}q = U_qU_{(w^2)^{-1}}q = U_qU_{p^{-1} + q^{-1}}q = U_qU_{q^{-1}}q = U_qq^{-1}$, and since $U_q\:J_0 \to J_0$ is bijective, the second equation of \eqref{UVED} is proved. The first one follows analogously. In particular, $\vph$, viewed as a linear map $J_2 \to J_0^{(q)}$ preserves units. But it also preserves $U$-operators since for $x,y \in J_2$ we may apply \eqref{UVEP} to obtain
\begin{align*}
U_{\vph(x)}^{(q)}\vph(y) =\,\,&U_{U_vx}U_qU_vy = U_vU_xU_vU_qU_vy = U_vU_x U_{U_vq}y = \vph(U_xU_cy) = \vph(U_xy),
\end{align*}
as claimed. Thus $\vph\:J_2 \to J_0^{(q)}$ is a homomorphism. Interchanging the role of $c$ and $d$, we also obtain a homotopy $\psi\:J_0 \to J_2$ induced by $U_w = U_v^{-1}$. Hence $\vph$ is bijective with inverse $\psi$, and we have established the first part of the exercise.

As to the second, $v^2 = 1_J$ implies $w = v$, and \eqref{WESQU} yields $p = c$, $q = d$. Hence $\vph\:J_2 \to J_0$ is an isomorphism. Conversely, let this be so. By what we have just shown, $1_{J_0} = \vph(1_{J_2}) = q^{-1}$, forcing $q = d$. On the other hand, $\vph^{-1}\:J_0 \to J_2$ is an isomorphism induced by $U_{v^{-1}}$, and from \eqref{WESQU} we conclude $p = U_{v^{-1}}d = \vph^{-1}(d) = c$. Now \eqref{WESQU} again implies first $w^2 = 1_J$ and then $v^2 = (w^2)^{-1} = 1_J$ as well.
\end{sol}

\begin{sol}{pr.JTWOSP} \label{sol.JTWOSP}  From (\ref{c.PEIBAS}.\ref{PEICIR}) we deduce for $x \in J_2$ that $V_x$ stabilizes $J_1$, so the map $\vph$ is well dfined and linear. Moreover, $\vph(1_{J_2}) = \vph(c) = \Eins_{J_1}$ by (\ref{t.PEDESI}.\ref{JONE}). Now let $x,y \in J_2$, $z \in J_1$. Then (\ref{ss.JABAS.fig}.\ref{UXYV}) yields
\[
U_{\vph(x)}\vph(y)z = V_xV_yV_xz = U_{x,y}V_xz + V_{x,y}V_xz, 
\]
where $U_{x,y}V_xz \in \{J_2J_1J_2\} = \{0\}$ by (\ref{t.PEDESI}.\ref{JIJJJL}).Applying (\ref{ss.JABAS.fig}.\ref{VUXYVX}), we therefore conclude 
\[
U_{\vph(x)}\vph(y)z = V_{x,y}V_xz = V_{U_xy}z + U_xV_yz = V_{U_xy}z
\]
since $U_xV_yz \in U_{J_2}J_1 = \{0\}$ by (\ref{t.PEDESI}.\ref{UJIJJ}). Thus $\vph$ is a homomorphism of Jordan algebras, and the first part of the problem is solved. Moreover, if $J_2$ is simple and $J_1 \neq \{0\}$, then $\Ker(\vph) \subseteq J$ is an ideal $\neq J_2$ since $\vph(c) = \Eins_{J_1} \neq 0$. Thus $\Ker(\vph) = \{0\}$ and $\vph$ is injective.
\end{sol} 

\begin{sol}{pr.PETRCO} \label{sol.PETRCO}  (ii) $\Rightarrow$ (i). Clear.

(i) $\Rightarrow$ (ii). We may assume
\begin{align}
\label{NODIS} \{i,j\} \cap \{l,m\} \neq \emptyset
\end{align}
and then consider the following cases. 

\case{1}
$\{i,j\} \cap \{l,m\} \cap \{n,p\} = \emptyset$. We may assume $l = j$ by \eqref{NODIS}, and since $T$ is not connected, we have $m \neq n,p$ but also $j \neq n,p$ since we are in Case~$1$. Thus $T = (np,jm,ij)$ with $\{n,p\} \cap \{j,m\} = \emptyset$, as in the first alternative of (ii). 

\case{2}
$\{i,j\} \cap \{l,m\} \cap \{n,p\} \neq \emptyset$. Then we may assume $n = l = j$, and since $T = (ij,jm,jp)$ is not connected, we necessarily have $m \neq i,j,p$. The first part of the problem is thereby solved.

Now let $T = (ij,jl,ij)$ be as in the second part of the problem. If $i \neq l\neq j$, the first part of the problem shows that $T$ is not connected. Conversely, if $l = i$ or $l = j$, then $T$ is clearly connected and there is a unique positive integer $m$ having $ij = \{i,j\} = \{l,m\} = lm$.
\end{sol}

\begin{sol}{pr.MUPEAL} \label{sol.MUPEAL}  Let $R := k[\bft_1^{\pm 1},\dots,\bft_r^{\pm 1}]$ be the ring of Laurent polynomials in the variable $\bft_1,\dots,\bft_r$. As in the proof of Thm.~\ref{t.PEDECOM}~(c), we have
\begin{align*}
A \subseteq A_R =\,\,&\bigoplus_{i_1,\dots,i_r\in \IZ}(\bft_1^{i_1}\cdots\bft_r^{i_r}A), \\
\End_k(A) \subseteq \End_k(A_R) =\,\,&\bigoplus_{i_1,\dots,i_r\in \IZ}\big(\bft_1^{i_1}\cdots\bft_r^{i_r}\End_k(A)\big) \subseteq \End_R(A_R).
\end{align*}
(a) First of all, we have
\[
\Eins_A = L_{1_A}R_{1_A} = L_{\sum c_i}R_{\sum c_j} = \sum_{i,j=1}^r L_{c_i}R_{c_j} = \sum_{i,j=1}^r E_{ij}.
\]
It remains to show that the $E_{ij}$ are orthogonal projections. Since every base change of $A$ is flexible, we conclude
\[
\sum_{i,j=1}^r\bft_i\bft_jL_{c_i}R_{c_j} = L_{\sum \bft_ic_i}R_{\sum\bft_jc_j} = R_{\sum\bft_jc_j}L_{\sum \bft_ic_i} = \sum_{i,j=1}^r\bft_i\bft_jR_{c_j}L_{c_i},
\]
and comparing coefficients yields
\begin{align}
\label{ELIARJ} E_{ij} = L_{c_i}R_{c_j} = R_{c_j}L_{c_i} &&(1 \leq i,j \leq r).
\end{align}
Similarly,
\[
\sum_{i=1}^ r\bft_i^ 2L_{c_i} = L_{\sum\bft_i^ 2c_i} = L_{(\sum\bft_ic_i)^ 2} = L_{\sum\bft_ic_i}^2 = (\sum_{i=1}^r\bft_iL_{c_i})^ 2 = \sum_{i,j=1}^ r\bft_i\bft_jL_{c_i}L_{c_j}
\]
and the analogous relation for the right multiplication imply
\begin{align}
\label{ELIELJ} L_{c_i}L_{c_j} = \delta_{ij}L_{c_i}, \quad R_{c_i}R_{c_j} = \delta_{ij}R_{c_i} &&(1 \leq i,j \leq r).
\end{align}
Given $i,j,l,m = 1,\dots,r$, we now apply \eqref{ELIARJ}, \eqref{ELIELJ} and obtain 
\[
E_{ij}E_{lm} = L_{c_i}R_{c_j}R_{c_m}L_{c_l} = \delta_{jm}L_{c_i}R_{c_j}L_{c_l} = \delta_{jm}L_{c_i}L_{c_l}R_{c_j} = \delta_{il}\delta_{jm}L_{c_i}R_{c_j} = \delta_{il}\delta_{jm}E_{ij}.
\]
Hence the $E_{ij}$ are orthogonal projections of $A$, as claimed.

(b) Let $1 \leq i,j \leq r$ and $x \in A_{ij}$. For $1 \leq l \leq r$, we apply \eqref{ELIELJ} to obtain
\[
c_lx = L_{c_l}E_{ij}x = L_{c_l}L_{c_i}R_{c_j}x = \delta_{il}E_{ij}x = \delta_{il}x
\]
and, similarly, $xc_l = \delta_{jl}x$. Thus the left-hand side of \eqref{DEPECO} is contained in the right. Conversely, let $x \in A$ and assume $c_lx = \delta_{il}x$, $xc_l = \delta_{jl}x$ for $l = 1,\dots,r$. Then this holds, in particular, for $l = i$ and $l = j$, which implies $E_{ij}x = c_i(xc_j) = x$, hence $x \in A_{ij}$. This completes the proof of (b).

(c) Since $\bft_\lambda \in R^\times$ for $1 \leq \lambda \leq r$, we have
\begin{align}
\label{WETLA} w := \sum_{\lambda=1}^r\bft_\lambda c_\lambda \in J_R^\times, \quad w^{-1} = \sum_{\lambda=1}^r\bft_\lambda^{-1}c_\lambda. 
\end{align} 
For $1 \leq i,j \leq r$, we now claim
\begin{align}
\label{AIJEL} A_{ij} = \{x \in A\mid L_wx = \bft_ix,\;\;R_wx = \bft_jx\} = \{x \in A\mid L_{w^{-1}}x = \bft_i^{-1}x,\;\;R_{w^{-1}}x = \bft_j^{-1}x\},
\end{align}
where the second equation follows from the first since $L_{w^{-1}} = L_w^{-1}$, $R_{w^{-1}} = R_w^{-1}$ by (\ref{p.INVALT}.\ref{OPINVALT}). In order to prove the first, we decompose $x$ as $x = \sum x_{lm}$, $x_{lm} \in A_{lm}$ and apply \eqref{DEPECO} to obtain
\begin{align}
\label{ELWEX} L_wx = \sum_{\lambda,l,m}\bft_\lambda c_\lambda x_{lm} = \sum_{l,m}\bft_lx_{lm}, \quad R_wx = \sum_{\lambda,l,m}\bft_\lambda x_{lm}c_\lambda = \sum_{l,m}\bft_mx_{lm}.
\end{align}
Thus $L_wx = \bft_ix$ and $R_wx = \bft_jx$ if and only if $x = x_{ij}$, which completes the proof \eqref{AIJEL}. In analogy to the proof of Thm.~\ref{t.PEDECOM}, we now put
\[
T := \{\bft_1^{i_1}\cdots\bft_r^{i_r}\mid i_1,\dots,i_r \in \IZ\}, \quad T_1 := \{\bft_1,\dots,\bft_r\}
\]
and claim:
\begin{enumerate}[label=($\ast$)]
\item \label{sol.MUPEAL.ast}  If $s,t \in T$ admit an element $x \neq 0$ in $A$ such that $L_wx = sx$ and $R_wx = tx$, then $s,t \in T_1$.
\end{enumerate}
In order to see this, write $x = \sum x_{lm}$, $x_{lm}\in A_{lm}$. Then \eqref{ELWEX} implies
\[
\sum_{l,m} sx_{lm} = sx = L_wx = \sum_{l,m} \bft_lx_{lm}, \quad \sum_{l,m} tx_{lm} = tx = R_wx = \sum_{l,m}\bft_mx_{lm},
\]
and comparing components yields unique indices $i,j = 1,\dots,r$ satisfying $s= \bft_i$, $t = \bft_j$ and $x = x_{ij}$.

For $1 \leq i,j,l,m \leq r$, $x_{ij} \in A_{ij}$, $y_{lm} \in A_{lm}$, we now obtain
\begin{align}
\label{ELWEXI} L_w(x_{ij}y_{lm}) = \bft_i\bft_j\bft_l^{-1}x_{ij}y_{lm}, \quad R_w(x_{ij}y_{lm}) = \bft_j^{-1}\bft_l\bft_mx_{ij}y_{lm}
\end{align}
since the left and right Moufang identities (\ref{ss.MOUF}.\ref{LMOUF}) and (\ref{ss.MOUF}.\ref{RMOUF}) combined with (\ref{p.INVALT}.\ref{OPINVALT}) and \eqref{AIJEL} imply
\begin{align*}
L_w(x_{ij}y_{lm}) =\,\,&w\Big(x_{ij}\big(w(w^{-1}y_{lm})\big)\Big) = (wx_{ij}w)(w^{-1}y_{lm}) = \bft_i\bft_j\bft_l^{-1}x_{ij}y_{lm}, \\
R_w(x_{ij}y_{lm}) =\,\,&\Big(\big((x_{ij}w^{-1})w\big)y_{lm}\Big)w = (x_{ij}w^{-1})(wy_{lm}w) = \bft_j^{-1}\bft_l\bft_mx_{ij}y_{lm},
\end{align*}
as claimed. Specializing $l \mapsto j$, $m \mapsto l$ in \eqref{ELWEXI} yields $L_w(x_{ij}y_{jl}) = \bft_ix_{ij}y_{jl}$, $R_w(x_{ij}y_{jl}) = \bft_lx_{ij}y_{lj}$, hence \eqref{AIJAJM}. Similarly, \eqref{ELWEXI} for $l = i$, $m = j$ yields \eqref{AIJAIJ}, while \eqref{AIJALM} follows from \ref{sol.MUPEAL.ast} and \eqref{ELWEXI} since $j \neq l$ and $(i,j) \neq (l,m)$ imply that $\bft_i\bft_j\bft_l^{-1}$ or $\bft_j^{-1}\bft_l\bft_m$ does not belong to $T_1$. Finally, for $1 \leq i,j \leq r$, $i \neq j$, and $x \in A_{ij}$, we apply \eqref{AIJAIJ} to obtain $x^2 \in A_{ji}$. But then \eqref{ELIARJ}, \eqref{DEPECO} and right alternativity yield $x^2= E_{ji}x^2 = R_{c_i}L_{c_j}x^2 = ((c_jx)x)c_i = 0$, as claimed.

(d) Write $E_{ij}^{(+)}$, $1 \leq i \leq j \leq r$, for the Peirce projections of $\Omega$, viewed as a complete orthogonal system of idempotents in $A^{(+)}$. Comparing  (\ref{p.COORPE}.\ref{COORPE}) with \eqref{ELIARJ}, we let $1 \leq i \leq j \leq r$ and conclude
\begin{align*}
E_{ii}^{(+)} = U_{c_i} = L_{c_i}R_{c_i} = E_{ii}, \quad E_{ij}^{(+)} = U_{c_i,c_j} = L_{c_i}R_{c_j} + L_{c_j}R_{c_i} = E_{ij} + E_{ji}.
\end{align*}
The assertion follows. 
\end{sol}

\begin{sol}{pr.CONCOM} \label{sol.CONCOM}  Since the Peirce triples $(ij,ii,ij)$ and $(ij,jj,ij)$ are connected, (\ref{t.PEDECOM}.\ref{UIJJL}) implies $U_{v_{ij}}J_{ii} \subseteq J_{jj}$ and $U_{v_{ij}}J_{jj} \subseteq J_{ii}$.

(a) (i) $\Rightarrow$ (ii). Applying Exc.~\ref{pr.ISTTWOZER} to $J^\prime := J_2(c_i + c_j)$, we conclude that $U_{v_{ij}}\:J_{ii} \to J_{jj}$ is an isotopy. The first inclusion of (ii) now follows from Prop.~\ref{p.HOIN} combined with Thm.~\ref{t.ISTJIT}~(a). The second one follows by symmetry.

(a) (ii) $\Rightarrow$ (iii). Clear since $c_i \in J_{ii}^\times$ and $c_j \in J_{jj}^\times$.

(a) (iii) $\Rightarrow$ (iv). Since the Peirce triples $(ij,ll,ij)$ are not connected for $l = 1,\dots,r$, $i \neq l \neq j$, we have $v_{ij}^2 = \sum_{l=1}^rU_{v_{ij}}c_l = U_{v_{ij}}c_i + U_{v_{ij}}c_j \in (J_{ii} \oplus J_{jj})^\times$ (by (iii)) $\subseteq J^{\prime\times}$, and this is (a) (iv). Basically the same argument also yields (b) (i) $\Rightarrow$ (b) (ii), while (b)(ii) $\Rightarrow$ (b) (i) is obvious. 

(a) (iv) $\Rightarrow$ (i).  By definition (\ref{ss.COIN}), invertibility in a Jordan subalgebra implies invertibility in the ambient Jordan algebra, while by Prop.~\ref{p.PROIN}~(b), invertibility of $x^2 = U_x1_J$ implies invertibility of $x$. 

(c) (\ref{ss.JABAS.fig}.\ref{UUXYUZ}) for $x = v_{ij}$, $y = v_{jl}$, $z = 1_J$ yields
\begin{align}
\label{UVIJVJL} U_{v_{ij} \circ v_{jl}}c_i = U_{v_{ij}}U_{v_{jl}}c_i + U_{v_{jl}}U_{v_{ij}}c_i + V_{v_{ij}}U_{v_{jl}}V_{v_{ij}}c_i - U_{U_{v_{ij}}v_{jl},v_{jl}}c_i.
\end{align}
Since $i,j,l$ are mutually distinct, the Peirce triples $(jl,ii,jl)$, $(jl,ij,jl)$ and $(ij,jl,ij)$ are not connected. Moreover, $V_{v_{ij}}c_i = c_i \circ v_{ij} = v_{ij}$ by (\ref{t.PEDECOM}.\ref{JIJ}) and (\ref{t.PEDESI}.\ref{JI}). Summing up, therefore, \eqref{UVIJVJL} collapses to $U_{v_{ij} \circ v_{jl}}c_i = U_{v_{jl}}U_{v_{ij}}c_i$ and, similarly, $U_{v_{ij} \circ v_{jl}} c_l = U_{v_{ij}}U_{v_{jl}}c_l$. Combining this with the characterization of (strong) connectedness in (a) (resp. (b)), the assertion follows.
\end{sol}

\begin{sol}{pr.ISTCOM} \label{sol.ISTCOM}  (a) Free use will be made of the Peirce multiplication rules assembled in the present section. Let $1 \leq i,j \leq r$. Setting $d_\lambda := c_\lambda^{(p)}$ for $\lambda = 1,\dots,r$, we apply (\ref{ss.JOISOT}.\ref{SQUAP}) and deduce $d_i^{(2,p)} = U_{p_i^{-1}}\sum_{\lambda =1}^r p_\lambda = p_i^{-1} = d_i$, so $d_i$ is an idempotents in $J^{(p)}$. Now assume $i \neq j$ Since \eqref{DIAGOM} is a direct sum of ideals, we obtain $U_{d_i}^{(p)}d_j = U_{d_i}U_{\sum p_\lambda}p_j^{-1} = U_{d_i}p_j = 0$, while (\ref{ss.JOISOT}.\ref{CIRP}) yields $d_i \circ^{(p)} d_j = \{d_ipd_j\} = 0$. Thus the idempotents $d_i,d_j$ are orthogonal in $J^{(p)}$ ((\ref{t.PEDESI}.\ref{JZERO}) and Prop.~\ref{p.ORIDJO}). The equation $\sum_{\lambda = 1}^r d_\lambda = p^{-1} = 1_{J^{(p)}}$ shows, therefore, that $\Omega^{(p)}$ is a complete orthogonal system of idempotents in $J^{(p)}$.  Write $E_{lm}$ (resp. $E_{lm}^{(p)}$) for the Peirce projections of $\Omega$ in $J$ (resp. $\Omega^{(p)}$ in $J^{(p)}$). Then $E_{lm}^{(p)} = E_{lm}U_p$ by (\ref{p.COORPE}.\ref{COORPE}), hence $J_{ij}^{(p)} = E_{ij}U_pJ = E_{ij}J = J_{ij}$, and the proof of (a) is complete.

(b) For $1 \leq i \leq r$ we have $(c_i^ {(p)})^ {(q)} = q_i^{(-1,p_i)} = U_{p_i^{-1}}q_i^{-1} = (U_{p_i}q_i)^{-1}$, and the assertion follows.

(c) Suppose $v_{1i} \in J_{1i}$ connects $c_1$ and $c_i$ for $2 \leq i \leq r$. By Exc.~\ref{pr.CONCOM}~(a) we have $p_i := U_{v_{1i}}c_1 \in J_{ii}^\times$. Put $p_1 := c_1 \in J_{11}^\times$ and $p := \sum_{\lambda = 1}^rp_\lambda \in \Diag_\Omega(J)^\times$. By the solution to Exc.~\ref{pr.ISTTWOZER}, $w_{1i} := v_{1i}^{-1}$ (the inverse of $v_{1i}$ in $J_2(c_1 + c_i)$), belongs to $J_{1i} = J_{1i}^{(p)}$ and satisfies $U_{w_{1i}}p_i = c_1$. Note that $v_{1i}^2 = q_{i1} + p_i$, for some $q_{i1} \in J_{11}^\times$ and hence $w_{1i}^2 = q_{i1}^{-1} + p_i^{-1}$. Hence (\ref{ss.JOISOT}.\ref{SQUAP}) implies
\begin{align*}
w_{1i}^{(2,p)} =\,\,& U_{w_{1i}}p = U_{w_{1i}}(p_1 + p_i) = U_{w_{1i}}p_i + U_{w_{1i}}c_1 = c_1 + U_{w_{1i}^2}p_i = c_1 + U_{p_i^{-1}}p_i \\
=\,\,&c_1 + p_i^{-1} = c_1^{(p)} + c_i^{(p)}.
\end{align*}
Thus $c_1^{(p)}$ and $c_i^{(p)}$ are strongly connected by $w_{1i} \in J_{1i}^{(p)}$.
\end{sol}

\begin{sol}{pr.LIFCOM} \label{sol.LIFCOM}  We argue by induction on $r$. For $r = 1$, the assertion agrees with Exc.~\ref{pr.PARALIFT}~(c). Now let $r > 1$ and assume the statement holds for $r - 1$ in place of $r$. Put $c^\prime := c_r^\prime$ and use Exc.~\ref{pr.PARALIFT}~(c) again to find an idempotent $c = c_r \in J$ satisfying $\vph(c) = c^\prime$. Write $E_i$ (resp. $J_i$) for the Peirce projections (resp. components) of $J$ relative to $c$ and $E_i^\prime$ (resp. $J_i^\prime$) for the Peirce projections (resp. components) of $J^\prime$ relative to $c^\prime$, where $i = 0,1,2$. We have $\vph \circ E_i = E_i^\prime \circ \vph$ and hence $\vph(J_i) = (\vph \circ E_i)(J) = E_i^\prime \circ \vph(J) = E_i^\prime(J^\prime) = J_i^\prime$. In particular, $\vph_0 := \vph\vert_{J_0}\:J_0 \to J_0^\prime$ is a surjective homomorphism of Jordan algebras whose kernel, $\Ker(\vph_0) = J_0 \cap \Ker(\vph) \subseteq \Ker(\vph)$ is a nil ideal in $J_0$. Moreover, $J_0^\prime$ contains $(c_1^\prime,\dots,c_{r-1}^\prime)$ as a complete orthogonal system of idempotents. Hence the induction hypothesis yields a complete orthogonal system $(c_1,\dots,c_{r-1})$ of idempotents in $J_0$ such that $\vph(c_i) = c_i^\prime$ for $1 \leq i < r$. By Prop.~\ref{p.ORIDJO} and Cor.~\ref{c.ORCOJO}, therefore, $(c_1,\dots,c_{r-1},c_r)$ is a complete orthogonal system of idempotents in $J$ with the desired properties.
\end{sol}

\begin{sol}{pr.CLIFSIM} \label{sol.CLIFSIM}  Write $t$ for the trace and $x \mapsto \bar x$ for the conjugation of $(M,q,e)$. Let $I \subset J := J(M,q,e)$ be an ideal that remains stable under conjugation. We claim $I = \{0\}$. In order to see this, let $x \in I$. Then $I$ contains $U_x\bar x^2 = q(x)^21_J$ (by (\ref{ss.JOPOID.fig}.\ref{JOCONJC})), and from $I \neq J$ we deduce $q(x) = 0$. Linearizing this, we obtai $q(x,y) = 0$ for all $y \in I$. Now let $y \in J \setminus I$. Since $\bar x \in I$ and $I \subset J$ is an ideal, $U_y\bar x = q(x,y)y - q(y)x$ belongs to $I$. But so does the second summand on the right of this expression, and we conclude $q(x,y)y \in I$, hence $q(x,y) = 0$. Summing up, we have shown $x \in \Rad(q) = \{0\}$, and our intermediate claim $I = \{0\}$ is proved.

Now assume that $J$ is not simple and let $I$ be any proper ideal of $J$. Then so is $\bar I$, and the ideals $I \cap \bar I$, $I + \bar I$ are stable under conjugation. This and our intermediate assertion imply $J = I \oplus \bar I$ as a direct sum of ideals. As observed in \ref{ss.DISUID}, $I$ and $\bar I$ are Jordan algebras in their own right such that the conjugation induces an isomorphism from $I$ onto $\bar I$. Put $c := 1_I$. Then $\bar c = 1_{\bar I}$ and $1_J = c + \bar c$. Moreover, $c$ and $\bar c$ are idempotents in $J$, and the assumption $q(c) = 1$ would imply $c \in J^\times$ (Exc.~\ref{pr.INVPOI}), hence $c = 1_J$, $\bar c = 0$, a contradiction. Hence $q(c) = 0$ (since $q(c) \in F$ by (\ref{ss.JOPOID.fig}.\ref{JOQUU}) is an idempotent), and  (\ref{ss.JOPOID.fig}.\ref{JOQUAD}) reduces to $0 \neq c = c^2 = t(c)c$. Taking traces, we conclude $0 \neq t(c) = t(c)^2$, hence $t(c) = 1$. Summing, we have shown that $c$ and $\bar c$ are elementary idempotents of $J$, whence the corresponding Peirce decomposition by Prop.~\ref{p.PEPOQUA} implies $I = J_2(c) = Fc$, $\bar I = J_0(c) = F\bar c$. Thus the assignment $(\alpha,\beta) \mapsto \alpha c + \beta\bar c$ defines an isomorphism from $(F \times F)^{(+)}$ onto $J$.
\end{sol}



\solnchap{Solutions for Chapter~\ref{c.CUNOS}}

\solnsec{Section~\ref{s.CUNOJA}}

\begin{sol}{pr.NSPCF} \label{sol.NSPCF}  (a) For $x \in X$, the assignment $y \mapsto N(x,y)$ defines linear form on $X$. By regularity, therefore, we find a unique element $x^\sharp \in X$ such that $T(x^\sharp,y) = N(x,y)$ for all $y \in X$. On the other hand, given $y \in X$, the assignment $x \mapsto N(x,y_R)$ for $R \in \kalg$, $x \in X_R$ defines a polynomial law $N(\emptyslot,y)\:X \to k$ which is homogeneous of degree $2$, hence a quadratic form (Exc.~\ref{pr.MULI}~(b)). But then $x \mapsto x^\sharp$ must be a quadratic map since $T$ is regular, and the gradient identity holds in all scalar extensions. Next we prove that $X$ together with $1,\sharp,N$ is a cubic array. By Euler's differential equation and Exc.~\ref{pr.EXOCUB}, \eqref{EFEXPLY}, we have $T(1,y) = N(1,1)N(1,y) - N(1,1,y) = 3N(1,y) - 2N(1,y) = N(1,y)$, which implies $1^\sharp = 1$, as claimed. To finish the proof of (a), it remains to establish the unit identity. Linearizing the gradient identity, we conclude $T(x \times y,z) = N(x,y,z)$ is totally symmetric in $x,y,z \in X$, where as usual $x \times y$ is the bilinearization of the adjoint. Hence
\begin{align*}
T(x,1 \times y) =\,\,&T(x \times y,1) = N(1,x,y) = N(1,x)N(1,y) - T(x,y) = T(x)T(y) - T(x,y) \\
=\,\,&T\big(x,T(y)1 - y\big),
\end{align*} 
and the unit identity follows again from  regularity of $T$.

(b) Let $1^\prime$, $\pN$, $\pT$ be the base point, norm, bilinear trace, respectively, of $\pX$. We first show that $\pX$ is regular. Since $\vph$ is linear, the chain rule in the form (\ref{ss.DICA}.\ref{DICL}) implies $\pN(\vph(x),\vph(y)) = N(x,y)$ and then $\pN(\vph(x),\vph(y),\vph(z)) = N(x,y,z)$ in all scalar extensions of $X$. Since $\vph$ also preserves base points, it therefore preserves bilinear traces as well: $\pT(\vph(x),\vph(y)) = T(x,y)$. Now let $\px \in \Rad(\pT)$ and write $\px = \vph(x)$ for some $x \in X$. Then $T(x,y) = \pT(\vph(x),\vph(y)) = 0$ for all $y \in X$, hence $x = 0$ since $T$ is regular. On the other hand, let $\plambda\:\pX \to k$ be a linear form. Then so is $\lambda := \plambda \circ \vph\:X \to k$, and regularity of $T$ yields an element $x \in X$ such that $\lambda = T(x,\emptyslot)$. Since $\vph$ is surjective, this implies $\plambda = \pT(\vph(x),\emptyslot)$, and we have shown that $\pX$ is indeed regular.  

For $x \in \Ker(\vph)$, $y \in X$ we have $T(x,y) = T^\prime(\vph(x),\vph(y)) = 0$, hence $x = 0$ since $T$ is regular. Thus, $\vph$ is injective and, altogether, bijective. On the other hand, the gradient identity gives
\begin{align*}
T^\prime\big(\vph(x^\sharp),\vph(y)\big) =\,\,&T(x^\sharp,y) = N(x,y) = N^\prime\big(\vph(x),\vph(y)\big) =T^\prime\big(\vph(x)^{\sharp^\prime},\vph(y)\big),
\end{align*} 
and regularity of $T^\prime$ combines with the surjectivity of $\vph$ to yield $\vph(x^\sharp) = \vph(x)^{\sharp\prime}$, as claimed.\end{sol}

\begin{sol}{pr.CHACUN} \label{sol.CHACUN}  If $X$ is a cubic norm structure over $k$, then
\eqref{RANI}--\eqref{RADJI} hold by definition, while
\eqref{QUAPO}, \eqref{CUPO} have been recorded in (\ref{ss.BACUNO.fig}.\ref{SADJ}),
(\ref{ss.BACUNO.fig}.\ref{DADJ}), respectively. Conversely, suppose
\eqref{RANI}--\eqref{CUPO} hold for all $x,y,z \in X$ and let
$R \in \kalg$. Since \eqref{RANI}, \eqref{RAGRA}, \eqref{QUAPO} are homogeneous of degree at most $2$ in each variable, they continue to hold in $X_R$, so it remains to verify
\eqref{RADJI}, \eqref{CUPO} over $R$. We first linearize \eqref{RAGRA} to conclude that $T(x \times y,z) = N(x,y,z)$ is totally symmetric in $x,y,z \in X$. Hence linearizing \eqref{QUAPO} with respect to $y$ implies
\begin{align}
\label{ITCUPO} x^\sharp \times (y \times z) + (x \times y) \times (x \times z) = T(x^\sharp,y)z
+ T(x^\sharp,z)y + T(x \times y,z)x
\end{align}
for all $x,y,z \in X$. Turning now to the proof of \eqref{CUPO}, linearity allows us to assume $y = v_R$ for some $v \in X$; we then
put
\[
M := \{ x \in X_R \mid x^\sharp \times (x \times y) = T(x^\sharp,y)x
+ N(x)y\}.
\]
By homogeneity, we have $rx \in M$ for all $x \in M$, $r \in R$.
Moreover, for $u \in X$, the validity of \eqref{CUPO} over $X$
yields
\begin{align*}
(u_R)^\sharp \times (u_R \times v_R) =\,\,&(u^\sharp)_R \times (u_R
\times v_R) = \big(u^\sharp \times (u \times v)\big)_R \\
=\,\,&\big(T(u^\sharp,v)u + N(u)v\big)_R =
T\big((u_R)^\sharp,v_R\big)u_R + N_R(u_R)v_R,
\end{align*}
hence $u_R \in M$. But the elements $u_R$, with $u \in X$, span $X_R$ as
an $R$-module. Therefore \eqref{CUPO} will hold in all of $X_R$ once
we have shown that $M$ is closed under addition. In order to do so,
let $x,z \in M$. Then the definition of $M$ and \eqref{ITCUPO} yield
\begin{align*}
(x + z)^\sharp \times \big((x + z) \times y\big) =\,\,&(x^\sharp + x
\times z + z^\sharp) \times (x \times y + z \times y) \\
=\,\,&x^\sharp \times (x \times y) + x^\sharp \times (z \times y) +
(x \times z) \times (x \times y) + \\
\,\,&(x \times z) \times (z \times y) + z^\sharp \times (x \times y)
+ z^\sharp \times (z \times y) \\
=\,\,&T(x^\sharp,y)x + N(x)y + T(x^\sharp,y)z + \\
\,\,&T(x^\sharp,z)y + T(x \times y,z)x + T(z^\sharp,y)x +
T(z^\sharp,x)y + \\
\,\,&T(x \times y,z)z + T(z^\sharp,y)z + N(z)y \\
=\,\,&T(x^\sharp + x \times z + z^\sharp,y)x + T(x^\sharp + x \times
z + z^\sharp,y)z + \\
\,\,&\big(N(x) + T(x^\sharp,z) + T(x,z^\sharp) + N(z)\big)y \\
=\,\,&T\big((x + z)^\sharp,y\big)(x + z) + N(x + z)y,
\end{align*}
forcing $x + z \in M$, and the proof of \eqref{CUPO} is complete.
Now \eqref{RADJI} will be treated in a similar manner by putting
\[
M^\prime := \{x \in X_R \mid x^{\sharp\sharp} = N(x)x\} \subseteq
X_R.
\]
Since both sides of \eqref{RADJI} are homogeneous of the same
degree, $N$ is stable under scalar multiplication: $x \in M^\prime$
and $r \in R$ implies $rx \in M^\prime$. Moreover, for $x \in X$ we
conclude, using the fact that $N$ is a polynomial law,
\begin{align*}
(x_R)^{\sharp\sharp} =\,\,&(x^{\sharp\sharp})_R = \big(N(x)x\big)_R
= \big(N(x)1_R\big)x_R = N(x_R)x_R,
\end{align*}
hence $x_R \in M^\prime$. But the elements $x_R$, $x \in X$, span
$X_R$ as an $R$-module. Therefore \eqref{RADJI} will hold in all of
$X_R$ once we have shown that $M^\prime$ is closed under addition.
In order to do so, let $x,y \in M^\prime$. Then \eqref{RAGRA},
\eqref{QUAPO}, \eqref{CUPO} and the definition of $M^\prime$ imply
\begin{align*}
(x + y)^{\sharp\sharp} =\,\,&(x^\sharp + x \times y +
y^\sharp)^\sharp \\
=\,\,&x^{\sharp\sharp} + x^\sharp \times (x \times y) + x^\sharp
\times y^\sharp + (x \times y)^\sharp + (x \times y) \times y^\sharp
+ y^{\sharp\sharp} \\
=\,\,&N(x)x + T(x^\sharp,y)x + N(x)y + T(x^\sharp,y)y + \\
\,\,&T(x,y^\sharp)x + T(x,y^\sharp)y + N(y)x + N(y)y \\
=\,\,&\big(N(x) + T(x^\sharp,y) + T(x,y^\sharp) + N(y)\big)\big(x +
y) = N(x + y)(x + y),
\end{align*}
forcing $x + y \in M^\prime$, as desired. 
\end{sol}

\begin{sol}{pr.GENTWO} \label{sol.GENTWO}  We first linearize (\ref{ss.BACUNO.fig}.\ref{LCADJ}) and obtain
\begin{align}
\label{LLADJ} u \times (v \times w) +\,\,&v \times (w \times u) + w
\times (u \times v) \\
=\,\,&\big(T(u)S(v,w) + T(v)S(w,u) + T(w)S(u,v) - T(u \times v,w)\big)1 - \notag \\
\,\,&\big(S(u,v)w + S(v,w)u + S(w,u)v\big) - \notag \\
\,\,&\big(T(u)v \times w + T(v)w \times u + T(w)u \times v\big)
\notag
\end{align}
for all $u,v,w \in X$. We must show that the submodule $M \subseteq
X$ spanned by the elements assembled in \eqref{TWOGEN} is stable
under the adjoint map. First of all, we have $1^\sharp = 1 \in M$ by
the base point identities (\ref{ss.BACUNO.fig}.\ref{BAPO}) and $1 \times M
\subseteq M$ by the unit identity (\ref{ss.BACUNO.fig}.\ref{UNI}).
Moreover, $M$ is spanned by the elements $1,s,t,s \times t$, $s \in
\{x,x^\sharp\}$, $t \in \{y,y^\sharp\}$. Hence it will be enough to
show that, for all $s,s^\prime \in \{x,x^\sharp\}$, $t,t^\prime \in
\{y,y^\sharp\}$,
\begin{align}
\label{ESTSHA} s^\sharp,t^\sharp,(s \times t)^\sharp \in\,\,&M, \\
\label{ESESPR} s \times s^\prime,s \times (s^\prime \times t) \in \,\,&M, \\
\label{ESESTE} t \times t^\prime,t \times (s \times t^\prime) \in \,\,&M, \\
\label{ESTEES} (s \times t) \times (s^\prime \times t^\prime) \in
\,\,&M.
\end{align}
The adjoint identity (\ref{ss.BACUNO.fig}.\ref{ADJI}) implies $s^\sharp
\in \{N(x)x,x^\sharp\}$, $t^\sharp \in \{N(y)y,y^\sharp\}$, hence
not only the first two inclusions of \eqref{ESTSHA} but also
$s^\sharp \times t^\sharp \in M$, while the third one now follows
from (\ref{ss.BACUNO.fig}.\ref{SADJ}). In the first inclusion of
\eqref{ESESPR} we may assume $s = x$, $s^\prime = x^\sharp$, whence
the assertion follows from \eqref{CADJ}. The second relation of
\eqref{ESESPR} follows from (\ref{ss.BACUNO.fig}.\ref{LCADJ}) for $s =
s^\prime$ and from (\ref{ss.BACUNO.fig}.\ref{ADJD})), (\ref{ss.BACUNO.fig}.\ref{DADJ}) for $s \neq s^\prime$.
Next, \eqref{ESESTE} is \eqref{ESESPR} with $x$ and $y$
interchanged, so we are left with \eqref{ESTEES}. We first note that
$1 \times M \subseteq M$ and \eqref{ESESPR}, \eqref{ESESTE} imply
\begin{align}
\label{ESTETI} s \times M + t \times M \subseteq M.
\end{align}
Then we combine \eqref{LLADJ} with \eqref{ESTSHA}--\eqref{ESESTE}
to conclude
\[
(s \times t) \times (s^\prime \times t^\prime) \equiv -s^\prime
\times \big(t^\prime \times (s \times t)\big) - t^\prime \times
\big((s \times t) \times s^\prime\big) \bmod M,
\]
and \eqref{ESTETI} implies \eqref{ESTEES}.
\end{sol}


\solnsec{Section~\ref{s.BAPRCU}}

 \begin{sol}{pr.HOCUJO} \label{sol.HOCUJO}  (a) Since cubic Jordan algebras and cubic norm structures are identified via Cor.~\ref{c.ISNOJO}, $\vph$ will be a homomorphism of cubic Jordan algebras once we have shown $N_{J^\prime} \circ \vph = N_J$ as polynomial laws over $k$. Moreover, since the entire set-up is stable under base change, applying Prop.~\ref{p.ZADE} to the polynomial law $g := N_{J^\prime} \circ \vph\:J \to k$, we see that it suffices to prove $N_{J^\prime}(\vph(x)) = N_J(x)$ for all $x \in J$ having $\vph(x) \in J^{\prime\times}$. But then the adjoint identity combined with the fact that $\vph$ preserves adjoints implies $N_J(x)\vph(x) = \vph(x^{\sharp\sharp}) = \vph(x)^{\sharp\sharp} = N_{J^\prime}(\vph(x))\vph(x)$, hence the assertion since $\vph(x)$ by Thm.~\ref{t.CHACUJ}~\ref{c.CUJUNIMO} is unimodular.

(b) Since $\vph$ preserves unit elements and norms, it is a surjective homomorphism of pointed cubic forms in the sense of Exc.~\ref{pr.NSPCF}. By part (b) of that exercise, therefore, $\vph$ is an isomorphism of cubic arrays, hence of cubic Jordan algebras.
\end{sol}

\begin{sol}{pr.TRACUJO} \label{sol.TRACUJO}  Write $N = N_J$, $T = T_J$, $S = S_J$, $1 = 1_J$. By (\ref{ss.BACUNO.fig}.\ref{TSONE}), we have $3 = T(1) = 0$ in $k$, hence $2 = -1 \in k^\times$. Thus $J$ may be viewed as a linear Jordan algebra. Since not only the linear trace but also the  quadratic one (by (\ref{ss.BACUNO.fig}.\ref{ESTESH})) and the bilinear one (by (\ref{ss.BACUNO.fig}.\ref{BQUAT}))  are identically zero,  the bilinear multiplication of $J$ by \ref{ss.COLIJO} and (\ref{ss.COCUJO}.\ref{JBIADJ}) is given by
\[
xy = \frac{1}{2}x \circ y = -x \circ y = -x \times y.
\]
We conclude $x^2 = x^\sharp$ (which follows also from (\ref{ss.COCUJO}.\ref{JADJ})), while (\ref{ss.UOPLIJ}.\ref{UOPLIX}) and (\ref{ss.COPAR}.\ref{CUNOP}) imply
\[
2x(xy) - x^2y = U_xy = T(x,y)x - x^\sharp \times y = -x^\sharp \times y = x^2y,
\]
hence the left alternative law $x(xy) =x^2y$. Being also commutative, $J$ is in fact alternative. Since the entire set-up is stable under base change, it remains to show $N(xy) = N(x)N(y)$ for all $x,y \in J$. By Artin's theorem (Cor.~\ref{c.ARTIN}), the unital subalgebra of $J$ generated by $x,y$ is (commutative) associative. Hence (\ref{ss.COCUAL}.\ref{CAHALT}) implies  $N(xy)1 = (xy)^3 = x^3y^3 = N(x)N(y)1$, and since $1 \in J$ is unimodular, the assertion follows.
\end{sol}

\begin{sol}{pr.RACUNO} \label{sol.RACUNO} Let $X$ be a rational cubic norm structure over $k$, with base point $1$, adjoint $x \mapsto x^\sharp$, bilinearized adjoint $(x,y) \mapsto x \times y$, bilinear trace $T$ and norm $N$. Combining \eqref{HOAD} and \eqref{ABESH} with the bilinearity of $\times$, we conclude that $\sharp\:X \to X$ is indeed a quadratic map with bilinearization $\times$. Beside the map $N\:X \to k$, we now consider the map $g\:X \times X \to k$ defined by $g(x,y) := T(x^\sharp,y)$ for $x,y \in X$. Since the adjoint is a quadratic map with bilinearization $\times$, we deduce from \eqref{HONO}, \eqref{NOAB}, \eqref{TEASH} that conditions (i)--(iv) of Exc.~\ref{pr.CUMA} hold for $(N,g)$, so $(N,g)\:X \to k$ is a cubic map. By parts (d), (e) of that exercise, therefore, we find a unique cubic form $\tin := N \ast g\:X \to k$ such that $\tin(x) = N(x)$ and $\tin(x,y) = T(x^\sharp,y)$ for all $x,y \in X$. Here \eqref{ASH} and \eqref{ONESH} imply $1 = 1^{\sharp\sharp} = N(1)1$, hence $N(1) = 1$ since $1$ is unimodular. Thus $X$ together with the base point $1$, the adjoint $\sharp$ and the norm $\tin$ is a cubic array over $k$, denoted by $\tiX$. Write $\tit$ (resp. $\tis$) for the (bi-)linear (resp. the quadratic) trace of $X$. In view of \eqref{ONESH}, we have $\tit(x) := \tin(1,x) = T(1,x)$ and $\tis(x) = \tin(x,1) = T(x^\sharp,1)$ for $x \in X$, while Exc.~\ref{pr.CUMA}~\ref{CUMA.c} implies $T(x \times y,z) = g(x,y,z) = T(x,y \times z)$ for all $x,y,z \in X$. Putting $x = 1$ and combining (\ref{ss.TRACU}.\ref{TEES}) with \eqref{ONESHB}, we conclude $\tit(y,z) = \tit(y)\tit(z) - \tis(y,z) = T(1,y)T(1,z) - T(1,y \times z) = T(y,z)$. Thus $T$ is the bilinear trace of $\tiX$. Hence equation \eqref{RAGRA} of Exc.~\ref{pr.CHACUN} holds over $k$ But so do \eqref{RANI}, \eqref{RADJI}--\eqref{CUPO} of that exercise, by \eqref{ONESHB}, \eqref{ASH}--\eqref{ASHT}, respectively. Thus all equations of Exc.~\ref{pr.CHACUN} hold over $k$, forcing the cubic array $\tiX$ to be, in fact, a cubic norm structure.

If $\vph\:X \to Y$ is a homomorphism of rational cubic norm structures, then $\vph:\tiX \to \tiY$ is a linear map preserving base points and adjoints. By Exc.~\ref{pr.HOCUJO}~(a), therefore,  it is a homomorphism of cubic norm structures. Summing up, we have constructed a functor from $\kracuno$ to $\kcuno$. Conversely, let $X$ be a cubic norm structure over $k$, with base point $1$, adjoint $x \mapsto x^\sharp$, bilinear trace $T$ and norm $N$. We claim that $1$, $x \mapsto x^\sharp$, $(x,y) \mapsto x \times y$, $T$ and the set map $N_k\:X \to k$ make $X$ a rational cubic norm structure over $k$, denoted by $\Rat(X)$, i.e., equations \eqref{HOAD}--\eqref{ONESHB} hold. Indeed, \eqref{HOAD}--\eqref{ABESH} are obvious, while \eqref{NOAB}--\eqref{ONESHB} have been established in this order in (\ref{ss.BACUNO.fig}.\ref{EXPF}), (\ref{ss.BACUNO.fig}.\ref{EUL}), (\ref{ss.BACUNO.fig}.\ref{ADJI}), (\ref{ss.BACUNO.fig}.\ref{SADJ}), (\ref{ss.BACUNO.fig}.\ref{DADJ}), (\ref{ss.BACUNO.fig}.\ref{BAPO}), (\ref{ss.BACUNO.fig}.\ref{UNI}), respectively. If $\vph\:X \to Y$ is a homomorphism of cubic norm structures, then $\vph\:\Rat(X) \to \Rat(Y)$ is trivially one of rational cubic norm structures. Thus we obtain a functor from $\kcuno$ to $\kracuno$ which by Exc.~\ref{pr.CUMA}~\ref{CUMA.d} is easily seen to be inverse to the functor $\kracuno \to \kcuno$ constructed before. This completes the proof.
\end{sol}

\begin{sol}{pr.NORID} \label{sol.NORID}  (a) Let $x \in I$ and $y \in J$. If $(\mfa,I)$ is a cubic ideal in $J$, then $I$ is an ordinary one, and (\ref{ss.BACUNO.fig}.\ref{ADJJA}), (\ref{ss.BACUNO.fig}.\ref{LADJJA}) and (\ref{ss.BACUNO.fig}.\ref{BQUAT}) imply
\[
x^\sharp = x^2 - T_J(x,1_J)x + T_J(x^\sharp,1_J)1_J, \quad x \times y = x \circ y - T_J(x)y - T_J(y)x + \big(T_J(x)T_J(y) - T_J(x,y)\big)1_J,
\]
which both belong to $I$ by (i), (ii). Hence (iv) holds. Conversely, if this is so, then (i), (ii), (iv) imply
\[
U_xy = T_J(x,y)x - x^\sharp \times y \in I, \quad U_yx = T_J(y,x)y - y^\sharp \times x \in I,
\]
so $I \subseteq J$ is an ordinary ideal and, therefore, $(\mfa,I)$ is a cubic one.

(b) The first part follows immediately from \ref{ss.SELICUG}. As to the second, let us assume that $\sigma$ is surjective. Since $1_{\pJ_1} \in \pJ_1$ is unimodular, there exists a linear form $\plambda$ on $\pJ_1$ over $\pK$ such that $\plambda(1_{\pJ_1}) = 1_{\pK}$. Put $\mfa := \Ker(\sigma)$, $I := \Ker(\vph)$, write $\bar{\sigma}\:K/\mfa \overset{\sim}\to \pK$ for the isomorphism in $\kalg$ induced by $\sigma$ and $\lambda\:J_1 \to K/\mfa$ for the $k$-linear map rendering the diagram
\[
\xymatrix{J_1 \ar[r]_{\vph} \ar@{.>}[d]_{\lambda} & \pJ_1 \ar[d]^{\plambda} \\
K/\mfa \ar[r]_{\bar{\sigma}}^{\cong} & \pK}
\] 
commutative. Then one checks that $\lambda$ is, in fact, $K$-linear, $\lambda(1_{J_1}) = 1_{K/\mfa}$ and $\lambda(I) = \{0\}$. Hence the cubic ideal $(\mfa,I)$ in $J_1$ is separated. 

(c) For the rest of this exercise, we systematically abbreviate $1 = 1_J$, $T = T_J$, $S = S_J$, $\sharp = \sharp_J$, $N = N_J$. By (i), there is a unique way of viewing $J_0$ as a Jordan algebra over $k_0$ such that $\pi$ becomes a $\sigma$-semi-linear homomorphism of (abstract) Jordan algebras. The crux is to show that the Jordan algebra $J_0$ is derived from an appropriate cubic norm structure over $k_0$ making $\pi$ a $\sigma$-semi-linear homomorphism of \emph{cubic} Jordan algebras. For this purpose, we consult (iv) to obtain
a unique map $\sharp_0\:X_0 \to X_0$ producing in
\begin{align}
\vcenter{\label{EXEXZE} \xymatrix{J \ar[r]_{\pi} \ar[d]_{\sharp} & J_0 \ar@{.>}[d]^{\exists!\,\sharp_0} \\
J \ar[r]_{\pi} & J_0}}
\end{align}
a commutative diagram of set maps. Since $\sigma$ and $\pi$ are surjective, one checks that $\sharp_0$ is, in fact, a quadratic map over $k_0$, so \eqref{EXEXZE} by Cor.~\ref{c.BAQUA} commutes as a $\sigma$-semi-linear polynomial square. Next we we consult (ii) again to find a unique symmetric bilinear form $T_0\:X_0 \times X_0 \to k_0$ such that the diagram
\begin{align}
\vcenter{\label{NORID.EXTIEX} \xymatrix{X \times X \ar[r]_{\pi \times \pi} \ar[d]_{T} & X_0 \times X_0 \ar@{.>}[d]^{\exists!\,T_0} \\
k \ar[r]_{\sigma} & k_0}} 
\end{align}  
commutes. Defining set maps 
\[
f := N_k\:J \to k, \quad g := (DN)_k\:J \times J \to k,
\] 
Exc.~\ref{pr.CUMA}~\ref{CUMA.e} implies $f \ast g = N$ as cubic forms acting on $J$. Moreover, the adjoint identity yields $g(x,y) = T(x^\sharp,y)$ for all $x,y \in X$. Now define $g_0\:J_0 \times J_0 \to k_0$ by $g_0(x_0,y_0) := T_0(x_0^{\sharp_0},y_0)$ for $x_0,y_0 \in J_0$. By \eqref{EXEXZE}, \eqref{NORID.EXTIEX}, the diagram
\begin{align}
\vcenter{\label{EXEXPI} \xymatrix{J \times J \ar[r]_{\pi \times \pi} \ar[d]_{g} & J_0 \times J_0 \ar[d]^{g_0} \\
k \ar[r]_{\sigma} & k_0}}
\end{align}
commutes, and using condition (iii) one checks that there is a unique set map $f_0\:J_0 \to k_0$ rendering the diagram
\begin{align}
\vcenter{\label{EXEFKAY} \xymatrix{J \ar[r]_{\pi} \ar[d]_{f} & J_0 \ar@{.>}[d]^{\exists!\,f_0} \\k \ar[r]_{\sigma} & k_0}}
\end{align}
commutative. By \eqref{EXEXPI}, \eqref{EXEFKAY}, we obtain a cubic map $(f_0,g_0)\:J_0 \to k_0$, giving rise, by Exc.~\ref{pr.CUMA} to a cubic form $N_0 := f_0 \ast g_0\:J_0 \to k_0$. Here Exc.~\ref{pr.CUMASE} combined with \eqref{EXEXPI}, \eqref{EXEFKAY} shows that the $\sigma$-semi-linear polynomial square
\begin{align}
\vcenter{\label{EXENKAY} \xymatrix{X \ar[r]_{\pi} \ar[d]_{N} & X_0 \ar[d]^{N_0} \\
k \ar[r]_{\sigma} & k_0}}
\end{align}
commutes. We claim that the $k_0$-module $J_0$ together with $1_0 := \pi(1) = 1_{J_0}$, $\sharp_0$, $N_0$ is a cubic norm structure $X_0$ over $k_0$. Since $\sigma$, $\pi$ are surjective and equations \eqref{RANI}--\eqref{CUPO} of Exc.~\ref{pr.CHACUN} hold in $J$, they hold mutatis mutandis in $X_0$. Hence we need only show that $1_0 \in X_0$ is unimodular. Applying \eqref{LAUNX}, we see that $\lambda$ factors uniquely through $\pi$ to a linear form $\lambda_0\:X_0 \to k_0$ satisfying $\lambda_0(1_0) = 1_{k_0}$. This proves the assertion. Combining \eqref{EXEXZE}, \eqref{EXENKAY} with \ref{ss.SELICUG}~(a), it follows that $\pi\:X \to X_0$, with $X := X(J)$, is a $\sigma$-semi-linear homomorphism of cubic norm structures. By \ref{ss.SELICUG}~(c), therefore, $\pi\:J = J(X) \to J(X_0)$ is a $\sigma$-semi-linear homomorphism of cubic Jordan algebras, which forces $J_0 = J(X_0)$ as abstract Jordan algebras, and the solution of (b) is complete.

(d) If the linear trace $T$ of $J$ is surjective, some $u \in J$ has $T(u) = 1$. Now let $(\mfa,I)$ be a cubic ideal in $J$ and write $\sigma\:k \to k/\mfa$ for the canonical projection. Then
\[
\lambda\:J \longrightarrow k/\mfa, \quad x \longmapsto \lambda(x) := \sigma\big(T(u,x)\big),
\]
is a $k$-linear map satisfying \eqref{LAUNX} since $T(u,x) \in \mfa$ for all $x \in I$ by (ii). Hence $(\mfa,I)$ is separated.

(e) Setting $I := \mfa J$, condition (i) trivially holds, and one checks that so does (ii). Hence $(\mfa,I)$ is a cubic ideal in $J$. Let $\plambda$ be a linear form on $J$ satisfying $\plambda(1) = 1_k$ and write $\sigma\:k \to k/\mfa$ for the canonical projection. Then the $k$-linear map $\lambda := \sigma \circ \plambda\:J \to k/\mfa$ satisfies \eqref{LAUNX}, forcing $(\mfa,I)$ to be separated. From (\ref{ss.REDID}.\ref{ALPHEX}) we conclude that $\can_{J,k/\mfa}\:J \to J_{k/\mfa} = J/\mfa J$ is the canonical projection, which by \ref{ss.SELICUG}~(c), (d) is a $\sigma$-semi-linear homomorphism of cubic Jordan algebras. In the sense of (b), therefore, $J/\mfa J = J_0$ as cubic Jordan algebras over $k_0 = k/\mfa$. The answer to the final question of (e) is no. In order to see this, assume $k \ne \{0\}$, consider the multiplicative cubic alternative algebra $(k \times k)_{\cub}$ of \ref{e.SPLIQUJO} and pass to the associated cubic Jordan algebra $J = (k \times k)_{\cub}^{(+)}$. Then $I := \{0\} \times k$ is an ideal in $J$ satisfying $\{0\} \ne I^\sharp \subseteq k \times \{0\}$ by (\ref{e.SPLIQUJO}.\ref{KAKASH}), hence violating (iv). In particular, $I$ cannot be extended to a cubic ideal in $J$.

(f) We begin by establishing (ii), so let $x \in I$, $y \in J$ and write $\sigma\:k \to k/\mfa$ for the canonical projection. Then $I$ by (iv) contains $1 \times x =T(x)1 - x$, and applying the map $\lambda$ of (iii), we deduce $\sigma(T(x)) = T(x)1_{k/\mfa} = \lambda(1 \times x) = 0$, hence $T(x) \in \mfa$. Since $x \times y \in I$, again by (iv), we now conclude $T(x,y) = T(x)T(y) - T(x \times y) \in \mfa$. Thus the first inclusion of (ii) holds, while the remaining ones by (iv), (v) are now obvious. By (a), therefore, $(\mfa,I)$ is a separated cubic ideal in $J$.
\end{sol}

\begin{sol}{pr.SELIBA} \label{sol.SELIBA}  (a) By \ref{ss.SCEM}, the assignment in question gives an additive bijection from the $\pK$-linear maps $X_{\pK} \to \pX$ to the $\sigma$-semi-linear maps $X \to \pX$ (i.e., the $K$-linear maps $X \to\,_K{\pX}$). For a $\pK$-linear map $\pvph\:X_{\pK} \to \pX$ and $\vph := \pvph \circ \can_{X,\pK}$, it therefore remains to show that $\pvph$ is a homomorphism of cubic arrays over $\pK$ if and only if $\vph$ is a $\sigma$-semi-linear homomorphism of cubic arrays.
	
Let us first assume that $\pvph$ is a homomorphism of cubic arrays over $\pK$. It is a general fact following trivially from the definitions that if $\tau\:L \to \pL$, $\ptau\:\pL \to \ppL$ are morphisms in $\kalg$, $Y$ (resp. $\pY$, $\ppY$) is a cubic array over $L$ (resp. $\pL$, $\ppL$) and $\psi\:Y \to \pY$ (resp. $\ppsi\:\pY \to \ppY$) is a $\tau$- (resp. $\ptau$-)semi-linear homomorphism of cubic arrays, then $\ppsi \circ \psi\:Y \to \ppY$ is a $(\ptau \circ \tau)$-semi-linear homomorphism of cubic arrays. Thus, since $\pvph\:X_{\pK} \to \pX$ is a homomorphism of cubic arrays and $\can_{X,\pK}\:X \to X_{\pK}$ by \ref{ss.CONOJO}~(d) is a $\sigma$-semi-linear one, so is $\vph = \pvph \circ \can_{X,\pK}$. 

Conversely, assume $\vph\:X \to \pX$ is a $\sigma$-semi-linear homomorphism of cubic arrays. We abbreviate $\sharp = \sharp_X$, $\psharp = \sharp_{\pX}$, $N = N_X$, $\pN= N_{\pX}$, specialize (\ref{ss.SELPOS}.\ref{TRIKALG})  to
\begin{align}
\vcenter{\label{IDSISI} \xymatrix{& k \ar[d] & \\
& K \ar[ld]_{\Eins_K} \ar[rd]^{\sigma} & \\
K \ar[rr]_{\sigma} && \pK.}}
\end{align}
Since the $\sigma$-semi-linear polynomial square (\ref{ss.SELICUG}.\ref{POSH}) commutes, so does
\begin{align}
\vcenter{\label{EXKAYF} \xymatrix{X \ar[r]_{_K\vph} \ar[d]_{\sharp} & _K\pX \ar[d]^{_K\psharp} \\
X \ar[r]_{_K\vph} & _K{\pX},}}
\end{align} 
by \eqref{IDSISI} and (\ref{ss.SELPOS}.\ref{TRIPALG}). Now consider the diagram
\begin{align}
\vcenter{\label{EXCANEX} \xymatrix{X \ar[r]_{\can_{X,\pK}} \ar[d]_{\sharp} 
\ar@/^1.5pc/[rr]^{_K\vph} & _K(X_{\pK}) \ar[r]_{_K\pvph} \ar[d]^{_K(\sharp_{\pK})} & _K\pX \ar[d]^{_K\psharp} \\ 
X \ar[r]_{\can_{X,\pK}} \ar@/_1.5pc/[rr]_{_K\vph} & _K(X_{\pK}) \ar[r]_{_K\pvph} & _K\pX.}}
\end{align}
Diagram chasing, (\ref{ss.SELPOS}.\ref{TRIPALG}) and \eqref{EXKAYF} imply that $_K(\pvph \circ \sharp_{\pK})$ and $_K(\psharp \circ \pvph)$, being equalized by $\can_{X,\pK}$, are the same, so the diagram
\begin{align}
\vcenter{\label{EXKAP} \xymatrix{X_{\pK} \ar[r]_{\pvph} \ar[d]_{\sharp_{\pK}} & \pX \ar[d]^{\psharp} \\
X_{\pK} \ar[r]_{\pvph} & \pX}}
\end{align}
of polynomial laws over $\pK$ commutes by Exc.~\ref{pr.RESCAIN}. Similarly, since (\ref{ss.SELICUG}.\ref{PON}) commutes, so does
\begin{align}
\vcenter{\label{EXENK} \xymatrix{X \ar[r]_{_K\vph} \ar[d]_{N} & _\pX \ar[d]^{_K\pN} \\
K \ar[r]_{_K\sigma} & _K\pK}}
\end{align}
by \eqref{IDSISI} and (\ref{ss.SELPOS}.\ref{TRIKALG}). Using \eqref{EXENK}, (\ref{ss.SELPOS}.\ref{TRIPALG}) and diagram chasing in
\begin{align}
\vcenter{\label{KAXKAP} \xymatrix{X \ar[r]_{\can_{X,\pK}} \ar[d]_{N} \ar@/^1.5pc/[rr]^{_K\vph} & _K(X_{\pK}) \ar[r]_{_K\pvph} \ar[d]^{_K(N \otimes_K \pK)} & _K\pX \ar@/^1.5pc/[ld]^{_K\pN} \\
X \ar[r]_{_K\sigma} & _K\pK}}
\end{align}
conclude that $_K(N \otimes _K \pK)$ and $_K(\pN \circ \pvph)$ are equalized by $\can_{X;\pK}$, hence are the same, and the diagram
\begin{align}
\vcenter{\label{ENTEK} \xymatrix{X_{\pK} \ar[r]_{\pvph} \ar[d]_{N \otimes_K\pK} & \pX \ar[ld]^{\pN} \\
\pK}}
\end{align}
of polynomial laws over $\pK$ commutes. Summing up, therefore, the $\pK$-linear map $\pvph\:X_{\pK} \to \pX$ is, in fact, a homomorphism of cubic arrays over $\pK$.

(b) Let $\vph\:X \to \pX$ be a $\sigma$-semi-linear map preserving base points and making \eqref{EXFIX} a commutative diagram of set maps. Then it is a commutative $\sigma$-semi-linear polynomial square. Let $\pvph\:X_{\pK} \to \pX$ be the $\pK$-liner map such that $\pvph \circ \can_{X,\pK} = \vph$. Arguing as in (a), we conclude that $\pvph$ preserves adjoints. By Exc.~\ref{pr.HOCUJO}~(a) and Cor.~\ref{c.ISNOJO}, therefore, $\pvph$ is a homomorphism of cubic norm structures over $\pK$. But then, by (a), $\vph = \pvph \circ \can_{X,\pK}$ is a $\sigma$-semi-linear homomorphism of cubic norm structures.
\end{sol}

\begin{sol}{pr.CUBNIL} \label{sol.CUBNIL}  For the first part of the problem, if $T(x),S(x),N(x) \in k$ are nilpotent, then so are $T(x)x^2, S(x)x, N(x)1 \in J$. By (\ref{t.CUNOJO}.\ref{UNIVT}) and Exc.~\ref{pr.PAPONI}, therefore, $x^3$ is nilpotent, whence $x$ is nilpotent. Conversely, assume that $x$ is nilpotent. Since $N$ permits Jordan composition, it commutes with taking powers, and we conclude $N(x) \in \Nil(k)$. Moreover, applying Exc.~\ref{pr.PAPONI} again, there exists a positive integer $m$ such that $x^n = 0$ for all integers $n \geq m$. We now show $T(x),S(x) \in \Nil(k)$ by induction on $m$. For $m = 1$ there is nothing to prove. Next assume $m > 1$ and that the assertion holds for all positive integers $< m$. Since $(x^2)^n = x^{2n} = 0$ for all integers $n \geq m - 1$, the induction hypothesis implies that $T(x^2)$ and $S(x^2)$ are nilpotent. On the other hand, from (\ref{ss.BACUNO.fig}.\ref{ADJJA}) and the adjoint identity (\ref{ss.BACUNO.fig}.\ref{ADJI}) we deduce
\begin{align*}
T(x)^2 = 2S(x) + T(x^2), \quad S(x)^2 = 2T(x)N(x) + S(x^2).
\end{align*}
Since $N(x)$ is nilpotent, the induction hypothesis and the second equation imply that $S(x)$ is nilpotent, which combines with the first equation to show that $T(x)$ is nilpotent as well. This completes the proof of the first part.

In order to establish the second part, we put
\begin{align}
\label{RIHANI} I := \{x \in J\mid\;\forall y \in J:\;T(x,y),T(x^\sharp,y),N(x) \in \Nil(k)\}
\end{align}
and have to show $I = \Nil(J)$. For $x,y \in I$ and $z \in J$, we apply (\ref{ss.BACUNO.fig}.\ref{EXPF}), (\ref{ss.BACUNO.fig}.\ref{DGRAD}) and obtain
\begin{align*}
N(x + y) =\,\,&N(x) + T(x^\sharp,y) + T(x,y^\sharp) + N(y) \in \Nil(k), \\
T(x + y,z) =\,\,&T(x,z) + T(y,z) \in \Nil(k), \\
T\big((x + y)^\sharp,z\big) =\,\,&T(x^\sharp,z) + T(x,y \times z) + T(y^\sharp,z) \in \Nil(k),
\end{align*} 
hence $x + y \in I$. Thus $I \subseteq J$ is a $k$-submodule. Next we claim
\begin{align}
\label{ISHA} I^\sharp + I \times J \subseteq I.
\end{align}
The inclusion $I^\sharp \subseteq I$ follows immediately from \eqref{RIHANI}, the adjoint identity (\ref{ss.BACUNO.fig}.\ref{ADJI}) and (\ref{ss.BACUNO.fig}.\ref{NADJ}). It remains to show $x \times y \in I$ for all $x \in I$, $y \in J$. Given $z \in J$, we apply (\ref{ss.BACUNO.fig}.\ref{LNADJ}), (\ref{ss.BACUNO.fig}.\ref{SADJ}), (\ref{ss.BACUNO.fig}.\ref{DGRAD})
\begin{align*}
N(x \times y) =\,\,&T(x^ \sharp,y)T(x,y^\sharp) - N(x)N(y) \in \Nil(k), \\
T(x \times y,z) =\,\,&T(x,y \times z) \in \Nil(k), \\
T\big((x \times y)^\sharp,z\big) =\,\,&T(x^ \sharp,y)T(y,z) + T(x,y^\sharp)T(x,z) - T(x^\sharp,y^\sharp \times z) \in \Nil(k).
\end{align*}
Hence $x \times y \in I$ and \eqref{ISHA} holds. Now let $x \in I$ and $y \in J$. Then $U_xy = T(x,y)x - x^\sharp \times y \in I$ by \eqref{ISHA}, so $I \subseteq J$ is an inner ideal. But since
\begin{align*}
T(U_yx,z) =\,\,&T(x,U_yz) \in \Nil(k), \\
T\big((U_yx)^\sharp,z\big) =\,\,&T(U_{y^\sharp}x^\sharp,z) = T(x^\sharp,U_{y^\sharp}z) \in \Nil(k), \\
N(U_yx) =\,\,&N(y)^2N(x) \in \Nil(k) 
\end{align*} 
for all $z \in J$ by (\ref{ss.BACUNO.fig}.\ref{JCOMP}), (\ref{ss.BACUNO.fig}.\ref{TRAU}), (\ref{ss.BACUNO.fig}.\ref{UADJ}) and \eqref{ISHA}, we also have $U_yx \in I$ by \eqref{RIHANI}, so $I \subseteq J$ is an outer ideal and altogether, therefore, an ideal in $J$ which by the first part of the problem consists entirely of nilpotent elements. This proves $I \subseteq \Nil(J)$.

Conversely, let $x \in \Nil(J)$. Then $T(x),S(x),N(x) \in \Nil(k)$ by the first part of the problem, and we have $\Nil(k)J \subseteq \Nil(J)$ by Exc.~\ref{pr.PARNIL}~(c). Furthermore, for all $y \in J$, $\Nil(J) \ni U_xy = T(x,y)x - x^\sharp \times y$, hence $x^\sharp \times y \in \Nil(J)$, which by (\ref{ss.BACUNO.fig}.\ref{ESTESH}), (\ref{ss.BACUNO.fig}.\ref{BQUAT}) implies $S(x)T(y) - T(x^\sharp,y) = T(x^\sharp \times y) \in \Nil(k)$. Thus $T(x^\sharp,y) \in \Nil(k)$. Since $z := U_yx = T(y,x)y^\sharp - y^\sharp \times x$ belongs to the ideal $\Nil(J)$, so does $x + z$, which by (\ref{ss.BACUNO.fig}.\ref{DGRAD}),  (\ref{ss.BACUNO.fig}.\ref{BQUAT}), implies
\begin{align*}
T(x,y)^2 - 2T(x^\sharp,y^\sharp) =\,\,&T(x,y)^2 - T(x \times x,y^\sharp) = T(x,y)^2 - T(x,y^\sharp \times x) = T(x,U_yx) \\
=\,\,&T(x,z) = T(x)T(z) - S(x,z) \\
=\,\,&T(x)T(z) + S(x) + S(z) - S(x + z) \in \Nil(k).
\end{align*}
Thus $T(x,y) \in \Nil(k)$ and summing up we have shown $x \in I$. This completes the proof of \eqref{NILARB}. 

Finally, let $x,y \in J$. If $\frac{1}{2} \in k$, then $T(x^\sharp,y) = \frac{1}{2}T(x \times x,y) = \frac{1}{2}T(x,x \times y)$ by (\ref{ss.BACUNO.fig}.\ref{DGRAD}), so \eqref{NILARB} implies \eqref{NILTWO}. Similarly, if $\frac{1}{3} \in k$, then $N(x) = \frac{1}{3}T(x,x^\sharp)$ by (\ref{ss.BACUNO.fig}.\ref{EUL}), so \eqref{NILARB} implies \eqref{NILTHR}. Finally, \eqref{NILSIX} follows by combining \eqref{NILTWO} and \eqref{NILTHR}.
\end{sol}

\begin{sol}{pr.CUADUN} \label{sol.CUADUN}  (a) $N$ as defined by \eqref{NOHAJ} is clearly a cubic form such that
\[
N\big((r,u),(s,v)\big) = sq(u) + rq(u,v)
\]
for $R \in \kalg$, $r,s \in R$, $u,v \in J_R$. Defining $T$ (resp. $S$) as a linear (resp. quadratic) form on $J$ by $T(x) := N(1_{\haj},x)$ (resp. $S(x) := N(x,1_{\haj})$) for $x \in \haj$, we therefore conclude that \eqref{LITHAJ}, \eqref{QITHAJ} hold. For $x = (r,u) \in \haj_R$, this and (\ref{ss.JOPOID.fig}.\ref{JOQUAD}), (\ref{ss.JOPOID.fig}.\ref{JOQUAT}) imply
\begin{align*}
x^3 - T(x)x^2 + S(x)x - N(x)1_{\haj} =\,\,&\Big(r^3 - \big(r + t(u)\big)r^2 + \big(rt(u) + q(u)\big)r - rq(u), \\
\,\,&u^3 - \big(r + t(u)\big)u^2 + \big(rt(u)  + q(u)\big)u - rq(u)e\Big) \\
=\,\,&\Big(0,u^3 - t(u)u^2 + q(u)u - r\big(u^2 - t(u)u + q(u)e\big)\Big) \\
=\,\,&0, \\
x^4 - T(x)x^3 + S(x)x^2 - N(x)x =\,\,&\Big(r^4 - \big(r + t(u)\big)r^3 + \big(rt(u + q(u)\big)r^2 - rq(u)r, \\
\,\,&u^4 - \big(r + t(u)\big)u^3 + \big(rt(u) + q(u)\big)u^2 -rq(u)u\Big) \\
=\,\,&\Big(0,U_u\big(u^2 - t(u)u + q(u)e\big) - r\big(u^3 - t(u)u^2 + q(u)u\big)\Big) \\
=\,\,&0.
\end{align*} 
Since $q$ permits Jordan composition by (\ref{ss.JOPOID.fig}.\ref{JOQUU}), so does $N$ by \eqref{NOHAJ}, and we have shown that $\haj$ is a cubic Jordan algebra with norm $N$ over $k$ such that \eqref{NOHAJ}, \eqref{LITHAJ}, \eqref{QITHAJ} hold. Linearizing \eqref{QITHAJ} and applying (\ref{ss.TRACU}.\ref{TEES}), we also obtain \eqref{BITHAJ}. Finally, (\ref{ss.BACUNO.fig}.\ref{ADJJA}) implies, for $x = (\alpha,u)$, $\alpha \in k$, $u \in J$,
\begin{align*}
x^\sharp = x^2 - T(x)x + S(x)1_{\haj} =\,\,&\Big(\alpha^2 - \big(\alpha + t(u)\big)\alpha + \big(\alpha t(u)+ q(u)\big), \\
\,\,&u^2 -\big(\alpha + t(u)\big)u + \big(\alpha t(u) + q(u)\big)e\Big) \\
=\,\,&\Big(q(u),\alpha\big(t(u)e - u\big)\Big) \\
=\,\,&\big(q(u),\alpha\bar u\big).
\end{align*}
This proves \eqref{ADHAJ}, and (a) is solved.

(b) Applying Exc.~\ref{pr.COBRPL}, we see that $C^{(+)}$ agrees with the Jordan algebra, $J := J(C,n_C,1_C)$, of the pointed quadratic module $(C,n_C,1_C)$. Thus
\[
\hac^{(+)} = k^{(+)} \times C^{(+)} = k^{(+)} \times J = \haj 
\]
 is a cubic Jordan algebra over $k$ with norm $N\:\hac \to k$ given by
 \begin{align}
\label{NOHAC} N\big((r,u)\big) = rn_C(u) &&(R \in \kalg,\;\;r \in R\;\;u \in C_R).
 \end{align}
 Hence $\hac$ is a cubic alternative $k$-algebra with norm $N$ if and only if $C$ is multiplicative. Finally, Example~\ref{e.SPLIQUJO} fits into this picture by setting $C := k$.
 \end{sol}

\begin{sol}{pr.SJ0} \label{sol.SJ0}
\ref{pr.SJ0.a}: We write $\theta$ for the canonical image of $\bft$ in $K$ and use the formula (\ref{pr.CUADUN}.\ref{QITHAJ}) for $S_J$, which says that for $u = c\theta + d$, $c,d \in F$, we have
\[
S_J\big((\alpha, u)\big) = \alpha t_K(c\theta + d) + \beta c^2 + cd + d^2.
\]
Now the trace on $K$ is given by 
\[
t_K(c\theta + d) = ct_K(\theta) + dt_K(1_K) = c + 2d = c,
\]
so 
\begin{equation} \label{SJ0.1}
S_J\big((\alpha, u)\big) = \alpha c + \beta c^2 + cd + d^2.
\end{equation}
That is the equation with respect to the basis $(1,0)$, $(0, \theta)$, $(0, 1)$ of $J$.  Using instead the basis $(1,1)$, $(0,\theta+1)$, $(0,1)$, we have for $r,s,d \in F$:
\[
S_J\big((r(1,1) + s(0,\theta+1) + (0,d)\big) = S_J\big((r, s\theta + r+s+d)\big) = \beta s^2 + sd + d^2 + r^2,
\]
a quadratic form of the claimed isomorphism class.

\ref{pr.SJ0.b}: 
By formula (\ref{pr.CUADUN}.\ref{LITHAJ}), 
the kernel of $T_J$ consists of elements $(\alpha, u)$ with $u = \alpha \theta + d$.  For such an element, \eqref{SJ0.1} gives that $S_J((\alpha, u)) =(1+ \beta) \alpha^2  + \alpha d + d^2$, which is the claim.
\end{sol}

\begin{sol}{pr.SJ0.3}\label{sol.SJ0.3}
\ref{pr.SJ0.3a}: Again, we write $\theta$ for the canonical image of $\bft$ in $K$. Since $\theta \in K$ has trace zero, we conclude
\[
t_K(c\theta + d) = 2d,
\]
and use the formula (\ref{pr.CUADUN}.\ref{QITHAJ}) for $S_J$, which says that for $u = c\theta + d$ we have
\begin{equation} \label{SJ0.3.1}
S_J\big((\alpha, u)\big) = 2 \alpha d +  d^2 - \beta c^2.
\end{equation}
Re-writing this as $2(\alpha d+ d^2 / 2) - \beta c^2$, we see that it has the same isomorphism class as $\la 2 \ra \otimes [0, \frac12] \perp \la -\beta \raq$.  Since $[0, \frac12]$ has an isotropic vector (namely, $(1, 0)$), it follows that it is isomorphic to $\bfh$, proving the claim.

\ref{pr.SJ0.3b}: By the preceding formula for the trace on $K$ and (\ref{pr.CUADUN}.\ref{LITHAJ}),  the kernel of $T_J$ consists of elements $(-2d, c\theta + d)$.  For such an element, \eqref{SJ0.3.1} gives that $S_J((-2d, c\theta + d)) = -3d^2 - \beta c^2$, which is the claim.

\ref{pr.SJ0.3c}: While it is possible to approach this by building on \ref{pr.SJ0.3b}, it is perhaps simpler to argue directly, regardless of $\car(F)$, as follows.  For $J := \plalg{(F \times F \times F)}$, $J^0$ consists of triples $(x, y, -x-y) \in F^3$ and $S_J^0$ sends such an element to $y(-x-y) + x(-x-y) + xy = -(x^2 + xy+y^2)$.
\end{sol}

\begin{sol}{pr.CUMOPO} \label{sol.CUMOPO}  (i) $\Rightarrow$ (ii) Let $X$ be a cubic norm structure over $k$ as in (i). Then $X = (M,e,\sharp,N)$ for some quadratic map $\sharp\:M \to M$, $x \mapsto x^\sharp$, and some cubic form $N\:M \to k$. Write $T$ (resp. $S$) for the (bi-)linear (resp. quadratic) trace of $X$ and put $J := J(X) = J(M,q,e)$. For $x \in J^\times$, we combine \eqref{INVPOI} of Exc.~\ref{pr.INVPOI} with (\ref{c.REGINV}.\ref{CUINV}) to conclude $q(x), N(x) \in k^\times$ and $q(x)^{-1}\bar x = x^{-1} = N(x)^{-1}x^\sharp$. Hence Prop.~\ref{p.ZADE} implies that the identity
 \begin{align}
 \label{ENBAQU} N(x)\bar x = q(x)x^\sharp
 \end{align}
 holds strictly. We now put
 \begin{align}
 \label{CATMIT} \lambda := T - t\:M \longrightarrow k
 \end{align}
 as a linear form on $M$ and have $\lambda(e) = 1$ by (\ref{ss.JOPOID.fig}.\ref{JOA1}), (\ref{ss.BACUNO.fig}.\ref{TSONE}). Eventually, we will show that $\lambda$ satisfies (ii). Noting 
 \begin{align}
 \label{LAKE} M = ke \oplus L, \quad L := \Ker(\lambda),
 \end{align}
 we begin by using (\ref{ss.JOPOID.fig}.\ref{JOCONJ}), (\ref{ss.JOPOID.fig}.\ref{JOQUAD}), (\ref{ss.BACUNO.fig}.\ref{ADJJA}) to expand \eqref{ENBAQU}: we obtain $N(x)\bar x = t(x)N(x)e - N(x)x$ and 
 \begin{align*}
 q(x)x^\sharp =\,\,&q(x)x^2 - q(x)T(x)x + q(x)S(x)e \\
 =\,\,&q(x)t(x)x - q(x)^2e - q(x)T(x)x + q(x)S(x)e \\
 =\,\,&q(x)\big(S(x) - q(x)\big)e - \lambda(x)q(x)x &&(\text{by} \eqref{CATMIT}).
 \end{align*}
 Thus
 \[
 \Big(N(x) - \lambda(x)q(x)\Big)x = \Big(t(x)N(x) - q(x)\big(S(x) - q(x)\big)\Big)e.
 \]
 Differentiating at $x$ in the direction $y$ yields
 \begin{align*}
 \big(N(x,y) -& \lambda(x)q(x,y) -\lambda(y)q(x)\big)x + \big(N(x) - \lambda(x)q(x)\big)y \\
 =& \Big(t(x)N(x,y) + t(y)N(x) - q(x)\big(S(x,y) - q(x,y)\big) - q(x,y)\big(S(x) - q(x)\big)\Big)e.
 \end{align*}
 Setting $y = e$ and consulting (\ref{ss.BACUNO.fig}.\ref{TSONE}), we conclude
 \begin{align*}
 \big(S(x) -& \lambda(x)t(x) - q(x)\big)x + \big(N(x) - \lambda(x)q(x)\big)e \\
 =& \Big(t(x)S(x) + 2N(x) - q(x)\big(2T(x) - t(x)\big) -t(x)\big(S(x) - q(x)\big)\Big)e,
 \end{align*}
 which after a short computation simplifies to
 \[
 \big(N(x) - \lambda(x)q(x)\big)e = \big(S(x) - q(x) - \lambda(x)t(x)\big)x.
 \]
 We claim that the identity
 \begin{align}
 \label{ENEXLA} N(x) = \lambda(x)q(x)
 \end{align}
 holds strictly, equivalently, that the cubic form
 \begin{align}
 \label{EFENKA} F\:M \longrightarrow k, \quad x \longmapsto F(x) := N(x) - \lambda(x)q(x),
 \end{align}
 is zero. By what we have just proved we know
 \begin{align}
 \label{EFEXE} F(x)e = \big(S(x) - q(x) - \lambda(x)t(x)\big)x.
 \end{align}
 Since $F(e) = 0$, the expansion formula \eqref{EFEXTE} of Exc.~\ref{pr.EXOCUB} combined with \eqref{LAKE} implies that it suffices to show
 \begin{align}
 \label{EFUFE} F(u) = F(e,u) = F(u,e) = 0
 \end{align}
 strictly for all $ u \in L$. Specializing $x$ to $u$ in \eqref{EFEXE} and applying $\lambda$, we obtain $F(u) = 0$. Moreover, $F(e,u) = N(e,u) - \lambda(e)q(e,u) - \lambda(u) q(e) = T(u) - t(u) = \lambda(u)$ (by \eqref{CATMIT}) $= 0$, so it remains to show $F(u,e) = 0$. But
 \[
 F(u,e) = N(u,e) - \lambda(u)q(u,e) - \lambda(e)q(u) = S(u) - q(u),
 \]
 so we must show $S = q$ on $L$. Since $F(u) = 0$ for $u \in L$, applying \eqref{EFEXE} to $x = u$ yields
 \begin{align}
 \label{ESUQU} \big(S(u)- q(u)\big)u = 0,
 \end{align}
 hence after linearizing
 \begin{align}
 \label{ESUQV} \big(S(u) - q(u)\big)v + \big(S(u,v) - q(u,v)\big)u = 0.
 \end{align}
Note by \eqref{LAKE} that $L$ is a projective $k$-module, so localizing if necessary, we may assume that $L$ is free. Pick a basis $(u_i)_{i\in I}$ of $L$. Then \eqref{ESUQU}, \eqref{ESUQV} imply $S(u_i) = q(u_i)$ ($i \in I$), $S(u_i,u_j) = q(u_i,u_j)$ ($i,j \in I$, $i \neq j$). This shows $S = q$ on $L$ and completes the proof of \eqref{EFUFE}. Thus \eqref{ENEXLA} holds strictly. Returning to \eqref{ENBAQU}, we therefore obtain $\lambda(x)q(x)\bar x = N(x)\bar x = q(x)x^\sharp$, and Prop.~\ref{p.ZADE} implies \eqref{SHALA}. Now the adjoint identity yields $N(x)x = x^{\sharp\sharp} = \lambda(x)^2\bar x^\sharp = \lambda(x)^2\lambda(\bar x)x$. Since invertible elements are unimodular by Thm.~\ref{t.CHACUJ}~\ref{c.CUJUNIMO}, we end up with \eqref{NOLA}, which linearizes to 
\[
N(x,y) = 2\lambda(x)\lambda(\bar x)\lambda(y) + \lambda(x)^2\lambda(\bar y);
\]
putting $x = e$ gives \eqref{CATLA}, while for $y = e$ we obtain \eqref{ESLA}. This in turn linearizes to $S(x,y) = 2(\lambda(x)\lambda(y) + \lambda(x)\lambda(\bar y) + \lambda(\bar x)\lambda(y))$, and a short computation involving \eqref{CATLA} and (\ref{ss.BACUNO.fig}.\ref{BQUAT}) yields \eqref{BICAT}.
 
(ii) $\Rightarrow$ (i). If \eqref{QULA} holds, so obviously does \eqref{TRALA}. Define $X := (M,e,\sharp,N)$, where $\sharp$, $N$ are given by \eqref{SHALA}, \eqref{NOLA}, respectively. From \eqref{SHALA}, \eqref{NOLA} we conclude that $X$ is a cubic array, and the preceding computations can be repeated verbatim to show that the (bi-)linear and quadratic traces of $X$satisfy \eqref{CATLA}--\eqref{BICAT}. Now we show that $X$ is a cubic norm structure. Linearizing \eqref{SHALA} yields
\begin{align}
\label{BISHAL} x \times y = \lambda(x)\bar y + \lambda(y)\bar x,
\end{align}
hence $e \times x = \lambda(x)e + \bar x = (\lambda(x) + t(x))e - x = T(x)e - x$ by \eqref{TRALA}, \eqref{CATLA}. Thus the unit identity holds. From $x^{\sharp\sharp} = \lambda(x)^2\bar x^\sharp = \lambda(x)^2\lambda(\bar x)x$ and \eqref{NOLA} we read off the adjoint identity, while comparing $N(x,y) = 2\lambda(x)\lambda(\bar x)\lambda(y)
+ \lambda(x)^2\lambda(\bar y)$ with 
\[
T(x^\sharp,y) = 2\lambda(x^\sharp)\lambda(y) + \lambda(\overline{x^\sharp})\lambda(\bar y) = 2\lambda(x)\lambda(\bar x)\lambda(y) + \lambda(x)^2\lambda(\bar y)
\]
proves also the gradient identity. Summing up, therefore, $X$ is a cubic norm structure. It remains to show $J(X) = J(M,q,e)$, i.e., 
\begin{align}
\label{CUCLI} q(x,\bar y)x - q(x)\bar y = T(x,y)x - x^\sharp \times y. 
\end{align}
From \eqref{CATLA}, \eqref{BISHAL}, \eqref{SHALA} we deduce
\begin{align*}
T(x,y)x - x^\sharp \times y =\,\,&\big(2\lambda(x)\lambda(y) + \lambda(\bar x)\lambda(\bar y)\big)x - \lambda(x^\sharp)\bar y  - \lambda(y)\overline{x^\sharp} \\
=\,\,&\big(2\lambda(x)\lambda(y) + \lambda(\bar x)\lambda(\bar y)\big)x - \lambda(x)\lambda(\bar x)\bar y - \lambda(x)\lambda(y)x \\
=\,\,&\big(\lambda(x)\lambda(y) + \lambda(\bar x)\lambda(\bar y)\big)x - \lambda(x)\lambda(\bar x)\bar y,
\end{align*}
and by \eqref{QULA} (possibly linearized), this is the left-hand side of \eqref{CUCLI}. Thus (i) holds.

In order to prove \eqref{NILLA}, we put
\[
I_1 := \{x \in J\mid x\;\text{is nilpotent}\}, \quad I_2 := \{x \in J\mid \lambda(x),\lambda(\bar x) \in k\;\text{are nilpotent}\}.
\]
It suffices to show $\Nil(J) \subseteq I_1 \subseteq I_2 \subseteq \Nil(J)$. The first inclusion is obvious. For $x \in I_1$, we apply Exc.~\ref{pr.NIPOQUA} to conclude that $q(x) = \lambda(x)\lambda(\bar x)$ and $t(x) = \lambda(x) + \lambda(\bar x)$ are nilpotent. Hence so is $\lambda(x)t(x) = \lambda(x)^2 + q(x)$, which implies $x \in I_2$. This proves the second inclusion. For $x \in I_2$, we conclude that $q(x) = \lambda(x)\lambda(\bar x)$ is nilpotent, as is $q(x,y) = \lambda(x)\lambda(\bar y) + \lambda(y)\lambda(\bar x)$ for all $y \in J$. Thus $x \in \Nil(J)$, again by Exc.~\ref{pr.NIPOQUA}, and we have established the third inclusion.

For the remainder of this exercise, we may assume that $(M,q)$ is a quadratic space of rank $r > 1$. Since (iii) implies (i) by \ref{e.SPLIQUJO}, we need only show

(ii) $\Rightarrow$ (iii). By \eqref{LAKE}, the $k$-module $L$ is projective of rank $r - 1$. Let us first assume that $L$ is free, and let $(u_i)_{1\leq i<r}$ be a $k$-basis of $L$. With $u_0 := e$, therefore, $(u_i)_{0\leq i<r}$ is a $k$-basis of $M$, and by \eqref{QULA}, the matrix of $q$ relative to this basis is
\[
S := \left(\begin{matrix}
1 & t(u_1) & \cdots & t(u_{r-1}) \\
0 & 0      & \cdots & 0          \\
\vdots & \vdots & \ddots & \vdots \\
0 & 0      & \cdots & 0 
\end{matrix}\right),
\]
whence the matrix of $Dq$ is
\[
\pS = S + S^\trans  = \left(\begin{matrix}
2 & t(u_1) & \cdots & t(u_{r-1}) \\
t(u_1) & 0      & \cdots & 0          \\
\vdots & \vdots & \ddots & \vdots \\
t(u_{r-1}) & 0      & \cdots & 0 
\end{matrix}\right).
\]
Since $Dq$ is regular by hypothesis, the determinant of $\pS$ is invertible in $k$. On the other hand, expanding with respect to the last row yields
\[
\det(\pS) = (-1)^{r+1}t(u_{r-1})\det\left(\begin{matrix}
t(u_1) & \cdots & t(u_{r-1}) \\
0 & \cdots & 0 \\
\vdots & \ddots & \vdots \\
0 & \cdots & 0
\end{matrix}\right),
\]
which would be zero for $r > 2$. Hence $r = 2$ and $k^\times \ni \det(\pS) = -t(u_1)^2$, forcing $t\:L \to k$ to be surjective.

Now let $L$ be arbitrary. Localizing if necessary, the preceding considerations show that we still have $r = 2$ and $t\:L \to k$ surjective. In particular, some $c \in L$ satisfies $t(c) = 1$. Moreover $q(c) = 0$ by the definitions of $q$ and $L$, so $c$ is an elementary idempotent in $J = J(M,q,e)$, and with $d := e - c$ we conclude $J = kc \oplus kd$ as a direct sum of ideals such that $\lambda(c) = 0$ and $\lambda(d) = 1$. By \eqref{NOLA}, therefore, the relation
\[
N(\xi_1c + \xi_2d) = \lambda(\xi_1c + \xi_2d)^2\lambda(\overline{\xi_1c + \xi_2d}) = \xi_1\xi_2^2
\] 
holds strictly for all $\xi_1,\xi_2 \in k$, whence the assignment $\xi_1c + \xi_2d \mapsto (\xi_1,\xi_2)$ gives an isomorphism $J \overset{\sim}\to (k \times k)^{(+)}_\cub$ of cubic Jordan algebras.

In view of (iii), the final assertion about quadratic \'etale algebras is obvious. 
\end{sol}

\begin{sol}{pr.PSUNOS} \label{sol.PSUNOS}  (a) The arguments solving Exc.~\ref{pr.RACUNO}~(a) go through for rational cubic norm pseudo-structures without change and prove (a).

(b) (i) $\Rightarrow$ (ii). Let $1 \leq i \leq n$. Then (i) and Euler's differential equation imply $\sigma(e_i^\sharp,e_i) = \mu(e_i,e_i) = 3\mu(e_i) = 3w_i$. Moreover, linearizing the relation $\mu(v,v^\prime) = \sigma(v^\sharp,v^\prime)$, we conclude $\sigma(v_1 \times v_2,v_3) = \mu(v_1,v_2,v_3)$, and this is totally symmetric in $v_1,v_2,v_3 \in V$. 

(ii) $\Rightarrow$ (i). For $R \in \kalg$, we define a map $\mu_R\:V_R \to W_R$ by setting
\begin{align}
\label{MUAR} \mu_R(\sum_{i=1}^n r_ie_{iR}) := \sum_{i=1}^n r_i^3w_{iR} + \sum_{1\leq i,j\leq n,i\neq j}r_i^2r_j\sigma(e_i^\sharp,e_j)_R + \sum_{1\leq i<j<l\leq n} r_ir_jr_l\sigma(e_i \times e_j,e_l)_R 
\end{align}
for $r_1,\dots,r_n \in R$. Then $\mu := (\mu_R)_{R\in\kalg}\:V \to W$ is obviously a homogeneous polynomial law of degree $3$ satisfying $\mu(e_i) = w_i$ for $1 \leq i \leq n$. Now let $v = \sum_{i=1}^nr_ie_{iR}, v^\prime = \sum_{i=1}^nr_i^\prime e_{iR} \in V_R$ with $r_i,r_i^\prime \in R$ for $i = 1,\dots,n$. Then \eqref{MUAR} implies
\begin{align*}
\mu(v,v^\prime) =\,\,&\sum_i3r_i^2r_i^\prime w_{iR} + \sum_{i\neq j}(2r_ir_i^\prime r_j + r_i^2r_j^\prime)\sigma(e_i^\sharp,e_j)_R \\
\,\,&+ \sum_{i<j<l}(r_i^\prime r_jr_l + r_ir_j^\prime r_l + r_ir_jr_l^\prime)\sigma(e_i \times e_j,e_l)_R.
\end{align*}
On the other hand
\[
v^\sharp = \sum_i r_i^2e_{iR}^\sharp + \sum_{i<j}r_ir_je_{iR} \times e_{jR},
\]
and hence (ii) gives
\begin{align*}
\sigma(v^\sharp,v^\prime) =\,\,&\sum_{i,l}r_i^2r_l^\prime\sigma(e_i^\sharp,e_l)_R + \sum_{i,j,l,i<j}r_ir_jr_l^\prime\sigma(e_i \times e_j,e_l)_R \\
=\,\,&\sum_ir_i^2r_i^\prime\sigma(e_i^\sharp,e_i)_R + \sum_{i\neq j}r_i^2r_j^\prime\sigma(e_i^\sharp,e_j)_R + \sum_{i<j}r_ir_i^\prime r_j\sigma(e_i \times e_j,e_i)_R \\
\,\,&+ \sum_{i<j}r_ir_jr_j^\prime\sigma(e_i \times e_j,e_j)_R +\Big(\sum_{l<i<j} + \sum_{i<l<j} + \sum_{i<j<l}\Big)r_ir_jr_l^\prime\sigma(e_i \times e_j,e_l)_R \\
=\,\,&\sum_i3r_i^2r_i^\prime w_{iR} + \sum_{i\neq j}(r_i^2r_j^\prime + 2r_ir_i^\prime r_j)\sigma(e_i^\sharp,e_j)_R \\
\,\,&+ \sum_{i<j<l}(r_i^\prime r_jr_l + r_ir_j^\prime r_l + r_ir_jr_l^\prime)\sigma(e_i \times e_j,e_l)_R \\
=\,\,&\mu(v,v^\prime),
\end{align*}
and (i) holds. Finally, uniqueness of $\mu$ follows immediately from \eqref{EFLICO} of Exc.~\ref{pr.EXOCUB}.

(c) Since a pseudo cubic norm pseudo-structure $X$ over $k$ is automatically a cubic norm structure if $\frac{1}{3} \in k$ ($\lambda := \frac{1}{3}T$ has $\lambda(1) = 1$), any example of the kind demanded in (c) must display pathologies of characteristic $3$. With this in mind, we let $K/F$ be a purely inseparable field extension of characteristic $3$ and exponent $1$ and fix an element $a \in K\setminus F$. Then $F$ is infinite and $J_0 := K^{(+)}$ is a cubic Jordan algebra over $F$, with norm $N_0\:K \to F$ (by Cor.~\ref{c.POMALA} an honest-to-goodness polynomial function) given by $N_0(v)1 = v^3$ for all $v \in J_0$. For $R \in \Falg$ and $v,v^\prime \in J_{0R}$, we conclude $N_0(v,v^\prime)1 = 3v^2v^\prime = 0$, so the (bi-)linear and quadratic traces of $J_0$ are both zero, and $v^\sharp = v^2$ by (\ref{ss.BACUNO.fig}.\ref{ADJJA}) as well as $v \times v^\prime = 2vv^\prime = -vv^\prime$ for all $v,v^\prime \in J_0$. Now let $W$ be a vector space over $F$ of finite dimension at least $2$, pick linearly independent vectors $w_1,w_2 \in W$ and let $\lambda\:J_0 \to W$ be a non-zero linear map satisfying $\lambda(1) = 0$. Define a symmetric bilinear map $\sigma\:J_0 \times J_0 \to W$ by $\sigma(v,v^\prime) := \lambda(vv^\prime)$ for $v,v^\prime \in J_0$. Then $\sigma(v^\sharp,v) = \lambda(v^2v) = \lambda(v^3) = N_0(v)\lambda(1) = 0$ for all $v \in J_0$, and $\sigma(v_1 \times v_2,v_3) = -\lambda(v_1v_2v_3)$ is totally symmetric in $v_1,v_2,v_3 \in J_0$. Thus (b) yields a homogeneous polynomial law $\mu\:J_0 \to W$ of degree $3$ such that $\mu(1) = 0$, $\mu(a) = w_1$, $\mu(a^\sharp) = \mu(a^2) = w_2$.

Now let $k := F \times W \in \Falg$ be the ``split-null extension'' of $F$ by $W$, so $W \subseteq k$ is an ideal that squares to zero and $F$ acts canonically on $W$. In other words, the multiplication of $k$ obeys the rule
\[
(\alpha,w)(\alpha^\prime,w^\prime) := (\alpha\alpha^\prime,\alpha w^\prime + \alpha^\prime w) 
\] 
for $\alpha,\alpha^\prime \in F$, $w,w^\prime \in W$. By letting $W$ act trivially on $J_0$, we may and always will view $J_0$ as a Jordan algebra $J$ over $k$. 

Next define a symmetric $F$-bilinear map $T\:J \times J \to k$
\[
T(v,v^\prime) := \big(0,\sigma(v,v^\prime)\big) = \big(0,\lambda(vv^\prime)\big)
\]
for $v,v^\prime \in J$, which is in fact a $k$-bilinear form since
\[
T\big((\alpha,w)v,v^\prime\big) = T(\alpha v,v^\prime) = \alpha T(v,v^\prime) = (\alpha,w)\big(0,\lambda(vv^\prime)\big) = (\alpha,w)T(v,v^\prime)
\]
for all $\alpha \in F$, $w \in W$, $v,v^\prime \in J$. Similarly, we obtain a map
\[
N := N_0 \times \mu\:J \longrightarrow k, \quad v \longmapsto N(v) := \big(N_0(v),\mu(v)\big)
\]
and claim that \emph{the identity element of $J$ as base point, the (bilinearized) adjoint of $J$ as (bilinearized) adjoint, the symmetric bilinear form $T$, and the map $N$ satisfy equations \eqref{HOAD}--\eqref{ONESHB} of Exc.}~\ref{pr.RACUNO}. To show this, let $\alpha \in F$, $w \in W$ and $v,v^\prime \in J$. Then
\begin{align}
\label{SHARK} \big((\alpha,w)v\big)^\sharp = (\alpha v)^\sharp = \alpha^2v^\sharp = (\alpha,2\alpha w)v^\sharp = (\alpha,w)^2v^\sharp
\end{align}
proves \eqref{HOAD} and
\[
N\big((\alpha,w)v\big) = N(\alpha v) = \alpha^3\big(N_0(v),\mu(v)\big) = (\alpha,w)^3\big(N_0(v),\mu(v)\big) = (\alpha,w)^3N(v)
\]
proves \eqref{HONO}. The equations \eqref{SHARK} and
\begin{align*}
\big((\alpha,w)v\big) \times v^\prime =\,\,&(\alpha v) \times v^\prime = \alpha(v \times v^\prime) = (\alpha,w)(v \times v^\prime)
\end{align*}
show that the adjoint $v \mapsto v^\sharp$ is not only $F$-quadratic but, in fact, $k$-quadratic. Moreover, the bilinearized adjoint is the $k$-bilinearization of the adjoint, and we have \eqref{ABESH}. Next we compute
\begin{align*}
N(v + v^\prime) =\,\,&\big(N_0(v + v^\prime),\mu(v + v^\prime)\big) = \big(N_0(v) + N_0(v^\prime),\mu(v) + \mu(v,v^\prime) + \mu(v^\prime,v) + \mu(v^\prime)\big) \\
=\,\,&\big(N_0(v),\mu(v)\big) + \big(0, \sigma(v^\sharp,v^\prime)\big) + \big(0,\sigma(v^{\prime\sharp},v)\big) + \big(N_0(v^\prime),\mu(v^\prime)\big) \\
=\,\,&N(v) + T(v^ \sharp,v^\prime) + T(v^{\prime\sharp},v) + N(v^\prime),
\end{align*}
and \eqref{NOAB} holds. Since $3 = 0$ in $k$, we have $T(v^\sharp,v) = \big(0,\sigma(v^ \sharp, v)\big) = 0 = 3N(v)$, hence \eqref{TEASH}. Since the adjoint identity holds in $J_0$, we obtain
\[
v^{\sharp\sharp} = N_0(v)v = \big(N_0(v), \mu(v)\big)v = N(v)v,
\]
hence \eqref{ASH}. Similarly,
\begin{align*}
v^\sharp \times v^{\prime\sharp} + (v \times v^\prime)^\sharp =\,\,&-v^2v^{\prime 2} + (vv^\prime)^2 = 0 = \big(0,\lambda(v^\sharp v^\prime)\big)v^\prime + \big(0, \lambda(vv^{\prime\sharp})\big)v \\
=\,\,&T(v^\sharp,v^\prime)v^\prime + T(v,v^{\prime\sharp})v 
\end{align*}
gives \eqref{ASHBE} and
\begin{align*}
v^\sharp \times (v \times v^\prime) =\,\,&v^\sharp(vv^\prime) = (v^\sharp v)v^\prime = \big(0 \oplus \lambda(v^\sharp v^\prime)\big)v + N_0(v)v^\prime = \big(N_0(v),\mu(v)\big)v^\prime \\
=\,\,&T(v^\sharp,v^\prime)v + N(v)v^\prime  
\end{align*}
gives \eqref{ASHT}. While \eqref{ONESH} is obvious, the unit identity for $J_0$ implies $1 \times v = -v = \big(0,\lambda(v)\big)1 - v = T(1,v)1 - v$, hence \eqref{ONESHB}.

By Exc.~\ref{pr.RACUNO}, therefore, we find a cubic form $\tin\:X \to k$ making $X$ a cubic norm pseudo-structure over $k$. Now consider the element $a \in J$. By construction we have
\begin{align*}
N(a^\sharp) =\,\,&\big(N_0(a^2),\mu(a^2)\big) = \big(N_0(a)^2,w_2\big), \\ 
N(a)^2 =\,\,&\big(N_0(a),\mu(a)\big)^2 = \big(N_0(a)^2,2N_0(a)\mu(a)\big) = \big(N_0(a)^2,- N_0(a)w_1\big),
\end{align*}
hence $N(a^\sharp) \neq N(a)^2$ since $w_1,w_2$ are linearly independent over $F$.
\end{sol}

\solnsec{Section~\ref{s.BUCUNO}}

\begin{sol}{pr.EICO} \label{sol.EICO}  (i) $k$ is a cubic norm structure with base point $1$, adjoint $\alpha \mapsto \alpha^2$, norm $r \mapsto r^3$, and bilinear (resp. quadratic) trace respectively given by $(\alpha,\beta) \mapsto 3\alpha\beta$ (resp. $\alpha \mapsto 3\alpha^2$). Scalar multiplication gives a bilinear action $k \times V \to V$, while the eikonal triple structure provides us with quadratic maps $Q\:V \to k$, $H\:V \to V$. Moreover, setting $v = u$ in \eqref{HADJ}, Euler's differential equation implies
\begin{align}
\label{QUHAEN} Q\big(u,H(u)\big) = N(u) &&(R \in \kalg,\;u \in V_R).
\end{align}
By \ref{ss.BUP}, therefore, formulas \eqref{EICBA}--\eqref{EICNO} define a cubic array $X$ over $k$ whose bilinear trace by (\ref{ss.BUP}.\ref{BITREX}) is regular. Thus the solution to (i) will be complete once we have shown that $X$ is, in fact, a cubic norm structure. In order to do so, it suffices to show, by Prop.~\ref{p.EXCUNO}, that the identities (\ref{p.EXCUNO.fig}.\ref{CUNISH})--(\ref{p.EXCUNO.fig}.\ref{CHXVT}) hold strictly in $X$. Here (\ref{p.EXCUNO.fig}.\ref{CQHV}) is just the eikonal equation, while (\ref{p.EXCUNO.fig}.\ref{CQXVH}) has been established in \eqref{QUHAEN} with $\han = N$.  Since (\ref{p.EXCUNO.fig}.\ref{CUNISH})--(\ref{p.EXCUNO.fig}.\ref{CHXV}), (\ref{p.EXCUNO.fig}.\ref{CTXQV}) are obvious in the special case at hand, we are left with (\ref{p.EXCUNO.fig}.\ref{CHHV}), (\ref{p.EXCUNO.fig}.\ref{CQHVWT}), (\ref{p.EXCUNO.fig}.\ref{CHXVT}), i.e., with
\begin{align}
\label{HAHAEN} H\big(H(u)\big) =\,\,&N(u)u - Q(u)H(u), \\
\label{QUHAQU} Q\big(u,H(u,v)\big) =\,\,&2Q\big(H(u),v\big), \\
\label{HAHAQU} H\big(u,H(u)\big) =\,\,&2Q(u)u.
\end{align}
We begin with \eqref{HAHAQU} by linearizing \eqref{HADJ} to conclude that
\begin{align}
\label{LIHADJ} Q\big(H(u,v),w\big) = \frac{1}{3}N(u,v,w)
\end{align}
is totally symmetric in $u,v,w \in V$. Combining this with the linearization of \eqref{EICEQ}, i.e., with
\begin{align}
\label{QUHAHA} Q\big(H(u),H(u,v)\big) = 2Q(u)Q(u,v),
\end{align}
we conclude $Q(H(u,H(u)),v) = Q(H(u,v),H(u)) =2Q(u)Q(u,v) = Q(2Q(u)u,v)$, and \eqref{HAHAQU} drops out since $Q$ is regular. Setting $w = u$ and using symmetry of \eqref{LIHADJ} again, we obtain $Q(H(u,v),u) = Q(H(u,u),v) = 2Q(H(u),v)$, hence \eqref{QUHAQU}. And finally, linearizing \eqref{HAHAQU}, we deduce
\[
H\big(v,H(u)\big) + H\big(u,H(u,v)\big) = 2Q(u,v)u + 2Q(u)v.
\]
Here we put $v = H(u)$ apply \eqref{HAHAQU}, \eqref{QUHAEN} and obtain
\begin{align*}
2H\big(H(u)\big) + 4Q(u)H(u) =\,\,&H\big(H(u),H(u)\big) + 2H\big(u,Q(u)u\big) \\
=\,\,& H\big(H(u),H(u)\big) + H\Big(u,H\big(u,H(u)\big)\Big) = 2Q\big(u,H(u)\big)u + 2Q(u)H(u) \\
=\,\,&2N(u)u + 2Q(u)H(u),
\end{align*}
hence \eqref{HAHAEN}.

(ii) The cubic norm substructure $X_0 \subseteq X$ is regular since $\frac{1}{3} \in k$. Hence $(X_0,V)$ is a complemented cubic norm substructure of $X$. Defining a quadratic map $H\:V \to V$ by $H(u) := Q(u) + u^\sharp$ for $u \in V$, we are in the situation \ref{ss.SPLAD} and conclude that the identities (\ref{ss.SPLAD}.\ref{QUH})--(\ref{ss.SPLAD}.\ref{SPQATR}), (\ref{ss.IDQH.fig}.\ref{QXV})--(\ref{ss.IDQH.fig}.\ref{NHV}) hold strictly in $X$. In particular, by (\ref{ss.SPLAD}.\ref{SPLIBILT}) we have $3Q(u,v) = T(u,v)$, whence $(V,Q)$ is a quadratic space over $k$. Moreover, for $u,v \in V$ the gradient identity implies $3Q(H(u),v) = T(H(u),v) = T(u^\sharp,v) = N(u,v)$, so the quadratic map $H$ is connected with the cubic form $N$ on $V$ by \eqref{HADJ}. Finally, (\ref{ss.IDQH.fig}.\ref{QHV}) shows that the eikonal equation \eqref{EICEQ} holds strictly in $V$. Summing up we have proved that $(V,Q,N\vert_V)$ is an eikonal triple over $k$.

(iii) follows by a straightforward verification.

(iv) The gradient of $N$ at $u \in \IR^n$ is defined by
\[
\big(\grad(N)\big)(u) = \left(\begin{matrix}
\frac{\partial N}{\partial x_1}(u) \\
\vdots \\
\frac{\partial N}{\partial x_n}(u)
\end{matrix}
\right) \in \IR^n
\]
and can be characterized by
\begin{align*}
\big(\grad(N)\big)(u)^\trans v = N(u,v) &&(v \in \IR^n).
\end{align*}
This and \eqref{HADJ} imply
\[
H(u)^\trans Qv =Q\big(H(u),v\big) = \frac{1}{3}\big(\grad(N)\big)(u)^\trans v,
\]
and we conclude $H(u)^\trans Q = \frac{1}{3}(\grad(N))(u)^\trans$, hence
\[
H(u)^\trans  = \frac{1}{3}\big(\grad(N)\big)(u)^\trans Q^{-1}, \quad H(u) = \frac{1}{3}Q^{-1}\big(\grad(N)\big)(u).  
\]
Now the left-hand side of \eqref{EICEQ} becomes
\begin{align*}
H(u)^\trans QH(u) =\,\,&\frac{1}{9}\big(\grad(N)\big)(u)^\trans Q^{-1}QQ^{-1}\big(\grad(N)\big)(u) = \sum_{i,j=1}^n q^{ij}\frac{\partial N}{\partial x_i}(u)\frac{\partial N}{\partial x_j}(u),
\end{align*}
while the right-hand side of \eqref{HADJ} is
\[
\Big(\sum_{i,j=1}^n q_{ij}u_iu_j\Big)^2.
\]
Hence \eqref{EICEQ} follows.
\end{sol}


\solnsec{Section~\ref{s.CUJOMA}}

\begin{sol}{pr.HERTWO} \label{sol.HERTWO} Write $\Gamma = \diag(\gamma_1,\gamma_2)$, $\gamma_1,\gamma_2 \in k^\times$, and put
\begin{align}
\label{UONTWO} u[12] := \gamma_2ue_{12} + \gamma_1\bar ue_{21} &&(u \in C)   
\end{align}
in terms of the usual matrix units in $\Mat_2(k)$. One checks that
\begin{align}
\label{EMHERT} M := \Her_2(C,\Gamma) = ke_{11} \oplus ke_{22} \oplus C[12], \quad C[12] := \{u[12] \mid u \in C\}.    
\end{align}
We define a quadratic form $q\:M \to k$ by
\[
q(\xi_1e_{11} + \xi_2e_{22} + u[12]) := \xi_1\xi_2 - \gamma_1\gamma_2n_C(u)
\]
for $\xi_1,\xi_2 \in k$, $u \in C$, which bilinearizes to
\begin{align}
 \label{QUXIET} q(\xi_1e_{11} + \xi_2e_{22} + u[12],\eta_1e_{11} + \eta_2e_{22} + v[12]) = \xi_1\eta_2 + \xi_2\eta_1 - \gamma_1\gamma_2n_C(u,v)   
\end{align}
for $\xi_i,\eta_i \in k$, $i = 1,2$, $u,v \in C$. Then $(M,q,e)$, $e := \Eins_2 = e_{11} + e_{22}$, is a pointed quadratic module over $k$ whose linear trace by \eqref{QUXIET} is given by
\[
t\:M \longrightarrow k, \quad t(\xi_1e_{11} + \xi_2e_{22} + u[12]) = \xi_1 + \xi_2
\]
for $\xi_1,\xi_2 \in k$, $u \in C$. It follows that $J := J(M,q,e)$ is a Jordan algebra of Clifford type over $k$ and, writing $x^n$, $n \in \IN$, for the $n$-th power of $x \in J$, a straightforward verification shows
\begin{align}
 \label{POMAPR} xx = x^2 = t(x)x - q(x)e,   
\end{align}
where the expression on the very left is the matrix product of $x$ with itself. Linearizing shows that the circle peoduct of $x,y \in J$ is $x \circ y = xy + yx$. We now claim
\begin{align}
\label{EXEXCU} (xx)x = x^3 = x(xx). 
\end{align}
Indeed, \eqref{POMAPR} implies $(xx)x = t(x)xx - q(x)ex = t(x)x^2 - q(x)x = x^3$ (by (\ref{ss.JOPOID.fig}.\ref{JOQUAT})), hence the first equation of \eqref{EXEXCU}; the second one follows analogously. Actually, \eqref{EXEXCU} holds strictly, so we are allowed to differentiate at $x$ in the direction $y$. Doing so with $x^3 = U_xx$ gives $U_xy + U_{x,y}x = U_xy + x^2 \circ y$ (by (\ref{ss.JABAS.fig}.\ref{UXYX})), while the same procedure for $(xx)x$ and \eqref{POMAPR} give
\[
(yx)x + (xy)x + (xx)y = (xy)x + x^2 \circ y + [y,x,x],
\]
hence the first equation of \eqref{UHEMAT}. The same argument applied to $x(xx)$ yields the second equation of \eqref{UHEMAT}.

Finally, assume $C$ is projective as a $k$-module and let $(\pM,\pq,\pe)$ be any pointed quadratic module over $k$ satisfying the conditions of the exercise. Then $\pM = \Her_2(C,\Gamma) = M$ and $\pe = \Eins_2 = e$. Moreover, \eqref{UHEMAT} shows that $J = J(M,q,e)$ and $\pJ = J(\pM,\pq,\pe)$ have the same $U$-operator, i.e., $\Eins_M\:J \to J$ is an isomorphism of Jordan algebras. Hence $\pq = \pq \circ \Eins_M = q$ by Exc.~\ref{pr.POJOCAT}, as desired.
\end{sol}

\begin{sol}{pr.UOPMA} \label{sol.UOPMA}  We begin by showing \eqref{SQUAMA}. To this end, we put
\begin{align}
\label{EXSQUA} x^2 =\,\,&\sum(\eta_ie_{ii} + v_i[jl]) &&(\eta_i \in k,\;v_i \in C,\;1 \le i \le3), \\
\label{UOPMA.EXTIEX} xx =\,\,&(v_{ij}) &&(v_{ij} \in C,\;1 \le i,j \le 3).
\end{align}
From (\ref{ss.BACUNO.fig}.\ref{ADJJA}) we deduce $x^2 = x^\sharp + T(x)x - S(x)\Eins_3$, which combines with (\ref{ss.THERCU}.\ref{MADJ}), (\ref{ss.THERCU}.\ref{MALIT}), (\ref{ss.THERCU}.\ref{MAQIT}) to yield
\begin{align*}
\eta_i =\,\,&\xi_j\xi_l - \gamma_j\gamma_ln_C(u_i) + (\xi_i + \xi_j + \xi_l)\xi_i - \big(\xi_j\xi_l - \gamma_j\gamma_ln_C(u_i)\big) \\
\,\,&-\big(\xi_l\xi_i - \gamma_l\gamma_in_C(u_j)\big) - \big(\xi_i\xi_j - \gamma_i\gamma_jn_C(u_l)\big)
\end{align*}
and then simplifies to
\begin{align}
\label{ETIEXP} \eta_i = \xi_i^2 + \gamma_l\gamma_in_C(u_j) + \gamma_i\gamma_jn_C(u_l).
\end{align}
Even simpler,
\[
v_i = -\xi_iu_i + \gamma_i\overline{u_ju_l} + (\xi_i + \xi_j + \xi_l)u_i,
\]
hence
\begin{align}
\label{VEIEXP} v_i = (\xi_j + \xi_l)u_i + \gamma_i\overline{u_ju_l}.    
\end{align}
On the other hand, inspecting (\ref{ss.COCUJM}.\ref{ELCUJM})and consulting \eqref{ETIEXP}, we obtain
\[
v_{11} = \xi_1^2 +\gamma_1\gamma_2n_C(u_3) + \gamma_3\gamma_1n_C(u_2) = \eta_1
\] and similarly $v_{ii} = \eta_i$ for $i = 1,2,3$. Moreover, by \eqref{VEIEXP}, 
\begin{align*}
v_{12} =\,\,&\gamma_2\xi_1u_3 + \gamma_2\xi_2u_3 + \gamma_2\gamma_3\bar u_2\bar u_1 = \gamma_2\big((\xi_1 + \xi_2)u_3 + \gamma_3\overline{u_1u_2}\big) = \gamma_2v_3, \\
v_{21} =\,\,&\gamma_1\xi_1\bar u_3 + \gamma_1\xi_2\bar u_3 + \gamma_3\gamma_1u_1u_2 = \gamma_1\big((\xi_1 + \xi_2)\bar u_3 + \gamma_3u_1u_2\big) = \gamma_1\bar v_3
\end{align*}
and, similarly, $v_{jl} = \gamma_lv_i$, $v_{lj} = \gamma_j\bar v_i$. This completes the proof of \eqref{SQUAMA}.

Next we deal with the first equation of \eqref{CUBMA} by setting
\begin{align}
\label{EXCUBE} x^3 =\,\,&\sum(\zeta_ie_{ii} + w_i[jl]) &&(\zeta_i \in k,\; w_i \in C,\;1 \le i \le 3), \\
\label{EXEXSQ} xx^2 =\,\,&(w_{ij}) &&(w_{ij} \in C,\;1 \le i,j \le 3). 
\end{align}
From (\ref{t.CUNOJO}.\ref{UNIVT}) we deduce $x^3 = T(x)x^2 - S(x)x + N(x)\Eins_3$, and \eqref{ETIEXP} implies
\begin{align*}
\zeta_i =\,\,&(\xi_i + \xi_j + \xi_l)\big(\xi_i^2 + \gamma_l\gamma_in_C(u_j) + \gamma_i\gamma_jn_C(u_l)\big) \\
\,\,&-\Big(\big(\xi_j\xi_l - \gamma_j\gamma_ln_C(u_i)\big) + \big(\xi_l\xi_i - \gamma_l\gamma_in_C(u_j)\big) + \big(\xi_i\xi_j - \gamma_i\gamma_jn_C(u_l)\big)\Big)\xi_i \\
\,\,&+\xi_i\xi_j\xi_l - \gamma_j\gamma_l\xi_in_C(u_i) - \gamma_l\gamma_i\xi_jn_C(u_j) - \gamma_i\gamma_j\xi_ln_C(u_l) + \gamma_1\gamma_2\gamma_3t_C(u_1u_2u_3).
\end{align*}
Expanding and collecting terms gives
\begin{align}
\label{ZETIEXP} \zeta_i =\,\,&\xi_i^3 + \gamma_l\gamma_i(\xi_l + 2\xi_i)n_C(u_j) + \gamma_i\gamma_j(2\xi_i + \xi_j)n_C(u_l) + \gamma_1\gamma_2\gamma_3t_C(u_1u_2u_3).    
\end{align}
Much simpler,
\begin{align*}
w_i =\,\,&(\xi_i + \xi_j + \xi_l)\big((\xi_j + \xi_l)u_i + \gamma_i\overline{u_ju_l}\big) \\
\,\,&-\big(\xi_j\xi_l - \gamma_j\gamma_ln_C(u_i) + \xi_l\xi_i - \gamma_l\gamma_in_C(u_j) + \xi_i\xi_j - \gamma_i\gamma_jn_C(u_l)\big)u_i,
\end{align*}
which after a short computation reduces to 
\begin{align}
\label{WEIEXP} w_i =\,\,&\big(\xi_j^2 + \xi_j\xi_l + \xi_l^2 + \gamma_j\gamma_ln_C(u_i) \\
\,\,&+ \gamma_l\gamma_in_C(u_j) + \gamma_i\gamma_jn_C(u_l)\big)u_i + \gamma_i(\xi_1 + \xi_2 + \xi_3)\overline{u_ju_l}. \notag
\end{align}
Again inspecting (\ref{ss.COCUJM}.\ref{ELCUJM}) and using \eqref{ETIEXP}, \eqref{VEIEXP}, we now compute
\begin{align*}
w_{11} =\,\,&\xi_1\big(\xi_1^2 + \gamma_3\gamma_1n_C(u_2) + \gamma_1\gamma_2n_C(u_3)\big) + \gamma_1\gamma_2u_3\big((\xi_1 + \xi_2)\bar u_3 + \gamma_3u_1u_2\big) \\
\,\,&+ \gamma_3\gamma_1\bar u_2\big((\xi_3 + \xi_1)u_2 + \gamma_2\overline{u_3u_1}\big),
\end{align*}
which by \eqref{ZETIEXP} simplifies to
\begin{align*}
w_{11} =\,\,&\zeta_1 + \gamma_1\gamma_2\gamma_3\big(u_3(u_1u_2) + \overline{(u_3u_1)u_2} - (u_3u_1)u_2 - \overline{(u_3u_1)u_2}\big) \\
=\,\,&\zeta_1 - \gamma_1\gamma_2\gamma_3[u_1,u_2,u_3]
\end{align*}
 and, similarly, $w_{ii} = \zeta_i - \gamma_1\gamma_2\gamma_3[u_1,u_2,u_3]$ for $i = 1,2,3$. On the other hand,
\begin{align*}
w_{23} =\,\,&\gamma_3\gamma_1\bar u_3\big((\xi_3 + \xi_1)\bar u_2 + \gamma_2u_3u_1\big) + \gamma_3\xi_2\big((\xi_2 + \xi_3)u_1 + \gamma_1\overline{u_2u_3}\big) \\
\,\,&+ \gamma_3\big(\xi_3^2 + \gamma_2\gamma_3n_C(u_1) + \gamma_3\gamma_1n_C(u_2)\big)u_1,
\end{align*}
which by \eqref{WEIEXP} simplifies to $w_{23} = \gamma_3w_1$. Similarly,
\begin{align*}
w_{32} =\,\,&\gamma_1\gamma_2u_2\big((\xi_1 + \xi_2)u_3 + \gamma_3\overline{u_1u_2}\big) + \gamma_2\big(\xi_2^2 + \gamma_3\gamma_1n_C(u_2) + \gamma_1\gamma_2n_C(u_3)\big)\bar u_1 \\
\,\,&+\gamma_2\xi_3\big((\xi_2 + \xi_3)\bar u_1 + \gamma_1u_2u_3\big) = \gamma_2\bar w_1.
\end{align*}
The same computation yields, more generally, $w_{jl} = \gamma_lw_i$, $w_{lj} = \gamma_j\bar w_i$ which completes the proof of the first part of \eqref{CUBMA}.

Equation \eqref{SQUAMA} amounts to saying that the squaring of $\Mat_3(C)$ restricts to the squaring of $J$. After linearizing, therefore, the symmetric matrix product $xy + yx$ of $\Mat_3(C)$ restricts to the circle product $x \circ y$ of $J$. Hence (\ref{p.JOAPOA}.\ref{TRIXMN}), \eqref{SQUAMA} and the first equation of \eqref{CUBMA} imply
\[
2x^3 = x \circ x^2 = x(xx) + (xx)x = x^3 - \gamma_1\gamma_2\gamma_3[u_1,u_2,u_3] + (xx)x,
\]
which yields the second equation of \eqref{CUBMA}. In order to derive \eqref{UOPMA}, we note that \eqref{CUBMA} holds strictly, so we are allowed to differentiate it at $x$ in the direction $y$. Doing so for the left-hand side,i.e., for $x^3 = U_xx$, we apply (\ref{ss.JABAS.fig}.\ref{UXYX}) and obtain $U_xy + U_{x,y}x = U_xy + x^2 \circ y = U_xy + (xx)y + yx^2$. Since the associator of $C$ is alternating, the same procedure applied to the middle term of \eqref{CUBMA} yields the expression
\[
yx^2 + x(yx) + x(xy)+ \gamma_1\gamma_2\gamma_3\big([u_1,u_2,v_3] + [u_2,u_3,v_1] + [u_3,u_1,v_2]\big)\Eins_3,
\]
and comparing, we end up with the first formula of \eqref{UOPMA}. The second one follows by the same argument, using the second formula of \eqref{CUBMA} instead of the first. Alternatively, we may invoke (\ref{ss.JABAS.fig}.\ref{TWOUX}) to conclude
\begin{align*}
2U_xy =\,\,&x \circ (x \circ y) - x^2 \circ y = x(xy + yx) + (xy + yx)x - (xx)y - y(xx) \\
=\,\,&x(yx) + (xy)x - [x,x,y] + [y,x,x],  
\end{align*}
so the second formula of \eqref{UOPMA} follows from this and the first. Subtracting the middle term of \eqref{CUBMA} from the right one, we obtain \eqref{ASXXX}. Invoking (\ref{ss.BACUNO.fig}.\ref{ADJJA}), (\ref{t.CUNOJO}.\ref{UNIVT}) and \eqref{SQUAMA}, \eqref{CUBMA}, we deduce
\begin{align*}
x^\sharp x =\,\,&(xx)x - T(x)xx + S(x)x = x^3 - T(x)x^2 + S(x)x + \gamma_1\gamma_2\gamma_3[u_1,u_2,u_3]\Eins_3 \\
=\,\,&\big(N(x)1_C + \gamma_1\gamma_2\gamma_3[u_1,u_2,u_3]1_C\big)\Eins_3,
\end{align*}
hence \eqref{SHAMA}, while \eqref{MASHA} follows from \eqref{SHAMA} and $xx^\sharp + x^\sharp x = x \circ x^\sharp = 2N(x)\Eins_3$. Finally, turning to \eqref{USHAI}, we replace $x$ by $x^\sharp$ in \eqref{MASHA} and apply (\ref{ss.BACUNO.fig}.\ref{NADJ}) to conclude
\[
x^\sharp x^{\sharp\sharp} = \big(N(x)^21_C - \gamma_1\gamma_2\gamma_3[u_1^\sharp,u_2^\sharp,u_3^\sharp]\big)\Eins_3.
\]
On the other hand, the adjoint identity and \eqref{SHAMA} yield
\[
x^\sharp x^{\sharp\sharp} = N(x)x^\sharp x = \big(N(x)^21_C + \gamma_1\gamma_2\gamma_3N(x)[u_1,u_2,u_3]\big)\Eins_3,
\]
hence \eqref{USHAI}.
\end{sol}


\solnsec{Section~\ref{s.ELIDMA}}

\begin{sol}{pr.IDCUJO} \label{sol.IDCUJO} (a) We provide the solution to Exc.~\ref{pr.DECID} with the following supplement. Let $f\:M \to N$ be a polynomial law over $k$, $\vep \in k$ an idempotent and $x \in M$. Then the identifications 
\begin{align*}
M_R = \vep M, \quad x_R = x \otimes \vep = \vep x &&(R = \vep k \in \kalg,\;\;x \in M),  
\end{align*}
ditto for $N$, of \ref{ss.PROPID} imply $f_R(\vep x) = f_R(x_R) = f_k(x)_R = \vep f_k(x)$, hence
\begin{align}
\label{EFAREP} f_R(\vep x) = \vep f_k(x) &&(x \in M).
\end{align}

\smallskip

(b) Before dealing with the problem itself, we characterize the various types of idempotents $e \in J$ in the following way.
\begin{align}
\label{CHAZER} e = 0 \quad\,\,&\Longleftrightarrow \quad T(e) = S(e) = N(e) = 0, \\
\label{CHAELM} e \;\; \text{is elementary} \quad\,\,&\Longleftrightarrow \quad T(e) = 1,\;\;S(e) = N(e) = 0 \\ 
\,\,&\Longleftrightarrow \quad T(e) = 1,\;\;S(e) = 0, \notag \\
\label{CHACOL} e \;\; \text{is co-elementary} \quad\,\,&\Longleftrightarrow \quad T(e) = 2,\;\;S(e) = 1,\;\;N(e) = 0 \\
\,\,&\Longleftrightarrow \quad T(e) = 2,\;\;S(e) = 1, \notag \\
\label{CHAUN} e = 1 \quad\,\,&\Longleftrightarrow \quad T(e) = S(e) = 3,\;\;N(e) = 1 \quad \Longleftrightarrow \quad N(e) = 1.
\end{align}
In \eqref{CHAZER}, the implication from left to right is obvious. Conversely, assume $T(e) = S(e) = N(e) = 0$. Then (\ref{t.CUNOJO}.\ref{UNIVT}) implies $e = e^3 = 0$. If $e$ is elementary, then $T(e) = 1$, $e^\sharp = 0$, hence $S(e) = T(e^\sharp) = 0$, while $N(e) = 0$ has been noted in \ref{ss.COELID}. Conversely, assume $T(e) = 1$, $S(e) = 0$. Then (\ref{ss.BACUNO.fig}.\ref{ADJJA}) yields $e^\sharp = e^2 - e = 0$, and $e$ is elementary, proving \eqref{CHAELM}. In \eqref{CHACOL}, we have
\begin{align*}
e \;\; \text{is co-elementary} \quad\,\,&\Longleftrightarrow \quad 1 - e \;\; \text{is elementary} \\ 
\quad\,\,&\Longleftrightarrow \quad T(1 - e) = 1,\;\;S(1 - e) = N(1 - e) = 0.
\end{align*}
But $T(1 - e) = 3 - T(e)$ and $S(1 - e) = S(1) - S(1,e) + S(e) = 3 - 2T(e) + S(e)$ by (\ref{ss.BACUNO.fig}.\ref{TSONE}), while $N(1 - e) = 1 - T(e) + S(e) - N(e)$ by (\ref{ss.BACUNO.fig}.\ref{EXPF}), (\ref{ss.BACUNO.fig}.\ref{BAPO}). This proves \eqref{CHACOL}. Finally, if $e = 1$ in \eqref{CHAUN}, then $T(e) = S(e) = 3$ and $N(e) = 1$ by (\ref{ss.BACUNO.fig}.\ref{BAPO}), (\ref{ss.BACUNO.fig}.\ref{TSONE}). Conversely, if $N(e) = 1$, then $e \in J^\times$ by Cor.~\ref{c.REGINV}. Hence $U_e$ is a bijective projection by Prop.~\ref{p.CHAIN.JOAL} and Thm.~\ref{t.PEDESI}, which implies $U_e = \Eins_J$, hence $e = e^2 = U_e1 = 1$. 

\smallskip

(c) Turning, finally, to the solution of the problem, we begin by showing uniqueness, so let $(\vep^{(i)})_{0\leq i\leq 3}$ be a complete orthogonal system of idempotents in $k$ with the desired properties. For the sake of clarity, we write $T^{(i)}, S^{(i)}, N^{(i)}$ for the trace, quadratic trace, norm, respectively, of $J^{(i)}$ over $k^{(i)}$. Invoking \eqref{EFAREP}--\eqref{CHAUN}, we obtain
\begin{align*}
T(e) =\,\,&\big(\vep^{(i)}T(e)\big)_{0\le i\le 3} = \big(T^{(i)}(\vep^{(i)}e)\big)_{0\le i\le 3} = \big(T^{(i)}(e^{(i)})\big)_{0\le i\le 3} \\
=\,\,&(0,\vep^{(1)},2\vep^{(2)},3\vep^{(3)})
\end{align*} 
and, similarly,
\[
S(e) = \big(S^{(i)}(e^{(i)})\big)_{0\le i\le 3} = (0,0,\vep^{(2)},3\vep^{(3)}), \quad N(e) = \big(N^{(i)}(e^{(i)})\big)_{0\le i\le 3} = (0,0,0,\vep^{(3)}).
\] 
Solving this system of linear equations in the unknowns $\vep^{(1)}, \vep^{(2)}, \vep^{(3)}$, we obtain \eqref{DECIONE}--\eqref{DECITHREE}. Now \eqref{DECIZERO} follows from the completeness of the orthogonal system $(\vep^{(i)})$. 

\smallskip

(d) In order to prove existence, we define the $\vep^{(i)}$, $0 \leq i \leq 3$ by \eqref{DECIZERO}--\eqref{DECITHREE} and must show that they have the desired properties. First of all, they clearly add up to $1$. Next we have to show that they are orthogonal idempotents, equivalently, that $\vep^{(i)}\vep^{(j)} = 0$ for $0 \leq i,j \leq 3$, $i \neq j$. To begin with, since $N$ preserves powers by (\ref{ss.BACUNO.fig}.\ref{JCOMP}), and $e,1 - e$ are idempotents, so are $\vep^{(0)}, \vep^{(3)}$. Moreover, $\vep^{(0)}\vep^{(3)} = N(e)^2N(1 - e) = N(U_e(1 - e)) = 0$. Summing up, we have proved
\begin{align}
\label{VEPZON} \vep^{(0)2} = \vep^{(0)}, \quad \vep^{(3)2} = \vep^{(3)}, \quad \vep^{(0)}\vep^{(3)} = 0.
\end{align}
Next we apply (\ref{t.CUNOJO}.\ref{UNIVF}) and obtain $e = e^4 = (T(e) - S(e) + N(e))e$, hence $N(1 - e)e = 0$. Since $N(1 - e) = \vep^{(3)}$ by \eqref{VEPZON} is an idempotent, applying $T$ and $S$ to the preceding equation, we end up with
\begin{align}
\label{ENONEE} N(1 - e)e = 0, \quad T(e)N(1 - e) = S(e)N(1 - e) = 0.
\end{align}
Combining with \eqref{DECIONE}, \eqref{DECITWO}, \eqref{VEPZON}, we conclude
\begin{align}
\label{VEZEON} \vep^{(0)}\vep^{(1)} = \vep^{(0)}\vep^{(2)} = 0.
\end{align}
On the other hand, replacing $e$ by $1 - e$ in \eqref{ENONEE} yields $0 = T(1 - e)N(e) = 3N(e) - T(e)N(e)$ and $0 = S(1 - e)N(e) = (3 - 2T(e) + S(e))N(e) = 3N(e) - 6N(e) + S(e)N(e)$, hence
\begin{align}
\label{TEENEE} T(e)N(e) = S(e)N(e) = 3N(e).
\end{align}
This implies $\vep^{(1)}\vep^{(3)} = 3N(e) - 6N(e) + 3N(e) = 0$, $\vep^{(2)}\vep^{(3)} = 3N(e) - 3N(e) = 0$, and we have
\begin{align}
\label{VEONTH} \vep^{(1)}\vep^{(3)} = \vep^{(2)}\vep^{(3)} = 0.
\end{align}
In view of \eqref{VEPZON}, \eqref{VEZEON}, \eqref{VEONTH}, the $\vep^{(i)}$ will form a complete orthogonal system of idempotents in $k$ once we have shown $\vep^{(1)}\vep^{(2)} = 0$. To this end, we first note $T(e) = T(e,e) + T(1 - e,U_ee) = T(e,e) + T(U_e(1 - e),e) = T(e,e)$ by (\ref{ss.BACUNO.fig}.\ref{TRAU}) and then apply (\ref{ss.BACUNO.fig}.\ref{BQUAT}) to conclude $2S(e) = S(e,e) = T(e)^2 - T(e,e)$, hence
\begin{align}
\label{TWOESE} 2S(e) = T(e)^2 - T(e).
\end{align}
Now we expand the last two equations of \eqref{ENONEE} by using \eqref{TEENEE} and \eqref{TWOESE}. We obtain $0 = T(e) - T(e)^2 + T(e)S(e) -3N(e)$, hence
\begin{align}
\label{TEEESE} T(e)S(e) = 2S(e) + 3N(e),
\end{align}
and $0 = S(e) -2S(e) - 3N(e) + S(e)^2 - 3N(e)$, hence
\begin{align}
\label{ESESQ} S(e)^2 = S(e) + 6N(e).
\end{align}
Using \eqref{TEENEE}, \eqref{TEEESE}, \eqref{ESESQ}, we can now compute
\begin{align*}
\vep^{(1)}\vep^{(2)} =\,\,&T(e)S(e) - 3T(e)N(e) -2S(e)^2 + 6S(e)N(e) + 3S(e)N(e) - 9N(e) \\
=\,\,&2S(e) + 3N(e) - 9N(e) - 2S(e) - 12N(e) + 27N(e) - 9N(e)  \\
=\,\,& 0.
\end{align*}
Thus $(\vep^{(i)})_{0\leq i\leq 3}$ is a complete orthogonal system of idempotents in $k$. In particular, we have $e = (e^{(i)})_{0\le i\le 3}$, $e^{(i)} = \vep^{(i)}e \in J^{(i)}$ for $0 \leq i \leq 3$. From \eqref{ENONEE} we deduce $e^{(0)} = N(1 - e)e = 0$, while \eqref{EFAREP} implies $N^{(3)}(e^{(3)}) = N^{(3)}(\vep^{(3)}e) = \vep^{(3)}N(e) =\vep^{(3)}$, hence $e^{(3)} = 1_{J^{(3)}}$ by \eqref{CHAUN}. It remains to show that $e^{(1)} \in J^{(1)}$ is elementary and $e^{(2)} \in J^{(2)}$ is co-elementary, which follows from 
\begin{align*}
T^{(1)}(e^{(1)}) =\,\,&T^{(1)}(\vep^{(1)}e) = \vep^{(1)}T(e) = T(e)^2 - 2T(e)S(e) + 3T(e)N(e) \\
=\,\,&T(e)^2 - 4S(e) -6N(e) + 9N(e) = T(e) -2S(e) + 3N(e) = \vep^{(1)} = 1_{k^{(1)}}, \\
S^{(1)}(e^{(1)}) =\,\,&S^{(1)}(\vep^{(1)}e) = \vep^{(1)}S(e) = T(e)S(e) - 2S(e)^2 + 3S(e)N(e) \\
=\,\,&2S(e) + 3N(e) -2S(e) -12N(e) + 9N(e) = 0
\end{align*}
and \eqref{CHAELM} in the first case, and from
\begin{align*}
T^{(2)}(e^{(2)}) =\,\,&T^{(2)}(\vep^{(2)}e) = \vep^{(2)}T(e) = T(e)S(e) - 3T(e)N(e) \\
=\,\,&2S(e) + 3N(e) - 9N(e) = 2\big(S(e) - 3N(e)\big) = 2\vep^{(2)} = 2\cdot 1_{k^{(2)}}, \\
S^{(2)}(e^{(2)}) =\,\,&S^{(2)}(\vep^{(2)}e) = \vep^{(2)}S(e) = S(e)^2 - 3S(e)N(e) = S(e) + 6N(e) - 9N(e) \\
=\,\,&S(e) - 3N(e) = \vep^{(2)} = 1_{k^{(2)}}
\end{align*}
and \eqref{CHACOL} in the second.
\end{sol}

\begin{sol}{pr.FERLE} \label{sol.FERLE}  Applying equation \eqref{EFLICO} of Exc.~\ref{pr.EXOCUB} and the (bilinearized) gradient identity, we expand $N(q) = N(\sum u_i) = \sum N(u_i) + \sum_{i\neq j}T(u_i^\sharp,u_j) + T(u_1 \times u_2,u_3) = \sum N(u_i) + T(u_1 \times u_2,u_3)$, where (\ref{ss.BACUNO.fig}.\ref{NADJ}) implies $N(u_i)^2 = N(u_i^\sharp) = 0$. Since, therefore, the scalars $N(q)$ and $T(u_1 \times u_2,u_3)$ differ by a nilpotent element in $k$, one of them is invertible if and only if so is the other. This proves the first assertion. Now assume $q \in J^\times$ and $T(u_1 \times u_2,u_3) \in k^\times$. From (\ref{ss.BACUNO.fig}.\ref{DGRAD}) we deduce $T(u_i \times u_j,u_l) = T(u_1 \times u_2,u_3)$. Hence the linear form $x \mapsto T(u_i \times u_j,x)$ from $J$ to $k$ takes $u_l$ to an invertible element in $k$, forcing $u_l \in J$ to be unimodular. The
adjoint identity (\ref{ss.BACUNO.fig}.\ref{ADJI}) in the special form
$N(u_l)u_l = u_l^{\sharp\sharp} = 0$ therefore yields $N(u_l) = 0$.
so we have
\begin{align}
\label{NORQ} N(q) = T(u_1 \times u_2,u_3);
\end{align}
in particular, 
\begin{align}
\label{QINVP} p = q^{-1} = N(q)^{-1}q^\sharp = T(u_1 \times
u_2,u_3)^{-1}\sum (u_j \times u_l)
\end{align}
since the $u_i$, $1 \leq i \leq 3$, have $u_i^\sharp = 0$. From (\ref{ss.BACUNO.fig}.\ref{DGRAD}) we deduce $T(u_j \times u_l,u_j) = 2T(u_j^\sharp,u_l) = 0$ and, similarly, $T(u_j \times u_l,u_l) = 0$. Applying (\ref{ss.ISCUNO}.\ref{ISOAD}), (\ref{ss.ISCUNO}.\ref{ISOLINT}) and \eqref{QINVP},
we not only conclude $u_i^{(\sharp,p)} = N(p)U_qu_i^\sharp = 0$ but
also $T^{(p)}(u_i) = T(p,u_i) = N(q)^{-1}T(\sum u_m \times u_n,u_i)
= N(q)^{-1}T(u_j \times u_l,u_i) = 1$ by \eqref{NORQ}. Thus
$u_1,u_2,u_3$ are elementary idempotents in $J^{(p)}$. That they do,
in fact, form an elementary frame in that algebra will follow from
Prop.~\ref{p.COMPOR} once we have shown $u_1 \times^{(p)} u_2 =
1^{(p)} - u_1 - u_2 = p^{-1} - u_1 - u_2 = q - u_1 - u_2$. In order
to do so, we first linearize (\ref{ss.ISCUNO}.\ref{ISOAD}) and apply
\eqref{NORQ}, \eqref{QINVP} to obtain
\begin{align}
\label{ISORTH} u_1 \times^{(p)} u_2 =\,\,&N(q)^{-1}U_q(u_1 \times
u_2) \notag \\
=\,\,&N(q)^{-1}T(q,u_1 \times u_2)q - N(q)^{-1}q^\sharp \times (u_1
\times u_2) \\
=\,\,&q - T(u_1 \times u_2,u_3)^{-1}\sum (u_j \times u_l) \times
(u_1 \times u_2) \notag
\end{align}
The term for $i = 3$ in the sum on the very right of \eqref{ISORTH}
by (\ref{ss.BACUNO.fig}.\ref{SADJ}) gives $(u_1 \times u_2) \times (u_1
\times u_2) = 2(u_1 \times u_2)^\sharp = 2(T(u_1^\sharp,u_2)u_2 +
T(u_1,u_2^\sharp)u_1 - u_1^\sharp \times u_2^\sharp) = 0$. On the
other hand, invoking (\ref{ss.BACUNO.fig}.\ref{DDADJ}) for $x = u_i$, $y
= u_j$, $z = u_l$, we deduce
\begin{align*}
(u_i \times u_j) \times (u_i \times u_l) =\,\,&T(u_i^\sharp, u_j)u_l
+ T(u_i^\sharp,u_l)u_j + T(u_i \times u_j,u_l)u_i - u_i^\sharp
\times (u_j \times u_l) \\
=\,\,&T(u_1 \times u_2,u_3)u_i
\end{align*}
Plugging all this this into \eqref{ISORTH} we end up with $u_1
\times^{(p)} u_2 = q - u_2 - u_1$, as claimed.

Finally, in order to derive the very last statement of the problem, it suffices to note $\sum e_i = 1$ and apply the previous parts of the problem.
\end{sol}

\begin{sol}{pr.RENIRA} \label{sol.RENIRA}  (a) follows immediately from  Exc.~\ref{pr.PARNIL}~(c),  Exc.~\ref{pr.CUBNIL} and the definition of cubic ideals.

(b) Let $e \in \pi^{-1}(e_0)$ be an idempotent. Then $\pi(e) = e_0$, and with the usual notational conventions, (\ref{ss.SELICUG}.\ref{PHISH}) implies $\pi(e^\sharp) =e_0^\sharp = 0$, so $e^\sharp \in I$ is nilpotent. But it is also an idempotent since (\ref{ss.BACUNO.fig}.\ref{UADJ}) implies $e^{\sharp 2} = e^{2\sharp} = e^\sharp$. Summing up, therefore, $e^\sharp = 0$. It remains to show $T(e) = 1$. From (\ref{ss.SELICUG}.\ref{NOSELI}) we deduce $\sigma(T(e)) = T_0(e_0) = 1_{k_0}$, hence $T(e) \in 1_k + \mfa \subseteq k^\times$. Thus $e \in J$ is unimodular. On the other hand, (\ref{ss.BACUNO.fig}.\ref{ESTESH}) and (\ref{ss.BACUNO.fig}.\ref{ADJJA}) yield $(1 - T(e))e = 0$, hence $T(e) = 1$, as claimed.
\end{sol}

\begin{sol}{pr.DIAFIX} \label{sol.DIAFIX}  (a) (i) $\Leftrightarrow$ (ii). Condition (i) implies $\Gamma^\prime\Delta^2\Gamma^\sharp =
\Gamma^\prime\Gamma^{\prime\sharp} = N(\Gamma^\prime)\Eins_3 =
N(\Delta)N(\Gamma)\Eins_3
=\Delta\Delta^\sharp\Gamma\Gamma^\sharp$. Canceling yields
$\Gamma^\prime\Delta = \Delta^\sharp\Gamma$, hence (ii). Thus (i)
implies (ii). Conversely, suppose (ii) holds. Then
$\Gamma^\prime\Delta = \Delta^\sharp\Gamma$, and taking norms we conclude
$N(\Gamma^\prime)N(\Delta) = N(\Delta)^2N(\Gamma)$, hence the second
relation of (i). Moreover,
$\Gamma^{\prime\sharp}\Delta^\sharp = (\Gamma^\prime\Delta)^\sharp =
(\Delta^\sharp\Gamma)^\sharp = \Delta^{\sharp\sharp}\Gamma^\sharp =
N(\Delta)\Delta\Gamma^\sharp = \Delta^2\Delta^\sharp\Gamma^\sharp$,
which implies $\Gamma^{\prime\sharp} = \Delta^2\Gamma^\sharp$, hence
the first relation of (i). Thus (i) holds.

(ii) $\Rightarrow$ (iii). Suppose (ii) holds. Hence so does (i).
Since $\vph := \vph_{C,\Delta}$ preserves base points, it suffices
to show that it preserves adjoints as well (Exc.~\ref{pr.HOCUJO}).
Accordingly, let
\begin{align*}
x = \sum(\xi_ie_{ii} + u_i[jl]) \in \Her_3(C,\Gamma)
&&(\xi_i \in k,\;u_i \in C,\;i = 1,2,3).
\end{align*}
Then (\ref{ss.THERCU}.\ref{MADJ}), \eqref{DIAFEL} imply
\begin{align}
\label{EXSH} \varphi(x^\sharp) = \sum\big(\xi_j\xi_l -
\gamma_j\gamma_ln_C(u_i)\big)e_{ii} +
\sum\big(-\delta_i^{-1}\xi_iu_i +
\delta_i^{-1}\gamma_i\overline{u_ju_l}\big)[jl].
\end{align}
On the other hand, \eqref{DIAFEL}, (\ref{ss.THERCU}.\ref{MADJ}) and
(i), (ii) yield
\begin{align*}
\varphi(x)^\sharp =\,\,&\sum\big(\xi_j\xi_l -
\gamma_j^\prime\gamma_l^\prime n_C(\delta_i^{-1}u_i)\big)e_{ii} +
\sum\big(-\xi_i\delta_i^{-1}u_i +
\gamma_i^\prime\overline{(\delta_j^{-1}u_j)(\delta_l^{-1}u_l)}\big)[jl]
\\
=\,\,&\sum\big(\xi_j\xi_l -
\delta_i^{-2}\gamma_j^\prime\gamma_l^\prime n_C(u_i)\big)e_{ii} +
\sum\big(-\delta_i^{-1}\xi_iu_i +
\delta_j^{-1}\delta_l^{-1}\gamma_i^\prime\overline{u_ju_l}\big)[jl]
\\
=\,\,&\sum\big(\xi_j\xi_l - \gamma_j\gamma_ln_C(u_i)\big)e_{ii}
+ \sum\big(-\delta_i^{-1}\xi_iu_i +
\delta_i^{-1}\gamma_i\overline{u_ju_l}\big)[jl],
\end{align*}
and a comparison with \eqref{EXSH} gives the assertion. 

Finally, let us assume that (iii) holds and $1_C \in C$ is unimodular. We wish to establish (ii), apply (\ref{ss.EICJM}.\ref{OFFTI}) to compute $\vph(1_C[li] \times 1_C[ij]) = \vph(1_C[li]) \times \vph(1_C[ij])$ and obtain $\delta_i^{-1}\gamma_i1_C[jl] = \delta_j^{-1}\delta_l^{-1}\gamma_i^\prime 1_C[jl]$. Hence $\gamma_i^\prime 1_C = \delta_j\delta_l\delta_i^{-1}\gamma_i1_C$, and (ii) follows.

(b) Let $\vep, \vep_1, \vep_2, \vep_3 \in k^\times$. For (i), it suffices to put $\delta_1 = \delta_2 = \delta_3 = \vep$ in (a), while (ii) follows by setting $\delta_i = \vep_j\vep_l$ in (a) since this implies $\delta_j\delta_l\delta_i^{-1} = \vep_l\vep_i\vep_i\vep_j\vep_j^{-1}\vep_l^{-1} = \vep_i^2$. Finally, in (iii), we put $\delta_i = \gamma_i^{-1}$ and have $\gamma_i^\prime = \gamma_j^{-1}\gamma_l^{-1}\gamma_i^2$.
\end{sol}
 
\begin{sol}{pr.DIISJO} \label{sol.DIISJO}  $\vph$ is a linear bijection which obviously preserves units. By Exc.~\ref{pr.HOCUJO}~(a), therefore, it suffices to show that $\vph$ preserves adjoints. For
\[
x = \sum(\xi_ie_{ii} + u_i[jl]) \in J := \Her_3(C,\Gamma),
\]
we note $p^{-1} = \sum \gamma_i^{-1}e_{ii}$ and apply (\ref{ss.ISCUNO}.\ref{ISOAD}),  (\ref{ss.THERCU}.\ref{MADJ}), Prop.~\ref{p.COORPE}, (\ref{p.PEDEDI}.\ref{PEDEDI}) to conclude
\begin{align*}
x^{(\sharp,p)} =\,\,&N(p)U_{p^{-1}}x^\sharp = \gamma_1\gamma_2\gamma_3U_{\sum \gamma_i^{-1}e_{ii}}x^\sharp = \gamma_1\gamma_2\gamma_3\sum(\gamma_i^{-2}U_{e_{ii}}x^\sharp + \gamma_j^{-1}\gamma_l^{-1}U_{e_{jj},e_{ll}}x^\sharp) \\
=\,\,&\gamma_1\gamma_2\gamma_3\sum\Big(\gamma_i^{-2}\big(\xi_j\xi_l - \gamma_j\gamma_ln_C(u_i)\big)e_{ii} + \big(\gamma_j^{-1}\gamma_l^{-1}(-\xi_iu_i + \gamma_i\overline{u_ju_l})\big)[jl]\Big) \\
=\,\,&\sum\Big(\gamma_i^{-1}\big(\gamma_j\gamma_l\xi_j\xi_l - \gamma_j^2\gamma_l^2n_C(u_i)\big)e_{ii} + \gamma_i(-\xi_iu_i + \gamma_i\overline{u_ju_l})[jl]\Big),
\end{align*}
hence
\begin{align*}
\vph(x^{(\sharp,p)}) =\,\,&\sum\Big(\big((\gamma_j\xi_j)(\gamma_l\xi_l) - n_C(\gamma_j\gamma_lu_i)\big)e_{ii} + \big(-(\gamma_i\xi_i)(\gamma_j\gamma_lu_i) + \overline{(\gamma_l\gamma_iu_j)(\gamma_i\gamma_ju_l)}\big)[jl]\Big) \\
=\,\,&\Big(\sum\big((\gamma_i\xi_i)e_{ii} + (\gamma_j\gamma_lu_i)[jl]\big)\Big)^\sharp = \vph(x)^\sharp,
\end{align*}
as claimed.
\end{sol}

\begin{sol}{pr.ISCOPA} \label{sol.ISCOPA}  Put
\[
(D,\Delta) := (C,\Gamma)^{(p,q)}, \quad \Delta = \diag(\delta_1,\delta_2,\delta_3), \quad \delta_i = \gamma_i^{(p,q)}
\] 
for $i = 1,2,3$. Then $C = D$ as $k$-modules, and we have
\begin{align}
\label{DEPEQU} \delta_1 = n_C(q)\gamma_1, \quad \delta_2 = n_C(p)\gamma_2, \quad \delta_3 = n_C(pq)^{-1}\gamma_3.
\end{align}
Moreover, from Exc.~\ref{pr.ISTQALT} we recall
\begin{align}
\label{NOCOPE} n_D(u) = n_C(pq)n_C(u), \quad \bar u^{(p,q)} = n_C(pq)^{-1}\overline{pq}\,\bar u\,\overline{pq} = n_C(pq)^{-1}\overline{(pq)u(pq)} &&(u \in C).
\end{align}
Now let
\begin{align}
\label{EXDEDE} x = \sum(\xi_ie_{ii} + u_i[jl]) \in \Her_3(D,\Delta) &&(\xi_i \in k,\;\;u_i \in D,\;\;i = 1,2,3)
\end{align}
and write $\sharp_{p,q}$ for the adjoint of the cubic Jordan matrix algebra $\Her_3(D,\Delta)$. Furthermore, denote by $\eta_i,\zeta_i \in k$, $v_i,w_i \in C$ the quantities satisfying
\begin{align}
\label{EXSHPE} x^{\sharp_{p,q}} =\,\,&\sum(\eta_ie_{ii} + v_i[jl]) \in \Her_3(D,\Delta), \\ 
\vph(x)^\sharp =\,\,&\sum(\zeta_ie_{ii} + w_i[jl]) \in \Her_3(C,\Gamma). \notag
\end{align}
Since
\begin{align}
\label{FIEX.ISCOPA} \vph(x) =\,\,&\sum(\xi_ie_{ii} + u_i^\prime[jl]) \in \Her_3(C,\Gamma), \\ 
\vph(x^{\sharp_{p,q}}) =\,\,&\sum(\eta_ie_{ii} + v_i^\prime[jl]) \in \Her_3(C,\Gamma), \notag
\end{align}
where the $u_i^\prime,v_i^\prime$ are to be determined by $u_i,v_i$, respectively, via \eqref{OFFDI}, we have to show
\begin{align}
\label{EZEVE} \eta_i = \zeta_i, \quad v_i^\prime = w_i &&(i = 1,2,3).
\end{align}
Reading (\ref{ss.THERCU}.\ref{MADJ}) in $\Her_3(D,\Delta)$ and invoking \eqref{DEPEQU}, \eqref{NOCOPE}, we compute
\begin{align*}
\eta_1 =\,\,&\xi_2\xi_3 - \delta_2\delta_3n_C(pq)n_C(u_1) = \xi_2\xi_3 - \gamma_2\gamma_3n_C(p)n_C(u_1) \\
=\,\,&\xi_2\xi_3 - \gamma_2\gamma_3n_C(u_1^\prime) = \zeta_1, \\
\eta_2 =\,\,&\xi_3\xi_1 - \delta_3\delta_1n_C(pq)n_C(u_2) = \xi_3\xi_1 - \gamma_3\gamma_1n_C(q)n_C(u_2) \\
=\,\,&\xi_3\xi_1 - \gamma_3\gamma_1n_C(u_2^\prime) = \zeta_2, \\
\eta_3 =\,\,&\xi_1\xi_2 - \delta_1\delta_2n_C(pq)n_C(u_3) = \xi_1\xi_2 - \gamma_1\gamma_2n_C(pq)^2n_C(u_3) \\
=\,\,&\xi_1\xi_2 - \gamma_1\gamma_2n_C(u_3^\prime) = \zeta_3.
\end{align*}
This proves the first set of equations in \eqref{EZEVE}. In the second one, we have to show $v_i^\prime = w_i$ for $i = 1,2,3$. The case $i = 3$ is comparatively harmless since 
\begin{align*}
v_3^\prime =\,\,&(pq)v_3(pq) = (pq)\big(-\xi_3u_3 + \delta_3\overline{(u_1p)(qu_2)}^{(p,q)}\big)(pq) \\
=\,\,&-\xi_3(pq)u_3(pq) + \gamma_3n_C(pq)^{-1}n_C(pq)^{-1}(pq)\big(\overline{pq}\,\overline{(u_1p)(qu_2)}\,\overline{pq}\big)(pq) \\
=\,\,&-\xi_3u_3^\prime + \gamma_3\overline{u_1^\prime u_2^\prime} = w_3,
\end{align*}
as desired. But the remaining cases $i = 1,2$ are considerably more involved. At a crucial stage, they require an application of (\ref{ss.UOPALT}.\ref{ULR}) combined with the Moufang identities:
\begin{align*}
v_1^\prime = v_1p =\,\,&(-\xi_1u_1 + \delta_1\overline{(u_2p)(qu_3)}^{(p,q)})p = \\
=\,\,&-\xi_1u_1p + \gamma_1n_C(q)n_C(pq)^ {-1}\overline{(pq)\big((u_2p)(qu_3)\big)(pq)}p,
\end{align*}
where
\begin{align*}
(pq)\big((u_2p)(qu_3)\big)(pq) =\,\,&\big((pq)(u_2p)\big)\big((qu_3)(pq)\big) = \big(p(qu_2)p\big)\big(q(u_3p)q\big) \\
=\,\,&p\Big[(qu_2)\Big(p\big(q(u_3p)q\big)\Big)\Big] = p\Big((qu_2)\big((pq)u_3(pq)\big)\Big) \\
=\,\,&p(u_2^ \prime u_3^ \prime).
\end{align*}

Hence Kirmse's identities (\ref{ss.IDCO}.\ref{KID}) yield
\begin{align*}
v_1^\prime =\,\,&-\xi_1u_1^\prime + \gamma_1n_C(p)^{-1}\overline{p(u_2^\prime u_3^\prime)}p = -\xi_1u_1^\prime + \gamma_1n_C(p)^{-1}(\overline{u_2^\prime u_3^\prime}\bar p)p \\
=\,\,&-\xi_1u_1^\prime + \gamma_1\overline{u_2^\prime u_3^\prime} = w_1.
\end{align*}
Similarly,
\begin{align*}
v_2^\prime = qv_2 =\,\,&q\big(-\xi_2u_2 + \delta_2\overline{(u_3p)(qu_1)}^{(p,q)}\big) \\
=\,\,&-\xi_2qu_2 + \gamma_2n_C(p)n_C(pq)^{-1}q\overline{(pq)\big((u_3p)(qu_1)\big)(pq)},
\end{align*}
where
\begin{align*}
(pq)\big((u_3p)(qu_1)\big)(pq) =\,\,&\big((pq)(u_3p)\big)\big((qu_1)(pq)\big) = \big(p(qu_3)p\big)\big(q(u_1q)q\big) \\
=\,\,&\Big[\Big(\big(p(qu_3)p\big)q\Big)(u_1p)\Big]q = \Big(\big((pq)u_3(pq)\big)(u_1p)\Big)q \\
=\,\,&(u_3^\prime u_1^\prime)q.
\end{align*}
Thus
\begin{align*}
v_2^\prime =\,\,&-\xi_2u_2^\prime + \gamma_2n_C(q)^{-1}q\overline{(u_3^\prime u_1^\prime)q} = -\xi_2u_2^\prime + \gamma_2n_C(q)^{-1}q(\bar q\overline{u_3^\prime u_1^\prime}) \\
=\,\,&-\xi_2u_2^\prime + \gamma_2\overline{u_3^\prime u_1^\prime} = w_2. 
\end{align*}
Now \eqref{NOCHA} follows since multiplying each entry of $\Gamma^{(p,q)}$ by $n_C(pq)$ and clearing squares converts $\Gamma^{(p,q)}$ into $\Gamma^\prime$, so Exc.~\ref{pr.DIAFIX}~(b) applies.
\end{sol}

\begin{sol}{pr.OUTCENCU} \label{sol.OUTCENCU}  Let $\alpha \in \Cenout(J)$. We have $[\alpha,U_x] = [\alpha,V_{x,y}] = 0$ for all $x,y \in J$ and must show $\alpha = \xi\Eins_J$ for some $\xi \in k$. Since $\alpha$ commutes with the Peirce projections relative to the diagonal frame of $J$, it stabilizes the corresponding Peirce components, so there are $\xi_i \in k$ and $k$-linear maps $\vph_i\:C \to C$ such that
\begin{align}
\label{ALPHEI} \alpha e_{ii} = \xi_ie_{ii}, \quad \alpha u_i[jl] = \vph_i(u_i)[jl] &&(\xi_i \in k,\;\;u_i \in C,\;\;i = 1,2,3).
\end{align} 
Applying (\ref{ss.EICJM}.\ref{UJLJJ}) and \eqref{ALPHEI}, we obtain $\gamma_j\gamma_l\xi_je_{ll} = U_{1_C[jl]}\xi_je_{jj} = U_{1_C[jl]}\alpha e_{jj} = \alpha U_{1_C[jl]}e_{jj} = \gamma_j\gamma_l\alpha e_{ll} = \gamma_j\gamma_l\xi_le_{ll}$, hence
\begin{align}
\label{ALEQ} \xi_1 = \xi_2 = \xi_3 = :\xi.
\end{align}
Similarly, since $u_i[jl] \circ 1_C[li] = \gamma_l\bar u_i[ij]$ by  (\ref{ss.EICJM}.\ref{OFFTI}), applying \eqref{ALPHEI} again gives $\gamma_l\vph_l(\bar u_i)[ij] = \alpha\gamma_l\bar u_i[ij] = \alpha V_{1_C[li]}u_i[jl] = V_{1_C[li]}\alpha u_i[jl] = \vph_i(u_i)[jl] \circ 1_C[li] = \gamma_l\overline{\vph_i(u_i)}[ij]$, hence
\begin{align}
\label{PHIL} \overline{\vph_i(u_i)} = \vph_l(\bar u_i) &&(u_i \in C,\;\;i = 1,2,3)
\end{align}
On the other hand, \eqref{ALPHEI} and (\ref{ss.EICJM}.\ref{UJLJL}) imply $\gamma_j\gamma_l\overline{\vph_i(v_i)}[jl] = U_{1_C[jl]}\vph_i(v_i)[jl] = U_{1_C[jl]}\alpha v_i[jl] = \alpha U_{1_C[jl]}v_i[jl] = \gamma_j\gamma_l\alpha\bar v_i[jl] = \gamma_j\gamma_l\vph_i(\bar v_i)[jl]$, hence
\begin{align}
\label{PHIBAR} \overline{\vph_i(v_i)} = \vph_i(\bar v_i) &&(v_i \in C,\;\;i = 1,2,3).
\end{align}
Combining \eqref{PHIL}, \eqref{PHIBAR}, we have
\begin{align}
\label{PHIEQ} \vph_1 = \vph_2 = \vph_3 = :\vph.
\end{align}
By \eqref{ALEQ} and \eqref{PHIEQ}, therefore, $\vph(u_1)[23] = \alpha u_1[23] = \alpha V_{u_1[23]}e_{22} = V_{u_1[23]}\alpha e_{22} = \xi u_1[23] \circ e_{22} = \xi u_1[23]$. Thus $\vph = \xi\Eins_C$, and we have shown $\alpha = \xi\Eins_J$, as claimed.
 \end{sol}

\begin{sol}{pr.IDCUJM} \label{sol.IDCUJM}  (a) We consider arbitrary elements
\[
x =\sum(\xi_ie_{ii} + u_i[jl]), \quad y = \sum(\eta_ie_{ii} + v_i[jl])
\]
in $J$, with $\xi_i,\eta_i \in k$, $u_i,v_i \in C$ for $i = 1,2,3$. If $I$ is an ideal in $(C,\iota_C)$ and $I_0 \subseteq I \cap k$ is a weakly $I$-ample ideal in $k$, we must show that $H := H_3(I_0,I,\Gamma)$ is an outer ideal in $J$. For $u \in I$ and $v \in C$ we first note $n_C(u,v) = t_C(u\bar v) \in I_0$ since $I \subseteq C$ is an ideal and $I_0$ is weakly $I$-ample. Now assume $y \in H$, i.e., $\eta_i \in I_0$ and $v_i \in I$ for $i = 1,2,3$. Inspecting (\ref{ss.THERCU}.\ref{MABIT}) and applying the preceding observation, we conclude $T(x,y) \in I_0$, hence $T(x,y)\xi_i \in I_0$ and $T(x,y)u_i \in I_0C \subseteq IC \subseteq I$. Thus $T(x,y)x \in H$. Moreover, since $I$ is stabilized by $\iota_C$, and by (\ref{ss.THERCU}.\ref{MABADJ}), we have $x \times y \in H$, which after replacing $x$ by $x^\sharp$ implies $U_xy\in H$. Thus $H$ is an outer ideal in $J$. Conversely, let this be so. Then the Peirce projections relative to the diagonal frame of $J$, i.e., $U_{e_{ii}}$ and $U_{e_{jj},e_{ll}}$ for $i = 1,2,3$ stabilize $H$ and Prop.~\ref{p.PEDEDI} yields
\begin{align}
\label{HEII} H = \sum\big((H \cap ke_{ii}) + (H \cap C[jl])\big).
\end{align}
We now put
\begin{align}
\label{IUCE} I := \{u \in C \mid u[23] \in H\}, \quad I_0 := \{\xi \in k \mid \xi e_{11} \in H\},
\end{align}
which are $k$-submodules of $C,k$, respectively. For $u \in I$, we apply (\ref{ss.EICJM}.\ref{UJLJL}) and obtain $\gamma_2\gamma_3\bar u[23] = U_{1_C[23]}u[23] \in H$, hence $\bar u \in I$. Thus $\bar I = I$. From (\ref{ss.EICJM}.\ref{OFFTI}) we conclude that $u[jl] \times v[li] = u[jl] \circ v[li]$ belongs to $H$ if one of the ``factors'' does. We now claim
\begin{align}
\label{HACEJL} H \cap C[jl] = I[jl].  
\end{align}
for $i = 1,2,3$. This follows from \eqref{IUCE} for $i = 1$. Arguing by ``cyclic induction mod $3$'', suppose \eqref{HACEJL} holds for some $i \in \{1,2,3\}$ and let $u \in C$. If $u \in I$, then $\bar u \in I$, and (\ref{ss.EICJM}.\ref{OFFTI}) yields $\gamma_ju[li] = 1_C[ij] \times \bar u[jl] \in H$, hence $u[li] \in H$. Conversely, if $u[li] \in H$, the $H$ contains $u[li] \times 1_C[ij] = \gamma_i\bar u[jl]$, which implies $\bar u \in I$, hence $u \in I$. Thus \eqref{HACEJL} holds for $j$ instead of $i$, which completes its proof. Next we claim
\begin{align}
\label{HACAPE} H \cap ke_{ii} = I_0e_{ii}
\end{align} 
for $i = 1,2,3$. Thanks to \eqref{IUCE}, the case $i = 1$ is again obvious. Next suppose \eqref{HACAPE} holds for some $i \in \{1,2,3\}$; we must prove it for $j$. Let $\xi \in k$. For $\xi \in I_0$ we obtain $\xi e_{ii} \in H$, so $H$ by (\ref{ss.EICJM}.\ref{UJLJJ}) contains $U_{1_C[ij]}\xi e_{ii} = \gamma_j\gamma_l\xi e_{jj}$, hence $\xi e_{jj}$. Conversely, if $\xi e_{jj} \in H$, then $H$ contains $U_{1_C[ij]}\xi e_{jj} = \gamma_i\gamma_j\xi e_{ii}$, and we conclude $\xi \in I_0$. This completes the proof of \eqref{HACAPE}. Combining \eqref{HEII} with \eqref{HACEJL}, \eqref{HACAPE}, we obtain
\begin{align}
\label{HAPEIR} H = \sum(I_0e_{ii} + I[jl]),
\end{align}
and it remains to show that (i) $I$ is an ideal in $(C,\iota_C)$, (ii) $I_0 \subseteq I \cap k$, and (iii) $I_0$ is weakly $I$-ample. For $u \in I$, $v \in J$, we combine (\ref{ss.EICJM}.\ref{OFFTI}) with \eqref{HACEJL} to obtain $\gamma_1uv[23] = \bar v[31] \circ \bar u[12] \in H$, hence $uv \in I$. Thus $I$, being stable under conjugation, is an ideal in $(C,\iota_C)$, and (i) is proved. Turning to (ii), let $\xi \in I_0$. Then (\ref{p.PEDEDI}.\ref{PEDEDI}), (\ref{t.PEDECOM}.\ref{JIJ}) imply $C[23] = J_{23} \subseteq J_1(e_{22})$, and (\ref{t.PEDESI}.\ref{JONE}) yields $\xi[23] =(\xi 1_C)[23] = e_{22} \circ (\xi 1_C)[23] = (\xi e_{22}) \circ 1_C[23] \in H$, hence $\xi \in I$. This proves (ii). Finally, turning to weak $I$-ampleness, let $u \in I$. Then $H$ contains $\{u[23]e_{22}1_C[23]\} =\gamma_2\gamma_3t_C(u)e_{33}$ (by (\ref{ss.EICJM}.\ref{UJLJJ}) linearized), which implies $t_C(u) \in I_0$ by \eqref{HACAPE} and completes the proof (iii).

(b) Let $I$ be an ideal in $(C,\iota_C)$ and $I_0$ an ideal in $k$ that is weakly $I$-ample and contained in $I \cap k$. By (a) it will be enough to show that $H := H(I_0,I,\Gamma)$ is an inner ideal of $J$ if and only if $I_0$ is $I$-ample. Suppose first that $H \subseteq J$ is an inner ideal and let $u \in I$. Then $H$ contains  $U_{u[23]}e_{22}$, which by (\ref{ss.EICJM}.\ref{UJLJJ}) agrees with $\gamma_2\gamma_3n_C(u)e_{33}$, and we conclude $n_C(u) \in I_0$. Thus $I_0$ is $I$-ample. Conversely, let this be so. For $x \in H$, i.e., $\xi_i \in I_0$ and $u_i \in I$, $i = 1,2,3$, we obtain $T(x,y)\xi_i \in I_0$, $T(x,y)u_i \in I$, hence $T(x,y)x \in H$. Moreover, $x^\sharp \in H$ by (\ref{ss.THERCU}.\ref{MADJ}) since $I_0$ is $I$-ample. But then $x^\sharp \times y \in H$ after inspecting (\ref{ss.THERCU}.\ref{MABADJ}), and we conclude $U_xy \in H$. Summing up, we have shown that $H \subseteq J$ is an inner ideal.

(c) For $\xi \in I \cap k$ we obtain $2\xi = t_C(\xi 1_C) \in I_0$. Thus $2(I \cap k) \subseteq I_0$. The rest is clear.

(d) Let $H$ be an outer ideal in $J$. The (a) implies $H = H(I_0,I,\Gamma)$ for some ideal $I \subseteq (C,\iota_C)$ and some weakly $I$-ample ideal $I_0$ of $k$ contained in $I \cap k$. By Exc.~\ref{pr.IDEALSCOAL}, there exists an ideal $\mfa \subseteq k$ such that $I = \mfa C$ and $I_0 \subseteq I \cap k = \mfa$. Since $C$ is regular, its trace form is surjective (Lemma~\ref{l.COMTRI}), so some $u \in C$ has $t_C(u) = 1$. For $\xi \in \mfa$, we obtain $\xi u \in I$ and therefore $\xi = t_C(\xi u) \in t_C(I) \subseteq I_0$ by weak $I$-ampleness. This proves $I_0 = \mfa$, hence $H = \mfa J$ by \eqref{HIIG}.

Finally, we must show that the homomorphism $\phi$ is injective. By what we have just seen, $\Ker(\phi) = \mfa J$ for some ideal $\mfa \subseteq k$, so it will be enough to prove $\mfa = \{0\}$. But this is clear since $\alpha \in \mfa$ implies $\alpha 1_J \in \Ker(\phi)$, hence $\alpha 1_A = \phi(\alpha 1_J) = 0$. But $1_A \in A$ is unimodular, which implies $\alpha =  0$, as claimed.

(e) By Kaplansky's theorem (Prop.~\ref{p.PRECOMFIELD}), $C$ is either a regular composition algebra or a purely inseparable field extension $K/F$ of characteristic $2$ and exponent at most $1$. In the former case, the assertion follows immediately from (d). In the latter case, we first note that the bilinearized norm of $C = K$ is zero and hence $\Rad(T) = H(\{0\},K,\Gamma)$ is an outer ideal in $J$. Conversely, suppose $H \subseteq J$ is a non-zero outer ideal and write $H = H(I_0,I,\Gamma)$ with $I_0,I$ as in (a). Then $H \neq \{0\}$ implies $I \neq \{0\}$, hence $I = K$, while $I_0 \subseteq I \cap F = F$ is an ideal in $F$, weak $I$-ampleness being automatic. If $I_0 = \{0\}$, then $H = \Rad(T)$, and if $I_0 = F$, then $H = J$. Thus $J$ is simple but nor outer simple.
\end{sol}

\begin{sol}{pr.ABZECU} \label{sol.ABZECU}  (a) We have $U_x = 0$ and must show $U_{x^\sharp} = 0$, i.e., $U_{x^\sharp}y = 0$ for all $y \in J$. Since cubic Jordan algebras are invariant under base change, viewing $U_{x^\sharp}$ as a polynomial law over $k$ and applying Prop.~\ref{p.ZADE} allow us to assume $y \in J^\times$. Then (\ref{c.REGINV}.\ref{CUINV}) for $x = y^{-1}$ and (\ref{ss.BACUNO.fig}.\ref{UADJ}) imply $U_{x^\sharp}y = N(y)U_{x^\sharp}y^{-1\sharp} = N(y)(U_xy^{-1})^\sharp = 0$.

(b) (i) $\Rightarrow$ (ii). By (i) we have $x^2 = U_x1 = 0$, and since $k$ is reduced, Exc.~\ref{pr.CUBNIL} shows $T(x) = S(x) = N(x) = 0$. By (\ref{ss.BACUNO.fig}.\ref{ADJJA}), therefore,  $x^\sharp = x^2 - T(x)x + S(x)1 = 0$. Thus $T(x,y)x = U_xy = 0$, which implies $T(x,y)^2 = T(T(x,y)x,y)) = 0$, hence $T(x,y) = 0$ by our hypothesis on $k$. Summing up, (ii) holds.

(ii) $\Rightarrow$ (i). Obvious, by the formula for the $U$-operator.

(ii) $\Rightarrow$ (iii). The adjoint identity yields $N(x)x = x^{\sharp\sharp} = 0$, and applying the norm we obtain $N(x)^4 = 0$, hence $N(x) = 0$. 

(c) Let $x \in J$ be an absolute zero divisor and assume first that $k$ is reduced. Then (b) implies $T(x,y) = T(x^\sharp,y) = N(x) = 0$ for all $y \in J$, and Exc.~\ref{pr.CUBNIL} shows $x \in \Nil(J)$.

Next assume that $k$ is arbitrary. Then $\bar k := k/\Nil(k) \in \kalg$ is reduced, and we have a canonical identification $J_{\bar k} = \bar J := J/\Nil(k)J$ via \ref{ss.REDID}, matching $z_{\bar k}$ for $z \in J$ with $\bar z$, the image of $z$ under the natural map $J \to \bar J $. Since $x$ is an absolute zero divisor in $J$, $\bar x$ is one in $\bar J$, so by the special case just treated we have $\bar T(\bar x,\bar y) = \bar T(\bar x^\sharp,\bar y) = \bar N(\bar x) = 0$ for all $y \in J$, where $\bar T = T_{\bar k}$ (resp. $\bar N = N \otimes \bar k$) is the bilinear trace (resp, norm) of $\bar J$ over $\bar k$. But this means that $T(x,y)$, $T(x^\sharp,y)$, $N(x)$ are nilpotent elements of $k$, for all $y \in J$, which means $x \in \Nil(J)$.

(d) Let $F$ be a field of characteristic $2$, $K/F$ a purely inseparable field extension of exponent at most $1$, $\Gamma = \diag(\gamma_1,\gamma_2,\gamma_3) \in \GL_3(F)$ and $J = \Her_3(K,\Gamma)$. By Exc.~\ref{pr.IDCUJM}, $J$ is a simple Jordan algebra over $F$. If $u_1,u_2,u_3 \in K$ are not all zero, neither is $x := \sum u_i[jl] \in J$. Moreover, since the bilinearized norm of $K/F$ is zero, (\ref{ss.THERCU}.\ref{MANO}) and (\ref{ss.THERCU}.\ref{MABIT}) show $N(x) = 0$ and $T(x,y) = 0$ for all $y \in J$, so $x$ satisfies condition (iii) of (b). On the other hand,
\[
x^\sharp = \sum(\gamma_j\gamma_ln_K(u_i)e_{ii} + \gamma_iu_ju_l[jl]) \neq 0,
\]
whence (b) shows that $x$ is \emph{not} an absolute zero divisor of $J$.

(e) First assume $\Nil(J) = \{0\}$. We have seen in Exc.~\ref{pr.PARNIL}~(c) that $\Nil(k)J$ is a nil ideal in $J$. This proves $\Nil(k)J = \{0\}$, and in particular $\Nil(k)1 = \{0\}$. But $1 \in J$ is unimodular, and we conclude that $k$ is reduced. Now let $x \in J$ be an absolute zero divisor. Then (b) implies $N(x) = T(x,y) = T(x^\sharp,y) = 0$ for all $y \in J$. From Exc.~\ref{pr.CUBNIL} we therefore deduce $x \in \Nil(J)$, hence $x = 0$, and $J$ has no absolute zero divisors. Conversely, let $k$ be reduced and suppose $J$ has no absolute zero divisors. For $x \in \Nil(J)$ and $y \in J$, we have $T(x,y) = T(x^\sharp,y) = N(x) = 0$. This implies $x^{\sharp\sharp} = N(x)x = 0$, so $x^\sharp$ is an absolute zero divisor by (b), forcing $x^\sharp = 0$ by hypothesis. But then, again by (b), $x$ is an absolute zero divisor, and we conclude $x = 0$, as desired.  
\end{sol}

\begin{sol}{pr.NIRAPE} \label{sol.NIRAPE}  Let $i = 0,2$. Then $J_i(e) \cap \Nil(J)$ is a nil ideal in $J_i(e)$ and hence contained in $\Nil(J_i(e))$. It therefore remains to prove the converse, i.e., 
\begin{align}
\label{NILRAPE} \Nil\big(J_i(e)\big) \subseteq \Nil(J) &&(i = 0,2).
\end{align}
In order to do so, we first reduce to the case that $e$ is elementary. Indeed, assume that this case has been settled and let $e$ be arbitrary. Clearly, \eqref{NILRAPE} holds for $e$ if and only if it holds for the complementary idempotent $1_J - e$. By assumption, therefore, it holds for $e = 0,1_J$, $e$ elementary and $e$ co-elementary. Applying now Exc.~\ref{pr.IDCUJO} and setting $e^{(0)} = 0$, $e^{(1)} = 1_{J^{(3)}}$, we observe that taking nil radicals commutes with direct products of ideals. Hence $J_i(e) = \prod_{l=0}^3 J^{(l)}_i(e^{(l)})$ implies
\begin{align*}
\Nil\big(J_i(e)\big) =\,\,&\prod_{l=0}^3 \Nil\big(J^{(l)}_i(e^{(l)})\big) \subseteq \prod_{l=0}^3 \Nil(J_i^{(l)}) = \Nil(J)
\end{align*}
and completes the reduction. We may thus assume that $e$ is elementary. Let us first consider the case $i = 2$. Since $\Nil(k) = \Nil(k^{(+)})$ and the map $\alpha \mapsto \alpha e$ is an isomorphism from $k^{(+)}$ to $J_2(e) = ke$, Exc.~\ref{pr.PARNIL}~(c) implies $\Nil(J_2(e)) = \Nil(k)e \subseteq \Nil(k)J \subseteq \Nil(J)$, as claimed. We are left with the case $i = 0$. Let $x \in \Nil(J_0(e))$ and put $f := 1 - e$. Combining Cor.~\ref{c.FAULE} with Exc.~\ref{pr.NIPOQUA}, we conclude $S(x),S(x,y) \in \Nil(k)$ for all $y \in J_0(e)$, hence in particular $T(x) = S(x,f) \in \Nil(k)$. Thus (\ref{ss.BACUNO.fig}.\ref{BQUAT}) implies $T(x,y) = T(x)T(y) - S(x,y) \in \Nil(k)$ for all $y \in J_0(e)$, and since the Peirce components of $J$ relative to $e$ are orthogonal with respect to the bilinear trace, by Prop.~\ref{p.PELE}~(b), this amounts to $T(x,J) \subseteq \Nil(k)$. On the other hand, from Prop.~\ref{p.PELE}~(c) we deduce $x^\sharp = S(x)e$, hence in particular, $T(x^\sharp,J) = S(x)T(e,J) \subseteq \Nil(k)$. And finally, $N(f) = N(1 - e) = 1 - T(e) + T(e^\sharp) - N(e) = 0$ and (\ref{ss.BACUNO.fig}.\ref{JCOMP}) yield $N(x) = N(U_fx) = N(f)^2N(x) = 0$. Summing up, Exc.~\ref{pr.CUBNIL} now shows $x \in \Nil(J)$, and the proof is complete.
\end{sol}

\begin{sol}{pr.HEMONI} \label{sol.HEMONI}  We put $I_0 := \Nil(k)$, $ I := \Nil(C)$. In \eqref{NILCEK}, the inclusion follows from Exc.~\ref{pr.NILRADCON}~(c). As to the rest, we clearly have $I \cap k \subseteq I_0$, while Exc.~\ref{pr.NILRAD} implies $I_0 \subseteq I_0C \subseteq I$. By Exc.~\ref{pr.IDCUJM}~(b), therefore,
\[
H_3(I_0,I,\Gamma) = \sum(I_0e_{ii} + I[jl])
\]
is an ideal in $J := \Her_3(C,\Gamma)$. Now let $\xi_i \in I_0$, $u_i \in I$, $\eta_i \in k$, $v_i \in C$ for $i = 1,2,3$ and put 
\[
x := \sum(\xi_ie_{ii} + u_i[jl]) \in H_3(I_0,I,\Gamma), \quad y := \sum(\eta_ie_{ii} + v_i[jl]) \in J.
\]
Using Exc.~\ref{pr.NILRADCON}~(c), (\ref{ss.THERCU}.\ref{MANO}), (\ref{ss.THERCU}.\ref{MABIT}), (\ref{ss.THERCU}.\ref{DIMANO}) and the relation $n_C(u_i,v_i) = t_C(u_i\bar v_i)$, one checks that $T(x,y)$, $T(x^\sharp,y)$, $N(x)$ all belong to $I_0$. Hence Exc.~\ref{pr.CUBNIL} implies $H_3(I_0,I,\Gamma) \subseteq \Nil(J)$. Conversely, let $x \in \Nil(J)$ and $i = 1,2,3$. Then Exc.~\ref{pr.NIRAPE} and the Peirce rules show $\xi_ie_{ii} = U_{e_{ii}}x \in J_2(e_{ii}) \cap \Nil(J) = \Nil(J_2(e_{ii})) = I_0e_{ii}$, and we conclude $\xi_i \in I_0$. On the other hand,
\[
u_i[jl] = U_{e_{jj},e_{ll}}x \in \Nil(J) \cap J_0(e_{ii}) = \Nil\big(J_0(e_{ii})\big),
\]
again by Exc.~\ref{pr.NIRAPE}. But, by Cor.~\ref{c.FAULE}, $J_0(e_{ii})$ is the Jordan algebra of the pointed quadratic module $(M_0,S_0,f)$, where $M_0 = J_0(e_{ii})$ as $k$-modules, $S_0 = S\vert_{M_0}$ and $f = 1_J - e_{ii} = e_{jj} + e_{ll}$. Here Exc.~\ref{pr.NIPOQUA} and (\ref{ss.THERCU}.\ref{MAQIT}), (\ref{ss.THERCU}.\ref{MABQIT}) imply that $S(u_i[jl]) = -\gamma_j\gamma_ln_C(u_i)$ and $S(u_i[jl],v_i[jl]) = -\gamma_j\gamma_ln_C(u_i,v_i)$ belong to $I_0$ for all $v_i \in C$. From Exc.~\ref{pr.NIRA} we therefore conclude $u_i \in I$, hence $x \in H_3(I_0,I,\Gamma)$, and \eqref{NILHER} follows.

Since $(I_0,I)$ is a conic ideal in $C$ by Exc.~\ref{pr.CONILID}~(b), we deduce from Exc.~\ref{pr.CONID}~(c) that there exists a unique way of viewing $C_0$ as a conic $k_0$-algebra making $\pi$ a $\sigma$-semi-linear homomorphism of conic algebras, and since $C$ is multiplicative alternative, so is $C_0$. Now put $\pJ_0 := \Her_3(C_0,\Gamma_0)$ as a cubic Jordan matrix algebra over $k_0$ and consider the map $\Phi\:J \to \pJ_0$ defined by
\begin{align}
\label{PHISI} \Phi\big(\sum(\xi_ie_{ii} + u_i[jl])\big) := \sum\big(\sigma(\xi_i)e_{ii} + \pi(u_i)[jl]\big)
\end{align}
for $\xi_i \in k$, $u_i \in C$, $i = 1,2,3$. Obviously, $\Phi$ is $\sigma$-semi-linear preserving base points and, using (\ref{ss.SELICON}.\ref{ENCEL}), one checks that it makes a commutative diagram
\[
\xymatrix{J \ar[r]_{\Phi} \ar[d]_{\sharp_J} & \pJ_0 \ar[d]^{\sharp_{\pJ_0}} \\
J \ar[r]_{\Phi} & \pJ_0}
\]
of set maps. Hence, by Exc.~\ref{pr.SELIBA}~(b) and \eqref{NILHER}, $\Phi$ is a surjective $\sigma$-semi-linear homomorphism of cubic Jordan algebras with kernel $\Nil(J)$. By Exec.~\ref{pr.RENIRA}, on the other hand, $J_0 := J/\Nil(J)$ carries the unique structure of a cubic Jordan algebra over $k_0$ making the canonical projection $\Pi\:J \to J_0$ a $\sigma$-semi-linear homomorphism of cubic Jordan algebras. Since $\Phi$ and $\Pi$ are $\sigma$-semi-linear surjective with the same kernel, there is a unique isomorphism $\Psi\:J_0 \overset{\sim} \to \pJ_0$ of Jordan algebras over $k_0$ such that the diagram
\[
\xymatrix{& J \ar[ld]_{\Pi} \ar[rd]^{\Phi} & \\
J_0 \ar@{.>}[rr]_{\exists! \Psi}^{\cong} && \pJ_0}
\]
commutes. It follows that
\[
\xymatrix{J \ar[r]_{\Pi} \ar[d]_{N_J} & J_0 \ar@<0.5ex>[d]^{N_{\pJ_0} \circ \Psi}_{}\ar@<-0.5ex>[d]_{N_{J_0}}^{} \\
k \ar[r]_{\sigma} & k_0}
\]
with \emph{both} vertical arrows on the right is a commutative $\sigma$-semi-linear polynomial square. Since $\Pi$ is surjective, we conclude $N_{\pJ_0} \circ \Psi  = N_{J_0}$ as cubic forms over $k_0$, whence $\Psi$ may be used to identify $J_0 = \pJ_0$ as cubic Jordan algebras over $k_0$. this completes the proof.
\end{sol}

\begin{sol}{pr.EXRACO} \label{sol.EXRACO}  We put
\begin{align*}
M :=\,\,&\{\sum \xi_ie_i \mid \xi_i \in k,\;\;\xi_iJ_{ij} = \{0\}\;\;\text{for}\;\;i = 1,2,3\}, \\
N :=\,\,&\{\sum \xi_ie_i \mid \xi_i^2 = 2\xi_i = 0\;\;\text{for}\;\;i = 1,2,3\}.
\end{align*}
By Exc.~\ref{pr.OUTIJO}~(c), $\Rex(J) \subseteq J$ is an ideal. Thus, given $x = \sum(\xi_ie_i + v_{jl}) \in \Rex(J)$, $\xi_i \in k$, $v_{jl} \in J_{jl}$, its Peirce components relative to $(e_1,e_2,e_3)$ by Prop.~\ref{p.COORPE} belong to $\Rex(J)$ as well, so we have $\xi_ie_i,v_{jl} \in \Rex(J)$ for $1 \leq i \leq 3$. Here (\ref{t.PEDESI}.\ref{JONE}) and (\ref{t.PEDECOM}.\ref{JCI}) imply $v_{jl} = e_j \circ v_{jl} = U_{v_{jl,1}}e_j = 0$, hence $x = \sum \xi_ie_i$. Similarly, $\xi_iJ_{ij} = \xi_ie_i \circ J_{ij} = U_{\xi_ie_i,1}J_{ij} = \{0\}$, and we have shown $\Rex(J) \subseteq M$. Now suppose $x = \sum \xi_ie_i$, $\xi_i \in k$, satisfies $x \in M$, i. e., $\xi_iJ_{ij} = \{0\}$ for $1 \leq i \leq 3$. Setting $u_{12} := u_{23} \times u_{31}$, we have $u_{jl} \in J_{jl}$ by Prop.~\ref{p.DECNOAD} and $S(u_{jl}) \in k^\times$ by Lemma~\ref{l.ESPE} for $i = 1,2,3$. In particular, $\xi_i^2S(u_{ij}) = S(\xi_iu_{ij}) = 0$ and $2\xi_iS(u_{ij}) = S(\xi_iu_{ij},u_{ij}) = 0$ imply $M \subseteq N$, and we have $U_{\xi_ie_i} = 0$. Moreover, the assignment $x \mapsto x \times u_{jl}$ by Lemma~\ref{l.ESPE} gives a linear bijection from $J_{ij}$ to $J_{li}$. Hence $\xi_iJ_{li} = \xi_iJ_{ij} \times u_{jl} = \{0\}$ and then $\xi_iJ_{jl}= \xi_iJ_{li} \times u_{ij} = \{0\}$. Summing up, we have shown
\[
\xi_iJ = \sum_m (k\xi_i)e_m,
\]
and from the Peirce rules we conclude
\begin{align*}
\{(\xi_ie_i)JJ\} =\,\,&\{e_i(\xi_iJ)J\} = \sum_m\{e_i(k\xi_i)e_mJ\} = \{e_ie_i(\xi_iJ)\} = \sum_m \{e_ie_i(k\xi_i)e_m\} \\
=\,\,&k\xi_i\{e_ie_ie_i\} = 2k\xi_ie_i = \{0\}.
\end{align*}
But this means $U_{\xi_ie_i,J} = \{0\}$, and we have shown $x = \sum \xi_ie_i \in \Rex(J)$. This completes the proof of \eqref{CHAREX}. Now let $(C,\Gamma)$ be a co-ordinate pair over $k$ and put $(J,\mfS) :=\bfHer_3(C,\Gamma)$ as a co-ordinated cubic Jordan algebra over $k$. By \eqref{CHAREX}, and with the notation used before, the elements of $\Rex(J)$ have the form $\sum \xi_ie_{ii}$, $\xi_i \in k$, $\xi_iJ_{ij} = 0$ for $1 \leq i \leq 3$. But $J_{ij} = C[ij]$ by Prop.~\ref{p.PEDEDI}, and we conclude $\xi_iC = \{0\}$. In particular, $\xi_i1_C = 0$, and if $1_C \in C$ is unimodular, we deduce $\xi_i = 0$. Thus the extreme radical of $J$ is zero.

Finally, let $k_0$ be any commutative ring in which $2 = 0$ and let $k := k_0[\vep]$ be the $k_0$-algebra of dual numbers. Then $k_0$ may be viewed  as an algebra $C \in \kalg$ under the homomorphism $k \to k_0$ satisfying $\vep \mapsto 0$. The squaring $\alpha_0 \mapsto \alpha_0^2$ may be regarded as a $k_0$-quadratic map $n_C\:C \to k$ with zero bilinearization, which is in fact $k$-quadratic since
\[
n_C\big((\beta_0 + \gamma_0\vep)\alpha_0\big) = n_C(\beta_0\alpha_0) = \beta_0^2n_C(\alpha_0) = (\beta_0 + \gamma_0\vep)^2n_C(\alpha_0)
\]
for all $\alpha_0,\beta_0,\gamma_0 \in C$. Thus $C$ together with $n_C$ is a multiplicative conic commutative associative $k$-algebra such that $\vep C = \{0\}$, and the extreme radical of $J$ is different from zero.
\end{sol}

\begin{sol}{pr.JACOCAT} \label{sol.JACOCAT}  (a) If
\[
\xymatrix{ (C,\Gamma) \ar[rr]_{(\eta,\Delta)} &&
(C^\prime,\Gamma^\prime) \ar[rr]_{(\eta^\prime,\Delta^\prime)} &&
(C^{\prime\prime},\Gamma^{\prime\prime}) }
\]
are morphisms of co-ordinate pairs, so obviously is
\begin{align}
\label{COMCOTR} (\eta^\prime,\Delta^\prime) \circ (\eta,\Delta) :=
(\eta^\prime \circ \eta,\Delta^\prime\Delta)\:(C,\Gamma)
\longrightarrow
(C^{\prime\prime},\Gamma^{\prime\prime}).
\end{align}
It is straightforward to check that we obtain a category in this way, denoted by $\kcopa$. For example, the identity morphism of $(C,\Gamma)$ is $\Eins_{(C,\Gamma)} := (\Eins_C,\Eins_3)$.

(b) If
\[
\xymatrix{ (J,\mfS) \ar[rr]_{(\varphi;\delta_1,\delta_2)} &&
(J^\prime,\mfS^\prime)
\ar[rr]_{(\varphi^\prime;\delta_1^\prime,\delta_2^\prime)} &&
(J^{\prime\prime},\mfS^{\prime\prime}) }
\]
with $\mfS^{\prime\prime} =
(e_1^{\prime\prime},e_2^{\prime\prime},e_3^{\prime\prime},u_{23}^{\prime\prime},u_{31}^{\prime\prime})$
are morphisms of co-ordinated cubic Jordan algebras, so is
\begin{align}
\label{COMCOC} (\varphi^\prime;\delta_1^\prime,\delta_2^\prime)
\circ (\varphi;\delta_1,\delta_2) := (\varphi^\prime
\circ\varphi;\delta_1^\prime\delta_1,\delta_2^\prime\delta_2)\:(J,\mfS)
\longrightarrow (J^{\prime\prime},\mfS^{\prime\prime}).
\end{align}
since $(\varphi^\prime\circ \varphi)(e_i) =
\varphi^\prime(e_i^\prime) = e_i^{\prime\prime}$ for $i = 1,2,3$ by \ref{JACOCAT.3} and $(\varphi^\prime \circ \varphi)(u_{jl}) =
\varphi^\prime(\delta_i^{-1}u_{jl}^\prime) = \delta_i^{\prime
-1}\delta_i^{-1}u_{jl}^{\prime\prime}$ for $i = 1,2$ by \ref{JACOCAT.4}.

It is straightforward to check that we obtain a category in this way, denoted by $\kcocujo$. For example, the identity morphism of $(J,\mfS)$ is $\Eins_{(J,\mfS)} := (\Eins_J;1,1)$.

(c) Let $\eta\:C \to C^\prime$ be a homomorphism of multiplicative conic alternative $k$-algebras. We define
\[
\Her_3(\eta)\:\Her_3(C,\Gamma) \longrightarrow \Her_3(C^\prime,\Gamma)
\]
component-wise by
\begin{align}
\label{FUHER} \Her_3\big(\eta)\big(\sum(\xi_ie_{ii} +
u_i[jl])\big) := \sum\big(\xi_ie_{ii} + \eta(u_i)[jl]\big)
\end{align}
for $\xi_i \in k$, $u_i \in C$, $i = 1,2,3$. Strictly speaking, $\Her_3(\eta)$ depends not only on $\eta$ (and thus on $C$ and $\pC$) but also on $\Gamma$. However, it will always be clear from the context which $\Gamma$ we have in mind. Since $\eta$ by
definition preserves units, norms, traces and conjugations, we easily deduce from (\ref{ss.THERCU}.\ref{MABA})--(\ref{ss.THERCU}.\ref{MANO})
that $\Her_3(\eta)$ is a homomorphism of cubic Jordan algebras.
Given homomorphisms
\[
\xymatrix{C \ar[rr]_{\eta} && C^\prime \ar[rr]_{\eta^\prime} && C^{\prime\prime}}
\]
of multiplicative conic alternative algebras, \eqref{FUHER} easily implies
$\Her_3(\eta^\prime \circ \eta) = \Her_3(\eta^\prime) \circ \Her_3(\eta)$.

Now let $(\eta,\Delta)\:(C,\Gamma) \to (C^\prime,\Gamma^\prime)$ with $\Delta =
\diag(\delta_1,\delta_2,\delta_3) \in \Diag_3(k)^\times$ be a
homomorphism of co-ordinate pairs. We claim that
\begin{align}
\label{COHE} \bfHer_3(\eta,\Delta) :=
\big(\Her_3(\eta,\Delta);\delta_1,\delta_2)\:\bfHer_3(C,\Gamma)
\longrightarrow \bfHer_3(C^\prime,\Gamma^\prime)
\end{align}
is a homomorphism of co-ordinated cubic Jordan algebras. Indeed, since $\Gamma^\prime = \Delta^\sharp\Delta^{-1}\Gamma$ by \ref{JACOCAT.2}, Exc.~\ref{pr.DIAFIX} shows that $\Her_3(\eta,\Delta) = \Her_3(\eta) \circ \vph_{C,\Delta} = \vph_{C^\prime,\Delta} \circ \Her_3(\eta)$ is a homomorphism of the underlying cubic Jordan algebras and by
\eqref{HERET} satisfies $\Her_3(\eta,\Delta)(e_{ii}) = e_{ii}$ for $i = 1,2,3$ and
$\Her_3(\eta,\Delta)(1_C[jl]) = (\delta_i^{-1}1_{C^\prime})[jl]$ for
$i = 1,2$. 

It is straightforward to check using \eqref{HERET} that in this way we obtain a functor
\begin{align}
\label{FUBFHER} \bfHer_3\:\kcopa \longrightarrow \kcocujo.
\end{align}

(d) Writing $T^\prime,S^\prime$ for the trace, quadratic trace of $J^\prime$, we put
\begin{align*}
\mfS =\,\,&(e_1,e_2,e_3,u_{23},u_{31}), \quad \omega := \omega_{J,\mfS}, \\
(C,\Gamma) :=\,\,&\Cot(J,\mfS) = (C_{J,\mfS},\Gamma_{J,\mfS}), \quad \Gamma = \diag(\gamma_1,\gamma_2,\gamma_3), \\
\mfS^\prime =\,\,&(e_1^\prime,e_2^\prime,e_3^\prime,u_{23}^\prime,u_{31}^\prime),
\quad \omega^\prime := \omega_{J^\prime,\mfS^\prime},
\\
(C^\prime,\Gamma^\prime) :=\,\,&\Cot(J^\prime,\mfS^\prime) =
(C_{J^\prime,\mfS^\prime},\Gamma_{J^\prime,\mfS^\prime}),
\quad \Gamma^\prime = \diag(\gamma_1^\prime,\gamma_2^\prime,\gamma_3^\prime).
\end{align*}
to deduce from (\ref{p.CONALT}.\ref{CONOM}) and (iv)
\begin{align*}
\omega^\prime =\,\,&S^\prime(u_{23}^\prime)^{-1}S^\prime(u_{31}^\prime)^{-1} =
S^\prime(\delta_1\varphi(u_{23}))^{-1}S^\prime(\delta_2\varphi(u_{31}))^{-1} \\
=\,\,&\delta_1^{-2}\delta_2^{-2}S(u_{23})^{-1}S(u_{31})^{-1},
\end{align*}
which amounts to
\begin{align}
\label{OMPROM} \omega^\prime = \delta_1^{-2}\delta_2^{-2}\omega.
\end{align}
Now let $v,w \in C = J_{12}$. Then \eqref{OMPROM}
combines with (\ref{p.CONALT}.\ref{CONPRO}) to imply
\begin{align*}
\varphi_{12}(v)\varphi_{12}(w) =\,\,&\omega^\prime\big(\varphi(v)
\times u_{23}^\prime\big) \times \big(u_{31}^\prime \times
\varphi(w)\big) \\
=\,\,&\delta_1^{-2}\delta_2^{-2}\delta_1\delta_2\omega\big(\varphi(v)
\times \varphi(u_{23})\big) \times \big(\varphi(u_{31}) \times
\varphi(w)\big) \\
=\,\,&\delta_1^{-1}\delta_2^{-1}\omega\varphi\big((v \times u_{23})
\times (u_{31} \times w)\big) =
\delta_1^{-1}\delta_2^{-1}\varphi_{12}(vw),
\end{align*}
and we conclude
\begin{align}
\label{MORMUL} \varphi_{12}(vw) =
\delta_1\delta_2\varphi_{12}(v)\varphi_{12}(w) &&(v,w \in C =
J_{12}).
\end{align}
In accordance with \eqref{COTFI}, we now define $\eta :=
\delta_1\delta_2\varphi_{12}\:C \to C^\prime$. Given $v,w \in C$
and applying \eqref{MORMUL}, we deduce $\eta(vw) =
\delta_1\delta_2\varphi_{12}(vw) =
\delta_1\delta_2\varphi_{12}(v)\delta_1\delta_2\varphi_{12}(w) =
\eta(v)\eta(w)$, so $\eta$ is an algebra homomorphism which by
(\ref{p.CONALT}.\ref{CONUN}), (\ref{p.CONALT}.\ref{CONNO}), \eqref{OMPROM} satisfies $\eta(1_C) =
\delta_1\delta_2\varphi(u_{23} \times u_{31}) =
\delta_1\delta_2\varphi(u_{23}) \times \varphi(u_{31}) =
u_{23}^\prime \times u_{31}^\prime = 1_{C^\prime}$,
$n_{C^\prime}(\eta(v)) = -\omega^\prime S^\prime(\delta_1\delta_2\varphi(v)) =
-\delta_1^2\delta_2^2\omega^\prime S^\prime(\varphi(v)) = -\omega
S(v) = n_C(v)$. Summing up we have thus shown that $\eta\:C \to
C^\prime$ is a homomorphism of conic algebras, so condition (i) holds for $(\eta,\Delta)$. We now verify condition \ref{JACOCAT.2}. An application of (\ref{p.CONALT}.\ref{CONMA}) yields $\gamma_1^\prime = -S^\prime(u_{31}^\prime) = -\delta_2^2S^\prime(\varphi(u_{31})) =
-\delta_2^2S(u_{31})$, $\gamma_2^\prime = -S^\prime(u_{23}^\prime) =
-\delta_1^2S^\prime(\varphi(u_{23})) = -\delta_1^2S(u_{23})$, and we
conclude
\begin{align}
\label{GAPRGA} \gamma_1^\prime = \delta_2^2\gamma_1, \quad
\gamma_2^\prime = \delta_1^2\gamma_2, \quad \gamma_3^\prime =
\gamma_3 = 1.
\end{align}
In view of \eqref{COTFI}, \eqref{GAPRGA}, we therefore have
$\delta_2\delta_3\delta_1^{-1}\gamma_1 =
\delta_2\delta_1\delta_2\delta_1^{-1}\gamma_1 =
\delta_2^2\gamma_1 = \gamma_1^\prime$,
$\delta_3\delta_1\delta_2^{-1}\gamma_2 =
\delta_1\delta_2\delta_1\delta_2^{-1}\gamma_2 = \delta_1^2\gamma_2 =
\gamma_2^\prime$, $\delta_1\delta_2\delta_3^{-1}\gamma_3 = \gamma_3
= \gamma_3^\prime$, which comes down to
$\delta_j\delta_l\delta_i^{-1}\gamma_i = \gamma_i^\prime$ for $i = 1,2,3$, hence to $\Gamma^\prime = \Delta^\sharp\Delta^{-1}\Gamma$, and we have shown that condition \ref{JACOCAT.2} holds as well. Thus
$\Cot(\varphi;\delta_1,\delta_2)$ is indeed a morphism of
co-ordinate pairs. It remains to show compatibility with
compositions, so let
\[
\xymatrix{(J,\mfS) \ar[rr]_{(\varphi;\delta_1,\delta_2)} &&
(J^\prime,\mfS^\prime)
\ar[rr]_{(\varphi^\prime;\delta_1^\prime,\delta_2^\prime)} &&
(J^{\prime\prime},\mfS^{\prime\prime}), }
\]
be morphisms of co-ordinated cubic Jordan algebras. Setting
$\Delta^\prime :=
\diag(\delta_1^\prime,\delta_2^\prime,\delta_3^\prime)$,
$\delta_3^\prime := \delta_1^\prime\delta_2^\prime$, gives
\[
\Delta^\prime\Delta =
\diag\big(\delta_1^\prime\delta_1,\delta_2^\prime\delta_2,(\delta_1^\prime
\delta_1)(\delta_2^\prime\delta_2)\big)
\]
Applying \eqref{COTFI}, we therefore conclude
\begin{align*}
\Cot(\varphi_1^\prime;\delta_1^\prime,\delta_2^\prime) \circ
\Cot(\varphi;\delta_1,\delta_2)
=\,\,&(\delta_1^\prime\delta_2^\prime\varphi_{12}^\prime,\Delta^\prime)
\circ (\delta_1\delta_2\varphi_{12},\Delta) \\
=\,\,&(\delta_1^\prime\delta_2^\prime\delta_1\delta_2\varphi_{12}^\prime
\circ \varphi_{12}, \Delta^\prime\Delta) \\
=\,\,&\Big(\big(\delta_1^\prime\delta_1\big)\big(\delta_2^\prime\delta_2\big)
\big(\varphi^\prime \circ
\varphi\big)_{12},\diag\big(\delta_1^\prime\delta_1,\delta_2^\prime\delta_2,
(\delta_1^\prime\delta_1)(\delta_2^\prime\delta_2)\big)\Big) \\
=\,\,&\Cot(\varphi^\prime \circ
\varphi;\delta_1^\prime\delta_1,\delta_2^\prime\delta_2) =
\Cot\big((\varphi^\prime;\delta_1^\prime,\delta_2^\prime) \circ
(\varphi;\delta_1,\delta_2)\big),
\end{align*}
and the assertion follows. Thus we have indeed a functor
\[
\Cot\:\kcocujo \longrightarrow \kcopa.
\]

(e) We employ the standard notation of \ref{ss.CUCOSY} for the
co-ordinated cubic Jordan algebra $(J,\mfS) :=
\bfHer_3(C,\Gamma)$ and then combine
(\ref{p.CONALT}.\ref{ONETWOCO}) with Prop.~\ref{p.PEDEDI} to
conclude $C^\prime = J_{12} = C[12]$ as $k$-modules,
giving the first assertion. Applying (\ref{e.COHER}.\ref{DIACOS}),
(\ref{e.COHER}.\ref{COHER}), (\ref{p.CONALT}.\ref{CONOM}),
(\ref{ss.THERCU}.\ref{MAQIT}), we obtain
\[
\omega = S(1_C[23])^{-1}S(1_C[31])^{-1} =
(-\gamma_2\gamma_3n(1_C))^{-1}(-\gamma_3\gamma_1n(1_C))^{-1},
\]
hence
\begin{align}
\label{OMDIAG} \omega = \gamma_1^{-1}\gamma_2^{-1}\gamma_3^{-2}.
\end{align}
Now, for $v,w \in C$, the multiplication in $C^\prime$ by
(\ref{p.CONALT}.\ref{CONPRO}), \eqref{OMDIAG},
(\ref{ss.EICJM}.\ref{OFFTI}) may be expressed as
\begin{align*}
v[12]w[12] =\,\,&\omega(v[12] \times 1_C[23]) \times (1_C[31] \times
w[12]) = \gamma_1^{-1}\gamma_2^{-1}\gamma_3^{-2}\gamma_2\gamma_1\bar
v[31] \times \bar w[23] \\
=\,\,&\gamma_3^{-2}\gamma_3\overline{\bar w\bar v}[12] =
\gamma_3^{-1}(vw)[12],
\end{align*}
which amounts to
\begin{align}
\label{PROPRI} (vw)[12] = \gamma_3v[12]w[12] &&(v,w \in C).
\end{align}
Similarly, combining (\ref{p.CONALT}.\ref{CONUN}) and
(\ref{e.COHER}.\ref{DIACOS}) with (\ref{ss.EICJM}.\ref{OFFTI})
yields $1_{C^\prime} = 1_C[23] \times 1_C[31] = \gamma_3\bar
1_C[12]$, hence
\begin{align}
\label{UNPRI} 1_{C^\prime} = \gamma_31_C[12].
\end{align}
Also, using (\ref{p.CONALT}.\ref{CONNO}),
(\ref{ss.THERCU}.\ref{MAQIT}), we obtain
\[
n_{C^\prime}(v[12]) = -\omega S(v[12]) =
-\gamma_1^{-1}\gamma_2^{-1}\gamma_3^{-2}(-\gamma_1\gamma_2n_C(v)),
\]
which may be summarized to
\begin{align}
\label{NOPRI} n_{C^\prime}(v[12]) = \gamma_3^{-2}n_C(v) &&(v \in C).
\end{align}
Invoking \eqref{PROPRI}--\eqref{NOPRI}, the map $\psi :=
\psi_{C,\Gamma}$ satisfies
\begin{align*}
\psi(vw) =\,\,&\gamma_3(vw)[12] =\gamma_3v[12]\gamma_3w[12] =
\psi(v)\psi(w), \\
\psi(1_C) =\,\,&\gamma_31_C[12] = 1_{C^\prime}, \\
n_{C^\prime}\big(\psi(v)\big) =\,\,&n_{C^\prime}(\gamma_3v[12]) =
\gamma_3^2\gamma_3^{-2}n_C(v) = n_C(v)
\end{align*}
for all $v,w \in C$. Hence $\psi\:C \to C^\prime$ is
an isomorphism of conic algebras. Finally, setting $\Gamma^\prime =
\diag(\gamma_1^\prime,\gamma_2^\prime,\gamma_3^\prime)$ and
observing (\ref{p.CONALT}.\ref{CONMA}),
(\ref{e.COHER}.\ref{DIACOS}), we deduce $\gamma_1^\prime =
-S(1_C[31]) = \gamma_3\gamma_1n_C(1_C)$, $\gamma_2^\prime =
-S(1_C[23]) = \gamma_2\gamma_3n_C(1_C)$, which yields
\begin{align}
\label{GAPRI} \gamma_1^\prime = \gamma_3\gamma_1, \quad
\gamma_2^\prime = \gamma_2\gamma_3, \quad \gamma_3^\prime = 1.
\end{align}
Combining \eqref{GAPRI} with \eqref{LAMCO}, we obtain
$\lambda_2\lambda_3\lambda_1^{-1}\gamma_1 = \gamma_3\gamma_1 =
\gamma_1^\prime$, $\lambda_3\lambda_1\lambda_2^{-1}\gamma_2 =
\gamma_3\gamma_2 = \gamma_2^\prime$,
$\lambda_1\lambda_2\lambda_3^{-1}\gamma_3 = \gamma_3^{-1}\gamma_3 =
1 = \gamma_3^\prime$, which may be unified to
$\lambda_j\lambda_l\lambda_i^{-1}\gamma_i = \gamma_i^\prime$. Thus,
setting $\Lambda := \Lambda_{C,\Gamma}$, we see $\Gamma^\prime =
\Lambda^\sharp\Lambda^{-1}\Gamma$. Hence $\Psi_{C,\Gamma}$ is an isomorphism of co-ordinate pairs,
as claimed.

(f) This is clear since $\Phi = \Phi_{J,\mfS}$ as defined in Thm.~\ref{t.JACO} matches \emph{exactly} the diagonal co-ordinate system of $\Her_3(C,\Gamma)$ with the co-ordinate system $\mfS$ of $J$.

(g) It suffices to show  
\begin{enumerate}[(i),start=5]
\item \label{JACOCAT.5} that the $\Phi_{J;\mfS}$ for all co-ordinated cubic Jordan algebras $(J,\mfS)$ over $k$ determine an isomorphism 
\[
\bfHer_3 \circ \Cot \overset{\sim} \longrightarrow \mathrm{Id}_{\kcocujo}
\]
of functors from $\kcocujo$ to itself and, similarly, 

\item \label{JACOCAT.6} that the $\Psi_{C,\Gamma}$ for all co-ordinate pairs $(C,\Gamma)$ over $k$ determine an isomorphism
\[
\mathrm{Id}_{\kcopa} \overset{\sim} \longrightarrow \Cot \circ \bfHer_3
\]
of functors from $\kcopa$ to itself.
\end{enumerate}
We begin with \ref{JACOCAT.5}. Given
a co-ordinated cubic norm structure $(J,\mfS)$ over $k$,
Thm.~\ref{t.JACO} shows that
\[
\Phi_{J,\mfS}\:(\bfHer_3 \circ \Cot)(J,\mfS) \overset{\sim}
\longrightarrow (J,\mfS)
\]
is an isomorphism in $\kcocujo$. Thus, by \cite[Def.~1.3,
p.~23]{MR1009787}, the proof will be complete once we have shown that
any homomorphism $(\varphi;\delta_1,\delta_2)\:(J,\mfS)
\to (J^\prime,\mfS^\prime)$ of co-ordinated cubic Jordan algebras
gives rise to a commutative diagram
\begin{align}
\vcenter{ \label{d.HERCOT} \xymatrix{ (\bfHer_3 \circ
\Cot)(J,\mfS) \ar[rr]_(0.65){\Phi_{J,\mfS}}^(0.65){\cong}
\ar[dd]_{(\bfHer_3 \circ
\Cot)(\varphi;\delta_1,\delta_2)} && (J,\mfS) \ar[dd]^{(\varphi;\delta_1,\delta_2)} \\ \\
(\bfHer_3 \circ \Cot)(J^\prime,\mfS^\prime)
\ar[rr]_(0.65){\Phi_{J^\prime,\mfS^\prime}}^(0.65){\cong} &&
(J^\prime,\mfS^\prime). }}
\end{align}
In order to establish the commutativity of \eqref{d.HERCOT}, we
write $S$ (resp. $S^\prime$) for the quadratic trace of $J$
(resp. $J^\prime)$; moreover, we put $\mfS =
(e_1,e_2,e_3,u_{23},u_{31})$, $\mfS^\prime =
(e_1^\prime,e_2^\prime,e_3^\prime,u_{23}^\prime,u_{31}^\prime)$.

Applying \eqref{PHIJAY}, we first obtain $(\varphi;\delta_1,\delta_2) \circ \Phi_{J,\mfS}
= (\varphi;\delta_1,\delta_2) \circ (\phi_{J,\mfS};1,1)$, hence
\begin{align}
\label{VARPHIPHI} (\varphi;\delta_1,\delta_2) \circ \Phi_{J,\mfS}
= (\varphi \circ \phi_{J,\mfS};\delta_1,\delta_2).
\end{align}
On the other hand, setting
\begin{align}
\label{VARPHIDEL} \eta := \delta_1\delta_2\varphi_{12}, \quad \Delta
:= \diag(\delta_1,\delta_2,\delta_1\delta_2) \in \Diag_3(k)^\times,
\end{align}
we deduce from \eqref{COHE}), \eqref{VARPHIDEL} that
\begin{align*}
\Phi_{J^\prime,\mfS^\prime} \circ
\bfHer_3\big(\Cot(\varphi;\delta_1,\delta_2)\big)
=\,\,&(\phi_{J^\prime,\mfS^\prime};1,1) \circ
\bfHer_3(\eta,\Delta) \\
=\,\,&(\phi_{J^\prime,\mfS^\prime};1,1) \circ
\big(\Her_3(\eta,\Delta);\delta_1,\delta_2\big) \\
=\,\,&\big(\phi_{J^\prime,\mfS^\prime} \circ
\Her_3(\eta,\Delta);\delta_1,\delta_2\big),
\end{align*}
so in view of \eqref{VARPHIPHI}, it suffices to show
\begin{align}
\label{FIFIHER} \varphi \circ \phi_{J,\mfS} =
\phi_{J^\prime,\mfS^\prime} \circ \Her_3(\eta,\Delta).
\end{align}
In order to do so, we write $(C,\Gamma) := \Cot(J,\mfS)$ (resp.
$(C^\prime,\Gamma^\prime) =
\Cot(J^\prime,\mfS^\prime)$). Combining
(\ref{t.JACO}.\ref{CONPHI}) with condition (a), (i) of Exc.~\ref{pr.DIAFIX}, we obtain $\varphi \circ \phi_{J,\mfS}(e_{ii}) = \varphi(e_i) = e_i^\prime =
\phi_{J^\prime,\mfS^\prime}(e_{ii}) =
\phi_{J^\prime,\mfS^\prime} \circ \Her_3(\eta,\Delta)(e_{ii})$
for $i = 1,2,3$, so it remains to show that both sides of
\eqref{FIFIHER} agree on $v_i[jl]$, for all $v_i \in C =
J_{12}$, $i = 1,2,3$. Since $\eta$ as defined in
\eqref{VARPHIDEL}, being an isomorphism of conic algebras, preserves conjugations, so does
every scalar multiple of it, in particular $\varphi_{12}$. With this
in mind, we apply (\ref{t.JACO}.\ref{CONPE}), \eqref{VARPHIDEL} and obtain
\begin{align*}
\varphi \circ \phi_{J,\mfS}(v_1[23])
=\,\,&-S(u_{31})^{-1}\varphi(u_{31} \times \bar v_1) =
-S^\prime\big(\varphi(u_{31})\big)^{-1}\varphi(u_{31}) \times
\varphi(\bar v_1) \\
=\,\,&-S^\prime(\delta_2^{-1}u_{31}^\prime)^{-1}\delta_2^{-1}u_{31}^\prime
\times \varphi_{12}(\bar v_1) =
-\delta_2^2\delta_2^{-1}S^\prime(u_{31}^\prime)^{-1}u_{31}^\prime
\times \overline{\varphi_{12}(v_1)} \\
=\,\,&-\delta_2S^\prime(u_{31}^\prime)^{-1}u_{31}^\prime \times
\overline{\varphi_{12}(v_1)} =
\delta_2\phi_{J^\prime,\mfS^\prime}\big(\varphi_{12}(v_1)[23]\big)
\\
=\,\,&\phi_{J^\prime,\mfS^\prime}\big(\delta_1^{-1}(\delta_1\delta_2\varphi_{12})(v_1)[23]\big)
= \phi_{J^\prime,\mfS^\prime} \circ \Her_3(\eta,\Delta)(v_1[23]).
\end{align*}
Similarly,
\begin{align*}
\varphi \circ \phi_{J,\mfS}(v_2[31])
=\,\,&-S(u_{23})^{-1}\varphi(u_{23} \times \bar v_2) =
-S^\prime\big(\varphi(u_{23})\big)^{-1}\varphi(u_{23}) \times
\varphi(\bar v_2) \\
=\,\,&-S^\prime(\delta_1^{-1}u_{23}^\prime)^{-1}\delta_1^{-1}u_{23}^\prime
\times \varphi_{12}(\bar v_2) =
-\delta_1^2\delta_1^{-1}S^\prime(u_{23}^\prime)^{-1}u_{23}^\prime
\times \overline{\varphi_{12}(v_2)} \\
=\,\,&-\delta_1S^\prime(u_{23}^\prime)^{-1}u_{23}^\prime \times
\overline{\varphi_{12}(v_2)} =
\delta_1\phi_{J^\prime,\mfS^\prime}\big(\varphi_{12}(v_2[31])\big)
\\
=\,\,&\phi_{J^\prime,\mfS^\prime}\big(\delta_2^{-1}(\delta_1\delta_2\varphi_{12})(v_2[31])\big)
= \phi_{J^\prime,\mfS^\prime} \circ \Her_3(\eta,\Delta)(v_2[31]).
\end{align*}
And, finally,
\begin{align*}
\varphi \circ \phi_{J,\mfS}(v_3[12]) =\,\,&\varphi(v_3) =
\varphi_{12}(v_3) =
\phi_{J^\prime,\mfS^\prime}\big(\varphi_{12}(v_3)[12]\big) \\
=\,\,&\phi_{J^\prime,\mfS^\prime}\big((\delta_1\delta_2)^{-1}(\delta_1\delta_2\varphi_{12})(v_3)[12]\big)
= \phi_{J^\prime,\mfS^\prime} \circ \Her_3(\eta,\Delta)(v_3[12]).
\end{align*}
Summing up, this completes the proof of \eqref{FIFIHER}, and we have
shown that the diagram \eqref{d.HERCOT} commutes. Thus \ref{JACOCAT.5} holds.

It remains to establish \ref{JACOCAT.6}. To this end, we let $(\eta,\Delta)\:(C,\Gamma)
\to (C^\prime,\Gamma^\prime))$ be a morphism of co-ordinate
triples and must show that the diagram
\begin{align}
\vcenter{ \label{d.COTHER} \xymatrix{ (C,\Gamma)
\ar[rr]_(0.45){\Psi_{C,\Gamma}}^(0.45){\cong}
\ar[dd]_{(\eta,\Delta)} &&
(\Cot \circ\,\bfHer_3)(C,\Gamma) \ar[dd]^{(\Cot \circ\,\bfHer_3)(\eta,\Delta)} \\ \\
(C^\prime,\Gamma^\prime)
\ar[rr]_(0.45){\Psi_{C^\prime,\Gamma^\prime}}^(0.45){\cong}
&& (\Cot \circ\,\bfHer_3)(C^\prime,\Gamma^\prime). }}
\end{align}
commutes. We begin by making the specifications
\begin{align}
\label{GAPR} \Gamma =\,\,&\diag(\gamma_1,\gamma_2,\gamma_3), \quad
\Gamma^\prime =
\diag(\gamma_1^\prime,\gamma_2^\prime,\gamma_3^\prime), \\
\label{CEHA} (\hC,\hGamma)
:=\,\,&\Cot\big(\bfHer_3(C,\Gamma)\big), \quad
\hGamma = \diag(\hgamma_1,\hgamma_2,\hgamma_3), \\
\label{CEHAPR} (\hC^\prime,\hGamma^\prime)
:=\,\,&\Cot\big(\bfHer_3(C^\prime,\Gamma^\prime)\big),
\quad \hGamma^\prime =
\diag(\hgamma_1^\prime,\hgamma_2^\prime,\hgamma_3^\prime).
\end{align}
We abbreviate $\Psi_{C,\Gamma} = (\psi,\Lambda)$,
$\Psi_{C^\prime,\Gamma^\prime} =
(\psi^\prime,\Lambda^\prime)$ and deduce from
\eqref{LAMCO}) that
\begin{align}
\label{LAM} \Lambda =\,\,&\diag(\lambda_1,\lambda_2,\lambda_3),\quad
\lambda_1 = \lambda_2 = 1, \quad \lambda_3 = \gamma_3, \\
\label{LAMPR} \Lambda^\prime
=\,\,&\diag(\lambda_1^\prime,\lambda_2^\prime,\lambda_3^\prime),
\quad \lambda_1^\prime = \lambda_2^\prime = 1, \quad
\lambda_3^\prime = \gamma_3^\prime.
\end{align}
Hence $\Psi_{C^\prime,\Gamma^\prime} \circ (\eta,\Delta) =
(\psi^\prime,\Lambda^\prime) \circ (\eta,\Delta)$, and
\eqref{COMCOTR}) implies
\begin{align}
\label{DORI} \Psi_{C^\prime,\Gamma^\prime} \circ
(\eta,\Delta) = (\psi^\prime \circ \eta,\Lambda^\prime\Delta),
\end{align}
where \eqref{PSICO} yields
\begin{align}
\label{PSET} (\psi^\prime \circ \eta)(v) = \gamma_3^\prime
\eta(v)[12] &&(v \in C).
\end{align}
On the other hand, writing $\Delta =
\diag(\delta_1,\delta_2,\delta_3)$ and invoking \eqref{LAMPR}, we
obtain
\begin{align}
\label{LAPRDEL}\Lambda^\prime\Delta =
\diag(\delta_1,\delta_2,\gamma_3^\prime\delta_3).
\end{align}
Next we use \eqref{COHE} to compute
\begin{align*}
\Cot\big(\bfHer_3(\eta,\Delta)\big)
=\,\,&\Cot\Big(\big(\Her_3(\eta,\Delta);\delta_1,\delta_2\big)\Big),
\end{align*}
where 
\begin{align}
\label{VERRI} \Cot\big(\bfHer_3(\eta,\Delta)\big) = (\heta,\hDelta),
\quad \heta := \delta_1\delta_2\Her_3(\eta,\Delta)_{12}, \quad
\hDelta := \diag (\delta_1,\delta_2,\delta_1\delta_2).
\end{align}
From $(\heta,\hDelta) \circ \Psi_{C,\Gamma} = (\heta,\hDelta)
\circ (\psi,\Lambda)$ we therefore deduce
\begin{align}
\label{RIDO} (\Cot \circ\,\bfHer_3)(\eta,\Delta) \circ
\Psi_{C,\Gamma} = (\heta \circ \psi,\hDelta\Lambda).
\end{align}
Given $v \in C$, we observe $\Gamma^\prime =
\Delta^\sharp\Delta^{-1}\Gamma$ (since $(\eta,\Delta)$ is a morphism
in $\kcopa$), which yields
\begin{align*}
(\heta \circ \psi)(v) =\,\,&\heta\big(\psi_{C,\Gamma}(v)\big) =
\gamma_3\heta(v[12]) =
\delta_1\delta_2\gamma_3\Her_3(\eta,\Delta)_{12}(v[12]) \\
=\,\,&\delta_1\delta_2\delta_3^{-1}\gamma_3\eta(v)[12] =
\gamma_3^\prime\eta(v)[12] = (\psi^\prime \circ \eta)(v).
\end{align*}
Thus $\heta \circ \psi = \psi^\prime \circ \eta$. Moreover, by
\eqref{VERRI}, \eqref{LAM}, \eqref{LAMPR},
\begin{align*}
\hDelta\Lambda
=\,\,&\diag(\delta_1,\delta_2,\delta_1\delta_2)\diag(1,1,\gamma_3) =
\diag(\delta_1,\delta_2,\delta_1\delta_2\gamma_3) \\
=\,\,&\diag(\delta_1,\delta_2,\delta_3\gamma_3^\prime) =
\diag(1,1,\gamma_3^\prime)\diag(\delta_1,\delta_2,\delta_3) =
\Lambda^\prime\Delta.
\end{align*}
Comparing now \eqref{DORI} with \eqref{RIDO}, we see that the
diagram \eqref{d.COTHER} commutes, which completes the proof of \ref{JACOCAT.6}.
\end{sol}


\solnsec{Section~\ref{s.JOALTH}}

\begin{sol}{pr.WEIMOD} \label{sol.WEIMOD}  (a) Define polynomial laws $f\:M \to k$ and $g\:M \to M$ by $f(u) := q(u)^2$, $g(u) := q(u)u$ for all $u \in M_R$, $R \in \kalg$. Note that $f$ is a bi-quadratic form, while $g$ is homogeneous of degree $3$. Combining the first and second order product rules (\ref{ss.DICA}.\ref{DIPR}), (\ref{ss.DICA}.\ref{SEDIPR}) with (\ref{ss.DICA}.\ref{DIQU}), (\ref{ss.DICA}.\ref{DISQ}), we conclude that
\begin{align*}
(Df)(u,v) =\,\,&2q(u)q(u,v) \\
(D^2f)(u,v) =\,\,&q(u,v)^2 + 2q(u)q(v), \\
(Dg)(u,v) =\,\,&q(u,v)u + q(u)v 
\end{align*}
hold strictly for $u,v \in M$. Fixing $v \in M$, the cubic form $h_v\:M \to k$ defined by $h_v(u) := 2q(u)q(u,v)$ for $u \in M_R$, $R \in \kalg$ satisfies
\[
(Dh_v)(u,w) = 2\big(q(u,v)q(u,w) + q(u)q(v,w)\big)
\]
strictly for $u,w \in M$. It now follows from a repeated application of Exc.~\ref{pr.POLALOC}~(b) that \eqref{WEIR} holds strictly if and only if \eqref{WEIR}--\eqref{QUVEU} hold over $k$.

Non-trivial weird quadratic modules exist, as the following simple example shows. Let $k_0$ be a commutative ring, $I$ a $k_0$-module and $k := k_0 \oplus I$ the split-null extension of $k_0$ by $I$. Regard any $k_0$-module $M_0$ as a $k$-module $M$ by letting $I$ act trivially on $M_0$. We claim that \emph{any quadratic form $q\:M \to k$ taking values in $I$ makes $(M,q)$ a weird quadratic module over $k$.} Indeed, since $I$ kills $M$ and $I^2 = \{0\}$ in $k$, not only \eqref{WEIR} but also \eqref{TWOQU}--\eqref{QUVEU} hold over $k$, which proves the assertion.

(b) Localizing if necessary and invoking Exc.~\ref{pr.POLALOC}~(a), we may assume that $M$ is free of rank at least $2$. Let $(e_i)_{i \in I}$ be a basis of $M$ over $k$. Given $i,j \in I$ distinct, \eqref{QUVEU} implies $q(e_i,e_j)e_i + q(e_i)e_j = 0$, hence $q(e_i) = q(e_i,e_j) = 0$. Thus $q = 0$.

(c) Put $\bfM_0 := (M_0,q_0,e_0)$ with $M_0 := ke_0$ and $q_0\:M_0 \to k$ being the quadratic form given by $q_0(\xi_0e_0) := \xi_0^2$ for $\xi_0 \in k$. Then $\bfM_0$ is a pointed quadratic $k$-module whose (bi-)linear trace and conjugation have the form $t_0(\xi_0e_0) = 2\xi_0$, $t_0(\xi_0e_0,\eta_0e_0) = 2\xi_0\eta_0$ and $\overline{\xi_0e_0} = \xi_0e_0$ for $\xi_0,\eta_0 \in k$. Let $M_0$ act bilinearly on $M_1 := M$ by $e_0\pt u := u$ for all $u \in M$ and define a quadratic map $Q\:M_1 \to M_0$ by $Q(u) := q(u)e_0$ for all $u \in M$. Then $\bfM_1 := (M_1,\pt,Q)$ is a Peirce-one extension of $\bfM_0$, and comparing (\ref{ss.BUEL}.\ref{BUMO})--(\ref{ss.BUEL}.\ref{BUNO}) with \eqref{EMEZ}--\eqref{EMEN} we see that $X = X(\bfM_0,\bfM_1)$ is indeed a cubic array such that \eqref{EMEB}--\eqref{EMELB} hold. By Prop.~\ref{p.ADPEMO}, therefore, $X$ is a cubic norm structure if and only if (\ref{p.ADPEMO}.\ref{ADEZ})--(\ref{p.ADPEMO}.\ref{ADQUX}) are strictly fulfilled. After the specifications described above (\ref{p.ADPEMO}.\ref{ADEZ})--(\ref{p.ADPEMO}.\ref{ADQU}), (\ref{p.ADPEMO}.\ref{ADTEX}) are trivial, while (\ref{p.ADPEMO}.\ref{QUQUX}), (\ref{p.ADPEMO}.\ref{ADQUX}) amount to the strict validity of \eqref{WEIR}. Thus $X$ is a cubic norm structure if and only if $(M,q)$ is weird. The remaining assertions in (c) now follow from the final statement of Prop.~\ref{p.ADPEMO}.

(d) Basically, the arguments in the solution to (c) will be read backwards. Write $X$ for the cubic norm structure underlying $J$. By definition, $e := 1 - c$ is an elementary idempotent in $J$ having $e_0 := 1 - e = c$. With $M := J_1(e)$, the Peirce decomposition of $J$ relative to $e$ therefore attains the form
\begin{align}
\label{WEIRX} X = ke \oplus M \oplus ke_0,
\end{align} 
where $ke$, $ke_0$ are free $k$-modules of rank $1$. We put $M_0 := ke_0$ and define $q_0\:M_0 \to k$ by $q(\xi_0e_0) := \xi_0^2$ for $\xi_0 \in k$ to obtain a pointed quadratic module $\bfM_0 = (M_0,q_0,e_0)$ whose (bi-)linear trace and conjugation are respectively given by $t_0(\xi_0e_0) = 2\xi_0$, $t_0(\xi_0e_0,\eta_0e_0) = 2\xi_0\eta_0$ and $\overline{\xi_0e_0} = \xi_0e_0$ for $\xi_0,\eta_0 \in k$. Applying the formalism of \ref{ss.CUPE}, there exists a quadratic form $q\:M \to k$ such that the quadratic map $Q$ of (\ref{ss.CUPE}.\ref{EXQU}) may be written as $Q(u) = q(u)e_0$ for all $u \in M$. With these specifications, and in view of (\ref{t.PEDESI}.\ref{JONE}), formulas (\ref{ss.CUPE}.\ref{ALPHEXX}), (\ref{ss.CUPE}.\ref{ENEXX}) are converted to \eqref{EMESH}, \eqref{EMEN}, which therefore hold strictly in $X$. But $X$ is a cubic norm structure . Hence (c) implies that $(M,q)$ is a weird quadratic module. By its very definition, the cubic norm structure underlying $\pJ := \Jcub(M,q)$  identifies canonically with $X$. Thus $J \cong \pJ$ under an isomorphism matching $e_0 = c \in J$ with $e_0 \in \pJ$. 
\end{sol} 

\begin{sol}{pr.UNUNORED} \label{sol.UNUNORED}  Let $X,\pX$ be two cubic norm structures over $k$ such that $J(X) = J = J(\pX)$ as abstract Jordan algebras. We must show $X = \pX$ and clearly have $X = J = \pX$ as $k$-modules as well as $1_X = 1_J = 1_{\pX}$. Put $N := N_X - N_{\pX}$, $T := T_X - T_{\pX}$, $S := S_X - S_{\pX}$. Then (\ref{t.CUNOJO}.\ref{UNIVT}) implies that
\begin{align}
\label{TEXX} T(x)x^2 - S(x)x + N(x)1_J = 0
\end{align} 
holds strictly for all $x \in J$. Let $\mfp \in \Spec(k)$ and $K$ an algebraic closure of $k(\mfp)$. Then the composite map
\[
\xymatrix{\vart:k \ar[r] & k/\mfp \ar[r] & k(\mfp) \ar[r] & K}
\]
has kernel $\mfp$ and makes $K$ a $k$-algebra. Let $\vph \in \{T,S,N\}$ and $x \in J_K$. If $\Xi_K^ J(x) \neq 0$, then \eqref{TEXX} shows $\vph_K(x) = 0$. Thus $((\vph \otimes K)(\Xi^J \otimes K))_K = 0$, and Exc.~\ref{pr.POLINF}~(b) implies $\vph \otimes K = 0$, in particular $\vph_K = 0$ as a set map $J_K \to K$. Given $x \in J$, we therefore conclude $0 = \vph_K(x_K) = \vph(x)_K = \vart (\vph(x))$, hence $\vph(x) \in \mfp$. This is true for all $\mfp \in \Spec(k)$, i.e.,
\[
\vph(x) \in \bigcap_{\mfp \in \Spec(k)} \mfp = \Nil(k) = \{0\}, 
\]
and we have shown
\begin{align}
\label{TEXNIL} T_X(x) = T_{\pX}(x), \quad S_X(x) = S_{\pX}(x), \quad N_X(x) = N_{\pX}(x) &&(x \in J),
\end{align} 
which by (\ref{ss.BACUNO.fig}.\ref{ADJJA}) implies
\begin{align}
\label{EXSHAX} x^{\sharp_X} = x^{\sharp_{\pX}} &&(x \in J)
\end{align}
Writing $\times_X$, $\times_{\pX}$ for the bilinearizations of $\sharp_X$, $\sharp_{\pX}$, respectively, we deduce from \eqref{TEXNIL}, \eqref{EXSHAX} that
\[
(X,1_X,\sharp_X,\times_X,T_X,N_X) = (\pX,1_{\pX},\sharp_{\pX}, \times_{\pX}, T_{\pX},N_{\pX})
\]
as rational cubic norm structures in the sense of Exc.~\ref{pr.RACUNO}. This exercise therefore implies $X = \pX$ as ordinary cubic norm structures.
\end{sol}

\begin{sol}{pr.NILEMZER} \label{sol.NILEMZER}  (a) Put $\mfa := \Nil(k)$, $I := \Nil(J)$. If $x \in I$, then $T(x,J) + T(x^\sharp,J) \subseteq \mfa$ by Exc.~\ref{pr.CUBNIL}. In particular, $\xi = T(x,e) \in \mfa$ by (\ref{ss.BUEL}.\ref{BUTBIL}). Moreover, $x_0 = U_{e_0}x \in J_0(e) \cap I = \Nil(J_0(e))$ by Exc.~\ref{pr.NIRAPE}, and Exc.~\ref{pr.NIPOQUA} implies $q_0(x_0),q_0(x_0,y_0) \in \mfa$ for all $y_0 \in M_0$. Next, by (\ref{ss.BUEL}.\ref{BUTBIL}) we have $t_0(Q(x_1,y_1)) = T(x,y_1) \in \mfa$ for all $y_1 \in M_1$ and, finally, 
\[
q_0\big(Q(x_1),y_0\big) = t_0\big(Q(x_1),\bar y_0\big) = -T(x^\sharp,\bar y_0) + \xi q_0(x_0,\bar y_0) \in \mfa 
\]
for all $y_0 \in M_0$ by (\ref{ss.BUEL}.\ref{BUSHA}), (\ref{ss.BUEL}.\ref{BUTBIL}) and by what we have seen before. Hence all the quantities of \eqref{XIQUZ} belong to $\mfa$. Conversely, let this be so. Then $t_0(x_0,y_0) = q_0(x_0,\bar y_0) \in \mfa$ for all $y_0 \in M_0$ by (\ref{ss.JOPOID.fig}.\ref{JONOCO})  and $T(x,J) \subseteq \mfa$ by (\ref{ss.BUEL}.\ref{BUTBIL}). For $y = (\eta e,y_1,y_0)$, $\eta \in k$, $y_i \in M_i$, $i = 0,1$, we apply (\ref{ss.BUEL}.\ref{BUSHA}), (\ref{ss.BUEL}.\ref{BUTBIL}) to conclude
\[
T(x^\sharp,y) = q_0(x_0)\eta + \xi q_0(x_0,y_0) - q_0\big(Q(x_1),\bar y_0\big) + t_0\big(Q(x_0\pt x_1,y_1)\big).
\]
By hypothesis, all summands on the right with the possible exception of the last belong to $\mfa$. But since (\ref{p.ADPEMO}.\ref{ADTEX}) implies that 
\[
t_0\big(Q(x_0\pt x_1,y_1)\big) = t_0\big(x_0,Q(x_1,y_1)\big) = q_0\big(x_0,\overline{Q(x_1,y_1)}\big) 
\]
belongs to $\mfa$ as well, so does $T(x^\sharp,y)$. Finally, (\ref{ss.BUEL}.\ref{BUNO}) combined with \eqref{XIQUZ} yields $N(x) \in \mfa$, and summing up, we have shown $x \in I$. 

(b) (i) $\Rightarrow$ (ii). By Exc.~\ref{pr.NIRAPE}, the nil radical of $J_2(e) \cong k^{(+)}$ is zero. Hence $k$ is reduced and, by (a), $q_0$ is non-degenerate.

(ii) $\Rightarrow$ (iii). By (a), (\ref{p.ADPEMO}.\ref{QUQUX}) and (ii), an element $x$ as in \eqref{EXXI} belongs to $\Nil(J)$ if and only if $\xi = 0$, $x_0 = Q(x_1) = 0$, and $t_0(Q(x_1,y_1)) = 0$ for all $y_1 \in M_1$. Hence (iii) holds.

(iii) $\Rightarrow$ (i). Obvious.
\end{sol}

\begin{sol}{pr.EXTEL} \label{sol.EXTEL}  From condition (iii) of (b) in Exc.~\ref{pr.NILEMZER} we conclude $Q(x_1) \neq 0$. We clearly have $T(\pe) = 1$, while (\ref{ss.BUEL}.\ref{BUSHA}) implies
\[
e^{\prime\sharp} = \Big(q_0\big(Q(x_1)\big),-Q(x_1)\pt x_1,-\overline{Q(x_1)} - Q(x_1)\Big),
\]
where $q_0(Q(x_1)) = 0$ by (\ref{p.ADPEMO}.\ref{QUQUX}) and $Q(x_1)\pt x_1 = t_0(Q(x_1))x_1$ (by (\ref{p.ADPEMO}.\ref{ADQUX})) $= 0$. Hence $e^{\prime\sharp} = 0$, i.e., $\pe$ is an elementary idempotent. Put $\pe_0 := 1_J - \pe = \big(0,-x_1,e_0 + Q(x_1)\big)$. For $y = (\eta e,y_1,y_0)$, $\eta \in k$, $y_i \in M_i$, $i = 0,1$, we apply (\ref{ss.BUEL}.\ref{BUTBIL}) and obtain
\begin{align*}
T(\pe_0,y) =\,\,&T\Big(\big(0,-x_1,e_0 + Q(x_1)\big),(\eta e,y_1,y_0)\Big) \\
=\,\,&t_0\big(Q(x_1),y_0\big) = q_0\big(Q(x_1),\bar y_0\big).
\end{align*} 
Since $Q(x_1) \neq 0$, $q_0(Q(x_1)) = 0$ and $q_0$ is non-degenerate by Exc.~\ref{pr.NILEMZER}~(b), some $y_0 \in M_0$ has $T(\pe_0,y) = q_0(Q(x_1),\bar y_0) \neq 0$. Now write $y = \peta\pe + \py_1 + \py_0$, $\peta \in k$, $\py_i \in J_i(\pe)$, $i = 0,1$. Since the Peirce components of $J$ relative to $\pe$ are orthogonal with respect to the bilinear trace of $J$, we may apply Cor.~\ref{c.FAULE} to $\pe$ to obtain $0 \neq T(\pe_0,y) = T(\pe_0,\py_0) = t^\prime_0(\py_0)$, and the problem is solved. 
\end{sol}

\begin{sol}{pr.SEPCUB} \label{sol.SEPCUB}  For $r > 2$, we are in case (i), by Thm.~\ref{t.SEPDETH}.  For $r = 1$, we are in case (iii), so we may assume $r = 2$. Thanks to Exc.~\ref{pr.CUMOPO} it suffices to show that \emph{there exists a pointed quadratic space $(M,q,e)$ over $k$ satisfying $J = J(M,q,e)$ as abstract Jordan algebras.} By faithfully flat descent, combining Thm.~\ref{t.FAITEL} with Exercises.~\ref{pr.FFPOLA} and \ref{pr.POJOCAT}, there is no harm in assuming that $J$ contains an elementary idempotent $c$. Put $d := 1_J - c$. For all $\mfp \in \Spec(k) $, the idempotent $c(\mfp) \in J(\mfp)$ continues to be elementary, forcing $d(\mfp) \neq 0$ by \ref{ss.COELID}. By Lemma~\ref{l.FAITH}, therefore, $d \in J$ is unimodular. Thus, by the Peirce rules, $kc \times kd$ is not only free of rank $2$ but also a direct summand of $J$ as a $k$-module. This implies $J = kc \times kd$ as a direct product of ideals, and one checks that the hyperbolic quadratic form $q\:J \to k$, $\xi c + \eta d \mapsto \xi\eta$, gives rise to a pointed quadratic space $(M,q,e)$, $M := J$, $e := 1_J$, satisfying $J = J(M,q,e)$ as abstract Jordan algebras.
\end{sol}


\solnsec{Section~\ref{s.FREAL}}

\begin{sol}{pr.CRITREG} \label{sol.CRITREG}  If $\vph_K\:M_K \to N_K$ is bijective for all fields $K \in \kalg$, then in particular so is $\vph(\mfp)\:M(\mfp) \to N(\mfp)$, for all $\mfp \in \Spec(k)$. Then a standard consequence of Nakayama's lemma (as in the solution to \ref{pr.LOCSIALG}) implies that $\vph_\mfp\:M_\mfp \to N_\mfp$ is bijective for all $\mfp \in \Spec(k)$. Thus $\vph\:M \to N$ is bijective. The converse is obvious.
	
Now suppose $J$ is a cubic Jordan algebra over $k$ that is finitely generated projective as a $k$-module and makes $T_K\:M_K \times M_K \to K$ a non-degenerate symmeetric bilinear form over $K$, for all fields $K \in \kalg$. Since $J_K$ has finite dimension over $K$, the natural map $J \to J^\ast$ (the dual module of $J$) induced by $T$ becomes an isomorphism when extending scalars from $k$ to $K$. By the first part of the exercise, therefore, the natural map in question must be an isomorphism to begin with, i.e., $J$ is regular.
\end{sol}

\begin{sol}{pr.CRITSEP} \label{sol.CRITSEP}  Since passing from $(C,\Gamma)$ to $J$ is compatible with base change, the specific nature of the conditions imposed on $J$ (resp. $C$) allows us to assume that $k = K$ is a field, in which case we have to show that $J$ is semi-simple, i.e., $\Nil(J) = \{0\}$, if and only if the quadratic form $n_C$ is non-degenerate. But this follows immediately from Exercises~\ref{pr.HEMONI} and \ref{pr.NIRA} combined:
\[
\Nil(J) = \{0\} \; \Leftrightarrow \;\Nil(C) = \{0\} \; \Leftrightarrow \;\text{$n_C$ is non-degenerate.}
\]
\end{sol}

 \begin{sol}{pr.JOMARK} \label{sol.JOMARK}  (i) $\Rightarrow$ (ii). By (\ref{ss.THERCU}.\ref{MADJ}), condition (i) is equivalent to
\begin{align}
\label{RKONE} \alpha_j\alpha_l = \gamma_j\gamma_ln_C(v_i), \quad \alpha_iv_i = \gamma_i\overline{v_jv_l} &&(i = 1,2,3),
\end{align}	
which in turn implies 
\[
\gamma_1\gamma_2\gamma_3v_1(v_2v_3) = (\gamma_2\gamma_3v_1)(\gamma_1v_2v_3) = \gamma_2\gamma_3\alpha_1v_1\bar v_1 = \alpha_1\gamma_2\gamma_3n_C(v_1)1_C = \alpha_1\alpha_2\alpha_31_C, 
\]
and (ii) holds.
	
(ii) $\Rightarrow$ (iii). Assume that (ii) holds. We first claim
\begin{align}
\label{TRAS} \gamma_1\gamma_2\gamma_3(v_1v_2)v_3 = \alpha_1\alpha_2\alpha_31_C.
\end{align}
If $v_1 = 0$, then $\alpha_2\alpha_3 = \gamma_2\gamma_3n_C(v_1) = 0$, forcing $\alpha_1\alpha_2\alpha_3 = 0$, and \eqref{TRAS} holds. Similarly, if $v_2 = 0$, then $\alpha_3\alpha_1 = \gamma_3\gamma_1n_C(v_2) = 0$, and we have again $\alpha_1\alpha_2 \alpha_3 = 0$, hence \eqref{TRAS}. Thus we may assume $v_1 \neq 0 \neq v_2$. Multiplying the second equation of (ii) from the left first by $\bar v_1$, then by $\bar v_2$ and invoking Krimse's identities as well as the fact that $C$ is multiplicative (Prop.~\ref{p.ALMU}) and its conjugation is an involution (Prop.~\ref{p.MULAD}~(a)), we obtain
\begin{align*}
\gamma_1\gamma_2\gamma_3n_C(v_1)v_2v_3 =\,\,&\alpha_1\alpha_2\alpha_3\bar v_1, \\
\gamma_1\gamma_2\gamma_3n_C(v_1v_2)v_3 =\,\,&\gamma_1\gamma_2\gamma_3n_C(v_1)n_C(v_2)v_3 = \alpha_1\alpha_2\alpha_3\bar v_2\bar v_1 = \alpha_1\alpha_2\alpha_3\overline{v_1v_2}. 
\end{align*}
On the other hand, multiplying the left-hand side of \eqref{TRAS} by $\overline{v_1v_2}$ from the left gives
\[
\gamma_1\gamma_2\gamma_3n_C(v_1v_2)v_3 = \gamma_1\gamma_2\gamma_3\overline{v_1v_2}\big((v_1v_2)v_3\big).
\]
In the expressions on the very right of the last two displayed equations, we may cancel the factor $\overline{v_1v_2}$ since $C$ has no zero divisors and arrive at \eqref{TRAS}. Applying \eqref{TRAS}, we now obtain
\begin{align*}
\gamma_1\gamma_2\gamma_3\bar v_3\big(v_3(v_1v_2)\big) = \,\,&\gamma_1\gamma_2\gamma_3n_C(v_3)v_1v_2 = \gamma_1\gamma_2\gamma_3\big((v_1v_2)v_3\big)\bar v_3 = \alpha_1\alpha_2\alpha_3\bar v_3.
\end{align*}
Hence
\begin{align}
\label{SHIFT} \gamma_1\gamma_2\gamma_3v_3(v_1v_2) = \alpha_1\alpha_2\alpha_31_C
\end{align}
provided $v_3 \neq 0$. But if $v_3 = 0$, then the first equation of (ii) yields $\alpha_1\alpha_2 = \gamma_1\gamma_2n_C(v_3) = 0$, hence $\alpha_1\alpha_2\alpha_3 = 0$, and \eqref{SHIFT} holds in full generality. We have thus arrived at both displayed equations of (iii) for $(ijl) = (123)$ and at the second one for $(ijl) = (312)$, i.e., for the right shift of $(123)$. Repeating the preceding arguments if necessary, we may therefore conclude that (iii) holds. 
	
Finally, suppose $k$ is an integral domain and $\gamma_i \neq 0$  for $i = 1,2,3$. Then we must show (iii) $\Rightarrow$ (i), so assuming (iii), we must establish \eqref{RKONE}. For $i = 1,2,3$ we obtain
\begin{align*}
\gamma_j\gamma_lv_i(\gamma_iv_jv_l) =\,\,&\gamma_1\gamma_2\gamma_3v_i(v_jv_l) = \alpha_1\alpha_2\alpha_31_C \\
=\,\,&\alpha_i\alpha_j\alpha_l1_C = \alpha_i\gamma_j\gamma_ln_C(v_i)1_C = \gamma_j\gamma_lv_i(\alpha_i\bar v_i). 
\end{align*} 
This proves $\alpha_i\bar v_i = \gamma_iv_jv_l$, hence the second equation of \eqref{RKONE}, if $v_i \neq 0$. But of $v_i = 0$, then $\alpha_j\alpha_l = \gamma_j\gamma_ln_C(v_i) = 0$, forcing $\alpha_j = 0$ or $\alpha_l = 0$ since $k$ is an integral domain, which in turn gives $0 = \alpha_i\alpha_j = \gamma_i\gamma_jn_C(v_l)$ in the former case and $0 = \alpha_l\alpha_i = \gamma_l\gamma_in_C(v_j)$ in the latter. Since $C$ has no zero divisors, its norm represents zero only trivially (Exc.~\ref{pr.ZERDIV}). Therefore $v_j = 0$ or $v_l = 0$, and the proof of \eqref{RKONE} is complete.
\end{sol}

\begin{sol}{pr.FLASET} \label{sol.FLASET} 
The first (resp. second) part of (a) follows by applying \ref{ss.FLAFALG} to $f\vert_N:N \rightarrow
M^\prime$ (resp. to $\pi \circ f:M \rightarrow M^\prime/N^\prime$,
$\pi:M^\prime \rightarrow M^\prime/N^\prime$ being the canonical
projection). In (b) we apply (a) with $N^\prime = P$ to the
natural embedding $i:N \rightarrow M$. For (c), we consider the
canonical map $\bigoplus_{\alpha \in I} N_\alpha \rightarrow M$
determined by the inclusions $N_\alpha \rightarrow M$ and apply
\ref{ss.FLAFALG}. Finally, in (d), we let $M^0$ be a free $k$-module with basis $(e_\alpha)_{\alpha \in I}$ and apply \ref{ss.FLAFALG} to the $k$-linear map $M^0 \rightarrow M$, $e_\alpha \mapsto x_\alpha$, $\alpha \in I$.
\end{sol}

\begin{sol}{pr.IDFREU} \label{sol.IDFREU} 
(a) Let $R \in \kalg$ be a flat $k$-algebra. As usual, for a submodule $M \subseteq J$, we identify $M_R \subseteq J_R$ canonically. Any $x \in \Sq(J)$ may be written as $x = \sum\alpha_jx_j^2$, $\alpha_j \in k$, $x_j \in J$, and we conclude $x_R = \sum\alpha_{jR}(x_{jR})^2 \in \Sq(J_R)$. This show $\Sq(J)_R \subseteq \Sq(J_R)$. Conversely, let $x \in \Sq(J_R)$. Then $x = \sum r_jx_j^2$, $r_j \in R$, $x_j \in J_R$, and writing $x_j = \sum_m x_{mj} \otimes r_{mj}$, $x_{mj} \in J$, $r_{mj} \in R$, we obtain
\[
x = \sum_{j,m} x_{mj}^2 \otimes r_jr_{mj}^2 + \sum_{j,m<n}(x_{mj} \circ x_{nj}) \otimes (r_jr_{mj}r_{nj}) \in \Sq(J)_R.
\]
Summing up, we have shown $\Sq(J)_R = \Sq(J_R)$.

Next assume that $J$ is split. The assertion $\Sq(J) = J$ is trivial for $n = 1$ and obvious for $n = 3$ since $\Sq(J)$ contains all idempotents. Let us therefore assume $n \geq 6$. By definition (\ref{ss.SPLIFRE}), we may assume $J = \Her_3(C)$, where $C$ is a split composition algebra of rank $\frac{n-3}{3}$ over $k$. The diagonal of $\Her_3(C)$ clearly belongs to $\Sq(J)$. Furthermore, for $i = 1,2,3$ and $u_i \in C$, we note $u_i[jl] \in J_1(e_{jj})$ by Prop.~\ref{p.PEDEDI} and (\ref{t.PEDECOM}.\ref{JIJ}). Hence (\ref{t.PEDESI}.\ref{JONE}), (\ref{ss.BACUNO.fig}.\ref{ADJJA}), (\ref{ss.THERCU}.\ref{MALIT}), (\ref{ss.THERCU}.\ref{MAQIT}) imply
\begin{align*}
(e_{jj} + u_i[jl])^2 =\,\,&e_{jj}^2 + e_{jj} \circ u_i[jl] + u_i[jl]^2 \\
=\,\,&e_{jj} + u_i[jl] + u_i[jl]^\sharp + T(u_i[jl])u_i[jl] - S(u_i[jl])1_J \\
=\,\,&e_{jj} + u_i[jl] - n_C(u_i)e_{ii} + n_C(u_i)1_J,
\end{align*}
which in turn shows $C[jl] \subseteq \Sq(J)$. Thus $\Sq(J) = J$ if $J$ is split.

Finally, let $J$ be arbitrary. Cor.~\ref{c.CHARF} yields a faithfully flat $k$-algebra $R \in \kalg$ making $J_R$ a split Freudenthal algebra of rank $n$ over $R$. By what we have just seen, $\Sq(J)_R = \Sq(J_R) = J_R$. Hence the inclusion $\Sq(J) \to J$ becomes an isomorphism after changing scalars from $k$ to $R$ and thus, by faithful flatness, must have been one all along. In other words, $\Sq(J) = J$.

\smallskip

(b) It suffices to prove the following implications. 

\step{1}
 \emph{If $\mfa \subseteq k$ is an ideal, then so is $\mfa J \subseteq J$.} Obvious. 
 
 \step{2}
 \emph{If $I \subseteq J$ is an outer ideal, then $I \cap k \subseteq k$ is an ideal.} This is clear since outer ideals of $J$ are, in particular, $k$-submodules. 
 
 \step{3}
 \emph{If $\mfa \subseteq k$ is an ideal, then $\mfa = (\mfa J) \cap k$.} For $\alpha \in \mfa$ we have $\alpha = \alpha 1_J \subseteq (\mfa J) \cap k$, so the left-hand side is contained in the right. Conversely, let $\alpha \in (\mfa J) \cap k$. Then $\alpha \in k$ satisfies $\alpha 1_J = \sum \alpha_ix_i$, $\alpha_i \in \mfa$, $x_i \in J$. Since $1_J \in J$ is unimodular, some linear form $\lambda\:J \to k$ has $\lambda(1_J) = 1$, hence $\alpha = \sum \alpha_i\lambda(x_i) \in \mfa$. 
 
 \step{4}
  \emph{If $I \subseteq J$ is an outer ideal, then $(I \cap k)J = I$.} Let $\alpha \in I \cap k$ and $x \in J$. Then $\alpha \in k$ and $\alpha 1_J \in I$. Since $I \subseteq J$ is an outer ideal, we conclude $\alpha x^2 = U_x(\alpha 1_J) \in I$, and (a) implies $\alpha J = \alpha\Sq(J) \subseteq I$. Hence $(I \cap k)J \subseteq I$.

In order to prove equality, we first assume that $J = \Her_3(C)$ is split, $C$ a split composition algebra of rank $\frac{n-3}{3}$ over $k$. Since regularity of $J$ is inherited by $C$, Exc.~\ref{pr.IDCUJM} yields an ideal $\mfb \subseteq k$ such that
\[
I = \mfb\Her_3(C) = H_3(\mfb,\mfb C,\Eins_3) = \sum(\mfb e_{ii} + (\mfb C)[jl]),
\]
hence $I \cap k = \mfb$ and $(I \cap k)J = \mfb J = I$.

Finally, let $J$ be arbitrary and apply Cor.~\ref{c.CHARF} to find a faithfully flat $k$-algebra $R \in \kalg$ making $J_R$ split of rank $n$ over $R$. Let $y \in I$ and write $x \in J_R$ as $x = \sum x_i \otimes r_i$, $x_i \in J$, $r_i \in R$. Then
\[
U_xy = \sum (U_{x_i}y) \otimes r_i^2 + \sum_{i<j} \{x_iyx_j\} \otimes (r_ir_j) \in I_R
\]
since $I \subseteq J$ is an outer ideal. Thus $U_xI_R = R(U_xI) \subseteq RI_R = I_R$, and we have shown that $I_R \subseteq J_R$ is an outer ideal. Now we combine Exc.~\ref{pr.FLASET}~(c),(d) to compute
\[
\big((I \cap k)J\big)_R = (I \cap k)_RJ_R = (I_R \cap k_R)J_R = (I_R \cap R)J_R = I_R.
\]
But this means that the inclusion $(I \cap k)J \to I$ becomes an isomorphism when extending scalars from $k$ to $R$, which by arguing as before completes the proof of $4^\circ$.
\end{sol}

\begin{sol}{pr.ELMFRA} \label{sol.ELMFRA}
(a) The second statement follows immediately from the first combined with \ref{ss.THREEFF}~(ii). Since elementary frames are invariant under base change, the defining equation for $\bfElf(J)$ does indeed define a subfunctor of $J_{\bfa}^3$. Let $u_1^\ast, \ldots, u_m^\ast$ be a finite set of generators of $J^\ast$, the dual of the $k$-module $J$. Then
\begin{align*}
\bfElf(J)(R) =\,\,&\{(e_1,e_2,e_3) \in J_R \mid T_J(e_j) = 1, \\
\,\,&\la u_{iR}^\ast,e_j^\sharp \ra = \la u_{jR},e_1 \times e_2 - e_3\ra = 0\;\; (0 \leq i \leq m,\,j = 1,2)\}
\end{align*}
for all $R \in \kalg$, and we conclude from \ref{ss.CLOSU}, \ref{ss.CONVEN},  Exercises~\ref{pr.EMMA} and \ref{pr.IDFIPR}~(b) that $\bfX:= \bfElf(J)$ is a finitely presented closed subscheme of $J_{\bfa}^3$; in particular, it is affine. It remains to show that $\bfX$ is fppf and smooth. 

Beginning with smoothness, let $R \in \kalg$ and $\mfa \subseteq R$ be an ideal such that $\mfa^2 = \{0\}$. By \ref{ss.SMOAFF}, we must prove that every elementary frame of $J_{R/\mfa} = J_R/\mfa J_R$ can be lifted to an elementary frame of $J_R$. But this follows from Exercises~\ref{pr.LIFCOM}, \ref{pr.RENIRA}~(b) and \ref{pr.NORID}~(e). 

In order to prove that $\bfX$ is fppf, it suffices to show by Prop.~\ref{p.SMOFPPF} that $\bfX$ has non-empty geometric fibers, so let $K \in \kalg$ be an algebraically closed field and put $I := \Nil(J_K)$. Then $\bar J := J_K/\Nil(J_K)$ is semi-simple cubic Jordan $K$-algebra of dimension at least $3$, and by Exercises~\ref{pr.LIFCOM} and \ref{pr.RENIRA}~(b) again, we will be through once we have shown that $\bar J$ contains two orthogonal elementary idempotents. Thanks to Thm.~\ref{t.RACUSE}, we may assume $J = F^{(+)} \times J(M,q,e)$, where $(M,q,e)$ is a non-degenerate pointed quadratic module of dimension at leas $2$ over $K$. We must show that $J(M,q,e)$ contains an elementary idempotent, i.e., by \ref{ss.ELICL}, an element $c$ having $q(c) = 0$ and $q(c,e) = 1$. If not, then $q(x) = 0$ implies $q(x,e) = 0$, for all $x \in M$, and from the Hilbert Nullstellensatz we deduce $q(x) = q(x,e)^2$ for all $x \in M$. But this contradicts non-degeneracy of $q$, so $\bfX(K) \ne \emptyset$, and the solution of (a) is complete.

(b) Since $J$ is regular, it will be enough to show that there is a \emph{unique} bilinear multiplication $J \times J \to J$, $(x,y) \mapsto xy$, making $J$ a cubic commutative associative $k$-algebra $E$ such that $E^{(+)} = J$ as cubic Jordan algebras. Since $J$ is regular, it is separable (Exc.~\ref{pr.CUBNIL}), so 
by (a) and faithfully flat descent, we may assume that $J$ contains an elementary frame $\Omega = (e_1,e_2,e_3)$. Since the corresponding Peirce components of $J$ are orthogonal relative to the bilinear trace (Prop.~\ref{p.COMPOR}~(c)), $J_0 := \sum ke_i \subseteq J$ is a regular cubic Jordan subalgebra, which implies $J = J_0 \oplus J_0^\perp$, and comparing ranks, we conclude that $J$ is free of rank $3$ as a $k$-module, with basis $(e_1,e_2,e_3)$. Introducing the multiplication rules $e_ie_j = \delta_{ij}e_i$ on the basis vectors and extending it bilinearly to all of $J$, we obtain the structure of a split cubic \'etale $k$-algebra $E$ on the $k$-module $J$. From (\ref{p.DECNOAD}.\ref{EXPNOR}) and (\ref{e.CUBET}.\ref{NXXI}) we deduce $N_J(\sum \xi_ie_I) = \xi_1\xi_2\xi_3 = N_E(\sum \xi_ie_i)$ in all scalar extensions, and combining with Prop.~\ref{p.CUALJO}, we conclude $J = E^{(+)}$ as cubic Jordan algebras. It remains to prove uniqueness, so let $E'$ be any cubic commutative associative $k$-algebra having $E'^{(+)} = J$ as cubic Jordan algebras. For $i,j = 1,2,3$ distinct, $e_ie_j = e_i^2e_j = U_{e_i}e_j = 0$, so $E' = E$ as cubic commutative associative $k$-algebras.
\end{sol}

\begin{sol}{pr.UOPID} \label{sol.UOPID}
(a) Expanding the $U$-operator, we obtain
\begin{align*}
U_xe =\,\,&U_{\xi e + x_1 + x_0}e = \xi^2U_ee + \xi\{eex_1\} + \xi\{eex_0\} + U_{x_1}e + \{x_1ex_0\} + U_{x_0}e,
\end{align*}
where the Peirce rules of Thm.~\ref{t.PEDESI} imply $\{x_1ex_0\} = U_{x_0}e = 0$, $U_{x_1}e \in J_0(e)$, $\{eex_i\} = e \circ x_i - \{e(1_J - e)x_i\} = ix_i$ for $i = 0,1$, $U_ee = e$, and (a) is proved. 

(b) Assume first that $J$ is split, so there is an identification $J = \Her_3(C)$, for some composition $k$-algebra $C$. We write
\begin{align*}
x = \sum(\xi_ie_{ii} + u_i[jl]) &&(\xi_i \in k,\,u_i \in C,\, i = 1,2,3),
\end{align*}
let $(mnp)$ be a cyclic permutation of $(123)$ and apply (a) to $e := e_{mm}$. Since $U_xe = e$, we conclude $\xi_m^2 = 1$, hence $\xi_m \in k^\times$ and  $x_1 = 0$.  Thus
\[
x = \sum \xi_ie_{ii} + x_0, \quad x_0 \in J_0(e_{mm}) = C[mn] \oplus C[pm].
\]
Since this is true for all $m = 1,2,3$, we actually have $x = \sum \xi_ie_{ii}$. Now the hypothesis and (\ref{ss.EICJM}.\ref{UJLJJ}), (\ref{ss.EICJM}.\ref{UJLLL}) imply $1_C[np] = U_x1_C[np] = \xi_n\xi_p1_C[np]$, hence $\xi_n\xi_p = 1$ and then $\xi_p = \xi_n\xi_p^2 = \xi_n$: the assertion is proved.

If $J$ is arbitrary, we apply Cor.~\ref{c.CHARF} to find a faithfully flat algebra $R \in \kalg$ making $J_R$ a split Freudenthal algebra over $R$. The special case just treated implies $x_R \in R1_{JR}$. On the other hand, since $1_J \in J$ is  unimodular, we have $J = k1_J \oplus M$, for some submodule $M \subseteq J$. In particular, $x = \xi 1_J + m$, for some $m \in M$. Changing scalars from $k$ to $R$, we conclude $m_R = 0$, hence $m = 0$ since $R$ is faithfully flat. Thus $x = \xi 1_J$, and $U_x = \Eins_J$ implies $\xi^2 = 1$.
\end{sol}

\solnsec{Section~\ref{s.ISNOAI}}

\begin{sol}{pr.SPAINV} \label{sol.SPAINV}  (a) Write $W$ for the subspace of $J^0$ spanned by the invertible trace-zero elements of $J$. For $v_i,v_l \in C^\times$, the quantities
\begin{align*}
x_1 :=\,\,&e_{ii} - e_{jj} + v_i[jl], \\
x_2 :=\,\,&e_{ll} - e_{ii} + v_l[ij], \\
x_3 :=\,\,&e_{jj} - e_{ll} - v_l[ij] 
\end{align*}
are invertible of trace zero and thus all belong to $W$. Hence so does $x_1 + x_2 + x_3 = v_i[jl]$. This not only implies $e_{ii} - e_{jj} = x_1 - v_i[jl] \in W$ but also $C[ij] \subseteq W$ by the hypothesis on $C$. Summing up we conclude $J^0 = W$. 

\smallskip
(b) Assume that $J$ is simple. By Thm.~\ref{t.RACUSE}, $J$ is either a division algebra, in which case the assertion is obvious, or it has the form $J \cong \Her_3(C,\Gamma)$, where $\Gamma \in \GL_3(F)$ is a diagonal matrix and $C$ is a pre-composition algebra over $F$. We check when the hypotheses of (a) are fulfilled. 

\begin{lem*}
A pre-composition algebra $C$ over $F$ either is spanned as a vector space over $F$ by $C^\times$, or $F = \IF_2$ and $C \cong F \times F$ is split quadratic \'etale.
\end{lem*}

\begin{proof}
If $F = \IF_2$ and $C \cong F \times F$, then $C^\times = \{1_C\}$ cannot possibly span $C$ as an $F$-vector space. Conversely, assume $F \ncong \IF_2$ or $C$ is not split quadratic \'etale. We must show $C = FC^\times$ and may assume that $C$ is not a division algebra, so we are reduced to the case that $C$ is a split composition algebra of dimension $> 1$. If $C = F \times F$ is split quadratic \'etale, then $F$ contains an element $\alpha \ne 0,1$, and since $\left(\begin{smallmatrix}
1 \\
0
\end{smallmatrix}\right) = (1 - \alpha)^{-1}(\left(\begin{smallmatrix}
1 \\
\alpha
\end{smallmatrix}\right) - \left(\begin{smallmatrix}
\alpha \\
\alpha
\end{smallmatrix}\right))$, $\left(\begin{smallmatrix}
0 \\
1
\end{smallmatrix}\right) = (1 - \alpha)^{-1}(\left(\begin{smallmatrix}
\alpha \\
1
\end{smallmatrix}\right) - \left(\begin{smallmatrix}
\alpha \\
\alpha
\end{smallmatrix}\right))$, the assertion follows. If $C = \Mat_2(F)$ is the split quaternions, the automorphism group of $C$ acts transitively on its elementary idempotents, but also on its non-zero nilpotent elements, so the equations
\[
\left(\begin{matrix}
1 & 0 \\
0 & 0
\end{matrix}\right) = \left(\begin{matrix}
1 & 1 \\
1 & 0
\end{matrix}\right) - \left(\begin{matrix}
0 & 1 \\
1 & 0
\end{matrix}\right), \;\left(\begin{matrix}
0 & 1 \\
0 & 0
\end{matrix}\right) = \left(\begin{matrix}
1 & 1 \\
0 & 1
\end{matrix}\right) - \left(\begin{matrix}
1 & 0 \\
0 & 1
\end{matrix}\right)
\]
give what we want. Finally, the case that $C = \Zor(F)$ is the split octonions is reduced to the quaternionic case since every element of $C$ is contained in a quaternion subalgebra (Exc.~\ref{pr.COVQUATALG}).
\end{proof}

In view of the lemma, the first part of (b) follows from (a) unless $F = \IF_2$ and $J \cong \Her_3(F \times F,\Gamma) \cong \Her_3(F \times F)$ (by Prop.~\ref{p.FRSPCOM}) $\cong A^{(+)}$, $A := \Mat_3(F)$ (by Prop.~\ref{p.HERKEI}). The subspace $W$ of $J$ spanned by $J^0 \cap J^\times$ is stable under the action of $\GL_3(F)$ on $A$ by conjugation. Since
\[
e_{12} = \left(\begin{matrix}
0 & 1 & 0 \\
0 & 0 & 0 \\
0 & 0 & 0
\end{matrix}\right) = \left(\begin{matrix}
1 & 1 & 0 \\
0 & 1 & 1 \\
0 & 1 & 0
\end{matrix}\right) + \left(\begin{matrix}
1 & 0 & 0 \\
0 & 1 & 1 \\
0 & 1 & 0
\end{matrix}\right) \in W,
\]
we therefore conclude $e_{ij} \in W$ ($1\le i,j\le 3,\,i \ne j$). On the other hand,
\[
e_{11} + e_{22} = \left(\begin{matrix}
1 & 0 & 0 \\
0 & 1 & 0 \\
0 & 0 & 0
\end{matrix}\right) = \left(\begin{matrix}
1 & 1 & 0 \\
1 & 1 & 0 \\
0 & 0 & 0
\end{matrix}\right) + e_{12} + e_{21},
\]
where the first summand on the very right is nilpotent of index $2$, hence conjugate to $e_{12}$ and thus belongs to $W$. This implies $e_{11} + e_{22} \in W$ and, similarly, $e_{11} + e_{33} \in W$. Summing up, therefore, we have $W = J^0$.

In order to prove the second part of (b), it suffices to note that $\diag(1,-1,1) \in J$ is invertible of trace $1$, hence together with $J^0$ spans all of $J$ as a vector space over $F$.
\end{sol}

\begin{sol}{pr.RANONE} \label{sol.RANONE}  (a) (\ref{ss.ISCUNO}.\ref{ISOAD}), the adjoint identity (\ref{ss.BACUNO.fig}.\ref{ADJI}) and (\ref{ss.ISCUNO}.\ref{ISONO}).

\smallskip

(b) If $u$ is an elementary idempotent in some isotope of $J$, then $u$ has rank $1$ by (a). Conversely, assume $u$ has rank $1$. By hypothesis, $\Nil(J) = \{0\}$ and hence $J$ has no absolute zero divisors (Exc.~\ref{pr.ABZECU}~(e)), while part (b) of that exercise implies $T(q,u) \neq 0$ for some $q \in J$. By Exc.~\ref{pr.SPAINV}~(b) we may assume $q \in J^\times$, whence (\ref{ss.ISCUNO}.\ref{ISOLINT}) implies $T^{(p)}(u) = 1$ for $p := T(q,u)^{-1}q \in J^\times$. Thus $u$ is an elementary idempotent in $J^{(p)}$.

Next assume $u$ is a co-elementary idempotent in some isotope $\pJ$ of $J$. By (a), we may assume $\pJ= J$. Then $u = 1_J - e$, for some elementary idempotent $e \in J$. Then $u^\sharp = e$ (by Prop.~\ref{p.PELE}~(a)) has rank $1$, so $u$ has rank $2$. Conversely, if $u$ has rank $2$, arguing as before with $u^\sharp$ in place of $u$, we find an invertible element $p \in J$ having $T(u^\sharp,p^{-1}) \neq 0$, hence $T(p^\sharp,u^\sharp) = N(p)T(u^\sharp,p^{-1}) \neq 0$. Passing to the isotope $J^{(p)}$ if necessary, we may therefore assume $T(u^\sharp) = S(u) \neq 0$. But this means that the element $e := T(u^\sharp)^{-1}u^\sharp \in J$ is an element of rank $1$ and trace $1$, hence an elementary idempotent. Now (\ref{ss.BACUNO.fig}.\ref{CADJ}) yields
\begin{align*}
e \times u =\,\,&T(u^\sharp)^{-1}u^\sharp \times u =T(u^\sharp)^{-1}\Big(\big(T(u)T(u^\sharp) - N(u)\big)1 - T(u^\sharp)u - T(u)u^\sharp\Big) \\
=\,\,&T(u)1 - u - T(u)e = T(u)(1 - e) - u.
\end{align*} 
This implies $u \in J_0(e)$ by Prop.~\ref{p.PELE}~(b), while Cor.~\ref{c.FAULE} shows that $u$ is invertible in $J_0(e)$. Write $v \in J_0(e)$ for the corresponding inverse and pass to the isotope $\pJ := J^{(p)}$, $p := e + v  \in J^\times$. Since $e$ continues to have rank $1$ in $\pJ$ and satisfies $T^{(p)}(e) = T(e,p) = T(e,e) = T(e) = 1$ by Prop.~\ref{p.PELE}~(b), it is, in fact, an elementary idempotent in $\pJ$. Moreover, $1_{\pJ} = p^{-1} = e + u$, forcing $u = 1_{\pJ} - e$, and (ii) holds. 

\smallskip

(c) $e$ can be extended to an elementary frame of $J$ if and only if $J_0(e) = J(M_0,q_0,e_0)$ as in Cor.~\ref{c.FAULE} contains an elementary idempotent (Prop.~\ref{p.COMPOR}). Since semi-simplicity of $J$ is inherited by $J_0(e)$ (Exc.~\ref{pr.NIRAPE}), the quadratic form $q_0$ is non-degenerate (Exc.~\ref{pr.NIPOQUA}). Moreover, by (\ref{p.ADPEMO}.\ref{QUQUX}) and Exc.~\ref{pr.NILEMZER}~(b), it is also isotropic. By Prop.~\ref{p.JAMAPO}, therefore, $J_0(e)$ contains an elementary idempotent if and only if the linear trace of $J_0(e)$, which by Cor.~\ref{c.FAULE} agrees with $T$ on $J_0(e)$, is different from zero. 

\smallskip

(d) Since 
\[
e = \sum(e_{ii} + 1_K[jl])
\]
and $2 = 0$ in $k$, we have $T(e) = 3 = 1$. On the other hand, since the the conjugation of $K$ is the identity, (\ref{ss.THERCU}.\ref{MADJ}) shows $e^\sharp = 0$. Thus $e$ is an elementary idempotent in $J$ and we have $e_0 := 1_J - e = \sum 1_K[jl]$. Now let $x = \sum(\xi_ie_{ii} + u_i[jl]) \in J_0(e)$. Then (\ref{p.PELE}.\ref{PEZERO}) implies $x = T(x)e_0 - e \times x = T(x)\sum 1_K[jl] - e \times x$. But (\ref{ss.THERCU}.\ref{MABADJ}) gives
\begin{align*}
e \times x =\,\,&\sum\big((\xi_j + \xi_l)e_{ii} + (-u_i - \xi_i1_K + u_l + u_j)[jl]\big) \\
=\,\,&\sum\big((\xi_j + \xi_l)e_{ii} + (\xi_i1_K + u_1 + u_2 + u_3)[jl]\big),
\end{align*}
and a comparison shows $\xi_1 = \xi_2 + \xi_3$, hence $T(x) = \sum\xi_i = 0$. By (c), therefore, $e$ cannot be extended to an elementary frame in $J$.

As to the second part of (d), we argue indirectly and assume $u = u_1 + u_2$ for some rank-one elements $u_1,u_2 \in J$. In particular, $u_j$ for $j = 1,2$ is a symmetric 3-by-3 matrix of rank $1$ with entries in $K$, and hence there are $a_j \in K^\times$, $0 \neq x_j \in K^3$ such that $u_j = a_jx_jx_j^\trans$. Since 
\[
u = a_1x_1x_1^\trans  + a_2x_2x_2^\trans 
\]
has rank $2$, the quantities $x_1,x_2$ are linearly independent over $K$, and $a_1,a_2$ are both different from zero. Moreover, they are contained in $\Img(u) = Ke_1 + Ke_2$, where $(e_1,e_2,e_3)$ stands for the canonical basis of $K^3$ over $K$. Thus $x_j = b_{1j}e_1 + b_{2j}e_2$ for $j = 1,2$ and some matrix $(b_{ij}) \in \GL_2(K)$. Writing $e_{ij}$ for the ordinary matrix units in $\Mat_3(K)$, we now obtain
\begin{align*}
e_{12} + e_{21} =\,\,&u = a_1(b_{11}e_1 + b_{21}e_2)(b_{11}e_1^\trans  + b_{21}e_2^\trans) + a_2(b_{12}e_1 + b_{22}e_2)(b_{12}e_1^\trans  + b_{22}e_2^\trans) \\
=\,\,&a_1(b_{11}^2e_{11} + b_{11}b_{21}e_{12} + b_{21}b_{11}e_{21} + b_{21}^2e_{22}) \\
\,\,&+a_2(b_{12}^2e_{11} + b_{12}b_{22}e_{12} + b_{22}b_{12}e_{21} + b_{22}^2e_{22}).
\end{align*}
Since the characteristic is $2$, this implies $a_1b_{11}^2 = a_2b_{12}^2$, $a_1b_{21}^2 = a_2b_{22}^2$, hence
\begin{align*}
a_1a_2\det(b_{ij})^2 =\,\,&a_1b_{11}^2a_2b_{22}^2 - a_2b_{12}^2a_1b_{21}^2 = 0,
\end{align*}
a contradiction. This completes the solution of (d).

\smallskip

(e) Since two vectors in $F^3$ are linearly independent if and only if they are distinct, one checks that that $J := \Her_3(F)$ contains precisely seven elements of rank $1$, namely,
\[
e_{11}, e_{22}, e_{33}, \left(\begin{matrix}
1 & 1 & 0 \\
1 & 1 & 0 \\
0 & 0 & 0 
\end{matrix}\right), \left(\begin{matrix}
0 & 0 & 0 \\
0 & 1 & 1 \\
0 & 1 & 1 
\end{matrix}\right), \left(\begin{matrix}
1 & 0 & 1 \\
0 & 0 & 0 \\
1 & 0 & 1 
\end{matrix}\right), \left(\begin{matrix}
1 & 1 & 1 \\
1 & 1 & 1 \\
1 & 1 & 1 
\end{matrix}\right), 
\]
and precisely four elementary idempotents, namely,
\[
e_{11}, e_{22}, e_{33}, e :=  \left(\begin{matrix}
1 & 1 & 1 \\
1 & 1 & 1 \\
1 & 1 & 1 
\end{matrix}\right). 
\]
It follows (again) that $e$, which agrees with the idempotent of (d), cannot be completed to an elementary frame in $J$.

\smallskip

(f) (i) $\Rightarrow$ (ii). Let $e_0$ be a co-elementary idempotent in some isotope $\pJ$ of $J$. Replacing $J$ by $\pJ$ if necessary, we may assume $\pJ= J$. Then $e := 1_J - e_0$ is an elementary idempotent in $J$, which by hypothesis can be extended to an elementary frame $(e = e_1,e_2,e_3)$ of $J$. Hence $e_0 = e_2 + e_3$, as claimed.

(ii) $\Rightarrow$ (iii). (b) and (ii).

(iii) $\Rightarrow$ (iv). If $J$ were not regular, then, by Thm.~\ref{t.RACUSE} and up to isotopy, it would be as in (d), hence would contain an element of rank $2$ that violates the condition described in (iii).

(iv) $\Rightarrow$ (i). By hypothesis, the bilinear trace $T$ of $J$ is regular, and the Peirce components of $J$ relative to an elementary idempotent $e \in J$ are orthogonal with respect to $T$. Hence $T$ stays regular on $J_0(e)$ and since $T(x_0) = T(1_J,x_0) = T(1_J - e,x_0)$ for all $x_0 \in J_0(e)$, the linear trace of $J$ cannot be zero on $J_0(e)$. The assertion now follows from (c).
\end{sol}

\begin{sol}{pr.INVORB}
Given $x \in J^\times$ such that $N_J(x) = 1$, we wish to find $g \in \Inv(J)$ such that $gx = 1_J$.  By Lemma \ref{l.KRUT}, we may assume that $x$ is diagonal, i.e., $x = \sum x_i e_{ii}$ for some $x_i \in \kx$ such that $\prod x_i = 1$.  Pick $y_i \in k^\times$ such that $y_i^2 = x_i$ and define $g := U_y$ for  $y = \sum y_i e_{ii}$, so $gx = 1_J$ as desired.  Furthermore, for $z \in J$, we have 
\[
N(gz) = N(U_yz) = N(y)^2N(z) = N(y^2)N(z) = N(x)N(z) = N(z), 
\]
so $U_y \in \Inv(J)$.
\end{sol}

\begin{sol}{pr.COFREN} \label{sol.COFREN}  
(a): Let $J$ be a Freudenthal algebra of rank $n > 6$ over $F$.  The norm of $J$ is a cubic form in more than 3 variables and is therefore isotropic by the Chevalley-Warning Theorem.  Then by Prop.~\ref{p.REDDEF}, $J \cong \Her_3(C, \Gamma)$ for a composition algebra $C$ over $F$ of rank $> 1$.  If the rank is $> 2$, then, again since $F$ is finite, $C$ is split, and therefore $J$ is itself split (Prop.~\ref{p.FRSPCOM}). If the rank is $2$, then $C$ is either split quadratic \'etale or isomorphic to $K$. In the former case, $J$ is split while in the latter case, the norm of $C  \cong K$ is surjective (since, again by Chevalley-Warning, any quadratic form in more than three variables over $F$ is isotropic), and (\ref{pr.ISCOPA}.\ref{NOCHA}) combined with Exc.~\ref{pr.ISTIS} implies $J \cong \Her_3(K)$.

\smallskip

(b): Let $k$ be a finite commutative ring, let $C$ be a composition algebra of rank $r > 2$ and $J$ a Freudenthal algebra of rank $n > 9$ over $k$. We first assume that $k$ is reduced. Since $k$, being finite, is artinian, we conclude from \cite[p.~203]{MR780184} that $k = k_1 \times \cdots \times k_m$ is a finite direct product of finite fields. This implies $C = C_{k_1} \times \cdots \times C_{k_m}$, where each $C_{k_i}$, for $1 \leq i \leq m$, is composition algebra of rank $r$ over the finite field $k_i$, hence split, by \ref{ss.FINFIE}. Thus $C$ is split over $k$. Similarly, $J = J_{k_1} \times \cdots \times J_{k_m}$, where each $J_i := J_{k_i}$ is a simple Freudenthal algebra of rank $n$  over the finite field $k_i$, and therefore is split by part (a). 

If $k$ is arbitrary, we only consider the case of a Freudenthal algebra, the composition algebra case being completely analogous. We put $\mfa := \Nil(k)$, $\bar k := k/\mfa$ and $\bar J := J/\mfa J = J_{\bar k}$. Since $\bar k$ is (finite and) reduced, the special case just treated implies that $\bar J$ is split over $\bar k$. Since $k$ is artinian, the ideal $\mfa \subseteq k$ is nilpotent, forcing $\mfa^m = \{0\}$ for some positive integer $m$. We show by induction on $i = 1,\dots,m$ that $J_{k/\mfa^i}$ is split over $k/\mfa^i$. By the special case treated at the beginning, this is true for $i = 1$. Now let $i > 1$ and suppose the assertion holds for $i - 1$. We have $k/\mfa^{i-1} = (k/\mfa^i)/(\mfa^{i-1}/\mfa^i)$ and ($\mfa^{i-1}/\mfa^i)^2 = \{0\}$. By Prop.~\ref{p.AFSPLIF}, the affine scheme $\bfSplid(J)$ is smooth, so the natural map
\begin{align}
\label{FINSURJ} \bfSplid(J)(k/\mfa^i) \longrightarrow \bfSplid(J)(k/\mfa^{i-1})
\end{align}
is surjective. By the induction hypothesis, $J_{k/\mfa^{i-1}}$ is split and hence there exists a splitting datum for $J_{k/\mfa^{i-1}}$ (Prop.~\ref{p.FASPLI}). We therefore conclude from the surjectivity of the map in \eqref{FINSURJ} that there exists a splitting datum for $J_{k/\mfa^i}$, which is therefore split over $k/\mfa^i$. This completes the induction, and for $i = m$ we conclude that $J = J_{k/\mfa^m}$ is split over $k/\mfa^m = k$.
\end{sol}

\begin{sol}{pr.str.outer} \label{sol.str.outer}
\ref{pr.str.outer.6}: We use the group homomorphism $\mu \!: \Str(J) \to \kx$ from Lemma \ref{l.ISOTSIM} and note that 
\[
\mu(\phi \eta \phi^{-1}) = \mu(\eta) \quad \text{and} \quad \mu(\eta^{\sharp - 1}) = \mu(\eta)^{-1}.
\]
Therefore, it suffices to find an $\eta$ with
$\mu(\eta)^2 \ne 1$.

Taking now $\eta$ to be multiplication by $c$, which is a similarity of the norm form, so belongs to $\Str(J)$.  Since $\mu(\eta)^2 = c^6 \ne 1$, this provides the desired example.

\ref{pr.str.outer.fix}: $\eta = \eta^{\sharp - 1}$ if and only if $\eta^{-1} = \eta^\sharp = \eta^{-1} U_{\eta(1_J)}$, i.e., $U_{\eta(1_J)} = \Eins_J$.  Exercise \ref{pr.UOPID} shows that this is equivalent to $\eta(1_J) = \zeta 1_J$ for some $\zeta \in \bfmu_2(k)$.

Given $\eta$ fixed by the automorphism, then, $\eta = \zeta \phi$ for some $\zeta \in \bfmu_2(k)$ and $\phi \in \Str(J)$ such that $\phi(1_J) = 1_J$, i.e., $\phi \in \Aut(J)$ by Thm.~\ref{t.STRINV.JOAL}~(c).  Conversely, given $\zeta \in \bfmu_2(k)$ and $\phi \in \Aut(J)$, $\eta := \zeta \phi$ is in $\Str(J)$ and satisfies $\eta(1_J) = \zeta 1_J$, so it is fixed by the automorphism, proving (1).  

For part (2), we aim to show that $\Inv(J) \cap \bfmu_2(k) \Aut(J) = \Aut(J)$.  For $\zeta \in \bfmu_2(k)$, $\phi \in \Aut(J)$, we have $N(\zeta\phi(1_J)) = \zeta^3 = \zeta$.  Consequently. $\zeta \phi$ belongs to $\Inv(J)$ if and only if $\zeta = 1_k$.
\end{sol}

\solnsec{Section~\ref{s.RFRALGF}}

\begin{sol}{pr.FREFISIX} \label{sol.FREFISIX}
As in the solution of Exc.~\ref{pr.COFREN}, every Freudenthal algebra over a finite field $F$ is of the form $J = \Her_3(F,\Gamma)$ for some $\Gamma = \diag(\gamma_1,\gamma_2,1) \in \GL_3(F)$. If $F$ has characteristic 2, then everything is a square, and $J$ is split (Exc. \ref{pr.DIAFIX}~(b)).  If $F$ has characteristic different from 2, then $J$ is regular, and $Q_J = \qform{\gamma_2,\gamma_1,\gamma_2 \gamma_1}$ is a classifying invariant (Thm.~\ref{t.SPR}). Note that over $F$ we may identify quadratic forms and symmetric bilinear forms, allowing us to use notation and terminology of O'Meara \cite{MR0152507}. Since the determinant of $Q_J$ is a square, \cite[62:1a, p.~157]{MR0152507} implies $Q_J \cong \bform{1,1,1}$, and $J$ is split.
\end{sol}

\begin{sol}{pr.CLAFRA} \label{sol.CLAFRA} Let $\Omega = (e_1,e_2,e_3)$, $\Omega^\prime = (\pe_1,\pe_2,\pe_3)$ be two elementary frames in $J$. Writing $C$ for the co-ordinate algebra of $J$, let
\[
\Phi\:J \overset{\sim} \longrightarrow \Her_3(C,\Gamma), \quad \Gamma = \diag(\gamma_1,\gamma_2,\gamma_3) \in \GL_3(F)
\]
be any $\Omega$-co-ordinatization of $J$. If $\eta$ is an automorphism of $J$ sending $\pe_i$ to $e_i$ for $1 \le i \le 3$, then $\Phi^\prime := \Phi \circ \eta$ is an $\Omega^\prime$-co-ordinatization of $J$, and (\ref{ss.NOCELID}.\ref{NOCLA}) combined with \eqref{NOCLOM} shows $\kappa(\Omega^\prime) = \kappa(\Omega)$. Conversely, let this be so and suppose
\[
\Phi^\prime\:J \overset{\sim} \longrightarrow \Her_3(C,\pGamma), \quad \pGamma = \diag(\pgamma_1,\pgamma_2,\pgamma_3) \in \GL_3(F)
\]
is any $\Omega^\prime$-co-ordinatization of $J$. By Exc.~\ref{pr.DIAFIX}~(b), we may assume $\gamma_3 = \pgamma_3 = 1$. Then (\ref{ss.NOCELID}.\ref{NOCLA}) and \eqref{NOCLOM} show $\cl(-\gamma_i) = \cl(-\pgamma_i)$ for $i = 1,2$, so there are $p,q \in C^\times$ having $\gamma_1 = n_C(p)\pgamma_1$, $\gamma_2 = n_C(q)\pgamma_2$. and we deduce $C^{(p,q)} \cong C$ from Exc.~\ref{pr.CDISOT}. By Exc.~\ref{pr.ISCOPA}, therefore, we find a diagonal isomorphism $\eta\:\Her_3(C,\pGamma) \overset{\sim} \to \Her_3(C,\Gamma)$, whence $\eta \circ \Phi^\prime\:J \to \Her_3(C,\Gamma)$ is an $\Omega^\prime$-co-ordinatization of $J$. In particular $\Phi^{-1} \circ \eta \circ \Phi^\prime \in \Aut(J)$ sends $\pe_i$ to $e_i$ for $1 \le i \le3$. The final statement of the exercise follows from the fact that the class group of $J$ is trivial under the hypotheses stated.
\end{sol}

\begin{sol}{pr.HYPERLINE}\label{sol.HYPERLINE}
Recall from (\ref{e.CINNER}.\ref{e.CINNER.1}) that $(x \times y)^\sharp = T(x, y^\sharp) x$ for all $y \in J$.  So it suffices to find a $y \in J$ with $T(x, y^\sharp) \ne 0$.  And, since $J$ is regular, it suffices to argue that the elements $y^\sharp$ span $J$.  If $J$ has rank 1 this is obvious and if $J$ has rank $\ge 6$ it follows from Prop.~\ref{p.REDDEF}, Exc.~\ref{pr.SPAINV}~(b) and $Fy^\sharp = Fy^{-1}$ for all $y \in J^\times$.

It remains to consider the case where $J$ has rank 3. If $J$ is a division algebra, the claim is obvious so suppose $J = \plalg{(F \times K)}$ for $K$ a quadratic 
\'etale $F$-algebra. Applying (\ref{pr.CUADUN}.\ref{ADHAJ}), one checks
\[
(\alpha,u) = \big(\alpha - n_K(u)\big)(0,1_K)^\sharp + (1,\bar u)^\sharp
\]
for $\alpha \in F$, $u \in K$, and the assertion follows.
\end{sol}

\begin{sol}{pr.POP}\label{sol.POP}
\ref{POP.nz}: Suppose first that $y = e_{11}$.  For sake of contradiction, suppose $x \times y = 0$, so $x = \xi_1 e_{11} + u_2[31]+ u_3[12]$ for some $\xi_1 \in F$ and $u_2, u_3 \in C$.  Put $x_0 = u_2[31] + u_3[12]$.  Then $x_0^\sharp = (x - \xi_1 e_{11})^\sharp = x^\sharp + \xi^2_1 e_{11}^\sharp = 0$, and $x_0 \times y = (x - \xi_1 e_{11}) \times y = 0$.  However,
\[
x_0^\sharp = {-\gamma_3} \gamma_1 n_C(u_2)e_{22} - \gamma_1 \gamma_2 n_C(u_3) e_{33} + \gamma_1 \overline{u_2 u_3}[23].
\]
Since $n_C$ does not represent 0, $x_0 = 0$.  Thus $x = \xi_1 e_{11} \in Fy$, a contradiction, verifying that $x \times y \ne 0$.

In the general case, since $C$ is regular, Thm.~\ref{t.TRANSEL1} provides a $g \in \Inv(J)$ such that $gy \in \Fx e_{11}$.  Then 
\[
0 \ne gx \times gy = g^{\sharp - 1} (x \times y)
\]
by Lemma \ref{l.ISOTSIM}~\ref{ISOTSIM.sharp2},
so $x \times y \ne 0$.

The fact that $(x\times y)^\sharp = 0$ follows immediately from (\ref{ss.BACUNO.fig}.\ref{SADJ}).

\smallskip

\ref{POP.tr}:  
The easy direction is when $x = y \times z$ for some $z \in J$.  Then by (\ref{ss.JOPOID.fig}.\ref{DGRAD}) we have 
\[
T(x,y) = T(y\times z, y) = T(z, y \times y) = 0.
\]
So suppose $T(x,y) = 0$.

In the case when $y = e_{11}$, the fact that $T(x,y) = 0$ implies that $x = \xi_2 e_{22} + \xi_3 e_{33} + \sum_i u_i[jl]$.  The diagonal entries of $x^\sharp$ are in particular zero, so for each $i$ we have
\[
0 = \xi_j \xi_l - \gamma_j \gamma_l n_C(u_i).
\]
For $i = 2, 3$, $j$ or $l$ is 1, so the $\xi_j \xi_l$ term vanishes and we find that $n_C(u_i) = 0$.  Since $C$ is division, $u_i = 0$.  Thus $x = \xi_2 e_{22} + \xi_3 e_{33} + u_1[23]$, i.e., $x$ is in $y \times J$, compare Example \ref{e.MATINNER}.

In general, Thm.~\ref{t.TRANSEL1} provides a $g \in \Inv(J)$ such that $gy \in \Fx e_{11}$.  Then
\[
0 = T(x,y) = T(g^{\sharp-1} x, gy)
\]
by Lemma \ref{l.ISOTSIM}~\ref{ISOTSIM.T},
so $g^{\sharp-1} x = gy \times gz$ for some $z \in J$ and we conclude that $x = y \times z$.

\smallskip

\ref{POP.line}: By \ref{POP.nz}, $x \times y$ has rank 1, so $(x \times y) \times J$ is a line.  It contains $Fx$ because $T(x, x \times y) = T(x \times x, y) = 0$ and similarly for $Fy$.

So suppose $z \times J$ is a line containing $Fx$ and $Fy$.  If $z = e_{11}$, then by hypothesis $x$ and $y$ have zeros in their top row and left column.  It immediately follows that $x \times y$ is in $Fe_{11}$.  Moreover, $x \times y \ne 0$, so $x \times y \in \Fx e_{11}$.  For general $z$, Thm.~\ref{t.TRANSEL1} provides a $g \in \Inv(J)$ such that $gz \in \Fx e_{11}$.  Then $g^{\sharp -1}x$ and $g^{\sharp - 1}y$ are rank 1 elements contained in $gz \times J$, so $g^{\sharp -1}x \times g^{\sharp -1}y \in \Fx gz$, i.e., $x \times y$ is contained in $\Fx z$.  That is, $Fz = F(x \times y)$, proving uniqueness.

\ref{POP.PROP}: Among the three axioms, the first was verified in \ref{POP.line}.  Two lines $x \times J$, $y \times J$ pass through a point $Fz$ if and only if $Fz \subseteq (x \times J) \cap (y \times J)$ if and only if $T(z, x) = T(z,y) = 0$ if and only if $Fx, Fy \subseteq z \times J$.  It follows now from \ref{POP.line} that $Fz = F(x \times y)$, proving the second axiom.  

Finally, consider the four points $Fe_{11}$, $Fe_{22}$, $Fe_{33}$, and $Fx$ for
\[
x = \sum_i \gamma_i e_{ii} + 1[jl].
\]
(One checks that $x^\sharp = 0$.)  It remains to verify that no three of them are collinear.  Any subset of three of these points contains, say, $e_{jj}$ and $e_{ll}$, so the line they lie on is $(e_{jj} \times e_{ll}) \times J = e_{ii}\times J$, which contains neither $e_{ii}$ nor $x$.
\end{sol}

\begin{sol}{pr.SKONOFR} \label{sol.SKONOFR}
Let $\vph\:J_1 \to \pJ_1$ be any isomorphism. By definition (\ref{ss.REDFRE}), the algebra $J_1$ contains an elementary frame $\Omega$, which by Prop.~\ref{p.ELFRAFI} can be extended to a co-ordinate system
\[
\mfS = (\Omega,u_{23},u_{31})\;\text{of $J_1$}, \quad \Omega = (e_1,e_2,e_3).
\]
Applying $\vph$, we obtain a co-ordinate system
\[
\pmfS = (\pOmega,\pu_{23},\pu_{31})\;\text{of $\pJ_1$}, \quad \pOmega = (\pe_1,\pe_2,\pe_3),
\]
where $\pe_i = \vph(e_i)$ for $i = 1,2,3$ and $\pu_{jl} = \vph(u_{jl})$ for $i = 1,2$. Note that $\mfS$ and $\pmfS$ are both co-ordinate systems of $J$. Since $\vph$ preserves cubic norm structures, Prop.~\ref{p.CONALT} implies
\begin{align*}
\omega := \omega_{J_1,\mfS} = \omega_{J,\mfS} &= \omega_{\pJ_1,\pmfS} = \omega_{J,\pmfS}, \\
C_1 := C_{\pJ_1,\mfS} \subseteq C_{J,\mfS} = :C, \quad &\pC_1 := C_{\pJ_1,\pmfS} \subseteq C_{J,\pmfS} = :\pC
\end{align*}
as composition subalgebras and
\[
\Gamma := \Gamma_{\pJ_1,\mfS} = \Gamma_{J,\mfS} = \Gamma_{\pJ_1,\pmfS} = \Gamma_{J,\pmfS}.
\]
Moreover, $\vph$ restricts to a linear bijection from $(J_1)_{12}(\Omega)$ to $(\pJ_1)_{12}(\pOmega)$, which actually is an isomorphism $\vph_{12}\:C_1 \to \pC_1$. Now the Jacobson co-ordinatization theorem \ref{t.JACO} produces isomorphisms $\phi_1,\pphi_1,\phi,\pphi$ making the two top rows of the diagram
\[
\xymatrix{J & J_1 \ar[l]_{\iota} \ar[rr]_{\vph}^{\cong} && \pJ_1 \ar[r]_{\piota} & J \\
\Her_3(C,\Gamma) \ar[u]^{\phi}_{\cong} \ar@/_3pc/[rrrr]_{\Her_3(\psi,\Gamma)} & \Her_3(C_1,\Gamma) \ar[l]_{\iota_\Gamma}  \ar[u]^{\phi_1}_{\cong} \ar[rr]_{\Her_3(\vph_{12},\Gamma)} && \Her_3(\pC_1,\Gamma) \ar[u]^{\pphi_1}_{\cong} \ar[r]_{\piota_\Gamma} & \Her_3(\pC,\Gamma) \ar[u]^{\pphi}_{\cong}}
\]
commutative. By the Jacobson-Faulkner theorem \ref{t.ALBJAC}, $C$ and $\pC$ are isomorphic. Hence the Skolem-Noether theorem for composition algebras (Exc.~\ref{pr.SKONOCO}) implies that there exists an isomorphism $\psi\:C \to \pC$ extending $\vph_{12}$. Setting
\[
\Psi := \pphi \circ \Her_3(\psi,\Gamma) \circ \phi^{-1} \in \Aut(J),
\]
diagram chasing shows $\Psi \circ \iota = \piota \circ \vph$, which solves the problem.
\end{sol}

\begin{sol}{pr.INVSIX} \label{sol.INVSIX}
Regularity of Freudenthal $F$-algebras in dimension $6$ implies that $F$ has characteristic not $2$, so there is a natural identification of quadratic forms over $F$ with symmetric bilinear ones. We may write $J = \Her_3(F,\Gamma)$, $\pJ = \Her_3(F,\pGamma)$ for some invertible matrices $\Gamma = \diag(\gamma_1,\gamma_2,\gamma_3)$, $\pGamma = \diag(\pgamma_1,\pgamma_2,\pgamma_3)$ over $F$. By Lemma~\ref{l.SPLICQU}, the quadratic forms $Q_J$ and $Q_{\pJ}$ are isometric, and from (\ref{SJ0.Her3}.\ref{QUAMAT}) we deduce
\[
\la\gamma_2\gamma_3,\gamma_3\gamma_1,\gamma_1\gamma_2\ra \cong \la\pgamma_2\pgamma_3,\pgamma_3\pgamma_1,\pgamma_1\pgamma_2\ra.
\]
Maltiplying the left-hand side with $\gamma_1\gamma_2\gamma_3$ and the right-hand side with $\pgamma_1\pgamma_2\pgamma_3$, we see that $\la\Gamma\ra \cong \alpha\la\pGamma\ra$ for some $\alpha \in F^\times$. But $\pJ$ does not change when multiplying $\pGamma$ by an invertible scalar. Hence we may assume $\alpha = 1$ and thus find a $P \in \GL_3(F)$ such that $\pGamma = P^\trans\Gamma P$. Now one checks that the assignment $x \mapsto P^{-1}xP$ is an isomorphism from $J$ onto $\pJ$. 
\end{sol}


\solnchap{Solutions for Chapter~\ref{c.TICOS}}

\solnsec{Section~\ref{s.KUFIT}}

\begin{sol}{pr.ISCUAL} \label{sol.ISCUAL} Put $N := N_A$. Then $\pN := N(pq)N$ is a cubic form on $\pA := A^{(p,q)}$. From (\ref{ss.ISTJALT}.\ref{ISTJALT}) we deduce $A^{\prime(+)} = A^{(+)(pq)}$. By Prop.~\ref{p.ITIS} and Cor.~\ref{c.ISNOJO}, therefore, $A^{\prime(+)}$ is a cubic Jordan algebra over $k$ with norm $\pN$. Since $\pN$ permits composition relative to the multiplication of $\pA$, it follows that $\pA$ is, in fact, a cubic alternative algebra with norm $\pN$ over $k$. Now let $x \in \pA$. Then Prop.~\ref{p.CUALJO}, (\ref{ss.ISCUNO}.\ref{ISOQUAT}) and Prop.~\ref{p.FOCUALT} imply
\[
S_{\pA}(x) = S_{A^{\prime(+)}}(x) = S_{A^{(+)(pq)}}(x) = T_{A^{(+)}}\big((pq)^\sharp,x^\sharp\big) = T_A(q^\sharp p^\sharp x^\sharp),
\]
hence the second equation of \eqref{ISCUAL}. Similarly,
\[
T_{\pA}(x) = T_{{\pA}^{(+)}}(x) = T_{A^{(+)(pq)}}(x) = T_{A^{(+)}}(pq,x) = T_A(pqx).
\]
In particular, $A$ and $A^p$ have the same linear and quadratic trace. Thanks to Prop.~\ref{p.FOCUALT}, the bilinear trace of $A^p$ satisfies $T_{A^p}(x,y) = T_{A^p}((xp^{-1})(py)) = T_A((xp^{-1})(py))) = T_A(([xp^{-1}]p)y) = T_A(x,y)$ and thus agrees with the bilinear trace of $A$.

\end{sol} 
	
\begin{sol}{pr.CUIDRE} \label{sol.CUIDRE}  (a) Since $(\mfa,I)$ is a cubic ideal in $A^{(+)}$, we have $\mfa A \subseteq I$ and $T_A(x,y), T_A(x^\sharp,y), N_A(x) \in \mfa$ for all $x \in I$, $y \in A$. Now let $x_i \in I$, $y_i \in A$, $i = 0,1,2$, and $x := x_0 + x_1j_1 + x_2j_2 \in J(I,\mu)$, $y := y_0 + y_1j_1 + y_2 j_2 \in J$. Then 
\begin{align}
\label{TEJAX}T_J(x,y) =\,\,&T_A(x_0,y_0) + \mu T_A(x_1,y_2) + \mu T_A(x_2,y_1) \\
\label{TEJASH} T_J(x^\sharp,y) =\,\,&T_A(x_0 - \mu x_1x_2,y_0) + \mu T_A(\mu x_2^\sharp - x_0x_1,y_2) + \mu T_A(x_1^\sharp - x_2x_0,y_1) \notag \\
=\,\,&T_A(x_0^\sharp,y_0 - \mu T_A(x_1,x_2y_0) + \mu^2T_A(x_2^\sharp,y_2) \\
\,\,&- \mu T_A(x_0,x_1y_2) + \mu T_A(x_1^\sharp,y_1) - \mu T_A(x_2,x_0y_1), \notag \\
\label{NJAX} N_J(x) =\,\,& N_A(x_0) + \mu N_A(x_1) + \mu^2 N_A(x_2) - \mu T_A(x_0,x_1x_2)
\end{align} 
all belong to $\mfa$. Moreover, $I^\sharp + I \times A \subseteq I$ by Exc.~\ref{pr.NORID}~(a), which by (\ref{t.EXTITS}.\ref{ADXFITI}), (\ref{t.EXTITS}.\ref{BADFITI}) implies $J(I,\mu)^\sharp + J(I,\mu) \times J \subseteq J(I,\mu)$. Hence $(\mfa,J(I,\mu))$ is a cubic ideal in $J$ (loc. cit.). Next suppose $(\mfa,I)$ is a cubic nil ideal in $A$, i.e., $\mfa$ and $I$ are both nil (Exc.~\ref{pr.RENIRA}). Then \eqref{TEJAX}\emptyslot\eqref{NJAX} combined with Exc.~\ref{pr.CUBNIL} show that the ideal $J(I,\mu) \subseteq J$ is nil. Thus $(\mfa,J(I,\mu))$ is a cubic nil ideal in $J$. On the other hand, assume $(\mfa,I)$ is separated, so some $k$-linear map $\lambda\:A \to k/\mfa$ has $\lambda(1_A) = 1_{k/\mfa}$ and $\lambda(I) = \{0\}$. Let $\plambda\:J \to k/\mfa$ be the projection from $J$ onto the initial summand followed by $\lambda$. Then $\plambda$ is $k$-linear, sends $1_J$ to $1_{k/\mfa}$ and kills $J(I,\mu)$. Hence $(\mfa,J(I,\mu))$ is separated.  Moreover, $A/I$ is an alternative algebra over $k/\mfa$, whence Exc.~\ref{pr.NORID}~(b) implies that $A^{(+)}/I = (A/I)^{(+)}$ carries a unique cubic Jordan algebra structure over $k/\mfa$ making the canonical projection $\pi\:A^{(+)} \to (A/I)^{(+)}$ a $\sigma$-semi-linear homomorphism of cubic Jordan algebras. Hence $A/I$ carries a unique structure of a cubic alternative algebra over $k/\mfa$ making $\pi$ a $\sigma$-semi-linear homomorphism of cubic alternative algebras. The final assertion of (a) is obvious.

(b) The first statement is clear while \eqref{JAYNIL} follows from the fact that $J(\Nil(A),\mu)$ by (a) is a nil ideal in $J$. Now suppose $k \neq \{0\}$ and $\mu$ is nilpotent. Setting $I_1 := \Nil(A) \oplus Aj_1 \oplus Aj_2$ and combining (\ref{t.EXTITS}.\ref{ADXFITI}), (\ref{t.EXTITS}.\ref{BADFITI}) with Exc.~\ref{pr.PARNIL}~(c), we conclude $I_1^\sharp + I_1 \times J \subseteq I_1$. By Exc.~\ref{pr.NORID}~(a) and  \eqref{TEJAX}\emptyslot\eqref{NJAX}, therefore, $I_1$ is a nil ideal in $J$ properly containing $J(\Nil(A),\mu)$. Hence we can't have equality in \eqref{JAYNIL}. Next let $F$ be a field of characteristic $3$ and $J = J(F,1)$. Then $\Nil(F) = \{0\} \neq \Nil(J)$ since, for example, $(j_1 - 1_J)^3 = j_1^3 - 1_J = 0$.

(c) By \eqref{JAYNIL}, the right-hand side of \eqref{NILJAY} is contained in the left. Conversely, put $\mfa := \Nil(k)$, $I := \Nil(A)$. By \eqref{JAYAIM},
\[
J(A,\mu)/J(I,\mu) \cong J\big(A/I,\sigma(\mu)\big),
\]
 where the right-hand side is regular by Cor.~\ref{c.REFITI}. But the nil radical of $J(A/I,\sigma(\mu))$ is zero by regularity and since $k_0$ is reduced. Hence $\Nil(J(A,\mu)) \subseteq J(I,\mu)$, and \eqref{NILJAY} is proved.
\end{sol}

\begin{sol}{pr.SESIAL} \label{sol.SESIAL} (a) Let $A$ be a semi-simple cubic alternative $F$-algebra. Then $J := A^{(+)}$ by Cor.~\ref{c.NICUAL} is a semi-simple cubic Jordan algebra over $F$, hence falls into precisely one of the categories listed in Thm.~\ref{t.RACUSE}. We now proceed in several steps. 

\step{1}
Suppose first $J$ is as in (ii) of \ref{t.RACUSE}, so up to isomorphism we have
\begin{align}
\label{JAYEFP} J = F^{(+)} \times J(M,q,e)
\end{align}
as a direct product of ideals, $(M,q,e)$ being a non-degnenerate pointed quadratic module over $F$. Then $c := (1,0) \in J$ by \eqref{ADHAJ}, \eqref{LITHAJ} of Exc.~\ref{pr.CUADUN} is an elementary idempotent satisfying $J_2(c) = F^{(+)}$, $J_1(c) = \{0\}$, $J_0(c) = J(M,q,e)$. Hence Example~\ref{e.PEALJO} implies $A_{11}(c) = F$, $A_{12}(c) = A_{21}(c) = \{0\}$, $A_{22}(c)^{(+)} = J(M,q,e)$, and the Peirce rules for alternative algebras (Exc.~\ref{pr.PEIRCEALT}) show that \eqref{JAYEFP} is in fact a decomposition into a direct product of ideals in $A$:
\begin{align}
\label{AEF} A = F \times C,
\end{align}
where $C := A_{22}(c)$ is a conic alternative $F$-algebra with identity element $1_C = e$ and norm $n_C = q$. By \eqref{NOHAJ} of Exc.~\ref{pr.CUADUN} we have $N_A((1,u)) = n_C(u)$ for all $u \in C$, and since $N_A$ permits composition, so does $n_C$. Summing up, we have shown that $C$ is a pre-composition algebra over $F$, forcing $A = \hC$ to be as in (iv) or (v) depending on whether $C$ is split quadratic \'etale or not. 

\step{2}
Next suppose $J$ is as in (iii) of \ref{t.RACUSE}, so $J = \Her_3(C,\Gamma)$ for some pre-composition algebra $C$ over $F$ and some $\Gamma = \diag(\gamma_1,\gamma_2,\gamma_3) \in \GL_3(F)$; in particular, $\dim_F(C) = 2^n$, $n \in \IN$. Arguing as in $1^\circ$,
\[
\pJ := Fe_{11} \oplus J_0(e_{11}) = Fe_{11} \oplus (Fe_{22} \oplus C[23] \oplus Fe_{33}) \cong F^{(+)} \times J(M,q,e),
\]
where $M = Fe_{22} \oplus C[23] \oplus Fe_{33}$ as a $k$-module, $e = e_{22} + e_{33}$ and $q\:M \to F$ is defined by
\begin{align}
\label{QUALPH} q(\alpha_2e_{22} + u_1[23] + \alpha_3e_{33}) := \alpha_2\alpha_3 - \gamma_2\gamma_3n_C(u_1)
\end{align}
for $\alpha_2,\alpha_3 \in F$, $u_1 \in C$, carries the structure of a cubic alternative $F$-algebra $\pA$ having ${\pA}^{(+)} = \pJ$ as cubic Jordan algebras. By $1^\circ$, therefore, $M$ becomes a pre-composition algebra $\pC$ over $F$ such that $n_{\pC} = q$ and $\pA = \hat{\pC}$; in particular, $2^n + 2 = \dim_F(\pC) = 2^m$ for some $m \in \IN$, which is impossible unless $n = 1$. But the norm of a pre-composition algebra over a field is either anisotropic or hyperbolic. Hence \eqref{QUALPH} implies that $C$ is split quadratic \'etale, and we deduce from Prop.~\ref{p.FRSPCOM} that
\[
A^{(+)} \cong \Her_3(F \times F,\Gamma) \cong \Her_3(F \times F) \cong \Mat_3(F)^{(+)}.
\] This isomorphism extends to one over the algebraic closure $\aclosF$ of $F$, and since $A_{\aclosF}^{(+)} \cong \Mat_3(\aclosF)^{(+)}$ is simple, so is $A_{\aclosF}$. From Cor.~\ref{ss.c-CENTSIM4} we therefore conclude that $A$ is central simple, and Thm.~\ref{t.CESIM} implies $A \cong \Mat_3(F)$, first as abstract alternative algebras and then as cubic ones, thanks to Cor.~\ref{c.ALUNO}. Thus we are in case (vi). 

\step{3}
Finally, we are able to deal with case (i) of \ref{t.RACUSE}, so $J$ is a cubic Jordan division algebra, forcing $A$ to be an alternative division algebra. Since the radical of the bilinear trace of $A$ by Prop.~\ref{p.FOCUALT}~(a) is a two-sided ideal in $A$, we either have $T_A = 0$ or $A$ is regular.

Assume first $T_A = 0$. Then $F$ has characteristic $3$ and (\ref{ss.BACUNO.fig}.\ref{EXPF}) implies that $N_A\:A \to F$ is a homomorphism of algebras over $\IZ$. Moreover, $A$ being a division algebra, this homomorphism is injective, forcing $A/F$ to be a purely inseparable field extension of exponent at most $1$. Thus we are either in case (i) or in case (ii).

We are left with the case that $A$ is regular. The $\pA:= A_{\aclosF}$ is a semi-simple cubic alternative algebra over $\aclosF$, so $\pJ:= {\pA}^{(+)}$ falls into one of the categories (i), (ii), (iii) of \ref{t.RACUSE}. If $\pJ$ satisfies (i), then $\pA\cong \aclosF$ by Exc.~\ref{pr.FIDIA}, hence $A \cong F$, and we are in case (i). If $\pJ$ satisfies (ii) of \ref{t.RACUSE}, so $\pJ= \aclosF^{(+)} \times J(\pM,\pq,\pe)$ is a direct product of ideals, where $(\pM,\pq,\pe)$ is a non-degenerate pointed quadratic module over $\aclosF$, then we may assume $\dim_{\aclosF} (\pM) \geq 3$ and $c := (1_{\aclosF},0)$ by Prop.~\ref{p.PEPOQUA} is the only elementary idempotent of $\pJ$ having $\pJ_1(c) = \{0\}$. Thus, by descent, $c$ belongs to $J$, and $J$ satisfies condition (ii) of \ref{t.RACUSE}. Now $2^\circ$ implies $A = \hC$ for some pre-composition algebra $C$ over $F$, which is actually  a composition algebra since $A$ is regular. Finally, if $\pJ$satisfies (iii) of \ref{t.RACUSE}, then $\pA$ is simple of dimension $9$ over $\aclosF$ by $2^\circ$, forcing $A$ to be central simple of dimension $9$ over $F$. By Thm.~\ref{t.CESIM}, therefore, we are in case (vi) or (vii) of our Exc.

\smallskip

(b) (i) $\Rightarrow$ (ii). Regular cubic alternative algebras are separable by Exc.~\ref{pr.CUBNIL} and Cor.~\ref{c.NICUAL}. Moreover, $k$ (resp. $(k \times k)_{\cub}$) are both separable but not regular unless $\frac{1}{3}$ (resp. $\frac{1}{2}$) belong to $k$.

(ii) $\Rightarrow$ (i). In view of (a), (b) of (ii) combined with Exc.~\ref{pr.CUMOPO}, we may assume $n \geq 3$. Note that $A$ is regular (resp. separable) if and only if so is $A_F$, for any algebraically closed field $F \in \kalg$ (for the regular case, see Exc.~\ref{pr.CRITREG}). Hence we may assume that $F$ is an algebraically closed field, so $A$ satisfies one of the conditions (i)\emptyslot(vii) of (a). Since $F$ is algebraically closed, only cases (iv)\emptyslot(vi) are possible, where $C$ in (v) has dimension at least $2$ and hence is regular. But then so is $A$ and we are done.
\end{sol}

\begin{sol}{pr.CHAMU} \label{sol.CHAMU}  We will repeatedly make use of the fact that a linear map between cubic Jordan algebras is a homomorphism if and only if it preserves unit elements and adjoints (Exc.~\ref{pr.HOCUJO}~(a)). 
	
(a) $\vph := \vph_{A,\mu,p}$ is the identity on $A^{(+)}$ and, in particular, preserves units. Furthermore, since $A$ is associative,
\begin{align*}
\vph(x)^\sharp =\,\,&\big(x_0 + (x_1p)j_1 + (p^\sharp x_2)j_2\big)^\sharp \\
=\,\,&\big(x_0^\sharp - N_A(p)\mu x_1x_2\big) + \big(N_A(p)\mu x_2^\sharp p - (x_0x_1)p\big)j_1 + \big(p^\sharp x_1^\sharp - p^\sharp x_2x_0\big)j_2 \\
=\,\,&\vph(x^\sharp),
\end{align*}
so $\vph\:J(A,N_A(p)\mu) \to J(A,\mu)$ is a homomorphism of cubic Jordan algebras. The final assertion is clear. 

(b) In order to show for $x \in A$ that $x[A,A] = \{0\}$ implies $x = 0$, we may assume that $k$ is a local ring with maximal ideal $\mfm$. Since $[A,A](\mfm) = [A(\mfm),A(\mfm)]$ is the kernel of the reduced trace of the central simple associative algebra $A(\mfm)$ of degree $3$ over the field $k(\mfm)$, it contains invertible elements. Hence so does $[A,A]$ over $k$, and we conclude $x = 0$, as desired.  

(ii) $\Rightarrow$ (i). Follows immediately from (a).

(i) $\Rightarrow$ (ii). Observe $A^{(+)\perp} = Aj_1 + Aj_2$ in $J(A,\nu)$ and $J(A,\mu)$, respectively, by (\ref{t.EXTITS}.\ref{BTRFITI}) since $\mu,\nu \in k$ are invertible. Moreover, given $x,y,z_1,z_2 \in A$ and setting $z := z_1j_1 + z_2j_2$, we may apply  (\ref{ss.FOFITI}.\ref{EXNAPT}) to obtain
\begin{align*}
x\pt (y\pt z) = x\pt \big((yz_1)j_1 + (z_2y)j_2\big) = (xyz_1)j_1 + (z_2yx)j_2 = (xy)\pt z - \big(z_2[x,y]\big)j_2.
\end{align*}
By what we have just shown, therefore, $x\pt (y\pt z) = (xy)\pt z$ for all $x,y \in A$ if and only if $z_2 = 0$, and we conclude
\begin{align}
\label{AJAYONE} Aj_1 = \{z \in A^{(+)\perp} \mid x\pt (y\pt z) = (xy)\pt z\;\text{for all $x,y \in A$}\}.
\end{align}
Since $\vph$ preserves traces, it maps $A^{(+)\perp} \subseteq J(A,\nu)$ to $A^{(+)\perp} \subseteq J(A,\mu)$, and since $\vph$ preserves bilinearized adjoints, we have $\vph(x\pt z) = x\pt \vph(z)$ for all $x \in A^{(+)}$, $z \in A^{(+)\perp}$. Hence \eqref{AJAYONE} implies that $\vph$ sends $Aj_1 \subseteq J(A,\nu)$ to $Aj_1 \subseteq J(A,\mu)$. In particular, $\vph(j_1) = pj_1$ for some $p \in A$, which implies
\[
N_A(p)\mu = N_{J(A,\mu)}(pj_1) = N_{J(A,\mu)}\big(\vph(j_1)\big) = N_{J(A,\nu)}(j_1) = \nu,
\]
proving the first part of (b). But we also have $\vph(j_2) = \vph(j_1^\sharp) = \vph(j_1)^\sharp = (pj_1)^\sharp = p^\sharp j_2$, and for $z_1,z_2 \in A$, we conclude
\[
\vph(z_1j_1 + z_2j_2) = z_1\pt (p\pt j_1 + z_2\pt (p^\sharp\pt j_2)) = (z_1p)j_1 + (p^\sharp z_2)j_2.
\]
Hence $\vph = \vph_{A,\mu,p}$, which is an isomorphism by (a) since $p \in A^\times$. 

(c) Put $\psi := \psi_{A,\mu,p}$. Then
\begin{align*}
\psi(x^\sharp) =\,\,&\psi\Big(\big(x_0^\sharp - \mu(x_1p^{-1})(px_2)\big) + \big(\mu x_2^\sharp - (x_0p^{-1})(px_1)\big)j_1 \\
\,\,&+ \big(x_1^\sharp - (x_2p^{-1})(px_0)\big)j_2\Big) \\
=\,\,&\big(x_0^\sharp - \mu(x_1p^{-1})(px_2)\big) + \big(\mu x_2^\sharp p^{-1} - [(x_0p^{-1})(px_1)]p^{-1}\big)j_1 \\
\,\,&+ \big(N_A(p)^{-1}px_1^\sharp - N_A(p)^{-1}p[(x_2p^{-1})(px_0)]\big)j_2 \\
=\,\,&\Big(x_0^\sharp - N_A(p)\mu\big(x_1p^{-1}\big)\big(N_A(p)^{-1}px_2\big)\Big) \\
\,\,&+ \Big(N_A(p)\mu\big(N_A(p)^{-1}px_2\big)^\sharp - x_0(x_1p^{-1})\Big)j_1 \\
\,\,&+ \Big((x_1p^{-1})^\sharp - \big(N_A(p)^{-1}px_2\big)x_0\Big)j_2 \\
=\,\,&\Big(x_0 + (x_1p^{-1})j_1  + \big(N_A(p)^{-1}px_2\big)j_2\Big)^\sharp = \psi(x)^\sharp.
\end{align*}
Thus $\psi$ is an isomorphism of cubic Jordan algebras. It remains to show that the diagram \eqref{COMPSI} commutes:
\begin{align*}
\psi_{A,N_A(q)\mu,p} \circ\,\,& \psi_{A^p,\mu,q}(x) \\
=\,\,&\psi_{A,N_A(q)\mu,p}\big(x_0 + [(x_1p^{-1})(pq^{-1})]j_1 + N_A(q)^{-1}[(qp^{-1})(px_2)]j_2\big) \\
=\,\,&x_0 + \big([(x_1p^{-1})(pq^{-1})]p^{-1}\big)j_1 + N_A(p)^{-1}N_A(q)^{-1}\big(p[(qp^{-1})(px_2)]\big)j_2 \\
=\,\,&x_0 + \big(x_1(pq)^{-1}\big)j_1 + N_A(pq)^{-1}\big((pq)x_2\big)j_2 = \psi_{A^{pq},\mu,pq}(x),
\end{align*} 
as desired.

(d) We apply Cor.~\ref{c.UPFITI} to $J := J(A^{\op},\mu^{-1})$, $J_0 := A^{\op (+)} \subseteq J$, $V := A^{\op}j_1 \oplus A^{\op}j_2 \subseteq J$ and $l := \mu j_2 = \mu j_1^\sharp \in J$. Then $l$ is a Kummer element of $J$ relative to $(J_0,V)$ satisfying $l^\sharp = \mu^2j_2^\sharp = \mu^2\mu^{-1}j_1 = \mu j_1$, and Thm.~\ref{t.INTITS}~(b) combined with (\ref{ss.FOFITI}.\ref{JAYONE}) implies
\[
A^{\op} = A_{j_1}(J,J_0,V) = A_{l^\sharp}(J,J_0,V) = A_l(J,J_0,V)^{\op},
\]
hence $A = A_l(J,J_0,V)$. We also have $N_J(l) = \mu^3N_J(j_2) = \mu^3\mu^{-2} = \mu$, and Cor.~\ref{c.UPFITI} yields a unique homomorphism
\[
\vph\:J(A,\mu) \longrightarrow J = J(A^{\op},\mu^{-1}) 
\]
of cubic Jordan algebras inducing the identity on $A^{(+)} = A^{\op(+)} = J_0$ and satisfying
\[
\vph(x_0 + x_1j_1 + x_2j_2) = x_0 + x_1\pt l + x_2\pt l^\sharp = x_0 + \mu x_2j_1 + \mu x_1j_2
\]
for $x_0,x_1,x_2 \in A$. Clearly, $\vph$ is an isomorphism.
\end{sol}

\begin{sol}{pr.NONKUM} \label{sol.NONKUM}  We begin with recalling some standard notation from linear algebra. Let $n$ be a positive integer. We denote by $k^n$ rank-$n$ column space over $k$ and by $(e_i)_{1\leq i\leq n}$ its canonical basis. Analogously, we denote by $k_n$ rank-$n$ row space over $k$ and have $(e_i^\trans)_{1\leq i\leq n}$ as its canonical basis. Matrices $x,y \in \Mat_n(k)$ can be written as columns of their row vectors, i.e.,
\begin{align}
\label{COLROW} x = \left(\begin{matrix}
x_1 \\
\vdots \\
x_n
\end{matrix}\right), \quad x_i = e_i^\trans x \in k_n, \quad y = \left(\begin{matrix}
y_1 \\
\vdots \\
y_n
\end{matrix}\right), \quad y_i = e_i^\trans y \in k_n &&(1 \leq i \leq n),
\end{align}
ditto for $xy \in \Mat_n(k)$. If $x = (\xi_{ij})_{1\leq i,j\leq n}$, then \eqref{COLROW} implies
\begin{align}
\label{XYROW} (xy)_i = e_i^\trans xy = x_iy = \sum_{j=1}^n\xi_{ij}e_j^\trans y = \sum_{j=1}^n\xi_{ij}y_j &&(1 \leq i \leq n).
\end{align}

Let us now get back to to the exercise proper. Note first that the off-diagonal components in
\[
C = \Zor(k) = \left(\begin{matrix}
k & k^3 \\
k^3 & k
\end{matrix}\right)
\]
are rank-$3$ \emph{column} vectors. Given $y \in \Mat_3(k)$ as in \eqref{COLROW} for $n = 3$, we put
\begin{align}
\label{ELY} l_y = \sum \left(\begin{matrix}
0 & 0 \\
y_i^\trans  & 0
\end{matrix}\right)[jl] \in \Her_3(C) = J.
\end{align}
We now identify $\Mat_3(k)^{(+)} = \Her_3(k \times k)$ via Prop.~\ref{p.HERKEI} and the split quadratic \'etale $k \times k$ with the diagonal of $C = \Zor(k)$. Then another matrix $x = (\xi_{ij})_{1\leq i,j\leq 3} \in \Mat_3(k)$ may be viewed canonically as an element of $J$, and we claim
\begin{align}
\label{EXTIY} x \times l_y = -l_{xy}.
\end{align}
To see this, we combine (\ref{ss.THERCU}.\ref{MABADJ}) with \eqref{XYROW} and obtain
\begin{align*}
-x \times l_y =\,\,&-\sum\Big(\xi_ie_{ii} + \left(\begin{matrix}
\xi_{jl} & 0 \\
0 & \xi_{lj} 
\end{matrix}\right)[jl]\Big) \times \sum \left(\begin{matrix}
0 & 0 \\
y_i^\trans  & 0
\end{matrix}\right)[jl] \\
=\,\,&\sum\Big(\xi_{ii}\left(\begin{matrix}
0 & 0 \\
y_i^\trans  & 0
\end{matrix}\right) - \overline{\left(\begin{matrix}
\xi_{li} & 0 \\
0 & \xi_{il}
\end{matrix}\right)\left(\begin{matrix}
0 & 0 \\
y_l^\trans  & 0
\end{matrix}\right)} - \overline{\left(\begin{matrix}
0 & 0 \\
y_j^\trans  & 0
\end{matrix}\right)\left(\begin{matrix}
\xi_{ij} & 0 \\
0 & \xi_{ji}
\end{matrix}\right)}\Big)[jl] \\
=\,\,&\sum\left(\begin{matrix}
0 & 0 \\
\xi_{ii}y_i^\trans  + \xi_{il}y_l^\trans  + \xi_{ij}y_j^\trans  & 0
\end{matrix}\right)[jl] \\
=\,\,&\sum\left(\begin{matrix}
0 & 0 \\
(\sum_{\nu=1}^3 \xi_{i\nu}y_\nu)^\trans  & 0
\end{matrix}\right)[jl] = \sum\left(\begin{matrix}
0 & 0 \\
(xy)_i^\trans  & 0
\end{matrix}\right)[jl] = l_{xy},
\end{align*}
as claimed.

We now put
\begin{align*}
u_i := \left(\begin{matrix}
0 & 0 \\
e_i & 0
\end{matrix}\right) \in C^0 \ni \left(\begin{matrix}
0 & e_i \\
0 & 0
\end{matrix}\right) =: v_i &&(1 \leq i \leq 3),
\end{align*}
where $C^0$ stands for the $k$-module of trace-zero-elements in $C$. Adopting the ternary cyclicity convention \ref{ss.TERCYC} and applying (\ref{ss.ZOVE}.\ref{ZOMU}), we conclude
\begin{align}
\label{UIUJVL} u_iu_j = v_l, \quad v_iv_j = u_l &&(1 \leq i \leq 3).
\end{align}
Now put
\begin{align}
\label{ELONE} l := l_{\Eins} =\sum u_i[jl] 
\end{align}
where $n_C(u_i) = 0$ and \eqref{UIUJVL} combined with (\ref{ss.THERCU}.\ref{MADJ}) imply
\begin{align}
\label{ELSHA} l^\sharp = \sum \bar v_i[jl] = -\sum v_i[jl].
\end{align}
Moreover, by (\ref{ss.THERCU}.\ref{MANO}), \eqref{UIUJVL}, (\ref{ss.ZOVE}.\ref{ZOMU}),
\begin{align}
\label{ENEL} N_J(l) = t_C(u_1u_2u_3) = t_C(u_1v_1) = t_C(\left(\begin{matrix}
0 & 0 \\
0 & -1
\end{matrix}\right)) = -1,
\end{align}
so $l$ is invertible in $J$.

(a) We put $J_0 := \Mat_3(k)^{(+)} = \Her_3(k \times k) \subseteq J$ as a regular cubic Jordan subalgebra and have
\[
J_0^\perp = \sum\left(\begin{matrix}
0 & k^3 \\
k^3 & 0
\end{matrix}\right)[jl],
\] 
whence \eqref{ELONE}, \eqref{ELSHA} imply that $l$ is strongly orthogonal to $J_0$. It is also invertible in $J$ and from \eqref{EXTIY} we deduce $x\pt (y\pt l) = x\pt (y\pt l_{\Eins}) = x\pt l_y = l_{xy} = (xy)\pt l$ for all $x,y \in \Mat_3(k)$, i.e., the stability condition for $l$. In other words, $l$ is a Kummer element of $J$ relative to $J_0$ and $A_l(J,J_0) = \Mat_3(k)$. Since the norm of $\Mat_3(k)$ is surjective, (a) follows from Cor.~\ref{c.UPFITI} combined with Exc.~\ref{pr.CHAMU}~(a).

(b) The regularity of $J_0 := \Sym_3(k) \subseteq J$ is clear, and since $l$ is an ivertible element of $J$ strongly orthogonal to $\Mat_3(k)^{(+)}$, it is so to $J_0$. But if $l$ were a Kummer element relative to $J_0$, then what we have seen in (a) would imply that $\Sym_3(k) \subseteq \Mat_3(k)$ were an \emph{associative} subalgebra, a contradiction.
\end{sol}

\begin{sol}{pr.NILTWO} \label{sol.NILTWO}
By Lemma~\ref{l.NIJONU}, there exists an elementary idempotent $e \in J$ such that $u \in J_0(e)$. Since $J$ is split, so the class group of $J$ is trivial (Lemma~\ref{l.NCLEID}), the automorphism group of $J$ acts transitively on the elementary idempotents (Cor.~\ref{c.ORBID}). Hence we may assume $J = \Her_3(C)$, $C := \Zor(F)$, $e = e_{11}$. Consulting the Peirce rules (\ref{t.PEDECOM}.\ref{PEDECOM}) and (\ref{p.PEDEDI}.\ref{PEDEDI}) combined with the fact that $u$ has trace zero (Exc.~\ref{pr.CUBNIL}), we deduce from (\ref{ss.THERCU}.\ref{MALIT}) that $u = \xi(e_{22} - e_{33}) + v[23]$, for some $\xi \in F$ and some $v \in C$.
\begin{itemize} 
\item [($\ast$)] We claim that \emph{there exists a $w \in C^\times$ such that $w^{-1}v$ is contained in a (split) quadratic \'etale subalgebra of $C$.}
\end{itemize}
If $v$ is invertible, then $w := v$ does the job. If not, then Zariski density yields a $w \in C^\times$ having $t_C(w^{-1}v) = 1$. Since $n_C(w^{-1}v) = n_C(w)^{-1}n_C(v) = 0$, we conclude that $w^{-1}v$ is an elementary idempotent in $C$, and ($\ast$) follows.

Replacing $w$ by $\sqrt{w}^{-1}w$ if necessary, we may assume $n_C(w) = 1$. Consulting (\ref{ss.THERCU}.\ref{MANO}), we see that the assignment
\[
\sum (\xi_ie_{ii} + v_i[jl]) \longmapsto \sum \xi_ie_{ii} + (w^{-1}v_1)[23] + (v_2w^{-1}[31] + (wv_3w)[12]
\]
gives a linear bijection $J \to J$ that preserves unit element and norm, hence is an automorphism of $J$ (Exc.~\ref{pr.HOCUJO}). By ($\ast$), we may therefore assume that $v$ itself belongs to a split quadratic \'etale subalgebra of $J$. But this means that $u = \xi(e_{22} - e_{33}) + v[23]$ belongs to a subalgebra of $J$ isomorphic to $\Her_3(F \times F) \cong \Mat_3(F)^{(+)}$. By Exercises \ref{pr.NONKUM} and \ref{pr.SKONOFR}, therefore, we may assume $J = J(\Mat_3(F),1)$ as a first Tits construction such that $u \in \Mat_3(F)^{(+)}$. Since nilpotent elements of $A := \Mat_3(F)$ having index $2$ (i.e., rank $1$) are conjugate under (inner) automorphisms of $A$, and every automorphism $\vph$ of $A$ extends to the automorphism $J(\vph,1)$ of $J$, the problem is solved.  
\end{sol}

\begin{sol}{pr.ELERTI} \label{sol.ELERTI}  We put $J := J(A,\mu)$. Then Thm.~\ref{t.EXTITS} yields
\[
U_px = T_J(p,x)p - p^\sharp \times x = T_A(p,x_0)p - p^\sharp \times x_0 + (p^\sharp x_1)j_1 + (x_2p^\sharp)j_2,
\]
and \eqref{UOPA} follows from Prop.~\ref{p.CUALJO}. Combining (\ref{ss.ISCUNO}.\ref{ISOAD}) with \eqref{UOPA}, we conclude
\begin{align*}
x^{(\sharp,p)} =\,\,&N_A(p)U_{p^{-1}}x^\sharp \\
=\,\,&N_A(p)\Big(p^{-1}(x_0^\sharp - \mu x_1x_2)p^{-1} + \big(p^{-1\sharp}(\mu x_2^\sharp - x_0x_1)\big)j_1 \\
\,\,&+ \big((x_1^\sharp - x_2x_0)p^{-1\sharp}\big)j_2\Big),
\end{align*}
and \eqref{ADJPA} has been proved. That $\til_p$, $\tir_p$ are isomorphisms of cubic Jordan algebras as indicated could now be proved by using \eqref{ADJPA} and showing that they preserve adjoints. Leaving this to the reader, we prefer to adopt a different approach and begin by deriving \eqref{UTILTIR}. The first two equations follow immediately from the left (resp. right) Moufang identity. For the third, we apply (\ref{ss.UOPALT}.\ref{ULR}) and obtain
\begin{align*}
\til_pU_q\tir_p(x) =\,\,&\til_pU_q\big(x_0p + (p^{-1}x_1p^\sharp)j_1 + (px_2)j_2\big) \\
=\,\,&\til_p\Big(q(x_0p)q + \big(q^\sharp(p^{-1}x_1p^\sharp)\big)j_1 + \big((px_2)q^\sharp\big)j_2\Big) \\
=\,\,&p\big(q(x_0p)q\big) + \Big(\big(q^\sharp(p^{-1}x_1p^\sharp)\big)p\Big)j_1 + \Big(p^\sharp\big((px_2)q^\sharp\big)p^{-1}\Big)j_2 \\
=\,\,&(pq)x_0(pq) + \big((pq)^\sharp x_1\big)j_1 + \big(x_2(pq)^\sharp\big)j_2 = U_{pq}(x)
\end{align*}
and, similarly, 
\begin{align*}
\tir_qU_p\til_q(x) =\,\,&\tir_qU_p\big(qx_0 + (x_1q)j_1 + (q^\sharp x_2q^{-1})j_2\big) \\
=\,\,&\tir_q\Big(p(qx_0)p + \big(p^\sharp(x_1q)\big)j_1 + \big((q^\sharp x_2q^{-1})p^\sharp\big)j_2\Big) \\
=\,\,&\big(p(qx_0)p\big)q + \Big(q^{-1}\big(p^\sharp(x_1q)\big)q^\sharp\Big)j_1 + \Big(q\big((q^\sharp x_2q^{-1})p^\sharp\big)\Big)j_2 \\
=\,\,&(pq)x_0(pq) + \big((pq)^\sharp x_1\big)j_1 + \big(x_2(pq)^\sharp\big)j_2 = U_{pq}(x).
\end{align*}
Thus \eqref{UTILTIR} holds. Its last equation combined with Prop.~\ref{p.ZADE} and \ref{ss.STRUG} implies $\til_p,\tir_p \in \Str(J)$ with $\til_p^\sharp = \tir_p$, $\tir_p^\sharp = \til_p$. Since $\til_p^\sharp(1_J) = p = \tir_p^\sharp(1_J)$, therefore, we deduce from Thm.~\ref{t.STRINV.JOAL}~(d) that $\til_p$, $\tir_p$ are indeed isomorphisms of abstract Jordan algebras from $J^{(p)}$ to $J$. It remains to prove that they are, in fact, homomorphisms of cubic ones, i.e., that they preserve norms. Putting $J := J(A,\mu)$, we compute
\begin{align*}
N_J \circ \til_p(x) =\,\,&N_J\big(px_0 + (x_1p)j_1 + (p^\sharp x_2p^{-1})j_2\big) \\
=\,\,&N_A(p)N_A(x_0) + \mu N_A(p)N_A(x_1) + \mu^2N_A(p)N_A(x_2) \\
\,\,&- \mu T_A\big((px_0)(x_1p)(p^\sharp x_2p^{-1}\big),
\end{align*}
where the last summand, up to the factor $-\mu$, agrees with
\begin{align*}
N_A(p)T_A\Big(\big(p(x_0x_1)p\big)(p^{-1}x_2p^{-1})\Big) =\,\,&N_A(p)T_A\Big(p\big((x_0x_1)(x_2p^{-1})\big)\Big) \\
=\,\,&N_A(p)T_A\big((x_0x_1)(x_2p^{-1})p\big) \\
=\,\,&N_A(p)T_A(x_0x_1x_2).
\end{align*}
Hence the relation $N_J \circ \til_p = N_{J^{(p)}}$ holds strictly, so $\til_p\:J^{(p)} \to J$ is indeed a homomorphism of cubic Jordan algebras. The map $\tir_p$ can be treated similarly. The rest is clear.
\end{sol}

\begin{sol}{pr.SPRIFO} \label{sol.SPRIFO}  (a) We begin by showing that \eqref{ETEPER} makes $E^\perp$ a left $E$-module. Since $1\pt u = u$, we must show
\begin{align}
\label{EXSHAY} x \times (y \times u) = -(xy) \times u
\end{align}
for all $x,y \in E$, $u \in E^\perp$. After a faithfully flat base change that splits $E$ (Cor.~\ref{c.CHARF}), and consulting Prop.~\ref{p.CHAFF}, we may assume that $E$ itself is split, so there is an elementary frame $\Omega = (e_1,e_2,e_3)$ in $J$ satisfying $E = \sum ke_i$. Let $J = \sum(ke_i + J_{jl})$ be the Peirce decomposition of $J$ relative to $\Omega$. Then $E^\perp = \sum J_{jl}$ by (\ref{p.DECNOAD}.\ref{EXPBT}), and we have $x = \sum\xi_ie_i$, $y = \sum\eta_ie_i$, $u = \sum u_{jl}$, for some $\xi_i,\eta_i \in k$, $u_{jl} \in J_{jl}$, $1 \leq i \leq 3$. Now (\ref{p.DECNOAD}.\ref{EXPBAT}) implies
\[
x \times (y \times u) = -x \times \sum\eta_iu_{jl} = \sum\xi_i\eta_iu_{jl} = (xy) \times u_{jl}, 
\]
and this is \eqref{EXSHAY}.

(b) The property of $q_E$ to be a quadratic form over $E$ can be proved in a similar manner as in (a), but we prefer an approach that is ``rational'' by avoiding scalar extensions once (a) is taken for granted. Note first that \eqref{TEBZED} is a viable definition of $q_E$ as a map since $T_E$ is a regular symmetric bilinear form on $E$. Actually, $q_E$ is nothing else than the $k$-quadratic map belonging to the complemented cubic Jordan subalgebra $(E,E^\perp)$ of $J$ via (\ref{ss.SPLAD}.\ref{QUH}) and denoted there by $Q$. In particular, by (\ref{ss.IDQH.fig}.\ref{QXV}), $q_E$ is homogeneous of degree $2$ over $E$: $q_E(x\pt u) = x^2q_E(u)$ for all $x \in E$, $u \in E^\perp$. Now let $x \in E$, $u,v \in E^\perp$ and $z \in E$. Then the equations
\begin{align*}
T_E\big(z,q_E(x\pt u,v)\big) =\,\,&T_J\big(z,(x \times u) \times v\big) = T_J\big(x \times (z \times v), u\big) \\
=\,\,&-T_J\big((xz) \times v,u\big) = -T_J(xz,u \times v) = T_E\big(xz,q_E(u,v)\big) \\
=\,\,&T_E\big(z,xq_E(u,v)\Big) 
\end{align*}
show that the $k$-bilinearization of $q_E$ is, in fact, $E$-bilinear, as desired.

(c) The first part follows from the decomposition $J = E \oplus E^\perp$. In \eqref{EARPER} we first note that $q_{E_R}$ is nothing else than the $R$-quadratic extension of $q_E$ viewed as a quadratic map over $k$, so adopting the notational conventions of \ref{ss.RESC}, the top row of \eqref{EARPER} becomes
\[
\big(_k(q_E)\big)_R\:\big(_k(E^\perp)\big)_R \longrightarrow (_kE)_R = E \otimes R,
\]
and it makes sense to regard the vertical arrow on the left of \eqref{EARPER} as the identification (\ref{ss.RESC}.\ref{EMAR}). On the other hand, the bottom row of \eqref{EARPER} should be written, more accurately, as
\begin{align}
\label{QUERE} (q_E)_{R_E}\:(E^\perp)_{R_E} \longrightarrow E_{R_E} = E \otimes_E R_E = R \otimes E. 
\end{align}
But $(_kE)_R = E_{R_E}$, again under the identification (\ref{ss.RESC}.\ref{EMAR}), and combining \eqref{QUERE} with (\ref{ss.RESC}.\ref{EXTER}), we obtain $r \otimes x = 1_E \otimes_E (r \otimes x) = x \otimes r$ for $x \in E$, $r \in R$. Thus the commutativity of \eqref{EARPER} is equivalent to $(_k(q_E))_R = (q_E)_{R_E}$, which follows from (\ref{ss.RESC}.\ref{QURK}). 

(d)  If $R \in \kalg$ and $S \in \Ralg$ are both faithfully flat, then so is $S \in \kalg$ (\ref{ss.FLAFALG}~(i)). The same conclusion holds for the property of being finitely presented in place of being faithfully flat (Exc.~\ref{pr.TRAFIPR}~(a)). By Cor.~\ref{c.CHARF}, therefore, we may assume that $E$ is split. Accordingly, let $\Omega = (e_1,e_2,e_3)$ be an elementary frame in $J$, with Peirce decomposition $J = \sum(ke_i + J_{jl})$, such that $E = \sum ke_i$. By Prop.~\ref{p.SPLISIX}, we may assume that $J$ is regular of rank $n \ge 9$. Changing scalars to any algebraically close field $K \in \kalg$ makes $J_K$ split of dimension $9$, $15$, or $27$, and Cor.~\ref{c.CONFRA} shows that the off-diagonal Peirce spces of $J_K$ relative to $\Omega_K$ all have the same dimension $2$, $4$, or $8$, respectively. From (\ref{p.DECNOAD}.\ref{EXPQTB}) combined with Prop.~\ref{p.COMPOR}~(c) we therefore deduce that the restriction of $S_J$ to $J_{jl}$ is a regular quadratic form of even rank. After an appropriate \'etale cover of $k$, Exc.~\ref{pr.FFSPLIQ}~(b) allows us to assume that this restriction is split hyperbolic. In particular $S(u_{jl}) = -1$ for some $u_{jl} \in J_{jl}$, and we have found a strong co-ordinate system in $J$ extending $\Omega$. 
Combining the Jacobson co-ordinatization theorem~\ref{t.JACO} with Exc.~\ref{pr.CRITSEP}, we therefore obtain an identification $J = \Her_3(C)$, for some composition algebra $C$ over $k$, such that $E = \sum ke_{ii}$ is the diagonal of $J$. After yet another fppf base change, we may assume that $C = C_{0r}$ is split of rank $r = \frac{1}{3}(n - 3)$, proving the first part of (d). Note that $E_R = R_E$ (by (b)) is faithfully flat over $E$ since faithful flatness is stable under base change (\ref{ss.FFBAS}). Combining \eqref{EARPER} with (\ref{ss.THERCU}.\ref{MABIT}), (\ref{ss.THERCU}.\ref{MABADJ}), we therefore conclude that
\[
(E^\perp)_{R_E} = E_R^\perp = \sum C[jl]
\]
is a free $E_R$-module of rank $r$. Hence $E^\perp$ by \cite[I, Lemme~3.6~(c)]{MR0417149} is a finitely generated projective $E$-module of rank $r$. It remains to show that $q_E$ is non-singular, and even regular unless $n = 6$ and $2 \notin k^\times$. Since a quadratic form on a projective module that becomes regular (resp. non-singular) after a faithfully flat base change must have been so all along (which follows immediately from the definitions in the regular case and from Exc.~\ref{pr.FFSEP} in the non-singular one), it suffices to prove the assertion for $q_{E_R}$ rather than $q_E$. But (\ref{ss.THERCU}.\ref{MADJ}) shows $q_{E_R} \cong (n_C)_{E_R}$, and $n_C$ by Cor.~\ref{c.COMALGLOC} is regular unless $r = 1$ and $\frac{1}{2} \notin k$, in which case it is at least non-singular. The assertion follows.

(e) Since $E^\perp$ is an $E$-module under the action \eqref{ETEPER}, the stability condition (iii) of \ref{ss.KUMELT} holds automatically for any $l \in E^\perp$. Thus $l$ is a Kummer element relative to $E$ if and only if $l$ is invertible in $J$ and both $l$ and $l^\sharp$ belong to $E^\perp$; the latter condition, in turn, is equivalent to $q_E(l) = 0$, proving the first part of (e). In this case, we clearly have $A_l(J,E) = E$, and Cor.~\ref{c.UPFITI} yields the indicated embedding from $J(E,\lambda)$ to $J$. The final assertion will follow once we have shown that $l$ is unimodular as an element of $E^\perp$ \emph{as an $E$-module}. But this is clear since $q_E(l,l^\sharp) = N_J(l)1_E$ by (\ref{ss.IDQH.fig}.\ref{QXVH}).
\end{sol}

\begin{sol}{pr.COHACE} \label{sol.COHACE}  (a) We begin with the reduction to the case $\mu = 1$, $p = 1_C$. For any $\alpha \in k^\times$, $p \in C^\times$, we have $\hp := (\alpha,p) \in A^\times$ and
\begin{align}
\label{AALPE} A^{\hp} = k \times C^p, \quad N_A(\hp) = \alpha n_C(p).
\end{align}
Now let $\alpha := \mu$ and $p := 1_C$. Then \eqref{AALPE} implies $A^{\hp} = A$, $N_A(\hp) = \mu$ and applying Exc.~\ref{pr.CHAMU} yields $J(A,\mu) \cong J(A,1)$. Next assume $\alpha = 1$ in \eqref{AALPE}. Assuming the case $\mu = 1$, $p = 1_C$ has been settled, we conclude
\[
\Her_3(C^p,\Gamma_0) \cong J(A^{\hp},1) \cong J(A,N_A(\hp)) \cong J(A,1) \cong \Her_3(C,\Gamma_0).
\]
From now on, we may therefore assume $\mu = 1$ and $p = 1_C$. 

(i) We abbreviate $t := t_C$, $n := n_C$, $T := T_J$, $S := S_J$, $N := N_J$ and write $C^0$ for the module of trace-zero elements in $C$. Moreover, we put $e_1 := (1,0) \in A$, have $A = ke_1 \oplus C$ as a direct sum of submodules, and recall from Exc.~\ref{pr.CUADUN}~(b) the relations
\begin{align}
\label{CUUN} 1 := 1_A =\,\,&e_1 + 1_C, \\
\label{CUNO} N(\alpha e_1 + u) =\,\,&\alpha n(u), \\
\label{CUAD} (\alpha e_1 + u)^\sharp =\,\,&n(u)e_1 + \alpha\bar u,
\\
\label{CUBIAD} (\alpha e_1 + u) \times (\beta e_1 + v)
=\,\,&n(u,v)e_1 +
\alpha\bar v + \beta\bar u, \\
\label{CUBILT} T(\alpha e_1 + u,\beta e_1 + v) =\,\,&\alpha\beta +
t(uv), \\
\label{CULIT} T(\alpha e_1 + u) =\,\,&\alpha + t(u), \\
\label{CUQUAT} S(\alpha e_1 + u) =\,\,&\alpha t(u) + n(u), \\
\label{CUQUAB} S(\alpha e_1 + u,\beta e_1 + v) =\,\,&\alpha t(v) + \beta t(u) + n(u,v)
\end{align}
for all $\alpha,\beta \in $, $u,v \in C$. From \eqref{CUAD} and \eqref{CULIT} we conclude that $e_1$ is an elementary idempotent in $J$. It remains to establish the corresponding Peirce rules. While \eqref{PEEONETWO} follows immediately from
(\ref{p.PELE}.\ref{PETWO}), the remaining assertions require a
proof. Let
\begin{align}
\label{EXDEC} x = x_0 + x_1j_1 + x_2j_2, \quad x_i = \alpha_ie_1 +
v_i \in A, \quad \alpha_i \in k, \quad v_i \in C, \quad 0 \leq i
\leq 2.
\end{align}
Then (\ref{t.EXTITS}.\ref{TRFITI}) and \eqref{CULIT} imply
\begin{align}
\label{EXTR} T(x) = T(x_0) = \alpha_0 + t(v_0).
\end{align}
On the other hand, (\ref{t.EXTITS}.\ref{BADFITI}) yields $e_1 \times x
= e_1 \times x_0 - (e_1x_1)j_1 - (x_2e_1)j_2$, and applying
\eqref{CUBIAD} we conclude
\begin{align}
\label{EXBIAD} e_1 \times x = \bar v_0 - (\alpha_1e_1)j_1 -
(\alpha_2e_1)j_2.
\end{align}
Now \eqref{PEEONEONE} follows by combining
(\ref{p.PELE}.\ref{PEONE}) with \eqref{EXTR} and \eqref{EXBIAD}.
Observing
\begin{align}
\label{COEONE} e_0 := 1 - e_1 = 1_C,
\end{align}
we obtain $T(x)e_0 - x = (\alpha_0 + t(v_0))1_C - x_0 - x_1j_1 -
x_2j_2 = \alpha_0(1_C - e_1) + \bar v_0 - x_1j_1 - x_2j_2$, and
comparing this with \eqref{EXBIAD}, we deduce \eqref{PEEONEZERO}
from (\ref{p.PELE}.\ref{PEZERO}). 

(ii) We clearly have $\sum e_i = 1$. Suppose we know that $e_2$ is an elementary idempotent of $J$. Then \eqref{PEEONEZERO} shows that it is orthogonal to $e_1$, and the assertion follows from Prop.~\ref{p.COMPOR}. Thus we only need to show that $e_2$ is an elementary idempotent.  From
\eqref{CULIT}, (\ref{t.EXTITS}.\ref{TRFITI}) and \eqref{CUAD} we deduce $T(e_2)
= t(u_0) = 1$, while (\ref{t.EXTITS}.\ref{BADFITI}) and \eqref{CUBIAD} yield $e_2^\sharp = (u_0^\sharp - n(u_0)e_1) +
(n(u_0)^2e_1^\sharp - n(u_0)u_0e_1)j_1 + (e_1^\sharp - n(u_0)e_1u_0)j_2 =
n(u_0)e_1 - n(u_0)e_1 = 0$, so $e_2$ is indeed an elementary
idempotent. 

(iii) Decomposing an arbitrary element $x \in J$ as in
\eqref{EXDEC}), we claim ,
\begin{align}
\label{TWOBIAD} e_2 \times x =\,\,&\Big(\big(n(u_0,v_0) -
n(u_0)\alpha_1 - \alpha_2\big)e_1 + \alpha_0\bar u_0\Big) + \\
\,\,&\Big(-\alpha_0e_1 + \big(n(u_0)\bar v_2 - u_0v_1\big)\Big)j_1
+ \Big(-n(u_0)\alpha_0e_1 + (\bar v_1 - v_2u_0)\Big)j_2, \notag \\
\label{THREEBIAD} e_3 \times x =\,\,&\Big(\big(n(\bar u_0,v_0) +
n(u_0)\alpha_1 + \alpha_2\big)e_1 + \alpha_0u_0\Big) + \\
\,\,&\Big(\alpha_0e_1 - \big(n(u_0)\bar v_2 + \bar
u_0v_1\big)\Big)j_1 + \Big(n(u_0)\alpha_0e_1 - (\bar v_1 + v_2\bar
u_0)\Big)j_2. \notag
\end{align}
To establish \eqref{TWOBIAD}, we expand the left-hand side by
(\ref{t.EXTITS}.\ref{BADFITI}) and apply \eqref{CUBIAD}
to conclude 
\begin{align*}
e_2 \times x =\,\,&\Big(u_0 + e_1j_1 + \big(n(u_0)e_1\big)j_2\Big)
\times \Big(x_0 + x_1j_1 + x_2j_2\Big) \\
=\,\,&\big(u_0 \times x_0 - e_1x_2 - n(u_0)x_1e_1\big) +
\big(n(u_0)e_1 \times x_2 - u_0x_1 - x_0e_1\big)j_1 + \\
\,\,&\big(e_1 \times x_1 - n(u_0)e_1x_0 - x_2u_0\big)j_2 \\
=\,\,&\big(n(u_0,v_0)e_1 + \alpha_0\bar u_0 - \alpha_2e_1 -
n(u_0)\alpha_1e_1\big) + \\
\,\,&\big(n(u_0)\bar v_2 - u_0v_1 - \alpha_0e_1\big)j_1 + (\bar v_1
- n(u_0)\alpha_0e_1 - v_2u_0\big)j_2.
\end{align*}
Hence \eqref{TWOBIAD} follows. Similarly,
\begin{align*}
e_3 \times x =\,\,&\Big(\bar u_0 - e_1j_1 -
\big(n(u_0)e_1\big)j_2\Big) \times \Big(x_0 + x_1j_1 +
x_2j_2\Big) \\
=\,\,&\big(\bar u_0 \times x_0 + e_1x_2 + n(u_0)x_1e_1\big) +
\big( -n(u_0)e_1 \times x_2 - \bar u_0x_1 + x_0e_1\big)j_1 + \\
\,\,&\big(-e_1 \times x_1 + n(u_0)e_1x_0 - x_2\bar u_0\big)j_2 \\
=\,\,&\big(n(\bar u_0,v_0)e_1 + \alpha_0u_0 + \alpha_2e_1 +
n(u_0)\alpha_1e_1\big) + \\
\,\,&\big(-n(u_0)\bar v_2 - \bar u_0v_1 + \alpha_0e_1\big)j_1 +
\big(-\bar v_1 + n(u_0)\alpha_0e_1 - v_2\bar u_0\big)j_2
\end{align*}
gives \eqref{THREEBIAD}.

We now prove \eqref{ETWOTH} by applying (\ref{ss.COELF}.\ref{PECOEL}), \eqref{TWOBIAD}, \eqref{THREEBIAD} and \eqref{EXTR} to obtain the following chain of equivalent statements. 
\begin{align*}
x \in J_{23} \Longleftrightarrow\,\,&T(x) = 0,\quad e_2 \times x = e_3 \times x = 0 \\
\Longleftrightarrow\,\,&\alpha_0 + t(v_0) = 0, \quad n(u_0,v_0) = n(u_0)\alpha_1 + \alpha_2, \quad \alpha_0 = 0, \\
\,\,&n(u_0)\bar v_2 = u_0v_1, \quad \bar v_1 = v_2u_0, \\
\,\,&n(\bar u_0,v_0) + n(u_0)\alpha_1 + \alpha_2 = 0, \quad n(u_0)\bar v_2 + \bar u_0v_1 = 0 = \bar v_1 + v_2\bar u_0; 
\end{align*}
if this holds, then $0 = \bar v_1 - \bar v_1 = v_2u_0 + v_2\bar u_0 = v_2$, hence $v_1 = 0$ as well, and \eqref{ETWOTH} is proved. Similarly, also using \eqref{EXBIAD},
\begin{align*}
x \in J_{31} \Longleftrightarrow\,\,&T(x) = 0, \quad e_1 \times x = 0 = e_3 \times x \\
\Longleftrightarrow\,\,& \alpha_0 + t(v_0) = 0, \quad v_0 = 0, \quad \alpha_1 = \alpha_2 = 0, \\ 
\,\,&\alpha_0 = 0, \quad n(u_0)\bar v_2 + \bar u_0v_1 = 0 = \bar v_1 + v_2\bar u_0 \\
\Longleftrightarrow\,\,&\alpha_0 = \alpha_1 = \alpha_2 = 0, \quad v_0 = 0, \quad v_1 = -u_0\bar v_2,
\end{align*} 
and this amounts to \eqref{ETHONE}. Finally,
\begin{align*}
x \in J_{12} \Longleftrightarrow\,\,&T(x) = 0, \quad e_1 \times x = 0 = e_2 \times x \\
\Longleftrightarrow\,\,&v_0 = 0, \quad \alpha_1 = \alpha_2 = 0, \\
\,\,&\alpha_0 = 0, \quad n(u_0)\bar v_2 - u_0v_1 = 0 = \bar v_1 - v_2u_0 \\
\Longleftrightarrow\,\,&v_0 = 0, \quad \alpha_0 = \alpha_1 = \alpha_2 = 0, \quad v_1 = \bar u_0\bar v_2,
\end{align*}
which proves \eqref{EONETW}. 

(iv) From \eqref{EONETW} we deduce that $u_{31}$ belongs to $J_{31}$. Since $2u_0 - 1_C$ has trace zero and
\[
u_{23} = (2u_0 - 1_C) +2e_1j_1 + \big(n_C(u_0,2u_0 - 1_C) - 2n_C(u_0)\big)e_1j_2,
\]
$u_{23}$ by \eqref{ETWOTH} belongs to $J_{23}$. It remains to compute $S(u_{jl})$ for $i = 1,2$. We do so by applying (\ref{t.EXTITS}.\ref{QRFITI}), \eqref{CUQUAT} and \eqref{CUBILT} to conclude
\begin{align*}
S(u_{23}) =\,\,&S(2u_0 - 1_C) - T\Big(2e_1,\big(2n_C(u_0) - 1\big)e_1\Big) = n(2u_0 - 1_C) - 2\big(2n(u_0) - 1\big) \\
=\,\,&4n(u_0) - 2 + 1 - 4n(u_0) + 2 = 1, \\
S(u_{31}) =\,\,&T(u_0) = t(u_0) = 1,
\end{align*} 
and the proof of (iv) is complete. 

(v) In the terminology of Prop.~\ref{p.CONALT} we have
\begin{align}
\label{OMEGON} \omega = 1, \quad \gamma_1 = \gamma_2 = -1, \quad \gamma_3 = 1.
\end{align}
We begin by verifying \eqref{UVEBAR}. By Prop.~\ref{p.ZADE}, we may assume that $u$ is invertible. Then the left Moufang identity yields
\begin{align*}
u\big((uv - \bar u\bar v)(uw)\big) =\,\,&\big(u(uv - \bar u\bar
v)u\big)w = \big(u(uvu)\big)w - n(u)(\bar vu)w \\
=\,\,&n(u,\bar v)u^2w - n(u)(u\bar v + \bar vu)w \\
=\,\,&n(u,\bar v)u^2w - n(u)\big(t(u)\bar v + t(v)u -
n(u,\bar v)1_C\big)w \\
=\,\,&u\big(n(u,\bar v)uw - t(u)\bar u(\bar vw) - n(u)t(v)w
+ n(u,\bar v)\bar uw\big) \\
=\,\,&u\Big(\big(t(u)n(u,\bar v) - n(u)t(v)\big)w - t(u)\bar u(\bar
vw)\Big).
\end{align*}
Cancelling $u$, the assertion follows.

Now let
\begin{align}
\label{EXBARU} x = (\bar u_0\bar v)j_1 + vj_2, \quad y = (\bar u_0\bar w)j_1 + wj_2, &&(v,w \in C)
\end{align}
be arbitrary elements of $J_{12} = C_{J,\mfS}$. Then (\ref{t.EXTITS}.\ref{BADFITI}) yields
\begin{align*}
x \times u_{23} =\,\,&\big((\bar u_0\bar v)j_1 + vj_2\big) \times \Big((2u_0 - 1_C) + 2e_1j_1 + \big(2n(u_0) - 1\big)e_1j_2\Big) \\
=\,\,&\Big(-\big(2n(u_0) - 1\big)(\bar u_0\bar v)e_1 - 2e_1v\Big) + \Big(\big(2n(u_0) - 1\big)v \times e_1 - (2u_0 - 1_C)\bar u_0\bar v\Big)j_1 \\
\,\,&+ \big(2(\bar u_0\bar v) \times e_1 - v(2u_0 - 1_C)\big)j_2 \\
=\,\,&\Big(\big(2n(u_0) - 1\big)\bar v - 2n(u_0)\bar v + \bar u_0\bar v\Big)j_1 + (2vu_0 - 2vu_0 + v)j_2 \\
=\,\,&(-\bar v+ \bar u_0\bar v)j_1 + vj_2 = (-u_0\bar v)j_1 + vj_2,
\end{align*}
hence
\begin{align}
\label{XUTWOTHREE} x \times u_{23} = -(u_0\bar v)j_1 + vj_2.
\end{align}
Similarly,
\begin{align*}
u_{31} \times y =\,\,&\big(-u_0j_1 + 1_Cj_2\big) \times \big((\bar
u_0\bar w)j_1 + wj_2\big) \\
=\,\,&(u_0w - \bar u_0\bar w) + (1_C \times w)j_1 - \big(u_0 \times
(\bar u_0\bar w)\big)j_2,
\end{align*}
hence
\begin{align}
\label{UONETHREEY} u_{31} \times y = \big(u_0w - \bar u_0\bar w\big)
+ \big(t(w)e_1\big)j_1 - \big(n(u_0,\bar u_0\bar w)e_1\big)j_2.
\end{align}
As a side remark we note $u_0w - \bar u_0\bar w \in C^0$ and $n(u_0,u_0w - \bar u_o\bar w) - n(u_0)t(w)
= n(\bar u_0u_0,w) - n(u_0,\bar u_0\bar w) - n(u_0)t(w) =
-n(u_0,\bar u_0\bar w)$, so
\[
u_{31} \times y = \big(u_0w - \bar u_0\bar w\big) +
\big(t(w)e_1\big)j_1 + \big([n(u_0,u_0w - \bar u_0\bar w) -
n(u_0)t(w)]e_1\big)j_2
\]
by \eqref{ETWOTH} does indeed belong to $J_{23}$,
as it should. Using \eqref{XUTWOTHREE}, \eqref{UONETHREEY} and
(\ref{t.EXTITS}.\ref{BADFITI}), \eqref{CUBIAD}, we can now
compute
\begin{align*}
(x \times u_{23}) \times (u_{31} \times y) =\,\,&\Big(-(u_0\bar
v)j_1 + vj_2\Big) \times \\
\,\,&\Big(\big(u_0w - \bar u_0\bar w \big) + \big(t(w)e_1\big)j_1 -
\big(n(u_0,\bar u_0\bar w)e_1\big)j_2\Big) \\
=\,\,&\big(-n(u_0,\bar u_0\bar w)(u_0\bar v)e_1 - t(w)e_1v\big) + \\
\,\,&\big(-n(u_0,\bar u_0\bar w)v \times e_1 + (u_0w - \bar u_0\bar
w)(u_0\bar v)\big)j_1 + \\
\,\,&\big(-t(w)(u_0\bar v) \times e_1 - v(u_0w - \bar u_0\bar
w)\big)j_2 \\
=\,\,&\big((u_0w - \bar u_0\bar w)(u_0\bar v) - n(u_0,\bar u_0\bar
w)\bar v\big)j_1 + \\
\,\,&\big(-t(w)v\bar u_0 - v(u_0w) + v(\bar u_0\bar w) \big)j_2.
\end{align*}
Here we note
\begin{align*}
-t(w)v\bar u_0 - v(u_0w) + v(\bar u_0\bar w) = -v(\bar u_0w) -
v(\bar u_0\bar w) - v(u_0w) + v(\bar u_0\bar w) = -vw
\end{align*}
since $u_0$ has trace $1$ and apply \eqref{UVEBAR} to conclude
\begin{align*}
(u_0w - \bar u_0\bar w)(u_0\bar v) - n(u_0,\bar u_0\bar w)\bar v
=\,\,&\big(n(u_0,\bar w) - n(u_0)t(w)\big)\bar v - \bar u_0(\bar
w\bar v) - n(u_0^2,\bar w)\bar v \\
=\,\,&-\bar u_0\overline{vw} + \big(n(u_0,\bar w) - n(u_0)t(w) \\
\,\,&-n(u_0,\bar w) + n(u_0)t(w)\big)\bar v \\
=\,\,&-\bar u_0\overline{vw}.
\end{align*}
Hence
\begin{align}
\label{COORPROD} (x \times u_{23}) \times (u_{31} \times y)
=\,\,&(-\bar u_0\overline{vw})j_1 + (-vw)j_2.
\end{align}
We put $\pC:= C_{J,\mfS}$ and $n^\prime := n_{C_{J,\mfS}}$. Since $\omega = 1$ by \eqref{OMEGON}, we combine
\eqref{COORPROD} with (\ref{p.CONALT}.\ref{CONPRO}) to conclude
that
\[
\varphi\:C \longrightarrow C^\prime, \quad v \longmapsto (- \bar u_0\bar
v)j_1 - vj_2
\]
is an isomorphism of $k$-algebras; in particular, it is unital. Moreover, by
(\ref{p.CONALT}.\ref{CONNO}),
$(n^\prime \circ \varphi)(v) = -S(-(\bar u_0\bar v)j_1 - vj_2) = n(v)$,
so $\varphi\:(C,n) \to (C^\prime,n^\prime)$ is actually an
isomorphism of \emph{conic} $k$-algebras.

(b) $C$ is free of rank $8$ as a $k$-module. By Exc.~\ref{pr.DIAFIX}, we may assume $\Gamma = \diag(\gamma_1,\gamma_2,\gamma_3)$, $\prod\gamma_i = 1$. Since $n_C$ is universal, we find $p,q \in C^\times$ satisfying $n_C(p) = \gamma_3^{-1}$, $n_C(q) = -\gamma_2^{-1}$, hence $n_C(pq) = -\gamma_1$. By Exc.~\ref{pr.ISCOPA}, therefore, we may assume $\Gamma = \Gamma_0$.  On the other hand, the Alsaody-Gille theorem \ref{t.ALSGIL} yields some $p \in C^\times$ making $C^p$ split. By (a), therefore, $\Her_3(C,\Gamma) \cong \Her_3(C^p,\Gamma_0)$ is split by Prop.~\ref{p.FRSPCOM}.
\end{sol}


\solnsec{Section~\ref{s.ISINV}}

\begin{sol}{pr.CAPI} \label{sol.CAPI} In order to get the terminology straight, let us denote by $\bfpainv$ the category of pointed alternative $k$-algebras with involution as defined in Exc.~\ref{pr.CAPI}. By contrast, let us denote by $\bfaist$ the category of alternative $k$-algebras with isotopy involution as defined in \ref{ss.HOIDIS}. 

We begin by defining a functor $\Phi\:\bfpainv \to \bfaist$. Let $((B,\tau),q)$ be an object of $\bfpainv$. Then Cor.~\ref{c.ENISINV} shows that
\[
\Phi\Big(\big((B,\tau),q\big)\Big) := (B^q,\tau^q,q)
\]
is an object of $\bfaist$. If $\vph\:((B,\tau),q) \to ((B^\prime,\tau^\prime),q^\prime)$ is a morphism in $\bfpainv$, then it is straightforward to check that
\[
\vph\:(B^q,\tau^q,q) \longrightarrow (B^{\prime q^\prime},\tau^{\prime q^\prime},q^\prime)
\]
is a morphism in $\bfaist$. Hence we obtain a functor $\Phi$ of the desired kind.

Next we define a functor $\Psi\:\bfaist \to \bfpainv$ in the opposite direction. Let $(B,\tau,q)$ be an object of $\bfaist$, i.e., an alternative $k$-algebra with isotopy involution. Then Prop.~\ref{p.TWISINV} and Lemma~\ref{l.TAUPO} show with $p := q^{-1}$ that
\[
\Psi\big((B,\tau,q)\big) := \big((B^p,\tau^p),q\big)
\]
is an object of $\bfpainv$. If $\psi\:(B,\tau,q) \to (B^\prime,\tau^\prime,q^\prime)$ is a morphism in $\bfaist$ and $p := q^{-1}$, $p^\prime := q^{\prime -1}$, then $\psi(p) = p^\prime$, and it is straightforward to check that
\[
\psi\:\big((B^p,\tau^p),p\big) \longrightarrow \big((B^{\prime p^\prime},\tau^{\prime p^\prime}),q^\prime\big)
\]
is a morphism in $\bfpainv$. Hence we obtain a functor $\Psi$ of the desired kind.

Now we simply apply (\ref{p.TWISINV}.\ref{TWISIT}) to conclude that $\Psi \circ \Phi$ is the identity on $\bfpainv$ and $\Phi \circ \Psi$ is the identity on $\bfaist$. Summing up
\[
\Phi\:\bfpainv \overset{\sim} \longrightarrow \bfaist
\]
is an isomorphism of categories with inverse $\Psi$.
\end{sol}

\begin{sol}{pr.SWISINV} \label{sol.SWISINV}
 The map $\epsilon_A$ has clearly period
two and fixes $p$. For the first part of the problem, it therefore
suffices to show (\ref{ss.COISINV}.\ref{ISINV}) for $\epsilon_A$ in
place of $\tau$. To this end, let $x,x^\prime,y,y^\prime \in A$.
Then
\begin{align*}
\Big(\epsilon_A\big((y,y^\prime)\big)p^{-1}\Big)\Big(p\epsilon_A\big((x,
x^\prime)\big)\Big) =\,\,&\big((y^\prime,y)(q,q)\big)\big((q^{-1},q^{-1})(x^\prime, x)\big) \\
=\,\,&\big((\py q^{-1})(qq),qy\big)\big((q^{-1}q^{-1})(q\px),xq^{-1}\big) \\
=\,\,&(y^\prime q,qy)( q^{-1}x^\prime,xq^{-1}) \\
=\,\,&\Big(\big((y^\prime q)q^{-1}\big)\big(q(q^{-1}x^\prime)\big),(xq^{-1})(qy)\Big) = \big(y^\prime x^\prime,(xq^{-1})(qy)\big) \\
=\,\,&\epsilon_A\Big(\big((xq^{-1})(qy)\big),y^\prime x^\prime\big)\Big) =
\epsilon_A\big((x,x^\prime)(y,y^\prime)\big),
\end{align*}
as desired. Turning to the second part of the problem, we first note
\[
H(A^q \times A^{\op},\epsilon_A) = \{(x,x) \mid x \in A\}.
\]
On the other hand, given $x,x^\prime \in A$, we have
\begin{align*}
\Big(\epsilon_A\big((x,x^\prime)\big)p\Big)(x,x^\prime) =\,\,&\big((x^\prime,x)(q^{-1},q^{-1})\big)(x,x^\prime) \\
=\,\,&\big((x^\prime q^{-1})(qq^{-1}),q^{-1}x)(x,x^\prime)\big) \\
=\,\,&(x^\prime q^{-1},q^{-1}x)(x,x^\prime)
= \Big(\big((x^\prime q^{-1})q^{-1}\big)(qx),x^\prime(q^{-1}x)\Big)
\\
=\,\,&\big((x^\prime q^{-2})(qx),x^\prime(q^{-1}x)\big).
\end{align*}
Combining, and using (\ref{ss.UNTISOTALT}.\ref{ITUNTISOTALT}), we
obtain the following chain of equivalent conditions:
\begin{align*}
\forall x,x^\prime \in A:\;\Big(\epsilon_A\big((x,x^\prime)\big)p\Big)(x,x^\prime) \in\,\,& H(A^q \times A^{\op},\epsilon_A)\\
\Longleftrightarrow\,\,&\forall x,x^\prime \in
A:\; (x^\prime q^{-2})(qx) = x^\prime(q^{-1}x) \\
\Longleftrightarrow\,\,& \forall x,x^\prime \in A:\;
(x^\prime q^{-2})(q^2x) = x^\prime x \\
\Longleftrightarrow\,\,& A = A^{q^2} \\
\Longleftrightarrow\,\,& q^2 \in \Nuc(A).
\end{align*}
This solves the problem.
\end{sol}


\solnsec{Section~\ref{s.SETINV}}

\begin{sol}{pr.DESELI} \label{sol.DESELI} (a) $M_0$ is clearly a $k$-submodule of $M$ and (\ref{ss.SCEM}.\ref{REX}) implies $\Phi(x_0 \otimes a) = ax_0$ for all $x_0 \in M_0$ and $a \in K$. In particular, \eqref{EMZEK} commutes. It remains to show that $\Phi$ is bijective. To this end, we note for any prime ideal $\mfp \subseteq k$ that $k_\mfp$ is a flat $k$-algebra (Bourbaki \cite[II, \S2, Thm.~1]{MR0360549}), which implies $H(M_\mfp,\tau_\mfp) = H(M,\tau)_\mfp$ by Exc.~\ref{pr.FLASET}. Hence we may assume that $k$ is a local ring and apply Prop.~\ref{p.SMASU} to pick $\theta \in K$ as in the hint. The
$k$-linear map $M_0^2 \to M_0 \otimes K$, $(x_0,y_0) \mapsto x_0 \otimes
1_K + y_0 \otimes \theta$ is bijective and $\Phi(x_0 \otimes 1_K + y_0
\otimes \theta) = x_0 + \theta y_0$ for all $x_0,y_0 \in M_0$. Given $z \in
M$, $x_0,y_0 \in M_0$, it therefore suffices to show
\begin{align}
\label{TAULOC} z = x_0 + \theta y_0 \Longleftrightarrow
\begin{cases}
x_0 = \big(1 - 4n_K(\theta)\big)^{-1}\big((1 - 2\theta)\bar\theta
z + (1 - 2\bar\theta)\theta\tau(z)\big), \\
y_0 = \big(1 - 4n_K(\theta)\big)^{-1}\big((1 - 2\bar\theta)z + (1 -
2\theta)\tau(z)\big).
\end{cases}
\end{align}
If $x_0,y_0$ have the form indicated in \eqref{TAULOC}, then it is
straightforward to check that they remain fixed under $\tau$, hence
belong to $M_0$, and $t_K(\theta) = 1$ implies
\begin{align*}
x_0 + \theta y_0 =\,\,&\big(1 - 4n_K(\theta)\big)^{-1}\Big(\big((1 -
2\theta)\bar\theta + (1 - 2\bar\theta)\theta\big)z + \big((1 -
2\bar\theta)\theta + (1 - 2\theta)\theta\big)\tau(z)\Big) \\
=\,\,&\big(1 - 4n_K(\theta)\big)^{-1}\Big(\big(1 -
4n_K(\theta)\big)z - 2\big(\theta^2 - \theta +
n_K(\theta)1_K\big)\tau(z)\Big) = z.
\end{align*}
Conversely, suppose $z \in M$ has $z = x_0 + \theta y_0$, $x_0,y_0 \in M_0$.
Then $\tau(z) = x_0 + \bar\theta y_0 = (x_0 + y_0) - \theta y_0$, forcing
\begin{align}
\label{TRAZET} 2x_0 + y_0 = z + \tau(z).
\end{align}
Applying \eqref{TRAZET} to $\bar\theta z = \bar\theta x_0 +
n_K(\theta)y_0 = (x_0 + n_K(\theta)y_0) - \theta x_0$ in place of $z$, we
conclude
\begin{align}
\label{TRATHE} x_0 + 2n_K(\theta)y_0 = 2\big(x_0 + n_K(\theta)y_0\big) - x_0 =
\bar\theta z + \tau(\bar\theta z).
\end{align}
Multiplying \eqref{TRATHE} by $2$ and subtracting the result from
\eqref{TRAZET}, we obtain $(1 - 4n_K(\theta))y_0 = (1 - 2\bar\theta)z
+ \tau((1 - 2\bar\theta)z)$. Similarly, multiplying \eqref{TRAZET}
by $2n_K(\theta)$ and subtracting the result from \eqref{TRATHE}, we
obtain$(1 - 4n_K(\theta))x_0 = (1 - 2\theta)\bar\theta z + \tau((1 -
2\theta)\bar\theta z)$, so $x_0,y_0$ have the form indicated in
\eqref{TAULOC}. This shows that $\Phi$ is bijective.

(b) We begin with the hint. Since $K$ is faithfully flat as a $k$-module, Prop.~\ref{p.CHAFF} shows that the natural map $M_0 \to M_{0K}$ is injective, hence gives an isomorphism $M_0 \overset{\sim} \to M_0 \otimes 1_K$. In order to prove $M_0 \otimes 1_K = H(M_{0K},\Eins_{M_0} \otimes \iota)$, we may again assume that $k$ is a local ring and pick $\theta$ as in (a). Again the map $M_0^2 \to M_0 \otimes K$, $(x_0,y_0) \mapsto x_0 \otimes 1_K + y_0 \otimes \theta$, is a
$k$-isomorphism and $(\Eins_{M_0} \otimes \iota)(x_0 \otimes 1_K + y_0
\otimes \theta) = x_0 \otimes 1_K + y_0 \otimes \bar\theta = (x_0 + y_0)
\otimes 1_K - y_0 \otimes \theta$ for all $x_0,y_0 \in M_0$. For $x_0
\otimes 1_K + y_0 \otimes \theta$ to remain fixed under $\Eins_{M_0}
\otimes \iota$ it is therefore necessary and sufficient that $x_0 + y_0
= x_0$ and $y_0 = -y_0$, so the assertion follows. In (b), everything is clear except for the final statement. Changing scalars from $k$ to $R$ in \eqref{EMZEK}, we obtaina commutative diagram
\[
\xymatrix{M_{0R} \otimes_RK_R \ar[rr]_{\Eins_{M_{0R}} \otimes_R \iota_{K_R}} \ar[d]_{\Phi_R}^{\cong} && M_{0R} \otimes_R K_R \ar[d]_{\cong}^{\Phi_R} \\
M_R \ar[rr]_{\tau_R} && M_R}
\]
of $R$-modules. Combined with the hint on the level of $R$ rather than $k$, this yields a chain of
isomorphisms
\[
\xymatrix{ M_{0R} \ar[r]_(0.2){- \otimes_R 1_{K_R}}^(0.2){\cong} &
M_{0R} \otimes_R 1_{K_R} = H(M_{0R} \otimes_R K_R,\Eins_{M_{0R}}
\otimes_R \iota_{K_R}) \ar[r]_(0.66){\Phi_R}^(0.65){\cong} &
H(M_R,\tau_R) = M_{R0}}
\]
sending $x_0 \otimes r \in M_{0R}$ ($x_0 \in M_0$, $r \in R$) after natural identifications\footnote{Given $k$-modules $M,N$ and $R \in \kalg$, we identify $(M \otimes N)_R = M_R \otimes_R N_R$ via $(x \otimes y) \otimes r = (x \otimes r) \otimes_R (y \otimes 1_R) = (x \otimes 1_R) \otimes_R (y \otimes r)$ and $(x \otimes r_1) \otimes_R (y \otimes r_2) = (x \otimes y) \otimes (r_1r_2)$ for $x \in M$, $y \in N$, $r,r_1,r_2 \in R$.} to
\begin{align*}
\Phi_R\big((x_0 \otimes r) \otimes_R1_{K_R}\big) =\,\,&\Phi_R\big((x_0
\otimes r) \otimes_R (1_K \otimes 1_R)\big) = \Phi_R\big((x_0
\otimes 1_K) \otimes r\big) \\
=\,\,&\Phi(x_0 \otimes 1_K) \otimes r = x_0 \otimes r \in M_{R0}
\end{align*}
and showing that $\Phi_R \circ(-- \otimes 1_{K_R})$ is indeed an isomorphism from $M_{0R}$ to $M_{R0}$.
\end{sol}

\begin{sol}{pr.MOVCOR} \label{sol.MOVCOR}
$\vrh\:K \to\, \rhK$ is $K$-linear and, after identifying $B = B \otimes_KK$ canonically, $\can := \can_{\msB} = \Eins_B \otimes_K \vrh$ is $\vrh$-semi-linear bijective. Thus $\rhtau\:\,\rhB \to\,\rhB$ exists as a $\rhiota$-semi-linear map, where $\rhiota = \iota_{\rhK}$ is the conjugation of $\rhK$. Write $\rhsharp$ for the adjoint of $\rhB$, which agrees with the $\rhK$-quadratic extensions of the adjoint of $B$. This implies $\rhtau \circ\,\rhsharp = \,\rhsharp \circ\,\rhtau$, and we deduce from \ref{ss.COINSY}~(b) that $\rhtau\:\,\rhB^{(+)} \to\,\rhB^{(+)}$ is a $\rhiota$-semi-linear automorphism of cubic Jordan algebras. Thus $\rhmsB$ is an involutorial $k$-system, and arguing as before, it follows that $(\vrh,\can)\:\msB \overset{\sim}\to\,\rhmsB$ is an isomorphism.
\end{sol}

\begin{sol}{pr.POCUAL} \label{sol.POCUAL}
 (a) For the first part, it suffices to show that $\msB(A,q)$ is an involutorial system, core splitness being obvious. Since $\vep_A$ is a $K/k$-involution of $B$ \emph{as a module} fixing $p$, it remains to show
\begin{align}
\label{EPSA} \varepsilon_A\big((u_1,u_1^\prime)(u_2,u_2^\prime)\big) =
\Big(\varepsilon_A\big((u_2,u_2^\prime)\big)p^{-1}\Big)\Big(p\varepsilon_A\big((u_1, u_1^\prime)\big)\Big)
\end{align}
for all $u_i,u_i^\prime \in B$, $i = 1,2$. Noting that the Jordan structures of $B$
and $B^{q\op}$ are the same by \ref{ss.UNISCU}, we manipulate the right-hand
side of \eqref{EPSA} to obtain
\begin{align*}
\Big(\varepsilon_A\big((u_2,u_2^\prime)\big)p^{-1}\Big)\Big(p\varepsilon_A&\big((u_1, u_1^\prime)\big)\Big) = \big((u_2^\prime,u_2)(q^{-1},q^{-1})\big)\big((q,q)(u_1^\prime,u_1)\big) \\
=\,\,&\Big(\big(u_2^\prime q^{-1},(q^{-1}q^{-1})(qu_2)\big)\Big)\Big(qu_1^\prime, (u_1q^{-1})(qq)\big)\Big) \\
=\,\,&(u_2^\prime q^{-1},q^{-1}u_2)(qu_1^\prime,u_1q) \\
=\,\,&\Big((u_2^\prime q^{-1})(qu_1^\prime),\big((u_1q)q^{-1}\big)\big(q(q^{-1}u_2)\big)\Big) \\
=\,\,&\big(u_2^\prime q^{-1})(qu_1^\prime),u_1u_2\big) \\
=\,\,&\varepsilon_A\Big(\big(u_1u_2,(u_2^\prime q^{-1})(qu_1^\prime)\big)\Big) \\
=\,\,&\varepsilon_A\big((u_1,u_1^\prime)(u_2,u_2^\prime)\big),
\end{align*}
as claimed. The second part of (a) follows by a straightforward computation.

(b) Since $(B,p)$ is a pointed cubic alternative algebra over $K$, its base change $(\rhB_+,{\rhp_+})$ is one over $k$, proving the first part of (b). As to the second, we have $\rhk_+ = \,\rhpk_+ = k$ as commutative rings, and since $\pvrh \circ \sigma = \vrh$ by definition, we conclude $a\alpha = \sigma(a)\alpha$ for all $a \in K$, $\alpha \in k$, i.e., $\Eins_k\:\,\rhk_+ \to \,\rhpk_+$ is $\sigma$-semi-linear. In the sense of \ref{ss.TESELI}, therefore, $\vph \otimes_\sigma \Eins_k\:\,\rhB_+ \to \,^{\pvrh}\hspace{-2pt}\pB_+$ exists as a $\sigma$-semi-linear map, which is obviously a unital homomorphism of alternative $k$-algebras and is easily checked to send $\rhp_+$ to $^{\pvrh}{\hspace{-2pt}}\pp_+$ and to  make the diagram
\begin{align*}
\xymatrix{\rhB_+ \ar[rr]_{\vph \otimes_\sigma \Eins_k} \ar[d]_{\rhsharp} && ^{\pvrh}\hspace{-2pt}\pB_+ \ar[d]^{\prhsharp} \\
\rhB_+ \ar[rr]_{\vph \otimes_\sigma \Eins_k} && ^{\pvrh}\hspace{-2pt}\pB_+ }
\end{align*}
of set maps commutative. Hence, by Exc.~\ref{pr.SELIBA}, it is a homomorphism of cubic alternative algebras.

(c) We must show that the composites 
\[
\Pocu \circ \Cosp\:\kpoc \to \kpoc \quad \text{and} \quad \Cosp \circ \Pocu\:\kosp \to \kosp 
\]
are (naturally) isomorphic to the respective identity functors. As to the first case, let $(A,q)$ be a pointed cubic alternative $k$-algebra. Then $\Cosp(A,q) = (\msB(A,q), \vrh)$, where $\msB(A,q) = (k \times k,B,\tau,p)$, $B = A \times A^{q\op}$, $\tau = \vep_B$, $p = (q,q)$ and $\vrh = \Eins_{k \times k}$. Writing $\pi_+\:K \to k$ and $\Pi_+\:B \to A$ for the respective projections onto the first factor, we obtain
\[
\Pi_+\big((\alpha_+,\alpha_-)(u_+,u_-)\big) = \Pi_+\big((\alpha_+u_+,\alpha_-u_-)\big) =\alpha_+u_+ = \pi_+\big((\alpha_+,\alpha_-)\big)\Pi_+(\big((u_+,u_-)\big)
\]
for $\alpha_\pm \in k$, $u_\pm \in A$. Moreover $\Pi_+$ commutes with the respective adjoints:
\[
\xymatrix{B \ar[r]_{\Pi_+} \ar[d]_{\sharp_B = \sharp_A \times \sharp_A} & A \ar[d]^{\sharp_A} \\
B \ar[r]_{\Pi_+} & A.}
\]
By Exc.~\ref{pr.SELIBA}, therefore, $\Pi_+$ is a $\pi_+$-semi-linear homomorphism of cubic alternative algebras, and part (a) of the same exercise yields a unique homomorphism $\Pi_{A,q}\:\rhB_+ = B \otimes_K \,\rhk_+ \to A$ of cubic alternative $\rhk_+$- algebras satisfying $\Pi_{A,q} \circ \can_{B,\,\rhk_+} = \Pi_+$, hence
\[
\Pi_{A,q}\big((u_+,u_-) \otimes_K \alpha_+\big) = \alpha_+(\Pi_{A,q} \circ \can_{B,\,\rhk_+})\big((u_+,u_-)\big)= \alpha_+\Pi_+\big((u_+,u_-)\big) = \alpha_+u_+
\]
for $u_\pm \in A$, $\alpha_+ \in \,\rhk_+$. On the other hand, we have a natural map $\vph\:A \to \,\rhB_+$ given by $u_+ \mapsto ((u_+,0) \otimes_K 1_{\rhk_+})$, and one checks that $\Pi_{A,q}$ and $\vph$ are inverse to one another. In particular, $\Pi_{A,q}$ is a bijection sending $\rhp_+ = (q,q) \otimes_K 1_{\rhk_+}$ to $q$. Summing up, therefore, we have found an isomorphism
\[
\Pi_{A,q}\:\Pocu(\Cosp(A,q)) \overset{\sim} \longrightarrow (A,q)
\]
of pointed cubic alternative $k$-algebras such that
\begin{align}
\label{PIAQ} \Pi_{A,q}\big((u_+,u_-) \otimes_K \alpha_+\big)= \alpha_+u_+
\end{align}
for $u_\pm \in A$, $\alpha_+ \in \,\rhk_+$. Now let $\vph\:(A,q) \to (\pA,\pq)$ be a homomorphism of pointed cubic alternative $k$-algebras. Then we have $\Cosp(\vph) = (\Eins_K,\vph \times \vph)$ and $\Pocu(\Cosp(\vph)) = (\vph \times \vph) \otimes_{\Eins_K} \Eins_k = (\vph \times \vph) \otimes_K \Eins_k$. We therefore have to show that the diagram
\[
\xymatrix{(\Pocu \circ \Cosp)(A,q) \ar[r]_(.65){\Pi_{A,q}} \ar[d]_{(\vph \times \vph) \otimes_K \Eins_k} & (A,q) \ar[d]^{\vph} \\
(\Pocu \circ \Cosp)(\pA,\pq) \ar[r]_(.65){\Pi_{\pA,\pq}} & (\pA,\pq)}
\]
commutes, which follows immediately from \eqref{PIAQ}.

Turning to the second case, let $(\msB,\vrh)$ with $\msB= (K,B,\tau,p)$ be a core-split involutorial system over $k$. We first treat the special case $K = k \times k$ and $\vrh = \Eins_{k \times k}$. As before, let $\pi_\pm\:K \to k$ be the two canonical projections making $k$ a $K$-algebra $k_\pm$. Then $B_\pm := B_{k_\pm}$ are cubic alternative $k$-algebras with norms $N_\pm = N_B \otimes_K k_\pm$, giving rise to natural identifications $B = B_+ \times B_-$ as a cubic alternative $K$-algebra with norm $N = N_+ \times N_-$. Put $\vep_+ := (1,0) \in K$, $\vep_- := (0,1) = \bar\vep_+ \in K$. Since $\tau$ is semi-linear with respect to the conjugation of $K$, there are $k$-linear maps $\sigma_\pm\:B_\pm \to B_\mp$ such that $\tau\big((u_+,u_-)\big) = \big(\sigma_-(u_-),\sigma_+(u_+)\big)$ for all $u_\pm \in B_\pm$, and since $\tau^2 = \Eins_B$, the map $\sigma_+$ is bijective with inverse $\sigma_-$. Thus
\begin{align}
\label{TAUUP} \tau\big((u_+,u_-)\big) = \big(\sigma_+^{-1}(u_-),\sigma_+(u_+)\big)&&(u_\pm \in B_\pm).
\end{align}
By (\ref{ss.COINSY}.\ref{ENTAUIO}), the relation $N_B(\tau(u)) = \overline{N_B(u)}$ holds strictly for all $u \in B$, which by \eqref{TAUUP} implies
\begin{align}
\label{ENMISA} N_- \circ \sigma_+ = N_+
\end{align}
as polynomial laws over $k$. Since $p = (p_+,p_-)$ remains fixed under $\tau$, we conclude
\begin{align}
\label{PEMIS} p_- = \sigma_+(p_+).
\end{align} 
$\tau$ being an $\iota$-semi-linear isotopy involution of $B$, hence a $k$-isomorphism $B \to B^{p\op}$, it may equally well be viewed as a $k$-isomorphism $B^{p\op} \to B$, which amounts to $\tau((vp^{-1})(pu)) = \tau(u)\tau(v)$ for all $u,v \in B$. Expanding with the aid of \eqref{TAUUP}, \eqref{PEMIS}, and applying \eqref{ENMISA}, we conclude that
\[
\sigma_+\:B_+^{p_+\op} \overset{\sim} \longrightarrow B_-
\]
is an isomorphism of cubic alternative $k$-algebras. It follows that
\begin{align}
\label{CHIEIN} \chi := \chi_{(\msB,\Eins_{k \times k})} := \Eins_{B_+} \times \sigma_+\:B_+ \times B_+^{p_+\op} \overset{\sim} \longrightarrow B
\end{align}
is an isomorphism of cubic alternative $K$-algebras. We now claim that
\begin{align}
\label{EICHI} (\Eins_K,\chi):\Cosp\big(\Pocu(\msB,\Eins_{k \times k})\big) \overset{\sim} \longrightarrow (\msB,\Eins_{k \times k})
\end{align}
 is an isomorphism in $\kosp$. From \eqref{PEMIS}, \eqref{CHIEIN} we deduce $\chi((p_+,p_+)) = p$, so we only have to show that $\chi$ respects isotopy involutions: $\chi \circ \vep_{B_+} = \tau \circ \chi$, which follows from \eqref{TAUUP} and a straightforward computation. Next we show that any morphism in $\kosp$ having the form $(\Eins_{k \times k},\vph)\:(\msB,\Eins_{k \times k}) \to (\pmsB,\Eins_{k \times k})$ gives rise to a commutative diagram
\begin{align}
\vcenter{\label{COPOBE} \xymatrix{\Cosp\big(\Pocu(\msB,\Eins_{k \times k})\big) \ar[rrr]_{\Cosp(\Pocu(\Eins_{k \times k},\vph))} \ar[d]_{(\Eins_{k \times k},\chi_{(\msB,\Eins_{k \times k})})}^{\cong} &&& \Cosp\big(\Pocu(\pmsB,\Eins_{k \times k})\big) \ar[d]_{\cong}^{(\Eins_{k \times k},\chi_{(\pmsB,\Eins_{k \times k})})} \\
(\msB,\Eins_{k \times k}) \ar[rrr]_{(\Eins_{k \times k},\vph)} &&& (\pmsB,\Eins_{k \times k}).}}
\end{align}
Indeed, \eqref{POCUSI} implies $\Pocu(\Eins_{k \times k},\vph) = \vph \otimes_{k \times k} \Eins_k = \vph_+$, hence $\Cosp(\Pocu(\Eins_{k \times k},\vph)) = (\Eins_{k \times k},\vph_+ \times \vph_+)$ by \eqref{COPHI}. Combining with \eqref{TAUUP} and the relation $\ptau \circ \vph = \vph \circ \tau$, the assertion follows.

Finally, we turn to the general case. By Exc.~\ref{pr.MOVCOR}, $(\vrh,\can_{\msB})\:(\msB,\vrh) \to (\rhmsB,\Eins_{k \times k})$ is an isomorphism of core-split involutorial systems over $k$. Write $\rhK$ for $k \times k$ viewed as a $K$-algebra via $\vrh$ and let $(\sigma,\vph)\:(\msB,\vrh) \to (\pmsB,\pvrh)$ be a morphism in $\kosp$. Then $\vph$ and $\Eins_{k \times k}\:\rhK \to\, \prhpK$ are both $\sigma$-semi-linear. By \ref{ss.TESELI}, therefore, $\vph \otimes_\sigma \Eins_{k \times k}$ exists as a $\sigma$-semi-linear map $\rhB \to\,\pprhB$, and one checks that it is compatible with adjoints, hence a $\sigma$-semi-linear isomorphism of cubic alternative algebras (Exc.~\ref{pr.SELIBA}~(b)). Moreover,
\[
(\Eins_{k \times k},\vph \otimes_\sigma \Eins_{k \times k})\:(\rhmsB,\Eins_{k \times k}) \longrightarrow (\pprhmsB,\Eins_{k \times k})
\]
is a homomorphism of core-split involutorial systems over $k$ which is easily seen to make the diagram
\begin{equation}
\label{BEROPH} \xymatrix{(\msB,\vrh) \ar[rr]_{(\sigma,\vph)} \ar[d]_{(\vrh,\can_{\msB})}^{\cong} && (\pmsB,\pvrh) \ar[d]_{\cong}^{(\pvrh,\can_{\pmsB})} \\
(\rhmsB,\Eins_{k \times k}) \ar[rr]_{(\Eins_{k \times k},\vph \otimes_\sigma \Eins_{k \times k})}&& ({\pprhmsB},\Eins_{k \times k})}
\end{equation}
commutative. Combining \eqref{COPOBE}, \eqref{BEROPH} and applying the functor $\Cosp \circ \Pocu$ to the latter, we obtain the commutative diagram
\[
\xymatrix{\Cosp\big(\Pocu(\msB,\vrh)\big) \ar[rrr]_{\Cosp(\Pocu(\sigma,\vph))} \ar[d]_{\Cosp(\Pocu(\vrh,\can_{\msB}))}^{\cong} &&& \Cosp\big(\Pocu(\pmsB,\pvrh)\big) \ar[d]_{\cong}^{\Cosp(\Pocu (\pvrh,\can_{\pmsB})} \\
\Cosp\big(\Pocu(\rhmsB,\Eins_{k \times k})\big) \ar[rrr]_{\Cosp(\Pocu(\Eins_{k \times k},\vph \otimes_\sigma \Eins_{k \times k}))} \ar[d]_{(\Eins_{k \times k},\chi_{(\msB,\Eins_{k \times k})})}^{\cong} &&& \Cosp\big(\Pocu(\pprhmsB,\Eins_{k \times k}) \ar[d]_{\cong}^{(\Eins_{k \times k},\chi_{(\pmsB,\Eins_{k \times k})})} \\
(\rhmsB,\Eins_{k \times k}) \ar[rrr]_{(\Eins_{k \times k},\vph \otimes_\sigma\Eins_{k \times k})} \ar[d]_{(\vrh,\can_{\msB})^{-1}}^{\cong} &&& (\pprhmsB,\Eins_{k \times k}) \ar[d]_{\cong}^{(\pvrh,\can_{\pmsB})^{-1}} \\
(\msB,\vrh) \ar[rrr]_{(\sigma,\vph)} &&& (\pmsB,\pvrh),}
\]
from which an isomorphism from $\Cosp \circ \Pocu$ to the identity functor of $\kosp$ can be immediately read off.
\end{sol}

\begin{sol}{pr.ETSITI} \label{sol.ETSITI}
 (a) (\ref{t.EXSETI}.\ref{TRAEXSE}) implies $w \in J_0^\perp$, while (\ref{t.EXSETI}.\ref{SHAEXSE}) yields $Q(w) = u(p\tau(u))$. In conjunction with Prop.~\ref{p.CUJAZ} we conclude 
 \begin{align*}
 N_{J_0}(Q(w)) =\,\,&N_B(u(p\tau(u))) = N_B(p)N_B(u)\overline{N_B(u)} \\
 =\,\,&\mu\bar\mu N_B(u)\overline{N_B(u)} = n_K(\mu N_B(u)).
 \end{align*}
Moreover, $N_J(w) = t_K(\mu N_B(u))$ by (\ref{t.EXSETI}.\ref{NOEXSE}). Hence (i) and (ii) are equivalent. (ii) combined with Prop.~\ref{p.SMASU} shows that $k[\mu N_B(u)] \subseteq K$ is a quadratic \'etale subalgebra. Since $K$ has rank $2$ as a finitely generated projective $k$-module, we have equality and (iii) holds. Finally, if (iii) holds, the preceding formulas imply the final assertion of (a), whence $w$ must be an \'etale element relative to $J_0$.
	
(b) Write $\pJ$ for the Jordan subalgebra of $J$ generated by $J_0$ and $j$. Then $\pJ$ is stable under the (bilinearized) adjoint, and we deduce from (\ref{ss.FOSETI}.\ref{EXZERU}) that $J_0j = J_0\pt j =J_0 \times j$ is contained in $\pJ$. On the other hand, (\ref{t.EXSETI}.\ref{SHAEXSE}) implies that $(p^\sharp j)^\sharp = -N_J(p)p + N_J(p)\bar\mu j$ belongs to $\pJ$. Hence so does the second summand, and since $\bar\mu$ by (a) generates $K$ as a $k$-algebra, we conclude $Kj \subseteq \pJ$. Summing up, therefore, $Bj = (KJ_0 )j \subseteq \pJ$, and we have shown $\pJ= J$, i.e., $J_0$ and $j$ generate $J$ as a Jordan $k$-algebra.
\end{sol}

\begin{sol}{pr.MOPA} \label{sol.MOPA}
(a) By (\ref{l.IDTAUIS}.\ref{TAUTAU}) and \eqref{BEPOW}, the quantity $p\pt w$ remains fixed under $\tau$. Hence $\msB\pt w$ will be an involutorial system over $k$ once we have shown that $\tau$ is a $k$-homomor\-phism from $B^w$ to $(B^w)^{(p\pt w)\op}$. It is certainly one from $B$ to $B^{p\op}$. Using (\ref{ss.UNTISOTALT}.\ref{ITUNTISOTALT}) both for $B$ and $B^w$ in place of $A$, it follows that $\tau$ is also a $k$-homomorphism from$B^w$ to
\begin{align*}
\big((B^p)^{\op}\big)^{\tau(w)} =\,\,&\big((B^p)^{\tau(w^{-1})}\big)^{\op} = (B^{p\tau(w^{-1})})^{\op} = \Big(\big((B^w)^{w^{-1}}\big)^{p\tau(w^{-1})}\Big)^{\op} \\
=\,\,&\big((B^w)^{(w^{-1}w^{-1})(w[p\tau(w^{-1})])}\big)^{\op} = \big((B^w)^{w^{-1}(p\tau(w^{-1}))}\big)^{\op} = (B^w)^{(p\pt w)\op},
\end{align*}
as desired. The assertion about $\mu$ being obvious, we proceed to show that $\Phi := \Phi_{\msB,w}$ is an isomorphism of cubic Jordan algebras. By Exc.~\ref{pr.HOCUJO}~(a), it suffices to show that $\Phi$ preserves adjoints. Let $x_0 \in H(\msB)$, $u \in B$ and put $x := x_0 + uj \in J(\msB,N_B(w)\mu)$. Then (\ref{t.EXSETI}.\ref{SHAEXSE}) implies
\[
x^\sharp = \Big(x_0^\sharp - u\big(p\tau(u)\big)\Big) + \Big(\overline{N_B(w)}\bar\mu\tau(u^\sharp)p^{-1} - x_0u\Big)j,
\]
hence
\begin{align}
\label{PHIXSH} \Phi(x^\sharp) = \Big(x_0^\sharp - u\big(p\tau(u)\big)\Big) + \Big(\overline{N_B(w)}\bar\mu\big[\tau(u^\sharp)p^{-1}\big]w - (x_0u)w\Big)j.
\end{align}
On the other hand, $\Phi(x) = x_0 + (uw)j \in J(\msB\pt w,\mu)$, hence
\begin{align}
\label{PHIXSHA} \Phi(x)^\sharp =\,\,&\Big(x_0^\sharp - \big[(uw)w^{-1}\big]\big[w\big\{\big((p\pt w)w^{-1}\big)\big(w\tau(uw)\big)\big\}\big]\Big) \\
\,\,&+ \Big(\bar\mu\big[\tau\big((uw)^\sharp\big)w^{-1}\big]\big[w(p\pt w)^{-1}\big] - (x_0w^{-1})(wuw)\Big)j. \notag
\end{align}
By the Moufang identities, the very last summands on the right of \eqref{PHIXSH} and \eqref{PHIXSHA} are the same. It therefore suffices to show that the second (resp. third) summand on the right of \eqref{PHIXSH} agrees with the second (resp. third) summand on the right of \eqref{PHIXSHA}. Indeed, making free use of the Moufang identities combined with (\ref{l.IDTAUIS}.\ref{TAUP}), \eqref{BEPOW}, the second summand on the right of \eqref{PHIXSHA} is the negative of
\[
u\Big(\big(w(p\pt w)\big)\tau(uw)\Big)= u\Big(\big(p\tau(w^{-1})\big)\tau(uw)\Big) = u\Big(p\tau\big((uw)w^{-1}\big)\Big) = u\big(p\tau(u)\big),
\]
which is also the negative on the right of the second summand of \eqref{PHIXSH}. Similarly, the equations
\begin{align*}
\bar\mu\Big(\tau\big((uw)^\sharp\big)w^{-1}\Big)\big(w(p\pt w)^{-1}\big) =\,\,&\bar\mu\Big(\tau\big((uw)^\sharp\big)w^{-1}\Big)\Big(w\big(\tau(w)p^{-1}\big)w\Big) \\=\,\,&\bar\mu\Big(\tau(w^\sharp u^\sharp)\big(\tau(w)p^{-1}\big)\Big)w = \bar\mu\Big(\tau\big(w(w^\sharp u^\sharp)\big)p^{-1}\Big)w \\
=\,\,&\overline{N_B(w)}\bar\mu\big(\tau(u^\sharp)p^{-1}\big)w
\end{align*}
show that the third summands on the right of \eqref{PHIXSH} and \eqref{PHIXSHA} are the same. Summing up we have proved that $\Phi$ is an isomorphism of cubic Jordan algebras. 

(b) For the first part, we must show $(p\pt v)\pt w = p\pt(vw)$, where we have to keep in mind that the second dot on the left refers to the algebra $B^v$ rather than $B$. Hence (\ref{l.IDTAUIS}.\ref{TAUP}) implies
\begin{align*}
(p\pt v)\pt w =\,\,&(w^{-1}v^{-1})\Big(v\big[\big((p\pt v)v^{-1}\big)\big(v\tau(w^{-1})\big)\big]\Big) \\
=\,\,&(vw)^{-1}\Big(\big(v(p\pt v)\big)\tau(w^{-1})\Big) = (vw)^{-1}\Big(\big(p\tau(v^{-1})\big)\tau(w^{-1})\Big) \\
=\,\,&(vw)^{-1}\Big(p\tau\big((vw)^{-1}\big)\Big) = p\pt (vw).
\end{align*}
The assertion about $\mu$ (resp. $N_B(w)\mu$ or $N_B(vw)\mu$) being an admissible scalar for $\msB\pt(vw)$ (resp. $\msB\pt v$ or $\msB$) being obvious, we show that the diagram \eqref{JAYBEP} is commutative. For $u \in B$ we compute
\begin{align*}
\Phi_{\msB\pt v,w} \circ \Phi_{\msB,v}(uj) =\,\,&\Phi_{\msB\pt v,w}\big((uv)j\big) = \Big(\big((uv)v^{-1}\big)\big(vw\big)\Big)j \\
=\,\,&\big(u(vw)\big)j = \Phi_{\msB,vw}(uj),
\end{align*}
and the assertion follows.

(c) For any $w \in B^\times$, we conclude from (a), (b) that $\pmu := N_B(w)^{-1}\mu$ is an admissible scalar for $\msB\pt w$ and
\[
\Phi_{\msB\pt w,w^{-1}}\:J(\msB\pt w,\pmu) \overset{\sim} \longrightarrow J(\msB,\mu)
\] 
is an isomorphism. In the special case $w :=\mu p^{-1}$ we conclude $\pmu = \mu^{-3}\mu\bar\mu\mu = \bar\mu\mu^{-1}$, hence $N_{B^w}(p\pt w) = n_K(\pmu) = 1$.
\end{sol}

\begin{sol}{pr.NINETY} \label{sol.NINETY} For sake of contradiction, suppose such a $z$ exists and write it as $z = x + iy$, $x,y \in \IZ$.  Cross-multiplying the equation $i = \overline{z}z^{-1}$, we obtain
\[
x - iy = \bar z = iz = -y + ix.
\]
Thus $y = -x$ and $z = x(1 - i)$. Hence $x \ne 0$ and taking norms, we obtain $n_C(z) = 2\vert x\vert > 1$ so $z$ is not invertible in $C$ (Prop.~\ref{p.CHIN}), a contradiction.
\end{sol}

\begin{sol}{pr.COSPAZ} \label{sol.COSPAZ}
Since $K$ is split, it contains a complete orthogonal system $(\vep_+,\vep_-)$ of elementary idempotents. Hence $K = k_+ \times k_-$, $k \cong k_\pm := k\vep_\pm \in \Kalg$, which implies $A = A_+ \times A_-$, $A_\pm = A_{k_\pm}$; in particular, $A_\pm$ are Azumaya algebras of degree $n$ over $k$. Since $\tau$ is unitary, there are $k$-linear maps $\tau_\pm\:A_\pm \to A_\mp$ such that $\tau((x_+,x_-)) = (\tau_-(x_-),\tau_+(x_+))$ for all $x_\pm \in A_\pm$, and since $\tau$ has order $2$, the map $\tau_+$ must be bijective with inverse $\tau_-$: $\tau((x_+,x_-)) = (\tau_+^{-1}(x_-),\tau_+(x_+))$. Finally, since $\tau$ is an anti-isomorphism, so is $\tau_+$, i.e., $\tau_+\:A_+^{\op} \to A_-$ is an isomorphism. Now one checks that
\[
\vph\:(A_+ \times A_+^{\op},\vep_{A_+}) \overset{\sim} \longrightarrow (A,\tau), \quad (x_+,y_+) \longmapsto \big(x_+,\tau_+(y_+)\big)
\]
is an isomorphism of algebras with involution. This proves the first part with $B = A_{k_+}$.

If $K$ is arbitrary, then $K_K \cong K \times K$ by Exc.~\ref{pr.MISPLET}, and (\ref{ss.RESC}.\ref{EMAR}) yields $A_K = (_kA)_K = A_{K_K} \cong A_{K \times K} = A \times A$ as $K \times K$-algebras. Now the first part proves the second.
\end{sol}

\begin{sol}{pr.REGFIR} \label{sol.REGFIR}
 Put $J_0 := H(\msB)$ and write $T_0$ (resp. $T$) for the bilinear trace of $J_0$ (resp. $B$). By (\ref{t.EXSETI}.\ref{TRAEXSE}), we obtain
\begin{align}
\label{TEJAYX} T_J(x,y) = T_0(x_0,y_0) + 2T\Big(\big(p\tau(u)\big)v\Big) 
\end{align}
for $x = x_0 + uj$, $y = y_0 + vj$, $x_0,y_0 \in J_0$, $u,v \in B$. In standard notation we have $J^\ast = J_0^\ast \oplus B^\ast j$, and writing $T_J^\ast\:J \to J^\ast$ for the linear map induced by the symmetric bilinear form $T_J$ (ditto for $T_0^\ast$, $T^\ast$), equation \eqref{TEJAYX} implies 
\begin{align}
\label{TEJAST} T_J^\ast = T_0^\ast \oplus (2\vph)j, \quad \vph := T^\ast \circ L_p  \circ \tau\:B \longrightarrow B^\ast.
\end{align}
Assume first that $B$ is regular and $2 \in k^\times$. Then $B = J_0 \oplus B_-$, $B_- = \{u \in B\mid \tau(u) = -u\}$, and since $B$ is finitely generated projective as a $k$-module, so are $J_0,B_-$ and $J$. Moreover, $B^\ast = J_0^\ast \oplus B_-^\ast$, and $T(J_0,B_-) = \{0\}$ implies $T^\ast = T_0^\ast \oplus T_-^\ast$, $T_-$ being the restriction of $T$ to $B_- \times B_-$. From this we deduce that $T_0^\ast$ is bijective and hence, by \eqref{TEJAST}, so is $T_J^\ast$. Thus $J$ is regular. Conversely, assume that $J$ is regular. Then $J_0$ and $B$ are finitely generated projective as $k$-modules and we conclude from \eqref{TEJAST} that multiplication by $2$ is a bijection $B \to B$; in particular, some $u \in B$ has $2u = 1_B$. But $1_B \in B$ is unimodular, so there exists a $k$-linear form $\lambda\:B \to k$ satisfying $\lambda(1_B) = 1$. Thus $2\lambda(u) = 1$ makes $2$ invertible in $k$, as claimed.
\end{sol}

\begin{sol}{pr.ISSETI} \label{sol.ISSETI}
 By Exc.~\ref{pr.MOPA}~(c), up to isomorphism of $J$, we may assume $N_B(p) = 1_K$, $n_K(\mu) = 1$. Passing to the $p$-isotope of $J$ and invoking Thm.~\ref{t.ISOTSE}, we may indeed assume $p = 1_B$ up to isotopy of $J$.  Now assume that the norm of $B$ is surjective as a set map from $B$ to $K$. Then Exc.~\ref{pr.MOPA}~(a) shows that, up to isomorphism of $J$, we may assume $\mu = 1_K$, which in particular implies $N_B(p) = 1_K$. Arguing as before, we may therefore assume $p = 1_B$, $\mu = 1_K$ up to isotopy of $J$.
\end{sol}

\begin{sol}{pr.INVQUA} \label{sol.INVQUA}
(a) It suffices to verify that $\tau_B$ is an involution of $C$. Localizing if necessary we may assume $C = \Cay(B,\mu) = B \oplus Bj_0$ for some $\mu \in k^\times$. Then $\vph := \tau \circ \iota_C\:C \to C$ satisfies $\vph(v + wj_0) = \tau(\bar v - wj_0) = v - wj_0$ for $v,w \in B$, hence is an isomorphism (Exc.~\ref{pr.CDSCAPAR}), forcing $\tau = \vph \circ \iota_C$ to be an involution.

(b) is clear. 

(c) Arbitrary elements of $J$ have the form
\begin{align}
\label{ELEXZU} x = x_0 + uj, \quad x_0 = (\alpha_0,\beta_01_C + w_0), \quad u = (\alpha,v + w)
\end{align}
with $\alpha_0,\beta_0,\alpha \in k$, $v \in B$, $w_0,w \in B^\perp$. Writing $T_0,S_0,N_0$ (resp. $T,S,N$) for the trace, quadratic trace, norm of $J_0$ (resp, $J$), we begin by using Prop.~\ref{p.PELE} to determine the Peirce components of $J$ relative to $e_1$: 
\begin{align}
\label{TWOONE} J_2(e_1) =\,\,&ke_1, \\
\label{ONEONE} J_1(e_1) =\,\,&(\{0\} \times C)j, \\ 
\label{ZERONE} J_0(e_1) =\,\,&\big(\{0\} \times H(C,\tau_B)\big) \oplus (ke_1)j.
\end{align}
here \eqref{TWOONE} is clear, while for the rest (\ref{t.EXSETI}.\ref{LTRAEXSE}), (\ref{t.EXSETI}.\ref{BSHAEXSE}) and Exc.~\ref{pr.CUADUN} imply
\begin{align}
\label{TEAZER} T(x) =\,\,&\alpha_0 + 2\beta_0, \\
\label{EONSHA} e_1 \times x =\,\,&(0,\beta_01_C - w_0) - (\alpha,0)j,
\end{align}
which immediately yields \eqref{ONEONE} but also \eqref{ZERONE} since $e_0 := 1_J - e_1 = (0,1_C)$ and
\begin{align*}
T(x)e_0 - x =\,\,&(\alpha_0 + 2\beta_0)(0,1_C) - (\alpha_0,\beta_01_C + w_0) - (\alpha,v + w)j \\=\,\,&\big(-\alpha_0,(\alpha_0 + \beta_0)1_C - w_0\big) - (\alpha,v + w)j.
\end{align*}
Now assume $x \in (J_0)_0(e_1)$, i.e., $\alpha_0 = \alpha = 0$, $v = w = 0$. Then (\ref{t.EXSETI}.\ref{SHAEXSE}) yields
\begin{align}
\label{EXSHAB} x^\sharp = \big(\beta_0^2 + n_C(w_0)\big)e_1.
\end{align}
Hence an elementary frame of the desired kind exists if and only if some $e_2 = (0,\beta_01_C + w_0)$, $\beta_0 \in k$, $w_0 \in B^\perp$ has $2\beta_0 = T(e_2) = 1$ and $e_2^\sharp = 0$, i.e., $\beta_0 = \frac{1}{2}$ and $n_C(w_0) = -\frac{1}{4}$ by \eqref{EXSHAB}. Combining Prop.~\ref{p.CDUNV} with regularity, we find an identification of $C$ as in \eqref{CECABE} matching $j_0$ with $2w_0$. Thus $e_2$ satisfies the first equation of \eqref{COPLET}, forcing $e_3 := 1_J - e_1 - e_2$ to do the same with the second. Conversely, if \eqref{CECABE} and \eqref{COPLET} hold, then $\Omega$ is an elementary frame of the desired kind.

For the rest of the solution, we assume this to be the case. We proceed to determine the Peirce-one-spaces of
\begin{equation}
\label{EPM} e_\pm := \big(0,\frac{1}{2}(1_C \pm j_0)\big),\;\text{i.e., of $e_+ = e_2$, $e_- = e_3$}.
\end{equation}
Adjusting the notation of \eqref{ELEXZU} to the identification \eqref{CECABE} by writing
\begin{equation}
\label{ELEXZA} x = x_0 + uj, \quad x_0 = (\alpha_0,\beta_01_C + w_0j_0), \quad u = (\alpha,v + wj_0)
\end{equation}
for $\alpha_0,\beta_0,\alpha \in k$, $w_0,v,w \in B$, we obtain
\begin{align*}
e_\pm \times x =\,\,&\big(0,\frac{1}{2}(1_C \pm j_0)\big) \times (\alpha_0,\beta_01_C + w_0j_0) + \big(0,\frac{1}{2}(1_C \pm j_0)\big) \times \big((\alpha,v + wj_0)j\big) \\
=\,\,&\big(\frac{1}{2}n_C(1_C \pm j_0,\beta_01_C + w_0j_0),\alpha_0\frac{1}{2}(1_C \mp j_0)\big) - \Big(\big(0,\frac{1}{2}(1_C \pm j_0)\big)\big(\alpha,v + wj_0\big)\Big)j \\
=\,\,&\big(\beta_0 \mp \frac{1}{2}t_B(w_0),\frac{\alpha_0}{2}1_C \mp \frac{\alpha_0}{2}j_0\big) - \big(0,\frac{1}{2}(v \pm \bar w) + \frac{1}{2}(w \pm \bar v)j_0\big)j.
\end{align*}
For $x$ as in \eqref{ELEXZA} to belong to $J_1(e_\pm)$ it is necessary and sufficient that $T(x) = 0$ and $e_\pm \times x = 0$. Hence
\begin{align}
\label{ONETWO} J_1(e_2) =\,\,&\big\{(0,w_0j_0) + (\alpha,v - \bar vj_0)j \mid \alpha \in k\;w_0 \in B^0,\;v \in B\big\}, \\
\label{ONETHREE} J_1(e_3) =\,\,&\big\{(0,w_0j_0) + (\alpha,v + \bar vj_0)j \mid \alpha \in k,\;w_0 \in B^0,\;v \in B\} .
\end{align}
Combining with \eqref{ONEONE}, the off-diagonal Peirce components of $J$ relative to $\Omega$ attain the form
\begin{align}
\label{OFFOT} J_{12} =\,\,&\{(0,v - \bar vj_0)j \mid v \in B\}, \\
\label{OFFTT} J_{23} =\,\,&\{(0,w_0j_0) + (\alpha,0)j \mid w_0 \in B^0,\,\alpha \in k\}, \\
\label{OFFTO} J_{31} =\,\,&\{(0,v + \bar vj_0)j \mid v \in B\}.
\end{align}
Now put
\begin{align}
\label{UTTO} u_{23} = (1,0)j \in J_{23}, \quad u_{31} = \big(0,\frac{1}{2}(1_B + j_0)\big)j \in J_{31}.
\end{align}
Using (\ref{t.EXSETI}.\ref{QUAEXSE}) one checks $S(u_{23}) = S(u_{31}) = -1$, so $\mfS:=(\Omega,u_{23},u_{31})$ is a strong co-ordinate system of $J$ in the sense of \ref{ss.STROCO}.

Next we consider the multiplicative conic alternative algebra $\pB:= C_{J,\mfS}$ of ~\ref{p.CONALT} living uon the $k$-module $J_{12}$ under the multiplication (\ref{p.CONALT}.\ref{CONPRO}) with $\omega = 1$. In order to describe this multiplication explicitly, we let $v,w \in B$, put $x := (0,v - \bar vj_0)j$, $y := (0,w - \bar wj_0)j$, write $x\bu y$ for their product in $\pB$ and obtain
\begin{align*}
x \times u_{23} =\,\,&(0,v - \bar vj_0)j \times (1,0)j = \htau_B\big((0,v - \bar vj_0) \times (1,0)\big)j \\
=\,\,&\htau_B\big((0,\bar v + \bar vj_0)\big)j = (0,v + \bar vj_0)j, \\
u_{31} \times y =\,\,&\big(0,\frac{1}{2}(1_B + j_0)\big)j \times (0,w - \bar wj_0)j \\
=\,\,&-\frac{1}{2}\big((0,1_B + j_0)(0,\bar w - \bar wj_0) + (0,w - \bar wj_0)(0,1_B + j_0)\big) \\
\,\,&+ \htau_B\Big(\big(0,\frac{1}{2}(1_B + j_0)\big) \times (0,w - \bar wj_0)\Big)j \\
=\,\,&-\frac{1}{2}(0,\bar w - \bar wj_0 + wj_0 - w + w +wj_0 - \bar wj_0 - \bar w) \\
\,\,&+ \htau_B\Big(\big(\frac{1}{2}n_C(1_B + j_0,w - \bar wj_0),0\big)\Big)j \\
=\,\,&\big(0,(\bar w - w)j_0\big) + \big(t_B(w),0\big)j,
\end{align*}
hence
\begin{align*}
x\bu y =\,\,& (x \times u_{23}) \times (u_{31} \times y) \\
=\,\,&(0,v + \bar vj_0)j \times \big(0,(\bar w - w)j_0\big) + (0,v + \bar vj_0)j \times \big(t_B(w),0\big)j \\
=\,\,&\Big(\big(0,(w - \bar w)j_0\big)(0,v + \bar vj_0)\Big)j - (0,v + \bar vj_0)\big(t_B(w),0\big) \\
\,\,&- \big(t_B(w),0\big)(0,\bar v + \bar vj_0) + \htau_B\Big((0,v + \bar vj_0)\times \big(t_B(w),0\big)\Big)j \\
=\,\,&\Big(0,\big((w - \bar w)\bar v\big)j_0 + v(w - \bar w)\Big)j + \htau_B\Big(\big(0,t_B(w)(\bar v - \bar vj_0)\Big)j \\
=\,\,&\Big(0,v(w - \bar w) + \big((w - \bar w)\bar v\big)j_0 + v(w + \bar w) - \big((w + \bar w)\bar v\big)j_0\Big)j \\
=\,\,&\big(0,2vw - (2\overline{vw})j_0\big)j.
\end{align*}
Thus the assignment $v \mapsto \frac{1}{2}(v - \bar vj_0)j$ yields an isomorphism $B \overset{\sim} \to \pB$ of composition algebras, and the Jacobson co-ordinatization theorem \ref{t.JACO} yields \eqref{JAYHER}. For the final assertion, we put $\Gamma = \diag(\gamma_1,\gamma_2,\gamma_3)$ and, by Exc.~\ref{pr.DIAFIX}, may assume $\det(\Gamma) = 1$. Then
\[
p := \sum \gamma_ie_{ii} \in \Her_3(B)\;\text{and}\;q := \sum\gamma_i^{-1}e_i \in J_0
\]
have norm $1$, and the isomorphism \eqref{JAYHER} matches $p^{-1}$ with $q$. By Exc.e~\ref{pr.DIISJO}, the assignment
\[
\sum(\xi_ie_{ii} + u_i[jl]) \longmapsto \sum(\gamma_i\xi_ie_{ii} + \gamma_j\gamma_lu_i[jl])
\]
defines an isomorphism $\vph\:\Her_3(B,\Gamma)^{(p)} \to \Her_3(B)$, and since $\vph(p^{-2}) = p^{-1}$, Thm.~\ref{t.ISOTSE} implies
\[
\Her_3(B,\Gamma) \cong \Her_3(B)^{(p^{-1})} \cong J^{(q)} \cong J(\qmsB,1),
\]
as claimed.
\end{sol}

\begin{sol}{pr.CAINSY} \label{sol.CAINSY}
(a) We have $J_0 = k \times k1_C$, so the elements of $J$ have the form
\begin{align}
\label{EXAlZE} x = (\alpha_0,\beta_01_C) + (\alpha,v)j &&(\alpha_0,\beta_0,\alpha \in k,\;v \in C),
\end{align}
and, writing $T$ for the linear trace of $J$, we conclude
\begin{align}
\label{TRAEXA} T(x) = \alpha_0 + 2\beta_0.
\end{align}
While $e_1$ is clearly an elementary idempotent, $e_2$ has trace $1$ and satisfies
\begin{align*}
e_2^\sharp =\,\,&\frac{1}{4}\big((0,1_C) + (1,0)j\big)^\sharp \\
=\,\,&\frac{1}{4}\big((1,0) - (1,0)(1,0)\big) + \frac{1}{4}\Big(\hiota_C\big((1,0)^\sharp\big) - (0,1_C)(1,0)\Big)= 0,
\end{align*}
hence is an elementary idempotent as well. Furthermore,
\begin{align*}
e_1 \times e_2 =\,\,&\frac{1}{2}(1,0) \times (0,1_C) + \frac{1}{2}(1,0) \times (1,0)j = \frac{1}{2}(0,1_C) - \frac{1}{2}(1,0)j \\
=\,\,&\frac{1}{2}(e_0 - e_1j) = e_3.
\end{align*}
By Prop.~\ref{p.COMPOR}, therefore, $\Omega$ is an elementary frame of $J$. Now assume $x \in J_1(e_3)$. Then $T(x) = 0$, i.e., $\alpha_0 + 2\beta_0 = 0$ by \eqref{TRAEXA}, and $e_3 \times x = 0$. But, quite generally,
\begin{align*}
e_3 \times x =\,\,&\frac{1}{2}\big((0,1_C) - (1,0)j\big) \times \big((\alpha_0,\beta_01_C) + (\alpha,v)j\big) \\
=\,\,&\frac{1}{2}\big((2\beta_0,\alpha_01_C) - (0,v)j + (\alpha_0,0)j + (\alpha,0) + (\alpha,0) - (0,v)j\big) \\
=\,\,&(\beta_0 + \alpha,\frac{1}{2}\alpha_01_C) + (\frac{1}{2}\alpha_0, -v)j.
\end{align*}
This proves $J_1(e_3) = \{0\}$. Moreover, $x \in J_0(e_3)$ means $e_3 \times x = T(x)(1_J - e_3) - x$, i.e., 
\begin{align*}
(\beta_0 + \alpha,\frac{1}{2}\alpha_01_C) + (\frac{1}{2}\alpha_0, -v)j =\,\,&(\alpha_0 + 2\beta_0)\big((1,\frac{1}{2}1_C) + (\frac{1}{2},0)j\big) - (\alpha_0,\beta_01_C) - (\alpha,v)j \\
=\,\,&(2\beta_0,\frac{1}{2}\alpha_01_C) + (\frac{1}{2}\alpha_0 + \beta_0 - \alpha,-v)j,
\end{align*}
and we conclude 
\[
J_0(e_3) = \{(\alpha_0,\alpha1_C) + (\alpha,v)j \mid \alpha_0,\alpha \in k, v \in C\},
\]
i.e.,
\begin{align}
J_0(e_3) = M_0 := ke_1 + ke_2 + (\{0\} \times C)j.
\end{align}
By Faulkner's lemma~\ref{c.FAULE}, it remains to determine $q_0 := S\vert _{M_0}$. For $\alpha_1,\alpha_2 \in k$, $v \in C$ we compute
\begin{align*}
q_0\big(\alpha_1e_1 + \alpha_2e_2 + (0,v)j\big) =\,\,&S\big((\alpha_1,\frac{1}{2}\alpha_21_C) + (\frac{1}{2}\alpha_2,v)j\big) \\
=\,\,&T\big((\frac{1}{4}\alpha_2^2,\frac{1}{2}\alpha_1\alpha_21_C)\big) - T\big((\frac{1}{2}\alpha_2,v)(\frac{1}{2}\alpha_2,\bar v)\big) \\
=\,\,&\frac{1}{4}\alpha_2^2 + \alpha_1\alpha_2 - \frac{1}{4}\alpha_2^2 - 2n_C(v) = \alpha_1\alpha_2 - 2n_C(v),
\end{align*}
and the final assertion follows.

(b) Since $2 = 0$ in $k$, we have
\begin{align}
\label{HASCR} H(\hmsC) = k \times C^0, \quad C^0 = \Ker(t_C).
\end{align}
Note that $C$ is regular since $r > 1$, so by Lemma~\ref{l.COMTRI}, some $w \in C$ has $t_C(w) = 1$. Thus $C = kw \oplus C^0$, and $C^0$ is finitely generated projective of rank $r - 1$, proving the rank formula for $J$. Moreover, the second Tits construction using $\hmsC$ is stable not just under flat but actually under arbitrary base change.

(i) Arbitrary elements of $J$ have the form
\begin{align}
\label{TWOZER} x = x_0 + uj, \quad x_0 = (\alpha_0,v_0), \quad u = (\alpha,v) &&(\alpha_0,\alpha \in k,\;v_0 \in C^0,\;v \in C).
\end{align}
Hence
\begin{align}
\label{TEEXA} T(x) = \alpha_0 
\end{align}
and
\begin{align}
\label{EXTWOZ} x^\sharp =\,\,& \Big((\alpha_0,v_0)^\sharp - (\alpha,v)\hiota_C\big((\alpha,v)\big)\Big) + \Big(\hiota_C\big((\alpha,v)^\sharp\big) - (\alpha_0,v_0)(\alpha,v)\Big)j \\
=\,\,&\Big(\big(n_C(v_0),\alpha_0v_0\big) -\big(\alpha^2,n_C(v)1_C\big)\Big) + \Big(\big(n_C(v),\alpha v\big) - (\alpha_0\alpha,v_0v)\Big)j \notag \\
=\,\,&\big(n_C(v_0) - \alpha^2,\alpha_0v_0 - n_C(v)1_C\big) + \big(n_C(v) - \alpha_0\alpha,\alpha v - v_0v\big)j. \notag
\end{align}
Thus $x$ is an elementary idempotent if and only if $\alpha_0 = 1$, $v_0 = n_C(v)1_C$, $\alpha = n_C(v)$, and we have shown that the elementary idempotents of $J$ are precisely of the form $e_v$, $v \in C$. Since the assignment $v \mapsto e_v$ is compatible with base change, it yields an isomorphism from $C_{\bfa}$ to $\bfEld(J)$. Conversely, the projection onto the last factor produces a morphism from $\bfEld(J)$ to $C_{\bfa}$, and by what we have just seen, the two morphisms are inverse to one another.

(ii) For $v,w \in C$, we obtain
\begin{align*}
T(e_v,e_w) =\,\,&T_0\Big(\big(1,n_C(v)1_C\big),\big(1,n_C(w)1_C\big)\Big) + T\Big(\big(n_C(v),v\big)\big(n_C(w),\bar w\big)\Big) \\
\,\,&+ T\Big(\big(n_C(w),w\big)\big(n_C(v),\bar v\big)\Big) \\
=\,\,&1 + n_C(vw) - t_C(v\bar w) - n_C(vw) - t_C(w\bar v) = 1, 
\end{align*}
whence $e_v, e_w$ can never be orthogonal.

(iii) Since $2 = 0$ in $k$ and $k$ is reduced, Exc.~\ref{pr.CUBNIL} yields
\begin{align}
\label{NILTZ} \Nil(J) = \big\{x \in J \mid T(x,J) = T(x^\sharp,J) = \{0\}\big\}.
\end{align}
For $\beta_0,\beta \in k$, $w_0 \in C^0$, $w \in C$, we obtain
\begin{align*}
T\big(x,(\beta_0,w_0)\big) =\,\,&\alpha_0\beta_0 + t_C(v_0,w_0), \\T\big(x,(\beta,w)\big)=\,\,&T\big((\alpha,v)(\beta,\bar w)\big) + T\big((\beta,w)(\alpha,\bar v)\big)= 0,
\end{align*}
so $T(x,J) = \{0\}$ is equivalent to $\alpha_0 = 0$ and $v_0 \in k1_C$. By \eqref{EXTWOZ}, \eqref{NILTZ}, therefore, $x \in \Nil(J)$ if and only if $\alpha_0 = 0$, $v_0 = \gamma 1_C$ for some $\gamma \in k$ and $\gamma^2 = n_C(v_0) = \alpha^2$. But since $2 = 0$ in $k$ and $k$ is reduced, this implies $\gamma = \alpha$, and (iii) is proved. 
\end{sol}

\begin{sol}{pr.INSEKI} \label{sol.INSEKI}
(a) By (\ref{ss.SELICON}.\ref{ENCEL}) we have $t_C(\tau(x)) = \iota_K(t_C(x))$ for all $x \in C$, which implies that the $k$-linear maps $\iota_C$ and $\tau$ commute. Thus $\tau_0 := \tau \circ \iota_C\:C \to C$ is an $\iota_K$-semi-linear automorphism of order $2$. By Exc.~\ref{pr.DESELI}, therefore,
\[
C_0 := H(C,\tau_0) = \{x \in C \mid \tau(x) = \bar x\} \subseteq C
\]
is a unital $k$-subalgebra, and the inclusion $C_0 \hookrightarrow C$ induces a $K$-linear bijection $\Phi\:C_0 \otimes K \to C$ making the diagram
\[
\xymatrix{C_0 \otimes K \ar[r]_{\Eins_{C_0} \otimes \iota_K} \ar[d]_{\Phi}^{\cong} & C_0 \otimes K \ar[d]_{\cong}^{\Phi} \\
C \ar[r]_{\tau_0} & C}
\]
commutative. Now one checks that $\Phi\:(C_0 \otimes K,\iota_{C_0} \otimes \iota_K) \to (C,\tau)$ is an isomorphism of $k$-algebras with involution. 

(b) One checks that $\msB$ is a unitary involutorial system over $k$. Setting $J := J(\msB,\mu)$, Cor.~\ref{c.SEFITI} and Exc.~\ref{pr.COHACE}~(a) imply that $J_K \cong J(\hC,\mu) \cong \Her_3(C,\Gamma_0)$, $\Gamma_0 := \diag(-1,-1,1)$, is a Freudenthal algebra of rank $3(r + 1)$ over $K$.  Since $K$, being quadratic \'etale, is faithfully flat over $k$, the assertion now follows from Cor.~\ref{c.CHARF}.
\end{sol}

\begin{sol}{pr.OPINSY} \label{sol.OPINSY}
 Since $\tau\:B \to B^{p\op}$ is an $\iota_K$-semi-linear isomorphism of cubic alternative $K$-algebras, it is also one from $B^{\op}$ to $B^p = (B^{\op})^{{\op}p} = (B^{\op})^{p^{-1}{\op}}$. Hence $\msB^{\op}$ is an involutorial system over $k$. Let $\mu \in K$ be an admissible scalar for $\msB$. The $\mu^{-1}$ is clearly one for $\msB^{\op}$, and we only have to show that $\Phi$ is an isomorphism. As usual, it suffices to show that $\Phi$ preserves adjoints, so let $x_0 \in H(\msB)$ and $u \in B$. Writing $u\pt v$ for the product in $B^{\op}$, we apply Lemma~\ref{l.IDTAUIS} and obtain
\begin{align*}
\Phi\big((x_0 + uj)^\sharp\big)=\,\,&\Phi\Big(\big[x_0^\sharp - u\big(p\tau(u)\big)\big] + \big[\bar\mu\tau(u^\sharp)p^{-1} - x_0u\big]j\Big) \\
=\,\,&\big[x_0^\sharp - u\big(p\tau(u)\big)\big] +
\tau\Big(\big[\bar\mu\tau(u^\sharp)p^{-1} - x_0u\big]p\Big)j \\
=\,\,&\Big(x_0^\sharp - \tau(up)\pt\big(p^{-1}\pt(up)\big)\Big) + \Big(\mu u^\sharp - \big(p\tau(u)\big)x_0\Big)j.
\end{align*}
But $\mu u^\sharp = \bar\mu^{-1}N_B(p)u^\sharp = \bar\mu^{-1}p(up)^\sharp = \bar\mu^{-1}\tau(\tau(up)^\sharp)\pt (p^{-1})^{-1}$ and $(p\tau(u))x_0 = x_0\pt \tau(up)$. Hence
\begin{align*}
\Phi\big((x_0 + uj)^\sharp\big) =\,\,&\Big(x_0^\sharp - \tau(up)\pt\big(p^{-1}\pt(up)\big)\Big) + \Big(\bar\mu^{-1}\tau\big(\tau(up)^\sharp\big)\pt(p^{-1})^{-1} \\
\,\,&- x_0\pt\tau(up)\Big)j \\
=\,\,&\big(x_0 + \tau(up)j\big)^\sharp = \Phi(x_0 + uj)^\sharp,
\end{align*}
as claimed.
\end{sol}


\solnsec{Section~\ref{s.FREPASE}}

\begin{sol}{pr.CUETET} \label{sol.CUETET}
 (a), (b) The assertion of (a) is clear if $E$ is a field. Next assume $E = F \times K$, where $K/F$ is a separable quadratic field extension. Letting $u$ be a generator of $K/F$ and setting $a := (1,u) \in E$, one checks that $1_E = (1,1_K),a,a^2$ are linearly independent over $F$. Hence $a$ generates $E$ as a unital $F$-algebra. Finally let us assume that $E = F \times F \times F$ is split. If $F$ contains more than two elements, some quantity $\alpha \in F^\times$ is different from $1$. Since the components of the vector $a :=  (\alpha,1,0) \in E$ are mutually distinct, a Vandermonde argument shows that $a$ generates $E$ as a unital $F$-algebra. This completes the proof of (a). On the other hand, if $F = \IF_2$ is the field with two elements, then the elements of $E$ have the form $a = (\alpha,\beta,\gamma)$ where at least two of the quantities $\alpha,\beta,\gamma \in F$ are the same. Hence $a$ does not generate $E$ as a unital $F$-algebra. This completes the proof of (b).

(c) Obvious.

(d)(i) Suppose first that $a \in J$ is \'etale relative to $J_0$. Then (c) implies $T_J(a) = 0$, while we deduce from Exc.~\ref{pr.ETSITI}~(b) that $a$ generates $E$ as a unital $F$-algebra. Finally $a^\sharp = \frac{1}{3}S_J(a)\cdot 1_E + (a^\sharp - \frac{1}{3}S_J(a)\cdot 1_E)$, where the second summand belongs to $J_0^\perp$. In the terminology of \ref{ss.ETEL} , therefore, $Q(a) = -\frac{1}{3}S_J(a)$, which must be different from zero by the definition of an \'etale element.

Conversely, suppose $a \in J_0^\perp$ generates $E$ as a unital $F$-algebra and has $S_J(a) \ne 0$. We claim that the discriminant
\begin{align}
\label{DICUET} D := -4S_J(a)^3 - 27N_J(a)^2 \in F^\times.
\end{align}
This is well known, see for example \cite[\S3.1]{Bh:small}.  Here is a sketch of a proof.

There is clearly no harm in assuming that $F$ is algebraically closed. Then $E = F \times F \times F$ is split and for some $\alpha_1,\alpha_2,\alpha_3 \in F$ we have $a = (\alpha_1,\alpha_2,\alpha_3)$ as well as $\sum\alpha_i = 0$. Short computations give
\[
S_J(a) = -(\alpha_1^2 + \alpha_1\alpha_2 + \alpha_2^2), \quad N_J(a) = -(\alpha_1^2\alpha_2 + \alpha_1\alpha_2^2)
\]
and then
\begin{align}
\label{DISPLI} D =\,\,&4\alpha_1^6 + 12\alpha_1^5\alpha_2 - 3\alpha_1^4\alpha_2^2 - 26\alpha_1^3\alpha_2^3 - 3\alpha_1^2\alpha_2^4 + 12\alpha_1\alpha_2^5 + 4\alpha_2^6.
\end{align}
On the other hand, let $\vrh\:E \to E$ be the \emph{shift operator} sending $(\alpha_1,\alpha_2,\alpha_3)$ to $(\alpha_3,\alpha_1,\alpha_2)$. Then $\vrh \in \Aut(E)$ has order $3$ and fixed algebra $F\cdot 1_E$. We have
\begin{align}
\label{ARHA} a - \vrh(a) = (2\alpha_1 + \alpha_2,\alpha_2 - \alpha_1, -\alpha_1 -2\alpha_2),
\end{align} 
which implies
\[
N_J\big(a - \vrh(a)\big) = (2\alpha_1 + \alpha_2)(\alpha_1 - \alpha_2)(\alpha_1 + 2\alpha_2) = 2\alpha_1^3 + 3\alpha_1^2\alpha_2 - 3\alpha_1\alpha_2^2 - 2\alpha_2^3
\]
Now a straightforward computation combined with \eqref{DISPLI} yields
\begin{align}
\label{NORDIS} N_J\big(a - \vrh(a)\big)^2 = D
\end{align}
Since $a$ generates $E$ as a unital $F$-algebra, the quantities $\alpha_1,\alpha_2,-(\alpha_1 + \alpha_2)$ are mutually distinct, and we conclude from \eqref{ARHA} that $a - \vrh(a)$ is invertible in $E$. This in conjunction with \eqref{NORDIS} proves $D \in F^\times$, as claimed.

As before we have $Q(a) = -\frac{1}{3}S_J(a)\cdot 1_J$. By hypothesis, this belongs to $J_0^\times$, and since $E/F$ by what we have just shown has discriminant
\begin{align*}
0 \neq D =\,\,&-4S_J(a)^3 - 27N_J(a)^2 = -27\Big(N_J(a)^2 - 4N_{J_0}\big(Q(a)\big)\Big), 
\end{align*} 
$a$ is an \'etale element of $J$ relative to $J_0$.

(d)(ii) We must show that $E$ is split cubic \'etale and $J = J(F,1)$ over $F = \IF_4$ if and only if $(J,J_0)$ does not admit \'etale elements.

Suppose first $E = F \times F \times F$ is split cubic \'etale and $J = J(F,1)$ over $F = \IF_4$. Note that $\ch(F) = 2$ and $F^\times = \{1,\zeta,\zeta^2\}$, with $\zeta \in F$ a primitive cube root of $1$. Assuming $a = (\alpha,\beta,\gamma) \in J$ with $\alpha,\beta,\gamma \in F$ is \'etale relative to $J_0$, then (i) implies $0 = T_J(a) = \alpha + \beta + \gamma$ and that $a$ generates $E$ as a unital $F$-algebra, equivalently, that $\alpha,\beta,\gamma$ are mutually distinct. If $\gamma = 0$, then $\beta = \alpha$, hence $a = \alpha(1,1,0)$ and $a^\sharp = \alpha^2(0,0,1)$, so we have the relation $\alpha a+ a^\sharp = \alpha^2\cdot 1_J$, whence $a$ cannot generate $E$ as a unital $F$-algebra, a contradiction.By symmetry, therefore, we may assume $\alpha,\beta,\gamma \in F^\times$, so up to applications of the shift operator we have $a = (1,\zeta,\zeta^2) \in J_0^\perp$, hence $a^\sharp = (1,\zeta^2,\zeta) \in J_0^\perp$. In particular, $a$ cannot be \'etale relative to $J_0$.

Conversely, suppose $(J,J_0)$ does not admit \'etale elements and apply (a) to pick $a \in J$ that generates $E$ as a unital $F$-algebra. Replacing $a$ by $a - \frac{1}{3}T_J(a)\cdot 1_J$ if necessary, we may assume $T_J(a) = 0$. Then (i) implies $S_J(a) = 0$ since $a$ is not \'etale relative to $J_0$. Hence $a = N_J(a)\cdot 1_J$, and since $E$ is reduced as a commutative ring, we conclude $N_J(a) \ne 0$. Also, $b := a + a^\sharp \in J_0^\perp$, and (\ref{ss.BACUNO.fig}.\ref{CADJ}) implies $b^\sharp = a^\sharp + a \times a^\sharp + a^{\sharp\sharp} = -N_J(a)\cdot 1_J + N_J(a)a + a^\sharp$. Thus $S_J(b) = -3N_J(a) \ne 0$. But $b$ is not \'etale relative to $J_0$, so we deduce from (i) that $1_J,b,b^\sharp$ are linearly dependent, yielding $\alpha,\beta \in F$ which satisfy $b^\sharp = \alpha1_J + \beta b = \alpha1_J + \beta a + \beta a^\sharp$. Comparing coefficients gives $N_J(a) = \beta = 1$. Summing up, therefore, 
\begin{align}
\label{EFFTEE} E = F[\bft]/(\bft^3 - 1) = F \times K, \quad K := F[\bft]/(\bft^2 + \bft + 1).    
\end{align}
Since $a \in J$ is a Kummer element of norm $1$ relative to $J_0$ (Thm.~\ref{t.REMSTAB}), Cor.~\ref{c.UPFITI} yields a natural identification $(J,J_0) = (J(F,1),F\cdot 1_J)$. Combining Prop.~\ref{p.ETFITI} with the fact that this Freudenthal pair does not admit \'etale elements implies $\alpha^3 = 1$ for all $\alpha \in F^\times$. But this amounts to $F = \IF_4$ and that $K$ as defined in \eqref{EFFTEE} is split. Hence so is $E$, and the proof of (ii) is complete.
\end{sol}

\begin{sol}{pr.SPRIFOSP} \label{sol.SPRIFOSP}
 \cite[pp.~96-98]{MR2032452}.
\end{sol}

\begin{sol}{pr.TERNET} \label{sol.TERNET}
(a) The first statement is clear while the second one follows immediately from (\ref{ss.THERCU}.\ref{MABIT}).
	
(b) We apply (\ref{ss.THERCU}.\ref{MANO}) to obtain $N_J(w) =\gamma_1\gamma_2\gamma_3t_C(w_1w_2w_3)$, so we only have to worry about the last factor in this product, where we note that Thm.~\ref{t.TERHER} implies $w_3 \in M = D^\perp$, hence $\bar w_3 = -w_3$, and then
\begin{align*}
t_C(w_1w_2w_3) =\,\,&t_C\big((w_1w_2)w_3\big) = n_C(w_1w_2,\bar w_3) = -n_C(w_1w_2,w_3) \\
=\,\,& -n_C\big(-h(w_1,w_2) + w_1 \times_{h,\Delta} w_2,w_3\big) = -t_D\big(h(w_1 \times_{h,\Delta} w_2,w_3)\big).
\end{align*}
Now the definition of the hermitian vector product (\ref{ss.HERVE}) yields \eqref{TERTRA}. The verification of \eqref{TERNOR} is a bit more involved. From (\ref{ss.THERCU}.\ref{MADJ}) and Thm.~\ref{t.TERHER} we deduce 
\begin{align*}
Q(w) = \sum\big(\gamma_j\gamma_lh(w_i,w_i)e_{ii} + \gamma_i\overline{h(w_j,w_l)}[jl]\big), 	  
\end{align*}
and taking norms via (\ref{ss.THERCU}.\ref{MANO}) yields
\begin{align*}
N_{J_0}\big(Q(w)\big) = \,\,&(\gamma_1\gamma_2\gamma_3)^2\Big(h(w_,w_1)h(w_2,w_2)h(w_3,w_3) - \sum h(w_i,w_i)n_D\big(h(w_j,w_l)\big) \\
\,\,&+ t_D\big(h(w_3,w_2)h(w_1,w_3)h(w_2,w_1)\big)\Big) \\
=\,\,&(\gamma_1\gamma_2\gamma_3)^2\big(h(w_1,w_1)h(w_2,w_2)h(w_3,w_3)  \\
\,\,&- h(w_1,w_1)h(w_2,w_3)h(w_3,w_2) -h(w_2,w_2)h(w_3,w_1)h(w_1,w_3)\\
\,\,&- h(w_3,w_3)h(w_1,w_2)h(w_2,w_1) + h(w_3,w_2)h(w_1,w_3)h(w_2,w_1)\\
\,\,&+ h(w_1,w_2)h(w_3,w_1)h(w_2,w_3)\big) \\
=\,\,&(\gamma_1\gamma_2\gamma_3)^2\det\Big(\big(h(w_i,w_j)\big)_{1\le i,j\le 3}\Big) \\
=\,\,&(\gamma_1\gamma_2\gamma_3)^2(\bigwedge^3 h)(w_1 \wedge w_2, \wedge w_3,w_1 \wedge w_2, \wedge w_3)
\end{align*}
But since $\Delta\:\bigwedge^3(M,h) \overset{\sim} \to (D,\la\det_\Delta(h)\ras)$ is an isometry by \ref{ss.DET} and $\det_\Delta(h) = 1$ by hypothesis, the last expression equals $(\gamma_1\gamma_2\gamma_3)^2n_D(\Delta(w_1\wedge w_2\wedge w_3))$, and \eqref{TERNOR} is proved.

(c) Applying \eqref{TERTRA}, \eqref{TERNOR}, we see that (i) holds if and only if $\delta(w) \in D^\times$ and $t_D(\delta(w))^2 - 4n_D(\delta(w)) \in D^\times$. By Prop.~\ref{p.SMASU}, the latter condition is equivalent to $D =k[\delta(w)]$. Hence (i) and (ii) are equivalent.

(d) We begin with a simple lemma.
\begin{lem*}
Let $M$ be a finitely generated projective $k$-module of rank $n \in \IN$. For elements $v_1,\dots,v_n \in M$ to form a basis of $M$ it is necessary and sufficient that $v_1 \wedge \cdots \wedge v_n$ form a basis of $\det(M) = \bigwedge^n(M)$.
\end{lem*} 

\begin{proof} 
The condition is clearly necessary. Conversely, suppose $v_1 \wedge \cdots \wedge v_n$ is a basis of $\det(M)$. Localizing if necessary, we may assume that $M$ is free. Let $e_1,\dots,e_n$ be a basis of $M$ and write $v_j = \sum_{i=1}^n \alpha_{ij}e_i$, for some $\alpha_{ij} \in k$, $1\le i,j\le n$. Then $v_1 \wedge \cdots \wedge v_n = \det((\alpha_{ij}))e_1 \wedge \cdots \wedge e_n$, 
so $\det(\alpha_{ij}) \in k^\times$, and the assertion follows. 
\end{proof} 

If $w \in J$ is an \'etale element relative to $J_0$, it may be written as $w = \sum w_i[jl]$, $w_i \in M$, $1\le i\le 3$, and (c) combined with the lemma implies that $w_1,w_2,w_3$ is a $D$-basis of $M$, and $D = k[a]$, $a = \delta(w) \in D^\times$. Conversely, let $v_1,v_2,v_3$ be a $D$-basis of $M$ and $a \in D^\times$ such that $D = k[a]$. Again by the lemma, $b := \Delta(v_1 \wedge v_2 \wedge v_3) \in D^\times$, and setting $w_1 := v_1b^{-1}a$, $w_i := v_i$ for $i = 2,3$, we put $w := \sum w_i[jl]$ and conclude $\delta(w) =a$, forcing $w$ by (c) to be an \'etale element of $J$ relative to $J_0$.
\end{sol}

\begin{sol}{pr.SKONOF} \label{sol.SKONOF}
 By definition, $J_0$ contains an elementary frame, which by Prop.~\ref{p.ELFRAFI} can be extended to a co-ordinate system $\mfS = (e_1,e_2,e_3,u[23],u[31])$ of $J_0$. Note that $\mfS$ is clearly a co-ordinate system of $J$, and consulting Prop.~\ref{p.CONALT}, we see that $C_{J_0,\mfS} \subseteq C_{J,\mfS}$ is a composition subalgebra. Moreover, $\Gamma := \Gamma_{J_0,\mfS} = \Gamma_{J;\mfS}$, and the Jacobson co-ordinatization theorem~\ref{t.JACO} yields isomorphisms $C_0 \cong C_{J_0,\mfS}$, $C \cong C_{J,\mfS}$ as well as \eqref{ISOFRE}. That the existence of the latter isomorphism is, in fact, independent of the preceding ones follows from the Skolem-Noether theorem for composition algebras (Exc.~\ref{pr.SKONOCO}).
\end{sol}

\begin{sol}{pr.NOPER} \label{sol.NOPER}
 We consider the cubic form $g := N_J\vert_{J_0^\perp}\:J_0^\perp \to F$ and assume that the set map $g_F\:J_0^\perp \to F$ is (identically) zero. The polynomial law $g$ is certainly different from zero since passing to the algebraic closure $\aclosF$ of $F$ allows us to identify $J_{\aclosF} = \Her_3(C)$, $C = \aclosF \times \aclosF$ and $J_{0\aclosF} = \sum \aclosF e_{ii}$, hence $J_{0\aclosF}^\perp = \sum C[jl]$, whence $u := \sum c[jl] \in J_{0\aclosF}^\perp$, $c := (1,0) \in C$ by (\ref{ss.THERCU}.\ref{MANO}) satisfies $g_{\aclosF}(u) = N_J(u) = t_C(c) = 1$. From Exc.~\ref{pr.EXOCUB} combined with (\ref{ss.BACUNO.fig}.\ref{GRADI}) linearized we therefore conclude
\[
T_J(u \times v,w) = N_J(u,v,w) = g_F(u,v,w) = 0
\]
for all $u,v,w \in J_0^\perp$. Since the bilinear trace of $J$ is regular on $J_0^\perp$, this implies $u \times v \in J_0$ for all $u,v \in J_0^\perp$. Thanks to bilinearity, this condition is stable under base change, hence continues to hold over $\aclosF$. On the other hand, applying (\ref{ss.EICJM}.\ref{OFFTI}), we obtain $1_C[23] \times 1_C[31] = 1_C[12] \notin J_{0\aclosF}$,a contradiction. Thus $N_J(u) \ne 0$ for some $u \in J_0^\perp$.
\end{sol}

\begin{sol}{pr.NINTWO} \label{sol.NINTWO}  As in the proof of Thm.~\ref{t.ALPAIR}, we define polynomial laws $e_i\:J\times J \to J$ for $1 \le i \le 9$ by
\begin{align*}
e_1(x,y) := 1_{J_R}, \quad e_2(x,y) := x, \quad e_3(x,y) := x^\sharp, \quad e_4(x,y) := y, \quad e_5(x,y) := y^\sharp, \\
e_6(x,y) := x \times y, \quad e_7(x,y) := x^\sharp \times y, \quad e_8(x,y) := x \times y^\sharp, \quad e_9(x,y) := x^\sharp \times y^\sharp
\end{align*}
for all $x,y \in J_R$, $R \in \kalg$.
We claim that \emph{it suffices to show
\begin{align}
\label{DETTEI} \det\Big(\big(T_J(e_i(x,y),e_j(x,y))_{1\le i,j\le 9}\big)\Big) \in k^\times 
\end{align}
for some $x,y \in J$.} Indeed, once such elemnts have been exhibited, the  quantities $e_i(x,y)$, $1\le i\le 9$, generate a free submodule $M \subseteq J$ of rank $9$ on which the bilinear trace of $J$ is regular. By Lemma~\ref{l.NOSU}, this implies $J = M \oplus M^\perp$, and by comparing ranks, we conclude $J = M$. On the other hand, $M$ is the subalgebra of $J$ generated by $x,y$ (Exc.~\ref{pr.GENTWO}), and we are done.

Thus we need only show that elements $x,y \in J$ exist satisfying \eqref{DETEI}. Since $k$ is an LG ring, the $k$-module $J \times J$ is free of rank $54$, so the scalar polynomial law $\det(T_J(e_i,e_j))$ by Cor.~\ref{c.FREPOLA} may be regarded as a polynomial in $k[\bft_1,\dots,\bft_{54}]$, which by the LG property represents an invertible element of  $k$ if and only if it represents an invertible element of $k/\mfm$, for each maximal ideal $\mfm \subseteq k$. We may therefore assume that $k = F$ is a field.

If $F$ contains more than four elements, we let $E \subseteq J$ be any cubic \'etale subalgebra. Then $(J,E)$ is a Freudenthal pair over $F$, hence by Thm.~\ref{t.ELGEET} admits \'etale elements. By Cor.~\ref{c.PHISER}, therefore, $J = J(\msB,\mu) = H(\msB) \oplus Bj$ is a second Tits construction, for some involutorial system $\msB = (K,B,\tau,p)$ of the second kind over $k$ and some admissible scalar $\mu$ for $\msB$. such that $E = H(\msB)$ and $j \in J$ is \'etale relative to $E$. Now Exc.~\ref{pr.ETSITI}~(b) shows that $J$ is generated by $E$ and $j$ as a Jordan $F$-algebra.  Since $E$ has a single generator by Exc.~\ref{pr.CUETET}~(a), the assertion follows.

We are left with the case that $F = \IF_2$ and consider the following subcases. (i) $J = \Mat_3(F^{(+)}$. Let $E/F$ be the unique (cyclic) cubic field extension of $F$, $\vrh$ a generator of its Galois group and write $A$ for the cyclic algebra $(E/F,\vrh,1)$ over $F$ \cite[8.5, p.~484]{MR1009787}. By \cite[Thm.~8.7, p.~477]{MR1009787}, $A$ is central simple of degree $3$ over $F$, hence isomorphic to $\Mat_3(F)$. It therefore suffices to show that the Jordan algebra $A^{(+)}$ is generated by two elements. Recall that $A$ as a three-dimensional left vector space over $E$ has a basis $1_A,y,y^2$ for some $y \in A$ having $y^3 = 1_A$, subject to the relations $ya = \vrh(a)y$ for all $a \in A$. This implies $(ay)^3 = (ay^2)^3 = N_E(a)1_A$, hence $T_A(ay) = T_A(ay^2) = 0$. Thus $y$ and $y^\sharp = y^2$ (by (\ref{ss.BACUNO.fig}.\ref{ADJJA})) both belong to $E^\perp$, and we conclude that $y$ is a Kummer element of $A^{(+)}$ relative to $E$. By \ref{ss.FOFITI}~(a), therefore, $A^{(+)}$ is generated by $E$ and $y$, and since $E$ is generated by a single element, we are done. (ii) $J = \Her_3(K)$, where $K/F$ is the unique quadratic field extension. Consulting Exc.~\ref{pr.COHACE}, we see that $J$ is isomorphic to the first Tits construction $J(F \times K,1)$. Since $F \times K$ is generated by one element (Exc.~\ref{pr.CUETET}~(a),(b)), $A^{(+)}$ is generated by two elements.
\end{sol}

\begin{sol}{pr.FREPAN} \label{sol.FREPAN}
 Arguing as in the proof of Thm.~\ref{t.ALPAIR}, we define polynomial laws $e_i\:J \to J$, $1 \le i \le 3$, by $e_1(x) := 1_J$, $e_2(x) := x$, $e_3(x) := x^\sharp$ for all $x \in J_R$, $R \in \kalg$. We claim that it suffices to show
\[
\det\Big(\big(T_J(e_i(x),e_j(x))\big)_{1\le i\le 3}\Big) \in k^\times
\]
for some $x \in J$. Indeed, if such an element has been exhibited, it generates a subalgebra $E \subseteq J$ which, as a submodule, is spanned by $1,x,x^\sharp$ ((\ref{ss.BACUNO.fig}.\ref{CADJ}) and Exc.~\ref{pr.GENTWO}) and hence is cubic \'etale, making $(J,E^{(+)})$ a Freudenthal pair over $k$. To find an element  making the above determinant invertible, we are reduced in the usual way to the case that $F := k$ is a field. If $F$ is algebraically closed, then $J = \Her_3(F \times F)$ is split and $x := \sum\xi_ie_{ii} \in J$, $\xi_i \in F$ mutually distinct, satisfies the desired condition. By Zariski density, this also settles the problem if $F$ is infinite. On the other hand, if $F$ is finite, it suffices to find a cubic \'etale $E \subseteq J$ generated by a single element. If $J = \Mat_3(F)^{(+)}$ is split, any cubic subfield of $J$ will do. If $J = \Her_3(K)$, $K$ the unique quadratic field extension of $F$, the norm of $K$ is a universal quadratic form, and setting $\Gamma_0 := \diag(-1,-1,1)$, (\ref{pr.ISCOPA}.\ref{NOCHA}) and Exc.~\ref{pr.COHACE} yield 
\[
J = \Her_3(K) \cong \Her_3(K,\Gamma_0) \cong J(F \times K,1),
\]
whence $F \times K$ is a cubic \'etale subalgebra of $J$ generated by a single element (Exc.~\ref{pr.CUETET}~(a),(b)).
\end{sol}

\begin{sol}{pr.RECUALF} \label{sol.RECUALF}
If $A$ has rank $n$, then $J(A,\mu)$ is finitely generated projective of rank $3n$ as a $k$-module and regular as a cubic Jordan algebra, so $(J(A,\mu),A^{(+)})$ is a balanced pair of cubic Jordan algebras in the sense of \ref{ss.BACUJO}. It remains to show that $J(A,\mu)$ is a Freudenthal algebra. By \ref{ss.GEOSIM}, this will follow once we have shown for all fields $K \in \kalg$ that $J(A,\mu)_K \cong J(A_K,\mu_K)$ is simple or cubic \'etale. We are thus reduced to the case that $k = F$ is a field, which we may assume as well from the outset to be algebraically closed.  The latter hypothesis makes the norm of $A$ surjective as a set map from $A$ to $F$, and by Exercise~\ref{pr.CHAMU}~(c), we may assume $\mu = 1$. According to Exercise~\ref{pr.SESIAL}, this up to isomorphism leaves the following options for $A$.

\begin{holgerenum}
\item $A = F$, $\ch(F) \ne 3$: By \ref{ss.FOFITI}~(b), $J = J(F,1)$ has degree $3$, so some $x \in J$ makes $1,x,x^2$ an $F$-basis of $J$, and we conclude that $J = F[x]$ is a $3$-dimensional regular cubic Jordan $F$-algebra. As such, it has no absolute zero divisors (Exercise~\ref{pr.ABZECU}~(e)) and hence is locally linear (Theorem~\ref{t.AZLOC}). Summing up, therefore, $J$ is split cubic \'etale.

\item $A = F \times C$, $C$ a regular (split) composition $F$-algebra: Since $F$ is algebraically closed,  Exercises~\ref{pr.COHACE}~(a) and \ref{pr.DIAFIX}~(b) imply $J \cong \Her_3(C)$, which is simple by Exercise~\ref{pr.IDCUJM}~(e).

\item $A = \Mat_3(F)$: By Exercise~\ref{pr.NONKUM}~(a) combined with Corollary~\ref{c.ALFITI}, $J = J(\Mat_3(F),1)$ is a split Albert algebra over $F$, hence simple (Exercise~\ref{pr.IDCUJM}~(e)).
\end{holgerenum}
\end{sol}

\solnsec{Section~\ref{s.FREDI}}

\begin{sol}{pr.C2CSA}
Suppose first that $A$ is a division algebra and note that the algebra $A$ has dimension $d^2$ for some $d$. For $c \in F$, consider the polynomial law $A \times F \to F$ defined by $(a, x) \mapsto \Nrd_A(a) - cx^d$.  It is a homogenous form of degree $d$ in $d^2 + 1$ variables so it is isotropic.  That is, there is some $(a, x)$ such that $\Nrd_A(a) = cx^d$.  Since $A$ is division, $cx^d = \Nrd_A(a) \ne 0$, whence $x \ne 0$ and we find $\Nrd_A(a/x) = \Nrd_A(a)/x^d = c$.

In the general case, by the Wedderburn-Artin Theorem, $A \cong M_n(D)$ for some $n \ge 1$ and some division algebra $D$, both of which are  determined by $A$.  By the preceding paragraph, we may assume that $n > 1$ and the claim holds for $D$, i.e., there is an element $d \in D$ such that $\Nrd_D(d) = c$.  Define $a := \diag(d, 1, 1, \ldots, 1) \in M_n(D)$.  It suffices to show
\[
\Nrd_{\Mat_n(D)}(a) = \Nrd_D(d). 
\]
In order to see this, we extend scalars to the algebraic closure $K$ of $F$ and write $r$ for the degree of $D$. Then the left-hand side of the displayed equation is the determinant of the block diagonal matrix $\diag(d_K,\Eins_r,\ldots,\Eins_r)$, which obviously is the same as the determinant of $d_K$, hence to its right-hand side.
\end{sol}

\begin{sol}{pr.TWOINV} \label{sol.TWOINV}
(a) From Cor.~\ref{c.AZUNTI} and Thm.~\ref{t.TITDIV} we deduce that $J$ is a reduced Albert algebra over $F$.

(b) Following Cor.~\ref{c.ALSETI}, we may realize $J$ by means of the second Tits construction: $J = J(B,\tau,p,\mu)$ as in part (b) of the problem. We define $f_3(J)$ as the $3$-Pfister form of the left $p$-twist of $\tau$. In order to show that $f_3(J)$ thus defined satisfies the characterizing conditions of (b), we may assume that $J$ itself is reduced, so by Thm.~\ref{t.TITDIV} some $q \in B^\times$ has $N_B(q) = \mu$. Now apply Exc.~\ref{pr.MOPA} to $\msB := (K,B,\tau,p)$ and $w := q^{-1}$. Since $B$ is associative, we have $\msB\pt w = (K,B,\tau,p\pt w)$, where $p\pt w = qp\tau(q)$, and one checks that the assignment $u \mapsto quq^{-1}$ gives an isomorphism $(B,\lptau) \to (B,\pwtau)$ of algebras with involution. Hence passing from $\msB$ to $\msB\pt w$ doesn't change the isometry class of $f_3(J)$, and Exc.~\ref{pr.MOPA}~(a) allows us to assume $\mu = 1_K$. This implies $N_B(p) = 1$, hence $J^{(p)} \cong J(B,\lptau,1_B,1_K)$ by Thm.~\ref{t.ISOTSE}, and since, by the Jacobson-Faulkner theorem \ref{t.ALBJAC}, the coefficient algebra of $J$ does not change when passing to an isotope, we conclude that $f_3(J)$ is indeed the co-ordinate norm of $J$. This proves existence. For uniqueness, it suffices to combine Springer's theorem with Cor.~\ref{c.SEPSUF} and Prop.~\ref{p.CUFIJO}~(b)~(5).
\end{sol}

\begin{sol}{pr.LINQUADI} \label{sol.LINQUADI}
\emph{Reduction to the finite-dimensional case}. Assume the finite-dimensional case has been settled and let $0 \ne x \in J$. We must show that the linear map $V_x\:J \to J$, $y \mapsto x \circ y$, is bijective. Gven $y \in J$, the cubic Jordan subalgebra, $\pJ$, of $J$ generated by $x,y$ is finite-dimensional (Exc.~\ref{pr.GENTWO}), so by assumption $V_xz = y$ for some $z \in \pJ$, and $x \circ y = 0$ implies $y = 0$. Thus $V_x$ is indeed bijective. For the rest of the proof, we may assume that $J$ is finite-dimensional, whence it suffices to establish \eqref{LINQUADI}. Suppose on the contrary that some non-zero elements $x,y \in J$ have $x \circ y = 0$. Since $1_J \circ y = 2y \ne 0$, and by symmetry, we obtain $x,y \notin F1_J$. Hence $E := F[x]$ and $\pE := F[y]$ are cubic subfields of $J$ (Prop.~\ref{p.CUFIJO}~(b)~(1)). By (\ref{ss.BACUNO.fig}.\ref{TEXCY}), we also have $T_J(x,y) = \frac{1}{2}T_J(x \circ y) = 0$. We now consider the following cases.

\step{1} \emph{$E$ is separable over $F$.} Then $J = E \oplus E^\perp$ as a direct sum of subspaces and the action
\begin{align}
\label{ZEDDOU} E \times E^\perp \longrightarrow E^\perp, \quad (z,w) \longmapsto z\pt w := -z \times w,    
\end{align}
makes $E^\perp$ a left vector space over $E$ (Exc.~\ref{pr.SPRIFO}~(a)). Here (\ref{ss.BACUNO.fig}.\ref{TRACI}) implies $E \circ E^\perp \subseteq E^\perp$. Writing $y = u + v$, $u \in E$, $v \in E^\perp$, we obtain $0 = x \circ y = x \circ u + x \circ v$, where $x \circ u \in E$, $x \circ v \in E^\perp$, hence $x \circ u = x \circ v = 0$. But $x \circ u = 2xu$ in terms of the product in $E$, so we have $u = 0$. Thus, $y = v \in E^\perp$. Now (\ref{ss.BACUNO.fig}.\ref{LADJJA}) and \eqref{ZEDDOU} show $x\pt y = -x \times y = T_J(x)y$, hence $(x - T_J(x)1_E)\pt y = 0$. Thus $x = T_J(x)1_J$, and taking traces we obtain $2T_J(x) = 0$, hence $x = 0$, a contradiction.

\step{2} \emph{$E$ is purely inseparable over $F$.} By symmetry, we may also assume that $\pE$ is purely inseparable over $F$. Then (\ref{ss.BACUNO.fig}.\ref{LADJJA}) gives $x \times y = x \circ y = 0$, and from (\ref{ss.BACUNO.fig}.\ref{DADJ}) we deduce $0 = x^\sharp \times (x \times y) = N_J(x)y + T_J(x^\sharp,y)x$. Since $N_J(x) \ne 0$ (Cor.~\ref{c.REGINV}), this forces $y = \alpha x$, for some $\alpha \in F^\times$, and we arrive at the contradiction $0 = x \circ y = 2\alpha x^2$.
\end{sol}

\begin{sol}{pr.CYCSUB}\label{sol.CYCSUB}
It suffices to show that $S_J^0$, the restriction of the quadratic trace to the elements of trace zero in $J$, is isotropic because any isotropic vector in $J^0$ relative to $S_J^0$ generates a cubic subfield of $J$ (Prop.~\ref{p.CUFIJO}~(b)~(1)) that, by the preceding hint to this exercise, must be cyclic. Let $E/F$ be any cubic subfield of $J$ (Cor.~\ref{c.SEPSUF}). Then the base change $J_E$ by Prop.~\ref{p.CUFIJO}~(b)~(5) is a reduced regular Freudenthal algebra of dimension $\ge 9$ over $E$, hence identifies with $J_E = \Her_3(C,\Gamma)$, for some (regular) composition algebra $C$ over $E$ and some diagonal matrix $\Gamma \in \GL_3(E)$ (Prop.~\ref{p.REDDEF}). By Springer's theorem, it suffuces to show that the binary quadratic form $q$ derived from restricting the quadratic trace of $J_E$ to the \emph{diagonal} elements of trace zero is isotropic. Since the diagonal of $J_E$ is split cubic \'etale, we can apply Exc.~\ref{pr.SJ0.3}~\ref{pr.SJ0.3b} with $\beta = 1$ to conclude $q \cong \la -3,-1\ra_{\mathrm{quad}}$. But $-3$ is a square in $F$ since $F$ contains the cube roots of $1$, so $q \cong \la 1,-1\ra_{\mathrm{quad}}$ is hyperbolic  (\cite[Cor.~13.3]{MR2427530}).
\end{sol}



\solnchap{Solutions for Chapter~\ref{c.LIEALG}}

\solnsec{Section~\ref{s.LIEALG}}

\begin{sol}{pr.malcev}
The solution is modeled after pages 6--9 of \cite{Myung}.

(a):  The displayed identity holds if and only if the following expression is identically zero:
\begin{multline*}
J(x,y,[xz]) - [J(x,y,z),x] = [[xy][xz]] + [[y[xz]]x] + [[[[xz]x]y] \\ - [[[xy]z]x] - [[[yz]x]x] - [[[zx]y]x].
\end{multline*}
On the right side, the second and last terms cancel, leaving an expression which vanishes if and only if the Malcev identity holds.  

(b): The very definition of $S$ shows that it is unchanged by a cyclic permutation of its arguments.  On the other hand, linearizing the flexible law, replacing $x$ with $x+z$, gives the identity $[x,y,z] = -[z,y,x]$.  We find
\[
S(z,y,x) = [z,y,x] + [y,x,z] + [x,z,y] = -[x,y,z] - [z,x,y] -[y,z,x] = -S(x,y,z).
\]
That is, for a cyclic permutation and for the transposition interchanging 1 and 3, the claim about how $S$ changes under permuting its arguments holds.  Because these permutations generate the group of all permutations, this verifies the claim in full generality.

(c): We expand $[[x,y,z],x] = [x,y,xz] - [x,y,zx]$.  Using that the Kleinfeld function is alternating, we have
\[
0 = f(x,z,y,x) = [xz,y,x] - z[x,y,x] - [z,y,x]x
\]
and we find $[x,y,z]x=[x,y,xz]$.  Similarly, $f(z,x,y,x) = 0$ implies that $[x,y,zx] = x[x,y,z]$.  Putting these together, we obtain
\[
[[x,y,z],x] = [x,y,[x,z]] \quad \forall x, y, z \in A.
\]
Now $J(x,y,z) = 2S(x,y,z) = 6[x,y,z]$ for all $x, y, z \in A$, so $A^-$ is Malcev by (a).

For the second claim, $A^-$ is Lie if and only if $J(x,y,z) = 0$ for all $x, y, z$, if and only if $[x, y, z] = 0$ for all $x, y, z$.
\end{sol}

\begin{sol}{pr.malcev.simple}
Suppose first that $k$ has characteristic 2.  Then the kernel $I$ of the trace map $C \to k$ is a codimension 1 subspace of $C$ containing $k$.  For $x \in I$ and $c \in C$, we have 
\[
t_C([x, c]) = t_C(xc) - t_C(cx) = 0.
\]
That is, the image of $I$ in $C^-/k$ is a codimension 1 ideal which is not zero, so $C^-/k$ is not simple.

It remains to consider the case where $k$ has characteristic different from 2.
By Exercise \ref{skip.pr.simple}, we may assume that $k$ is algebraically closed and in particular that $C$ is split.  Let $I$ be a subspace of $C$ properly containing $k$ such that the image of $I$ in $C^-/k$ is an ideal; we aim to show that $I = C$.  Note that $I$ is an ideal in $C^-$.

\smallskip
We treat the case 
$C = \Mat_2(k)$.  Suppose first that $I$ contains $\stbtmat{-1}{0}{0}{1}$.  Since
\[
\left[ \stbtmat{a}{b}{c}{d}, \stbtmat{-1}{0}{0}{1} \right] = \stbtmat{-a}{b}{-c}{d} - \stbtmat{-a}{-b}{c}{d} = \stbtmat{0}{2b}{-2c}{0},
\]
we deduce that $I = C$.  If $I$ contains a trace zero element $x$ that is not nilpotent, then the Jordan form of $x$ is $\stbtmat{a}{0}{0}{-a}$ for some nonzero $a$, and the previous calculation shows that $I = C$.  

If $I$ contains a trace zero element $x$ that is nilpotent, then up to conjugacy $x$ is $\stbtmat{0}{1}{0}{0}$.  Since
\[
\left[ \stbtmat{a}{b}{c}{d}, \stbtmat{0}{1}{0}{0} \right] = \stbtmat{-c}{a-d}{0}{c},
\]
we deduce that $\stbtmat{-1}{0}{0}{1}$ is in $I$, whence $I = C$, completing the argument for the rank 4 case.

\smallskip
Now assume $C$ is an octonion algebra.  There is a nonzero element of $I \setminus k$, which by Exercise \dref{22.14}{pr.COVQUATALG} is contained in a quaternion subalgebra $B$ of $C$.  Because $I$ is an ideal in $C^-$, 
$I \cap B \supseteq [B, I \cap B]$, which is $B$ by the rank 4 case.

From the 4-dimensional subspace of $C$ perpendicular to $B$, pick a $j$ with $n_C(j) \ne 0$.  Then
$C$ is obtained by the internal Tits construction $C = \Cay(B, n_C(j))$ as in Prop.~\dref{20.1}{p.INTCD}.  For each $b$ in the space $B_0$ of trace zero elements of $B$, set $u = \frac12 b$, so $u - \bar{u} = b$.  In the notation of the external Tits construction \dref{20.2}{ss.EXTCD}, the calculations in \dref{20.12}{ss.CDCOMASS} give that $I$ contains $[u, j] = bj$ and therefore $I$ contains $B_0 j$.  
For $x, y \in B_0$, we have $[\frac12 x, yj] = (yx)j$, hence $(B_0 B_0) j \subseteq I$.  Since $B_0^2 = B$, we conclude that $I = C$ as desired.

(For the rank 8 case, one might compare p.~435 in \cite{Sagle1961}, which starts with an explicit multiplication table for a basis for $C$.)
\end{sol}

\begin{sol}{pr.sl.simple}
See for example \cite{Hogeweij}.
\end{sol}

\begin{sol}{pr.sp.char2}
\ref{sp.char2.der1}: It is straightforward to check that $q([A, B]x,x) = 0$ for all $A, B \in \mfsp_{2n}(F)$ and $x \in k^{2n}$, verifying $\subseteq$.  For the opposite containment, let us simplify notation by writing $L := [\mfsp_{2n}(F), \mfsp_{2n}(F)]$.  For any $B \in \Alt_{n}(F)$, pick $A \in \Mat_n(F)$ such that $B = A - A^\trans$ and note that
\[
\left[ \stbtmat{0}{\Eins_n}{0}{0}, \stbtmat{A}{0}{0}{-A^\trans} \right] = \stbtmat{0}{B}{0}{0} \quad \in L.
\]
An analogous computation shows that $\stbtmat{0}{0}{B}{0} \in L$.  For $A, A' \in \sl_{n}(F)$, we have
\[
\left[ \stbtmat{A}{0}{0}{-A^\trans}, \stbtmat{A'}{0}{0}{-(A')^\trans} \right] = \stbtmat{[A, A']}{0}{0}{-[A, A']^\trans}.
\]
By Exc.~\ref{pr.sl.simple}\ref{sl.simple.perfect} we conclude that 
\begin{equation} \label{sol.sp.char2.1}
\stbtmat{A}{0}{0}{-A^\trans} \in L
\end{equation}
for all $A \in \sl_{n}(F)$.  Take now $B \in \Sym_n(F)$ to have a 1 in the $(1,1)$-entry and zeros elsewhere.  Then
\[
\left[ \stbtmat{0}{B}{0}{0}, \stbtmat{0}{0}{B}{0} \right] = \stbtmat{B}{0}{0}{-B^\trans} \in L.
\]
In particular, matrices of the form \eqref{sol.sp.char2.1} include those with $A$ of trace zero and at least one of trace 1, ergo include all $A \in \Mat_n(F)$.

\ref{sp.char2.der2}: Put $L' := [L, L]$.  The same computations as in part \ref{sp.char2.der1} show that $L'$ contains the claimed subalgebra, i.e., those $\stbtmat{A}{B}{C}{-A^\trans}$ with $A \in \sl_n(F)$ and $B, C \in \Alt_n(F)$.  Conversely, if we take generic elements $\stbtmat{A}{B}{C}{-A^\trans}$, $\stbtmat{A'}{B'}{C'}{-(A')^\trans}$ of $L$, then the upper left corner of the commutator is
$[A, A'] + BC' - B'C$.  Therefore, we are done once we verify the following: \emph{If $B, C \in \Alt_n(F)$, then $\Tr(BC) = 0$.}

To verify that, note that $B = B_0 - B_0^\trans$ and $C = C_0 - C_0^\trans$ for some matrices $B_0, C_0$.  Plugging this in, we find
\begin{align*}
\Tr(BC) &= \Tr((B_0 - B_0^\trans)(C_0 - C_0^\trans)) \\
&= \Tr(B_0 C_0 + B_0^\trans C_0^\trans) - \Tr(B_0^\trans C_0 + B_0 C_0^\trans).
\end{align*}
Since $\Tr(XY) = \Tr(X^\trans Y^\trans)$ for all $X, Y \in \Mat_n(F)$ and $\car(F) = 2$, we conclude that $\Tr(BC) = 0$.
\end{sol}

\begin{sol}{pr.A2D2}
The action of $\phi \in \sl_{2n}(k)$ on $\wedge^d k^{2n}$ for any $d$ is given by
\begin{multline*}
\phi(x_1 \wedge \cdots \wedge x_d) = (\phi x_1) \wedge x_2 \wedge \cdots \wedge x_d + x_1 \wedge (\phi x_2) \wedge x_3 \wedge \cdots \wedge x_d \\
+ \cdots + x_1 \wedge \cdots \wedge x_{d-1} \wedge (\phi x_d).
\end{multline*}
Writing out the coordinates, we find that for a basis $e_1, \ldots, e_{2n}$ of $k^{2n}$, we find that the action of $\phi$ on $\wedge^{2n} k^{2n}$ is given by
\[
\phi(e_1 \wedge \cdots \wedge e_{2n}) = \Tr(\phi) \, (e_1 \wedge \cdots \wedge e_{2n}).
\]
Therefore, for $x, y \in \wedge^n k^{2n}$, we have
\begin{align*}
b(\rho(\phi) x, y) + b(x, \rho(\phi)y) &= \Delta(\rho(\phi)x \wedge y + x \wedge \rho(\phi)y) \\
&= \Delta(\phi(x \wedge y)) = \tr(\phi) \Delta(x \wedge y) = 0,
\end{align*}
so the image of $\phi$ is contained in $\mfg$, as required.  We omit the verification that $\rho$ is compatible with the Lie bracket.   Let $\phi \in \sl_{2n}(k)$ be such that $\phi(e_1) = e_{n+1}$ and $\phi(e_i) = 0$ for $i \ne 1$.  Then 
\[
\rho(\phi)(e_1 \wedge \cdots \wedge e_n) = e_{n+1} \wedge e_2 \wedge e_3 \wedge \cdots \wedge e_n \ne 0,
\]
i.e., $\rho(\phi) \ne 0$, completing the verification of \ref{pr.A2D2.hom}.

Note that $\rho(\Eins_{2n})$ is multiplication by $n$.  Therefore, if $n = 0$ in $k$, we have $\Eins_{2n} \subseteq \ker \rho$.

\ref{pr.A2D2.0}: Suppose now that $k$ is a field of characteristic not dividing $2n$.  Then $\sl_{2n}(k)$ is simple (Exc.~\ref{pr.sl.simple}), so $\ker \rho$ must be the zero ideal, i.e., $\rho$ is injective.  In the case $n = 2$, we note that $\sl_4(k)$ and $\mfo_6(k)$ both have dimension 15, so $\rho$ is an isomorphism.

\ref{pr.A2D2.2}: Finally assume $k$ is a field of characteristic 2.  Since $\rho$ is not the zero map and its kernel contains $\Eins_{2n}$, Exc.~\ref{pr.sl.simple} gives that $\ker(\rho) = k\Eins$.  Since $\Eins$ is the center of $\sl_{4}(k)$, $\rho$ induces an inclusion
\[
[\sl_4(k), \sl_4(k)] / k\Eins \hookrightarrow [\mfg, \mfg].
\]
But $[\sl_4(k), \sl_4(k)] = \sl_4(k)]$ by Exc.~\ref{pr.sl.simple}\ref{sl.simple.perfect}, so the domain is $\sl_4(k)/k\Eins$ of dimension 14.  The codomain has
$\dim [\mfg, \mfg] = \dim \mfg - 1$ by Exc.~\ref{pr.sp.char2}, which is $\binom{6}{2} - 1 = 14$.
\end{sol}

\solnsec{Section~\ref{s.DER}}

\begin{sol}{pr.LIEBC}
The linear maps $\Delta_A$ and $s_A$ commute with base change.  Recall the descriptions of ${\LMDer}(A)$, $\InDer(A)$ and $\ComDer(A)$ in (\ref{DEREXT}.\ref{LMDEREXT}), (\ref{ss.INDERALT}.\ref{INDEREXT}) and (\ref{ss.INDERALT}),
\begin{align*}
\LMDer(A)  &= \{\Delta_x\, |\, x \in W(A), s(x) \in {\rm Nuc}(A)\} = \Delta_A({s_A}^{-1}({\rm Nuc}(A))), \\
\InDer(A) &= \{\Delta_x\, |\, x \in W(A), s(x) = 0\} = \Delta_A({\rm Ker}(s_A)), \\
\ComDer(A) &= \{\Delta_a\, |\,  a \in A, s(a) = 0 \} = \Delta_A({\rm Ker}(s_A))\cap A).
\end{align*}
Using (\ref{ss.FLAFALG}.\ref{KERR}), (\ref{ss.FLAFALG}.\ref{IMF}), 
we obtain (a), (b) and (c). Finally $(D_{a, b})_R = D_{a_R,\, b_R}$ yields (d).
\end{sol}

\solnsec{Section~\ref{s.DEROCT}}

\begin{sol}{pr.LIEZOR}
If $C$ is a norm associative conic algebra then $t_C(cd) = t_C(dc)$ (\ref{ss.NORAS}.\ref{ADMASSBILT}) and $[C, C] \subseteq C^0$.

For the opposite inclusion, we are reduced by the usual argument to considering the case $C = \Zor(k)$.
Let $c = \left(\begin{matrix}
\alpha_1 & u_2 \\
u_1 & \alpha_2
\end{matrix}\right)$, $d = \left(\begin{matrix}
\beta_1 & v_2 \\
v_1 & \beta_2
\end{matrix}\right)\in C$. Using (\ref{ss.ZOVE}.\ref{ZOMU}), 
\[
[c,d] = \left(\begin{smallmatrix}
v_2^\trans u_1 - u_2^\trans v_1 & (\alpha_1 - \alpha_2)v_2 + (\beta_2 - \beta_1)u_2 + 2u_1 \times v_1 \\
(\beta_1 - \beta_2)u_1 + (\alpha_2 - \alpha_2)v_1 + 2u_2 \times v_2 & v_1^\trans u_2 - u_1^\trans v_2
\end{smallmatrix}\right)
\]
One sees easily that $[C, C] \supseteq C^0$. 
\end{sol}

\begin{sol}{pr.IDDERA}
We use the notation and results of \ref{ss.REDOCT}. Let $C = \Zor(k)$, $k$ a field of characteristic not 3, $e = e_1 = \left(\begin{matrix}1 &0 \\ 0  &0 \end{matrix} \right)$ and  $\Der(C) =
\mfg = \mfg_0 \oplus \mfg_1 \oplus \mfg_2$ the $e$-grading of $\Der(C)$ as in Prop. \ref{p.PEIRGRAD}. By Prop. \ref{p.GSL3}, the Lie subalgebra $\mfg_0$ of $\mfg$ is isomorphic to $\sl(M) := \{g \in \End_k(M)\mid  \tr(g) = 0\}$ via the map $\phi \: \sl(M) \rightarrow \mfg_0$ given by
\[
\phi(g)\left (
 \begin{matrix}
\alpha_1 & v^\ast \\
v        & \alpha_2
\end{matrix}
\right) := \left(
\begin{matrix}
0 & -g^\ast v^\ast \\
gv        & 0
\end{matrix}
\right) .
\]
That the two representations of $\sl_3(k)$ above are faithful and irreducible is clear. Since $k$ is a field, $M$ and $M^*$ are vector spaces of dimension 3 and we may write $M$ for both. Also $\langle x^*, y\rangle = x^{*\trans}y$ and if $g \in \End(M)$ then $g^* = g^\trans$.  
Assume there exists an isomorphism of  $\psi: M \rightarrow M$ of $\mfsl$-modules, i.e., $\psi(gv) = -g^\trans\psi(v)$. So we have $c^{-1}gc = -g^\trans$ for some invertible $c \in \Mat_3(K)$ and all $g \in\mfsl(M)$. Equivalently $gc = -cg^\trans$ for all $g \in \mfsl(M)$. Write $c = \sum_{i,j}\gamma_{ij}e_{ij}$ with $\gamma_{ij} \in k$ and $e_{ij}$ the usual matrix units in $\Mat_3(k)$. For distinct $p,q \in \{1,2,3\}$, we put $g := e_{pq} \in \mfsl_3(k)$ and obtain $g^\trans = e_{qp}$, hence
\[
gc = \sum_j\gamma_{qj}e_{pj}, \quad -cg^\trans = -\sum_i\gamma_{iq}e_{ip}.
\]
Comparing gives $\gamma_{qj} = \gamma_{iq} = 0$ for $j \neq p \ne i$, hence $c = 0$, a contradiction.
\end{sol}
 
\solnsec{Section~\ref{s.LIEALGJ}}

\noexsec

\solnsec{Section~\ref{s.DERALB}}

\begin{sol}{pr.VX0}
If $x = \alpha 1_J$, $\alpha \in k$ then $V_x y = 2\alpha y = 0$
 and $V_x = 0$. Let $x = \sum (\alpha_ie_i +  a_l[ij]) \in J$ and assume $V_x = 0$.  Then $V_x e_i = 2\alpha_i e_i + a_l[ij] + a_j[li] = 0$ and $a_l[ij] = 0 = a_j[li]$. Similarly $a_l[ij] = 0$ and $x = \sum \alpha_ie_i$. In that case $V_x 1_C[ij] = (\alpha_i + \alpha_j)1_C[ij] = 0$. Hence $\alpha_i = \alpha_j$. Similarly $\alpha_i = \alpha_l$ and $x = \alpha_i1_J$. 
\end{sol}

\begin{sol}{pr.DIHEDJ}
We need only verify non-zero products. Assume $\{i, j, l\} = \{1,2,3\}$. We have $J_0 \circ J_i \subseteq J_i$ by Th.\ref{t.PEDESI} (note conflicting notation) and $J_i \circ J_j \subseteq J_l$. Consider the triple  products. We use \ref{ss.EICJM}. If at least 2 of the entries are from $J_0$ the result of the product is in the submodule corresponding to the third entry. If exactly one of the entries is from $J_0$, the result of the product is in $J_0$ if the other two entries are from the same submodule, while if they are from two distinct submodules then the result belongs to the third. If two of the entries are from the same off-diagonal submodule then the result belongs to the submodule of the third entry.  Finally if each of the entries are from distinct off-diagonal submodules the result lands in $J_0$. So the triple products respect the grading. 
\end{sol}


\solnchap{Solutions for Chapter \ref{c.Groups}}

\solnsec{Section  \ref{skip.grp.back}}


\begin{sol}{skip.pr.SO.iso} 
The split hyperbolic plane $\bfh$ contains a hyperbolic pair $(u, v)$ as defined in \ref{ss.ISVE}.  For every $R \in \kalg$ and $t \in R^\times$, we have a linear transformation $g_t$ of the module underlying $\bfh \perp Q$ defined by setting $g_t(u) = tu$, $g_t(v) = t^{-1} v$, and $g_t$ to be the identity on the module underlying $Q$.
\end{sol}

\solnsec{Section~\ref{s.AUTGRP}}

\noexsec

\solnsec{Section  \ref{skip.co.sec}}

\begin{sol}{skip.pr.sets} \label{skip.sol.sets}The crux is to prove \ref{skip.pr.sets.1} implies \ref{skip.pr.sets.3}.  For sake of contradiction, suppose that the collection $\mcC$ is a set and $k$ is not the zero ring.  Then there is a field $K \in \kalg$ and we define 
\[
\mcD := \{ \dim (M \otimes K) \mid M \in \mcC \}
\]
where $\dim$ is the cardinality of a basis of the vector space $M \otimes K$.  (It is the same for all bases of $M \otimes K$ by \cite[{II.7.2}, Th.~3]{MR0354207}.)  Note that it is a set because $\mcC$ is a set.
 
Consider the case of $k$-modules.  For any cardinal $n$, we take $M$ to be the module $k^n$.  Since $M \otimes K \cong K^n$,  $\mcD$ contains every cardinal, contradicting that $\mcD$ is a set by \cite[{III.3.6}, Cor.]{MR2102219}.

Consider the case of $k$-algebras.  For each $n$, we consider the module $k^n$ as a nonassociative $k$-algebra with a multiplication that is identically zero and set $M \in \kalg$ to be the ``unital hull'' of $k^n$.  As a $k$-module, $M = k1_M \oplus k^n$, where $1_M$ is the identity of the multiplication on $M$.  Constructing the set $\mcD$ by the same recipe, we find that it in this case $\mcD$ contains $\dim (M \otimes K) = n + 1$.  Again we find a contradiction, so $\mcC$ is not a set.

For the rest of the implications, note that the only module under the zero ring is the module 0.
\end{sol}

\begin{sol}{skip.pr.closed}
Because $\sigma_i(M)$, $\tau_i(M)$ are finitely generated projective, by (\ref{skip.tensors}.\ref{skip.tensor.id}) we may view each $\alpha_i$ as an element of $\sigma_i(M)^* \otimes \tau_i(M)$.  Therefore, for every $R \in \kalg$, $\Aut(A_R)$ is the subgroup of $\GL(M_R)$ consisting of those $g$ stabilizing $\alpha_i$ for all $i$, and the claim follows from Exercise \ref{pr.skip.GLstab}.
\end{sol}

\begin{sol}{skip.pr.form.iso}
Write $i$ for the inclusion of $\bfGm$ in $\bfAut(f)$.
Trivially, in the weight space decomposition (\ref{skip.tori}.\ref{skip.wt.space}) of $V$ relative to $i(\bfGm)$, $V_\chi \ne 0$ for some nonzero $\chi$.  The composition $\chi \circ i$ is raising to a power $e \in \IZ$ where $e \ne 0$ because $\chi \ne 0$.  The statement about the weight space means that there is a nonzero $v \in V$ such that
for every $F \in \kalg$ and $t \in F^\times$, we have:
\[
f_F(v) = f_F(i(t)v) = f_F(t^e v) = t^{de} f_F(v).
\]
In particular, taking $F = k(\bft)$ and $t = \bft$, we find $(\bft^{de} - 1) f_{k(\bft)}(v \otimes 1) = 0$, ergo the second term in the product must be zero and $f_k(v) = 0$.
\end{sol}

\begin{sol}{pr.LANG}
Since $F$ is perfect, we may apply Prop.~\ref{skip.coh.gal}, reducing us to showing that the Galois cohomology set $H^1(K/F, \bfG)$ is zero for every finite extension $K \supseteq F$.

For such a $K$, the group $\Aut_{\Falg}(K)$ is cyclic of order $[K:F]$, generated by $\sigma$.  Following the notation in \ref{skip.gc}, a 1-cocycle is a map $g \: \Aut_{\Falg}(K) \to \bfG(K)$ satisfying (\ref{skip.gc}.\ref{skip.gc.1}).  By Lang, there is an $h \in \bfG(K)$ such that $g(\sigma) = h^{-1} \sigma(h)$.  By induction, we have
\[
g(\sigma^{i+1}) = g(\sigma) \cdot \sigma g(\sigma^i) = h^{-1} \sigma(h) \sigma(h^{-1}) \sigma^{i+1}(h) = h^{-1} \sigma^{i+1}(h),
\]
so the cocycle $g$ is equivalent to the trivial cocycle that maps all of $\Aut_{\Falg}(K)$ to $1 \in \bfG(K)$.  Since this holds for any 1-cocycle, we have shown that $H^1(K/F, \bfG) = 0$.
\end{sol}

\begin{sol}{skip.pr.finite}
Suppose $\bfX$ is a $\bfG$-torsor in the flat topology, i.e., a representative of a class in $H^1(k, \bfG)$.  Our aim is to show that $\bfX$ is the trivial torsor, i.e., $\bfX(k)$ is nonempty.  There is some fppf $R \in \kalg$ such that $\bfX(R)$ is nonempty.

If $k$ is reduced, then since it is artinian it is a finite product $k_1 \times \cdots \times k_m$ of finite fields $k_i$.  Write $R = \prod R_i$ where each $R_i$ is a $k_i$-algebra which is necessarily faithfully flat over $k_i$.  Then
$H^1(R/k, \bfG) = \prod_i H^1(R_i/k_i, \bfG_{k_i})$, where each term in the product is zero by Lang's Theorem (Exc.~\ref{pr.LANG}).

Drop the hypothesis that $k$ is reduced and put $\mfa := \Nil(k)$.  Because $k$ is finite, there is some minimal $m \ge 1$ such that $\mfa^m = 0$.  We proceed by induction on $m$, with the case $m = 1$ being settled by the previous paragraph, so suppose $m \ge 2$.  Put $I := \mfa^{m-1}$.  The ring $k/I$ has $\Nil(k/I)^{m-1} = (\Nil(k)/I)^{m-1} = 0$, so by induction $X(k/I)$ is nonempty.  On the other hand, $I^2 = \mfa^{2m-2} = \mfa^m \cdot \mfa^{m-2} = 0$ and $\bfX$ is smooth, so the natural map $\bfX(k) \to \bfX(k/I)$ is surjective.
\end{sol}

\begin{sol}{pr.const.tensor}
The first key point is an expression of elements of $\bfX_\Gamma(K)$ in terms of constant classes.  (This description has already appeared in the solutions to Exercises \ref{pr.ETAUT} and \ref{etale.aut} in the case where $\Gamma$ is a symmetric group.)
An element $\phi \in \bfX_\Gamma(K)$ is a homomorphism of rings $\IZ^\Gamma \to K$.  Under this homomorphism, the elements $\{ \phi(1_\gamma) \in K \mid \gamma \in \Gamma \}$ is a complete orthogonal system of idempotents in $K$.  This corresponds to an isomorphism $K \cong \prod_{\gamma \in \Gamma} K_\gamma$ for $K_\gamma \in \kalg$, some of which may be the zero ring.  The image of $\phi$ under the natural projection $\bfX_\Gamma(K) \to \bfX_\Gamma(K_\gamma)$ equals the image of 
\[
\gamma \in \Gamma = \bfX_\Gamma(\IZ) \to \bfX_\Gamma(K_\gamma).
\]

The second key point is that the description of a tensor system consists of terms that are $k$-linear, so that the isomorphism $K \cong \prod K_\gamma$ determined by $\phi$ induces isomorphisms $\bfAut(A_K) \cong \prod_{\gamma \in \Gamma} \bfAut(A_{K_\gamma})$.  The diagram
\[
\xymatrix{
\bfAut(A_K) \ar[r]^{f(K)} \ar[d]_{\cong} & \bfX_\Gamma(K) \ar[d]^{\cong} \\
\prod \bfAut(A_{K_\gamma}) \ar[r]^{\prod f(K_\gamma)} & \prod \bfX_\Gamma(K_\gamma)}
\]
commutes because of functoriality of $f$.
The image of $\phi$ under the composition $\bfX_\Gamma(K) \to \bfX_\Gamma(K_\gamma)$ is the constant class that is the image of $\gamma$, so by hypothesis it is in the image of $f(K_\gamma)$.  That is, the element $\phi \in \bfX_\Gamma(K)$ maps down to something in the image of the bottom arrow.  Commutativity of the diagram proves the claim.
\end{sol}

\solnsec{Section \ref{skip.app.sec}}

\begin{sol}{skip.pr.frddense}
By Corollary \ref{skip.Frd.3}, $J = E^{(+)}$ for some cubic \'etale $k$-algebra $E$.  Therefore, the claim is equivalent to the same statement for $E$, which is Exercise \ref{skip.pr.etale.gen}. 
\end{sol}

\begin{sol}{skip.unitary9}
The proof proceeds in the same manner as the proofs of Propositions \ref{skip.octonion.G2} and \ref{skip.Frd.A1}.  Note that for each central simple associative $F$-algebra $B$, we get an algebra with unitary involution $A = B \times \opgrp{B}$, $K = F \times F$, and $\tau$ the map that swaps the two factors, and in this case we abuse notation and write simply $B$ for the involutorial system.  The automorphism group of the algebra $\Mat_3(F)$ is the same as the automorphism group of the corresponding algebra with unitary involution, which is $\PGL_3 \rtimes \IZ/2$ as in \cite[pp.~346, 400]{MR2000a:16031}.  The automorphism group of $\Sym(\Mat_3(F))$ is the same by Proposition \ref{skip.frd.69}.\ref{skip.frd.9}, and the automorphism group of the split adjoint semisimple group of type $\rsA_2$, $\PGL_3$, is the same as described in \ref{skip.split.ring}.

The $F$-algebras with involution that are twisted forms of the algebra with unitary involution corresponding to $\Mat_3(F)$ are exactly the ones described in (i) by \cite[\S{29.D}]{MR2000a:16031}.  The rank 9 Freudenthal $F$-algebras are the $F$-forms of $\plalg{\Mat_3(F)}$ by Corollary \ref{c.CHARF}.  The adjoint semisimple groups of type $\rsA_2$ are by definition the $F$-forms of $\PGL_3$.
\end{sol}

\begin{sol}{skip.Frd.C3}
The solution is similar to that of Exc.~\ref{skip.unitary9}.

The automorphism group of $(\Mat_6(F), \tau_{\spl})$ is $\PGSp_6$ by \cite[pp.~347, 359]{MR2000a:16031}.  We already asserted that the split Freudenthal algebra of rank 15 had automorphism group $\PGSp_6$, see Theorem \ref{skip.aut.alb}.  Finally, the automorphism group of $\PGSp_6$ is itself $\PGSp_6$ by Lemma \ref{skip.ss.ff}.  

The $F$-algebras with involution $(A, \sigma)$ in \ref{skip.Frd.C3.1} are twisted forms of $(\Mat_6(F), \tau_\spl)$, because $A \otimes \aclosF \cong \Mat_6(\aclosF)$ and all non-degenerate skew-symmetric bilinear forms on $\aclosF^6$ are equivalent, so all symplectic involutions on $\Mat_6(\aclosF)$ are conjugate to $\tau_\spl$.  Conversely, if $(A, \sigma)$ is a twisted form of $(\Mat_6(F), \tau_\spl)$, then $A$ is a central simple $F$-algebra (Corollary \ref{ss.c-CENTSIM4}) and the involution $\sigma$ is evidently also symplectic.
\end{sol}

\solnsec{Section \ref{s.ALBINT}}

\begin{sol}{pr.CHARPOS} \label{sol.CHARPOS}  (a) Being euclidean, $J$ does not contain nilpotent elements other than zero. By Exc.~\ref{pr.PRIMID}~(c), therefore, the unital commutative associative subalgebra $\IR[x] \subseteq J$, which has dimension at most $3$ by Thm.~\ref{t.CUNOJO}, is a direct product of finitely many finite algebraic field extensions of $\IR$, i.e., of copies of $\IR$ and $\IC$. Copies of $\IC$ cannot occur since $1^2 + i^2 = 0$ and hence $\IC^{(+)}$ is not a \emph{euclidean} Jordan algebra. Thus, for some $n = 1,2,3$, there is a complete orthogonal system $(d_1,\dots,d_n)$ of non-zero idempotents in $J$ such that $\IR[x] = \IR d_1 \oplus \cdots \oplus \IR d_n$ as a direct sum of ideals. By Exc.~\ref{pr.IDCUJO}, each $d_i$, $1 \leq i \leq n$, is either elementary or co-elementary or equal to $1_J$ and hence, by Exc.~\ref{pr.RANONE}~(f), splits into the orthogonal sum of elementary idempotents. By \ref{ss.COELF}, therefore, an elementary frame $(c_1,c_2,c_3)$ of $J$ and scalars $\alpha_1,\alpha_2,\alpha_3$ exist such that \eqref{EXSUCE} holds. 

\smallskip

(b) For $x \in J$ as given in (a), we recall from (\ref{e.CUBET}.\ref{NXXI})--(\ref{e.CUBET}.\ref{QUATXXI}) that
\begin{align}
\label{NOEXAL} N(x) = \alpha_1\alpha_2\alpha_3, \quad T(x) = \alpha_1 + \alpha_2 + \alpha_3, \quad S(x) = \alpha_2\alpha_3 + \alpha_3\alpha_1 + \alpha_1\alpha_2.
\end{align}
(i) $\Rightarrow$ (ii). Pick an elementary frame $(c_1,c_2,c_3)$ in $J$ and positive real numbers $\alpha_1,\alpha_2,\alpha_3$ such that \eqref{EXSUCE} holds. Then $y := \sum\sqrt{\alpha_i}c_i \in J$ satisfies $N(y) = \sqrt{\alpha_1\alpha_2\alpha_3} \in \IR^\times$ by \eqref{NOEXAL}, hence $y \in J^\times$, and $y^2 = x$. 

(ii) $\Rightarrow$ (iii). Choose $y$ as in (ii). By (a) and \eqref{NOEXAL}, there exist an elementary frame $(c_1,c_2,c_3)$ of $J$ and $\beta_1,\beta_2,\beta_3 \in \IR^\times$ such that $y = \sum\beta_ic_i$, hence $x = \sum\alpha_ic_i$, $\alpha_i = \beta_i^2 > 0$ for $i = 1,2,3$. Since $T$ is an associative bilinear form, $L_x$ is self-adjoint relative to $T$, and it suffices to prove \eqref{TEXBUZ}. Before doing so, we show
\begin{align}
\label{TECEDE} T(c,d) \geq 0
\end{align}
for all idempotents $c,d \in J$. By Exc.~\ref{pr.IDCUJO}, there are four cases: $c = 0$ (trivial), $c$ elementary, $c$ co-elementary, and $c = 1_J$. The same distinction can be made about $d$; hence $T(d) \in \{0,1,2,3\}$, and the case $c = 1_J$ is obvious. Now let $c$ be elementary. Then the Peirce decomposition of $d$ with respect to $c$ attains the form $d = T(c,d)c + d_1 + d_0$, $d_i \in J_i(c)$, $i = 0,1$. The Peirce-$2$-component of $d^2$ with respect to $c$ may therefore be written as $\alpha c$,
\begin{align*}
\alpha =\,\,&T(c,d^2) = T(c,d)^2T(c,c) + T(c,d_1^2) = T(c,d)^2 + T(c \bu d_1,d_1) \\
=\,\,&T(c,d)^2 + \frac{1}{2}T(d_1,d_1) \geq T(c,d)^2
\end{align*}
since $T$ is positive definite. But $d^2 = d$ and we conclude $T(c,d) \geq T(c,d)^2$, i.e.,
\begin{align}
\label{ZETECE} 0 \leq T(c,d) \leq 1
\end{align}
if $c$ is elementary. Finally, let $c$ be co-elementary. Then $\pc := 1_J - c$ is elementary, and since \eqref{TECEDE} is trivial for $d = 0$, we may assume $d \neq 0$, in which case \eqref{ZETECE} yields $T(c,d) = T(d) - T(\pc,d) \geq T(d) - 1 \geq 0$. This completes the proof of \eqref{TECEDE}. Returning to \eqref{TEXBUZ}, we now write $0 \neq z \in J$ as $z =\sum\gamma_id_i$ for some elementary frame $(d_1,d_2,d_3)$ in $J$ and $\gamma_i \in \IR$, $1 \leq i \leq 3$, not all zero. From \eqref{TECEDE} we deduce $T(c_i,z^2) = \sum_m\gamma_m^2T(c_i,d_m) \geq 0$ for $1 \leq i \leq 3$, and setting $\alpha := \min_{1\leq i\leq 3}\alpha_i > 0$, we obtain
\begin{align*}
T(x \bu z,z) =\,\,&T(x,z^2) = \sum\alpha_iT(c_i,z^2) \geq \alpha\sum T(c_i,z^2) = \alpha T(z^2) = \alpha T(z,z) > 0, 
\end{align*}
as claimed.

(iii) $\Rightarrow$ (i). Writing $x$ as in (a), we obtain $\alpha_i = T(x \bu c_i,c_i) > 0$ for $1 \leq i \leq 3$, hence (i).

(iii) $\Rightarrow$ (iv). Being positive definite, the linear operator $L_x$ is bijective, and we see that $x$ is linearly invertible in the sense of Exc.~\ref{pr.LININV}, hence invertible. Thus
\begin{align}
\label{PDLEMU} X := \{\px \in J \mid L_{\px}\,\text{is positive definite relative to $T$}\} \subseteq J^\times.
\end{align}
Clearly, $X$ is convex, hence connected with $1_J \in X$, and we have shown
\begin{align}
\label{EXINPO} X \subseteq \Pos(J).
\end{align}
In particular, (iv) holds.

Before tackling to the final implication (iv) $\Rightarrow$ (iii) of (b), we prove:

(c) (i) $\Rightarrow$ (ii) $\Rightarrow$ (iii) $\Rightarrow$ (i). This follows from the corresponding implications in (b) by noting that (c) (i) (resp. (c) (ii), (c) (iii)) holds for $x$ if and only if (b) (i) (resp. (b) (ii), (b) (iii)) holds for for $x + \vep 1_J$ and all $\vep > 0$. By the same token, the solution to the entire exercise will be complete once we have established the implication 

\smallskip

(b) (iv) $\Rightarrow$ (iii). By homogeneity, \eqref{TEXBUZ} holds if and only if
\begin{align}
\label{TEXBUS} T(x \bu z,z) > 0 &&(z \in S),
\end{align}
where $S$ stands for the ``sphere''
\[
S := \{z \in J \mid T(z,z) = 1\}.
\]
Since $S$ is compact, \eqref{TEXBUS} and a standard continuity argument show that $X$ as defined in \eqref{PDLEMU} is open in $\Pos(J)$. On the other hand, $\Pos(J)$ is connected, so if $X$ were not all of $\Pos(J)$, it could not be closed in $\Pos(J)$, i.e., we would have $X \subset \bar X \cap \Pos(J)$. Write $\px \in (\bar X \cap \Pos(J)) \setminus X$ as in (a). The set
\[
Y := \{y \in J \mid L_y\,\text{is positive semi-definite relative to $T$}\}
\]
is closed in $J$ and contains $X$, hence $\bar X$, hence $\px$. The implication (iii) $\Rightarrow$ (i) in (c) therefore implies $\alpha_i \geq 0$ for $1 \leq i \leq 3$, but also $\alpha_i = 0$ for some $i = 1,2,3$ since, otherwise, $\px$ would belong to $X$. But then, by \eqref{NOEXAL}, $\px$ cannot be invertible, hence does not belong to $\Pos(J)$, a contradiction. Thus $X = \Pos(J)$, implying (iii), and the solution is complete.
\end{sol}

\begin{sol}{pr.PROPOS} \label{sol.PROPOS} 
In each of the four cases, it suffices for obvious reasons to establish only the first of the two claimed relations. 

(a) We put $Y := \Pos(J)$. Since $U_p\:J^\times \to J^\times$ is a topological map, it permutes the connected components of $J^\times$. In particular, $U_pY \subseteq J^\times$ is a connected component containing $U_p1_J = p^2$. But by Exc.~\ref{pr.CHARPOS}~(b), this element is also contained in $Y$, which implies $U_pY = Y$.

(b) By Exc.~\ref{pr.CHARPOS}~(b) we have $p = y^2 = U_y1_J$ for some $y \in J^\times$, hence $J^{(p)} \cong J$ by Cor.~\ref{c.ORSTR}. From (a) we deduce $U_p^{-1}Y = U_{p^{-1}}Y = Y$, whence $p^{-1} = U_p^{-1}p \in Y$. By definition, $\Pos(J^{(p)})$ is the connected component of $J^{(p)\times} = J^\times$ containing $1_{J^{(p)}} = p^{-1} \in Y$. Thus $\Pos(J^{(p)}) = Y$. 

(c) Assume first $\eta \in \Aut(J)$. Again, $\eta\:J^\times \to J^\times$ is a topological map, forcing $\eta(Y) \subseteq J^\times$ to be a connected component containg $\eta(1_J) = 1_J$. Thus $\eta(Y) = Y$.  Since $\Str^0(J)$ acts canonically on $J^\times$, the orbit $\Str^0(J)p \in J^\times$ of $p \in Y$ is a connected subset containing $p$ and hence belonging to $Y$. Thus $\eta(Y) \subseteq Y$ for $\eta \in \Str^0(J)$, and applying this also to $\eta^{-1}$ in place of $\eta$ gives the assertion. 

(d) $\Pos(\pJ) \subseteq {\pJ}^\times \subseteq J^\times$ is a connected subset containing $1_{\pJ} = 1_J$. Hence $\Pos(J^\prime) \subseteq \Pos(J) \cap {\pJ}^\times \subseteq \Pos(J) \cap \pJ$. Conversely, suppose $y \in \Pos(J) \cap \pJ$. Then $y \in \pJ$, being invertible in $J$, is so in $\pJ$ since we are dealing with \emph{finite-dimensional} Jordan algebras. Thus $y \in \Pos(J) \cap {\pJ}^\times$. Now Exc.~\ref{pr.CHARPOS}, \eqref{TEXBUZ} shows $T(y \bu z,z) > 0$ for all non-zero elements $z \in J$. Since the bilinear trace of $J$ restricts to the bilinear trace of $\pJ$, we conclude from $T(y \bu z,z) > 0$ for all non-zero elements $z \in \pJ$ that $y \in \Pos(\pJ)$. This proves the first two equations of (d). As to the last, $\bPos(J) \cap \pJ$ is a closed subset of $\pJ$ containing $\Pos(J) \cap \pJ = \Pos(\pJ)$. Thus $\bPos(\pJ) \subseteq \bPos(J) \cap \pJ$. Conversely, let $y \in \bPos(J) \cap \pJ$. Then Exc.~\ref{pr.CHARPOS}, \eqref{TEXBU}) shows $T(y \bu z,z) \geq 0$ for all $z \in J$, hence, in particular, for all $z \in \pJ$. Thus $y \in \bPos(J^\prime)$.
\end{sol}

\begin{sol}{pr.POSDEF} \label{sol.POSDEF}  (a) Let $y \in \GL_n(\ID)$ and $\eta := \Phi_y$. For $x,z \in J$, we obtain
\[
U_{\eta(x)}z = U_{\bar y^\trans xy}z = \bar y^\trans xyz\bar y^\trans xy = \eta U_x\eta^\sharp z,
\]
with $\eta^\sharp := \Phi_{\bar y^\trans }$. Hence (a) follows from (\ref{ss.STRUG}.\ref{STREL}). 

(b) We have $\Phi_{y_1y_2} = \Phi_{y_2} \circ \Phi_{y_1}$ for all $y_1,y_2 \in \GL_n(\ID)$, so (b) will follow once we have shown $\Phi_y \in \Str^0(J)$ for all $y \in \Tri_n(\ID)$. But $\Tri_n(\ID)$, being homeomorphic to $\ID^{\frac{n(n-1)}{2}}$ as a topological space, is connected, and the assertion follows from the fact that $\Phi$ is continuous. 

(c) (i) $\Rightarrow$ (ii). The case $n = 1$ being trivial, let us assume $n > 1$ and that the implication holds for $n - 1$ in place of $n$. Write
\begin{align}
\label{EXEXONE} x = \left(\begin{matrix}
x_1 & u \\
\bar u^\trans  & \xi 
\end{matrix}\right), \quad x_1 \in J_1 := \Her_{n-1}(\ID),\;u \in \Mat_{n-1,1}(\ID),\;\xi \in \IR.
\end{align}
Since $x_1$ along with $x$ is positive definite, hence invertible, \eqref{BLOMA} yields $y \in \Tri_n(\ID)$ such that
\[
\px := \bar y^\trans xy = \left(\begin{matrix}
x_1 & 0 \\
0 & \xi_n 
\end{matrix}\right) \in J.
\]
Since $\px$ is positive definite as well, we obtain $\xi_n > 0$, and the induction hypothesis leads to an element $y_1 \in \Tri_{n-1}(\ID)$ making $\bar y_1^\trans x_1y_1$ a diagonal matrix with positive diagonal entries. This property carries over to $\overline{y\py}^\trans xy\py = \overline{\py}^\trans \bar y^\trans  xy\py$ with $\py = \left(\begin{smallmatrix}
y_1 & 0 \\
0 & 1 
\end{smallmatrix}\right)$, and the induction is complete.

(ii) $\Rightarrow$ (i). By (ii) we find a matrix $z \in J^\times$ satisfying $ \px := \bar y^\trans xy = z^2 = \bar z^\trans z$. Hence $\px$ is positive definite and thus so is $x$.

Now let $n = 3$. (ii) $\Rightarrow$ (iii). Again some $z \in J^\times $ has $\px := \bar y^\trans xy = z^2$, which by Exc.~\ref{pr.CHARPOS}~(b) belongs to $\Pos(J)$. Hence so does $x = \Phi_{y^{-1}}\px$ by (b) and Exc.~\ref{pr.PROPOS}~(c).

(iii) $\Rightarrow$ (i). $x \in \Pos(J)$ leads to an element $z \in J^\times$ such that $x = z^2 = \bar z^\trans z$. Hence $x$ is positive definite.
\end{sol}

\begin{sol}{pr.MINPOS} \label{sol.MINPOS} 
(a) (i) $\Rightarrow$ (ii). From Exc.~\ref{pr.CHARPOS}~(b) we deduce $T(x,e_{ii}) = T(x,e_{ii}^2) = T(x \bu e_{ii},e_{ii}) > 0$ for $1 \leq i \leq 3$. The same exercise implies $x = y^2$ for some $y \in J^\times$, hence $x^\sharp = y^{2\sharp} = y^{\sharp 2} \in \Pos(J)$, and we also have $T(x^\sharp,e_{ii}) > 0$. Finally, $N(x) = N(y)^2 > 0$. 

(ii) $\Rightarrow$ (iii). Obvious.

(iii) $\Rightarrow$ (i).  Note first that the unital subalgebras of $\IO$ up to conjugation under $\Aut(\IO)$ are $\IR$, $\IC$, $\IH$ and $\IO$ itself (Exc.~\ref{pr.SKONOCO}). 

Next we claim that \emph{there is no harm in assuming $u_1 = 0$ or $u_1 = 1_{\IO}$.} Indeed, if $u_1 \neq 0$, we replace $x$ by $\alpha^{-1}x$, $\alpha := \sqrt{n_{\IO}(u_1)}$, if necessary to ensure $n_{\IO}(u_1) = 1$. By Exc.~\ref{pr.AURE}, the assignment
\[
 \sum(\eta_ie_{ii} + v_i[jl]) \longmapsto \sum\eta_ie_{ii} + (u_1^{-1}v_1)[23] + (v_2u_1^{-1})[31] + (u_1v_3u_1)[12]
\]
defines an automorphism $\vph$ of $J$ that leaves all minors of $x$ as well as $\Pos(J)$ (Exc.~\ref{pr.PROPOS}~(c)) invariant and sends $x$ to 
\[
\sum\xi_ie_{ii} + 1_{\IO}[23] + (u_2u_1^{-1})[31] + (u_1u_3u_1)[12],
\]
proving our claim.

From now on we therefore assume $u_1 \in \{0,1_{\IO}\}$. Then Exc.~\ref{pr.ARTCONALG} combined with our initial observation implies that there is no harm in assuming $x \in \pJ:= \Her_3(\IH)$. Write
\[
x = \left(\begin{matrix}
x_1 & u \\
\bar u^\trans  & \xi_3 
\end{matrix}\right), \quad x_1 \in \Her_2(\IH), \quad u \in \Mat_{2,1}(\IH), \quad \xi_3 \in \IR.
\]
More precisely, we have $x_1 \in \ppJ := \Her_2(\IR[u_3])$, and $\IR[u_3]$ is isomorphic to $\IR$ or $\IC$, hence commtative associative. Since the principal minors of $x_1$ are positive by (iii), we deduce from Bourbaki \cite[IX.7, Prop.~3]{MR0107661} that $x_1$ is positive definite, hence invertible. Invoking \eqref{BLOMA} of Exc.~\ref{pr.POSDEF}, we therefore find an element $y \in \Tri_3(\IH)$ such that
\begin{align}
\label{EXITH} \px := \bar y^\trans xy = \left(\begin{matrix}
x_1 & 0 \\0 & \pxi_3 
\end{matrix}\right) \in \ppJ 
\end{align}
for some $\pxi_3 \in \IR$. Since $\Phi_y$ by Exc.~\ref{pr.POSDEF}~(a) belongs to the structure group of $\pJ$, and $\pJ$ has degree $3$, it reproduces the norm up to the factor $N(\bar y^\trans y)$, which must be positive since $\Tri_3(\IH)$ is connected  and $y \mapsto N(\bar y^\trans y)$ is continuous. This proves $N(\px) > 0$, and we have shown that all principal minors of $\px$ are positive, forcing $\px$ itself to be positive definite (Bourbaki, loc. cit.). From Exercises~\ref{pr.POSDEF}~(b) and \ref{pr.PROPOS}~(d), we therefore deduce $\px \in \Pos(\ppJ) \subseteq \Pos(\pJ)$, hence $x \in \Pos(\pJ) \subseteq \Pos(J)$. 

(b) Apply (a) to $x + \vep \Eins_3$, $\vep > 0$, and let $\vep \to 0$.
\end{sol}

\solnsec{Section~\ref{s.E6}}

\noexsec

\addtocontents{toc}{\vspace{\baselineskip}}

\bibliographystyle{amsalpha}
\bibliography{b-albalg}
\printindex


\end{document}